\documentclass[reqno,12pt]{article}

\usepackage{amsmath,amssymb,amsthm,mathtools,bm,mathrsfs}
\usepackage{geometry}
\usepackage{microtype}
\usepackage{enumitem}
\usepackage{booktabs,array,aliascnt,fancyhdr}
\usepackage[numbers,sort&compress]{natbib}
\usepackage[hidelinks]{hyperref}
\usepackage[nameinlink,noabbrev]{cleveref}
\allowdisplaybreaks[3]
\numberwithin{equation}{section}
\setlist{itemsep=0.2em,topsep=0.35em}
\makeatletter
\renewcommand\section{\@startsection{section}{1}{\z@}%
  {-2.2ex \@plus -0.5ex \@minus -0.2ex}%
  {1.1ex \@plus 0.2ex}%
  {\normalfont\large\bfseries\raggedright}}
\makeatother

\newcommand{\RunningTitle}{Viscous Contact Waves under Large Perturbations}
\theoremstyle{plain}
\newtheorem{theorem}{Theorem}[section]
\newtheorem{lemma}{Lemma}[section]
\newaliascnt{proposition}{lemma}
\newtheorem{proposition}[proposition]{Proposition}
\aliascntresetthe{proposition}
\newaliascnt{corollary}{lemma}
\newtheorem{corollary}[corollary]{Corollary}
\aliascntresetthe{corollary}
\theoremstyle{remark}
\newaliascnt{remark}{lemma}
\newtheorem{remark}[remark]{Remark}
\aliascntresetthe{remark}
\crefname{theorem}{Theorem}{Theorems}
\Crefname{theorem}{Theorem}{Theorems}
\crefname{lemma}{Lemma}{Lemmas}
\Crefname{lemma}{Lemma}{Lemmas}
\crefname{proposition}{Proposition}{Propositions}
\Crefname{proposition}{Proposition}{Propositions}
\crefname{corollary}{Corollary}{Corollaries}
\Crefname{corollary}{Corollary}{Corollaries}
\crefname{remark}{Remark}{Remarks}
\Crefname{remark}{Remark}{Remarks}
\crefname{equation}{equation}{equations}
\Crefname{equation}{Equation}{Equations}
\crefname{section}{Section}{Sections}
\Crefname{section}{Section}{Sections}
\crefname{appendix}{Appendix}{Appendices}
\Crefname{appendix}{Appendix}{Appendices}

\newcommand{\dd}{\,\mathrm d}
\newcommand{\R}{\mathbb R}
\newcommand{\pa}{\partial}

\title{\bfseries Stability and Optimal Decay of Viscous Contact Waves under Large Perturbations\thanks{This work is partially supported by the National Natural Science Foundation of
China (12471227, 12494541).}}
\author{Jing Li$^a$; Ziyin Liu$^b$
\thanks{a. Institute of Applied Mathematics \& State Key Laboratory of Mathematical Sciences, AMSS,\\ \& Hua Loo-Keng Key Laboratory of Mathematics,
Chinese Academy of Sciences,\\ Beijing, 100190 People's Republic of China, \\ and School of Mathematics and Computer Sciences, \\ Institute of Mathematics and Interdisciplinary Sciences,\\ Nanchang
University, Nanchang, 330031, People's Republic of China (ajingli@gmail.com).\\
b. Institute of Applied Mathematics, Academy of Mathematics and Systems Science,\\
Chinese Academy of Sciences, Beijing 100190, People's Republic of China\\
(a1909462594@gmail.com).}}
\date{}
\hypersetup{
  pdftitle={Stability and Optimal Decay of Viscous Contact Waves under Large Perturbations},
  pdfauthor={Jing Li; Ziyin Liu},
  pdfkeywords={viscous contact wave; large perturbation; optimal decay; one dimension}
}

\begin{document}
\maketitle
\begin{abstract}
The stability of viscous contact waves is established for the
one-dimensional full compressible Navier--Stokes equations under large
integrable perturbations. A unique global positive solution is obtained,
which converges uniformly to a fixed Eulerian contact profile at the
optimal inverse square-root rate in time. The wave-strength threshold
depends only on bounds
for the initial density perturbation in the first-order Sobolev norm,
the energy of the initial velocity perturbation, the thermal relative
entropy and the total conserved excess, together with positive lower
bounds for the initial density and temperature. It is independent of
higher Sobolev norms and the common
far-field velocity. No square-integrability condition is imposed on the
initial anti-derivatives. Derivative and pressure decomposition
estimates are also established for the evolving solution.
\end{abstract}

\section{Introduction}\label{sec:intro}

We study the one-dimensional full compressible Navier--Stokes equations
in Eulerian coordinates $(x,t)\in\R\times[0,\infty)$.
The density $\rho$, velocity $u$ and absolute temperature $\theta$ satisfy
\begin{equation}\label{eq:bg-eqs}
\begin{cases}
 \rho_t+u\rho_x+\rho u_x=0,\\
 \rho(u_t+uu_x)+p_x=\nu u_{xx},\\
 c_v\rho(\theta_t+u\theta_x)+pu_x=\kappa\theta_{xx}+\nu u_x^2.
\end{cases}
\end{equation}
The viscosity coefficient $\nu$
and heat conductivity $\kappa$ are constants satisfying
\begin{equation}\label{eq:cf-hyp}
 \nu>0,\qquad \kappa>0.
\end{equation}
The gas is polytropic and ideal, with
\begin{equation}\label{eq:state}
 p=R\rho\theta,\qquad c_v=\frac R{\gamma-1},\qquad
 R>0,\quad\gamma>1.
\end{equation}
We write $c_p=\gamma c_v$ and $c=\kappa/c_p$.

The positive end states satisfy
\begin{equation}\label{eq:ct-st}
 u_-=u_+=:u_*,\qquad R\rho_-\theta_-=R\rho_+\theta_+,\qquad
 \theta_-\ne\theta_+.
\end{equation}
The corresponding Euler equations have the contact discontinuity
\begin{equation}\label{eq:ct}
 (\rho,u,\theta)(x,t)=
 \begin{cases}
 (\rho_-,u_*,\theta_-),&x<u_*t,\\
 (\rho_+,u_*,\theta_+),&x>u_*t,
 \end{cases}
\end{equation}
across which the velocity and pressure are continuous. We put
\begin{equation}\label{e:end}
 p_*=R\rho_\pm\theta_\pm,\qquad
 \delta=|\theta_+-\theta_-|>0,\qquad \tau=1+t,
\end{equation}
and call $(R,\gamma,\nu,\kappa,\rho_-,\theta_-)$ the fixed data.

Small-perturbation stability and convergence rates for viscous contact
waves were established in
\cite{HuangMatsumuraXin2006,HuangXinYang2008,HuangLiMatsumura2010};
the optimal uniform diffusion rate for nonzero excess mass was obtained
in \cite{LiuQiuWangXu2025,LiuWangXu2026}.
Hong~\cite{Hong2015} treated a class of large perturbations by choosing
the smoothing width sufficiently large. Huang and
Wang~\cite{HuangWang2016} proved stability for perturbations of any
prescribed $H^1$ size, with qualitative convergence and a wave-strength
threshold depending on that size. Peng, Shi and
Wu~\cite{PengShiWu2026} removed the initial velocity and temperature
gradient norms from this dependence. The result below gives the
optimal uniform diffusion rate for large integrable perturbations of a
contact profile of fixed initial width, with a threshold controlled by
the lower-order quantities in \eqref{e:dat}.

\subsection{Contact wave, initial data and conserved excess}

Following \cite{HuangMatsumuraXin2006}, the viscous contact wave is
\begin{equation}\label{eq:ct-w}
 \bar Q=(\bar\rho,\bar u,\bar\theta)
 =\left(\frac{p_*}{R\bar\theta},\,
 u_*+\frac{\kappa R}{c_pp_*}\partial_x\bar\theta,\,\bar\theta\right).
\end{equation}
Here $\bar\theta$ is the positive monotone diffusion profile satisfying
\begin{equation}\label{eq:ct-euler}
\begin{aligned}
 (\partial_t+u_*\partial_x)\bar\theta
 &=\frac{\kappa R}{c_pp_*}
   \bigl(\bar\theta\partial_x^2\bar\theta-(\partial_x\bar\theta)^2\bigr),\\
 \bar\theta(\pm\infty,t)&=\theta_\pm .
\end{aligned}
\end{equation}
Its initial profile and spatial phase are fixed by the standard
construction \eqref{eq:prof}, \eqref{e:Xc} in
Section~\ref{s:bg-coordinates}. In particular, $R\bar\rho\bar\theta=p_*$.

For $e_0=(e_\rho,e_u,e_\theta)$, set
\begin{equation}\label{eq:bg-dat}
 W_0=(\rho_0,u_0,\theta_0)=\bar Q(\cdot,0)+e_0.
\end{equation}
The initial and far-field conditions are
\begin{equation}\label{eq:data}
 (\rho,u,\theta)(x,0)=W_0(x),
\end{equation}
\begin{equation}\label{eq:far}
 (\rho,u,\theta)(x,t)\longrightarrow
 (\rho_\pm,u_*,\theta_\pm)\quad\text{as }x\to\pm\infty.
\end{equation}
For positive $\rho_0$, the initial mass labels are
\begin{equation}\label{e:Y0}
 Y_0(x)=\int_0^{x}\rho_0(z)\dd z,\qquad X_0=Y_0^{-1}.
\end{equation}
Let $(V_c,U_c,\Theta_c)$ be the mass-coordinate profile in
\eqref{eq:prof}, and write $H_c=U_c-u_*$. With
\[
 (V_0,U_0,\Theta_0)
 =\left(\frac1{\rho_0},u_0,\theta_0\right)\circ X_0,
\]
define the initial conserved excess $Q=(Q_v,Q_{\rm mom},Q_E)$ by
\begin{equation}\label{e:Q-intro}
\begin{aligned}
 Q_v&=\int_\R[V_0(y)-V_c(y,0)]\dd y,\\
 Q_{\rm mom}&=\int_\R[U_0(y)-u_*-H_c(y,0)]\dd y,\\
 Q_E&=\int_\R\left(c_v[\Theta_0(y)-\Theta_c(y,0)]
                      +\frac{[U_0(y)-u_*]^2}{2}\right)\dd y.
\end{aligned}
\end{equation}
Lemma~\ref{l:cmp} proves conservation of these three integrals;
\eqref{e:mass} gives the resulting contact phase and acoustic masses.

\subsection{Main result}

Our main result is stated as follows.

\begin{theorem}[One-dimensional stability and decay]\label{thm:planar}
Let $R>0$ and $\gamma>1$, let $\nu,\kappa$ satisfy
\eqref{eq:cf-hyp}, and let $\rho_-,\theta_->0$ and $A\ge1$.
Let $Q$ be the initial conserved excess vector defined in
\eqref{e:Q-intro}, and put $\Phi(z)=z-\log z-1$.
There is $\delta_0(A)>0$, depending only on $A$ and the fixed data and
independent of $u_*$, such that the following holds for every $u_*\in\R$, every right
state $(\rho_+,\theta_+)$ with \eqref{eq:ct-st} and
$0<\delta\le\delta_0(A)$, and every $e_0$ with
$e_\rho\in H^3(\R)$, $(e_u,e_\theta)\in H^2(\R)$ and
$e_0\in L^1(\R)$ satisfying
\begin{equation}\label{e:dat}
\begin{gathered}
 \Vert e_\rho\Vert_{H^1(\R)}+\Vert e_u\Vert_{L^2(\R)}^2+|Q|
 +\int_\R\Phi\left(\frac{\theta_0(x)}
                             {\bar\theta(x,0)}\right)\dd x\le A,\\
 \rho_0\ge A^{-1},\qquad \theta_0\ge A^{-1}.
\end{gathered}
\end{equation}

Problem \eqref{eq:bg-eqs} with the initial value $W_0$ of
\eqref{eq:bg-dat} has a unique global solution $(\rho,u,\theta)$
such that, for every $S<\infty$,
\begin{equation}\label{e:brg}
\begin{aligned}
 &\rho-\bar\rho\in C([0,S];H^3(\R)),\qquad
 (u-\bar u,\theta-\bar\theta)\in C([0,S];H^2(\R)),\\
 &\partial_x(u,\theta)\in L^2(0,S;H^2(\R)),\\
 &\rho_t\in C([0,S];H^1(\R))\cap L^2(0,S;H^2(\R)),\\
 &(u_t,\theta_t)\in C([0,S];L^2(\R))\cap L^2(0,S;H^1(\R)),
\end{aligned}
\end{equation}
and $\rho,\theta$ lie between two positive constants on
$\R\times[0,\infty)$.
Moreover, for $t\ge0$, with $\bar Q$ of \eqref{eq:ct-w},
\begin{equation}\label{e:r1}
 \Vert(\rho,u,\theta)(t)-\bar Q(t)\Vert_{L^\infty(\R)}
 +\Vert p(t)-p_*\Vert_{L^\infty(\R)}
 \le C_b(1+t)^{-1/2}.
\end{equation}
The $L^1$ assumption defines $Q$ and bounds the initial primitives
used for \eqref{e:r1}; it is not needed for global existence and
positivity (Proposition~\ref{l:bg-entropy}). The threshold $\delta_0(A)$
requires no quantitative bound on $\Vert e_\rho\Vert_{H^3}$,
$\Vert(e_u,e_\theta)\Vert_{H^2}$ or $\Vert e_0\Vert_{L^1}$ beyond
\eqref{e:dat}. These finite norms may affect $C_b$ and the entrance times.
\end{theorem}

\begin{remark}[Meaning of the initial bounds]\label{r:data-condition}
Condition~\eqref{e:dat} bounds the density $H^1$ size, velocity
$L^2$ energy, thermal relative entropy, conserved excess and the
reciprocals of the initial positive lower bounds by an arbitrary
$A\ge1$. The perturbation need not be small. Since
\[
 \Phi(z)\sim\tfrac12(z-1)^2\quad(z\to1),\qquad
 \Phi(z)\sim z\quad(z\to\infty),
\]
the thermal bound permits large temperature peaks and does not
control $\Vert e_\theta\Vert_{L^2}$. Likewise, $|Q|$ controls the
net conserved excess, not the full $L^1$ norm.
\end{remark}

\begin{remark}[Main improvements over previous results]\label{r:oned-comparison}
Large-data stability is known from
\cite{HuangWang2016,PengShiWu2026}. Theorem~\ref{thm:planar}
improves two quantitative aspects of that theory.

\emph{Wave-strength threshold.}
Huang and Wang~\cite[Theorem~1.1]{HuangWang2016} impose a wave-strength
restriction depending on the full initial $H^1$ size. Peng, Shi and
Wu~\cite[Theorem~1.1]{PengShiWu2026} remove the velocity and temperature
gradient norms from this dependence, but retain the velocity $L^4$
and temperature $L^2$ sizes. Here $\delta_0(A)$ depends only on
\eqref{e:dat} and the fixed data. For fixed $A$, the same threshold
therefore applies even when these two norms and the higher Sobolev
norms become arbitrarily large. They may affect $C_b$ and the
entrance times.

\emph{Optimal decay for large perturbations.}
The above large-data results give qualitative convergence.
Liu, Wang and Xu~\cite{LiuWangXu2026} obtain the optimal $L^\infty$
rate $(1+t)^{-1/2}$ for small Sobolev perturbations with small
corrected $L^2$ primitives. In our results, the rate holds for large integrable
perturbations satisfying \eqref{e:dat}; the corrected primitives need
only be bounded, without smallness or square integrability.
Nonzero conserved excess is allowed, and
Corollary~\ref{l:rate-sharp} gives a matching lower bound when an
acoustic mass is nonzero.

These improvements concern the present single-contact-wave data
class. The results in \cite{HuangWang2016,PengShiWu2026} also cover
composite waves and require no initial $L^1$ condition for global
stability; here $L^1$ is used for the quantitative decay.
\end{remark}

\begin{remark}[Eulerian formulation and comparison]\label{r:coordinates}
The initial perturbation and \eqref{e:r1} are specified in Eulerian
coordinates, while $Q$ and the proof use the mass labels
\eqref{e:Y0}, \eqref{eq:m-map}. Lemma~\ref{l:rate} bounds
$X-\bar X$ and $Y-\bar Y$ uniformly, so the solution is compared
with the fixed profile at the same physical position, without a
time-dependent translation. Thus \eqref{e:r1} is an Eulerian
stability estimate proved in mass coordinates. The results in
\cite{HuangWang2016,PengShiWu2026,LiuWangXu2026} are stated in
Lagrangian coordinates.
\end{remark}

\subsection{Notation and proof structure}

Sobolev spaces are inhomogeneous, vector norms are Euclidean, and
$\Vert\cdot\Vert$ denotes the $L^2$ norm. Functions of $y$ are measured
on $\R_y$, and Eulerian functions on $\R_x$. We use
\begin{equation}\label{e:der}
 D_t=\partial_t+u\partial_x,\qquad \tau_s=1+s.
\end{equation}
Constants $C_M$ depend only on an auxiliary bound $M$ for the datum,
the fixed data and any displayed fixed parameters. Constants $C_b$
may also depend on the complete solution datum, and $C_{b,S}$ on a
fixed starting time $S$; they are independent of the final time.
The notation $\omega_b(s)$ denotes a nonnegative quantity depending on
the solution with $\omega_b(s)\to0$ as $s\to\infty$; it may change
between estimates. Gaussian exponents are denoted by $c_G$ or $c_j$.
The proof of Theorem~\ref{thm:planar} specifies the order of choosing
the wave-strength threshold and the entrance times.

\paragraph{Technical difficulties and proof strategy.}
The small-perturbation argument of Liu, Wang and
Xu~\cite{LiuWangXu2026} uses the $L^2$ energy of corrected
primitives. This energy need not be finite here. Although the
physical perturbation eventually becomes small, this alone does
not provide the required $L^2$ condition on its primitive, so their
decay theorem cannot be applied directly after a late restart.

Proposition~\ref{l:bg-entropy} first controls the initial layer by
relative entropy and restarts the global existence theory of
Huang and Wang~\cite{HuangWang2016}. The entrance times may depend
on the complete datum, while the wave-strength threshold depends
only on the lower-order bounds and the fixed data. The initial
$L^1$ assumption gives bounded initial primitives. Their uniform
boundedness in time is proved separately in Lemma~\ref{l:bd}:
a preliminary sub-square-root growth bound is improved by scalar
comparison and an exponential transformation.

The conserved excess determines the contact phase and two acoustic
diffusion waves. Normal-form corrections separate the leading
nonlinear interactions from errors controlled by dissipation tails.
Qualitative convergence alone does not specify the decay rate.
On windows with duration proportional to the inverse square of the
candidate amplitude, local derivative estimates and parabolic kernel
bounds give the conditional contraction \eqref{e:cpt} in
Lemma~\ref{l:win}. Its feedback coefficient is $C\delta$, with $C$
depending only on the fixed data. A doubling argument and this
contraction yield the sharp remainder estimate of
Lemma~\ref{l:asy}, without a global $L^2$ bound on the primitives.
Finally, Lemma~\ref{l:rate} bounds the mass-map displacement and
transfers the decay estimate to the fixed Eulerian contact profile.

These arguments are developed in Sections~\ref{s:bg}--\ref{s:pot}.
Section~\ref{sec:proof} completes the proof of the main theorem.
Appendix~\ref{app:bg} contains the existence, regularity and auxiliary
balance arguments. Appendix~\ref{app:sc} gives the scalar parabolic
proofs. Additional derivative estimates and the weighted pressure
decomposition are collected in Appendices~\ref{s:coef} and~\ref{s:pres};
their use in the planar problem is specified at the start of
Appendix~\ref{s:coef}.
\section{Global solution and characteristic estimates}\label{s:bg}

The global $H^1$ solution follows from Huang and
Wang~\cite{HuangWang2016}. We record its mass-coordinate form and the
regularity and characteristic estimates used in the decay proof.

\subsection{Contact profile and mass coordinates}\label{s:bg-coordinates}

We record the standard mass-coordinate construction of the contact
profile \cite{HuangMatsumuraXin2006,HuangLiMatsumura2010}. The reduced
temperature equation is
\begin{equation}\label{eq:pr-eq}
 \pa_t\Theta_c=\frac{\kappa p_*}{c_pR}\pa_y\Bigl(\frac{\pa_y\Theta_c}{\Theta_c}\Bigr).
\end{equation}
For $\Theta_c(y,t)=\widehat\Theta(\xi)$, $\xi=\frac y{\sqrt{1+t}}$,
\eqref{eq:pr-eq} is equivalent to
\begin{equation}\label{eq:th-ode}
 (\log\widehat\Theta)''=-a\,\xi\widehat\Theta(\log\widehat\Theta)',\qquad
 a=\frac{c_pR}{2\kappa p_*},
\end{equation}
and by van Duyn and Peletier~\cite{vanDuynPeletier1977} this equation
has exactly one solution with $\widehat\Theta(\pm\infty)=\theta_\pm$.
By \eqref{eq:th-ode},
\[
 (\log\widehat\Theta)'(\xi)=(\log\widehat\Theta)'(0)
 \exp\Bigl\{-a\int_0^\xi s\widehat\Theta(s)\dd s\Bigr\}.
\]
Hence $(\log\widehat\Theta)'$ does not change sign, $\widehat\Theta$ is
monotone, and
$\min\{\theta_-,\theta_+\}\le\widehat\Theta\le\max\{\theta_-,\theta_+\}$.
Its mass-coordinate form is
\begin{equation}\label{eq:prof}
 \Theta_c(y,t)=\widehat\Theta\Bigl(\frac y{\sqrt{1+t}}\Bigr),\qquad
 V_c=\frac{R\Theta_c}{p_*},\qquad
 U_c=u_*+\frac{\kappa}{c_p}\frac{\pa_y\Theta_c}{\Theta_c},
\end{equation}
where $V_c$ is the specific volume. By \cite{HuangLiMatsumura2010}, if
$\delta$ is sufficiently small, then for $\pm y\ge0$ and $t\ge0$ (here
and below, either all upper or all lower signs are taken)
\begin{equation}\label{eq:ct-gau}
 (1+t)|\pa_y^2\Theta_c|+(1+t)^{1/2}|\pa_y\Theta_c|+|\Theta_c-\theta_\pm|
 \le C\delta e^{-\frac{c_G y^2}{1+t}},
\end{equation}
where $C,c_G>0$ are independent of $y$, $t$ and $\delta$.
Differentiating \eqref{eq:th-ode} repeatedly and using
$\min\{\theta_-,\theta_+\}\le\widehat\Theta\le\max\{\theta_-,\theta_+\}$
gives Gaussian bounds for every derivative. Absorbing the polynomial
factors by replacing $c_G$ with $c_G/2$, \eqref{eq:prof} yields, for
$y\in\R$ and $t\ge0$,
\begin{equation}\label{eq:pr-gau}
\begin{aligned}
 |\pa_y^k(V_c,\Theta_c)|
 &\le C\delta(1+t)^{-\frac k2}e^{-\frac{c_G y^2}{1+t}} &&(k\ge1),\\
 |\pa_y^k(U_c-u_*)|
 &\le C\delta(1+t)^{-\frac{k+1}2}e^{-\frac{c_G y^2}{1+t}} &&(k\ge0),
\end{aligned}
\end{equation}
where $C$ depends on $k$, and $C,c_G$ are independent of $y$, $t$ and
$\delta$. By \eqref{eq:pr-eq} and \eqref{eq:prof},
\begin{equation}\label{eq:pr-id}
 \pa_tV_c=\pa_yU_c,\qquad
 \pa_tU_c=\frac\kappa{c_p}\pa_y\Bigl(\frac{\pa_yU_c}{V_c}\Bigr),\qquad
 c_v\pa_t\Theta_c+p_*\pa_yU_c=\kappa\pa_y\Bigl(\frac{\pa_y\Theta_c}{V_c}\Bigr),
\end{equation}
and the self-similar form gives
\begin{equation}\label{e:pt}
 (V_c)_t=-\frac{y(V_c)_y}{2\tau},\qquad
 (U_c)_t=-\frac{U_c-u_*+y(U_c)_y}{2\tau},
\end{equation}
so that, by \eqref{eq:pr-gau}, for every integer $j\ge0$,
\begin{equation}\label{e:put}
 \lvert\partial_y^j(U_c)_t(y,t)\rvert
 \le C_j\delta\tau^{-\frac{j+3}2}e^{-\frac{c_jy^2}{\tau}},
\end{equation}
where $C_j,c_j>0$ depend only on $j$ and the fixed data. We put
$H_c=U_c-u_*$ and fix the Eulerian position of the wave by
\begin{equation}\label{e:Xc}
 \bar X(y,t)=\int_0^yV_c(z,t)\dd z+\int_0^tU_c(0,r)\dd r,
 \qquad \bar Y(\cdot,t)=\bar X(\cdot,t)^{-1};
\end{equation}
thus $\bar X_y=V_c$, $\bar X_t=U_c$ and $\bar X(0,0)=0$. The Eulerian
profile \eqref{eq:ct-w} is equivalently given by
\begin{equation}\label{e:ct-pullback}
 \bar Q=(\bar\rho,\bar u,\bar\theta)
 =\left(\frac1{V_c},U_c,\Theta_c\right)\circ\bar Y.
\end{equation}
Since $\bar X_y=V_c$ and $\bar X_t=U_c$, the chain rule in
\eqref{eq:prof} gives
$\bar u-u_*=\frac{\kappa R}{c_pp_*}\partial_x\bar\theta$, and
\eqref{eq:pr-eq} becomes \eqref{eq:ct-euler}.

In the mass labels \eqref{e:Y0}, the solution is represented by
$(V_b,U_b,\Theta_b)$ satisfying
\begin{align}
 (V_b)_t-(U_b)_y&=0,\label{eq:bg-m}\\
 (U_b)_t+(P_b)_y
 &=\nu\partial_y\left(\frac{(U_b)_y}{V_b}\right),\label{eq:bg-mom}\\
 \partial_t\left(c_v\Theta_b+\frac{U_b^2}{2}\right)
 +\partial_y(P_bU_b)
 &=\partial_y\left(\frac{\kappa(\Theta_b)_y+\nu U_b(U_b)_y}{V_b}\right),
 \qquad P_b=\frac{R\Theta_b}{V_b},\label{eq:bg-en}
\end{align}
with the initial value
$(\frac1{\rho_0},u_0,\theta_0)\circ X_0$. The map
\begin{equation}\label{eq:m-map}
 X(y,t)=\int_0^yV_b(z,t)\dd z+\int_0^tU_b(0,r)\dd r,
 \qquad Y(\cdot,t)=X(\cdot,t)^{-1}
\end{equation}
satisfies $X_y=V_b$, $X_t=U_b$ and $X(\cdot,0)=X_0$, and
\begin{equation}\label{eq:bg-eul}
 (\rho,u,\theta)=\left(\frac1{V_b},U_b,\Theta_b\right)\circ Y,
 \qquad p=P_b\circ Y.
\end{equation}

For a function $F(x,t)$ we write
$F\circ X=F(X(y,t),t)$; in particular,
\begin{equation}\label{e:chain}
 \partial_y(F\circ X)=V_b(\partial_xF)\circ X,\qquad
 \partial_t(F\circ X)=(D_tF)\circ X,\qquad
 \dd x=V_b\dd y.
\end{equation}
Mass-coordinate functions are composed with $Y$ when used in Eulerian
coordinates. The maps $X$ and $\bar X$ retain their initial
normalization at every later starting time.

Put
\begin{equation}\label{e:Um}
 b_m=U_b-u_*,\qquad
 Z_m=(\varphi_m,\upsilon_m,\vartheta_m)
    =(V_b-V_c,\,b_m-H_c,\,\Theta_b-\Theta_c).
\end{equation}
Throughout this section $M\ge1$,
\begin{equation}\label{e:nat}
 \Vert e_\rho\Vert_{H^3(\R)}
 +\Vert(e_u,e_\theta)\Vert_{H^2(\R)}\le M,\qquad
 \rho_0,\ \theta_0\ge M^{-1},
\end{equation}
and $0<\delta\le\delta_1(M)$, where $\delta_1(M)>0$ depends only on $M$
and $(R,\gamma,\nu,\kappa,\rho_-,\theta_-)$ and may be decreased
in each lemma. Subtracting $U_b$ times \eqref{eq:bg-mom} from
\eqref{eq:bg-en} gives
\begin{equation}\label{e:bgth}
 c_v(\Theta_b)_t+P_b(U_b)_y
 =\kappa\partial_y\frac{(\Theta_b)_y}{V_b}+\frac{\nu(U_b)_y^2}{V_b}.
\end{equation}

\subsection{Global existence}

\begin{lemma}[Global one-dimensional solution]\label{l:bg-known}
Problem \eqref{eq:bg-m}--\eqref{eq:bg-en} with initial value
$(\frac1{\rho_0},u_0,\theta_0)\circ X_0$ admits a positive global
solution with $Z_m\in C([0,\infty);H^1)$ and
$\inf_{\R\times[0,T_e]}V_b>0$ for every $T_e<\infty$. It satisfies
\begin{equation}\label{e:blow}
\begin{gathered}
 K(T)=\sup_{0\le t\le T}\Vert Z_m(t)\Vert_{H^1}^2
 +\int_0^T\bigl(\Vert(\varphi_m)_y\Vert_{L^2}^2
       +\Vert\partial_y(\upsilon_m,\vartheta_m)\Vert_{H^1}^2\bigr)\dd t,\\
 K(\infty)<\infty,\qquad \lim_{t\to\infty}\Vert Z_m(t)\Vert_{L^\infty}=0.
\end{gathered}
\end{equation}
\end{lemma}
\begin{proof}
Put $B=(V_c,U_c,\Theta_c)(\cdot,0)$. By \eqref{e:Y0}, \eqref{e:Xc} and
\eqref{eq:bg-dat}, $|Y_0(x)-\bar Y(x,0)|=|\int_0^{x}e_\rho\dd z|
\le\Vert e_\rho\Vert_{L^2}|x|^{1/2}$. With
$Z_m(\cdot,0)=(-\frac{e_\rho}{\rho_0\bar\rho(\cdot,0)},e_u,e_\theta)\circ X_0
+B\circ\bar Y(\cdot,0)\circ X_0-B$, the chain rule (the maps vanish at
$0$ and have derivatives $\rho_0,\bar\rho(\cdot,0)\in[C_M^{-1},C_M]$),
the mean value theorem, \eqref{eq:pr-gau} and \eqref{e:nat} give
\begin{equation}\label{e:bgin}
 \Vert\varphi_m(0)\Vert_{H^3}+\Vert(\upsilon_m,\vartheta_m)(0)\Vert_{H^2}
 \le C_M,\qquad
 \inf_\R V_b(\cdot,0),\ \inf_\R\Theta_b(\cdot,0)\ge C_M^{-1}.
\end{equation}

Apply \cite[Theorem~1.1 and Proposition~3.1]{HuangWang2016}, with the
viscosity there replaced by $\nu$ and the common velocity $u_*$
subtracted. Choosing $\delta_1(M)$ below their threshold for
$\Vert Z_m(0)\Vert_{H^1}\le2C_M$ and
$\inf V_b(\cdot,0),\inf\Theta_b(\cdot,0)\ge(2C_M)^{-1}$ gives the
assertions. The finite-interval positivity used in this continuation is
made explicit in Appendix~\ref{app:bg-positive}, using the estimates
of the cited proof before taking any large-time limit.
Proposition~\ref{l:bg-entropy} gives the alternative threshold
controlled by the initial relative entropy, the initial volume $H^1$ size and
the positive lower bounds. Its conversion to Eulerian coordinates is
given in the proof of Theorem~\ref{thm:planar}.
\end{proof}

\subsection{Regularity and coordinate conversion}

For $0<T_e<\infty$, the class \eqref{e:bgcl} on $[0,T_e]$ consists of
the $(V_b,U_b,\Theta_b)$ such that, for every $t_0\in(0,T_e)$,
\begin{equation}\label{e:bgcl}
\begin{aligned}
 &\varphi_m\in C([0,T_e];H^3),\qquad
 (\upsilon_m,\vartheta_m)\in C([0,T_e];H^2)\cap C([t_0,T_e];H^3),\\
 &\partial_y(\upsilon_m,\vartheta_m)
       \in L^2(0,T_e;H^2)\cap L^2(t_0,T_e;H^3),\\
 &(V_b)_t\in C([0,T_e];H^1)\cap L^2(0,T_e;H^2),\qquad
 ((U_b)_t,(\Theta_b)_t)\in C([0,T_e];L^2)\cap L^2(0,T_e;H^1).
\end{aligned}
\end{equation}

\begin{lemma}\label{l:bg}
The solution of Lemma~\ref{l:bg-known} belongs to the class
\eqref{e:bgcl} on every $[0,T_e]$. There are $c_b,C_b>0$ such that
$c_b\le V_b,\Theta_b\le C_b$ on $\R\times[0,\infty)$ and, for every
$t_0>0$,
\[
\begin{gathered}
 \sup_{t\ge0}\bigl(\Vert Z_m\Vert_{H^2}^2
       +\Vert\partial_y^3\varphi_m\Vert_{L^2}^2\bigr)
 +\int_0^\infty\bigl(\Vert\partial_y\varphi_m\Vert_{H^2}^2
       +\Vert\partial_y(\upsilon_m,\vartheta_m)\Vert_{H^2}^2\bigr)\dd t\\
 +\sup_{t\ge t_0}\Vert\partial_y^3(\upsilon_m,\vartheta_m)\Vert_{L^2}^2
 +\int_{t_0}^\infty\Vert\partial_y^4(\upsilon_m,\vartheta_m)\Vert_{L^2}^2\dd t
 <\infty,
\end{gathered}
\]
\[
 \lim_{t\to\infty}\Bigl(\Vert Z_m(t)\Vert_{L^\infty}
       +\sum_{j=1}^3\Vert\partial_y^jZ_m(t)\Vert_{L^2}\Bigr)=0.
\]
The map \eqref{eq:m-map} satisfies $X_y=V_b$, $X_t=U_b$ and
$X(\cdot,0)=X_0$, and the fields \eqref{eq:bg-eul} solve
\eqref{eq:bg-eqs} with initial value $W_0$ in \eqref{e:brg} for every
$S<\infty$. The solution is unique among the solutions with this
initial value such that, for every $T_e<\infty$,
$Z_m\in C([0,T_e];H^1)$, $\partial_y(U_b,\Theta_b)\in L^2(0,T_e;H^1)$
and $V_b,\Theta_b$ lie between positive constants on
$\R\times[0,T_e]$; \eqref{eq:bg-eul} is unique in the corresponding
Eulerian class relative to $\bar Q$.
\end{lemma}
We use the one-dimensional interpolation inequality
\begin{equation}\label{e:agm}
 \Vert f\Vert_{L^\infty}^2\le2\Vert f\Vert_{L^2}\Vert f_y\Vert_{L^2}.
\end{equation}

The proof of Lemma~\ref{l:bg} is given in Appendix~\ref{app:regularity}.

\subsection{Characteristic balances}

The next two auxiliary lemmas are proved in Appendix~\ref{app:balances}.

\begin{lemma}\label{l:osp}
Let $t_1<t_2$, and for $i=1,2$ let $E_i\ge0$, $F_i$, $S_i$ be functions
on $\R\times[t_1,t_2]$ with
\[
 E_i\in C([t_1,t_2];L^1(\R))\cap L^\infty(\R\times(t_1,t_2)),\qquad
 F_i,\ S_i\in L^1(\R\times(t_1,t_2)),
\]
and $\partial_tE_i+\partial_xF_i=S_i$ in $\mathcal D'(\R\times(t_1,t_2))$.
Put $P_i(x,t)=\int_x^\infty E_i(z,t)\dd z$. Then
$0\le P_i(x,t)\le\Vert E_i(t)\Vert_{L^1}$, $\partial_xP_i=-E_i$, and in
$\mathcal D'(\R\times(t_1,t_2))$
\begin{equation}\label{e:osp1}
 \partial_tP_i=F_i+\int_x^\infty S_i(z,t)\dd z,\qquad
 \partial_t\int_{-\infty}^xE_i(z,t)\dd z
 =-F_i+\int_{-\infty}^xS_i(z,t)\dd z .
\end{equation}
Moreover, $I(t)=\int_\R E_1(x,t)P_2(x,t)\dd x$ satisfies
\begin{equation}\label{e:osp2}
 I(t_2)-I(t_1)
 =\int_{t_1}^{t_2}\!\!\int_\R\Bigl(E_1F_2-F_1E_2+S_1P_2
   +E_1\int_x^\infty S_2(z,t)\dd z\Bigr)\dd x\dd t .
\end{equation}
\end{lemma}

Let $\chi_e$ be smooth and even, with $0\le\chi_e\le1$, $\chi_e=1$ on
$[-1,1]$ and $\chi_e=0$ outside $[-2,2]$. For $\varepsilon>0$ put
\begin{equation}\label{e:gep}
 g_\varepsilon(z)=\int_0^z\chi_e\Bigl(\frac a\varepsilon\Bigr)\dd a,
 \qquad \Phi_\varepsilon(z)=\int_0^zg_\varepsilon(a)\dd a .
\end{equation}
Then $g_\varepsilon(z)=z$ for $|z|\le\varepsilon$,
$|g_\varepsilon(z)|\le\min\{|z|,2\varepsilon\}$ and
$0\le g_\varepsilon'\le1$.

\begin{lemma}\label{l:trc}
Let $s\ge0$, $\bar z>0$, $w_1<w_2$, $K=[-\bar z,\bar z]\times[w_1,w_2]$,
$\Lambda\in C^2(K)$ with $\Lambda\ge\Lambda_0>0$, and $\alpha_0\in\R$.
Let $u,w,\alpha,F,f,R,G$ be functions on $\R\times[s,\infty)$ with
$(u,w)\in K$ and
\[
 u_t+\alpha u_y=F_y+R,\qquad w_t=f_y+G,
\]
such that, for every $T>s$, $u\in C([s,T];H^1)\cap C^1([s,T];L^2)$,
$w_y,w_t\in C([s,T];L^2)$, $\alpha\in L^\infty(\R\times(s,T))$,
$F,f\in C([s,T];H^1)$ and $R,G\in C([s,T];L^1)$. For $h(a)=a$ or
$h=g_\varepsilon$ put
$E_h(z,w)=\int_0^z\frac{h(a)}{\Lambda(a,w)}\dd a$,
$\Phi_h(z)=\int_0^zh(a)\dd a$, $J_h=(E_h)_zF+(E_h)_wf$ and
\[
\begin{aligned}
 S_h={}&(E_h)_zR+(E_h)_wG
   -\bigl[(E_h)_{zz}u_y+(E_h)_{zw}w_y\bigr]F
   -\bigl[(E_h)_{zw}u_y+(E_h)_{ww}w_y\bigr]f\\
 &-\bigl(\alpha-\alpha_0\Lambda(u,w)\bigr)(E_h)_zu_y,
\end{aligned}
\]
with the derivatives of $E_h$ evaluated at $(u,w)$. Then, for
$s\le t'\le t$,
\begin{equation}\label{e:trc1}
 \partial_t[E_h(u,w)]+\alpha_0\partial_y\Phi_h(u)=\partial_yJ_h+S_h,
 \qquad
 \int_\R E_h(u,w)\dd y\Big|_{t'}^{t}
 =\int_{t'}^t\!\!\int_\R S_h\dd y\dd r,
\end{equation}
and $|S_h|\le C\mathcal M$, with $C$ depending only on $K$, $\Lambda_0$ and
$\Vert\Lambda\Vert_{C^2(K)}$, where
\[
 \mathcal M=|u_yF|+|u|\bigl(|R|+|w_yF|+|u_yf|
   +|\alpha-\alpha_0\Lambda(u,w)||u_y|\bigr)+u^2\bigl(|G|+|w_yf|\bigr).
\]
If moreover $\mathcal M\in L^1(\R\times(s,\infty))$ and
$\lim_{t\to\infty}\Vert u(t)\Vert_{L^\infty}=0$, then
$\lim_{t\to\infty}\Vert u(t)\Vert_{L^2}=0$.
\end{lemma}

For the solution of Lemma~\ref{l:bg} and $\sigma=\pm1$ define
\begin{equation}\label{e:char}
\begin{gathered}
 m=V_b\Bigl(\frac{P_b}{p_*}\Bigr)^{1/\gamma},\qquad e_m=m-V_c,\qquad
 H_T=c\partial_y\log\Theta_b,\qquad
 f_m=\Bigl(\frac{P_b}{p_*}\Bigr)^{1/\gamma}H_T,\qquad f_e=f_m-H_c,\\
 J_a=\frac{2\sqrt{\gamma p_*}\sqrt m}{\gamma-1}
 \Bigl[\Bigl(\frac{P_b}{p_*}\Bigr)^{\frac{\gamma-1}{2\gamma}}-1\Bigr],\qquad
 q_\pm=m^{-1/4}(\upsilon_m\pm J_a),\\
 \Gamma=\bigl((b_m)_y,(P_b)_y,H_T-H_c,(e_m)_y,(q_+)_y,(q_-)_y,
          m_yq_+,m_yq_-,H_Tq_+,H_Tq_-\bigr),\\
 c_\sigma(z,m)=\sqrt{\gamma p_*}\,m^{-1/2}
    \Bigl(1+\frac{\sigma(\gamma-1)m^{-1/4}z}{4\sqrt{\gamma p_*}}\Bigr)^{\frac{\gamma+1}{\gamma-1}},
 \qquad
 E_\sigma(q,m)=\int_0^q\frac{z\dd z}{c_\sigma(z,m)}.
\end{gathered}
\end{equation}

\begin{lemma}\label{l:dep}
There are $S_0\ge1$ and functions $S_e$ on $\R\times(0,\infty)$ and
$J_\pm,S_\pm$ on $\R\times[S_0,\infty)$ with the following properties,
where $c_{\rm ac}=\frac14\sqrt{\gamma p_*\rho_-}$ and $C_*$ depends only
on the fixed data. First,
\begin{equation}\label{e:mbal}
 \partial_t\frac{e_m^2}2=\partial_y(e_mf_e)+S_e\quad
 \text{on }\R\times(0,\infty),\qquad
 S_e\in L^1(\R\times(0,\infty)).
\end{equation}
Second, for $\sigma=\pm1$, on $\R\times[S_0,\infty)$,
$c_\sigma(z,m)\ge2c_{\rm ac}$ for $z$ between $0$ and $q_\sigma$, and
\begin{equation}\label{e:hcur}
\begin{gathered}
 \partial_t[E_\sigma(q_\sigma,m)]
       +\sigma\partial_y\frac{q_\sigma^2}{2}=(J_\sigma)_y+S_\sigma,\qquad
 C_*^{-1}q_\sigma^2\le E_\sigma(q_\sigma,m)\le C_*q_\sigma^2,\\
 C_*^{-1}\bigl((P_b-p_*)^2+\upsilon_m^2\bigr)\le q_+^2+q_-^2
 \le C_*\bigl((P_b-p_*)^2+\upsilon_m^2\bigr),\\
 |J_\sigma|+|\partial_y[E_\sigma(q_\sigma,m)]|\le C_*|q_\sigma||\Gamma|,
 \qquad |f_e|+|(e_m)_y|\le C_*|\Gamma|,
\end{gathered}
\end{equation}
with $S_\pm\in L^1(\R\times(S_0,\infty))$. Third,
$\int_0^\infty\Vert\Gamma\Vert_{H^1}^2\dd t<\infty$,
\[
 \lim_{t\to\infty}\bigl(\Vert(e_m,q_+,q_-)(t)\Vert_{H^1}
       +\Vert f_m(t)\Vert_{L^\infty}\bigr)=0,
\]
there is a fixed $c_*>0$ such that, for every $0<c_G\le c_*$,
\[
 \int_0^\infty\!\!\int_\R
       \tau^{-1}e^{-\frac{2c_Gy^2}{\tau}}|Z_m|^2\dd y\dd t<\infty,
\]
and the following quantity is finite for $r\ge S_0$ and tends to zero
as $r\to\infty$:
\[
\begin{gathered}
 \sup_{t\ge r}\Vert Z_m(t)\Vert_{H^3}^2
 +\int_r^\infty\bigl(\Vert\partial_y\varphi_m\Vert_{H^2}^2
       +\Vert\partial_y(\upsilon_m,\vartheta_m)\Vert_{H^3}^2
       +\Vert\Gamma\Vert_{L^\infty}^2
       +\Vert(S_+,S_-,S_e)\Vert_{L^1}\bigr)\dd t\\
 +\int_r^\infty\!\!\int_\R
       \Bigl[q_+^2q_-^2+(q_+^2+q_-^2)
       \Bigl(e_m^2+\tau^{-1/2}\operatorname{sech}^2\frac y{\sqrt\tau}\Bigr)
       \Bigr]\dd y\dd t .
\end{gathered}
\]
\end{lemma}
\begin{proof}
By \cite[Lemma~2.3]{HuangWang2016}, applied to Lemma~\ref{l:bg},
there is a fixed $c_*>0$
such that, for $0<c_G\le c_*$ and every $T_e<\infty$,
\[
 \int_0^{T_e}\!\!\int_\R\tau^{-1}e^{-\frac{2c_Gy^2}{\tau}}|Z_m|^2\dd y\dd t
 \le C_{b,c_G}\Bigl(1+\int_0^{T_e}\Vert\partial_yZ_m\Vert_{L^2}^2\dd t\Bigr)
 \le C_{b,c_G} ,
\]
where $C_{b,c_G}$ may depend on the chosen exponent, but not on $T_e$.
Choose $c_G$ no larger than the exponent in \eqref{eq:pr-gau}; then
\begin{equation}\label{e:gprod}
 |Z_m|\bigl(|(V_c)_y|+|(\Theta_c)_y|+|H_c|\bigr)
 \in L^2((0,\infty)\times\R).
\end{equation}
By \eqref{eq:bg-m}, \eqref{e:bgth}, $P_b=\frac{R\Theta_b}{V_b}$ and
$R=(\gamma-1)c_v$,
\begin{equation}\label{e:pbt}
 (P_b)_t+\frac{\gamma P_b}{V_b}(b_m)_y
 =\frac{(\gamma-1)\kappa}{V_b}\partial_y\frac{(\Theta_b)_y}{V_b}
  +\frac{(\gamma-1)\nu}{V_b^2}(b_m)_y^2 .
\end{equation}
The identities
\[
\begin{aligned}
 (P_b)_y&=\frac{R(\vartheta_m)_y-p_*(\varphi_m)_y}{V_b}
                            -\frac{(P_b-p_*)(V_b)_y}{V_b},\qquad
 H_T-H_c=c\Bigl(\frac{(\vartheta_m)_y}{\Theta_b}
             -\frac{(\Theta_c)_y\vartheta_m}{\Theta_b\Theta_c}\Bigr),\\
 (e_m)_y&=\Bigl(\frac{P_b}{p_*}\Bigr)^{1/\gamma}(\varphi_m)_y
   +\Bigl[\Bigl(\frac{P_b}{p_*}\Bigr)^{1/\gamma}-1\Bigr](V_c)_y
   +\frac{m}{\gamma P_b}(P_b)_y ,
\end{aligned}
\]
which follow from $R\Theta_c=p_*V_c$ and $H_c=c\partial_y\log\Theta_c$,
together with \eqref{e:gprod}, Lemma~\ref{l:bg}, \eqref{eq:pr-gau} and
the equations \eqref{eq:bg-mom}, \eqref{e:bgth}, \eqref{e:pbt}, give
\begin{equation}\label{e:bprod}
\begin{gathered}
 (P_b)_y,\ H_T-H_c,\ (e_m)_y,\ b_mH_c,
     \Bigl(\frac1{V_b}-\frac1{V_c}\Bigr)(V_c)_y
                  \in L^2((0,\infty)\times\R),\\
 \int_0^\infty\bigl(\Vert(b_m)_y\Vert_{H^1}^2
       +\Vert((b_m)_t,(\Theta_b)_t,(P_b)_t)\Vert_{L^2}^2\bigr)\dd t<\infty .
\end{gathered}
\end{equation}

Since $\log m=\log V_b+\frac1\gamma\log\frac{P_b}{p_*}$, \eqref{eq:bg-m}
and \eqref{e:pbt} give, with
$A_e(p)=\frac{\gamma-1}{\gamma p}(\frac p{p_*})^{1/\gamma}$,
\[
 m_t=A_e(P_b)\Bigl(\kappa\partial_y\frac{(\Theta_b)_y}{V_b}
       +\frac{\nu(b_m)_y^2}{V_b}\Bigr).
\]
By $\kappa\frac{(\Theta_b)_y}{V_b}=\frac{c_pP_b}{R}H_T$,
$A_e(P_b)\frac{c_pP_b}{R}=(\frac{P_b}{p_*})^{1/\gamma}$ and
$-A_e'(p)\frac{c_pp}{R}=A_e(p)$, this is
\begin{equation}\label{e:mt}
 m_t=(f_m)_y+G_m,\qquad
 G_m=A_e(P_b)H_T(P_b)_y+\frac{\nu A_e(P_b)}{V_b}(b_m)_y^2 ,
\end{equation}
and $H_T=c\frac{m_y}{m}+\frac{c(\gamma-1)(P_b)_y}{\gamma P_b}$ gives
\begin{equation}\label{e:fm2}
 f_m=\frac c{V_b}m_y+cA_e(P_b)(P_b)_y .
\end{equation}
With $(V_c)_t=(H_c)_y$, \eqref{e:mt} gives $(e_m)_t=(f_e)_y+G_m$ and
the balance in \eqref{e:mbal} with $S_e=-(e_m)_yf_e+e_mG_m$.

Multiplying $(b_m)_t+(P_b)_y=\nu\partial_y\frac{(b_m)_y}{V_b}$ by
$A_e(P_b)H_c$, inserting the result into the term $A_e(P_b)H_c(P_b)_y$
of $G_m$, and using \eqref{e:fm2} and $H_c=\frac{c(V_c)_y}{V_c}$ in
$f_e$, we obtain, for $\widetilde e_m=e_m+A_e(P_b)H_cb_m$,
\[
 (\widetilde e_m)_t
   =\partial_y\Bigl(\frac c{V_b}(\widetilde e_m)_y+F_e\Bigr)
     +\widetilde S_e,
\]
\[
\begin{aligned}
 F_e&=c\Bigl(\frac1{V_b}-\frac1{V_c}\Bigr)(V_c)_y
       +cA_e(P_b)(P_b)_y
       +\frac{\nu A_e(P_b)H_c(b_m)_y}{V_b}
       -\frac c{V_b}\partial_y\bigl(A_e(P_b)H_cb_m\bigr),\\
 \widetilde S_e&=A_e(P_b)(H_T-H_c)(P_b)_y
       +\frac{\nu A_e(P_b)}{V_b}(b_m)_y^2
       +b_mH_c\partial_tA_e(P_b)+b_mA_e(P_b)(H_c)_t\\
    &\quad-\frac{\nu(b_m)_yH_c}{V_b}\partial_yA_e(P_b)
       -\frac{\nu A_e(P_b)(b_m)_y(H_c)_y}{V_b}.
\end{aligned}
\]
The estimates \eqref{e:bprod}, \eqref{e:gprod}, \eqref{eq:pr-gau} and,
through $\int_0^\infty\Vert b_m\Vert_{L^2}\Vert(H_c)_t\Vert_{L^2}\dd t<\infty$,
\eqref{e:put} yield, over $(0,\infty)\times\R$, $F_e\in L^2$,
$\widetilde S_e\in L^1$,
\[
 f_e=\Bigl(\frac{P_b}{p_*}\Bigr)^{1/\gamma}(H_T-H_c)
       +\Bigl[\Bigl(\frac{P_b}{p_*}\Bigr)^{1/\gamma}-1\Bigr]H_c
 \in L^2,
\]
and, since $e_m$ is bounded and $e_mH_c\in L^2$ by \eqref{e:gprod},
$S_e\in L^1$. Lemma~\ref{l:trc} with $s=1$, $u=\widetilde e_m$,
$w\equiv1$, $\Lambda\equiv1$, $\alpha=\alpha_0=0$,
$F=\frac c{V_b}(\widetilde e_m)_y+F_e$, $R=\widetilde S_e$ and $f=G=0$
applies by \eqref{e:bgcl}, with
$\mathcal M=|(\widetilde e_m)_yF|+|\widetilde e_m||\widetilde S_e|$
integrable on $(1,\infty)\times\R$ and, by $|e_m|\le C_b|Z_m|$ and
$|H_c|\le C\delta\tau^{-1/2}$, $\Vert\widetilde e_m(t)\Vert_{L^\infty}\to0$.
Hence $\Vert\widetilde e_m(t)\Vert_{L^2}\to0$, and
$\Vert\widetilde e_m-e_m\Vert_{L^2}\le C_b\delta\tau^{-1/2}\Vert b_m\Vert_{L^2}$
gives $\Vert e_m(t)\Vert_{L^2}\to0$.

Put $C_a=\sqrt{\gamma P_b/V_b}$. At fixed $m$, by
$V_b=m(\frac{P_b}{p_*})^{-1/\gamma}$, $\partial_{P_b}J_a=C_a^{-1}$; at
fixed $P_b$, $\partial_mJ_a=\frac{J_a}{2m}$. With
$\sigma m^{-1/4}J_a=\frac{q_\sigma-q_{-\sigma}}2$ this gives, for
$\partial=\partial_y$ and $\partial=\partial_t$,
\begin{equation}\label{e:dq}
 \partial q_\sigma=m^{-1/4}\partial\upsilon_m
       +\sigma m^{-1/4}C_a^{-1}\partial P_b-\frac{q_{-\sigma}}{4m}\partial m .
\end{equation}
We insert $(\upsilon_m)_t=-(P_b)_y+\nu\partial_y\frac{(b_m)_y}{V_b}-(H_c)_t$,
\eqref{e:pbt} with
$\frac{(\gamma-1)\kappa}{V_b}\partial_y\frac{(\Theta_b)_y}{V_b}
=\frac{\gamma}{V_b}(P_bH_T)_y$, and \eqref{e:mt} into \eqref{e:dq} with
$\partial=\partial_t$, and use \eqref{e:dq} with $\partial=\partial_y$ in
the form
\[
 -m^{-1/4}(P_b)_y-\sigma m^{-1/4}C_a(b_m)_y
 =-\sigma C_a(q_\sigma)_y-\frac{\sigma C_am_yq_{-\sigma}}{4m}
   -\sigma m^{-1/4}C_a(H_c)_y .
\]
With
$\frac{\sigma m^{-1/4}}{C_a}\frac{\gamma}{V_b}(P_bH_T)_y=\sigma m^{-1/4}C_a(H_T)_y
+\frac{\sigma m^{-1/4}}{C_a}\frac{\gamma H_T}{V_b}(P_b)_y$, this gives
\begin{equation}\label{e:aceq}
 (q_\sigma)_t+\sigma C_a(q_\sigma)_y
   =(F_\sigma)_y+R_\sigma-\frac{\sigma C_am_yq_{-\sigma}}{4m},
\end{equation}
\[
\begin{aligned}
 F_\sigma={}&\frac{\nu m^{-1/4}(b_m)_y}{V_b}
             +\sigma m^{-1/4}C_a(H_T-H_c)-\frac{q_{-\sigma}}{4m}f_m,\\
 R_\sigma={}&-\frac{\nu\partial_y(m^{-1/4})(b_m)_y}{V_b}
       -\sigma\partial_y(m^{-1/4}C_a)(H_T-H_c)
       +\partial_y\Bigl(\frac{q_{-\sigma}}{4m}\Bigr)f_m
       -\frac{q_{-\sigma}}{4m}G_m-m^{-1/4}(H_c)_t\\
 &+\frac{\sigma m^{-1/4}}{C_a}
       \Bigl(\frac{\gamma H_T}{V_b}(P_b)_y
                  +\frac{(\gamma-1)\nu}{V_b^2}(b_m)_y^2\Bigr).
\end{aligned}
\]

The function $V_c$ takes values between $\frac1{\rho_-}$ and
$\frac1{\rho_+}$, and $|\frac1{\rho_+}-\frac1{\rho_-}|\le C\delta$. Fix
$z_*,m_*>0$, depending only on the fixed data, such that
$c_\sigma\ge2c_{\rm ac}$ on
$K_*=[-z_*,z_*]\times[\rho_-^{-1}-m_*,\rho_-^{-1}+m_*]$ (note
$c_\sigma(0,\rho_-^{-1})=4c_{\rm ac}$), and decrease $\delta_1(M)$ so
that $C\delta_1(M)\le\frac{m_*}2$. Since $|q_\pm|+|e_m|\le C_b|Z_m|$ and
$\Vert Z_m(t)\Vert_{L^\infty}\to0$, there is $S_0\ge1$ with
$(q_\sigma,m)\in K_*$ on $\R\times[S_0,\infty)$; then $K_*$ contains the
segments from $(0,m)$ to $(q_\sigma,m)$. In this proof $C_*$ depends
only on $K_*$ and the fixed data. On $K_*$,
\[
\begin{gathered}
 C_*^{-1}q^2\le E_\sigma(q,m)\le\frac{q^2}{4c_{\rm ac}},\qquad
 |(E_\sigma)_q|+|(E_\sigma)_{qm}|\le C_*|q|,\\
 |(E_\sigma)_{qq}|\le C_*,\qquad
 |(E_\sigma)_m|+|(E_\sigma)_{mm}|\le C_*q^2 ,
\end{gathered}
\]
and, on $\R\times[S_0,\infty)$, $(\frac{P_b}{p_*})^{\frac{\gamma-1}{2\gamma}}
=1+\frac{(\gamma-1)m^{-1/4}(q_+-q_-)}{4\sqrt{\gamma p_*}}$ and
$C_a^2=\frac{\gamma p_*}m(\frac{P_b}{p_*})^{\frac{\gamma+1}\gamma}$ give
$C_a=c_\sigma(q_\sigma-q_{-\sigma},m)$, hence
$|c_\sigma(q_\sigma,m)-C_a|\le C_*|q_{-\sigma}|$, and
$C_*^{-1}|P_b-p_*|\le|J_a|\le C_*|P_b-p_*|$.
Lemma~\ref{l:trc} applies on $[S_0,\infty)$ with $K=K_*$, $h(a)=a$ and
\begin{equation}\label{e:acch}
 (u,w,\Lambda,\alpha,\alpha_0,F,R,f,G)
 =\Bigl(q_\sigma,m,c_\sigma,\sigma C_a,\sigma,F_\sigma,
        R_\sigma-\frac{\sigma C_am_yq_{-\sigma}}{4m},f_m,G_m\Bigr):
\end{equation}
its equations are \eqref{e:aceq} and \eqref{e:mt}, and its regularity
assumptions hold by \eqref{e:bgcl}. Then $E_h=E_\sigma$ and
$\Phi_h(z)=\frac{z^2}2$, and with $J_\sigma=J_h=(E_\sigma)_qF_\sigma
+(E_\sigma)_mf_m$ and $S_\sigma=S_h$ the identity \eqref{e:trc1} is the
balance in \eqref{e:hcur}. By \eqref{e:char} and
$|(\frac{P_b}{p_*})^{1/\gamma}-1|\le C_*|J_a|\le C_*(|q_+|+|q_-|)$,
\[
\begin{gathered}
 |F_\sigma|+|q_\sigma f_m|\le C_*|\Gamma|,\qquad
 |\partial_y[E_\sigma(q_\sigma,m)]|
       \le C_*|q_\sigma|\bigl(|(q_\sigma)_y|+|m_yq_\sigma|\bigr),\\
 |f_e|\le C_*\bigl(|H_T-H_c|+|H_Tq_+|+|H_Tq_-|\bigr),
\end{gathered}
\]
and these bounds, those on $K_*$ and
$q_+^2+q_-^2=2m^{-1/2}(\upsilon_m^2+J_a^2)$ give the remaining
assertions of \eqref{e:hcur}.

By \eqref{e:bprod}, \eqref{e:dq} with $\partial=\partial_y$,
$m_y=(e_m)_y+(V_c)_y$, $H_T=(H_T-H_c)+H_c$, $|q_\pm|\le C_b|Z_m|$ and
\eqref{e:gprod}, $\Gamma\in L^2((0,\infty)\times\R)$. Differentiating
once more, Lemma~\ref{l:bg}, \eqref{eq:pr-gau} and
$\Vert Z_m\partial_y^2(V_c,\Theta_c)\Vert_{L^2}
+\Vert Z_m|\partial_y(V_c,\Theta_c)|^2\Vert_{L^2}\le C_b\delta\tau^{-1}$
give $\int_0^\infty\Vert\Gamma_y\Vert_{L^2}^2\dd t<\infty$; with
\eqref{e:agm}, $\int_0^\infty\Vert\Gamma\Vert_{L^\infty}^2\dd t<\infty$.

Let $\mathcal M_\sigma$ be the function $\mathcal M$ of
Lemma~\ref{l:trc} for \eqref{e:acch}, so that
$|S_\sigma|\le C_*\mathcal M_\sigma$, and put
\[
 \mathcal M_\sigma'=|(q_\sigma)_yF_\sigma|
   +|q_\sigma|\bigl(|R_\sigma|+|m_yF_\sigma|+|(q_\sigma)_yf_m|\bigr)
   +q_\sigma^2\bigl(|G_m|+|m_yf_m|\bigr).
\]
By $|c_\sigma(q_\sigma,m)-C_a|\le C_*|q_{-\sigma}|$,
$\mathcal M_\sigma\le\mathcal M_\sigma'
+C_*|q_+q_-|(|m_y|+|(q_\sigma)_y|)$. Since $(q_\sigma)_y$, $m_yq_\sigma$,
$H_Tq_\sigma$ are components of $\Gamma$, $|F_\sigma|\le C_*|\Gamma|$ and
$|f_m|\le C_*|H_T|$, \eqref{e:bprod}, \eqref{e:put} and
$\Gamma\in L^2((0,\infty)\times\R)$ give
$\mathcal M_\sigma'\in L^1((S_0,\infty)\times\R)$. Thus
\[
 N=1+\sup_{t\ge S_0}\Vert(q_+,q_-)(t)\Vert_{H^1}^2
   +\int_{S_0}^\infty\Bigl(\Vert\Gamma\Vert_{L^2}^2
       +\sum_\sigma\Vert\mathcal M_\sigma'\Vert_{L^1}\Bigr)\dd t<\infty .
\]
For $T>S_0$ and $L>0$ put
\[
\begin{gathered}
 \omega=\tau^{-1/2}\operatorname{sech}^2\frac y{\sqrt\tau},\qquad
 B_T=\int_{S_0}^T\!\!\int_\R q_+^2q_-^2\dd y\dd t,\qquad
 X_T=\delta\int_{S_0}^T\!\!\int_\R\omega(q_+^2+q_-^2)\dd y\dd t,\\
 w_\sigma=1-\sigma L\delta\tanh\frac y{\sqrt\tau},\qquad
 \mathcal H(t)=\sum_\sigma\int_\R w_\sigma E_\sigma(q_\sigma,m)\dd y .
\end{gathered}
\]
Splitting $m_y=(V_c)_y+(e_m)_y$ and using $|(V_c)_y|\le C_*\delta\omega$,
\[
 \sum_\sigma\int_{S_0}^T\!\!\int_\R|S_\sigma|\dd y\dd t
 \le C_*N+C_*X_T+C_*\sqrt{NB_T}.
\]
Multiplying the balance in \eqref{e:hcur} by $w_\sigma$, with
$\partial_yw_\sigma=-\sigma L\delta\omega$, and integrating over $\R$,
we have
\[
 \mathcal H'(t)+\frac{L\delta}{2}\int_\R\omega(q_+^2+q_-^2)\dd y
 =\sum_\sigma\int_\R\Bigl[-\sigma L\delta\,
   \partial_t\tanh\frac y{\sqrt\tau}\,E_\sigma
 +\sigma L\delta\omega J_\sigma+w_\sigma S_\sigma\Bigr]\dd y .
\]
For every $\eta>0$,
$|\partial_t\tanh\frac y{\sqrt\tau}|
=\frac{|y|}{2\tau^{3/2}}\operatorname{sech}^2\frac y{\sqrt\tau}
\le\eta\omega+C_\eta\tau^{-1}e^{-\eta\sqrt\tau}$, the first term
bounding it on $|y|\le2\eta\tau$ and the second on the complement. By
$|J_\sigma|\le C_*|q_\sigma||\Gamma|$, $\omega\le1$ and the Cauchy
inequality,
\[
 L\delta\sum_\sigma\int_{S_0}^T\!\!\int_\R\omega|J_\sigma|\dd y\dd t
 \le\frac L8X_T+C_*L\delta\int_{S_0}^T\Vert\Gamma\Vert_{L^2}^2\dd t .
\]
Choose $L\ge12C_*$, then $\eta\le\frac1{8C_*}$, and decrease
$\delta_1(M)$ so that $L\delta_1(M)\le\frac12$; $L$ and $\eta$ depend
only on the fixed data. Then $\frac12\le w_\sigma\le\frac32$,
$0\le\mathcal H\le C_*N$, and integration over $(S_0,T)$ gives
\[
 X_T\le C_*\bigl(N+\sqrt{NB_T}\bigr).
\]

Apply Lemma~\ref{l:osp} on $[S_0,T]$ with $x=y$,
$(E_1,F_1,S_1)=(E_+(q_+,m),\frac{q_+^2}2-J_+,S_+)$ and
$(E_2,F_2,S_2)=(E_-(q_-,m),-\frac{q_-^2}2-J_-,S_-)$; these are the
balances in \eqref{e:hcur}, and on $[S_0,T]$ the hypotheses of the
lemma hold by \eqref{e:hcur}, \eqref{e:bgcl} and $|S_\sigma|\le
C_*\mathcal M_\sigma$. For
$I(t)=\int_\R E_+(x,t)\int_x^\infty E_-(z,t)\dd z\dd x$, where
$E_\pm(x,t)=E_\pm(q_\pm,m)(x,t)$, we have $0\le I\le C_*N^2$ and, by
\eqref{e:osp2},
\[
\begin{aligned}
 I(T)-I(S_0)+\frac12\int_{S_0}^T\!\!\int_\R(q_+^2E_-+E_+q_-^2)\dd x\dd t
 ={}&\int_{S_0}^T\!\!\int_\R\Bigl[J_+E_--E_+J_-
 +S_+\int_x^\infty E_-\dd z\\
 &\quad+E_+\int_x^\infty S_-\dd z\Bigr]\dd x\dd t .
\end{aligned}
\]
For each $\sigma$, by \eqref{e:agm},
\[
 \int_\R|J_\sigma|E_{-\sigma}\dd x
 \le C_*\Vert q_{-\sigma}\Vert_{L^\infty}^2
                  \Vert q_\sigma\Vert_{L^2}\Vert\Gamma\Vert_{L^2}
 \le C_*N\Vert(q_{-\sigma})_y\Vert_{L^2}\Vert\Gamma\Vert_{L^2},
\]
whose time integral is at most $C_*N^2$, and
$\int_x^\infty E_\mp\dd z\le C_*N$. With $E_\sigma\ge C_*^{-1}q_\sigma^2$,
$I\ge0$ and the bounds for $|S_\sigma|$ and $X_T$,
\[
 B_T\le C_*N^2+C_*N\bigl(N+X_T+\sqrt{NB_T}\bigr)
 \le C_*(N^2+N^3)+\frac12B_T .
\]
Since $B_T<\infty$, $B_T\le C_*(N^2+N^3)$ for every $T$. Hence
$B_\infty+X_\infty<\infty$, $\mathcal M_\sigma\in L^1((S_0,\infty)\times\R)$
and $S_\sigma\in L^1((S_0,\infty)\times\R)$.

In the coordinate $x=\sigma y$, $\partial_x=\sigma\partial_y$, and the
balances \eqref{e:hcur} and \eqref{e:mbal} become those of
Lemma~\ref{l:osp} with
$(E_1,F_1,S_1)=(E_\sigma(q_\sigma,m),\frac{q_\sigma^2}2-\sigma J_\sigma,S_\sigma)$
and $(E_2,F_2,S_2)=(\frac{e_m^2}2,-\sigma e_mf_e,S_e)$, evaluated at
$y=\sigma x$. For
$\mathcal I_\sigma(t)=\int_\R E_1(x,t)\int_x^\infty E_2(z,t)\dd z\dd x$,
\eqref{e:osp2} gives
\[
\begin{aligned}
 \mathcal I_\sigma(T)-\mathcal I_\sigma(S_0)
 +\frac14\int_{S_0}^T\!\!\int_\R q_\sigma^2e_m^2\dd y\dd t
 \le\int_{S_0}^T\Bigl(&\frac12\int_\R|J_\sigma|e_m^2\dd y
       +\int_\R E_\sigma|e_mf_e|\dd y\\
 &+\frac12\Vert e_m\Vert_{L^2}^2\Vert S_\sigma\Vert_{L^1}
       +C_*\Vert q_\sigma\Vert_{L^2}^2\Vert S_e\Vert_{L^1}\Bigr)\dd t .
\end{aligned}
\]
By \eqref{e:agm} and \eqref{e:hcur},
\[
\begin{gathered}
 \int_\R|J_\sigma|e_m^2\dd y
 \le C_*\Vert e_m\Vert_{L^2}\Vert(e_m)_y\Vert_{L^2}
             \Vert q_\sigma\Vert_{L^2}\Vert\Gamma\Vert_{L^2},\\
 \int_\R E_\sigma|e_mf_e|\dd y
 \le C_*\Vert q_\sigma\Vert_{L^2}\Vert(q_\sigma)_y\Vert_{L^2}
             \Vert e_m\Vert_{L^2}\Vert f_e\Vert_{L^2},\qquad
 0\le\mathcal I_\sigma\le C_*\Vert q_\sigma\Vert_{L^2}^2\Vert e_m\Vert_{L^2}^2,
\end{gathered}
\]
and the time integrals of the right sides are finite. Hence
$\int_{S_0}^\infty\int_\R(q_+^2+q_-^2)e_m^2\dd y\dd t<\infty$.

Lemma~\ref{l:trc} with \eqref{e:acch}, $s=S_0$ and $K=K_*$, now with
$\mathcal M_\sigma\in L^1((S_0,\infty)\times\R)$ and
$\Vert q_\sigma(t)\Vert_{L^\infty}\le C_b\Vert Z_m(t)\Vert_{L^\infty}\to0$,
gives $\Vert q_\sigma(t)\Vert_{L^2}\to0$. Since
\[
 \upsilon_m=\frac{m^{1/4}}2(q_++q_-),\qquad
 J_a=\frac{m^{1/4}}2(q_+-q_-),\qquad
 V_b=m\Bigl(\frac{P_b}{p_*}\Bigr)^{-1/\gamma},\qquad
 \Theta_b=\frac{P_bV_b}R,
\]
we have $|Z_m|\le C_b(|q_+|+|q_-|+|e_m|)$, so
$\Vert Z_m(t)\Vert_{L^2}\to0$. By \eqref{e:dq},
$|(q_\pm)_y|+|(e_m)_y|\le C_b(|\partial_yZ_m|+|(V_c)_y||Z_m|)$, so
$\Vert(e_m,q_+,q_-)(t)\Vert_{H^1}\to0$ by Lemma~\ref{l:bg}; and
$|f_m|\le C_b(|(\vartheta_m)_y|+|(\Theta_c)_y|)$, \eqref{e:agm},
Lemma~\ref{l:bg} and \eqref{eq:pr-gau} give
$\Vert f_m(t)\Vert_{L^\infty}\to0$. With the limit of
Lemma~\ref{l:bg}, $\sup_{t\ge r}\Vert Z_m(t)\Vert_{H^3}\to0$; the
remaining quantities are tails of integrals shown to be finite above
and in Lemma~\ref{l:bg} (with $t_0=1$), the one with
$\operatorname{sech}^2$ by $X_\infty<\infty$.
\end{proof}

\section{Conserved integrals, diffusion waves and scalar parabolic tools}\label{s:tools}

For each fixed datum choose an auxiliary $M\ge A$ such that
\begin{equation}\label{e:dat-aux}
\begin{gathered}
 \Vert e_\rho\Vert_{H^3(\R)}
 +\Vert(e_u,e_\theta)\Vert_{H^2(\R)}
 +\Vert e_0\Vert_{L^1(\R)}\le M,\\
 \rho_0\ge M^{-1},\qquad \theta_0\ge M^{-1}.
\end{gathered}
\end{equation}
Such an $M$ exists by \eqref{e:dat} and the prescribed regularity.
It need not be controlled by $A$.

Throughout this section $M\ge1$, \eqref{e:dat-aux} holds, $0<\delta\le\delta_1(M)$
with $\delta_1(M)$ as in Section~\ref{s:bg}, and $(V_b,U_b,\Theta_b)$ is the
solution of Lemma~\ref{l:bg}. The index $\sigma$ belongs to $\{-1,+1\}$, and
$v_\sigma,c_\sigma,s_\sigma,m_\sigma,d_\sigma,k_\sigma$ denote
$v_\pm,c_\pm,s_\pm,m_\pm,d_\pm,k_\pm$ for $\sigma=\pm1$.

\subsection{Conserved integrals and the phase}\label{s:mass}

Put
\[
 E_m=c_v\Theta_b+\frac{b_m^2}2,\qquad \mathcal U=(V_b,b_m,E_m),\qquad
 \widehat{\mathcal U}=(V_c,H_c,c_v\Theta_c).
\]
By \eqref{eq:bg-m}, \eqref{eq:bg-mom}, and \eqref{eq:bg-en} minus $u_*$ times
\eqref{eq:bg-mom},
\[
\begin{gathered}
 (V_b)_t-(b_m)_y=0,\qquad
 (b_m)_t+(P_b)_y-\nu\partial_y\frac{(b_m)_y}{V_b}=0,\\
 (E_m)_t+\partial_y(P_bb_m)
   -\partial_y\frac{\kappa(\Theta_b)_y+\nu b_m(b_m)_y}{V_b}=0,
\end{gathered}
\]
where $\Theta_b=\frac1{c_v}(E_m-\frac{b_m^2}2)$ and $P_b=\frac{R\Theta_b}{V_b}$.
With $\widehat\theta=\Theta_c-\frac{H_c^2}{2c_v}$ and
$\widehat p=\frac{R\widehat\theta}{V_c}=p_*-\frac{(\gamma-1)H_c^2}{2V_c}$,
\eqref{eq:pr-id} gives
\begin{equation}\label{e:Uhat}
\begin{gathered}
 (V_c)_t-(H_c)_y=0,\qquad
 (H_c)_t+\widehat p_y-\nu\partial_y\frac{(H_c)_y}{V_c}=\partial_y\widehat R_2,\\
 (c_v\Theta_c)_t+\partial_y(\widehat pH_c)
   -\partial_y\frac{\kappa\widehat\theta_y+\nu H_c(H_c)_y}{V_c}
 =\partial_y\widehat R_3,\\
 \widehat R_2=\frac{(c-\nu)(H_c)_y}{V_c}-\frac{(\gamma-1)H_c^2}{2V_c},\qquad
 \widehat R_3=\Bigl(\frac\kappa{c_v}-\nu\Bigr)\frac{H_c(H_c)_y}{V_c}
   -\frac{(\gamma-1)H_c^3}{2V_c}.
\end{gathered}
\end{equation}
Put
\begin{equation}\label{e:mass}
\begin{gathered}
 v_\pm=\frac1{\rho_\pm},\qquad \Delta_v=v_+-v_-,\qquad
 Q=(Q_v,Q_{\rm mom},Q_E)=\int_\R(\mathcal U-\widehat{\mathcal U})(y,0)\dd y,\\
 c_\pm=\sqrt{\frac{\gamma p_*}{v_\pm}},\qquad s_\pm=\pm c_\pm,\qquad
 r_0=\Bigl(1,0,\frac{p_*}{\gamma-1}\Bigr),\qquad r_\sigma=(-1,s_\sigma,p_*),\\
 m_0=\frac{\gamma-1}\gamma\Bigl(Q_v+\frac{Q_E}{p_*}\Bigr),\qquad h_*=\frac{m_0}{\Delta_v},\\
 m_-=\frac{c_+(m_0-Q_v)-Q_{\rm mom}}{c_-+c_+},\qquad
 m_+=\frac{c_-(m_0-Q_v)+Q_{\rm mom}}{c_-+c_+}.
\end{gathered}
\end{equation}

\begin{lemma}\label{l:cmp}
For every $T_e>0$,
\[
 Z_m,\quad b_m,\quad (\gamma-1)E_m-p_*V_b,\quad \mathcal U-\widehat{\mathcal U}
 \ \in C([0,T_e];L^1),
\]
and $\int_\R(\mathcal U-\widehat{\mathcal U})(y,t)\dd y=Q$ for every $t\ge0$.
The numbers $h_*,m_-,m_+$ of \eqref{e:mass} are the only real numbers with
\[
 Q=h_*\Delta_vr_0+m_-r_-+m_+r_+,
\]
and
\[
 \delta|h_*|+|m_-|+|m_+|\le C|Q|,
\]
where $C$ depends only on the fixed data. Under \eqref{e:dat},
the right side is at most $CA$.
\end{lemma}
\begin{proof}
Put $\eta=\bar Y(\cdot,0)\circ X_0$. By
\eqref{e:Y0}, \eqref{e:Xc} and \eqref{eq:bg-dat},
$Y_0(x)-\bar Y(x,0)=\int_0^{x}e_\rho\dd z$, hence
$|\eta(y)-y|\le\Vert e_\rho\Vert_{L^1}$ and
\[
 \int_\R|f(\eta(y))-f(y)|\dd y
 \le\int_\R\int_{y-\Vert e_\rho\Vert_{L^1}}^{y+\Vert e_\rho\Vert_{L^1}}|f_y(z)|\dd z\dd y
 =2\Vert e_\rho\Vert_{L^1}\Vert f_y\Vert_{L^1}.
\]
By \eqref{eq:bg-dat},
\[
 V_b(\cdot,0)-V_c(\eta,0)
   =-\Bigl(\frac{e_\rho}{\rho_0\bar\rho(\cdot,0)}\Bigr)\circ X_0,\qquad
 (U_b,\Theta_b)(\cdot,0)-(U_c,\Theta_c)(\eta,0)=(e_u,e_\theta)\circ X_0 .
\]
Changing variables by $y=Y_0(x)$ and applying \eqref{eq:pr-gau} to the profile
difference gives $\Vert Z_m(0)\Vert_{L^1}\le C_M\Vert e_0\Vert_{L^1(\R)}$.
Since $b_m(0)=\upsilon_m(0)+H_c(0)\in L^2$ and
\[
 \mathcal U-\widehat{\mathcal U}
   =\Bigl(\varphi_m,\ \upsilon_m,\ c_v\vartheta_m+\frac{b_m^2}2\Bigr),
\]
the initial conserved excess is integrable and $Q$ is well-defined.

To propagate integrability, fix $0<T_e<\infty$ and a coefficient $a$ on
$\R\times[0,T_e]$ with $a\ge a_*>0$ and bounded $a_y,a_{yy}$.
Let $S_a(t,s)$, $0\le s\le t\le
T_e$, be the solution operator of $f_t=\partial_y(af_y)$. Then, for $p=1,2$,
\[
 \Vert S_a(t,s)f\Vert_{L^p}\le\Vert f\Vert_{L^p},\qquad
 \Vert S_a(t,s)\partial_yf\Vert_{L^p}\le C_a(t-s)^{-\frac12}\Vert f\Vert_{L^p},
\]
where $C_a$ depends only on $T_e$, $a_*$, $\Vert a_y\Vert_{L^\infty}$ and
$\Vert a_{yy}\Vert_{L^\infty}$. The first bound follows from positivity
and mass conservation for $p=1$, and the energy identity for $p=2$.
By duality and time reversal, the second follows from
$\Vert g_y(t)\Vert_{L^q}\le C_a(t-s)^{-\frac12}\Vert g(s)\Vert_{L^q}$,
$q=\infty,2$, for solutions of $\mathcal L_ag=0$ on $(s,t)$, where
$\mathcal L_a=\partial_t-a\partial_y^2-a_y\partial_y$. For $z=g_y$,
\[
\begin{gathered}
 \mathcal L_a(g^2)=-2az^2,\qquad
 \mathcal L_a(z^2)=-2az_y^2+2a_yzz_y+2a_{yy}z^2,\\
 \mathcal L_a\bigl((t-s)z^2+Ag^2\bigr)
 \le-(t-s)az_y^2+\Bigl[1+(t-s)\Bigl(2\Vert a_{yy}\Vert_{L^\infty}
   +\frac{\Vert a_y\Vert_{L^\infty}^2}{a_*}\Bigr)-2Aa_*\Bigr]z^2\le0
\end{gathered}
\]
for $A=\frac1{2a_*}(1+T_e(2\Vert a_{yy}\Vert_{L^\infty}+\frac{\Vert
a_y\Vert_{L^\infty}^2}{a_*}))$, and the maximum principle gives
$(t-s)g_y^2\le A\Vert g(s)\Vert_{L^\infty}^2$. In $L^2$,
\[
 \frac12\frac{\dd}{\dd t}\Vert g\Vert_{L^2}^2+\int_\R ag_y^2\dd y=0,\qquad
 \frac12\frac{\dd}{\dd t}\Vert g_y\Vert_{L^2}^2+\int_\R ag_{yy}^2\dd y
 =\frac12\int_\R a_{yy}g_y^2\dd y;
\]
at a time $t_1\in(s,\frac{s+t}2)$ with $\Vert g_y(t_1)\Vert_{L^2}^2\le
\frac2{t-s}\int_s^{(s+t)/2}\Vert g_y\Vert_{L^2}^2\le\frac{\Vert
g(s)\Vert_{L^2}^2}{a_*(t-s)}$, Gronwall's inequality on $[t_1,t]$ gives the
$L^2$ bound. Galerkin approximation constructs $S_a$ and the adjoint
evolution.

Put $a_1=\frac\nu{V_b}$ and $a_2=\frac\kappa{c_vV_b}$. By \eqref{e:bgpar},
\eqref{eq:pr-id} and $P_b-p_*=\frac{R\vartheta_m-p_*\varphi_m}{V_b}$,
\[
\begin{gathered}
 (\varphi_m)_t=(\upsilon_m)_y,\qquad
 (\upsilon_m)_t=\partial_y(a_1(\upsilon_m)_y)-\mathfrak f_y,\qquad
 (\vartheta_m)_t=\partial_y(a_2(\vartheta_m)_y)-\mathfrak g_y+\mathfrak s,\\
 \mathfrak f=\frac{R\vartheta_m-p_*\varphi_m}{V_b}
   +\frac{\nu(H_c)_y\varphi_m}{V_bV_c}+\frac{(c-\nu)(H_c)_y}{V_c},\qquad
 \mathfrak g=\frac{P_b\upsilon_m}{c_v}
   +\frac{\kappa(\Theta_c)_y\varphi_m}{c_vV_bV_c},\\
 \mathfrak s=\frac{(P_b)_y\upsilon_m-(P_b-p_*)(H_c)_y}{c_v}
   +\frac{\nu(U_b)_y^2}{c_vV_b}.
\end{gathered}
\]
Let $\mathcal P_B=(1-\partial_y^2)^{-1}$ and $\mathcal J_B=\partial_y\mathcal
P_B$; their kernels $\frac12e^{-|y|}$ and $-\frac12\operatorname{sgn}(y)e^{-|y|}$
have $L^1$ norm $1$, and $\mathcal J_B\partial_y=\mathcal P_B-I$. With
$(a_1^{-1})_t=\frac{(U_b)_y}\nu$, the function $h_v=\varphi_m+\mathcal
J_B(a_1^{-1}\upsilon_m)$ satisfies
\[
 (h_v)_t=\mathcal J_B\upsilon_m
   +(\mathcal P_B-I)\Bigl(\frac{(a_1)_y}{a_1}\upsilon_m-\frac{\mathfrak f}{a_1}\Bigr)
   +\mathcal J_B\Bigl(-\partial_y\Bigl(\frac{(a_1)_y}{a_1}\Bigr)\upsilon_m
     +(a_1^{-1})_y\mathfrak f+\frac{(U_b)_y}\nu\upsilon_m\Bigr).
\]
By Lemma~\ref{l:bg}, \eqref{e:bgcl} and \eqref{eq:pr-gau}, the functions
$a_i,(a_i)_y,(a_i)_{yy},(U_b)_y,P_b,(P_b)_y$ and the coefficients of
$\varphi_m,\vartheta_m$ in $\mathfrak f$, $\mathfrak g$ are bounded on
$\R\times[0,T_e]$. Regard them as fixed and insert
$\varphi_m=h_v-\mathcal J_B(a_1^{-1}\upsilon_m)$. On an interval $[s,s+h]$, the
$h_v$ equation integrated in time, and
\[
 \upsilon_m(t)=S_{a_1}(t,s)\upsilon_m(s)-\int_s^tS_{a_1}(t,t')\mathfrak f_y\dd t',
 \qquad
 \vartheta_m(t)=S_{a_2}(t,s)\vartheta_m(s)
   +\int_s^tS_{a_2}(t,t')(\mathfrak s-\mathfrak g_y)\dd t',
\]
form an affine system for $(h_v,\upsilon_m,\vartheta_m)$. By the bounds for
$S_{a_1},S_{a_2}$, its linear part has norm at most $C_{b,T_e}(h+2\sqrt h)$ on
$C([s,s+h];L^1\cap L^2)$ and on $C([s,s+h];L^2)$; its source terms
$\frac{(c-\nu)(H_c)_y}{V_c}$ and $\frac{\nu(U_b)_y^2}{c_vV_b}$ belong to
$C([0,T_e];L^1\cap L^2)$ and $L^1(0,T_e;L^1\cap L^2)$. For
$C_{b,T_e}(h+2\sqrt h)\le\frac12$ and data at $s$ in $L^1\cap L^2$, the system
has a fixed point in
$C([s,s+h];L^1\cap L^2)$, which by uniqueness in $C([s,s+h];L^2)$ is
$(h_v,\upsilon_m,\vartheta_m)$. Finitely many such intervals cover $[0,T_e]$;
thus $h_v,\upsilon_m,\vartheta_m,\varphi_m\in C([0,T_e];L^1)$. The formula for
$\mathcal U-\widehat{\mathcal U}$, $b_m=\upsilon_m+H_c$,
\[
 (\gamma-1)E_m-p_*V_b=R\vartheta_m-p_*\varphi_m+\frac{\gamma-1}2b_m^2,
\]
$b_m\in C([0,T_e];L^2)$ and $\Vert H_c(t)\Vert_{L^1}=c|\log\frac{v_+}{v_-}|$,
by \eqref{eq:prof}, give the other memberships.

Subtracting \eqref{e:Uhat} from the equations of $\mathcal U$ gives
$(\mathcal U-\widehat{\mathcal U})_t=\partial_y\mathcal J$ with
\[
 \mathcal J=\Bigl(\upsilon_m,\
   \frac{\nu(b_m)_y}{V_b}-(P_b-p_*)-\frac{c(H_c)_y}{V_c},\
   \frac{\kappa(\Theta_b)_y+\nu b_m(b_m)_y}{V_b}-\frac{\kappa(\Theta_c)_y}{V_c}
   -(P_b-p_*)b_m-p_*\upsilon_m\Bigr).
\]
By \eqref{e:bgcl} and \eqref{eq:pr-gau},
$\sup_{0\le t\le T_e}\Vert\mathcal J(t)\Vert_{L^2}\le C_{b,T_e}$.
With $\chi_e$ from \eqref{e:gep} and $\chi_K=\chi_e(\frac\cdot K)$,
\[
 \Bigl|\int_\R\chi_K\bigl[(\mathcal U-\widehat{\mathcal U})(y,t)
   -(\mathcal U-\widehat{\mathcal U})(y,0)\bigr]\dd y\Bigr|
 \le\frac{\Vert\chi_e'\Vert_{L^2}}{\sqrt K}\int_0^t\Vert\mathcal J(t')\Vert_{L^2}\dd t'
 \le\frac{C_{b,T_e}}{\sqrt K}.
\]
Letting $K\to\infty$ gives $\int_\R(\mathcal U-\widehat{\mathcal U})(y,t)\dd y=Q$.

In components, $Q=h_*\Delta_vr_0+m_-r_-+m_+r_+$ with $m_0=h_*\Delta_v$ reads
\[
 m_0-m_--m_+=Q_v,\qquad -c_-m_-+c_+m_+=Q_{\rm mom},\qquad
 \frac{p_*m_0}{\gamma-1}+p_*(m_-+m_+)=Q_E,
\]
whose only solution is \eqref{e:mass}: $\det(r_0,r_-,r_+)=-\frac{\gamma
p_*(c_-+c_+)}{\gamma-1}\ne0$ and $|\Delta_v|=\frac{R\delta}{p_*}$. Hence
$\delta|h_*|+|m_-|+|m_+|\le C|Q|$, as asserted.
\end{proof}

\subsection{Diffusion waves and the conjugation}\label{s:wave}

For $n>0$ put
\begin{equation}\label{e:ac}
 d(n)=\frac{\nu+(\gamma-1)c}{2n},\qquad
 s_\sigma(n)=\sigma\sqrt{\frac{\gamma p_*}n},\qquad
 k(n)=\frac{\gamma+1}{4n},\qquad
 \alpha=\frac{\gamma+1}{2(\nu+(\gamma-1)c)},
\end{equation}
so that $\alpha d(n)=k(n)$, and let $d_\sigma=d(v_\sigma)$, $k_\sigma=k(v_\sigma)$;
then $s_\sigma(v_\sigma)=s_\sigma$. For real numbers $m_-,m_+$, use the
Cole--Hopf representation~\cite{Cole1951,Hopf1950}:
\begin{equation}\label{e:CH}
\begin{gathered}
 F_d(\xi)=\frac1{\sqrt{4\pi d}}\int_{-\infty}^\xi e^{-\frac{z^2}{4d}}\dd z,\qquad
 p_\sigma=1+\bigl(e^{-\alpha s_\sigma m_\sigma}-1\bigr)
   F_{d_\sigma}\Bigl(\frac{y-s_\sigma t}{\sqrt\tau}\Bigr),\\
 A_\sigma=-\frac1\alpha\log p_\sigma,\qquad a_\sigma=\partial_yA_\sigma,\qquad
 \vartheta_\sigma=\frac{a_\sigma}{s_\sigma},\qquad
 R_U=\mathcal U-\widehat{\mathcal U}(y+h_*,t)
   -\sum_{\sigma=\pm1}\vartheta_\sigma r_\sigma .
\end{gathered}
\end{equation}
In $R_U$, and wherever no other masses are fixed, $m_\pm$ and $h_*$ are those
of \eqref{e:mass}. The scalar operators are
\begin{equation}\label{e:ops}
 L_\sigma=\partial_t-d(V_c)\partial_y^2+s_\sigma(V_c)\partial_y,\qquad
 P_\sigma=\partial_t-d_\sigma\partial_y^2+s_\sigma\partial_y,\qquad
 \mathscr P_\sigma=L_\sigma+2k(V_c)a_\sigma\partial_y .
\end{equation}

\begin{lemma}\label{l:pr}
Let $M_a\ge1$ and $|m_-|+|m_+|\le M_a$. Then $p_\sigma$ takes values between
$1$ and $e^{-\alpha s_\sigma m_\sigma}$, $A_\sigma$ between $0$ and $s_\sigma m_\sigma$,
\[
\begin{gathered}
 A_\sigma(-\infty,t)=0,\qquad A_\sigma(+\infty,t)=s_\sigma m_\sigma,\qquad
 \int_\R\vartheta_\sigma(y,t)\dd y=m_\sigma,\\
 P_\sigma p_\sigma=0,\qquad P_\sigma A_\sigma=-k_\sigma a_\sigma^2,\qquad
 (\vartheta_\sigma)_t+s_\sigma(\vartheta_\sigma)_y
   +2k_\sigma s_\sigma\vartheta_\sigma(\vartheta_\sigma)_y
 =d_\sigma(\vartheta_\sigma)_{yy},\\
 |\partial_y^ja_\sigma|+|\partial_y^j\vartheta_\sigma|
 \le C_j\tau^{-\frac{j+1}2}e^{-\frac{(y-s_\sigma t)^2}{8d_\sigma\tau}}
 \qquad(j\ge0),
\end{gathered}
\]
where $C_j$ depends only on $j$, $M_a$ and the fixed data. For every smooth
$z$,
\[
 L_\sigma(p_\sigma z)=p_\sigma\bigl(\mathscr P_\sigma z-\alpha R^{\rm prof}_\sigma z\bigr),
 \qquad
 P_\sigma(p_\sigma z)=p_\sigma\bigl(P_\sigma z+2k_\sigma a_\sigma z_y\bigr),
\]
where
\[
 R^{\rm prof}_\sigma=L_\sigma A_\sigma+k(V_c)a_\sigma^2
 =(s_\sigma(V_c)-s_\sigma)a_\sigma-(d(V_c)-d_\sigma)(a_\sigma)_y
   +(k(V_c)-k_\sigma)a_\sigma^2 .
\]
With $c_G$ from \eqref{eq:ct-gau} and $C$ depending only on $M_a$ and the
fixed data,
\[
 \Vert R^{\rm prof}_\sigma(t)\Vert_{L^\infty}
   +\Vert(V_c)_ya_\sigma(t)\Vert_{L^\infty}\le C\delta e^{-c_et},\qquad
 c_e=\frac{\min\{c_-,c_+\}^2}4\min\Bigl\{c_G,\frac1{8d_-},\frac1{8d_+}\Bigr\}.
\]
The masses of \eqref{e:mass} satisfy $|m_-|+|m_+|\le C_M$; for them,
$R_U(t)\in L^1$ and $\int_\R R_U(y,t)\dd y=0$ for every $t\ge0$.
\end{lemma}
\begin{proof}
In $x=y-s_\sigma t$, $F_{d_\sigma}(x/\sqrt\tau)$ solves the heat equation,
so $P_\sigma p_\sigma=0$. Since $F_{d_\sigma}$ increases from $0$ to $1$,
\eqref{e:CH} gives the ranges and endpoint values, as well as
$\int_\R a_\sigma\dd y=s_\sigma m_\sigma$. By
$a_\sigma=-\frac{(p_\sigma)_y}{\alpha p_\sigma}$ and $\alpha d_\sigma=k_\sigma$,
\[
 P_\sigma A_\sigma
 =-\frac{P_\sigma p_\sigma}{\alpha p_\sigma}
   -\frac{d_\sigma(p_\sigma)_y^2}{\alpha p_\sigma^2}
 =-k_\sigma a_\sigma^2 ;
\]
its $y$-derivative divided by $s_\sigma$ is the equation for
$\vartheta_\sigma$. For $k\ge1$, $\partial_y^kp_\sigma=(e^{-\alpha s_\sigma
m_\sigma}-1)\tau^{-\frac k2}F_{d_\sigma}^{(k)}(\frac x{\sqrt\tau})$ and
$|F_d^{(k)}(\xi)|\le C_ke^{-\frac{\xi^2}{8d}}$; with $p_\sigma\ge
e^{-\alpha|s_\sigma m_\sigma|}$ and $a_\sigma=-\frac1\alpha\partial_y\log
p_\sigma$ this gives the bounds for $\partial_y^ja_\sigma$ and
$\partial_y^j\vartheta_\sigma$.

By the product rule,
\[
 L_\sigma(p_\sigma z)=p_\sigma L_\sigma z-2d(V_c)(p_\sigma)_yz_y+zL_\sigma p_\sigma,
\]
where $(p_\sigma)_y=-\alpha a_\sigma p_\sigma$, $\alpha d(V_c)=k(V_c)$ and
$L_\sigma p_\sigma=-\alpha p_\sigma(L_\sigma A_\sigma+\alpha d(V_c)a_\sigma^2)$;
this is the first identity. The second is the same computation for
$P_\sigma$, with $P_\sigma p_\sigma=0$. Inserting
$L_\sigma=P_\sigma-(d(V_c)-d_\sigma)\partial_y^2+(s_\sigma(V_c)-s_\sigma)\partial_y$
and $P_\sigma A_\sigma=-k_\sigma a_\sigma^2$ into $L_\sigma A_\sigma$ gives the
formula for $R^{\rm prof}_\sigma$.

The formula for $R^{\rm prof}_\sigma$ gives
\[
\begin{aligned}
 |R^{\rm prof}_\sigma|+|(V_c)_ya_\sigma|
 &\le C\bigl(|V_c-v_\sigma|+|(V_c)_y|\bigr)\\
 &\qquad\times\bigl(|a_\sigma|+|(a_\sigma)_y|+a_\sigma^2\bigr).
\end{aligned}
\]
On $\sigma y\ge0$, the first factor is bounded by
$C\delta e^{-c_Gy^2/\tau}$ by \eqref{eq:ct-gau}. Combine its Gaussian with that
of $a_\sigma$ using $y^2+(y-s_\sigma t)^2\ge s_\sigma^2t^2/2$.
On $\sigma y<0$, use $|V_c-v_\sigma|+|(V_c)_y|\le C\delta$ and
$|y-s_\sigma t|\ge|s_\sigma|t$. Since $t^2/\tau\ge t-1$, both bounds give
the stated rate $c_e$.

Lemma~\ref{l:cmp} and \eqref{e:dat-aux} give $|m_-|+|m_+|\le C_M$. Since
$\widehat{\mathcal U}(\pm\infty,t)=(v_\pm,0,c_v\theta_\pm)$ and
$c_v(\theta_+-\theta_-)=\frac{p_*\Delta_v}{\gamma-1}$, Fubini's theorem gives
\[
 \int_\R[\widehat{\mathcal U}(y+h_*,t)-\widehat{\mathcal U}(y,t)]\dd y
 =\int_0^{h_*}\!\!\int_\R\partial_y\widehat{\mathcal U}(y+z,t)\dd y\dd z
 =h_*\Delta_vr_0 .
\]
With Lemma~\ref{l:cmp} and $\int_\R\vartheta_\sigma\dd y=m_\sigma$, this gives
$R_U(t)\in L^1$ and $\int_\R R_U\dd y=Q-h_*\Delta_vr_0-m_-r_--m_+r_+=0$.
\end{proof}

\subsection{Positive responses and reception}\label{s:resp}

Let $0\le r<T_e<\infty$ and let $\ell_-<\ell_+$ be affine functions on
$[r,T_e]$. A finite domain is a set
\[
 \Omega=\{(y,t):r<t<T_e,\ \ell_-(t)<y<\ell_+(t)\};
\]
its parabolic boundary $\partial_p\Omega$ consists of the bottom $\{t=r\}$ and
the two lateral sides. On $\Omega$ let $P=\partial_t-d\partial_y^2+B\partial_y$,
where $d$ is Lipschitz on $\overline\Omega$, $d\ge d_{\min}>0$, and $B$ is
measurable with $\int_r^{T_e}\Vert B(t)\Vert_{L^\infty}^2\dd t<\infty$. For
$f\in L^2(\Omega)$, and for $z$ continuous on $\overline\Omega$ with $z_y$
bounded and $z_t,z_{yy}\in L^2(\Omega)$ (so that $Pz\in L^2(\Omega)$), define
\begin{equation}\label{e:Tdef}
\begin{gathered}
 P\,T_{P,\Omega}f=f\ \text{ in }\Omega,\qquad
 T_{P,\Omega}f=0\ \text{ on }\partial_p\Omega,\\
 H_{P,\Omega}[z|_{\partial_p\Omega}]=z-T_{P,\Omega}(Pz),\qquad
 T_{P,r}f=\lim_{n\to\infty}T_{P,(-n,n)\times(r,r+n)}f\quad(f\ge0).
\end{gathered}
\end{equation}
For other $f\ge0$, $T_{P,\Omega}f=\lim_{n\to\infty}T_{P,\Omega}\min\{f,n\}$. A
domain is a finite domain or the union $\Omega$ of an increasing sequence of
finite domains $\Omega_n$, with coefficients as above on each $\Omega_n$; then
$T_{P,\Omega}f=\lim_{n\to\infty}T_{P,\Omega_n}f$ for $f\ge0$. For real $f$,
$T_{P,\Omega}f=T_{P,\Omega}f^+-T_{P,\Omega}f^-$ when both are finite. By the
comparison below these limits are increasing, so a bound for nonnegative
sources that holds with one constant on every finite domain also holds for
the limits.

The first line of \eqref{e:Tdef} is meant in the weak sense. If $\ell_\pm$ are
constant and $I=(\ell_-,\ell_+)$, then $u=T_{P,\Omega}f$ is the function in
$C([r,T_e];L^2(I))\cap L^2(r,T_e;H^1_0(I))$ with $u(r)=0$ and
\[
 \iint_\Omega\Bigl(-u\,\partial_t\frac\phi d+u_y\phi_y+\frac Bdu_y\phi\Bigr)
   \dd y\dd t
 =\iint_\Omega\frac{f\phi}d\dd y\dd t
\]
for every $\phi\in H^1(\Omega)$ vanishing on the lateral sides and at $t=T_e$;
this is $Pu=f$ divided by $d$. The map $(y,t)\mapsto(\frac{y-\ell_-(t)}
{\ell_+(t)-\ell_-(t)},t)$ takes a finite domain onto $(0,1)\times(r,T_e)$ and
$P$ to an operator of the same form. Testing with $u$ gives
\[
 \frac12\frac{\dd}{\dd t}\int_I\frac{u^2}d\dd y+\int_Iu_y^2\dd y
 =\int_I\frac{fu}d\dd y-\int_I\frac Bdu_yu\dd y
   +\frac12\int_I\partial_t\Bigl(\frac1d\Bigr)u^2\dd y,
\]
hence $\frac{\dd}{\dd t}\int_I\frac{u^2}d+\int_Iu_y^2\le C_\Omega(1+\Vert
B(t)\Vert_{L^\infty}^2)\int_I\frac{u^2}d+C_\Omega\Vert f(t)\Vert_{L^2}^2$, with
$C_\Omega$ depending on $d_{\min}$, $\sup d$ and the Lipschitz constant of $d$.
Galerkin approximation gives existence, and the same estimate gives uniqueness.
For $f\ge0$, testing with $-u^-$ gives
$u\ge0$. Testing with $(u-\mathcal W)^+$ gives $u\le\mathcal W$ for every weak
supersolution $\mathcal W$ of $P\mathcal W\ge f$ with $\mathcal W\ge0$ on
$\partial_p\Omega$; every $\mathcal W$ with the regularity of $z$ above and
$P\mathcal W\ge f$ almost everywhere is one. Consequently
$|T_{P,\Omega}f|\le T_{P,\Omega}|f|$, $H_{P,\Omega}[g]$ depends only on $g$,
$H_{P,\Omega}[1]=1$ and $|H_{P,\Omega}[g]|\le\sup_{\partial_p\Omega}|g|$.

Let $\Omega=\R\times(r,T_e)$, let $d$ and $B$ be bounded, and let $u$ be
bounded, with the regularity of $z$ above on each
$\Omega_n=(-n,n)\times(r,T_e)$, $u=0$ at $t=r$ and $Pu=f$ almost everywhere,
where $f$ is bounded. Then $u=T_{P,\Omega}f$. Indeed,
\[
 W_n=\frac{\Vert u\Vert_{L^\infty}}{n^2}
   \bigl(y^2+2(\sup d+n\sup|B|)(t-r)\bigr)
\]
satisfies $PW_n\ge0$ in $\Omega_n$ and $W_n\ge|u-T_{P,\Omega_n}f|$ on
$\partial_p\Omega_n$, so testing with $(\pm(u-T_{P,\Omega_n}f)-W_n)^+$ gives
$|u-T_{P,\Omega_n}f|\le W_n$. As $n\to\infty$, $W_n\to0$ and
$T_{P,\Omega_n}f\to T_{P,\Omega}f$, since $T_{P,\Omega_n}f^\pm$ increase in $n$
and are at most $\Vert f\Vert_{L^\infty}(t-r)$.

\begin{lemma}\label{l:pos}
Let $\Omega\subset\R\times(r,\infty)$ be a finite domain and let $P,P'$ be
operators as in \eqref{e:Tdef}.
\begin{enumerate}[label=(\roman*)]
\item For $f,g\ge0$ and $h=h(t)\ge0$,
 $T_{P,\Omega}(fg)\le(T_{P,\Omega}f^2)^{\frac12}(T_{P,\Omega}g^2)^{\frac12}$ and
 $T_{P,\Omega}h\le\int_r^\infty h(t)\dd t$.
\item If $z$ has the regularity in \eqref{e:Tdef},
 $|z|\le K$, $\mu>0$, $F\ge0$ and $Pz+\mu z_y^2\le F$, then
 $\mu T_{P,\Omega}z_y^2\le\Vert T_{P,\Omega}F\Vert_{L^\infty(\Omega)}+2K$.
\item If $\int_r^{T_e}\Vert w(t)\Vert_{L^\infty}^2\dd t<\infty$ and
 $\Vert T_{P,\Omega}w^2\Vert_{L^\infty(\Omega)}\le\frac{d_{\min}}{16}$, then
 $\Vert T_{P+w\partial_y,\Omega}f\Vert_{L^\infty(\Omega)}
 \le2\Vert T_{P,\Omega}f\Vert_{L^\infty(\Omega)}$ for $f\ge0$.
\item Let $p>0$ on $\overline\Omega$ with $p,p_t,p_y,p_{yy}$ bounded and
 continuous, and let $V$ be measurable with
 $\int_r^{T_e}\Vert V(t)\Vert_{L^\infty}\dd t<\infty$ and $P(pz)=p(P'z+Vz)$ for
 every smooth $z$. Then, for $f\ge0$, each of
 $\Vert T_{P,\Omega}f\Vert_{L^\infty(\Omega)}$ and
 $\Vert T_{P',\Omega}f\Vert_{L^\infty(\Omega)}$ is at most
 $\frac{\sup_\Omega p}{\inf_\Omega p}\exp(\int_r^{T_e}\Vert V(t)\Vert_{L^\infty}\dd t)$
 times the other.
\item For the waves of Lemma~\ref{l:pr} and $f\ge0$, with $C$ depending only
 on $M_a$ and the fixed data,
 \[
\begin{gathered}
  \Vert T_{P_\sigma+2k_\sigma a_\sigma\partial_y,\Omega}f\Vert_{L^\infty(\Omega)}
  \le e^{\alpha|s_\sigma m_\sigma|}\Vert T_{P_\sigma,\Omega}f\Vert_{L^\infty(\Omega)},
  \\
  \Vert T_{\mathscr P_\sigma,\Omega}f\Vert_{L^\infty(\Omega)}
  \le e^{\alpha|s_\sigma m_\sigma|+C\delta}
    \Vert T_{L_\sigma,\Omega}f\Vert_{L^\infty(\Omega)},
 \end{gathered}
 \]
 and the same holds with the two operators of each inequality interchanged.
\end{enumerate}
\end{lemma}
\begin{proof}
For $\varepsilon>0$, $2fg\le\varepsilon f^2+\varepsilon^{-1}g^2$ and positivity
give $2T_{P,\Omega}(fg)\le\varepsilon T_{P,\Omega}f^2+\varepsilon^{-1}T_{P,\Omega}g^2$;
the infimum over rational $\varepsilon$ at each point gives the first bound in
(i). The second is comparison with $\int_r^th(t')\dd t'$. For (ii),
$z=H_{P,\Omega}[z|_{\partial_p\Omega}]+T_{P,\Omega}(Pz)$ and $Pz\le F-\mu z_y^2$ give
$\mu T_{P,\Omega}z_y^2\le T_{P,\Omega}F+H_{P,\Omega}[z|_{\partial_p\Omega}]-z
\le\Vert T_{P,\Omega}F\Vert_{L^\infty}+2K$.

For (iii) let $v=T_{P+w\partial_y,\Omega}f\ge0$ and $T=T_{P,\Omega}$.
By $Pv=f-wv_y$ and $P(v^2)=2vf-2vwv_y-2dv_y^2$, with zero values on
$\partial_p\Omega$, $v=T(f-wv_y)$ and $2d_{\min}Tv_y^2\le T(2vf-2vwv_y)$. By (i),
\[
 \max\bigl\{\Vert v\Vert_{L^\infty},
   \sqrt{d_{\min}}\Vert Tv_y^2\Vert_{L^\infty}^{1/2}\bigr\}
 \le\Vert Tf\Vert_{L^\infty}
   +\Vert Tw^2\Vert_{L^\infty}^{1/2}\Vert Tv_y^2\Vert_{L^\infty}^{1/2}.
\]
Since $\Vert Tw^2\Vert_{L^\infty}\le d_{\min}/16$, absorption gives
$\Vert v\Vert_{L^\infty}\le\frac43\Vert Tf\Vert_{L^\infty}$.
For bounded smooth $w$ and $f$ these quantities
are finite; the general case follows by approximation.

For (iv) put $E(t)=\exp\int_r^t\Vert V(t')\Vert_{L^\infty}\dd t'$. If
$w=T_{P,\Omega}f$, the identity with $z=\frac wp$ gives
$P'(\frac wp)=\frac fp-V\frac wp$, hence
\[
 P'\Bigl(\frac{Ew}p\Bigr)=\frac{Ef}p+\bigl(\Vert V\Vert_{L^\infty}-V\bigr)\frac{Ew}p
 \ge\frac f{\sup_\Omega p},
\]
and comparison gives $T_{P',\Omega}f\le(\sup_\Omega p)\frac{Ew}p$. If
$u=T_{P',\Omega}f$, then $P(pEu)=pE(f+(V+\Vert V\Vert_{L^\infty})u)\ge(\inf_\Omega
p)f$, so $T_{P,\Omega}f\le\frac{pEu}{\inf_\Omega p}$. By Lemma~\ref{l:pr},
(iv) applies to $(P,P',p,V)=(P_\sigma,P_\sigma+2k_\sigma a_\sigma\partial_y,
p_\sigma,0)$ and $(L_\sigma,\mathscr P_\sigma,p_\sigma,-\alpha R^{\rm prof}_\sigma)$,
with $\frac{\sup p_\sigma}{\inf p_\sigma}\le e^{\alpha|s_\sigma m_\sigma|}$ and
$\alpha\int_0^\infty\Vert R^{\rm prof}_\sigma\Vert_{L^\infty}\dd t\le
\frac{\alpha C\delta}{c_e}$; this is (v).
\end{proof}

\begin{lemma}\label{l:gsrc}
Let $v,h\in\R$, $b_0>0$, $0<d_{\min}\le d_{\max}$, and let $P$ be as in
\eqref{e:Tdef} on a domain $\Omega\subset\R\times(r,\infty)$, with
$d_{\min}\le d\le d_{\max}$ and $\sigma(B-v)\ge b_0$. Put $x=\sigma(y-vt+h)$,
$\xi=\frac x{\sqrt\tau}$ and $\tau_r=1+r$. For nonnegative $g\in L^1(\R)$,
\begin{equation}\label{e:stcmp}
\begin{gathered}
 p_g(x)=\frac1{d_{\min}}\int_x^\infty e^{-\frac{b_0}{d_{\max}}(z-x)}g(z)\dd z,
 \qquad
 \mathcal V_g(x)=\int_{-\infty}^xp_g(z)\dd z,\\
 P\mathcal V_g=\Bigl(\sigma(B-v)-\frac{db_0}{d_{\max}}\Bigr)p_g+\frac d{d_{\min}}g
 \ge g,\qquad
 0\le\mathcal V_g\le\frac{d_{\max}}{b_0d_{\min}}\Vert g\Vert_{L^1},
\end{gathered}
\end{equation}
and $T_{P,\Omega}[g(x)]\le\mathcal V_g(x)$. For $p_w>\frac12$, $a>0$ and
$N\ge0$,
\begin{equation}\label{e:gauss}
\begin{gathered}
 \Vert T_{P,\Omega}(\tau^{-\frac12}(1+\xi^2)^{-p_w})\Vert_{L^\infty(\Omega)}\le C,
 \qquad
 \Vert T_{P,\Omega}(\tau^{-\frac12}e^{-a\xi^2})\Vert_{L^\infty(\Omega)}\le C,\\
 \Vert T_{P,\Omega}(\tau^{-1}(1+|\xi|)^Ne^{-a\xi^2})\Vert_{L^\infty(\Omega)}
 \le C\tau_r^{-\frac12},
\end{gathered}
\end{equation}
where $C$ depends only on $p_w$ (respectively $a,N$), $b_0$, $d_{\min}$ and
$d_{\max}$.
\end{lemma}
\begin{proof}
In the variables $(x,t)$, $P=\partial_t-d\partial_x^2+\sigma(B-v)\partial_x$.
With $\partial_x\mathcal V_g=p_g$ and $\partial_xp_g=\frac{b_0}{d_{\max}}p_g
-\frac g{d_{\min}}$ this gives $P\mathcal V_g$; the inequality uses
$\sigma(B-v)\ge b_0\ge\frac{db_0}{d_{\max}}$ and $d\ge d_{\min}$, and
$\int_\R p_g\dd x=\frac{d_{\max}}{b_0d_{\min}}\Vert g\Vert_{L^1}$ bounds
$\mathcal V_g$. Comparison, first for bounded $g$, gives $T_{P,\Omega}g\le
\mathcal V_g$.

Let $G=\tau^{-\frac12}(1+\xi^2)^{-p_w}$ and $\mathcal I(x,t)=\int_{-\infty}^x
G(z,t)\dd z$, so that $0\le\mathcal I\le\int_\R(1+\xi^2)^{-p_w}\dd\xi$. Since
$G_t=-\frac1{2\tau}\partial_x(xG)$,
\[
 P\mathcal I=\sigma(B-v)G-\frac{xG}{2\tau}-dG_x
 \ge\frac{b_0}2G-\Bigl(\frac x{2\tau}-\frac{b_0}2\Bigr)_+G-d_{\max}|G_x| .
\]
The positive part vanishes unless $x>b_0\tau$, that is $\xi>b_0\sqrt\tau$, and
then $(\frac x{2\tau}-\frac{b_0}2)_+G\le\frac{\xi(1+\xi^2)^{-p_w}}{2\tau}
\le\frac{b_0^{1-2p_w}}2\tau^{-\frac12-p_w}$. Moreover
$|G_x|=\tau^{-1}|\frac{\dd}{\dd\xi}(1+\xi^2)^{-p_w}|\le\frac{2p_w}{\tau+x^2}
\le g(x)=2p_w\min\{\tau_r^{-1},x^{-2}\}$, with
$\Vert g\Vert_{L^1}=8p_w\tau_r^{-\frac12}$.
By \eqref{e:stcmp}, the nonnegative supersolution for the source $G$
\[
 \frac2{b_0}\Bigl(\mathcal I+d_{\max}\mathcal V_g
   +\frac{b_0^{1-2p_w}}2\int_r^t(1+t')^{-\frac12-p_w}\dd t'\Bigr)
\]
is bounded by $C$. Comparison gives the first bound
in \eqref{e:gauss}. The other two follow by positivity from
$e^{-a\xi^2}\le\max\{1,a^{-1}\}(1+\xi^2)^{-1}$,
$(1+|\xi|)^Ne^{-a\xi^2}\le C_{N,a}(1+\xi^2)^{-1}$ and
$\tau^{-1}\le\tau_r^{-\frac12}\tau^{-\frac12}$ on $\Omega$.
\end{proof}

\begin{lemma}\label{l:opp}
Let $b_0>0$, $0<d_{\min}\le d_{\max}$, and let $P$ be as in \eqref{e:Tdef} on
a domain $\Omega\subset\R\times(r,\infty)$, with $d_{\min}\le d\le d_{\max}$ and
$\sigma B\ge b_0$. If $\tau_r\ge4h_*^2$, then
\[
 \Bigl\Vert T_{P,\Omega}\Bigl(a_{-\sigma}^2+|V_c(\cdot+h_*,t)-V_c|^2
   +|H_c(\cdot+h_*,t)-H_c|^2+|H_c(\cdot+h_*,t)|^4+H_c^4+a_-^4+a_+^4\Bigr)
   \Bigr\Vert_{L^\infty(\Omega)}
 \le C_M\tau_r^{-\frac12}.
\]
Here $C_M$ may also depend on $b_0,d_{\min},d_{\max}$.
\end{lemma}
\begin{proof}
By Lemma~\ref{l:pr}, $a_{-\sigma}^2\le C_M\tau^{-1}e^{-\frac{(y-s_{-\sigma}t)^2}
{4d_{-\sigma}\tau}}$, and $\sigma(B-s_{-\sigma})=\sigma B+c_{-\sigma}\ge b_0$. For
$\tau\ge4h_*^2$ and $0\le\theta\le1$, $(y+\theta h_*)^2\ge\frac{y^2}2-h_*^2$;
the mean value theorem, \eqref{eq:pr-gau} and $\delta|h_*|\le C_M$
(Lemma~\ref{l:cmp}) give
\[
 |V_c(y+h_*,t)-V_c(y,t)|^2\le C\delta^2h_*^2\tau^{-1}e^{-\frac{c_Gy^2}\tau}
 \le C_M\tau^{-1}e^{-\frac{c_Gy^2}\tau}.
\]
The third bound in \eqref{e:gauss}, with $(v,h)=(s_{-\sigma},0)$ and $(0,0)$,
bounds the responses to these two terms by $C_M\tau_r^{-\frac12}$. By
\eqref{eq:pr-gau} and Lemma~\ref{l:pr}, the remaining terms are bounded by
$C_M\tau^{-2}$, and Lemma~\ref{l:pos}(i) bounds their response by
$C_M\tau_r^{-1}$.
\end{proof}

\begin{lemma}\label{l:rec}
Let $r\ge S_0$, $b_0>0$, $0<d_{\min}\le d_{\max}$, and let $P$ be as in
\eqref{e:Tdef} on a domain $\Omega\subset\R\times(r,\infty)$, with
$d_{\min}\le d\le d_{\max}$, and let $C$ depend only on $b_0$, $d_{\min}$,
$d_{\max}$ and the fixed data. If $\sigma B\ge b_0$ on $\Omega$, then
\[
\begin{gathered}
 \bigl\Vert T_{P,\Omega}\bigl(q_{-\sigma}^2+e_m^2+|\Gamma|^2
   +|(V_c)_y|^2+|(\Theta_c)_y|^2+H_c^2\bigr)\bigr\Vert_{L^\infty(\Omega)}
 \le C\bigl(\omega_b(r)+\delta^2\tau_r^{-\frac12}\bigr),\\
 \bigl\Vert T_{P,\Omega}\bigl(|(V_c)_y|+|(\Theta_c)_y|+|H_c|\bigr)
   \bigr\Vert_{L^\infty(\Omega)}\le C\delta .
\end{gathered}
\]
If $|B|\le c_{\rm ac}$ on $\Omega$, then
$\Vert T_{P,\Omega}(q_+^2+q_-^2+|\Gamma|^2)\Vert_{L^\infty(\Omega)}\le C\omega_b(r)$.
\end{lemma}
\begin{proof}
By exhaustion, it suffices to work on finite domains. Put
\[
 \omega_b(r)=\sup_{t\ge r}\Vert(q_+,q_-,e_m)(t)\Vert_{L^2}^2
 +\int_r^\infty\bigl(\Vert\Gamma\Vert_{L^\infty}^2
   +\Vert(S_+,S_-,S_e)\Vert_{L^1}\bigr)\dd t,
\]
which tends to zero by Lemma~\ref{l:dep}. For $\iota=\pm1$, set
$\mathcal I_\iota(y,t)=\int_{\iota z\le\iota y}E_\iota(q_\iota,m)(z,t)\dd z$.
These primitives have the regularity in \eqref{e:Tdef} by \eqref{e:bgcl}.
Writing $E_\iota=E_\iota(q_\iota,m)$, \eqref{e:hcur} and \eqref{e:osp1} give
\[
 P\mathcal I_\iota=-\frac{q_\iota^2}2+\iota BE_\iota
   +\iota J_\iota-\iota d(E_\iota)_y
   +\int_{\iota z\le\iota y}S_\iota(z,t)\dd z .
\]
If $\iota B\le c_{\rm ac}$, then
$E_\iota\le q_\iota^2/(4c_{\rm ac})$ and
$|J_\iota|+|(E_\iota)_y|\le C_*|q_\iota||\Gamma|$ yield
\[
 P\mathcal I_\iota\le-\frac{q_\iota^2}8
   +C|\Gamma|^2+\Vert S_\iota(t)\Vert_{L^1}.
\]
Since $0\le\mathcal I_\iota\le C\omega_b(r)$, comparison with
\[
 8\Bigl[C\omega_b(r)-\mathcal I_\iota
   +\int_r^t\bigl(C\Vert\Gamma\Vert_{L^\infty}^2
     +\Vert S_\iota\Vert_{L^1}\bigr)\dd t'\Bigr]
\]
gives $\Vert T_{P,\Omega}q_\iota^2\Vert_{L^\infty}\le C\omega_b(r)$.
This applies to $\iota=-\sigma$ when $\sigma B\ge b_0$, and to both
families when $|B|\le c_{\rm ac}$. Also, Lemma~\ref{l:pos}(i) gives
$\Vert T_{P,\Omega}|\Gamma|^2\Vert_{L^\infty}\le\omega_b(r)$.

For $\sigma B\ge b_0$, use $x=\sigma y$ and
$\mathcal I_e(x,t)=\int_x^\infty\frac{e_m^2}2(\sigma z,t)\dd z$.
By \eqref{e:mbal} and \eqref{e:osp1},
\[
 P\mathcal I_e=-\frac{\sigma Be_m^2}2+de_m\partial_xe_m-\sigma e_mf_e
   +\int_x^\infty S_e(\sigma z,t)\dd z
 \le-\frac{b_0e_m^2}4+C|\Gamma|^2+\Vert S_e(t)\Vert_{L^1},
\]
where \eqref{e:hcur} was used in the last step. The same comparison,
with $8$ replaced by $4/b_0$, gives
$\Vert T_{P,\Omega}e_m^2\Vert_{L^\infty}\le C\omega_b(r)$.
Finally, \eqref{eq:pr-gau} gives
\[
 |(V_c)_y|+|(\Theta_c)_y|+|H_c|
 \le C\delta\tau^{-\frac12}e^{-\frac{c_Gy^2}\tau}.
\]
The second and third bounds in \eqref{e:gauss}, with $v=h=0$, give
the remaining reference terms.
\end{proof}

\subsection{Scalar kernels}\label{s:kern}

Let
\begin{equation}\label{e:nco}
 F_c(y,t)=\int_0^yf_c(z,t)\dd z,\qquad f_c=\sqrt{\rho_-V_c},\qquad
 b(F_c(y,t),t)=f_cB+(F_c)_t-d(f_c)_y,\qquad D_0=df_c^2,
\end{equation}
where $d$ is $d(V_c)$ or $\frac c{V_c}$, and $B$ is a scalar function or a
matrix function; in the second case $(F_c)_t-d(f_c)_y$ is multiplied by the
identity matrix. By \eqref{e:ac}, $D_0$ is the constant $d(\rho_-^{-1})$ or
$c\rho_-$, and $f_cs_\sigma(V_c)=s_\sigma(\rho_-^{-1})$. For
$u(y,t)=U(F_c(y,t),t)$ and $q(y,t)=f_c(y,t)Q(F_c(y,t),t)$,
\[
\begin{aligned}
 u_t-du_{yy}+Bu_y&=\bigl[U_t-D_0U_{\xi\xi}+bU_\xi\bigr](F_c,t),\\
 q_t-\pa_y(dq_y-Bq)&=f_c\bigl[Q_t-\pa_\xi(D_0Q_\xi-bQ)\bigr](F_c,t).
\end{aligned}
\]
The chain rule gives the first identity; the second follows from
$dq_y-Bq=[D_0Q_\xi-bQ](F_c,t)+(F_c)_tQ(F_c,t)$ and $\pa_y(F_c)_t=(f_c)_t$.
By \eqref{e:pt}, $(f_c)_t=-\frac{\sqrt{\rho_-}\,y(V_c)_y}{4\tau\sqrt{V_c}}$, so
\eqref{eq:pr-gau} gives
$|(F_c)_t|\le C\delta\int_0^{|y|}\frac z{\tau^{3/2}}e^{-\frac{c_Gz^2}\tau}\dd z
\le C\delta\tau^{-\frac12}$ and
\[
 |\pa_y^j[(F_c)_t-d(f_c)_y]|+|\pa_y^{j+1}f_c|
 \le C\delta\tau^{-\frac{j+1}2}\qquad(0\le j\le2).
\]
The operators $L_\sigma$ of \eqref{e:ops},
$P_V=\pa_t-\pa_y\bigl(\frac c{V_c}\pa_y\bigr)$ and $\mathscr P_\sigma$ of
\cref{l:pr} have
$(d,B)=(d(V_c),s_\sigma(V_c))$, $\bigl(\frac c{V_c},\frac{c(V_c)_y}{V_c^2}\bigr)$ and
$(d(V_c),s_\sigma(V_c)+2k(V_c)a_\sigma)$. Since
$\pa_y[g(F_c(y,t),t)]=f_cg_\xi(F_c(y,t),t)$ and $f_c,f_c^{-1}\le C$, their
drifts satisfy, for $0\le j\le2$, respectively,
\[
 \bigl|\pa_\xi^j\bigl(b-s_\sigma(\rho_-^{-1})\bigr)\bigr|\le C\delta\tau^{-\frac{j+1}2},\qquad
 |\pa_\xi^jb|\le C\delta\tau^{-\frac{j+1}2},\qquad
 \bigl|\pa_\xi^j\bigl(b-s_\sigma(\rho_-^{-1})\bigr)\bigr|\le C_M\tau^{-\frac{j+1}2},
\]
where the bounds of \cref{l:pr} for $\pa_y^ja_\sigma$ are used in the last one.

Let $I\subset[0,\infty)$ be an interval, $D_0>0$, and let $b$ be a real function
on $\R\times I$ such that $b$ and $b_\xi$ are bounded and uniformly H\"older
continuous on $\R\times J$ for every compact interval $J\subset I$. For $s<t$ in
$I$, let $H(t,s;\xi,\eta)$ and $K(t,s;\xi,\eta)$ be the fundamental solutions
given by classical parabolic theory~\cite[Chapter~1]{Friedman1964} for
\[
 v_t=D_0v_{\xi\xi}-\pa_\xi(bv),\qquad u_t=D_0u_{\xi\xi}-bu_\xi,
\]
so that $v(\xi,t)=\int_\R H(t,s;\xi,\eta)v(\eta,s)\dd\eta$ and
$u(\xi,t)=\int_\R K(t,s;\xi,\eta)u(\eta,s)\dd\eta$. By the maximum principle
$H,K\ge0$. Moreover $\int_\R H\dd\xi=\int_\R K\dd\eta=1$, and $K_\xi=-H_\eta$,
since $u_\xi$ solves the first equation when $u$ solves the second. For $b$ given
by \eqref{e:nco} with scalar $B$, put $\xi=F_c(y,t)$ and $\eta=F_c(z,s)$; then the
kernels on $\R_y$ of $q_t=\pa_y(dq_y-Bq)$ and of $u_t=du_{yy}-Bu_y$, and the
conversion of densities, are
\begin{equation}\label{e:xder}
\begin{gathered}
 f_c(y,t)H(t,s;\xi,\eta),\qquad f_c(z,s)K(t,s;\xi,\eta),\qquad
 \dd y=\frac{\dd\xi}{f_c},\\
 \pa_y[f_cv(F_c,t)]=(f_c)_yv(F_c,t)+f_c^2v_\xi(F_c,t),\qquad
 (f_c)_y=\frac{\sqrt{\rho_-}\,(V_c)_y}{2\sqrt{V_c}}.
\end{gathered}
\end{equation}

Fix $\bar b\in\R$ and $A_0,A_1\ge0$. Write $h=t-s$, $\tau_s=1+s$, $\tau_t=1+t$,
$X=\xi-\eta-\bar bh$ and $G(h,x)=(4\pi D_0h)^{-\frac12}e^{-\frac{x^2}{4D_0h}}$, and
assume
\begin{equation}\label{e:kH}
 |b-\bar b|\le A_0\tau^{-\frac12},\qquad |b_\xi|\le A_1\tau^{-1}
 \qquad\text{on }\R\times I.
\end{equation}

\begin{lemma}\label{l:kern}
Assume \eqref{e:kH}. There are $\varpi\in(0,1]$, depending only on $D_0$ and an
upper bound for $A_0$, and $C\ge1$, depending only on $D_0$ and upper bounds for
$A_0$ and $A_1$, such that \textup{(a)}--\textup{(c)} hold for $s<t$ in $I$. In
the bounds for $H$, $K$, $H_\eta$, $K_\xi$, $H_\xi+H_\eta$, $K_\xi+K_\eta$ and
$K_{\xi\xi}+K_{\xi\eta}$, the constant $C$ depends only on $D_0$ and an upper
bound for $A_0$.
\begin{enumerate}
\item[(a)] If $h\le\varpi\tau_s$, then
\[
\begin{gathered}
 C^{-1}G(C^{-1}h,X)\le H,K\le CG(Ch,X),\\
 \sqrt h\bigl(|H_\xi|+|H_\eta|+|K_\xi|+|K_\eta|\bigr)
 +h\bigl(|H_{\xi\eta}|+|K_{\xi\eta}|\bigr)\le CG(Ch,X),\\
 |H_\xi+H_\eta|+|K_\xi+K_\eta|+\sqrt h\,|K_{\xi\xi}+K_{\xi\eta}|
 \le C\sqrt h\,\Vert b_\xi\Vert_{L^\infty(\R\times[s,t])}G(Ch,X).
\end{gathered}
\]
\item[(b)] For all $s<t$ in $I$,
\[
\begin{gathered}
 H+\sqrt h\,|H_\xi|\le CG(Ch,X),\qquad
 |H_\eta|\le C\bigl(h^{-\frac12}+\tau_s^{-\frac12}\bigr)G(Ch,X),\qquad
 \int_\R K\dd\xi\le\Bigl(\frac{\tau_t}{\tau_s}\Bigr)^{A_1},\\
 K+\sqrt h\,|K_\xi|\le C\Bigl(\frac{\tau_t}{\tau_s}\Bigr)^{A_1}G(Ch,X),\qquad
 |K_\eta|+\sqrt h\,|K_{\xi\eta}|
 \le C\Bigl(\frac{\tau_t}{\tau_s}\Bigr)^{A_1}\bigl(h^{-\frac12}+\tau_s^{-\frac12}\bigr)G(Ch,X).
\end{gathered}
\]
\item[(c)] Let $r<t_0$ in $I$ with $T_0=t_0-r\le\varpi\tau_r$, let $f$ be
measurable on $\R\times(r,t_0)$ with $|f|\le F_0$, and let
\[
 u(\xi,t)=\int_r^t\!\!\int_\R K(t,s;\xi,\eta)f(\eta,s)\dd\eta\dd s,\qquad
 \bar u(\xi,t)=\int_r^t\!\!\int_\R K(t,s;\xi,\eta)|f(\eta,s)|\dd\eta\dd s.
\]
Then $u=T_{P,\Omega}f$ and $\bar u=T_{P,\Omega}|f|$ in the notation of
\eqref{e:Tdef}, for $P=\pa_t-D_0\pa_\xi^2+b\pa_\xi$ and $\Omega=\R\times(r,t_0)$.
Let $\mathcal Sg=\Psi*g$, where $\Psi$ is a Schwartz function with
$\hat\Psi\in C_c^\infty(\R)$ and $\hat\Psi=1$ near $0$. With $\mathcal S$ acting in
$\xi$ and the norms taken in $\eta$,
\[
\begin{aligned}
 \int_0^{T_0}\sup_\xi\Vert\mathcal S\pa_\xi K(t_0,t_0-h;\xi,\cdot)\Vert_{L^2}\dd h
 &\le C(1+T_0)^{\frac14},\\
 \int_0^{T_0}\sup_\xi\Vert\mathcal S\pa_\xi\pa_\eta K(t_0,t_0-h;\xi,\cdot)\Vert_{L^2}\dd h
 &\le C,\\
 \int_0^{T_0}\sup_\xi\Vert\mathcal S\pa_\xi K(t_0,t_0-h;\xi,\cdot)\Vert_{L^\infty}\dd h
 &\le C\log(2+T_0),\\
 \sup_\xi\Vert K_\xi(t_0,r;\xi,\cdot)\Vert_{L^1}\le CT_0^{-\frac12},\qquad
 \Vert g-\mathcal Sg\Vert_{L^\infty}&\le C\Vert g_\xi\Vert_{L^2}\quad(g\in H^1(\R)),
\end{aligned}
\]
where $C$ may also depend on $\Psi$. If $\xi_0\in\R$, $\rho\in(0,\infty]$,
$\varepsilon_0\ge0$ and $\bar u(\cdot,t_0)\le\varepsilon_0$ on
$(\xi_0-\rho,\xi_0+\rho)$, then
\begin{equation}\label{e:kres}
 |u_\xi(\xi_0,t_0)|\le C\sqrt{F_0\bigl(\varepsilon_0+F_0T_0e^{-\frac{\rho^2}{CT_0}}\bigr)},
\end{equation}
with $e^{-\frac{\rho^2}{CT_0}}=0$ for $\rho=\infty$. If $\rho<\infty$,
$8|\bar b|T_0\le\rho$, $f(\eta,s)=0$ for $|\eta-\xi_0|>\rho$, and
$T_{P,J\times(r,t_0)}|f|\le\varepsilon_0$ with
$J=(\xi_0-\frac32\rho,\xi_0+\frac32\rho)$, then
$\bar u\le\varepsilon_0+CF_0T_0e^{-\frac{\rho^2}{CT_0}}$ on $J\times(r,t_0]$.
\end{enumerate}
If $b\equiv\bar b$, then $H=K=G(h,X)$, and \textup{(a)} and \textup{(c)} hold
without the conditions $h\le\varpi\tau_s$ and $T_0\le\varpi\tau_r$, with $C$
depending only on $D_0$ \textup{(}and $\Psi$\textup{)}.
\end{lemma}

The proof is given in Appendix~\ref{app:sc}.

\section{Potentials and asymptotics}\label{s:pot}

Throughout this section $M\ge1$, \eqref{e:dat-aux} holds, $0<\delta\le\delta_1(M)$,
$(V_b,U_b,\Theta_b)$ is the solution of Lemma~\ref{l:bg}, $m_\pm$ and $h_*$ are
given by \eqref{e:mass}, $\sigma\in\{-1,+1\}$ with the notation of
Section~\ref{s:tools}, and every restart time $r$ satisfies $r\ge S_0$, with $S_0$
from Lemma~\ref{l:dep}. In each lemma $\delta_1(M)$, which depends only on $M$ and
the fixed data, may be decreased, and $S_0$ may be increased.

\subsection{Potentials and the equation of \texorpdfstring{$D_\sigma$}{D}}\label{s:poteq}

Put
\begin{equation}\label{e:pot}
\begin{gathered}
 \mathcal F=(\gamma-1)E_m-p_*V_b,\qquad
 \chi_m=\frac{\gamma-1}\gamma\Bigl(V_b+\frac{E_m}{p_*}\Bigr),\qquad
 \mathfrak b_\sigma(n)=\frac1{2s_\sigma(n)n},\\
 \mathcal B(y,t)=\int_{-\infty}^yb_m(z,t)\dd z,\qquad
 \mathcal G(y,t)=\int_{-\infty}^y\mathcal F(z,t)\dd z,\qquad
 \mathcal N_m(y,t)=\int_{-\infty}^y(\chi_m-V_c)(z,t)\dd z,\\
 \Psi_\sigma=\frac{\mathcal B}2+\mathfrak b_\sigma(V_c)\mathcal G .
\end{gathered}
\end{equation}
Since $\chi_m-V_c=\frac{\gamma-1}\gamma(1,0,\frac1{p_*})\cdot(\mathcal U-\widehat{\mathcal U})$,
Lemma~\ref{l:cmp} gives $b_m,\mathcal F,\chi_m-V_c\in C([0,T_e];L^1)$, and
$\Vert\int_{-\infty}^yf\dd z\Vert_{L^\infty_y}\le\Vert f\Vert_{L^1}$ gives
$\mathcal B,\mathcal G,\mathcal N_m,\Psi_\sigma\in C([0,T_e];C_b(\R))$ for every
$T_e<\infty$.

For a state $U=(v,b,E)$ with $v>0$ and for $n>0$ put
\begin{equation}\label{e:sfun}
\begin{gathered}
 p(U)=\frac{\gamma-1}v\Bigl(E-\frac{b^2}2\Bigr),\qquad
 \mathsf F(U)=\bigl(-b,p(U),p(U)b\bigr)^{\!\top},\qquad
 \mathsf A(U)=D_U\mathsf F(U),\\
 \mathsf B(U)=
 \begin{pmatrix}
 0&0&0\\0&\frac\nu v&0\\0&\frac bv\bigl(\nu-\frac\kappa{c_v}\bigr)&\frac\kappa{c_vv}
 \end{pmatrix},\qquad
 U_*(n)=\Bigl(n,0,\frac{p_*n}{\gamma-1}\Bigr),\\
 \ell_\sigma(n)=\Bigl(-\frac{p_*}{2s_\sigma(n)n},\frac12,\frac{\gamma-1}{2s_\sigma(n)n}\Bigr),
 \qquad \ell_0=\frac{\gamma-1}\gamma\Bigl(1,0,\frac1{p_*}\Bigr),\\
 R_\sigma(n)=\Bigl(-\frac1{s_\sigma(n)},1,\frac{p_*}{s_\sigma(n)}\Bigr)^{\!\top},\qquad
 R_0=r_0^\top,\qquad f_0(n)=\frac c2\log\frac n{v_-},\\
 O_\sigma(n,U)=\ell_\sigma(n)\mathsf B(U)-d(n)\ell_\sigma(n).
\end{gathered}
\end{equation}
Substitution gives $\mathsf A(U_*(n))R_\rho(n)=s_\rho(n)R_\rho(n)$ for $\rho=\pm1$,
$\mathsf A(U_*(n))R_0=0$, and duality of the rows $\ell_\sigma(n),\ell_{-\sigma}(n),
\ell_0$ with the columns $R_\sigma(n),R_{-\sigma}(n),R_0$. Since $p(\mathcal U)=P_b$
and $E_m=c_v\Theta_b+\frac{b_m^2}2$, the equations of $\mathcal U$ in Section~\ref{s:mass}
read
\begin{equation}\label{e:nfflx}
 \mathcal U_t+\mathsf A(\mathcal U)\mathcal U_y=\partial_y\bigl(\mathsf B(\mathcal U)\mathcal U_y\bigr),
 \qquad
 \mathsf B(\mathcal U)\mathcal U_y=\Bigl(0,\frac{\nu(b_m)_y}{V_b},
   \frac{\kappa(\Theta_b)_y+\nu b_m(b_m)_y}{V_b}\Bigr)^{\!\top}.
\end{equation}
For a state $U$ and $j_0\in\R$ let $m(U)$, $J_a(U)$ and
$q_\pm(U;j_0)=m(U)^{-\frac14}(b-j_0\pm J_a(U))$ be given by \eqref{e:char} with
$(V_b,P_b,\upsilon_m)$ replaced by $(v,p(U),b-j_0)$; then $m(\mathcal U)=m$ and
$q_\pm(\mathcal U;H_c)=q_\pm$. Finally put
\[
 \mathfrak a_\sigma(n,U)=\ell_\sigma(n)\bigl(U-U_*(n)\bigr),\qquad
 \mathfrak j_\sigma(n,U)=\frac{(v-n)(p(U)-p_*)}{2n}+\frac{(\gamma-1)b^2}{4n}
   -\frac{(\gamma-1)(p(U)-p_*)b}{2s_\sigma(n)n},
\]
$\mathcal A_\sigma=L_\sigma[\mathfrak b_\sigma(V_c)]$ and, for a constant $n>0$,
$\Psi^n_\sigma=\frac{\mathcal B}2+\mathfrak b_\sigma(n)\mathcal G$.

\begin{lemma}\label{l:pot}
On $\R\times(0,\infty)$,
\begin{equation}\label{e:potq}
\begin{aligned}
 \mathcal B_t&=\frac{\nu(b_m)_y}{V_b}-(P_b-p_*),\\
 \mathcal G_t&=(\gamma-1)\frac{\kappa(\Theta_b)_y+\nu b_m(b_m)_y}{V_b}
   -\gamma p_*b_m-(\gamma-1)(P_b-p_*)b_m,\\
 (\Psi_\sigma)_y&=\mathfrak a_\sigma(V_c,\mathcal U)+\mathfrak b_\sigma'(V_c)(V_c)_y\mathcal G,\\
 L_\sigma\Psi_\sigma&=\mathfrak j_\sigma(V_c,\mathcal U)+O_\sigma(V_c,\mathcal U)\mathcal U_y
   +\mathcal A_\sigma\mathcal G-2d(V_c)\mathfrak b_\sigma'(V_c)(V_c)_y\mathcal F,
\end{aligned}
\end{equation}
and, for every constant $n>0$,
\[
 (\Psi^n_\sigma)_y=\mathfrak a_\sigma(n,\mathcal U),\qquad
 \bigl(\partial_t-d(n)\partial_y^2+s_\sigma(n)\partial_y\bigr)\Psi^n_\sigma
 =\mathfrak j_\sigma(n,\mathcal U)+O_\sigma(n,\mathcal U)\mathcal U_y .
\]
There are $\epsilon_*\in(0,1]$ and $C$, depending only on the fixed data, and a
smooth map $U_{\rm rec}(n,a,Z)$, defined for $n\in[\frac{v_-}2,2v_-]$, $a\in\R$ and
$Z=(q,e,j_0)\in\R^3$ with $|a|+|Z|<\epsilon_*$, such that $U_{\rm rec}(n,0,0)=U_*(n)$
and
\[
 U=U_{\rm rec}\bigl(n,\mathfrak a_\sigma(n,U),Z\bigr),\qquad
 Z=\bigl(q_{-\sigma}(U;j_0),\,m(U)-n,\,j_0\bigr),
\]
whenever $|U-U_*(n)|<\epsilon_*$ and $|\mathfrak a_\sigma(n,U)|+|Z|<\epsilon_*$. Put
$\mathcal H_\sigma(n,a,Z)=-\mathfrak j_\sigma(n,U_{\rm rec}(n,a,Z))$ and
$h_{\sigma,n}(a)=\mathcal H_\sigma(n,a,0)$. For $|a|+|Z|<\epsilon_*$,
\begin{equation}\label{e:hex}
\begin{gathered}
 h_{\sigma,n}(0)=h_{\sigma,n}'(0)=0,\qquad
 h_{\sigma,n}''(0)=\frac{\gamma+1}{2n},\qquad
 |\mathcal H_\sigma(n,a,Z)-k(n)a^2|\le C\bigl(|a|^3+|a||Z|+|Z|^2\bigr),\\
 C^{-1}a^2-C|Z|^2\le\mathcal H_\sigma(n,a,Z)\le Ca^2+C|Z|^2,
\end{gathered}
\end{equation}
and, with the linear function $\mathcal F(U)=(\gamma-1)E-p_*v$,
\begin{equation}\label{e:hsplit}
\begin{gathered}
 \mathcal H_\sigma(n,a,Z)=h_{\sigma,n}(a)+a\mathcal H^1_\sigma+\mathcal H_\sigma(n,0,Z),
 \qquad
 \mathcal F\bigl(U_{\rm rec}(n,a,Z)\bigr)=a\mathcal F^1_\sigma
   +\mathcal F\bigl(U_{\rm rec}(n,0,Z)\bigr),\\
 \mathcal H^1_\sigma=\int_0^1\bigl[\partial_a\mathcal H_\sigma(n,\theta a,Z)
   -\partial_a\mathcal H_\sigma(n,\theta a,0)\bigr]\dd\theta,\qquad
 \mathcal F^1_\sigma=\int_0^1(-p_*,0,\gamma-1)\,\partial_aU_{\rm rec}(n,\theta a,Z)\dd\theta,\\
 |\mathcal H^1_\sigma|\le C|Z|,\qquad |\mathcal H_\sigma(n,0,Z)|\le C|Z|^2,\qquad
 |\mathcal F^1_\sigma|\le C,\qquad
 \bigl|\mathcal F\bigl(U_{\rm rec}(n,0,Z)\bigr)\bigr|\le C|Z| .
\end{gathered}
\end{equation}
Each $h_{\sigma,n}$ extends to a smooth function on $\R$, denoted in the same way,
with $C^{-1}\le h_{\sigma,n}''\le C$ on $\R$.
\end{lemma}
\begin{proof}
Integrating the momentum equation in \eqref{e:nfflx}, and its energy equation
multiplied by $\gamma-1$ minus its volume equation multiplied by $p_*$, gives
the first two lines of \eqref{e:potq}. The integrals converge by
$b_m,\mathcal F\in C([0,T_e];L^1)$; the fluxes vanish at $-\infty$ uniformly
on compact time intervals by \eqref{e:bgcl} and \eqref{eq:pr-gau}.

Since $\mathcal F(U_*(n))=0$, the definitions give
$\frac b2+\mathfrak b_\sigma(n)\mathcal F(U)=\mathfrak a_\sigma(n,U)$; differentiation of
\eqref{e:pot} therefore gives the third line of \eqref{e:potq} and
\[
 (\Psi_\sigma)_{yy}=\ell_\sigma(V_c)\mathcal U_y+2\mathfrak b_\sigma'(V_c)(V_c)_y\mathcal F
   +\bigl(\mathfrak b_\sigma''(V_c)(V_c)_y^2+\mathfrak b_\sigma'(V_c)(V_c)_{yy}\bigr)\mathcal G .
\]
By the first two lines of \eqref{e:potq} and \eqref{e:nfflx},
\[
\begin{aligned}
 \frac{\mathcal B_t}2+\mathfrak b_\sigma(V_c)\mathcal G_t+s_\sigma(V_c)\mathfrak a_\sigma(V_c,\mathcal U)
 ={}&\ell_\sigma(V_c)\mathsf B(\mathcal U)\mathcal U_y-\frac{P_b-p_*}2
   -\gamma p_*\mathfrak b_\sigma(V_c)b_m\\
 &-(\gamma-1)\mathfrak b_\sigma(V_c)(P_b-p_*)b_m+\frac{s_\sigma(V_c)b_m}2+\frac{\mathcal F}{2V_c}.
\end{aligned}
\]
Here $\gamma p_*\mathfrak b_\sigma(n)=\frac{s_\sigma(n)}2$,
$s_\sigma(n)\mathfrak b_\sigma(n)=\frac1{2n}$ and
$\mathcal F=(P_b-p_*)V_b+\frac{\gamma-1}2b_m^2$, so the last five terms equal
$\mathfrak j_\sigma(V_c,\mathcal U)$:
\begin{equation}\label{e:potid}
 \frac{\mathcal B_t}2+\mathfrak b_\sigma(V_c)\mathcal G_t
   +s_\sigma(V_c)\mathfrak a_\sigma(V_c,\mathcal U)
 =\mathfrak j_\sigma(V_c,\mathcal U)+\ell_\sigma(V_c)\mathsf B(\mathcal U)\mathcal U_y .
\end{equation}
Add $\mathfrak b_\sigma'(V_c)(V_c)_t\mathcal G$, $s_\sigma(V_c)\mathfrak b_\sigma'(V_c)(V_c)_y\mathcal G$
and $-d(V_c)(\Psi_\sigma)_{yy}$ to both sides. The left side becomes
$L_\sigma\Psi_\sigma$, and $O_\sigma=\ell_\sigma\mathsf B-d\ell_\sigma$ and
$\mathcal A_\sigma=\mathfrak b_\sigma'(V_c)L_\sigma V_c-d(V_c)\mathfrak b_\sigma''(V_c)(V_c)_y^2$
turn the right side into that of the last line of \eqref{e:potq}. For constant $n$
the same computation, without the terms containing $\mathfrak b_\sigma'$ and
$\mathfrak b_\sigma''$, gives the identities for $\Psi^n_\sigma$.

At $U_*(n)$, $D_Up=\frac1n(-p_*,0,\gamma-1)$. Thus the differentials at $(U_*(n),0)$ of
$\mathfrak a_\sigma(n,\cdot)$, $q_{-\sigma}(\cdot\,;j_0)$ and $m$ with respect to $U$ are
$\ell_\sigma(n)$, $2n^{-\frac14}\ell_{-\sigma}(n)$ and $\ell_0$, and that of
$q_{-\sigma}$ with respect to $j_0$ is $-n^{-\frac14}$. By the duality, the inverse
function theorem applies to $(U,j_0)\mapsto(\mathfrak a_\sigma(n,U),q_{-\sigma}(U;j_0),
m(U)-n,j_0)$ uniformly for $n\in[\frac{v_-}2,2v_-]$; it gives $\epsilon_*$ and
$U_{\rm rec}$, with $\partial_aU_{\rm rec}(n,0,0)=R_\sigma(n)$. Along $Z=0$, therefore,
and since $D_Up(U_*(n))R_\sigma(n)=s_\sigma(n)$, $U=U_{\rm rec}(n,a,0)=(v,b,E)$ satisfies
\[
 v-n=-\frac a{s_\sigma(n)}+O(a^2),\qquad b=a+O(a^2),\qquad
 p(U)-p_*=s_\sigma(n)a+O(a^2).
\]
Substitution into the quadratic expression for $\mathfrak j_\sigma$ gives
\[
 h_{\sigma,n}(a)=\Bigl(\frac1{2n}-\frac{\gamma-1}{4n}+\frac{\gamma-1}{2n}\Bigr)a^2+O(a^3)
 =\frac{\gamma+1}{4n}a^2+O(a^3),
\]
and $\mathcal H_\sigma$ and its first derivatives in $(a,Z)$ vanish at $(a,Z)=0$.
Taylor's formula in $a$ at fixed $Z$ and in $Z$ at $a=0$ gives \eqref{e:hsplit}; for
$\mathcal F$ it uses that $\mathcal F$ is linear and $\mathcal F(U_*(n))=0$. With
$|h_{\sigma,n}(a)-k(n)a^2|\le C|a|^3$ this gives the bound for
$\mathcal H_\sigma-k(n)a^2$ in \eqref{e:hex}, and, after decreasing $\epsilon_*$,
Young's inequality gives its last line.

Decrease $\epsilon_*$ so that $\frac{\gamma+1}{4n}\le h_{\sigma,n}''\le\frac{\gamma+1}n$ on
$(-\epsilon_*,\epsilon_*)$, and let $\phi\in C^\infty(\R)$, $0\le\phi\le1$, equal $1$ on
$[-\frac{\epsilon_*}4,\frac{\epsilon_*}4]$ and $0$ outside $[-\frac{\epsilon_*}2,\frac{\epsilon_*}2]$.
The function with value and slope zero at $a=0$ and second derivative
$\phi h_{\sigma,n}''+(1-\phi)\frac{\gamma+1}{2n}$ agrees with $h_{\sigma,n}$ on
$[-\frac{\epsilon_*}4,\frac{\epsilon_*}4]$, and its second derivative lies between
$\frac{\gamma+1}{8v_-}$ and $\frac{2(\gamma+1)}{v_-}$. Replacing $\epsilon_*$ by
$\frac{\epsilon_*}4$ gives the extension.
\end{proof}

After decreasing $\delta_1(M)$, $V_c$ takes values in
$[\frac{v_-}2,2v_-]$, and for $n$ in this interval $|s_\sigma(n)|\ge\frac{c_-}2$ and
$\frac{d_-}2\le d(n)\le2d_-$. Fix $n\in[\frac{v_-}2,2v_-]$, write $s=s_\sigma(n)$ and
$U_*=U_*(n)$, and index the families by $i\in\{s,0,o\}$, in this order:
\[
\begin{gathered}
 (R_s,R_o)=\bigl(R_\sigma(n),R_{-\sigma}(n)\bigr),\qquad
 (\ell_s,\ell_o)=\bigl(\ell_\sigma(n),\ell_{-\sigma}(n)\bigr),\\
 (\lambda_s,\lambda_0,\lambda_o)=(s,0,-s),\qquad
 \mathfrak c_i=s-\lambda_i,\qquad z_i=\ell_iz,
\end{gathered}
\]
so that $z=\sum_iz_iR_i$ and $\mathsf A(U_*)R_i=\lambda_iR_i$. Put
\[
\begin{gathered}
 \ell^\perp_\sigma=\sum_{i=0,o}\frac{O_\sigma(n,U_*)R_i}{\mathfrak c_i}\ell_i,\qquad
 k_{\rm nf}=-\frac{2n}{\gamma+1}\bigl(D_UO_\sigma(n,U_*)[R_s]R_s
   +\ell^\perp_\sigma D_U\mathsf A(U_*)[R_s]R_s\bigr),\\
 \ell^0_\sigma=\ell^\perp_\sigma+k_{\rm nf}\ell_s,\qquad
 \Lambda_\sigma(n,U)=O_\sigma(n,U)-\ell^0_\sigma\bigl[s\mathbb I-\mathsf A(U)\bigr],\qquad
 \mathsf C_{ij}=D_U\Lambda_\sigma(n,U_*)[R_i]R_j,\\
 P_{2,\sigma}(n,z)=\sum_{i=0,o}\frac{\mathsf C_{ii}}{2\mathfrak c_i}z_i^2
   +\sum_{i<j}\frac{\mathsf C_{ij}-\mathsf C_{ji}}{\mathfrak c_j-\mathfrak c_i}z_iz_j,\qquad
 F_{2,\sigma}(n,z)=\sum_{i<j}\frac{\mathfrak c_j\mathsf C_{ji}-\mathfrak c_i\mathsf C_{ij}}
   {\mathfrak c_j-\mathfrak c_i}z_iz_j,\\
 \mathcal P_\sigma(n,z)=\ell^0_\sigma(n)z+P_{2,\sigma}(n,z).
\end{gathered}
\]
Since $(\mathfrak c_s,\mathfrak c_0,\mathfrak c_o)=(0,s,2s)$, every denominator has
absolute value at least $\frac{c_-}2$, and all coefficients are smooth functions of
$n$ determined by the fixed data. With $z=U-U_*(n)$ define
\begin{equation}\label{e:nf}
 Q_\sigma(n,U)=f_0(n)+\mathcal P_\sigma(n,z),\qquad
 \mathcal R_\sigma(n,U)=O_\sigma(n,U)-(Q_\sigma)_U(n,U)\bigl[s_\sigma(n)\mathbb I-\mathsf A(U)\bigr]
   -(F_{2,\sigma})_z(n,z).
\end{equation}
In the following, $n=V_c$, $U=\mathcal U$, $z=\mathcal U-U_*(V_c)$,
$\widetilde z=z-(0,H_c,\frac{H_c^2}2)$, $d=d(V_c)$, and the
derivatives in $n$ are taken at fixed $U$ and $U_y$:
\begin{equation}\label{e:nfs}
\begin{aligned}
 F_{\rm nf}&=F_{2,\sigma}(n,z)-(Q_\sigma)_U\bigl[\mathsf B(U)-d\mathbb I\bigr]U_y,\\
 G_{\rm nf}&=\mathcal R_\sigma(n,U)U_y+(Q_\sigma)_{UU}\bigl[U_y,\mathsf B(U)U_y\bigr],\\
 E_{\rm par}&=-(F_{\rm nf})_n(V_c)_y-(Q_\sigma)_nL_\sigma V_c+2d(Q_\sigma)_{Un}[U_y](V_c)_y+d(Q_\sigma)_{nn}(V_c)_y^2,\\
 F_{\rm bulk}&=F_{2,\sigma}(n,\widetilde z)-(Q_\sigma)_U\bigl[\mathsf B(U)-d\mathbb I\bigr]z_y,\\
 G_{\rm bulk}&=\mathcal R_\sigma(n,U)z_y+(Q_\sigma)_{UU}\bigl[z_y,\mathsf B(U)z_y\bigr],\\
 F_{\rm c}&=F_{2,\sigma}(n,z)-F_{2,\sigma}(n,\widetilde z)-(Q_\sigma)_U\bigl[\mathsf B(U)-d\mathbb I\bigr]r_0^\top(V_c)_y,\\
 R_{\rm c}&=G_{\rm nf}-G_{\rm bulk}+E_{\rm par}+\partial_yF_{\rm c}.
\end{aligned}
\end{equation}

\begin{lemma}\label{l:nf}
\begin{enumerate}[label=(\roman*)]
\item For $n\in[\frac{v_-}2,2v_-]$, $U_*=U_*(n)$ and $z=U-U_*$,
\[
\begin{gathered}
 \Lambda_\sigma(n,U_*)=0,\qquad \mathsf C_{ss}=0,\qquad
 D_zP_{2,\sigma}(n,z)\bigl[s\mathbb I-\mathsf A(U_*)\bigr]+D_zF_{2,\sigma}(n,z)
   =D_U\Lambda_\sigma(n,U_*)[z],\\
 \mathcal R_\sigma(n,U)=\Lambda_\sigma(n,U)-D_U\Lambda_\sigma(n,U_*)[z]
   +D_zP_{2,\sigma}(n,z)\bigl[\mathsf A(U)-\mathsf A(U_*)\bigr],\\
 \ell^0_\sigma(n)R_0=f_0'(n)=\frac c{2n},\qquad (Q_\sigma)_n(n,U_*)=0,
\end{gathered}
\]
and, for $|z|\le\epsilon_*$, $|\mathcal R_\sigma(n,U)|\le C|z|^2$ and
$|(Q_\sigma)_n(n,U)|\le C|z|$, with $C$ depending only on the fixed data.
\item On $\R\times(0,\infty)$,
\begin{equation}\label{e:nfb}
 O_\sigma(V_c,\mathcal U)\mathcal U_y
 =L_\sigma Q_\sigma(V_c,\mathcal U)+\partial_yF_{\rm bulk}+G_{\rm bulk}+R_{\rm c}.
\end{equation}
\item Let $g_r(y)=\frac{C_1\delta}{\tau_r+y^2}$, with $C_1\ge1$ depending only on the
fixed data such that, for $t\ge r$,
\[
 |(V_c)_t|+|(V_c)_{yy}|+(V_c)_y^2\le g_r,\qquad
 |g_r'|\le\tau_r^{-\frac12}g_r,\qquad g_r^2\le C_1\delta\tau_r^{-1}g_r,\qquad
 \int_\R g_r\dd y=\pi C_1\delta\tau_r^{-\frac12}.
\]
For $t\ge S_0$,
\begin{equation}\label{e:nfr}
 |R_{\rm c}|\le C_M\bigl[(|z|+|z|^2)|(V_c)_y|+|z_y||(V_c)_y|
   +|(V_c)_t|+|(V_c)_{yy}|+(V_c)_y^2\bigr],
\end{equation}
and, as $r\to\infty$,
\[
\begin{gathered}
 \Vert F_{\rm bulk}\Vert_{L^2(r,\infty;H^1)}+\Vert G_{\rm bulk}\Vert_{L^1(r,\infty;L^2)}
   +\Vert G_{\rm bulk}\Vert_{L^2(r,\infty;L^2)}\longrightarrow0,\\
 \sup_{t\ge r}\Vert\mathcal P_\sigma(V_c,z)(t)\Vert_{L^\infty}
   +\Vert\partial_y[\mathcal P_\sigma(V_c,z)]\Vert_{L^2(r,\infty;L^\infty)}
   +\sup_{t\ge r}\Vert R_{\rm c}(t)\Vert_{W^{1,\infty}}\longrightarrow0,
\end{gathered}
\]
\[
 |R_{\rm c}|\le C_M\omega_b(r)|(V_c)_y|+C_M\bigl(|(V_c)_t|+|(V_c)_{yy}|+(V_c)_y^2\bigr)
 \le C_M\bigl(\omega_b(r)|(V_c)_y|+g_r\bigr)\qquad(t\ge r).
\]
\item With $L\ge1$ depending only on the fixed data and $\tau_r\ge L^2$, put
\begin{equation}\label{e:nfwdef}
 a_r(y)=\exp\Bigl\{-\sigma L\int_{-\infty}^yg_r(x)\dd x\Bigr\},\qquad
 w_0=\frac{a_r}{|s_\sigma(V_c)|},\qquad w_1=a_r|s_\sigma(V_c)|.
\end{equation}
Then, with $s=s_\sigma(V_c)$, $d=d(V_c)$, and $C,c_1>0$ depending only on the fixed
data,
\begin{equation}\label{e:nfwts}
\begin{gathered}
 C^{-1}\le w_i\le C,\qquad |(w_i)_y|^2\le Cg_r\qquad(i=0,1),\\
 (w_0)_t+(w_0s)_y+(w_0d)_{yy}\le-c_1g_r,\qquad
 (w_1)_t+(w_1s)_y-2w_1s_y+\bigl((w_1)_yd\bigr)_y\le-c_1g_r .
\end{gathered}
\end{equation}
\item For $F\in L^2(r,\infty;H^1)$ and $G\in L^1(r,\infty;L^2)\cap L^2(r,\infty;L^2)$ the
problem $L_\sigma u=\partial_yF+G$, $u(r)=0$, has a unique solution
$u\in C([r,\infty);H^1)$ with $u_y,u_{yy}\in L^2((r,\infty)\times\R)$, and, with $C$
independent of $F,G,r$,
\begin{equation}\label{e:nfe}
\begin{aligned}
 \sup_{t\ge r}\Vert u(t)\Vert_{L^2}^2+\int_r^\infty\Vert u_y\Vert_{L^2}^2\dd t
 &\le C\bigl[\Vert F\Vert_{L^2(r,\infty;L^2)}^2+\Vert G\Vert_{L^1(r,\infty;L^2)}^2\bigr],\\
 \sup_{t\ge r}\Vert u_y(t)\Vert_{L^2}^2+\int_r^\infty\Vert u_{yy}\Vert_{L^2}^2\dd t
 &\le C\bigl[\Vert F_y\Vert_{L^2(r,\infty;L^2)}^2+\Vert G\Vert_{L^2(r,\infty;L^2)}^2\bigr].
\end{aligned}
\end{equation}
\item Let $u_{\sigma,r}$ be the solution in (v) with $(F,G)=(F_{\rm bulk},G_{\rm bulk})$:
\begin{equation}\label{e:ubC}
 L_\sigma u_{\sigma,r}=\partial_yF_{\rm bulk}+G_{\rm bulk},\qquad u_{\sigma,r}(r)=0,
 \qquad W_{\sigma,r}=\mathcal P_\sigma(V_c,z)+u_{\sigma,r}.
\end{equation}
Then $L_\sigma(f_0(V_c)+W_{\sigma,r})=O_\sigma(V_c,\mathcal U)\mathcal U_y-R_{\rm c}$ and
$\sup_{t\ge r}\Vert W_{\sigma,r}(t)\Vert_{L^\infty}
+\Vert\partial_yW_{\sigma,r}\Vert_{L^2(r,\infty;L^\infty)}\le\omega_b(r)$.
\item There is $B_{\sigma,r}\ge0$ such that, with $C$ depending only on the fixed data,
\[
\begin{gathered}
 L_\sigma B_{\sigma,r}\ge C^{-1}(|(V_c)_y|+g_r),\qquad B_{\sigma,r}\le C\delta,\\
 |(B_{\sigma,r})_y|\le C(|(V_c)_y|+g_r),\qquad
 |(B_{\sigma,r})_y|^2\le C\delta\tau_r^{-\frac12}(|(V_c)_y|+g_r).
\end{gathered}
\]
For every operator $P$ as in \eqref{e:Tdef} on $\R\times(r,\infty)$ with
$\sigma B\ge b_0>0$ and $d_{\min}\le d\le d_{\max}$,
$\Vert T_{P,r}|R_{\rm c}|\Vert_{L^\infty}\le C_M(\delta\omega_b(r)+\delta\tau_r^{-\frac12})$,
where $C_M$ depends also on $b_0,d_{\min},d_{\max}$. The function
$u^{\rm c}_{\sigma,r}=T_{L_\sigma,r}R_{\rm c}$ satisfies
\[
 \sup_{t\ge r}\bigl(\Vert u^{\rm c}_{\sigma,r}(t)\Vert_{L^\infty}
   +\Vert\partial_yu^{\rm c}_{\sigma,r}(t)\Vert_{L^\infty}\bigr)
 +\bigl\Vert T_{L_\sigma,r}|\partial_yu^{\rm c}_{\sigma,r}|^2\bigr\Vert_{L^\infty}^{\frac12}
 \le\omega_b(r).
\]
\item Let $n_0\in[\frac{v_-}2,2v_-]$, and let $F_{\rm nf},G_{\rm nf}$ be given by
\eqref{e:nfs} with $n=n_0$ and $z=\mathcal U-U_*(n_0)$. Then
\[
 O_\sigma(n_0,\mathcal U)\mathcal U_y
 =\bigl(\partial_t-d(n_0)\partial_y^2+s_\sigma(n_0)\partial_y\bigr)\mathcal P_\sigma(n_0,z)
   +\partial_yF_{\rm nf}+G_{\rm nf},
\]
and, where $|z|\le\epsilon_*$, $|F_{\rm nf}|\le C(|z|^2+|\mathcal U_y|)$ and
$|G_{\rm nf}|\le C(|z|^2|\mathcal U_y|+|\mathcal U_y|^2)$.
\end{enumerate}
\end{lemma}
\begin{proof}
(i) Since $\ell_\sigma(n)\mathsf B(U_*)=(0,\frac\nu{2n},\frac{(\gamma-1)\kappa}{2sc_vn^2})$,
$s^2n=\gamma p_*$ and $\frac\kappa{\gamma c_v}=c$,
\[
 \ell_\sigma(n)\mathsf B(U_*)R_s=\frac\nu{2n}+\frac{(\gamma-1)\kappa p_*}{2s^2c_vn^2}=d(n),
 \qquad
 \ell_\sigma(n)\mathsf B(U_*)R_0=\frac{\kappa p_*}{2sc_vn^2}=\frac{cs}{2n};
\]
thus $O_\sigma(n,U_*)R_s=0$ and $\ell^0_\sigma R_0=\frac{O_\sigma(n,U_*)R_0}{\mathfrak c_0}
=\frac c{2n}=f_0'(n)$. With $D_Up(U_*)R_s=s$, $D_U^2p(U_*)[R_s,R_s]=\frac{3-\gamma}n$
and $b=0$ at $U_*$,
\[
\begin{gathered}
 D_U\mathsf A(U_*)[R_s]R_s=D_U^2\mathsf F(U_*)[R_s,R_s]=\Bigl(0,\frac{3-\gamma}n,2s\Bigr)^{\!\top},\\
 \ell_sD_U\mathsf A(U_*)[R_s]R_s=\frac{3-\gamma}{2n}+\frac{\gamma-1}n=\frac{\gamma+1}{2n}.
\end{gathered}
\]
Since $\ell_iR_s=0$ for $i\ne s$ and $\ell^0_\sigma R_i=\frac{O_\sigma(n,U_*)R_i}{\mathfrak c_i}$
for $i=0,o$, $\Lambda_\sigma(n,U_*)R_i=O_\sigma(n,U_*)R_i-\mathfrak c_i\ell^0_\sigma R_i=0$
for every $i$. Since $D_U\Lambda_\sigma[X]=D_UO_\sigma[X]+\ell^0_\sigma D_U\mathsf A[X]$,
\[
 \mathsf C_{ss}=D_UO_\sigma(n,U_*)[R_s]R_s+\ell^\perp_\sigma D_U\mathsf A(U_*)[R_s]R_s
   +k_{\rm nf}\frac{\gamma+1}{2n}=0 .
\]
In the coordinates $z_i$,
$D_zP_{2,\sigma}[s\mathbb I-\mathsf A(U_*)]R_j=\mathfrak c_j\partial_{z_j}P_{2,\sigma}$,
$D_zF_{2,\sigma}R_j=\partial_{z_j}F_{2,\sigma}$ and
$D_U\Lambda_\sigma(n,U_*)[z]R_j=\sum_iz_i\mathsf C_{ij}$. The coefficient of $z_j$ in
$\mathfrak c_j\partial_{z_j}P_{2,\sigma}+\partial_{z_j}F_{2,\sigma}$ is $\mathsf C_{jj}$
(for $j=s$ by $\mathfrak c_s=0$ and $\mathsf C_{ss}=0$), and for $i<j$ the coefficients of
$z_i$ in the $j$-th and of $z_j$ in the $i$-th of these expressions are
\[
 \mathfrak c_j\frac{\mathsf C_{ij}-\mathsf C_{ji}}{\mathfrak c_j-\mathfrak c_i}
   +\frac{\mathfrak c_j\mathsf C_{ji}-\mathfrak c_i\mathsf C_{ij}}{\mathfrak c_j-\mathfrak c_i}
 =\mathsf C_{ij},\qquad
 \mathfrak c_i\frac{\mathsf C_{ij}-\mathsf C_{ji}}{\mathfrak c_j-\mathfrak c_i}
   +\frac{\mathfrak c_j\mathsf C_{ji}-\mathfrak c_i\mathsf C_{ij}}{\mathfrak c_j-\mathfrak c_i}
 =\mathsf C_{ji}.
\]
This is the identity for $D_zP_{2,\sigma}$ and $D_zF_{2,\sigma}$. Inserting
$(Q_\sigma)_U=\ell^0_\sigma+D_zP_{2,\sigma}$ into \eqref{e:nf} and using it gives the
expression of $\mathcal R_\sigma$, and $\Lambda_\sigma(n,U_*)=0$ gives
$|\mathcal R_\sigma|\le C|z|^2$. Since $\frac{\dd}{\dd n}U_*(n)=R_0$ and
$\mathcal P_\sigma(n,0)=0$, $(Q_\sigma)_n(n,U_*)=f_0'(n)-\ell^0_\sigma R_0=0$; Taylor's
formula gives $|(Q_\sigma)_n|\le C|z|$.

(ii) Let $\partial^U_y$ denote the derivative in $y$ at fixed $n$; all derivatives of
$Q_\sigma$ below are evaluated at $(V_c,\mathcal U)$. The chain rule gives
\[
\begin{aligned}
 L_\sigma Q_\sigma={}&(Q_\sigma)_U\bigl(\mathcal U_t-d\mathcal U_{yy}+s_\sigma(V_c)\mathcal U_y\bigr)\\
 &-d(Q_\sigma)_{UU}[\mathcal U_y,\mathcal U_y]+(Q_\sigma)_nL_\sigma V_c\\
 &-2d(Q_\sigma)_{Un}[\mathcal U_y](V_c)_y-d(Q_\sigma)_{nn}(V_c)_y^2,
\end{aligned}
\]
and, by \eqref{e:nfflx} and
$(Q_\sigma)_U\{\partial_y(\mathsf B\mathcal U_y)-d\mathcal U_{yy}\}
=\partial^U_y\{(Q_\sigma)_U[\mathsf B-d\mathbb I]\mathcal U_y\}-(Q_\sigma)_{UU}[\mathcal U_y,(\mathsf B-d\mathbb I)\mathcal U_y]$,
\[
\begin{aligned}
 L_\sigma Q_\sigma={}&(Q_\sigma)_U\bigl[s_\sigma(V_c)\mathbb I-\mathsf A\bigr]\mathcal U_y
   +\partial^U_y\bigl\{(Q_\sigma)_U[\mathsf B-d\mathbb I]\mathcal U_y\bigr\}
   -(Q_\sigma)_{UU}[\mathcal U_y,\mathsf B\mathcal U_y]\\
 &+(Q_\sigma)_nL_\sigma V_c-2d(Q_\sigma)_{Un}[\mathcal U_y](V_c)_y-d(Q_\sigma)_{nn}(V_c)_y^2 .
\end{aligned}
\]
By \eqref{e:nf}, $(Q_\sigma)_U[s_\sigma(V_c)\mathbb I-\mathsf A]\mathcal U_y
=O_\sigma\mathcal U_y-\mathcal R_\sigma\mathcal U_y-\partial^U_yF_{2,\sigma}(V_c,z)$. With
$\partial^U_yF_{\rm nf}=\partial_yF_{\rm nf}-(F_{\rm nf})_n(V_c)_y$ this gives
$O_\sigma(V_c,\mathcal U)\mathcal U_y=L_\sigma Q_\sigma+\partial_yF_{\rm nf}+G_{\rm nf}+E_{\rm par}$.
Since $\mathcal U_y=z_y+r_0^\top(V_c)_y$, $F_{\rm nf}=F_{\rm bulk}+F_{\rm c}$, and
$G_{\rm nf}+E_{\rm par}=G_{\rm bulk}+R_{\rm c}-\partial_yF_{\rm c}$ by definition; this is
\eqref{e:nfb}. For constant $n=n_0$, this also proves (viii), whose bounds
follow from (i) and the quadratic form of $F_{2,\sigma}$.

(iii) For $t\ge r$, \eqref{eq:pr-gau} and \eqref{e:pt} give
$|(V_c)_t|+|(V_c)_{yy}|+(V_c)_y^2\le C\delta\tau^{-1}e^{-\frac{c_Gy^2}{2\tau}}
\le\frac{C\delta}{\tau_r+y^2}$, which fixes $C_1$; moreover
$\frac{|g_r'|}{g_r}=\frac{2|y|}{\tau_r+y^2}\le\tau_r^{-\frac12}$ and
$g_r\le C_1\delta\tau_r^{-1}$. For $t\ge S_0$, by \eqref{e:char},
$z=(\varphi_m,b_m,c_v\vartheta_m+\frac{b_m^2}2)$ and
$\widetilde z=(\varphi_m,\upsilon_m,c_v\vartheta_m+\frac{\upsilon_m(b_m+H_c)}2)$; the
quantities $\ell_\pm(V_c)\widetilde z-\frac{V_c^{1/4}}2q_\pm$ and $\ell_0\widetilde z-e_m$
vanish where $Z_m=0$ and are bounded by $C|Z_m|(|Z_m|+|H_c|)$. Since $F_{2,\sigma}$ contains
only mixed products,
\[
\begin{aligned}
 |F_{2,\sigma}(V_c,\widetilde z)|^2
 &\le C\bigl\{q_+^2q_-^2+(q_+^2+q_-^2)e_m^2+|Z_m|^6+H_c^2|Z_m|^4\bigr\},\\
 |\partial_yF_{2,\sigma}(V_c,\widetilde z)|
 &\le C\bigl\{|Z_m||\partial_yZ_m|+(|(V_c)_y|+|(H_c)_y|)|Z_m|^2\bigr\},\\
 |G_{\rm bulk}|&\le C(|z|^2|z_y|+|z_y|^2),\qquad
 |\partial_y[\mathcal P_\sigma(V_c,z)]|\le C(|z_y|+|z||(V_c)_y|),
\end{aligned}
\]
and the other part of $F_{\rm bulk}$ and its first derivative are bounded by
$C(|z_y|+|z_{yy}|+|z_y|^2+|(V_c)_y||z_y|)$. Lemma~\ref{l:dep}, \eqref{e:agm},
\eqref{eq:pr-gau} and, for example,
\[
 \Vert G_{\rm bulk}\Vert_{L^1(r,\infty;L^2)}
 \le C\bigl(\Vert z\Vert_{L^4(r,\infty;L^\infty)}^2+\Vert z_y\Vert_{L^2(r,\infty;L^\infty)}\bigr)
   \Vert z_y\Vert_{L^2(r,\infty;L^2)}
\]
give the three source tails and the two polynomial tails. Expanding \eqref{e:nfs}
with $H_c=\frac{c(V_c)_y}{V_c}$ gives
\[
\begin{aligned}
 |G_{\rm nf}-G_{\rm bulk}|&\le C(|z|^2+|z_y|)|(V_c)_y|+C(V_c)_y^2,\\
 |E_{\rm par}|&\le C(|z|+|z_y|)|(V_c)_y|+C(V_c)_y^2+C|z|\bigl(|(V_c)_t|+|(V_c)_{yy}|\bigr),\\
 |\partial_yF_{\rm c}|&\le C\bigl\{|z_y||(V_c)_y|+(V_c)_y^2+|(V_c)_{yy}|+|H_c||\widetilde z_y|
   +\bigl(|(H_c)_y|+|(V_c)_yH_c|\bigr)\bigl(|\widetilde z|+|H_c|\bigr)\bigr\},
\end{aligned}
\]
which prove \eqref{e:nfr}; with
$\sup_{t\ge r}(\Vert z\Vert_{L^\infty}+\Vert z_y\Vert_{L^\infty})\le\omega_b(r)$ they give
the last display of (iii). One more derivative involves state derivatives of order
at most two, whose suprema tend to zero by Lemma~\ref{l:dep}; this gives the
$W^{1,\infty}$ tail.

(iv) Since $\sigma s=|s|$, $w_0s=\sigma a_r$ and $w_1s=\sigma a_rs^2$; with
$(a_r)_y=-\sigma Lg_ra_r$,
\[
 (w_0s)_y=-Lg_ra_r,\qquad (w_1s)_y-2w_1s_y=-Lg_ra_rs^2 .
\]
The other terms of \eqref{e:nfwts} are $(w_0)_t$, $(w_0d)_{yy}$ with
$w_0d=\frac{a_rd(V_c)}{|s_\sigma(V_c)|}$, $(w_1)_t$, and $((w_1)_yd)_y$ with
$(w_1)_y=-\sigma Lg_ra_r|s|+a_r\partial_y|s|$. By
$(a_r)_{yy}=(L^2g_r^2-\sigma Lg_r')a_r$, each of them is $a_r$ times a bounded function
times one of $L^2g_r^2$, $L|g_r'|$, $Lg_r|(V_c)_y|$, $|(V_c)_t|$, $|(V_c)_{yy}|$,
$(V_c)_y^2$. By (iii) and $|(V_c)_y|\le C\delta\tau^{-\frac12}$ they are bounded by
$C_2a_rg_r(1+L\tau_r^{-\frac12}+L^2\delta\tau_r^{-1})$, with $C_2$ depending only on the
fixed data. Take $L=12C_2\max\{1,\frac4{c_-^2}\}$. For $\tau_r\ge L^2$,
$L\tau_r^{-\frac12}\le1$, the error terms are at most
$3C_2a_rg_r=\frac L4\min\{1,\frac{c_-^2}4\}a_rg_r$, and
$e^{-\pi C_1}\le a_r\le e^{\pi C_1}$ by $\int_\R g_r=\pi C_1\delta\tau_r^{-\frac12}$; with
$s^2\ge\frac{c_-^2}4$ this gives the two inequalities of \eqref{e:nfwts} with
$c_1=\frac34\min\{1,\frac{c_-^2}4\}Le^{-\pi C_1}$. Finally
$|(w_i)_y|^2\le C(L^2g_r^2+(V_c)_y^2)\le Cg_r$.

(v) Testing the equation with $w_0u$, and its derivative in $y$ with $w_1u_y$, gives
\[
\begin{aligned}
 &\frac12\frac{\dd}{\dd t}\int_\R w_0u^2\dd y+\int_\R w_0du_y^2\dd y
   -\frac12\int_\R\bigl[(w_0)_t+(w_0s)_y+(w_0d)_{yy}\bigr]u^2\dd y\\
 &\qquad=-\int_\R(w_0u)_yF\dd y+\int_\R w_0uG\dd y,\\
 &\frac12\frac{\dd}{\dd t}\int_\R w_1u_y^2\dd y+\int_\R w_1du_{yy}^2\dd y\\
 &\qquad-\frac12\int_\R\bigl[(w_1)_t+(w_1s)_y-2w_1s_y+((w_1)_yd)_y\bigr]u_y^2\dd y
   =-\int_\R(w_1u_y)_y(F_y+G)\dd y .
\end{aligned}
\]
By \eqref{e:nfwts} the third terms on the left are at least $\frac{c_1}2\int_\R g_ru^2$
and $\frac{c_1}2\int_\R g_ru_y^2$; with $|(w_i)_y|^2\le Cg_r$, Young's inequality absorbs
the terms containing $(w_i)_y$ on the right. Integration over $(r,T)$, $u(r)=0$ and
$\int_r^T|\int_\R w_0uG\dd y|\dd t\le C\sup_{[r,T]}\Vert u\Vert_{L^2}
\Vert G\Vert_{L^1(r,T;L^2)}$ give \eqref{e:nfe} on $(r,T)$ with $C$ independent of $T$.
Galerkin approximation and the same estimates for differences give the unique
solution.

(vi) Combine \eqref{e:nfb}, \eqref{e:nfe}, (iii), \eqref{e:agm} and
$\int_r^\infty\Vert u_y\Vert_{L^\infty}^2\dd t\le2\Vert u_y\Vert_{L^2((r,\infty)\times\R)}
\Vert u_{yy}\Vert_{L^2((r,\infty)\times\R)}$.

(vii) Let $\mathcal V_{g_r}$ be the function of \eqref{e:stcmp} for $P=L_\sigma$,
$x=\sigma y$ ($v=h=0$), $b_0=\frac{c_-}2$, $d_{\min}=\frac{d_-}2$, $d_{\max}=2d_-$ and the
even function $g_r$. Then $L_\sigma\mathcal V_{g_r}\ge g_r$ and
$0\le\mathcal V_{g_r}\le\frac{8\pi C_1\delta}{c_-}\tau_r^{-\frac12}$. Since
$g_r(z)\le e^{\tau_r^{-1/2}|z-x|}g_r(x)$, for $\tau_r^{-\frac12}\le\frac{c_-}{8d_-}$,
\[
 |(\mathcal V_{g_r})_y|=p_{g_r}(\sigma y)
 \le\frac2{d_-}\int_0^\infty e^{-(\frac{c_-}{4d_-}-\tau_r^{-1/2})\zeta}\dd\zeta\,
   g_r(y)\le\frac{16}{c_-}g_r(y).
\]
Put $B_{\sigma,r}=(2+2d_-)\mathcal V_{g_r}(\sigma y)+\sigma\operatorname{sgn}(\Delta_v)(V_c-v_{-\sigma})$.
Since $V_c$ is monotone between $v_-$ and $v_+$, the second term equals
$|V_c-v_{-\sigma}|\le|\Delta_v|$, and $\sigma\operatorname{sgn}(\Delta_v)s_\sigma(V_c)(V_c)_y
=|s_\sigma(V_c)||(V_c)_y|$ gives
\[
 L_\sigma\bigl[\sigma\operatorname{sgn}(\Delta_v)V_c\bigr]
 \ge\frac{c_-}2|(V_c)_y|-|(V_c)_t|-2d_-|(V_c)_{yy}|
 \ge\frac{c_-}2|(V_c)_y|-(1+2d_-)g_r .
\]
Hence $L_\sigma B_{\sigma,r}\ge\frac{c_-}2|(V_c)_y|+g_r$, $0\le B_{\sigma,r}\le C\delta$
and $|(B_{\sigma,r})_y|\le C(|(V_c)_y|+g_r)$; with $g_r^2\le C_1\delta\tau_r^{-1}g_r$ and
$(V_c)_y^2\le C\delta\tau_r^{-\frac12}|(V_c)_y|$ this gives the bound for
$|(B_{\sigma,r})_y|^2$. For $P$ as in (vii), \eqref{e:stcmp} with the parameters of $P$
gives $\Vert T_{P,r}g_r\Vert_{L^\infty}\le\frac{d_{\max}}{b_0d_{\min}}\pi C_1\delta\tau_r^{-\frac12}$,
Lemma~\ref{l:rec} gives $\Vert T_{P,r}|(V_c)_y|\Vert_{L^\infty}\le C\delta$, and the last
display of (iii) gives the bound for $T_{P,r}|R_{\rm c}|$.

Let $u=u^{\rm c}_{\sigma,r}$. The preceding bound with $P=L_\sigma$ gives
$\sup_{t\ge r}\Vert u(t)\Vert_{L^\infty}\le\omega_b(r)$. Since
$L_\sigma(u^2)=2uR_{\rm c}-2du_y^2$ and $u(r)=0$,
$d_-T_{L_\sigma,r}u_y^2\le2\sup_{t\ge r}\Vert u(t)\Vert_{L^\infty}\,T_{L_\sigma,r}|R_{\rm c}|$,
which bounds the last term. For the gradient, apply Lemma~\ref{l:kern}(c)
after \eqref{e:nco} on $[\max\{r,t-1\},t]$, increasing $S_0$ so that
$1\le\varpi\tau_r$. The response to $|R_{\rm c}|$ on this interval is at most
$\sup_{l\ge r}\Vert R_{\rm c}(l)\Vert_{L^\infty}$, so \eqref{e:kres} bounds the
forced gradient by a constant times this supremum. The homogeneous term is zero for
$t\le r+1$ and has age $1$ otherwise. Thus the $L^1$ bound for $K_\xi$ gives
\[
 \sup_{t\ge r}\Vert u_y(t)\Vert_{L^\infty}
 \le C\sup_{t\ge r}\bigl(\Vert u(t)\Vert_{L^\infty}
       +\Vert R_{\rm c}(t)\Vert_{L^\infty}\bigr)\le\omega_b(r)
\]
by (iii) and the preceding bound for $u$.
\end{proof}

With $m_-$ from \eqref{e:mass} put
\begin{equation}\label{e:plat}
 B_L=s_-m_-,\qquad G_c=\gamma p_*m_-,\qquad
 K_\sigma(n)=\frac{B_L}2+\mathfrak b_\sigma(n)G_c .
\end{equation}

\begin{lemma}\label{l:Deq}
Let $\zeta$ be smooth on $\R\times[r,\infty)$ with $0\le\zeta\le1$ and bounded
derivatives, let $u$ be the solution of
$L_\sigma u=\partial_y(\zeta F_{\rm bulk})+\zeta G_{\rm bulk}$, $u(r)=0$, given by
Lemma~\ref{l:nf}(v), and let $\Omega\subset\R\times(r,\infty)$ be open with $\zeta=1$ on
$\Omega$; for $\zeta\equiv1$, $u=u_{\sigma,r}$. Put
\[
 W=\mathcal P_\sigma(V_c,z)+u,\qquad
 D_\sigma=\Psi_\sigma-f_0(V_c)-K_\sigma(V_c)-W,\qquad
 \Delta\mathcal G=\mathcal G-G_c,\qquad Z=(q_{-\sigma},e_m,H_c).
\]
Then on $\R\times(r,\infty)$
\[
\begin{gathered}
 \mathfrak a_\sigma(V_c,\mathcal U)=(D_\sigma)_y+j_\sigma,\qquad
 j_\sigma=W_y+f_0'(V_c)(V_c)_y-\mathfrak b_\sigma'(V_c)(V_c)_y\Delta\mathcal G,\\
 -\mathfrak j_\sigma(V_c,\mathcal U)=\mathcal H_\sigma\bigl(V_c,\mathfrak a_\sigma(V_c,\mathcal U),Z\bigr),
\end{gathered}
\]
with $|\mathcal U-U_*(V_c)|+|\mathfrak a_\sigma(V_c,\mathcal U)|+|Z|<\epsilon_*$, and on
$\Omega$
\[
 L_\sigma D_\sigma=\mathfrak j_\sigma(V_c,\mathcal U)+\mathcal A_\sigma\Delta\mathcal G
   -2d(V_c)\mathfrak b_\sigma'(V_c)(V_c)_y\mathcal F+R_{\rm c}.
\]
The same holds with $(B_L,G_c)$ replaced by $(0,0)$, and with $W$ replaced by
$W+u^{\rm c}_{\sigma,r}$ and $R_{\rm c}$ removed.
\end{lemma}
\begin{proof}
On $\Omega$, $\partial_y(\zeta F_{\rm bulk})+\zeta G_{\rm bulk}=\partial_yF_{\rm bulk}+G_{\rm bulk}$,
so \eqref{e:nfb} gives $L_\sigma(f_0(V_c)+W)=O_\sigma(V_c,\mathcal U)\mathcal U_y-R_{\rm c}$
there. Since $L_\sigma[K_\sigma(V_c)]=G_c\mathcal A_\sigma$, subtraction from the last line
of \eqref{e:potq} gives the equation for $D_\sigma$. The third line of \eqref{e:potq}
and $K_\sigma'(n)=\mathfrak b_\sigma'(n)G_c$ give $\mathfrak a_\sigma=(D_\sigma)_y+j_\sigma$.
By Lemma~\ref{l:bg} and \eqref{eq:pr-gau},
$\sup_y(|\mathcal U-U_*(V_c)|+|H_c|)\to0$ as $t\to\infty$; for $t\ge S_0$ the
smallness condition holds, and Lemma~\ref{l:pot} with $n=V_c$, $j_0=H_c$,
$q_{-\sigma}(\mathcal U;H_c)=q_{-\sigma}$ and $m(\mathcal U)-V_c=e_m$ gives the identity for
$\mathfrak j_\sigma$. The other two cases follow in the same way, with
$L_\sigma u^{\rm c}_{\sigma,r}=R_{\rm c}$.
\end{proof}

\subsection{Bounded primitives}\label{s:bdd}

\begin{lemma}\label{l:bd}
\begin{equation}\label{e:bdpot}
 \sup_{t\ge0}\bigl(\Vert\mathcal B(t)\Vert_{L^\infty}+\Vert\mathcal G(t)\Vert_{L^\infty}\bigr)\le C_b,
\end{equation}
equivalently $\sup_{t\ge0}\Vert\Psi_\sigma(t)\Vert_{L^\infty}\le C_b$ for $\sigma=\pm1$,
and for every $t\ge0$
\begin{equation}\label{e:bdends}
\begin{gathered}
 \mathcal B(-\infty,t)=\mathcal G(-\infty,t)=0,\qquad
 \mathcal B(+\infty,t)=Q_{\rm mom}+c\log\frac{v_+}{v_-},\\
 \mathcal G(+\infty,t)=(\gamma-1)Q_E-p_*Q_v=\gamma p_*(m_-+m_+).
\end{gathered}
\end{equation}
\end{lemma}
\begin{proof}
The momentum component of $\int_\R(\mathcal U-\widehat{\mathcal U})\dd y=Q$
(Lemma~\ref{l:cmp}) and $\int_\R H_c\dd y=c\log\frac{v_+}{v_-}$ give $\mathcal B(+\infty,t)$.
By $R\Theta_c=p_*V_c$, $(-p_*,0,\gamma-1)\cdot(\mathcal U-\widehat{\mathcal U})=\mathcal F$, and
$(-p_*,0,\gamma-1)\cdot r_0=0$, $(-p_*,0,\gamma-1)\cdot r_\pm=\gamma p_*$ give
$\mathcal G(+\infty,t)$. This proves \eqref{e:bdends}. The identities
\[
 \mathcal B=\Psi_++\Psi_-,\qquad \mathcal G=\sqrt{\gamma p_*V_c}\,(\Psi_+-\Psi_-),
\]
show the asserted equivalence. Put
$K(S,T)=\max_\sigma\sup_{S\le t\le T}\Vert\Psi_\sigma(t)\Vert_{L^\infty}$ and $K_S=K(S,S)$;
by \eqref{e:pot}, $K(S,T)<\infty$ for $T<\infty$.

By \eqref{e:potid}, the third line of \eqref{e:potq},
$\mathfrak b_\sigma'(n)=-\frac{\mathfrak b_\sigma(n)}{2n}$ and
$2\mathfrak b_\sigma(V_c)\mathcal G=\Psi_\sigma-\Psi_{-\sigma}$,
\[
\begin{aligned}
 \bigl(\partial_t+s_\sigma(V_c)\partial_y\bigr)\Psi_\sigma
 ={}&\mathfrak j_\sigma(V_c,\mathcal U)+\frac{\nu(b_m)_y}{2V_b}
   +\frac{\gamma-1}{2s_\sigma(V_c)V_c}\,\frac{\kappa(\Theta_b)_y+\nu b_m(b_m)_y}{V_b}\\
 &-\frac{(\partial_t+s_\sigma(V_c)\partial_y)V_c}{4V_c}\,(\Psi_\sigma-\Psi_{-\sigma}).
\end{aligned}
\]
Let $0\le S<T$, $S\le t\le T$, and let $y(l)$, $S\le l\le t$, solve
$\dot y=s_\sigma(V_c(y,l))$. Then $\sigma y(l)$ increases with speed at least
$\frac{c_-}2$, and $|(V_c)_t|\le\frac{C\delta}{\tau+y^2}$ by \eqref{e:pt} and
\eqref{eq:pr-gau}; hence
\[
 \int_S^t|(V_c)_t(y(l),l)|\dd l\le\frac{2C\delta}{c_-}\int_\R\frac{\dd x}{1+x^2}\le C\delta .
\]
Since $\frac{\dd}{\dd l}V_c(y(l),l)=((V_c)_t+s_\sigma(V_c)(V_c)_y)(y(l),l)$ and
$\sigma\operatorname{sgn}(\Delta_v)s_\sigma(V_c)(V_c)_y=|s_\sigma(V_c)(V_c)_y|$,
\[
 \int_S^t|s_\sigma(V_c)(V_c)_y|(y(l),l)\dd l
 \le|\Delta_v|+\int_S^t|(V_c)_t(y(l),l)|\dd l\le C\delta .
\]
By Lemma~\ref{l:bg}, the first three terms on the right of the equation for
$\Psi_\sigma$ are at most $C_b|(V_c)_y|$ plus $C_b$ times
$\Vert Z_m\Vert_{L^\infty}^2+\Vert\partial_y(\upsilon_m,\vartheta_m)\Vert_{L^\infty}
+\Vert(H_c)_y\Vert_{L^\infty}+\Vert H_c\Vert_{L^\infty}^2$, which belongs to
$L^2(0,\infty)$ by Lemma~\ref{l:bg}, \eqref{e:agm} and \eqref{eq:pr-gau}.
Integration along these curves gives
\[
 K(S,T)\le K_S+C_b\delta
   +C_b\int_S^T\bigl(\Vert Z_m\Vert_{L^\infty}^2+\Vert\partial_y(\upsilon_m,\vartheta_m)\Vert_{L^\infty}
   +\Vert(H_c)_y\Vert_{L^\infty}+\Vert H_c\Vert_{L^\infty}^2\bigr)\dd t+C\delta K(S,T).
\]
Let $C\delta\le\frac12$ and $S=0$. The Cauchy--Schwarz inequality on $(t_1,T)$,
$T\to\infty$ and then $t_1\to\infty$ give
\begin{equation}\label{e:bdgrow}
 \lim_{T\to\infty}\frac{K(0,T)}{\sqrt{1+T}}=0,\qquad
 \lim_{r\to\infty}\frac{K_r}{\sqrt{\tau_r}}=0 .
\end{equation}

Fix a late $r$, let $u=u^{\rm c}_{\sigma,r}$ and $\Xi=\Psi_\sigma-f_0(V_c)-W_{\sigma,r}-u$.
Lemma~\ref{l:Deq}, with $(B_L,G_c)$ replaced by $(0,0)$ and $W_{\sigma,r}+u$ in place of
$W$, gives, with $Z=(q_{-\sigma},e_m,H_c)$,
\[
\begin{gathered}
 L_\sigma\Xi=-\mathcal H_\sigma(V_c,\Xi_y+j,Z)+\mathcal A_\sigma\mathcal G
   -2d(V_c)\mathfrak b_\sigma'(V_c)(V_c)_y\mathcal F,\\
 j=(W_{\sigma,r})_y+u_y+f_0'(V_c)(V_c)_y-\mathfrak b_\sigma'(V_c)(V_c)_y\mathcal G,
\end{gathered}
\]
and $|\mathcal F|\le C(|\Xi_y+j|+|Z|)$ by \eqref{e:hsplit}. Let $T>r$ and $K\ge1$ with
$K(r,T)\le K$. On $\R\times(r,T)$, $|\mathcal G|\le CK$,
$|\mathcal A_\sigma|\le C(|(V_c)_y|+g_r)$ and
$j^2\le C(\Vert(W_{\sigma,r})_y\Vert_{L^\infty}^2+u_y^2+(K+1)^2(V_c)_y^2)$. With
\eqref{e:hex} and $(\Xi_y+j)^2\le2\Xi_y^2+2j^2$ this gives $C_4>0$, depending only on
the fixed data, with
\begin{equation}\label{e:bdtwo}
\begin{aligned}
 L_\sigma\Xi&\le C_4\bigl[|Z|^2+\Vert(W_{\sigma,r})_y\Vert_{L^\infty}^2+u_y^2
   +(K+1)(|(V_c)_y|+g_r)+(K+1)^2(V_c)_y^2\bigr],\\
 L_\sigma\Xi&\ge-C_4\Xi_y^2-C_4\bigl[|Z|^2+\Vert(W_{\sigma,r})_y\Vert_{L^\infty}^2+u_y^2
   +(K+1)(|(V_c)_y|+g_r)+(K+1)^2(V_c)_y^2\bigr].
\end{aligned}
\end{equation}
By Lemmas~\ref{l:rec}, \ref{l:pos}(i) and \ref{l:nf}(vi), (vii), and by \eqref{e:stcmp}
and \eqref{e:gauss},
\[
\begin{gathered}
 \bigl\Vert T_{L_\sigma,r}\bigl(|Z|^2+\Vert(W_{\sigma,r})_y\Vert_{L^\infty}^2+u_y^2\bigr)
   \bigr\Vert_{L^\infty}\le\omega_b(r),\\
 \Vert T_{L_\sigma,r}(|(V_c)_y|+g_r)\Vert_{L^\infty}\le C\delta,\qquad
 \Vert T_{L_\sigma,r}(V_c)_y^2\Vert_{L^\infty}\le C\delta^2\tau_r^{-\frac12}.
\end{gathered}
\]
For late $r$, $|f_0(V_c)|+|W_{\sigma,r}|+|u|\le1$ on $\R\times(r,\infty)$, so
$|\Xi(r)|\le K_r+1$, and the first line of \eqref{e:bdtwo} and comparison give
\[
 \sup_{\R\times(r,T)}\Xi\le K_r+2+C\delta(K+1)+C\delta^2(K+1)^2\tau_r^{-\frac12}.
\]

For the other side put $\eta=\frac{2C_4}{d_-}$, so that $\eta d(V_c)\ge C_4$, and, with
$A_0\ge1$ and $C_5=\frac{4\eta C_4}{d_-}$,
\[
 \mathcal V=\exp\Bigl\{A_0(K+1)B_{\sigma,r}-C_5u^2
   +\eta C_4\int_r^t\Vert(W_{\sigma,r})_y(l)\Vert_{L^\infty}^2\dd l\Bigr\}.
\]
By $L_\sigma(u^2)=2uR_{\rm c}-2du_y^2$, with $d=d(V_c)$ and $B=B_{\sigma,r}$,
\[
\begin{aligned}
 \frac{L_\sigma\mathcal V}{\mathcal V}
 ={}&A_0(K+1)L_\sigma B-dA_0^2(K+1)^2B_y^2-2C_5uR_{\rm c}
   +2C_5d(1-2C_5u^2)u_y^2\\
 &+4dA_0(K+1)C_5B_yuu_y+\eta C_4\Vert(W_{\sigma,r})_y\Vert_{L^\infty}^2,\\
 |4dA_0(K+1)C_5B_yuu_y|
 \le{}&\frac{C_5d_-}4u_y^2+64d_-A_0^2(K+1)^2C_5u^2B_y^2 .
\end{aligned}
\]
Let $64C_5u^2\le1$ on $\R\times(r,\infty)$. By Lemma~\ref{l:nf}(vii), with its constant $C$,
\[
\begin{gathered}
 A_0(K+1)L_\sigma B\ge\frac{A_0}C(K+1)(|(V_c)_y|+g_r),\\
 (d+64d_-C_5u^2)A_0^2(K+1)^2B_y^2\le3d_-CA_0^2(K+1)^2\delta\tau_r^{-\frac12}(|(V_c)_y|+g_r),\\
 2C_5d(1-2C_5u^2)u_y^2-\frac{C_5d_-}4u_y^2\ge\frac{C_5d_-}4u_y^2=\eta C_4u_y^2,
\end{gathered}
\]
by the last display of Lemma~\ref{l:nf}(iii),
$2C_5|uR_{\rm c}|\le2C_5C_M|u|(|(V_c)_y|+g_r)$, and
$(K+1)^2(V_c)_y^2\le C\delta(K+1)^2\tau_r^{-\frac12}|(V_c)_y|$. Take $A_0=4C\eta C_4$,
then $\varepsilon_1>0$ small, both depending only on the fixed data, and let
$(K+1)\tau_r^{-\frac12}\le\varepsilon_1$ and $2C_5C_M|u|\le\eta C_4$ on
$\R\times(r,\infty)$. Then
\[
 \frac{L_\sigma\mathcal V}{\mathcal V}\ge\eta C_4\bigl[\Vert(W_{\sigma,r})_y\Vert_{L^\infty}^2+u_y^2
   +(K+1)(|(V_c)_y|+g_r)+(K+1)^2(V_c)_y^2\bigr],
\]
and, since $2d|(\log\mathcal V)_y|\le4d_-(A_0(K+1)|B_y|+2C_5|uu_y|)\le CA_0\varepsilon_1+C|uu_y|$,
for late $r$
\[
 \sigma\bigl[s_\sigma(V_c)-2d(\log\mathcal V)_y\bigr]\ge\frac{c_-}4 .
\]
By the second line of \eqref{e:bdtwo} and $\eta d\ge C_4$,
$L_\sigma e^{-\eta\Xi}=-\eta e^{-\eta\Xi}(L_\sigma\Xi+\eta d\Xi_y^2)$ is at most $\eta C_4$
times the bracket in \eqref{e:bdtwo} times $e^{-\eta\Xi}$. With the product rule,
\[
 \bigl[L_\sigma-2d(\log\mathcal V)_y\partial_y\bigr]\bigl(e^{-\eta\Xi}\mathcal V^{-1}\bigr)
 =\frac{L_\sigma e^{-\eta\Xi}}{\mathcal V}-\frac{L_\sigma\mathcal V}{\mathcal V}\,
   e^{-\eta\Xi}\mathcal V^{-1}
 \le\eta C_4|Z|^2e^{-\eta\Xi}\mathcal V^{-1}.
\]
By Lemma~\ref{l:rec}, for late $r$ the response of
$L_\sigma-2d(\log\mathcal V)_y\partial_y$ on $\R\times(r,T)$ to $|Z|^2$ is at most
$\frac1{2\eta C_4}$. Since $e^{-\eta\Xi}\mathcal V^{-1}$ is bounded on $\R\times(r,T)$ and
at most $e^{\eta(K_r+1)}$ at $t=r$, where $\mathcal V\ge1$, comparison gives
$\sup e^{-\eta\Xi}\mathcal V^{-1}\le e^{\eta(K_r+1)}+\frac12\sup e^{-\eta\Xi}\mathcal V^{-1}$,
the suprema over $\R\times(r,T)$. With
$\mathcal V\le\exp\{CA_0(K+1)\delta+\eta C_4\omega_b(r)\}$,
\[
 \sup_{\R\times(r,T)}(-\Xi)\le K_r+1+\frac{\log2}\eta+C_4\omega_b(r)+C\delta(K+1).
\]

Thus there are $C_7\ge1$ and $C_8$, depending only on the fixed data, such that, if
$K(r,T)\le K$ and $(K+1)\tau_r^{-\frac12}\le\varepsilon_1$, then
\[
 K(r,T)\le K_r+C_7+C_8\delta(K+1)+C_8\delta^2(K+1)^2\tau_r^{-\frac12}.
\]
Decrease $\delta_1(M)$ so that $C_8\delta\le\frac1{16}$. By \eqref{e:bdgrow} there is a
late $r$ with the preceding properties for which $K_*=4(K_r+C_7)$ satisfies
$(K_*+1)\tau_r^{-\frac12}\le\varepsilon_1$ and $C_8\delta^2(K_*+1)\tau_r^{-\frac12}\le\frac1{16}$.
If $K(r,T)\le K_*$, the bound with $K=K_*$ gives
$K(r,T)\le\frac{K_*}4+\frac{K_*+1}8<\frac{K_*}2$, since $K_*\ge4$. As $T\mapsto K(r,T)$
is continuous and $K(r,r)=K_r<K_*$, $K(r,T)\le K_*$ for every $T>r$. With
$K(0,r)<\infty$ this proves \eqref{e:bdpot}.
\end{proof}

\subsection{Limits of the potentials}\label{s:lim}

\begin{lemma}\label{l:ham}
Fix $c_0$ with $0<c_0<\frac14\min\{c_-,c_+\}$; $S_0$ may depend on $c_0$. For $T_e>r$
and $\ell>c_0T_e$ put
\[
\begin{gathered}
 I_\sigma(t)=\{y:\sigma y>c_0t\},\qquad
 \Omega_{\sigma,r}=\{(y,t):t>r,\ y\in I_\sigma(t)\},\\
 \Omega_\sigma(r,T_e,\ell)=\{(y,t):r<t<T_e,\ c_0t<\sigma y<\ell\}.
\end{gathered}
\]
There is a function $W^e_{\sigma,r}$ on $\Omega_{\sigma,r}$ such that
$W^e_{\sigma,r}-\mathcal P_\sigma(v_\sigma,\mathcal U-U_*(v_\sigma))$ vanishes at $t=r$ and on
$\sigma y=c_0t$, and
\[
 P_\sigma W^e_{\sigma,r}=O_\sigma(v_\sigma,\mathcal U)\mathcal U_y,\qquad
 \sup_{t\ge r}\Vert W^e_{\sigma,r}(t)\Vert_{L^\infty(I_\sigma(t))}
 +\Bigl(\int_r^\infty\Vert\partial_yW^e_{\sigma,r}(t)\Vert_{L^\infty(I_\sigma(t))}^2\dd t
   \Bigr)^{\frac12}\le\omega_b(r).
\]
Put $\Xi_{\sigma,r}=\Psi^{v_\sigma}_\sigma-W^e_{\sigma,r}$ and
$Z=(q_{-\sigma},m-v_\sigma,H_c)$. Then on $\Omega_{\sigma,r}$
$|\mathcal U-U_*(v_\sigma)|+|\mathfrak a_\sigma(v_\sigma,\mathcal U)|+|Z|<\epsilon_*$,
\[
 \mathfrak a_\sigma(v_\sigma,\mathcal U)=\partial_y\Xi_{\sigma,r}+\partial_yW^e_{\sigma,r},\qquad
 P_\sigma\Xi_{\sigma,r}=\mathfrak j_\sigma(v_\sigma,\mathcal U)
 =-\mathcal H_\sigma(v_\sigma,\partial_y\Xi_{\sigma,r}+\partial_yW^e_{\sigma,r},Z),
\]
and there are $\mu_H>0$, depending only on the fixed data, and $C_M$ such that
\[
\begin{gathered}
 R_H=C_M\bigl(q_{-\sigma}^2+e_m^2+H_c^2+|V_c-v_\sigma|^2+|\partial_yW^e_{\sigma,r}|^2\bigr),\qquad
 P_\sigma\Xi_{\sigma,r}+\mu_H(\partial_y\Xi_{\sigma,r})^2\le R_H,\\
 \Vert T_{P_\sigma,\Omega_\sigma(r,T_e,\ell)}R_H\Vert_{L^\infty}\le\omega_b(r),
\end{gathered}
\]
with $\omega_b(r)$ independent of $T_e$ and $\ell$.
\end{lemma}
\begin{proof}
On $\Omega_{\sigma,r}$, $\frac{y^2}\tau\ge\frac{c_0^2t^2}{1+t}\ge c_0^2(t-1)$; hence
\eqref{eq:ct-gau}, \eqref{eq:pr-gau} and $V_c=\frac{R\Theta_c}{p_*}$ give
\[
 |V_c-v_\sigma|+|H_c|+|(V_c)_y|+|(H_c)_y|\le C\delta e^{-c_Gc_0^2t}
\]
there, and with Lemma~\ref{l:bg},
$\sup_{t\ge r}\Vert\mathcal U-U_*(v_\sigma)\Vert_{L^\infty(I_\sigma(t))}\to0$. Let
$z=\mathcal U-U_*(v_\sigma)$, and let $F_{\rm nf},G_{\rm nf}$ be those of
Lemma~\ref{l:nf}(viii) with $n_0=v_\sigma$. Since
$z=\widetilde z+(0,H_c,\frac{H_c^2}2)+U_*(V_c)-U_*(v_\sigma)$, with $\widetilde z$ as in the
proof of Lemma~\ref{l:nf}(iii), the differences $\ell_\pm(v_\sigma)z-\frac{V_c^{1/4}}2q_\pm$
and $\ell_0z-e_m$ are at most $C(|Z_m|^2+|V_c-v_\sigma|+|H_c|)$, and, for $i\ne j$,
\[
 |z_iz_j|^2\le C\bigl\{q_+^2q_-^2+(q_+^2+q_-^2)e_m^2+|Z_m|^6+|V_c-v_\sigma|^2+H_c^2\bigr\}.
\]
With the spatial norms taken on $I_\sigma(t)$, Lemma~\ref{l:dep} and the estimates in the
proof of Lemma~\ref{l:nf}(iii) give
\[
 \Vert F_{\rm nf}\Vert_{L^2(r,\infty;H^1)}+\Vert G_{\rm nf}\Vert_{L^1(r,\infty;L^2)}
 +\Vert G_{\rm nf}\Vert_{L^2(r,\infty;L^2)}\longrightarrow0 .
\]
In $x=\sigma y-c_0t$, $P_\sigma=\partial_t+(|s_\sigma|-c_0)\partial_x-d_\sigma\partial_x^2$. Let
$u$ solve
\[
 u_t+(|s_\sigma|-c_0)u_x-d_\sigma u_{xx}=\sigma\partial_xF_{\rm nf}+G_{\rm nf}
 \quad(x>0,\ t>r),\qquad u(x,r)=u(0,t)=0 .
\]
With $\Vert\cdot\Vert=\Vert\cdot\Vert_{L^2(0,\infty)}$, testing with $u$ and with $-u_{xx}$ gives
\[
\begin{aligned}
 \frac12\frac{\dd}{\dd t}\Vert u\Vert^2+d_\sigma\Vert u_x\Vert^2
 &=-\sigma\int_0^\infty u_xF_{\rm nf}\dd x+\int_0^\infty uG_{\rm nf}\dd x,\\
 \frac12\frac{\dd}{\dd t}\Vert u_x\Vert^2+d_\sigma\Vert u_{xx}\Vert^2
   +\frac{|s_\sigma|-c_0}2u_x(0,t)^2
 &=-\int_0^\infty u_{xx}\bigl(\sigma\partial_xF_{\rm nf}+G_{\rm nf}\bigr)\dd x .
\end{aligned}
\]
Young's inequality and $u(r)=0$ give
\[
 \sup_{t\ge r}\Vert u(t)\Vert_{H^1}^2+\int_r^\infty\bigl(\Vert u_x\Vert^2+\Vert u_{xx}\Vert^2\bigr)\dd t
 \le C\bigl(\Vert F_{\rm nf}\Vert_{L^2(r,\infty;H^1)}^2+\Vert G_{\rm nf}\Vert_{L^1(r,\infty;L^2)}^2
   +\Vert G_{\rm nf}\Vert_{L^2(r,\infty;L^2)}^2\bigr),
\]
and Galerkin approximation gives the solution. Put
$W^e_{\sigma,r}(y,t)=\mathcal P_\sigma(v_\sigma,z(y,t))+u(\sigma y-c_0t,t)$. By
Lemma~\ref{l:nf}(viii), $P_\sigma W^e_{\sigma,r}=O_\sigma(v_\sigma,\mathcal U)\mathcal U_y$; the
traces are those of $u$, and $\Vert u\Vert_{L^\infty(0,\infty)}^2\le2\Vert u\Vert\Vert u_x\Vert$,
$|\mathcal P_\sigma(v_\sigma,z)|\le C|z|$ and $|\partial_y\mathcal P_\sigma(v_\sigma,z)|\le C|\mathcal U_y|$
give the bounds.

Lemma~\ref{l:pot} with $n=v_\sigma$ and $j_0=H_c$, for late $r$, gives the box condition
and the identities for $\Xi_{\sigma,r}$. By \eqref{e:hex},
$\frac12(\partial_y\Xi_{\sigma,r})^2-|\partial_yW^e_{\sigma,r}|^2\le(\partial_y\Xi_{\sigma,r}+\partial_yW^e_{\sigma,r})^2$ and
$|Z|^2\le2(q_{-\sigma}^2+2e_m^2+2|V_c-v_\sigma|^2+H_c^2)$, there is $\mu_H$ with
$P_\sigma\Xi_{\sigma,r}+\mu_H(\partial_y\Xi_{\sigma,r})^2\le R_H$. Lemma~\ref{l:rec}, with
$\sigma s_\sigma=c_\sigma$, bounds the response to $q_{-\sigma}^2+e_m^2+H_c^2$, and
positivity and Lemma~\ref{l:pos}(i) bound those to $|V_c-v_\sigma|^2\le C\delta^2e^{-2c_Gc_0^2t}$
and to $|\partial_yW^e_{\sigma,r}|^2\le\Vert \partial_yW^e_{\sigma,r}(t)\Vert_{L^\infty(I_\sigma(t))}^2$.
\end{proof}

\begin{lemma}\label{l:neu}
\begin{equation}\label{e:npplim}
\begin{gathered}
 \mathcal N_m(-\infty,t)=0,\qquad \mathcal N_m(+\infty,t)=m_0\qquad(t\ge0),\\
 \sup_{t\ge0}\Vert\mathcal N_m(t)\Vert_{L^\infty}\le C_b,\qquad
 \lim_{t\to\infty}\Vert\mathcal N_m(t)-h_*(V_c(t)-v_-)\Vert_{L^\infty}=0 .
\end{gathered}
\end{equation}
Put
\begin{equation}\label{e:npr}
\begin{aligned}
 d_0&=\frac c{\chi_m},\qquad b_0=\frac{c(V_c)_y}{V_c\chi_m},\\
 g_\chi&=\partial_y\log\frac{\chi_m}{V_c},\qquad
 P_\chi=\partial_t-d_0\partial_y^2+b_0\partial_y,\\
 F_N&=\frac{(\gamma-1)c}{\gamma p_*}\Bigl(\frac\nu{V_b}-d_0\Bigr)(b_m)_y
       +\frac{(\gamma-1)(\nu-c)H_cb_m}{\gamma p_*V_b}\\
 &\quad+\frac{H_c}{\gamma p_*}\Bigl(\frac\nu{V_b}-d_0\Bigr)b_m,\\
 R_N&=-\frac{(\gamma-1)(P_b-p_*)\upsilon_m}{\gamma p_*}
       -\frac{(P_b-p_*)H_c}{p_*}\\
 &\quad+\frac{c\mathcal F(\chi_m)_y}{\gamma p_*V_b\chi_m}
       +\frac{(\gamma-1)c(P_b-p_*)(V_b)_y}{\gamma p_*V_b}\\
 &\quad+\frac{(\gamma-1)(\nu-c)\upsilon_m(b_m)_y}{\gamma p_*V_b}
       -\frac{(\gamma-1)(\nu-c)(H_c)_yb_m}{\gamma p_*V_b}\\
 &\quad+\frac{(\gamma-1)(\nu-c)H_cb_m(V_b)_y}{\gamma p_*V_b^2}
       -\frac{(\gamma-1)cd_0g_\chi(b_m)_y}{\gamma p_*}\\
 &\quad+\frac{(P_\chi H_c)\mathcal B}{\gamma p_*}
       -\frac{(H_c)_y}{\gamma p_*}\Bigl(\frac\nu{V_b}+d_0\Bigr)b_m\\
 &\quad+\frac{H_c}{\gamma p_*}
       \Bigl(\frac{\nu(V_b)_y}{V_b^2}-d_0g_\chi\Bigr)b_m .
\end{aligned}
\end{equation}
Then $P_\chi u=u_t-\partial_y(d_0u_y)-d_0g_\chi u_y$,
\[
 P_\chi\Bigl(\mathcal N_m+\frac{(\gamma-1)c}{\gamma p_*}b_m+\frac{H_c\mathcal B}{\gamma p_*}\Bigr)
 =\partial_yF_N+R_N,\qquad
 P_\chi\bigl[h_*(V_c-v_-)\bigr]=\Bigl(1-\frac{V_c}{\chi_m}\Bigr)h_*(V_c)_t ,
\]
and, for $t\ge S_0$ and as $r\to\infty$,
\[
 |b_0|\le C\delta\tau^{-\frac12},\qquad
 |R_N|\le C_M\bigl(q_+^2+q_-^2+|\Gamma|^2\bigr)+C_b\tau^{-\frac32},\qquad
 \Vert F_N\Vert_{L^2(r,\infty;H^1)}+\int_r^\infty\Vert g_\chi(t)\Vert_{L^\infty}^2\dd t\longrightarrow0 .
\]
\end{lemma}
\begin{proof}
By \eqref{e:pot},
\[
 \chi_m=V_b+\frac{\mathcal F}{\gamma p_*},\qquad
 E_m=\frac{p_*\chi_m}{\gamma-1}+\frac{\mathcal F}\gamma,\qquad
 \mathcal F=(P_b-p_*)V_b+\frac{\gamma-1}2b_m^2 ,
\]
and $\int_\R(\chi_m-V_c)\dd y=\frac{\gamma-1}\gamma(Q_v+\frac{Q_E}{p_*})=m_0$ by Lemma~\ref{l:cmp} and
\eqref{e:mass}; this gives the endpoints in \eqref{e:npplim}. Since
$d_0(\frac{(\chi_m)_y}{\chi_m}-\frac{(V_c)_y}{V_c})=-(d_0)_y-b_0$, $P_\chi$ has the stated
divergence form.

The volume and energy equations of $\mathcal U$ and $(V_c)_t=(H_c)_y$ give
\[
 (\mathcal N_m)_t=-\frac{(\gamma-1)(P_b-p_*)b_m}{\gamma p_*}
   +\frac{\gamma-1}{\gamma p_*}\,\frac{\kappa(\Theta_b)_y+\nu b_m(b_m)_y}{V_b}-H_c ,
\]
and, by $c_v\Theta_b=E_m-\frac{b_m^2}2$, the expression of $E_m$ and $\frac\kappa{c_v}=\gamma c$,
\[
 \frac{\gamma-1}{\gamma p_*}\,\frac{\kappa(\Theta_b)_y+\nu b_m(b_m)_y}{V_b}
 =\frac{c(\chi_m)_y}{V_b}+\frac{(\gamma-1)c\mathcal F_y}{\gamma p_*V_b}
   +\frac{(\gamma-1)(\nu-\gamma c)b_m(b_m)_y}{\gamma p_*V_b}.
\]
Since $b_0\chi_m=H_c$ by \eqref{eq:prof}, $b_0V_c=d_0(V_c)_y$ and
$\frac c{V_b}-d_0=\frac{c\mathcal F}{\gamma p_*V_b\chi_m}$,
\[
\begin{aligned}
 P_\chi\mathcal N_m
 &=(\mathcal N_m)_t-d_0(\chi_m)_y+H_c\\
 &=-\frac{(\gamma-1)(P_b-p_*)b_m}{\gamma p_*}
       +\frac{c\mathcal F(\chi_m)_y}{\gamma p_*V_b\chi_m}\\
 &\quad+\frac{(\gamma-1)c\mathcal F_y}{\gamma p_*V_b}
       +\frac{(\gamma-1)(\nu-\gamma c)b_m(b_m)_y}{\gamma p_*V_b}.
\end{aligned}
\]
The momentum equation, $(P_b)_y=\frac{\mathcal F_y}{V_b}-\frac{(\gamma-1)b_m(b_m)_y}{V_b}
-\frac{(P_b-p_*)(V_b)_y}{V_b}$ and $(d_0)_y+b_0=-d_0g_\chi$ give
\[
 P_\chi b_m=-\frac{\mathcal F_y}{V_b}+\frac{(\gamma-1)b_m(b_m)_y}{V_b}
   +\frac{(P_b-p_*)(V_b)_y}{V_b}+\partial_y\Bigl[\Bigl(\frac\nu{V_b}-d_0\Bigr)(b_m)_y\Bigr]
   -d_0g_\chi(b_m)_y ,
\]
and the first line of \eqref{e:potq}, $\mathcal B_y=b_m$ and
$\partial_y(\frac\nu{V_b}-d_0)=-\frac{\nu(V_b)_y}{V_b^2}+b_0+d_0g_\chi$ give
\[
 P_\chi(H_c\mathcal B)=(P_\chi H_c)\mathcal B-(P_b-p_*)H_c
   +\partial_y\Bigl[H_c\Bigl(\frac\nu{V_b}-d_0\Bigr)b_m\Bigr]
   -(H_c)_y\Bigl(\frac\nu{V_b}+d_0\Bigr)b_m
   +H_c\Bigl(\frac{\nu(V_b)_y}{V_b^2}-d_0g_\chi\Bigr)b_m .
\]
Add the last three identities with the factors $1$, $\frac{(\gamma-1)c}{\gamma p_*}$ and
$\frac1{\gamma p_*}$. The terms with $\mathcal F_y$ cancel; the coefficient of
$b_m(b_m)_y$ becomes $\frac{(\gamma-1)(\nu-c)}{\gamma p_*V_b}$; with $b_m=\upsilon_m+H_c$ the coefficient of
$(P_b-p_*)H_c$ becomes $-\frac{\gamma-1}{\gamma p_*}-\frac1{\gamma p_*}=-\frac1{p_*}$; and
\[
 \frac{(\gamma-1)(\nu-c)H_c(b_m)_y}{\gamma p_*V_b}
 =\partial_y\frac{(\gamma-1)(\nu-c)H_cb_m}{\gamma p_*V_b}-\frac{(\gamma-1)(\nu-c)(H_c)_yb_m}{\gamma p_*V_b}
   +\frac{(\gamma-1)(\nu-c)H_cb_m(V_b)_y}{\gamma p_*V_b^2}.
\]
This is the identity for $\mathcal N_m+\frac{(\gamma-1)c}{\gamma p_*}b_m+\frac{H_c\mathcal B}{\gamma p_*}$. By
$(V_c)_t=(H_c)_y=c\partial_y\frac{(V_c)_y}{V_c}$,
$-d_0(V_c)_{yy}+b_0(V_c)_y=-\frac{V_c}{\chi_m}(V_c)_t$, which gives the equation for
$h_*(V_c-v_-)$.

For $t\ge S_0$, $\frac{v_-}4\le\chi_m\le4v_-$, and $|b_0|\le C|(V_c)_y|\le C\delta\tau^{-\frac12}$.
With $q_\pm=m^{-\frac14}(\upsilon_m\pm J_a)$ and \eqref{e:char},
$|P_b-p_*|+|\upsilon_m|\le C(|q_+|+|q_-|)$, and $V_b=m(\frac{P_b}{p_*})^{-\frac1\gamma}$,
$P_b=\frac{R\Theta_b}{V_b}$, $H_T=c\partial_y\log\Theta_b$ and \eqref{eq:prof} give
\[
 \partial_y\log\frac m{V_c}=\frac{H_T-H_c}c-\frac{(\gamma-1)(P_b)_y}{\gamma P_b},\qquad
 \partial_y\log\frac{V_b}{V_c}=\frac{H_T-H_c}c-\frac{(P_b)_y}{P_b}.
\]
Hence, with $(V_c)_yq_\pm=\frac{V_c}c(H_Tq_\pm-(H_T-H_c)q_\pm)$,
$e_m(V_c)_y=V_c(e_m)_y-mV_c\partial_y\log\frac m{V_c}$,
$(V_b-V_c)(V_c)_y=V_c(V_b-V_c)_y-V_bV_c\partial_y\log\frac{V_b}{V_c}$ and
\[
 \chi_m-m=V_b\Bigl[1+\frac{P_b-p_*}{\gamma p_*}-\Bigl(\frac{P_b}{p_*}\Bigr)^{\frac1\gamma}\Bigr]
   +\frac{(\gamma-1) b_m^2}{2\gamma p_*},
\]
the entries of $\Gamma$ give
\[
 |(V_b-V_c)_y|+|(\chi_m-V_c)_y|+|e_m(V_c)_y|+|(V_b-V_c)(V_c)_y|\le C|\Gamma|,\qquad
 |g_\chi|\le C|\Gamma|+CH_c^2|(V_c)_y| .
\]
With $\int_r^\infty\Vert\Gamma\Vert_{L^\infty}^2\dd t\to0$ (Lemma~\ref{l:dep}) this gives the
tail of $g_\chi$. In $R_N$, by $(\chi_m)_y=(V_c)_y+(\chi_m-V_c)_y$,
$(V_b)_y=(V_c)_y+(V_b-V_c)_y$, $H_c=\frac{c(V_c)_y}{V_c}$ and the expression of $\mathcal F$,
the terms containing $(P_b-p_*)(V_c)_y$ have the coefficient
\[
 \frac c{\gamma p_*\chi_m}+\frac{(\gamma-1)c}{\gamma p_*V_b}-\frac c{p_*V_c}
 =-\frac c{p_*V_c}\Bigl(\frac{\chi_m-V_c}{\gamma\chi_m}+\frac{(\gamma-1)(V_b-V_c)}{\gamma V_b}\Bigr),
\]
so that, with $\chi_m-V_c=(\chi_m-m)+e_m$, they are bounded by
$C(q_+^2+q_-^2+|\Gamma|^2)+C(H_c^4+H_c^2|(V_c)_y|)$. The other terms of $R_N$ are bounded
by $C(q_+^2+q_-^2+|\Gamma|^2)+C\tau^{-\frac32}$, using \eqref{eq:pr-gau}, except
$|(P_\chi H_c)\mathcal B|\le C_b(|(H_c)_t|+|(H_c)_{yy}|+|b_0(H_c)_y|)\le C_b\tau^{-\frac32}$, by
Lemma~\ref{l:bd}, \eqref{e:put} and \eqref{eq:pr-gau}. The term with $(b_m)_y$ in $F_N$
has a vanishing $L^2(r,\infty;H^1)$ tail by Lemma~\ref{l:dep}, and so do the others,
by $H_cb_m=H_c\upsilon_m+H_c^2$, $H_cq_\pm=H_Tq_\pm-(H_T-H_c)q_\pm$ and \eqref{eq:pr-gau}.

For every solution of $P_\chi u=\partial_yF$, the divergence form gives
\begin{equation}\label{e:pchie}
\begin{aligned}
 \frac12\frac{\dd}{\dd t}\Vert u\Vert_{L^2}^2+\int_\R d_0u_y^2\dd y
 &=\int_\R d_0g_\chi uu_y\dd y-\int_\R Fu_y\dd y,\\
 \frac12\frac{\dd}{\dd t}\Vert u_y\Vert_{L^2}^2+\int_\R d_0u_{yy}^2\dd y
 &=\int_\R b_0u_yu_{yy}\dd y-\int_\R F_yu_{yy}\dd y .
\end{aligned}
\end{equation}
Let $u_r$ solve $P_\chi u_r=\partial_yF_N$ on $\R\times(r,\infty)$, $u_r(r)=0$. Young's
inequality, Gronwall's inequality with $\int_r^\infty\Vert g_\chi\Vert_{L^\infty}^2\dd t<\infty$
in the first line of \eqref{e:pchie}, and then the second line with $|b_0|\le C\delta$, give
\[
 \sup_{t\ge r}\Vert u_r(t)\Vert_{H^1}^2+\int_r^\infty\bigl(\Vert(u_r)_y\Vert_{L^2}^2
   +\Vert(u_r)_{yy}\Vert_{L^2}^2\bigr)\dd t\le C\Vert F_N\Vert_{L^2(r,\infty;H^1)}^2 ,
\]
and Galerkin approximation gives $u_r$. Put
\[
 D=\mathcal N_m+\frac{(\gamma-1)c}{\gamma p_*}b_m+\frac{H_c\mathcal B}{\gamma p_*}-u_r-h_*(V_c-v_-).
\]
By \eqref{e:agm}, Lemmas~\ref{l:bg} and \ref{l:bd} and \eqref{eq:pr-gau},
$\sup_{t\ge r}\Vert D-\mathcal N_m+h_*(V_c-v_-)\Vert_{L^\infty}\le\omega_b(r)$; moreover
$D(\cdot,r)$ is bounded, continuous and zero at $y=\pm\infty$, since $h_*\Delta_v=m_0$, and
\[
 P_\chi D=R_N-\Bigl(1-\frac{V_c}{\chi_m}\Bigr)h_*(V_c)_t .
\]
If $w(r)\in H^1$ has compact support and $P_\chi w=0$, then \eqref{e:pchie} with $F=0$ gives
$\sup_{t\ge r}\Vert w\Vert_{L^2}<\infty$,
$\int_r^\infty(\Vert w_y\Vert_{L^2}^2+\Vert w_{yy}\Vert_{L^2}^2)\dd t<\infty$ and
$\int_r^\infty|\frac{\dd}{\dd t}\Vert w_y\Vert_{L^2}^2|\dd t<\infty$; hence
$\Vert w_y(t)\Vert_{L^2}\to0$ and $\Vert w(t)\Vert_{L^\infty}^2\le2\Vert w\Vert_{L^2}\Vert w_y\Vert_{L^2}\to0$.
Uniform approximation and the maximum principle extend $\Vert w(t)\Vert_{L^\infty}\to0$ to
every bounded continuous $w(r)$ with zero limits at $y=\pm\infty$, in particular to
$w(r)=D(\cdot,r)$. Let $0<a\le\min\{c_G,\frac{v_-}{32c}\}$ and $\xi=\frac y{\sqrt\tau}$. Then
\[
 P_\chi e^{-a\xi^2}=\tau^{-1}\bigl\{2ad_0+(a-4a^2d_0)\xi^2-2ab_0\sqrt\tau\,\xi\bigr\}e^{-a\xi^2}
 \ge\frac a4\min\Bigl\{1,\frac c{v_-}\Bigr\}\tau^{-1}(1+\xi^2)e^{-a\xi^2}
\]
for $C\delta\le\frac14\min\{1,\frac c{v_-}\}$, and $|(V_c)_t|\le C\delta\tau^{-1}|\xi|e^{-c_G\xi^2}$
by \eqref{e:pt} and \eqref{eq:pr-gau}; comparison gives
$\Vert T_{P_\chi,r}|(V_c)_t|\Vert_{L^\infty}\le C\delta$. Since $|b_0|\le c_{\rm ac}$ for late
$r$, Lemma~\ref{l:rec} and Lemma~\ref{l:pos}(i) give
$\Vert T_{P_\chi,r}|R_N|\Vert_{L^\infty}\le C_M\omega_b(r)+C_b\tau_r^{-\frac12}$. Hence
$D-w=T_{P_\chi,r}(P_\chi D)$ satisfies
\[
 \sup_{t\ge r}\Vert D-w\Vert_{L^\infty}\le C_M\omega_b(r)+C_b\tau_r^{-\frac12}
   +C\delta|h_*|\sup_{t\ge r}\Vert\chi_m-V_c\Vert_{L^\infty}\le\omega_b(r).
\]
Letting $t\to\infty$ and then $r\to\infty$ gives the limit in \eqref{e:npplim}; with
$\mathcal N_m\in C([0,T_e];C_b(\R))$ it gives the uniform bound.
\end{proof}

Let $P$ be as in \eqref{e:Tdef} on a domain $\Omega\subset\R\times(r,\infty)$, with
$\sigma B\ge b_0>0$ and $d_{\min}\le d\le d_{\max}$. Then
\begin{equation}\label{e:Aresp}
 \bigl\Vert T_{P,\Omega}|\partial_y^i\mathcal A_\sigma|\bigr\Vert_{L^\infty(\Omega)}
 \le C\delta\tau_r^{-\frac i2}\qquad(i=0,1,2),
\end{equation}
with $C$ depending only on $b_0,d_{\min},d_{\max}$ and the fixed data. Indeed,
$\mathcal A_\sigma=\mathfrak b_\sigma'(V_c)[(V_c)_t-d(V_c)(V_c)_{yy}+s_\sigma(V_c)(V_c)_y]
-d(V_c)\mathfrak b_\sigma''(V_c)(V_c)_y^2$, so \eqref{e:pt} and \eqref{eq:pr-gau} give,
with $\xi=\frac y{\sqrt\tau}$,
\[
 |\partial_y^i\mathcal A_\sigma|\le C\delta\tau^{-\frac i2}
   \bigl[\tau^{-\frac12}e^{-c_G\xi^2}+\tau^{-1}(1+|\xi|)e^{-c_G\xi^2}\bigr],
\]
and the second and third bounds of \eqref{e:gauss}, with $v=h=0$, and $\tau\ge\tau_r$
on $\Omega$ give \eqref{e:Aresp}.

\begin{lemma}\label{l:cen}
Put $a_0=\frac{c_-}{128}$; then $a_0<\frac1{32}\min\{c_-,c_+\}$. For $\sigma=\pm1$,
\[
 \lim_{t\to\infty}\sup_{\sigma y\le a_0t}
 \bigl|\Psi_\sigma(y,t)-f_0(V_c(y,t))-K_\sigma(V_c(y,t))\bigr|=0,
\]
and
\[
 \lim_{t\to\infty}\sup_{|y|\le a_0t}\Bigl(\Bigl|\mathcal B(y,t)-B_L-c\log\frac{V_c(y,t)}{v_-}\Bigr|
   +|\mathcal G(y,t)-G_c|\Bigr)=0 .
\]
\end{lemma}
\begin{proof}
Since $v_+\le2v_-$, $c_+\ge\frac{c_-}{\sqrt2}$, which gives the bound for $a_0$. The momentum
component of the identity in Lemma~\ref{l:cmp} is $Q_{\rm mom}=-c_-m_-+c_+m_+$, and
$\frac{G_c}{\sqrt{\gamma p_*v_\pm}}=c_\pm m_-$, $B_L=-c_-m_-$. With \eqref{e:bdends},
\[
\begin{gathered}
 \Psi_+(-\infty,t)-f_0(v_-)-K_+(v_-)=-\frac{B_L}2-\frac{G_c}{2\sqrt{\gamma p_*v_-}}=0,\\
 \Psi_-(+\infty,t)-f_0(v_+)-K_-(v_+)=\frac{Q_{\rm mom}+c_-m_--c_+m_+}2=0 .
\end{gathered}
\]
Fix a late $r$ and let $D_\sigma$, $j_\sigma$ and $Z$ be those of Lemma~\ref{l:Deq} with
$\zeta\equiv1$. Since $W_{\sigma,r}(t)$ vanishes at $y=\pm\infty$, $D_+(-\infty,t)=0$ and
$D_-(+\infty,t)=0$ for $t\ge r$. By Lemma~\ref{l:Deq}, \eqref{e:hsplit},
$\mathfrak a_\sigma(V_c,\mathcal U)=(D_\sigma)_y+j_\sigma$ and Taylor's formula for the extension
of $h_{\sigma,V_c}$,
\[
 \bigl[L_\sigma+v_{\sigma,r}\partial_y\bigr]D_\sigma=\mathcal A_\sigma\Delta\mathcal G+f_{\sigma,r},
\]
where, with $\mathcal H^1_\sigma$ and $\mathcal F^1_\sigma$ taken at
$(V_c,\mathfrak a_\sigma(V_c,\mathcal U),Z)$,
\[
\begin{aligned}
 v_{\sigma,r}={}&(D_\sigma)_y\int_0^1(1-\theta)h_{\sigma,V_c}''(\theta(D_\sigma)_y)\dd\theta
   +j_\sigma\int_0^1h_{\sigma,V_c}''(\theta(D_\sigma)_y)\dd\theta
   +\mathcal H^1_\sigma+2d(V_c)\mathfrak b_\sigma'(V_c)(V_c)_y\mathcal F^1_\sigma,\\
 f_{\sigma,r}={}&R_{\rm c}-\mathcal H_\sigma(V_c,0,Z)
   -2d(V_c)\mathfrak b_\sigma'(V_c)(V_c)_y\mathcal F\bigl(U_{\rm rec}(V_c,0,Z)\bigr)\\
 &-\bigl(\mathcal H^1_\sigma+2d(V_c)\mathfrak b_\sigma'(V_c)(V_c)_y\mathcal F^1_\sigma\bigr)j_\sigma
   -j_\sigma^2\int_0^1(1-\theta)h_{\sigma,V_c}''\bigl((D_\sigma)_y+\theta j_\sigma\bigr)\dd\theta .
\end{aligned}
\]
By Lemma~\ref{l:bd}, $|j_\sigma|\le|(W_{\sigma,r})_y|+C_b|(V_c)_y|$, and the bounds of
\eqref{e:hsplit} and of $h_{\sigma,V_c}''$ give
\[
 |v_{\sigma,r}|\le C\bigl(|\mathfrak a_\sigma(V_c,\mathcal U)|+|Z|+C_b|(V_c)_y|\bigr)
   +C|(W_{\sigma,r})_y|,\qquad
 |f_{\sigma,r}|\le|R_{\rm c}|+C\bigl(|Z|^2+C_b^2(V_c)_y^2+|(W_{\sigma,r})_y|^2\bigr).
\]
Put $b_{\sigma,r}=s_\sigma(V_c)+\max\{-\frac{c_-}8,\min\{v_{\sigma,r},\frac{c_-}8\}\}$,
$w_{\sigma,r}=s_\sigma(V_c)+v_{\sigma,r}-b_{\sigma,r}$ and
$P^0_{\sigma,r}=\partial_t-d(V_c)\partial_y^2+b_{\sigma,r}\partial_y$, so that
$P^0_{\sigma,r}+w_{\sigma,r}\partial_y=L_\sigma+v_{\sigma,r}\partial_y$. Then
$\sigma b_{\sigma,r}\ge\frac{3c_-}8$, and for late $r$ the first part of the bound for
$v_{\sigma,r}$ is at most $\frac{c_-}8$, so that $|w_{\sigma,r}|\le C|(W_{\sigma,r})_y|$ and,
by Lemma~\ref{l:nf}(vi), $\int_r^\infty\Vert w_{\sigma,r}\Vert_{L^\infty}^2\dd t\le\omega_b(r)$.
By Lemma~\ref{l:pos}(i), (iii), for late $r$ the responses of $L_\sigma+v_{\sigma,r}\partial_y$
on $\R\times(r,\infty)$ are at most twice those of $P^0_{\sigma,r}$. Hence
\eqref{e:Aresp}, Lemma~\ref{l:nf}(vi), (vii), Lemma~\ref{l:rec} and \eqref{e:gauss} give
\[
 \Vert T_{L_\sigma+v_{\sigma,r}\partial_y,r}|\mathcal A_\sigma|\Vert_{L^\infty}\le C_9\delta,\qquad
 \Vert T_{L_\sigma+v_{\sigma,r}\partial_y,r}|f_{\sigma,r}|\Vert_{L^\infty}\le\omega_b(r),
\]
with $C_9$ depending only on the fixed data.

Let $r_1\ge r$, and let $Y$ be bounded with $(L_\sigma+v_{\sigma,r}\partial_y)Y=0$ on
$\R\times(r_1,\infty)$ and $Y(y,r_1)\to0$ as $\sigma y\to-\infty$. For
$0<\eta\le\frac{c_-}{32d_-}$ and $L\in\R$, the function
\[
 E_L=\exp\Bigl\{\eta\Bigl[\sigma y-\frac{c_-}8(t-r_1)+\int_{r_1}^t\Vert w_{\sigma,r}(l)\Vert_{L^\infty}\dd l
   -L\Bigr]\Bigr\}
\]
satisfies
$(L_\sigma+v_{\sigma,r}\partial_y)E_L=\eta\bigl[\sigma(b_{\sigma,r}+w_{\sigma,r})-\frac{c_-}8
+\Vert w_{\sigma,r}\Vert_{L^\infty}-d(V_c)\eta\bigr]E_L\ge\frac{\eta c_-}8E_L$. Given $\varepsilon>0$,
choose $L$ with $|Y(y,r_1)|\le\varepsilon$ for $\sigma y\le L$; comparison gives
$|Y|\le\varepsilon+\sup|Y(r_1)|\,E_L$, and on $\sigma y\le a_0t$
\[
 E_L\le\exp\Bigl\{\eta\Bigl[a_0t-\frac{c_-}8(t-r_1)+\sqrt{t-r_1}\,\omega_b(r)^{\frac12}-L\Bigr]\Bigr\}
 \longrightarrow0\qquad(t\to\infty).
\]
Hence $\sup_{\sigma y\le a_0t}|Y(y,t)|\to0$. If $g$ vanishes for $t>R_1$ and
$|g|\le C'e^{\eta\sigma y}$, then $\frac{\eta c_-}8T_{L_\sigma+v_{\sigma,r}\partial_y,r_1}|g|\le
C'e^{\frac{\eta c_-}8(R_1-r_1)}E_0$ on $[r_1,R_1]$, and the response of $g$ at $t=R_1$ is bounded
and tends to zero as $\sigma y\to-\infty$; this applies to
$g=\mathcal A_\sigma\Delta\mathcal G\,\mathbf 1_{t<R_1}$, since
$|\mathcal A_\sigma\Delta\mathcal G|\le C_b\delta e^{-\frac{c_Gy^2}{1+R_1}}$ for $t\le R_1$.

Put
\[
 \varrho=\limsup_{t\to\infty}\max_{\sigma=\pm1}\sup_{|y|\le a_0t}
   \bigl|\Psi_\sigma-f_0(V_c)-K_\sigma(V_c)\bigr|,
\]
which is finite by Lemma~\ref{l:bd}. Since
\[
 \Delta\mathcal G=\sqrt{\gamma p_*V_c}\,\bigl[(\Psi_+-f_0(V_c)-K_+(V_c))-(\Psi_--f_0(V_c)-K_-(V_c))\bigr],
\]
for every $\varepsilon>0$ there is $R_1\ge r$ with $|\Delta\mathcal G|\le C(\varrho+\varepsilon)$ on
$\{|y|\le a_0t,\ t\ge R_1\}$; on $|y|>a_0t$, $|\mathcal A_\sigma\Delta\mathcal G|\le C_b\delta
e^{-\frac12c_Ga_0^2t}$ by \eqref{e:pt} and \eqref{eq:pr-gau}. On $\R\times(r,\infty)$, $D_\sigma$ is the sum of the
solution $Y$ with $Y(r)=D_\sigma(r)$ and of the responses of $L_\sigma+v_{\sigma,r}\partial_y$ to
$f_{\sigma,r}$ and to $\mathcal A_\sigma\Delta\mathcal G$ times $\mathbf 1_{t<R_1}$,
$\mathbf 1_{t\ge R_1,|y|\le a_0t}$ and $\mathbf 1_{t\ge R_1,|y|>a_0t}$. By the two preceding
paragraphs and Lemma~\ref{l:pos}(i), on $\sigma y\le a_0t$ the first and third terms tend
to zero as $t\to\infty$, and the others are at most $\omega_b(r)$,
$CC_9\delta(\varrho+\varepsilon)$ and $C_b\delta\int_{R_1}^\infty e^{-\frac12c_Ga_0^2t}\dd t$. With
$|\Psi_\sigma-f_0(V_c)-K_\sigma(V_c)|\le|D_\sigma|+|W_{\sigma,r}|$, letting $R_1\to\infty$
and $\varepsilon\to0$ gives
\[
 \varrho\le\omega_b(r)+C_{10}\delta\varrho ,
\]
with $C_{10}$ depending only on the fixed data. Decrease $\delta_1(M)$ so that
$C_{10}\delta\le\frac12$; then $r\to\infty$ gives $\varrho=0$. Repeating the representation
with $\varrho=0$ on $\sigma y\le a_0t$ gives the first limit of the lemma. By
$K_++K_-=B_L$, $2f_0(V_c)=c\log\frac{V_c}{v_-}$ and $\mathcal B=\Psi_++\Psi_-$,
$\mathcal B-B_L-c\log\frac{V_c}{v_-}=\sum_{\sigma=\pm1}(\Psi_\sigma-f_0(V_c)-K_\sigma(V_c))$; with
the expression of $\Delta\mathcal G$ this gives the second limit.
\end{proof}

\begin{lemma}\label{l:all}
Put $\bar A_-=A_--B_L$, $\bar A_+=A_+$ and
$E^{\rm base}_\sigma=\Psi_\sigma-f_0(V_c)-K_\sigma(V_c)-\bar A_\sigma$. Then
$\lim_{t\to\infty}\Vert E^{\rm base}_\sigma(t)\Vert_{L^\infty}=0$ for $\sigma=\pm1$, and
\begin{equation}\label{e:applim}
 \lim_{t\to\infty}\Bigl(\Bigl\Vert\mathcal B(t)-c\log\frac{V_c(t)}{v_-}-A_-(t)-A_+(t)\Bigr\Vert_{L^\infty}
 +\Bigl\Vert\mathcal G(t)-\sum_{\sigma=\pm1}s_\sigma v_\sigma A_\sigma(t)\Bigr\Vert_{L^\infty}\Bigr)=0 .
\end{equation}
\end{lemma}
\begin{proof}
Fix $c_0$ with $0<c_0<a_0$ in Lemma~\ref{l:ham}, and put
$\Phi_-=A_-$, $\Phi_+=f_0(v_+)+K_+(v_+)+A_+$; by Lemma~\ref{l:pr},
$P_\sigma e^{-\alpha\Phi_\sigma}=0$. Let $\Omega=\Omega_\sigma(r,T_e,\ell)$,
$\Xi=\Xi_{\sigma,r}$, $T=T_{P_\sigma,\Omega}$,
$a=\mathfrak a_\sigma(v_\sigma,\mathcal U)=\Xi_y+\partial_yW^e_{\sigma,r}$
and $Z=(q_{-\sigma},m-v_\sigma,H_c)$.
Norms in this comparison are over $\Omega$, and the constants are independent of
$T_e$ and $\ell$. Lemmas~\ref{l:bd} and~\ref{l:ham} give
$|\Xi|\le C_b$, $P_\sigma\Xi=-\mathcal H_\sigma(v_\sigma,a,Z)$ and
\[
 \bigl\Vert T(R_H+|Z|^2+|\partial_yW^e_{\sigma,r}|^2)\bigr\Vert_\infty
 \le\omega_b(r).
\]
Since $P_\sigma\Xi+\mu_H\Xi_y^2\le R_H$, Lemma~\ref{l:pos}(ii) gives
$\mu_H\Vert T\Xi_y^2\Vert_\infty\le\omega_b(r)+2C_b$;
thus $\Vert T\Xi_y^2\Vert_\infty+\Vert Ta^2\Vert_\infty\le C_b$ for late $r$.
By \eqref{e:hex},
\[
 |P_\sigma\Xi+k_\sigma\Xi_y^2|
 \le C\bigl(|a|^3+|a||Z|+|Z|^2
       +|\Xi_y\partial_yW^e_{\sigma,r}|+|\partial_yW^e_{\sigma,r}|^2\bigr).
\]
Moreover, $\sup_{\Omega_{\sigma,r}}|a|\to0$ as $r\to\infty$ by
Lemma~\ref{l:bg} and \eqref{eq:ct-gau}. Positivity and
Lemma~\ref{l:pos}(i) therefore yield
$\Vert T|P_\sigma\Xi+k_\sigma\Xi_y^2|\Vert_\infty\le\omega_b(r)$.
Since \eqref{e:ac} gives
$P_\sigma e^{-\alpha\Xi}=-\alpha e^{-\alpha\Xi}(P_\sigma\Xi+k_\sigma\Xi_y^2)$,
\[
 \bigl\Vert e^{-\alpha\Xi}
   -H_{P_\sigma,\Omega}[e^{-\alpha\Xi}|_{\partial_p\Omega}]\bigr\Vert_\infty
 \le\alpha e^{\alpha C_b}\Vert T|P_\sigma\Xi+k_\sigma\Xi_y^2|\Vert_\infty
 \le\omega_b(r).
\]
Both functions lie in $[e^{-\alpha C_b},e^{\alpha C_b}]$ by comparison.
Taking logarithms gives
\[
 \Bigl\Vert\Xi+\frac1\alpha\log H_{P_\sigma,\Omega}
   \bigl[e^{-\alpha\Xi}|_{\partial_p\Omega}\bigr]\Bigr\Vert_\infty
 \le\omega_b(r).
\]
The outer limits of $\Xi(\cdot,t)$ and $\Phi_\sigma(\cdot,t)$ agree: for $\sigma=-$ both are
zero, and for $\sigma=+$, by \eqref{e:bdends}, \eqref{e:plat} and Lemma~\ref{l:pr},
\[
 \Psi^{v_+}_+(+\infty,t)-f_0(v_+)-K_+(v_+)-s_+m_+=\frac{Q_{\rm mom}+c_-m_--c_+m_+}2=0 .
\]
On $\sigma y=c_0t\le a_0t$, Lemma~\ref{l:cen}, $|V_c-v_\sigma|\le C\delta e^{-c_Gc_0^2t}$,
$f_0(v_-)+K_-(v_-)=B_L$, \eqref{e:CH} and the trace of $W^e_{\sigma,r}$ in
Lemma~\ref{l:ham} give $(\Xi-\Phi_\sigma)(\sigma c_0t,t)\to0$ as $t\to\infty$, and this trace
does not depend on $r$. In $x=\sigma y-c_0t$, with $b=|s_\sigma|-c_0$, the solution of
$\partial_tu+b\partial_xu-d_\sigma\partial_x^2u=0$ on $x>0$ with $u(0,t)=0$ and $u(x,r)=f(x)$ is
$\int_0^\infty\mathcal K(x,\xi,t-r)f(\xi)\dd\xi$, where, for $\theta>0$ and $x,\xi>0$,
\begin{equation}\label{e:halfker}
\begin{aligned}
 \mathcal K(x,\xi,\theta)&=\frac1{\sqrt{4\pi d_\sigma\theta}}
   \exp\Bigl\{\frac{b(x-\xi)}{2d_\sigma}-\frac{b^2\theta}{4d_\sigma}\Bigr\}
   \Bigl(e^{-\frac{(x-\xi)^2}{4d_\sigma\theta}}-e^{-\frac{(x+\xi)^2}{4d_\sigma\theta}}\Bigr),\\
 0\le\mathcal K(x,\xi,\theta)&\le\frac1{\sqrt{4\pi d_\sigma\theta}}e^{-\frac{(x-\xi-b\theta)^2}{4d_\sigma\theta}} ;
\end{aligned}
\end{equation}
the first line follows from the odd reflection of the heat kernel and the factor
$\exp\{\frac{bx}{2d_\sigma}-\frac{b^2t}{4d_\sigma}\}$, which removes the drift, and
$x\xi>0$ gives the second. Let
$f_r(x)=(e^{-\alpha\Xi}-e^{-\alpha\Phi_\sigma})(\sigma(x+c_0r),r)$, which is bounded and tends
to zero as $x\to\infty$. Since
$e^{-\alpha\Phi_\sigma}=H_{P_\sigma,\Omega}[e^{-\alpha\Phi_\sigma}|_{\partial_p\Omega}]$, comparison on
$\Omega$ with $\int_0^\infty\mathcal K(x,\xi,t-r)|f_r(\xi)|\dd\xi$, with constants and with
$e^{\alpha C_b}e^{\sigma y-\ell+d_\sigma(t-r)}$, which satisfies
$P_\sigma e^{\sigma y-\ell+d_\sigma(t-r)}=|s_\sigma|e^{\sigma y-\ell+d_\sigma(t-r)}\ge0$ and
is at least $e^{\alpha C_b}$ on $\sigma y=\ell$, gives, for every $L>0$,
\[
\begin{aligned}
 \bigl|H_{P_\sigma,\Omega}\bigl[e^{-\alpha\Xi}|_{\partial_p\Omega}\bigr]-e^{-\alpha\Phi_\sigma}\bigr|
 \le{}&\frac{\Vert f_r\mathbf 1_{[0,L]}\Vert_{L^1}}{\sqrt{4\pi d_\sigma(t-r)}}+\sup_{x>L}|f_r(x)|
   +\sup_{l\ge r}\bigl|(e^{-\alpha\Xi}-e^{-\alpha\Phi_\sigma})(\sigma c_0l,l)\bigr|\\
 &+e^{\alpha C_b}e^{\sigma y-\ell+d_\sigma(t-r)} ;
\end{aligned}
\]
here both exponentials lie between $e^{-\alpha C_b}$ and $e^{\alpha C_b}$. Taking logarithms,
letting $\ell\to\infty$ at fixed $(y,t)$, then $t\to\infty$, $L\to\infty$ and $r\to\infty$, and
using the bounds for $W^e_{\sigma,r}$, one obtains
\[
 \lim_{t\to\infty}\sup_{\sigma y>c_0t}\bigl|\Psi^{v_\sigma}_\sigma(y,t)-\Phi_\sigma(y,t)\bigr|=0 .
\]
On $\sigma y>c_0t$,
$|\Psi_\sigma-\Psi^{v_\sigma}_\sigma|=|\mathfrak b_\sigma(V_c)-\mathfrak b_\sigma(v_\sigma)||\mathcal G|\le
C_b\delta e^{-c_Gc_0^2t}$, $f_0(V_c)+K_\sigma(V_c)\to f_0(v_\sigma)+K_\sigma(v_\sigma)$ uniformly,
and $\Phi_\sigma-\bar A_\sigma=f_0(v_\sigma)+K_\sigma(v_\sigma)$, by $f_0(v_-)+K_-(v_-)=B_L$; hence
$E^{\rm base}_\sigma\to0$ uniformly on $\sigma y>c_0t$. On $\sigma y\le a_0t$,
$\Psi_\sigma-f_0(V_c)-K_\sigma(V_c)\to0$ uniformly by Lemma~\ref{l:cen}, and $\bar A_\sigma\to0$
uniformly by \eqref{e:CH}, since there $\sigma(y-s_\sigma t)\le-(c_\sigma-a_0)t$ and
$A_-(+\infty,t)=s_-m_-=B_L$. As $c_0<a_0$, this proves
$\Vert E^{\rm base}_\sigma(t)\Vert_{L^\infty}\to0$.

By $\mathcal B=\Psi_++\Psi_-$, $\mathcal G=\sqrt{\gamma p_*V_c}(\Psi_+-\Psi_-)$, $K_++K_-=B_L$,
$\sqrt{\gamma p_*V_c}(K_+-K_-)=G_c=-\sqrt{\gamma p_*v_-}B_L$ and $s_\pm v_\pm=\pm\sqrt{\gamma p_*v_\pm}$,
\[
\begin{aligned}
 \mathcal B-c\log\frac{V_c}{v_-}-A_--A_+&=E^{\rm base}_++E^{\rm base}_-,\\
 \mathcal G-\sum_{\sigma=\pm1}s_\sigma v_\sigma A_\sigma
 &=\sqrt{\gamma p_*V_c}\,(E^{\rm base}_+-E^{\rm base}_-)
   +\bigl(\sqrt{\gamma p_*V_c}-\sqrt{\gamma p_*v_-}\bigr)(B_L-A_-)\\
 &\quad+\bigl(\sqrt{\gamma p_*V_c}-\sqrt{\gamma p_*v_+}\bigr)A_+ ,
\end{aligned}
\]
and the last two terms tend to zero uniformly by \eqref{eq:ct-gau} and \eqref{e:CH}. This
is \eqref{e:applim}.
\end{proof}

\subsection{Square responses}\label{s:sq}

\begin{lemma}\label{l:sq}
For $\mathscr L\in\{L_\sigma,\mathscr P_\sigma\}$,
\[
 \max_{\sigma=\pm1}\bigl\Vert T_{\mathscr L,r}|R_U|^2\bigr\Vert_{L^\infty(\R\times(r,\infty))}
 \le\omega_b(r).
\]
\end{lemma}
\begin{proof}
Fix $\sigma$. Lemma~\ref{l:pr} allows $S_0$ to be chosen so that, for $t\ge S_0$,
\[
 \tau\ge4h_*^2,\qquad \frac{v_-}2\le V_c\le2v_-,\qquad |a_\sigma|\le\frac{\epsilon_*}2,\qquad
 \sigma\bigl(s_\sigma(V_c)+2k(V_c)a_\sigma\bigr)\ge\frac{c_-}4 .
\]
Thus Lemmas~\ref{l:opp}, \ref{l:rec}, \ref{l:nf}(vii) and \eqref{e:Aresp} apply to $L_\sigma$ and
$\mathscr P_\sigma$ with $b_0=\frac{c_-}4$, $d_{\min}=d(2v_-)$ and $d_{\max}=d(\frac{v_-}2)$. Let
$\Omega\subset\R\times(r,\infty)$ be a finite domain and $T=T_{\mathscr P_\sigma,\Omega}$; all bounds
below are uniform in $\Omega$.

By Lemma~\ref{l:pr}, $L_\sigma A_\sigma+k(V_c)a_\sigma^2=R^{\rm prof}_\sigma$,
$\mathscr P_\sigma(-A_\sigma)+k(V_c)a_\sigma^2=-R^{\rm prof}_\sigma$, $|A_\sigma|\le C_M$ and
$T_{\mathscr L,\Omega}|R^{\rm prof}_\sigma|\le\int_r^\infty C_M\delta e^{-c_et}\dd t\le C_M$. Lemma~\ref{l:pos}(ii)
with $z=\pm A_\sigma$ and $\mu=k(2v_-)$ gives, for $\mathscr L\in\{L_\sigma,\mathscr P_\sigma\}$,
\begin{equation}\label{e:Ta2}
 \Vert T_{\mathscr L,\Omega}a_\sigma^2\Vert_{L^\infty(\Omega)}\le C_M .
\end{equation}
For $\rho=\pm1$, \eqref{eq:ct-gau} and Lemma~\ref{l:pr}, with $\R$ split at $y=0$ and at
$y=\frac{s_\rho t}2$, give
\begin{equation}\label{e:Vca}
 |V_c-v_\rho||a_\rho|\le C_M\delta e^{-c_et}.
\end{equation}
Let $a_0$ be as in Lemma~\ref{l:cen}, $\Delta\mathcal G=\mathcal G-G_c$, and
$\varepsilon_c(t)=\sup_{s\ge t}\sup_{|y|\le a_0s}|\Delta\mathcal G(y,s)|$; then $\varepsilon_c(t)\to0$ by
Lemma~\ref{l:cen}, and $|\Delta\mathcal G|\le C_b$ by Lemma~\ref{l:bd}. For $|y|>a_0t$ and $t\ge1$, the
pointwise bound below \eqref{e:Aresp} and \eqref{eq:pr-gau} give
$|\mathcal A_\sigma|+|\partial_y\mathcal A_\sigma|+|(V_c)_y|\le C\delta e^{-\frac{c_Ga_0^2t}4}$; hence
\begin{equation}\label{e:Acen}
 |\Delta\mathcal G|\bigl(|\mathcal A_\sigma|+|\partial_y\mathcal A_\sigma|+|(V_c)_y|\bigr)
 \le\varepsilon_c(t)\bigl(|\mathcal A_\sigma|+|\partial_y\mathcal A_\sigma|+|(V_c)_y|\bigr)
   +C_be^{-\frac{c_Ga_0^2t}4}\qquad(t\ge1).
\end{equation}

Let $W_{\sigma,r}$, $D_\sigma$ and $j_\sigma$ be as in
Lemma~\ref{l:Deq} with $\zeta\equiv1$, and put
$Z=(q_{-\sigma},e_m,H_c)$, $\mathfrak a=\mathfrak a_\sigma(V_c,\mathcal U)$ and
$E=D_\sigma-\bar A_\sigma=E^{\rm base}_\sigma-W_{\sigma,r}$.
Lemma~\ref{l:Deq} gives
\[
 L_\sigma D_\sigma=-\mathcal H_\sigma(V_c,\mathfrak a,Z)
   +\mathcal A_\sigma\Delta\mathcal G
   -2d(V_c)\mathfrak b_\sigma'(V_c)(V_c)_y\mathcal F+R_{\rm c},
 \qquad \mathfrak a=E_y+a_\sigma+j_\sigma .
\]
By $|j_\sigma|\le|(W_{\sigma,r})_y|+C_b|(V_c)_y|$,
Lemmas~\ref{l:pos}(i), \ref{l:nf}(vi), (vii), \ref{l:rec}, \ref{l:pr},
\eqref{e:Aresp} and \eqref{e:Acen},
\[
 \Bigl\Vert T\bigl(j_\sigma^2+|Z|^2+|\mathcal A_\sigma\Delta\mathcal G|
   +|(V_c)_y||\mathcal F|+|R_{\rm c}|+|R^{\rm prof}_\sigma|\bigr)
 \Bigr\Vert_{L^\infty(\Omega)}\le\omega_b(r).
\]
Here $T(|(V_c)_y||\mathcal F|)\le
C\delta\sup_{t\ge r}\Vert\mathcal F(t)\Vert_{L^\infty}\to0$ by Lemmas~\ref{l:rec}
and~\ref{l:bg}, since
$\mathcal F=(\gamma-1)(E_m-c_v\Theta_c)-p_*\varphi_m$.
Moreover, \eqref{e:Ta2} gives $\Vert Ta_\sigma^2\Vert_{L^\infty}\le C_M$,
and Lemma~\ref{l:bg} gives $\sup_{t\ge r}\Vert\mathfrak a(t)\Vert_{L^\infty}\to0$.
On the finite domain $\Omega$, Lemma~\ref{l:pos}(i), \eqref{e:bgcl} and
Lemma~\ref{l:nf}(vi) give
$\Vert TE_y^2\Vert_{L^\infty}\le\int_r^{T_e}\Vert E_y(t)\Vert_{L^\infty}^2\dd t<\infty$.
By \eqref{e:hex}, $\mathfrak a=E_y+a_\sigma+j_\sigma$ and
$(D_\sigma)_y=E_y+a_\sigma$, Young's inequality for products containing $E_y$
and Lemma~\ref{l:pos}(i) for products with $a_\sigma$ give, for late $r$,
\[
 \bigl\Vert T|L_\sigma D_\sigma+k(V_c)(D_\sigma)_y^2|\bigr\Vert_{L^\infty}
 \le\frac{k(2v_-)}2\Vert TE_y^2\Vert_{L^\infty}+\omega_b(r).
\]
Since $L_\sigma\bar A_\sigma=R^{\rm prof}_\sigma-k(V_c)a_\sigma^2$,
\begin{equation}\label{e:bsq}
 T\bigl(k(V_c)E_y^2\bigr)=H_{\mathscr P_\sigma,\Omega}[E|_{\partial_p\Omega}]-E
   +T\bigl(L_\sigma D_\sigma+k(V_c)(D_\sigma)_y^2-R^{\rm prof}_\sigma\bigr).
\end{equation}
Lemmas~\ref{l:all} and~\ref{l:nf}(vi) give
$\sup_{t\ge r}\Vert E(t)\Vert_{L^\infty}\le\omega_b(r)$.
Thus \eqref{e:bsq} and $k(V_c)\ge k(2v_-)>0$ absorb the preceding gradient term and give
$\Vert TE_y^2\Vert_{L^\infty(\Omega)}\le\omega_b(r)$.

Using $\partial_aU_{\rm rec}(n,0,0)=R_\sigma(n)$ from the proof of Lemma~\ref{l:pot},
$\vartheta_\sigma r_\sigma=a_\sigma R_\sigma(v_\sigma)$ and $U_*(V_c)=\widehat{\mathcal U}-(0,H_c,0)$, we write
\[
\begin{aligned}
 R_U={}&\bigl[U_{\rm rec}(V_c,\mathfrak a,Z)-U_{\rm rec}(V_c,a_\sigma,0)\bigr]
   +\bigl[U_{\rm rec}(V_c,a_\sigma,0)-U_*(V_c)-a_\sigma R_\sigma(V_c)\bigr]\\
 &+a_\sigma\bigl(R_\sigma(V_c)-R_\sigma(v_\sigma)\bigr)-(0,H_c,0)
   -\bigl(\widehat{\mathcal U}(y+h_*,t)-\widehat{\mathcal U}\bigr)-\vartheta_{-\sigma}r_{-\sigma},
\end{aligned}
\]
and $\mathfrak a-a_\sigma=E_y+j_\sigma$; therefore
\[
 |R_U|^2\le C\bigl[E_y^2+j_\sigma^2+|Z|^2+a_\sigma^4+|V_c-v_\sigma|^2a_\sigma^2
   +|V_c(y+h_*,t)-V_c|^2+|H_c(y+h_*,t)-H_c|^2+a_{-\sigma}^2\bigr].
\]
The preceding estimates, Lemmas~\ref{l:rec} and~\ref{l:opp}, and \eqref{e:Vca}
bound every response by $\omega_b(r)$. This proves the assertion for
$\mathscr P_\sigma$; Lemma~\ref{l:pos}(v) gives it for $L_\sigma$.
\end{proof}

\subsection{Windows}\label{s:win}

Since $(\gamma-1)c_v\Theta_c=p_*V_c$, one has
$\ell_\sigma(n)(\widehat{\mathcal U}(y+h_*,t)-\widehat{\mathcal U})=\frac12(H_c(y+h_*,t)-H_c)$ and
$\ell_0(\widehat{\mathcal U}(y+h_*,t)-\widehat{\mathcal U})=V_c(y+h_*,t)-V_c$; with $r_\rho=s_\rho R_\rho(v_\rho)$,
$\ell_0r_\rho=0$, the duality in \eqref{e:sfun}, $f_0'(V_c)(V_c)_y=\frac{H_c}2$ and
$K_\sigma'=\mathfrak b_\sigma'G_c$, the definitions \eqref{e:CH}, \eqref{e:pot} and \eqref{e:potq} give,
on $\R\times(0,\infty)$,
\begin{equation}\label{e:rows}
\begin{aligned}
 \mathfrak a_\sigma(V_c,\mathcal U)-a_\sigma&=\ell_\sigma(V_c)R_U+\frac12H_c(y+h_*,t)
   +\sum_{\rho=\pm1}a_\rho\ell_\sigma(V_c)\bigl(R_\rho(v_\rho)-R_\rho(V_c)\bigr),\\
 \partial_y\bigl[\Psi_\sigma-f_0(V_c)-K_\sigma(V_c)-\bar A_\sigma\bigr]
 &=\mathfrak a_\sigma(V_c,\mathcal U)-a_\sigma-\frac{H_c}2+\mathfrak b_\sigma'(V_c)(V_c)_y\Delta\mathcal G,\\
 \partial_y\bigl[\mathcal N_m-h_*(V_c-v_-)\bigr]&=\ell_0R_U+V_c(y+h_*,t)-V_c-h_*(V_c)_y .
\end{aligned}
\end{equation}
For $n\in[\frac{v_-}2,2v_-]$ and $X\in\R^3$, the duality gives
$X=(\ell_+(n)X)R_+(n)+(\ell_-(n)X)R_-(n)+(\ell_0X)R_0$, hence, with $C$ depending only on the
fixed data,
\begin{equation}\label{e:rinv}
 |X|\le C\bigl(|\ell_+(n)X|+|\ell_-(n)X|+|\ell_0X|\bigr).
\end{equation}
Since $\mathcal U-\widehat{\mathcal U}=(\varphi_m,\upsilon_m,c_v\vartheta_m+\frac{b_m^2}2)$ and
$b_m=\upsilon_m+H_c$, Lemmas~\ref{l:pr}, \ref{l:cmp} and \eqref{eq:pr-gau} give, for $t\ge S_0$,
\begin{equation}\label{e:ZR}
 |R_U|\le C|Z_m|+C_M\tau^{-\frac12},\qquad |Z_m|\le C|R_U|+C_M\tau^{-\frac12}.
\end{equation}

Put $\bar c=10(1+c_-+c_+)$. For $z_n\in\R$, $t_n>0$ and $H_n\ge1$ put
\begin{equation}\label{e:win}
\begin{gathered}
 r_n=t_n-H_n,\qquad S_n=J_n\times[t_n-2H_n,t_n],\\
 J_n=[z_n-16\bar cH_n,z_n+16\bar cH_n],\qquad I_n=(z_n-12\bar cH_n,z_n+12\bar cH_n).
\end{gathered}
\end{equation}
Let $\varpi$ be the smaller of the numbers $\varpi$ given by Lemma~\ref{l:kern} for
$D_0=d(\rho_-^{-1})$ and for $D_0=c\rho_-$, with the upper bound $1$ for $A_0$.

\begin{lemma}\label{l:win}
Let $\Lambda,M_0>0$, $t_n\to\infty$, $\alpha_n\to0$ with $\alpha_n>0$, $z_n\in\R$ and
$H_n=\Lambda\alpha_n^{-2}$, and assume that, for all large $n$, $4H_n\le t_n$ and
$\sup_{S_n}|Z_m|\le M_0\alpha_n$. Then
\begin{enumerate}[label=\textup{(\alph*)}]
\item
\[
\begin{gathered}
 \sup_{r_n\le t\le t_n}\bigl(\Vert\partial_yZ_m(t)\Vert_{L^2(I_n)}+\Vert\Gamma(t)\Vert_{L^2(I_n)}\bigr)
 =o(\alpha_n),\\
 \int_{r_n}^{t_n}\!\!\int_{I_n}\Bigl[|\partial_y^2\upsilon_m|^2
   +\Bigl|\partial_y^2\log\frac{\Theta_b}{\Theta_c}\Bigr|^2\Bigr]\dd y\dd t=o(\alpha_n^2);
\end{gathered}
\]
\item if also $4H_n\le\varpi t_n$ for all large $n$, then, with $C$ depending only on the
fixed data,
\begin{equation}\label{e:cpt}
 |R_U(z_n,t_n)|\le C\delta\sup_{S_n}|R_U|+o(\alpha_n).
\end{equation}
\end{enumerate}
Here $o(\cdot)$ refers to $n\to\infty$ and may depend on $\Lambda$, $M_0$ and the initial datum.
\end{lemma}
\begin{proof}
Here $C_b$ may also depend on $\Lambda$ and $M_0$. For large $n$, on $S_n$,
$t\ge t_n-2H_n\ge\frac{t_n}2$, hence
\[
 \tau^{-1}\le\frac2{t_n}\le\frac{\alpha_n^2}{2\Lambda},\qquad \int_{t_n-2H_n}^{t_n}\tau^{-1}\dd t\le1 .
\]

Put
\[
 Q=\partial_y\log\frac{V_b}{V_c},\qquad K=\partial_y\log\frac{\Theta_b}{\Theta_c},\qquad
 h=\frac{(\upsilon_m)_y}{V_b},\qquad A=(H_c)_y\Bigl(\frac1{V_b}-\frac1{V_c}\Bigr).
\]
By \eqref{eq:prof}, $H_c=\frac{c(V_c)_y}{V_c}$ and $\Theta_c=\frac{p_*V_c}R$; by \eqref{eq:pr-id},
$(\log V_c)_t=(\log\Theta_c)_t=\frac{(H_c)_y}{V_c}$ and $(H_c)_t=c\,\partial_y\frac{(H_c)_y}{V_c}$. With
$\frac\kappa{cc_v}=\gamma$, the definitions and \eqref{eq:bg-m}--\eqref{eq:bg-en} give
\begin{equation}\label{e:bmat}
\begin{gathered}
 (V_b)_y=V_b\Bigl(Q+\frac{H_c}c\Bigr),\qquad (\Theta_b)_y=\Theta_b\Bigl(K+\frac{H_c}c\Bigr),\qquad
 (P_b)_y=P_b(K-Q),\qquad \Bigl(\log\frac{V_b}{V_c}\Bigr)_t=h+A,\\
 Q_t=h_y+A_y,\qquad
 (\upsilon_m)_t=\nu h_y-P_b(K-Q)+\nu A_y+(\nu-c)\,\partial_y\frac{(H_c)_y}{V_c},\qquad
 K_t=-(\gamma-1)h_y+F_y,\\
 (\log\Theta_b)_t=-(\gamma-1)h+F+\frac{(H_c)_y}{V_c},\qquad
 F=\frac\kappa{c_vV_b}\bigl[K_y+(K-Q)K\bigr]+\frac{\nu V_bh^2}{c_v\Theta_b}+B,\\
 B=\frac{\gamma H_c(K-Q)}{V_b}+\frac{2\nu(H_c)_yh}{c_v\Theta_b}
   +\frac{\nu(H_c)_y^2}{c_vV_b\Theta_b}+A .
\end{gathered}
\end{equation}

Let $\psi_n$ be a smooth function of $y$ with $0\le\psi_n\le1$, $\psi_n=1$ on $I_n$, $\psi_n=0$
outside $J_n$ and $(\psi_n)_y^2\le CH_n^{-2}\psi_n$ (the square of a cutoff). Put
$\lambda_n=\frac{(\psi_n)_y}{\psi_n}$ where $\psi_n>0$, $\lambda_n=0$ elsewhere,
$w_n=\frac{\psi_nV_c^2}{V_b^2}$, and
\[
 \mathcal E_n=\frac12\int_\R w_n\Bigl(RQ^2+c_vK^2+\frac{V_b^2h^2}{\Theta_b}\Bigr)\dd y,\qquad
 \mathcal D_n=\int_\R w_n\Bigl(\frac{\nu V_bh_y^2}{\Theta_b}+\frac{\kappa K_y^2}{V_b}\Bigr)\dd y .
\]
Then $(w_n)_t=-2w_n(h+A)$, $(w_n)_y=w_n(\lambda_n-2Q)$,
$\partial_y\frac{w_n(\upsilon_m)_y}{\Theta_b}=\frac{w_nV_b}{\Theta_b}[h_y+(\lambda_n-Q-K)h]$,
$|w_n\lambda_n|\le C_bH_n^{-1}$ and $w_n\lambda_n^2\le C_bH_n^{-2}$. Multiply the equations for $Q$ and
$K$ in \eqref{e:bmat} by $Rw_nQ$ and $c_vw_nK$, and the $y$-derivative of the equation for
$\upsilon_m$ by $\frac{w_n(\upsilon_m)_y}{\Theta_b}$, integrate over $\R$ and integrate by parts. The
terms containing $hQ_y$, $hK_y$ and $h_y(K-Q)$ add up to
$-R\int_\R w_n((K-Q)h)_y\dd y$, and
\[
 \mathcal E_n'+\mathcal D_n=I_1+I_2,
\]
where
\[
\begin{aligned}
 I_1={}&R\int_\R w_n\lambda_nh(K-Q)\dd y-c_p\int_\R w_nhK^2\dd y
   +\frac{\gamma-3}2\int_\R\frac{w_nV_b^2h^3}{\Theta_b}\dd y\\
 &+\nu\int_\R\frac{w_nV_b}{\Theta_b}(Q+K-\lambda_n)hh_y\dd y
   -\nu\int_\R\frac{w_nV_b}{\Theta_b}\bigl[K_y+(\lambda_n-2Q)K\bigr]h^2\dd y\\
 &-\kappa\int_\R\frac{w_n}{V_b}\bigl[(\lambda_n+K-3Q)KK_y+(\lambda_n-2Q)(K-Q)K^2\bigr]\dd y\\
 &-\frac\kappa{2c_v}\int_\R\frac{w_nV_bh^2}{\Theta_b}\bigl[K_y+(K-Q)K\bigr]\dd y
   -\frac\nu{2c_v}\int_\R\frac{w_nV_b^3h^4}{\Theta_b^2}\dd y,\\
 I_2={}&R\int_\R w_nQA_y\dd y
   -\int_\R\frac{w_nV_b}{\Theta_b}\bigl[h_y+(\lambda_n-Q-K)h\bigr]
     \Bigl(\nu A_y+(\nu-c)\,\partial_y\frac{(H_c)_y}{V_c}\Bigr)\dd y\\
 &-\int_\R w_n\Bigl[c_v\bigl(K_y+(\lambda_n-2Q)K\bigr)+\frac{V_b^2h^2}{2\Theta_b}\Bigr]B\dd y
   -\int_\R w_nA\Bigl(RQ^2+c_vK^2+\frac{V_b^2h^2}{\Theta_b}\Bigr)\dd y\\
 &-\frac12\int_\R\frac{w_n(H_c)_yV_b^2h^2}{V_c\Theta_b}\dd y .
\end{aligned}
\]
In the two cubic terms of $I_1$ integrate $h=\frac{(\upsilon_m)_y}{V_b}$ by parts:
\[
\begin{aligned}
 \int_\R w_nhK^2\dd y
 &=-\int_\R\frac{\upsilon_mw_n}{V_b}\Bigl[2KK_y+\Bigl(\lambda_n-3Q-\frac{H_c}c\Bigr)K^2\Bigr]\dd y,\\
 \int_\R\frac{w_nV_b^2h^3}{\Theta_b}\dd y
 &=-\int_\R\frac{\upsilon_mw_nV_b}{\Theta_b}\bigl[2hh_y+(\lambda_n-Q-K)h^2\bigr]\dd y,
\end{aligned}
\]
where $|\upsilon_m|\le M_0\alpha_n$ on $\operatorname{supp}w_n$ for $t_n-2H_n\le t\le t_n$. Young's
inequality, the bounds for $w_n\lambda_n$ and the sign of the quartic term give
\[
 I_1\le\frac14\mathcal D_n+C_b\bigl[\alpha_n^2+\tau^{-1}+\Vert(h,K)\Vert_{L^\infty}^2
   +\alpha_n\Vert(h,K)\Vert_{L^\infty}\bigr]\mathcal E_n
   +\frac{C_b}{H_n}\Vert(Q,K,h)\Vert_{L^2(J_n)}^2 .
\]
By \eqref{eq:pr-gau}, $|H_c|\le C\delta\tau^{-\frac12}$, $|(H_c)_y|+|A|\le C\delta\tau^{-1}$,
$\Vert\partial_y\frac{(H_c)_y}{V_c}\Vert_{L^2}+\Vert(H_c)_y^2\Vert_{L^2}\le C\delta\tau^{-\frac54}$ and
\[
 A_y=-\frac{(H_c)_yQ}{V_b}+\Bigl((H_c)_{yy}-\frac{(H_c)_yH_c}c\Bigr)\Bigl(\frac1{V_b}-\frac1{V_c}\Bigr),
 \qquad |A_y|\le C_b\delta\bigl(\tau^{-1}|Q|+\tau^{-\frac32}|\varphi_m|\bigr);
\]
Young's inequality gives
\[
 |I_2|\le\frac14\mathcal D_n+C_b\bigl[\tau^{-1}+\Vert(h,K)\Vert_{L^\infty}^2\bigr]\mathcal E_n
   +\frac{C_b}{H_n}\Vert(Q,K,h)\Vert_{L^2(J_n)}^2
   +C_b\delta^2\tau^{-2}\Vert\varphi_m\Vert_{L^2}^2+C_b\delta^2\tau^{-\frac52}.
\]
Hence
\begin{equation}\label{e:beng}
 \mathcal E_n'+\frac12\mathcal D_n\le X_n\mathcal E_n+\frac{C_b}{H_n}\Vert(Q,K,h)\Vert_{L^2(J_n)}^2
   +C_b\delta^2\tau^{-2}\Vert\varphi_m\Vert_{L^2}^2+C_b\delta^2\tau^{-\frac52},
\end{equation}
with $X_n=C_b[\alpha_n^2+\tau^{-1}+\Vert(h,K)\Vert_{L^\infty}^2+\alpha_n\Vert(h,K)\Vert_{L^\infty}]$.

By $(\varphi_m)_y=V_bQ+\frac{H_c}c\varphi_m$, $(\vartheta_m)_y=\Theta_bK+\frac{H_c}c\vartheta_m$,
$(\upsilon_m)_y=V_bh$, Lemma~\ref{l:dep}, $\alpha_n^2H_n=\Lambda$ and the Cauchy--Schwarz inequality,
\[
 \int_{t_n-2H_n}^{t_n}X_n\dd t\le C_b,\qquad
 \varepsilon_n=\int_{t_n-2H_n}^{t_n}\Vert(Q,K,h)(t)\Vert_{L^2(J_n)}^2\dd t\longrightarrow0 .
\]
Since $\int_{t_n-2H_n}^{r_n}\mathcal E_n\dd t\le C_b\varepsilon_n$, there is
$s_n\in[t_n-2H_n,r_n]$ with $\mathcal E_n(s_n)\le C_b\varepsilon_nH_n^{-1}$. Gronwall's inequality
for \eqref{e:beng} on $[s_n,t_n]$, $\sup_{t\ge t_n-2H_n}\Vert\varphi_m(t)\Vert_{L^2}\to0$
(Lemma~\ref{l:dep}) and $4H_n\le t_n$ give
\begin{equation}\label{e:btail}
 \sup_{r_n\le t\le t_n}\mathcal E_n(t)+\int_{r_n}^{t_n}\mathcal D_n\dd t
 \le C_b\Bigl[\frac{\varepsilon_n}{H_n}
   +\frac{\delta^2H_n}{t_n^2}\sup_{t\ge t_n-2H_n}\Vert\varphi_m(t)\Vert_{L^2}^2
   +\frac{\delta^2H_n}{t_n^{\frac52}}\Bigr]=o(\alpha_n^2).
\end{equation}

On $I_n$, $w_n\ge C_b^{-1}$. By the identities above, $(b_m)_y=V_bh+(H_c)_y$, $H_T-H_c=cK$,
$\frac{m_y}m=Q+\frac{K-Q}\gamma+\frac{H_c}c$, $|e_m|+|q_\pm|+|P_b-p_*|\le C|Z_m|$, and the smoothness
of $J_a$ with $J_a=0$ at $P_b=p_*$,
\[
 |\partial_yZ_m|+|\Gamma|\le C_b\bigl(|Q|+|K|+|h|+|H_c||Z_m|+|(H_c)_y|\bigr),\qquad
 (\upsilon_m)_{yy}=V_bh_y+V_b\Bigl(Q+\frac{H_c}c\Bigr)h .
\]
With \eqref{e:btail}, $\sup_{t\ge r_n}\Vert Z_m(t)\Vert_{L^2}\to0$,
$\Vert(H_c)_y\Vert_{L^2}\le C\delta\tau^{-\frac34}$, $\tau^{-\frac12}\le C_b\alpha_n$ and
$\int_{r_n}^{t_n}\Vert Q\Vert_{L^\infty}^2\dd t\le C_b$ (Lemma~\ref{l:dep}), this proves (a).

For (b) fix $\sigma$, let $S_0$ be as in the proof of Lemma~\ref{l:sq}, $r_n\ge S_0$, and put
$\mathfrak a=\mathfrak a_\sigma(V_c,\mathcal U)$, $Z=(q_{-\sigma},e_m,H_c)$ and $z=\mathcal U-U_*(V_c)$. By
\eqref{e:ZR}, \eqref{e:rows}, Lemma~\ref{l:pr} and $\tau^{-\frac12}\le C_b\alpha_n$,
\[
 |R_U|+|Z|+|\mathfrak a|+|\mathfrak a-a_\sigma|+|z|+|\chi_m-V_c|\le C_b\alpha_n\qquad\text{on }S_n .
\]
The drifts of $L_\sigma$ and $P_V$ in the coordinate $\xi=F_c(y,t)$ of \eqref{e:nco} satisfy
\eqref{e:kH} on $[r_n,t_n]$ with $A_0=A_1=C\delta\le1$ (after decreasing $\delta_1(M)$) and
$\bar b=s_\sigma(\rho_-^{-1})$,
respectively $\bar b=0$. Since $t_n-r_n=H_n\le\varpi\tau_{r_n}$, Lemma~\ref{l:kern}(a), (c) apply on
$[r_n,t_n]$, with constants depending only on the fixed data. Let $\mathcal S$ be the filter of
Lemma~\ref{l:kern}(c), $\xi_n=F_c(z_n,t_n)$, and, for a function $g$ on $\R$,
$\mathcal S_ng=(\mathcal S\hat g)(\xi_n)$ with $\hat g(F_c(y,t_n))=\frac{g(y)}{f_c(y,t_n)}$. For
$\mathscr L\in\{L_\sigma,P_V\}$, $T=T_{\mathscr L,\R\times(r_n,t_n)}$ and $u=Tf$, write
$u(y,t)=\tilde u(F_c(y,t),t)$ and $f(y,t)=\tilde f(F_c(y,t),t)$; Lemma~\ref{l:kern}(c) represents
$\tilde u$ by the kernel $K$ and $\tilde f$, $u_y=f_c\tilde u_\xi$ and
$\mathcal S_n[u_y(\cdot,t_n)]=(\mathcal S\tilde u_\xi)(\xi_n,t_n)$. Hence, for $j\ge1$, with $C$ depending only on
the fixed data, $C_j$ only on $j$ and the fixed data, and suprema over $r_n<s<t_n$,
\begin{equation}\label{e:Sn}
\begin{gathered}
 |u_y(y,t_n)|\le C\bigl(\Vert f\Vert_{L^\infty}\Vert T|f|\Vert_{L^\infty}\bigr)^{\frac12},\\
 |\mathcal S_n[u_y(\cdot,t_n)]|\le C\sup_s\bigl(H_n^{\frac14}\Vert f_1(s)\Vert_{L^2}
   +\log(2+H_n)\Vert f_2(s)\Vert_{L^1}+\Vert f_3(s)\Vert_{L^2}\bigr)
   \quad\text{if }f=f_1+f_2+\partial_yf_3,\\
 |\mathcal S_n[u_y(\cdot,t_n)]|\le C_jH_n^{-j}\sup_s\bigl(\Vert f^0(s)\Vert_{L^1}+\Vert f^1(s)\Vert_{L^\infty}\bigr)
   \quad\text{if }f=f^0+\partial_yf^1,\\
 |v_y(y,t_n)|\le CH_n^{-\frac12}\Vert v(r_n)\Vert_{L^\infty}\quad\text{if }\mathscr Lv=0
   \text{ on }\R\times(r_n,t_n),\\
 |g(z_n)-f_c(z_n,t_n)\mathcal S_ng|\le C\bigl(\Vert g_y\Vert_{L^2}+\delta t_n^{-\frac12}\Vert g\Vert_{L^2}\bigr),
 \qquad |\mathcal S_ng|\le C\min\bigl\{\Vert g\Vert_{L^\infty},\Vert g\Vert_{L^2}\bigr\},\\
 |\mathcal S_ng|+|\mathcal S_n(g_y)|\le C_jH_n^{-j}\Vert g\Vert_{L^\infty}
   \quad\text{if }g=0\text{ on }(z_n-2\bar cH_n,z_n+2\bar cH_n).
\end{gathered}
\end{equation}
where in the third line $f^0$ and $f^1$ vanish unless $2\bar cH_n\le|y-z_n|\le4\bar cH_n$.
The first line is \eqref{e:kres} with $\rho=\infty$. The second follows from the first three
kernel bounds of Lemma~\ref{l:kern}(c), $\Vert\tilde f(s)\Vert_{L^p}\le C\Vert f(s)\Vert_{L^p}$ and,
by \eqref{e:xder}, $\widetilde{\partial_yf_3}=\partial_\xi(f_c\tilde f_3)-(\partial_\xi f_c)\tilde f_3$
with $|\partial_\xi f_c|\le C\delta\tau^{-\frac12}$ and $H_n^{\frac14}\tau_{r_n}^{-\frac12}\le1$. For
the third: if $|\xi-\xi_n|\le\frac{\bar cH_n}2$ and $\eta$ lies in the support of $\tilde f(\cdot,s)$,
then $|\xi-\eta-\bar b(t_n-s)|\ge\frac{\bar cH_n}2$, since $f_c\in[\frac1{\sqrt2},\sqrt2]$,
$|(F_c)_t|\le C\delta$ and $|\bar b|\le\frac{\bar c}{10}$; so Lemma~\ref{l:kern}(a) gives
$|\tilde u_\xi(\xi,t_n)|\le C_jH_n^{-j}\sup_s(\Vert f^0(s)\Vert_{L^1}+\Vert f^1(s)\Vert_{L^\infty})$ there,
while $\Vert u(t_n)\Vert_{L^\infty}\le C\sqrt{H_n}\sup_s(\Vert f^0(s)\Vert_{L^1}+\Vert f^1(s)\Vert_{L^\infty})$ and
the kernel of $\mathcal S$ is $O(H_n^{-j})$ with its derivative for $|\xi-\xi_n|\ge\frac{\bar cH_n}4$. The
fourth follows from $\sup_\xi\Vert K_\xi(t_n,r_n;\xi,\cdot)\Vert_{L^1}\le CH_n^{-\frac12}$, the fifth from
$\Vert\hat g-\mathcal S\hat g\Vert_{L^\infty}\le C\Vert\partial_\xi\hat g\Vert_{L^2}$, and the last from the decay
of the kernel of $\mathcal S$ and $|F_c(y,t_n)-\xi_n|\ge\frac{|y-z_n|}{\sqrt2}$.

Let $\zeta_n,\eta_n$ be smooth functions of $y$ with values in $[0,1]$, $\zeta_n=1$ for
$|y-z_n|\le8\bar cH_n$, $\zeta_n=0$ for $|y-z_n|\ge11\bar cH_n$, $\eta_n=1$ for
$|y-z_n|\le2\bar cH_n$, $\eta_n=0$ for $|y-z_n|\ge4\bar cH_n$, and
$|\partial_y^j\zeta_n|+|\partial_y^j\eta_n|\le CH_n^{-j}$ for $j\le2$. Let $u_n$, $W_n$ and $D_\sigma$ be
given by Lemma~\ref{l:Deq} with $r=r_n$, $\zeta=\zeta_n$ and
$\Omega_n=\{|y-z_n|<8\bar cH_n\}\times(r_n,\infty)$, and put
$E_n=D_\sigma-\bar A_\sigma=E^{\rm base}_\sigma-W_n$. By \eqref{e:nfs} and Lemma~\ref{l:nf}(i),
$|F_{\rm bulk}|\le C(|\widetilde z|^2+|z_y|)$, $|G_{\rm bulk}|\le C(|z|^2|z_y|+|z_y|^2)$ and
$|\partial_y[\mathcal P_\sigma(V_c,z)]|\le C(|z||(V_c)_y|+|z_y|)$, where
$|\widetilde z|\le C|Z_m|$ and $|z_y|\le C(|\partial_yZ_m|+|(b_m)_y|)$. Part (a),
$\sup_{t\ge r_n}\Vert Z_m(t)\Vert_{L^2}\to0$ and $\Vert(V_c)_y\Vert_{L^2}\le C\delta\tau^{-\frac14}$ give
\[
 \sup_{r_n<t<t_n}\Bigl(\Vert\zeta_nF_{\rm bulk}\Vert_{L^2}+\frac{\Vert\zeta_nG_{\rm bulk}\Vert_{L^1}}{\alpha_n}
   +\Vert\partial_y[\mathcal P_\sigma(V_c,z)]\Vert_{L^2(I_n)}\Bigr)=o(\alpha_n).
\]
With the second line of \eqref{e:Sn} for $u_n$, the splitting
$\partial_y[\mathcal P_\sigma(V_c,z)]=\zeta_n\partial_y[\mathcal P_\sigma(V_c,z)]
+\partial_y[(1-\zeta_n)\mathcal P_\sigma(V_c,z)]+(\zeta_n)_y\mathcal P_\sigma(V_c,z)$ and the last line of
\eqref{e:Sn}, $|\mathcal S_n[(W_n)_y(\cdot,t_n)]|=o(\alpha_n)$. By \eqref{e:nfe} and the tails in
Lemma~\ref{l:nf}(iii), $\sup_{t\ge r_n}\Vert u_n(t)\Vert_{H^1}\to0$; with Lemmas~\ref{l:all} and
\ref{l:nf}(iii), $\sup_{t\ge r_n}\Vert E_n(t)\Vert_{L^\infty}\to0$ and
$\sup_{t\ge r_n}\Vert W_n(t)\Vert_{L^\infty}\to0$.

By Lemmas~\ref{l:Deq} and \ref{l:pr}, on $\Omega_n$
\begin{equation}\label{e:Erow}
 L_\sigma E_n=\mathfrak j_\sigma(V_c,\mathcal U)+k(V_c)a_\sigma^2-R^{\rm prof}_\sigma+\mathcal A_\sigma\Delta\mathcal G
   -2d(V_c)\mathfrak b_\sigma'(V_c)(V_c)_y\mathcal F+R_{\rm c},
\end{equation}
and, with $d=d(V_c)$ and $s=s_\sigma(V_c)$,
$L_\sigma(\eta_nE_n)=\eta_nL_\sigma E_n+(s(\eta_n)_y-d(\eta_n)_{yy}+2(d(\eta_n)_y)_y)E_n
-2\partial_y(d(\eta_n)_yE_n)$. On $\R\times[r_n,t_n]$, $\eta_nE_n$ is the sum of the solution of
$L_\sigma v=0$, $v(r_n)=\eta_nE_n(r_n)$, and the responses $T$ to the terms of
$L_\sigma(\eta_nE_n)$, with $\eta_nL_\sigma E_n$ given by \eqref{e:Erow}; these contribute to $\mathcal S_n[\partial_y(\eta_nE_n)(\cdot,t_n)]$ as follows.
The function $v$ contributes $CH_n^{-\frac12}\Vert E_n(r_n)\Vert_{L^\infty}=o(\alpha_n)$, by the fourth
line of \eqref{e:Sn}. The cutoff terms contribute $C_jH_n^{-j}\sup_s\Vert E_n(s)\Vert_{L^\infty}$, by the
third line. By $-\mathfrak j_\sigma=\mathcal H_\sigma(V_c,\mathfrak a,Z)$ and \eqref{e:hex},
\[
 |\mathfrak j_\sigma(V_c,\mathcal U)+k(V_c)a_\sigma^2|\le C\bigl(|\mathfrak a|^3+|\mathfrak a||Z|+|Z|^2\bigr)
   +k(V_c)|\mathfrak a-a_\sigma|\bigl(|\mathfrak a-a_\sigma|+2|a_\sigma|\bigr);
\]
on $\operatorname{supp}\eta_n$ this is at most $C_b\alpha_n^2$, and by \eqref{e:rows},
Lemmas~\ref{l:sq}, \ref{l:opp}, \ref{l:rec}, \ref{l:pos}(i), \eqref{e:Ta2} and \eqref{e:Vca} its
response is at most $\omega_b(r_n)$; the first line of \eqref{e:Sn} gives the contribution
$C_b\alpha_n\omega_b(r_n)^{\frac12}=o(\alpha_n)$. In the same way $\eta_nR^{\rm prof}_\sigma$ contributes
$C_M\delta e^{-c_er_n}$, and $\eta_nd(V_c)\mathfrak b_\sigma'(V_c)(V_c)_y\mathcal F$, which is at most
$C_b\alpha_n\tau^{-\frac12}\le C_b\alpha_n^2$ by \eqref{e:hsplit} and has response at most $C_b\delta\alpha_n$
by Lemma~\ref{l:rec}, contributes $o(\alpha_n)$. For $R_{\rm c}$ let $F_1$ and $F_2=C_M|z_y||(V_c)_y|$ be
the sum of the other terms and the term with $z_y$ on the right of \eqref{e:nfr}, and
$R_{\rm c}=R_{\rm c}^1+R_{\rm c}^2$ with $R_{\rm c}^i=\frac{F_i}{F_1+F_2}R_{\rm c}$. On
$\operatorname{supp}\eta_n$, $|R_{\rm c}^1|\le F_1\le C_b\alpha_n^2$, and by Lemma~\ref{l:rec}, \eqref{e:pt},
\eqref{eq:pr-gau} and the third bound of \eqref{e:gauss}, $T F_1\le C_M\delta\sup_{t\ge r_n}\Vert z(t)\Vert_{L^\infty}
+C_M\delta\tau_{r_n}^{-\frac12}$; by the second line of \eqref{e:Sn},
$\eta_nR_{\rm c}^2$ contributes $CH_n^{\frac14}C_M\delta\tau_{r_n}^{-\frac12}\sup_s\Vert z_y(s)\Vert_{L^2(I_n)}
=o(\alpha_n)$. Hence all these terms contribute $o(\alpha_n)$.

Let $V_n=T(\eta_n\mathcal A_\sigma\Delta\mathcal G)$ with $\mathscr L=L_\sigma$. In the coordinate $\xi$,
$(\partial_t-D_0\partial_\xi^2+b\partial_\xi+b_\xi)\partial_\xi\tilde V_n=\widetilde{f_c^{-1}f_y}$ with
$\partial_\xi\tilde V_n(r_n)=0$ and
\[
 f_y=\eta_n\bigl[(\partial_y\mathcal A_\sigma)\Delta\mathcal G+\mathcal A_\sigma\mathcal F\bigr]
   +(\eta_n)_y\mathcal A_\sigma\Delta\mathcal G .
\]
Since $|b_\xi|\le C\delta\tau^{-1}$ and $\tau_{t_n}\le2\tau_{r_n}$, comparison gives
$|\partial_yV_n(t_n)|\le4T|f_y|(t_n)$. By \eqref{e:CH} and $(\gamma-1)c_v\Theta_c=p_*V_c$,
\begin{equation}\label{e:FRU}
 \mathcal F=\gamma p_*\sum_{\rho=\pm1}\vartheta_\rho+(\gamma-1)(R_U)_3-p_*(R_U)_1 .
\end{equation}
By \eqref{e:Aresp} with $C$ depending only on the fixed data, \eqref{e:Acen}, Lemma~\ref{l:pos}(i),
and $|\mathcal A_\sigma\vartheta_\rho|\le C_M\delta e^{-c_et}$ (by the splitting used for \eqref{e:Vca}),
\[
 T|f_y|\le C\delta(\gamma-1+p_*)\sup_{S_n}|R_U|+C\delta\varepsilon_c(r_n)\tau_{r_n}^{-\frac12}
   +C_be^{-\frac{c_Ga_0^2r_n}4}+C_M\delta e^{-c_er_n},
\]
hence, with $|\mathcal S_ng|\le C\Vert g\Vert_{L^\infty}$,
$|\mathcal S_n[\partial_yV_n(\cdot,t_n)]|\le C\delta\sup_{S_n}|R_U|+o(\alpha_n)$, with $C$ depending
only on the fixed data. Altogether,
\[
 \bigl|\mathcal S_n[\partial_y(\eta_nE_n)(\cdot,t_n)]\bigr|\le C\delta\sup_{S_n}|R_U|+o(\alpha_n).
\]
By the first two lines of \eqref{e:rows},
\[
 \ell_\sigma(V_c)R_U=\partial_yE_n+\partial_yW_n-\frac12\bigl(H_c(y+h_*,t)-H_c\bigr)
   -\mathfrak b_\sigma'(V_c)(V_c)_y\Delta\mathcal G
   -\sum_{\rho=\pm1}a_\rho\ell_\sigma(V_c)\bigl(R_\rho(v_\rho)-R_\rho(V_c)\bigr).
\]
Let $g=\eta_n\ell_\sigma(V_c)R_U$ at $t=t_n$. With
$\eta_n\partial_yE_n=\partial_y(\eta_nE_n)-(\eta_n)_yE_n$,
$\eta_n\partial_yW_n=\partial_yW_n-\partial_y((1-\eta_n)W_n)-(\eta_n)_yW_n$, the last three lines of
\eqref{e:Sn}, $|H_c(y+h_*,t)-H_c|\le C_M\tau^{-1}$, \eqref{e:Acen} and \eqref{e:Vca},
$|\mathcal S_ng|\le C\delta\sup_{S_n}|R_U|+o(\alpha_n)$. Since $\Vert R_U(t)\Vert_{L^2}\to0$ and
$\Vert\partial_yR_U(t_n)\Vert_{L^2(I_n)}=o(\alpha_n)$, by Lemmas~\ref{l:dep}, \ref{l:pr}, \eqref{eq:pr-gau}
and (a), the fifth line of \eqref{e:Sn} gives
\[
 |\ell_\sigma(V_c)R_U(z_n,t_n)|\le C\delta\sup_{S_n}|R_U|+o(\alpha_n).
\]

For the third row let $u_N$ solve $P_\chi u_N=\partial_y(\zeta_nF_N)$ on $\R\times(r_n,\infty)$,
$u_N(r_n)=0$, with $P_\chi$, $F_N$, $R_N$, $d_0$, $b_0$, $g_\chi$ from \eqref{e:npr}, and put
$C_N=\frac{(\gamma-1)c}{\gamma p_*}b_m+\frac{H_c\mathcal B}{\gamma p_*}-u_N$, $D_N=\mathcal N_m+C_N-h_*(V_c-v_-)$,
$e_N=\chi_m-V_c-h_*(V_c)_y$, $\Delta_d=\frac c{\chi_m}-\frac c{V_c}$ and
$\Delta_b=\frac{c(V_c)_y}{V_c^2}-b_0$. By \eqref{e:pot} and Lemma~\ref{l:neu}, $(D_N)_y=e_N+(C_N)_y$
and, on $\Omega_n$,
\begin{equation}\label{e:Nrow}
 P_VD_N=R_N-\Bigl(1-\frac{V_c}{\chi_m}\Bigr)h_*(V_c)_t+\partial_y\bigl(\Delta_d(D_N)_y\bigr)
   +\bigl(\Delta_b-(\Delta_d)_y\bigr)(D_N)_y .
\end{equation}
By \eqref{e:npr}, \eqref{eq:pr-gau}, Lemma~\ref{l:dep} and (a),
$\sup_{r_n<t<t_n}\Vert\zeta_nF_N\Vert_{L^2}=o(\alpha_n)$ and
$\Vert\partial_y(\zeta_nF_N)\Vert_{L^2(\R\times(r_n,t_n))}=o(\alpha_n)$. The first line of
\eqref{e:pchie}, Gronwall's inequality with $\int_{r_n}^\infty\Vert g_\chi\Vert_{L^\infty}^2\dd t\le C_b$,
and then the second line with $|b_0|\le C\delta\tau^{-\frac12}$ give
\[
\begin{gathered}
 \sup_{r_n\le t\le t_n}\Vert u_N\Vert_{L^2}^2+\int_{r_n}^{t_n}\Vert(u_N)_y\Vert_{L^2}^2\dd t
 \le C_bH_n\,o(\alpha_n^2)=o(1),\\
 \sup_{r_n\le t\le t_n}\Vert(u_N)_y\Vert_{L^2}^2
 \le C\delta^2\tau_{r_n}^{-1}\int_{r_n}^{t_n}\Vert(u_N)_y\Vert_{L^2}^2\dd t+o(\alpha_n^2)=o(\alpha_n^2),
\end{gathered}
\]
hence, with \eqref{e:npplim}, Lemma~\ref{l:bd} and $\Vert b_m(t)\Vert_{L^\infty}\to0$,
$\sup_{r_n\le t\le t_n}\Vert D_N(t)\Vert_{L^\infty}=o(1)$ and
$\sup_{r_n\le t\le t_n}\Vert(C_N)_y\Vert_{L^2(I_n)}=o(\alpha_n)$. On $S_n$, $|\Delta_d|+|e_N|\le C_b\alpha_n$; moreover
$\Vert\Delta_b\Vert_{L^2}+\Vert(\Delta_d)_y\Vert_{L^2(I_n)}=o(\alpha_n)$ and $\Vert e_N\Vert_{L^2}=o(1)$ for
$r_n\le t\le t_n$, by (a), \eqref{eq:pr-gau}, $|h_*|\le C_M\delta^{-1}$ and
$\Vert\chi_m-V_c\Vert_{L^2}\to0$. Write $\eta_nD_N$ on
$\R\times[r_n,t_n]$ as for $\eta_nE_n$, now with $\mathscr L=P_V$ and \eqref{e:Nrow}. By Lemma~\ref{l:neu},
$|R_N|\le F_1+F_2$ with $F_1=C_M(q_+^2+q_-^2)+C_b\tau^{-\frac32}$ and $F_2=C_M|\Gamma|^2$; split $R_N$ as
$R_{\rm c}$ above. On $\operatorname{supp}\eta_n$, $F_1\le C_b\alpha_n^2$, and, since the drift of $P_V$
is at most $c_{\rm ac}$ in absolute value, Lemmas~\ref{l:rec} and \ref{l:pos}(i) give
$TF_1\le\omega_b(r_n)$; by (a), $\sup_s\Vert\eta_nF_2(s)\Vert_{L^1}=o(\alpha_n^2)$. The term with
$(V_c)_t$ is at most $C_b\alpha_n\tau^{-1}$ on $\operatorname{supp}\eta_n$, and $T$ of it at most
$C_b\alpha_nH_n\tau_{r_n}^{-1}$. In the flux $\partial_y(\eta_n\Delta_d(D_N)_y)$,
$\sup_s\Vert\eta_n\Delta_d(D_N)_y\Vert_{L^2}\le C_b\alpha_n\Vert e_N\Vert_{L^2}+C_b\sup_s\Vert(C_N)_y\Vert_{L^2(I_n)}=o(\alpha_n)$.
In $\eta_n(\Delta_b-(\Delta_d)_y)(D_N)_y$, the part with $e_N$ has $L^2$ norm at most
$C_b\alpha_n\,o(\alpha_n)$, and the part with $(C_N)_y$ has $L^1$ norm $o(\alpha_n^2)$; the term
$-(\eta_n)_y\Delta_d(D_N)_y$ has $L^1$ norm at most $C_b\alpha_nH_n^{-\frac12}$. With
$\alpha_nH_n^{\frac14}+\alpha_n\log(2+H_n)\to0$, \eqref{e:Sn} gives
$|\mathcal S_n[\partial_y(\eta_nD_N)(\cdot,t_n)]|=o(\alpha_n)$, the cutoff terms and the initial
value being treated as for $\eta_nE_n$. By the third line of \eqref{e:rows},
\[
 \ell_0R_U=(D_N)_y-(C_N)_y-\bigl(V_c(y+h_*,t)-V_c-h_*(V_c)_y\bigr),\qquad
 |V_c(y+h_*,t)-V_c-h_*(V_c)_y|\le C_b\tau^{-1},
\]
and \eqref{e:Sn} applied to $\eta_n\ell_0R_U$ at $t=t_n$ gives $|\ell_0R_U(z_n,t_n)|=o(\alpha_n)$.

By \eqref{e:rinv} with $n=V_c(z_n,t_n)$, the bounds for $\ell_\pm(V_c)R_U$ and $\ell_0R_U$ at
$(z_n,t_n)$ give $|R_U(z_n,t_n)|\le C\delta\sup_{S_n}|R_U|+o(\alpha_n)$, with $C$ depending only on the
fixed data, hence \eqref{e:cpt}.
\end{proof}

\subsection{Asymptotic expansion and Eulerian decay}\label{s:asy}

\begin{lemma}\label{l:dbl}
Let $k\ge1$, $p>0$, $c_1>0$, let $f:\R\times[0,\infty)\to\R^k$ be continuous and bounded on
$\R\times[0,T_e]$ for every $T_e<\infty$, with $\Vert f(t)\Vert_{L^\infty}\to0$ as $t\to\infty$, and let
$(z_n^0,t_n^0)$ satisfy
\[
 t_n^0\to\infty,\qquad a_n=|f(z_n^0,t_n^0)|>0,\qquad t_n^0a_n^p\to\infty .
\]
After deleting finitely many $n$, there are points $(z_n,t_n)$ such that, with
$\alpha_n=|f(z_n,t_n)|$, $H_n=\alpha_n^{-p}$ and $S_n$ given by \eqref{e:win} with $\bar c$ replaced
by $c_1$,
\[
\begin{gathered}
 \alpha_n\ge a_n,\qquad 0\le t_n^0-t_n\le\frac{2^{p+1}}{2^p-1}a_n^{-p},\qquad
 |z_n-z_n^0|\le\frac{2^{p+4}c_1}{2^p-1}a_n^{-p},\qquad \sup_{S_n}|f|\le2\alpha_n,\\
 \alpha_n\to0,\qquad t_n\alpha_n^p\to\infty,\qquad \frac{t_n}{t_n^0}\to1,\qquad
 \sup_{(y,t)\in S_n}\frac{|y-z_n^0|+t_n^0-t}{t_n^0}\to0 .
\end{gathered}
\]
In particular, if $c_2>0$ and $|z_n^0|\le c_2t_n^0$ for all $n$, then $|y|\le2c_2t$ on $S_n$ for all
large $n$.
\end{lemma}
\begin{proof}
Delete the $n$ with $t_n^0\le(\frac{2^{p+1}}{2^p-1}+2)a_n^{-p}$. Fix $n$ and put
$(y_0,s_0)=(z_n^0,t_n^0)$, $b_0=a_n$. If the box
$[y_j-16c_1b_j^{-p},y_j+16c_1b_j^{-p}]\times[s_j-2b_j^{-p},s_j]$ contains a point $(y_{j+1},s_{j+1})$
with $b_{j+1}=|f(y_{j+1},s_{j+1})|>2b_j$, choose it. Then $b_j\ge2^ja_n$, so
\[
 0\le t_n^0-s_j\le2a_n^{-p}\sum_{i\ge0}2^{-ip}=\frac{2^{p+1}}{2^p-1}a_n^{-p},\qquad
 |y_j-z_n^0|\le16c_1a_n^{-p}\sum_{i\ge0}2^{-ip}=\frac{2^{p+4}c_1}{2^p-1}a_n^{-p},
\]
all boxes lie in $\{t>0\}$, and, since $f$ is bounded on $\R\times[0,t_n^0]$ and $b_j\to\infty$ otherwise,
the choice stops at some $j$; put $(z_n,t_n)=(y_j,s_j)$. This gives the first line. For
$(y,t)\in S_n$, the first line gives
$|y-z_n^0|+t_n^0-t\le(\frac{2^{p+4}c_1+2^{p+1}}{2^p-1}+16c_1+2)a_n^{-p}=o(t_n^0)$, since
$t_n^0a_n^p\to\infty$; hence $\frac{t_n}{t_n^0}\to1$, $t_n\alpha_n^p\ge t_na_n^p\to\infty$ and
$\alpha_n\le\Vert f(t_n)\Vert_{L^\infty}\to0$. If $|z_n^0|\le c_2t_n^0$, then on $S_n$,
$|y|\le c_2t_n^0+o(t_n^0)$ and $t\ge(1-o(1))t_n^0$, so $|y|\le2c_2t$ for large $n$.
\end{proof}

\begin{lemma}\label{l:asy}
\[
 \lim_{t\to\infty}\sqrt{1+t}\,\Vert R_U(t)\Vert_{L^\infty}=0,\qquad
 \sup_{t\ge0}\sqrt{1+t}\,\Vert R_U(t)\Vert_{L^\infty}\le C_b .
\]
\end{lemma}
\begin{proof}
Decrease $\delta_1(M)$ so that $4C\delta_1(M)\le1$, where the constant $C$ in
\eqref{e:cpt} depends only on the fixed data. By
Lemmas~\ref{l:bg}, \ref{l:pr} and \eqref{e:ZR}, $R_U$ is continuous, bounded on $\R\times[0,T_e]$ for
every $T_e<\infty$, and $\Vert R_U(t)\Vert_{L^\infty}\to0$.

Suppose that $\sup_{t\ge0}\sqrt\tau\Vert R_U(t)\Vert_{L^\infty}=\infty$. Then there are
$(z_n^0,t_n^0)$ with $t_n^0\to\infty$ and $t_n^0|R_U(z_n^0,t_n^0)|^2\to\infty$. Lemma~\ref{l:dbl} with
$f=R_U$, $p=2$ and $c_1=\bar c$ gives $(z_n,t_n)$ with $\alpha_n=|R_U(z_n,t_n)|\to0$,
$t_n\alpha_n^2\to\infty$ and $\sup_{S_n}|R_U|\le2\alpha_n$, $S_n$ as in \eqref{e:win} with
$H_n=\alpha_n^{-2}$. On $S_n$, $\tau^{-\frac12}\le2t_n^{-\frac12}=o(\alpha_n)$; hence, for large $n$,
$4H_n\le\varpi t_n$ and, by \eqref{e:ZR}, $\sup_{S_n}|Z_m|\le3C\alpha_n$ with $C$ from \eqref{e:ZR}.
Lemma~\ref{l:win} with $\Lambda=1$ and $M_0=3C$ gives
$\alpha_n\le2C\delta\alpha_n+o(\alpha_n)\le\frac12\alpha_n+o(\alpha_n)$, a contradiction. This proves the
second assertion.

Let $A_\infty=\limsup_{t\to\infty}\sqrt\tau\Vert R_U(t)\Vert_{L^\infty}$, and suppose $A_\infty>0$. Choose
$(z_n,t_n)$ with $t_n\to\infty$ and $\sqrt{\tau_{t_n}}|R_U(z_n,t_n)|\to A_\infty$, and put
$\alpha_n=\tau_{t_n}^{-\frac12}$ and $\Lambda=\frac\varpi8$. Then $H_n=\Lambda\tau_{t_n}$, $4H_n\le\varpi t_n$,
and $\tau\ge(1-2\Lambda)\tau_{t_n}\ge\frac34\tau_{t_n}$ on $S_n$. By the second assertion and
\eqref{e:ZR}, $\sup_{S_n}|Z_m|\le C_b\alpha_n$, and
$\sup_{S_n}|R_U|\le(1-2\Lambda)^{-\frac12}(A_\infty+o(1))\alpha_n$. Lemma~\ref{l:win} with $M_0=C_b$ gives
$A_\infty\le C\delta(1-2\Lambda)^{-\frac12}A_\infty\le\frac14(\frac43)^{\frac12}A_\infty$, a contradiction.
\end{proof}

\begin{lemma}\label{l:rate}
For $t\ge0$,
\begin{equation}\label{e:bmrate}
 \Bigl\Vert\Bigl(\frac1{V_b},U_b,\Theta_b\Bigr)-\Bigl(\frac1{V_c},U_c,\Theta_c\Bigr)\Bigr\Vert_{L^\infty}
   +\Vert P_b-p_*\Vert_{L^\infty}\le C_b\tau^{-\frac12},
\end{equation}
the maps \eqref{eq:m-map} and \eqref{e:Xc} satisfy
\[
 \sup_{t\ge0}\bigl(\Vert X(\cdot,t)-\bar X(\cdot,t)\Vert_{L^\infty(\R_y)}
   +\Vert Y(\cdot,t)-\bar Y(\cdot,t)\Vert_{L^\infty(\R_x)}\bigr)\le C_b,
\]
and, with $\bar Q$ from \eqref{eq:ct-w},
\[
 \Vert (\rho,u,\theta)(t)-\bar Q(t)\Vert_{L^\infty(\R_x)}+\Vert p(t)-p_*\Vert_{L^\infty(\R_x)}
 \le C_b\tau^{-\frac12}.
\]
\end{lemma}
\begin{proof}
By \eqref{e:ZR} and Lemma~\ref{l:asy},
\[
 \Vert Z_m(t)\Vert_{L^\infty}
 \le C\Vert R_U(t)\Vert_{L^\infty}+C_M\tau^{-1/2}
 \le C_b\tau^{-1/2}\qquad(t\ge S_0).
\]
Since $V_b,V_c$ have positive lower bounds and
$P_b-p_*=(R\vartheta_m-p_*\varphi_m)/V_b$, this proves
\eqref{e:bmrate} for $t\ge S_0$; Lemma~\ref{l:bg} covers $[0,S_0]$.

Put $A_v=\int_{-\infty}^y(V_b-V_c)\dd z$. By \eqref{e:pot} and $\ell_0R_0=1$,
$\chi_m=V_b+\frac{\mathcal F}{\gamma p_*}$, hence $A_v=\mathcal N_m-\frac{\mathcal G}{\gamma p_*}$ and
$|A_v|\le C_b$ by Lemmas~\ref{l:bd} and \ref{l:neu}. By \eqref{eq:bg-m}, \eqref{eq:pr-id} and
Lemma~\ref{l:cmp}, $(A_v)_t=\int_{-\infty}^y(U_b-U_c)_y\dd z=U_b-U_c$, and \eqref{eq:m-map},
\eqref{e:Xc} give
\[
 X-\bar X=\int_0^y(V_b-V_c)\dd z+\int_0^t(U_b-U_c)(0,s)\dd s=A_v(y,t)-A_v(0,0).
\]
Since $\bar X_y=V_c\ge\min\{v_-,v_+\}$ and $\bar X(Y(x,t),t)-x=(\bar X-X)(Y(x,t),t)$,
$|Y-\bar Y|\le\frac{\Vert X-\bar X\Vert_{L^\infty}}{\min\{v_-,v_+\}}$. At fixed $x$, with $y=Y(x,t)$,
\eqref{eq:bg-eul} and \eqref{e:ct-pullback} give
\[
\begin{aligned}
 |(\rho,u,\theta)-\bar Q|(x,t)
 &\le\Bigl|\Bigl(\frac1{V_b},U_b,\Theta_b\Bigr)
       -\Bigl(\frac1{V_c},U_c,\Theta_c\Bigr)\Bigr|(y,t)\\
 &\quad+\Bigl\Vert\partial_y\Bigl(\frac1{V_c},U_c,\Theta_c\Bigr)(t)
       \Bigr\Vert_{L^\infty}|Y-\bar Y|(x,t).
\end{aligned}
\]
Also $p-p_*=(P_b-p_*)\circ Y$; \eqref{e:bmrate}, \eqref{eq:pr-gau} and the bound for $Y-\bar Y$ give the
last assertion.
\end{proof}

\begin{corollary}[Sharpness of the decay rate]\label{l:rate-sharp}
Under the assumptions of Lemma~\ref{l:rate}, if $m_\sigma\ne0$ for
some $\sigma\in\{-1,1\}$, then
\[
 \liminf_{t\to\infty}\sqrt{1+t}\,
       \Vert u(t)-\bar u(t)\Vert_{L^\infty(\R)}
 \ge\frac{2|\tanh(\alpha s_\sigma m_\sigma/2)|}
              {\alpha\sqrt{4\pi d_\sigma}}>0.
\]
Consequently the exponent $1/2$ in \eqref{e:r1} cannot be increased
for the full class of integrable perturbations considered here.
\end{corollary}
\begin{proof}
At $y=s_\sigma t$, the explicit formula \eqref{e:CH} gives
\[
 \sqrt\tau\,a_\sigma(s_\sigma t,t)
   =\frac{2\tanh(\alpha s_\sigma m_\sigma/2)}
                {\alpha\sqrt{4\pi d_\sigma}}.
\]
The opposite acoustic wave and the contact velocities at this point
are exponentially small. Since
$b_m=H_c(y+h_*,t)+a_-+a_++(R_U)_2$, Lemma~\ref{l:asy} therefore gives the
same limit for $\sqrt\tau(U_b-U_c)(s_\sigma t,t)$.
By Lemma~\ref{l:rate}, the two mass maps differ by a bounded amount;
by \eqref{eq:pr-gau}, $\Vert(U_c)_y\Vert_\infty\le C\delta\tau^{-1}$.
Thus replacing $U_c(y,t)$ by
$\bar u(X(y,t),t)=U_c(\bar Y(X(y,t),t),t)$ changes the velocity
difference by $O_b(\tau^{-1})$. Evaluate at
$x=X(s_\sigma t,t)$ to obtain the stated Eulerian lower bound.
\end{proof}

\section{Proof of the main theorem}\label{sec:proof}

For each fixed datum choose an auxiliary $M$ as in
\eqref{e:dat-aux}. We verify that every wave-strength restriction
depends only on $A$ and the fixed data. Throughout this section
$C_A$ has only this dependence, whereas $C_b$ and the entrance times
may depend on the complete datum.

\begin{proof}[Proof of Theorem~\ref{thm:planar}]
By \eqref{e:dat} and the one-dimensional Sobolev inequality,
$A^{-1}\le\rho_0\le C_A$, so $X_0,Y_0$ are bi-Lipschitz with
constants controlled by $A$. The chain rule, the Gaussian bounds and
\[
 |Y_0(x)-\bar Y(x,0)|
 \le\Vert e_\rho\Vert\,|x|^{1/2}
\]
give $\Vert\varphi_m(0)\Vert_{H^1}+\Vert\upsilon_m(0)\Vert\le C_A$ and
positive initial mass-coordinate volume and temperature lower bounds.
For positive $a,\widetilde a$ in a fixed compact interval, the identity
$\Phi(rs)=\Phi(r)+(s-1)(r-1)+\Phi(s)$ and
$|r-1|\le C(\Phi(r)^{1/2}+\Phi(r))$ give
\[
 \Phi\left(\frac\theta{\widetilde a}\right)
 \le C\Phi\left(\frac\theta a\right)+C|a-\widetilde a|^2.
\]
Take $a=\bar\theta(x,0)$ and
$\widetilde a=\Theta_c(Y_0(x),0)$.
The Gaussian derivative bound and the fact that $Y_0$ and
$\bar Y(\cdot,0)$ vanish at zero and have positive bounded derivatives
give
\[
 \Vert\bar\theta(\cdot,0)-\Theta_c(Y_0(\cdot),0)\Vert
 \le C_A\delta\Vert e_\rho\Vert.
\]
Changing variables by $\dd y=\rho_0\dd x$ bounds the initial
relative entropy of Proposition~\ref{l:bg-entropy} by $C_A$.
This proposition gives the global solution for a threshold
depending only on $A$. The estimate \eqref{e:bgch} and the
coordinate conversion in Lemma~\ref{l:bg} propagate its prescribed
finite higher regularity with constants $C_b$.
The $L^1$ membership makes $Q$ well-defined and justifies the
conservation argument of Lemma~\ref{l:cmp}. Lemma~\ref{l:cmp} gives
$|m_-|+|m_+|+\delta|h_*|\le C|Q|\le C_A$.
Consequently the diffusion waves, their Cole--Hopf factors and their
inverses in Lemma~\ref{l:pr} have $C_A$ bounds. Finite primitive
sizes and source budgets may depend on the complete datum, including
its possibly large $L^1$ norm.

We next check the large-time estimates. The characteristic box
in Lemma~\ref{l:dep} and the normal-form coefficients in
Lemma~\ref{l:pot} are fixed from the physical data. In
Lemmas~\ref{l:bd}, \ref{l:neu} and~\ref{l:cen}, the surviving
absorptions have coefficients $C\delta$ with fixed-data $C$;
the primitive sizes, source responses and high derivative bounds
occur only in errors that tend to zero with the entrance time.
In particular, the coefficient $C$ in \eqref{e:cpt} depends only on
the fixed data: the only surviving acoustic feedback is \eqref{e:FRU}, estimated
by \eqref{e:Aresp}, while the neutral contribution is $o(\alpha_n)$.
Hence $4C\delta\le1$ closes Lemma~\ref{l:asy} independently of
the higher initial norms. Choose $\delta_0(A)>0$ below these
thresholds. Proposition~\ref{l:bg-entropy} gives the uniform positive
upper and lower bounds. The regularity \eqref{e:brg} follows from
\eqref{e:bgch} and Lemma~\ref{l:bg}, and the decay estimate
\eqref{e:r1} follows from Lemma~\ref{l:rate}.
\end{proof}

\appendix

\section{Existence, regularity and auxiliary balances}\label{app:bg}

All spaces are over $\R_y$, and the system is \eqref{eq:bg-m},
\eqref{eq:bg-mom}, \eqref{e:bgth}. Let $B=(B_v,B_u,B_T)$ be a smooth
function of $y$ with $\partial_yB\in H^k$ for every integer $k\ge0$ and
with $B_v$, $B_T$ between positive constants. The energy identities
below are justified in the stated classes by mollification in $y$.

\begin{lemma}\label{l:a0}
Let $(V_0,U_0,\Theta_0)-B\in H^3$ and $\inf_\R V_0,\inf_\R\Theta_0>0$.
There are $T_0>0$, depending only on
$\Vert(V_0,U_0,\Theta_0)-B\Vert_{H^3}$, the lower bounds of $V_0$ and
$\Theta_0$, $B$ and the fixed data, and a solution with these initial
values such that
\[
 (V_b,U_b,\Theta_b)-B\in C([0,T_0];H^3),\qquad
 \partial_y(U_b,\Theta_b)\in L^2(0,T_0;H^3),
\]
with $(V_b)_t\in L^2(0,T_0;H^3)$ and
$((U_b)_t,(\Theta_b)_t)\in L^2(0,T_0;H^2)$; it is unique in this class.
If the maximal existence time $T_*$ is finite, then
\[
 \sup_{0\le t<T_*}\bigl(
 \Vert(V_b,U_b,\Theta_b)(t)-B\Vert_{H^3}
 +\Vert V_b(t)^{-1}\Vert_{L^\infty}
 +\Vert\Theta_b(t)^{-1}\Vert_{L^\infty}\bigr)=\infty .
\]
\end{lemma}
\begin{proof}
Put $(V^0,U^0,\Theta^0)=(V_0,U_0,\Theta_0)$ and, for $n\ge0$, let
\[
\begin{aligned}
 U^{n+1}_t&=\nu\partial_y\frac{U^{n+1}_y}{V^n}
                         -\partial_y\frac{R\Theta^n}{V^n},\\
 c_v\Theta^{n+1}_t&=\kappa\partial_y\frac{\Theta^{n+1}_y}{V^n}
       -\frac{R\Theta^n}{V^n}U^{n+1}_y
       +\frac{\nu(U^{n+1}_y)^2}{V^n},\\
 V^{n+1}(t)&=V_0+\int_0^tU^{n+1}_y(r)\dd r,
\end{aligned}
\]
with the initial values $(U_0,\Theta_0)$; the two linear parabolic
problems for $U^{n+1}-B_u$ and $\Theta^{n+1}-B_T$ are solved by Fourier
projections and passage to the projection limit. Suppose that
$(V^n,\Theta^n)$ lies in a fixed positive box and
$\Vert(V^n,U^n,\Theta^n)-B\Vert_{H^3}\le L$ on $[0,T]$. Multiplying
$\partial_y^j$ of the equation for $U^{n+1}$, which does not contain
$\Theta^{n+1}$, by $\partial_y^j(U^{n+1}-B_u)$ for $0\le j\le3$, with the
pressure flux paired with $\partial_y^4U^{n+1}$ for $j=3$, and then
doing the same for the equation for $\Theta^{n+1}$, where the algebra
property of $H^2(\R)$ and the bound for $U^{n+1}$ give
\[
 \Bigl\Vert\frac{R\Theta^n}{V^n}U^{n+1}_y\Bigr\Vert_{H^2}
 +\Bigl\Vert\frac{\nu(U^{n+1}_y)^2}{V^n}\Bigr\Vert_{H^2}
 \le C_{L,B},
\]
we obtain
\[
\begin{aligned}
 &\sup_{0\le t\le T}\Vert(U^{n+1}-B_u,\Theta^{n+1}-B_T)(t)\Vert_{H^3}^2
       +c_L\int_0^T\Vert\partial_y(U^{n+1},\Theta^{n+1})\Vert_{H^3}^2\dd t\\
 &\qquad\le C\Vert(U_0-B_u,\Theta_0-B_T)\Vert_{H^3}^2+TC_{L,B},\\
 &\sup_{0\le t\le T}\Vert V^{n+1}(t)-V_0\Vert_{H^3}^2
       \le T\int_0^T\Vert U^{n+1}_y\Vert_{H^3}^2\dd t,\qquad
 \Vert(V^{n+1},\Theta^{n+1})(t)-(V_0,\Theta_0)\Vert_{L^\infty}
 \le C_{L,B}\sqrt T .
\end{aligned}
\]
Choosing $L$ from the initial norm and then $T$ small, depending on $L$,
$B$ and the lower bounds, the box and the bound $L$ are preserved. For
the differences $\Delta_n$ of consecutive iterates, the same
multipliers through order two, with the pressure difference treated as
a flux in $H^2$ and the heating difference as a source in $H^1$, give
\[
 \sup_{[0,T]}\Vert\Delta_{n+1}(V,U,\Theta)\Vert_{H^2}^2
 +c_L\int_0^T\Vert\partial_y\Delta_{n+1}(U,\Theta)\Vert_{H^2}^2\dd t
 \le TC_{L,B}\sup_{[0,T]}\Vert\Delta_n(V,U,\Theta)\Vert_{H^2}^2 .
\]
For $TC_{L,B}<\frac12$ the iterates converge in $C([0,T];H^2)$ to a
solution; the bounds of order three pass to the limit, the equations
give the time derivatives, the energy identity gives the continuity of
$U_b,\Theta_b$ in $H^3$, and $V_b=V_0+\int_0^t(U_b)_y$ gives that of
$V_b$. The difference estimate gives uniqueness. Since $T_0$ depends
only on the quantities in the continuation criterion, the solution
extends beyond $T_*$ if they stay bounded.
\end{proof}

\begin{lemma}\label{l:au}
Let $(V_i,U_i,\Theta_i)$, $i=1,2$, be distributional solutions on
$[0,T_e]$ with
\[
 (V_i,U_i,\Theta_i)-B\in C([0,T_e];H^1),\qquad
 \partial_y(U_i,\Theta_i)\in L^2(0,T_e;H^1),\qquad
 0<c_{T_e}\le V_i,\Theta_i\le C_{T_e},
\]
and put
$E_\Delta(t)=\Vert(V_1-V_2,U_1-U_2,\Theta_1-\Theta_2)(t)\Vert_{L^2}^2$.
Then, for $0\le t\le T_e$,
\[
 E_\Delta(t)\le E_\Delta(0)\exp\int_0^tk(r)\dd r,\qquad
 k=C\bigl(1+\Vert(U_1)_y\Vert_{L^\infty}^2
                 +\Vert(U_2)_y\Vert_{L^\infty}^2
                 +\Vert(\Theta_2)_y\Vert_{L^\infty}^2\bigr)\in L^1(0,T_e),
\]
with $C$ depending only on $c_{T_e}$, $C_{T_e}$ and the fixed data. In
particular, equal initial values give equal solutions on $[0,T_e]$.
\end{lemma}
\begin{proof}
Put $\Delta V=V_1-V_2$, $\Delta U=U_1-U_2$ and
$\Delta\Theta=\Theta_1-\Theta_2$. Subtract the equations and test by
$\Delta V$, $\Delta U$ and $\Delta\Theta/c_v$. Using
\[
 |V_1^{-1}-V_2^{-1}|\le C|\Delta V|,\qquad
 |P_1-P_2|\le C(|\Delta V|+|\Delta\Theta|),
 \qquad P_i=R\Theta_i/V_i,
\]
Cauchy's inequality and absorption of the derivative terms give
\[
 E_\Delta'
 +\nu\int_\R\frac{(\Delta U)_y^2}{V_1}\dd y
 +\frac\kappa{c_v}\int_\R\frac{(\Delta\Theta)_y^2}{V_1}\dd y
 \le k(t)E_\Delta .
\]
The constant in $k$ has precisely the dependence stated above.
By \eqref{e:agm} and the $L^2(0,T_e;H^1)$ bounds for the gradients,
$k\in L^1(0,T_e)$. Gronwall's inequality proves the assertion.
\end{proof}

\subsection{Finite-time positivity}\label{app:bg-positive}

We record the positivity argument needed in the application of
\cite[Theorem~1.1 and Proposition~3.1]{HuangWang2016} in
Lemma~\ref{l:bg-known}. Work on a finite bootstrap interval $[0,T]$
of the cited proof. Constants below may depend on the complete initial
$H^1$ bound and the positive initial lower bounds, but are independent
of $k$ and $T$.

For a fixed unit interval $I_k=[k,k+1]$, the entropy and Jensen
estimates in \cite[Lemmas~3.2--3.4]{HuangWang2016} give
\[
 c\le\int_{I_k}V_b(y,t)\dd y\le C,\qquad
 c\le\int_{I_k}\Theta_b(y,t)\dd y\le C.
\]
The same estimates give $c\le B_k(y,t)\le C$ for the factors in the
localized representation~\cite{Jiang1999}
\[
 V_b(y,t)=B_k(y,t)Y_k(t)
 +\frac R\nu\int_0^t
 \frac{B_k(y,t)Y_k(t)}{B_k(y,r)Y_k(r)}\Theta_b(y,r)\dd r,
 \qquad y\in I_k,
\]
where $Y_k(0)=1$ and $B_k(y,0)=V_b(y,0)$. Put
\[
 F_k(t)=\int_{I_k}\frac{V_b}{B_k}(y,t)\dd y,\qquad
 a_k(t)=\frac R\nu\int_{I_k}\frac{\Theta_b}{B_k}(y,t)\dd y.
\]
Then $c\le F_k,a_k\le C$ and
\[
 \frac{F_k(t)}{Y_k(t)}=1+\int_0^t\frac{a_k(r)}{Y_k(r)}\dd r,
 \qquad
 \frac{Y_k(t)}{Y_k(s)}=\frac{F_k(t)}{F_k(s)}
 \exp\left(-\int_s^t\frac{a_k(r)}{F_k(r)}\dd r\right).
\]
The second identity follows by differentiating the first. It yields
$c e^{-C(t-s)}\le Y_k(t)/Y_k(s)\le C e^{-c(t-s)}$ for
$0\le s\le t\le T$. The logarithmic-temperature Jensen estimate
\cite[(3.30)]{HuangWang2016} gives, whenever $1\le t\le T$,
\[
 \int_{t-1}^t\inf_{y\in I_k}\Theta_b(y,r)\dd r\ge c.
\]
Keeping only $r\in[t-1,t]$ in the representation proves $V_b(y,t)\ge c$
for $t\ge1$; for $t\le1$, its positive term $B_kY_k$ gives the same
bound. The upper bound is the finite-interval estimate
\cite[(3.38)]{HuangWang2016}.

The subsequent finite-interval estimates
\cite[Lemma~3.5]{HuangWang2016} give
\[
 \sup_{0\le t\le T}\Vert\vartheta_m(t)\Vert^2
 +\int_0^T\Vert\partial_y\vartheta_m(t)\Vert^2\dd t\le C.
\]
Completing the square in \eqref{e:bgth} gives
\[
 (\Theta_b)_t-\frac\kappa{c_v}\partial_y
       \left(\frac{(\Theta_b)_y}{V_b}\right)
 =\frac{\nu(U_b)_y^2-R\Theta_b(U_b)_y}{c_vV_b}
 \ge-C\Theta_b^2.
\]
Hence the comparison principle yields, for $0\le s\le t\le T$,
\[
 \inf_y\Theta_b(y,t)\ge
 \left[(\inf_y\Theta_b(y,s))^{-1}+C(t-s)\right]^{-1}.
\]
Let $a=\frac12\min\{\theta_-,\theta_+\}$. If
$\inf_y\Theta_b(y,t)\le a$, then $\Vert\vartheta_m(t)\Vert_\infty\ge a$.
The inequality
$\Vert\vartheta_m\Vert_\infty^2\le
2\Vert\vartheta_m\Vert\Vert\partial_y\vartheta_m\Vert$
therefore gives $\Vert\partial_y\vartheta_m(t)\Vert^2\ge c>0$.
Consequently,
\[
 \left|\{t\in[0,T]:\inf_y\Theta_b(y,t)\le a\}\right|\le C.
\]
Choose a fixed $L>C+1$. Every interval $[t-L,t]\subset[0,T]$
contains a time $s$ with $\inf_y\Theta_b(y,s)>a$. Propagating its
lower bound to $t$, and using the initial lower bound for $t\le L$, gives
\[
 \inf_{\R\times[0,T]}\Theta_b\ge
 \left[\max\{(\inf_y\Theta_b(y,0))^{-1},a^{-1}\}+CL\right]^{-1}>0.
\]
Thus the positivity constants are independent of $T$. The bootstrap
bound can first be fixed from the full initial $H^1$ data, and the wave
strength can then be chosen small enough for the cited absorptions.
Continuation gives the global solution before asymptotic limits are used.

\subsection{Higher regularity and Eulerian coordinates}\label{app:regularity}

\begin{proof}[Proof of Lemma~\ref{l:bg}]
Put $B=(V_c,U_c,\Theta_c)(\cdot,0)$.
The finite-time positive state bounds of Appendix~\ref{app:bg-positive},
the positive contact states and \eqref{e:blow} give
$c_b\le V_b,\Theta_b\le C_b$ on $\R\times[0,\infty)$.

To propagate the higher regularity, we prove the following estimate.
Suppose a solution in \eqref{e:bgcl} on $[0,T]$ satisfies
\[
 \underline k\le V_b,\Theta_b\le\overline k\ \text{ on }\R\times[0,T],
 \qquad K(T)\le K_1,\qquad
 \Vert\varphi_m(0)\Vert_{H^3}+\Vert(\upsilon_m,\vartheta_m)(0)\Vert_{H^2}
 \le2C_M.
\]
Then, with $C$ depending only on $\underline k$, $\overline k$, $K_1$,
$M$ and the fixed data, for $0<t_0<T$,
\begin{equation}\label{e:bgch}
\begin{aligned}
 &\sup_{0\le t\le T}\bigl(\Vert Z_m(t)\Vert_{H^2}^2
        +\Vert\partial_y^3\varphi_m(t)\Vert_{L^2}^2\bigr)
 +\int_0^T\bigl(\Vert\partial_y\varphi_m\Vert_{H^2}^2
        +\Vert\partial_y(\upsilon_m,\vartheta_m)\Vert_{H^2}^2\bigr)\dd t\\
 &\quad+\sup_{t_0\le t\le T}\Vert\partial_y^3(\upsilon_m,\vartheta_m)(t)\Vert_{L^2}^2
 +\int_{t_0}^T\Vert\partial_y^4(\upsilon_m,\vartheta_m)\Vert_{L^2}^2\dd t
 \le C\Bigl(1+\frac1{t_0}\Bigr).
\end{aligned}
\end{equation}

We first prove \eqref{e:bgch} for smooth solutions. Put
$Q_d=\nu\partial_y\log\frac{V_b}{V_c}-\upsilon_m$. By \eqref{eq:bg-m}, \eqref{eq:bg-mom} and
\eqref{eq:pr-id},
$\nu \partial_t\partial_y\log\frac{V_b}{V_c}=(U_b)_t+(P_b)_y-\nu\partial_y\frac{(U_c)_y}{V_c}$; inserting
$(P_b)_y=\frac{R(\Theta_b)_y}{V_b}-P_b\partial_y\log V_b$,
$\partial_y\log V_b=\frac{Q_d+\upsilon_m}\nu+\frac{(V_c)_y}{V_c}$ and
$R(\Theta_c)_y=p_*(V_c)_y$, we obtain
\begin{equation}\label{e:qdamp}
 (Q_d)_t+\frac{P_b}{\nu}Q_d
 =\frac R{V_b}(\vartheta_m)_y-\frac{P_b}{\nu}\upsilon_m
       -\frac{R\vartheta_m(V_c)_y}{V_bV_c}
       +(U_c)_t-\nu\partial_y\frac{(U_c)_y}{V_c}.
\end{equation}
Differentiating \eqref{e:qdamp} in $y$, using
$\frac{(P_b)_y}{\nu}(Q_d+\upsilon_m)=(P_b)_y\partial_y\log\frac{V_b}{V_c}$, and multiplying by
$(Q_d)_y$, we have
\begin{equation}\label{e:qd1}
\begin{gathered}
 \frac12\partial_t(Q_d)_y^2+\frac{P_b}{\nu}(Q_d)_y^2=Q_1(Q_d)_y,\\
 Q_1=\partial_y\Bigl[\frac R{V_b}(\vartheta_m)_y
       -\frac{R\vartheta_m(V_c)_y}{V_bV_c}
       +(U_c)_t-\nu\partial_y\frac{(U_c)_y}{V_c}\Bigr]
 -\frac{P_b}{\nu}(\upsilon_m)_y-(P_b)_y\partial_y\log\frac{V_b}{V_c} .
\end{gathered}
\end{equation}
The identities
\[
\begin{aligned}
 \partial_y\log\frac{V_b}{V_c}&=\frac{(\varphi_m)_y}{V_b}-\frac{\varphi_m(V_c)_y}{V_bV_c},\\
 (P_b)_y&=\frac{R(\vartheta_m)_y-P_b(\varphi_m)_y}{V_b}
       -\frac{(R\vartheta_m-p_*\varphi_m)(V_c)_y}{V_b^2},
\end{aligned}
\]
$K(T)\le K_1$ and \eqref{eq:pr-gau} give, for
$D(t)=\Vert\partial_yZ_m(t)\Vert_{L^2}^2+\delta^2\tau^{-1}$,
\[
 \Vert(P_b)_y\Vert_{L^2}^2+\Vert \partial_y\log\frac{V_b}{V_c}\Vert_{L^2}^2\le CD,\qquad
 \sup_{0\le t\le T}D\le C,\qquad
 \int_0^TD^3\dd t\le C .
\]
The estimate \eqref{e:agm} and $\nu \partial_y^2\log\frac{V_b}{V_c}=(Q_d)_y+(\upsilon_m)_y$ yield,
for every $\epsilon\in(0,1]$,
\[
\begin{aligned}
 \Bigl\Vert(P_b)_y\partial_y\log\frac{V_b}{V_c}\Bigr\Vert_{L^2}^2
 &\le2\Vert(P_b)_y\Vert_{L^2}^2
   \Bigl\Vert\partial_y\log\frac{V_b}{V_c}\Bigr\Vert_{L^2}
   \Bigl\Vert\partial_y^2\log\frac{V_b}{V_c}\Bigr\Vert_{L^2}\\
 &\le\epsilon\Vert(Q_d)_y\Vert_{L^2}^2+C\Vert(\upsilon_m)_y\Vert_{L^2}^2
       +C_\epsilon D^3,
\end{aligned}
\]
and \eqref{eq:pr-gau}, \eqref{e:put} and $(V_b)_y=(\varphi_m)_y+(V_c)_y$
yield
\[
 \Vert Q_1+(P_b)_y\partial_y\log\frac{V_b}{V_c}\Vert_{L^2}^2
 \le C\bigl(\Vert(\vartheta_m)_{yy}\Vert_{L^2}^2
       +\Vert(\vartheta_m)_y\Vert_{L^\infty}^2
       +\Vert\partial_yZ_m\Vert_{L^2}^2+\delta^2\tau^{-3/2}\bigr).
\]
Integrating \eqref{e:qd1} over $\R\times(0,t)$, using the Cauchy
inequality, $P_b\ge\frac{R\underline k}{\overline k}$,
$\Vert(Q_d)_y(0)\Vert_{L^2}\le C$, $K(T)\le K_1$ and \eqref{e:agm}, and
choosing $\epsilon$ small, we obtain
\[
 \sup_{0\le t\le T}\Vert(Q_d)_y\Vert_{L^2}^2
 +\int_0^T\Vert(Q_d)_y\Vert_{L^2}^2\dd t\le C.
\]
Differentiating $\partial_y\log V_b=(Q_d+\upsilon_m)/\nu+(V_c)_y/V_c$ gives
\[
\begin{aligned}
 (\varphi_m)_{yy}
 ={}&\frac{V_b}{\nu}\bigl((Q_d)_y+(\upsilon_m)_y\bigr)
       +V_b\Bigl(\partial_y\log\frac{V_b}{V_c}\Bigr)^2\\
 &+2V_b\frac{(V_c)_y}{V_c}\partial_y\log\frac{V_b}{V_c}
       +\frac{\varphi_m(V_c)_{yy}}{V_c}.
\end{aligned}
\]
By \eqref{e:agm},
\[
\begin{aligned}
 \Bigl\Vert\Bigl(\partial_y\log\frac{V_b}{V_c}\Bigr)^2\Bigr\Vert_{L^2}^2
 &\le2\Bigl\Vert\partial_y\log\frac{V_b}{V_c}\Bigr\Vert_{L^2}^3
       \Bigl\Vert\partial_y^2\log\frac{V_b}{V_c}\Bigr\Vert_{L^2}\\
 &\le\Bigl\Vert\partial_y^2\log\frac{V_b}{V_c}\Bigr\Vert_{L^2}^2+CD^3,\\
 \Bigl\Vert\frac{(V_c)_y}{V_c}\partial_y\log\frac{V_b}{V_c}\Bigr\Vert_{L^2}^2
 &\le C\delta^2\tau^{-1}D,\qquad
 \Vert\varphi_m(V_c)_{yy}\Vert_{L^2}^2\le C\delta^2\tau^{-3/2}.
\end{aligned}
\]
Together with the preceding bound for $(Q_d)_y$, this yields
\[
 \sup_{0\le t\le T}\Vert(\varphi_m)_{yy}\Vert_{L^2}^2
 +\int_0^T\Vert(\varphi_m)_{yy}\Vert_{L^2}^2\dd t\le C .
\]

By \eqref{eq:bg-mom}, \eqref{e:bgth} and \eqref{eq:pr-id},
\begin{equation}\label{e:bgpar}
\begin{aligned}
 (\upsilon_m)_t
 &=\partial_y\Bigl[\frac{\nu((\upsilon_m)_y+(U_c)_y)}{V_b}-P_b\Bigr]
   -(U_c)_t,\\
 c_v(\vartheta_m)_t
 &=\partial_y\Bigl[\frac{\kappa(\vartheta_m)_y}{V_b}
       -\frac{\kappa\varphi_m(\Theta_c)_y}{V_bV_c}\Bigr]
   -P_b(\upsilon_m)_y-(P_b-p_*)(U_c)_y+\frac{\nu(U_b)_y^2}{V_b}.
\end{aligned}
\end{equation}
Differentiating the first equation in $y$ and multiplying by
$-\partial_y^3\upsilon_m$, we have
\begin{equation}\label{e:ep2}
\begin{gathered}
 \frac12\frac{\dd}{\dd t}\Vert(\upsilon_m)_{yy}\Vert_{L^2}^2
 +\int_\R\frac\nu{V_b}|\partial_y^3\upsilon_m|^2\dd y
 =\int_\R Q_2\,\partial_y^3\upsilon_m\dd y,\\
 Q_2=-2\nu(\upsilon_m)_{yy}\partial_y\frac1{V_b}
     -\nu(\upsilon_m)_y\partial_y^2\frac1{V_b}
     -\nu\partial_y^2\frac{(U_c)_y}{V_b}+(P_b)_{yy}+(U_c)_{yt}.
\end{gathered}
\end{equation}
The bound $\sup_{0\le t\le T}\Vert\varphi_m\Vert_{H^2}\le C$,
\eqref{e:agm}, \eqref{eq:pr-gau}, \eqref{e:put} and
$P_b-p_*=\frac{R\vartheta_m-p_*\varphi_m}{V_b}$ yield
\[
 \Vert Q_2\Vert_{L^2}
 \le C\bigl(\Vert\partial_yZ_m\Vert_{H^1}+\delta\tau^{-3/4}\bigr).
\]
Integrating \eqref{e:ep2} over $(0,t)$ and using the Cauchy inequality,
we obtain
\[
 \sup_{0\le t\le T}\Vert(\upsilon_m)_{yy}\Vert_{L^2}^2
 +\int_0^T\Vert\partial_y^3\upsilon_m\Vert_{L^2}^2\dd t\le C .
\]
Similarly, multiplying the $y$-derivative of the second equation of
\eqref{e:bgpar} by $-\partial_y^3\vartheta_m$, and using
$\Vert(U_b)_y\Vert_{L^\infty}\le C$, which follows from the preceding
bound and \eqref{e:agm}, in
\[
 \Bigl\Vert\partial_y\frac{\nu(U_b)_y^2}{V_b}\Bigr\Vert_{L^2}
 \le C\bigl(\Vert(\upsilon_m)_{yy}\Vert_{L^2}+\Vert(\upsilon_m)_y\Vert_{L^2}
       +\delta\tau^{-3/4}\bigr),
\]
we obtain
\[
 \sup_{0\le t\le T}\Vert(\vartheta_m)_{yy}\Vert_{L^2}^2
 +\int_0^T\Vert\partial_y^3\vartheta_m\Vert_{L^2}^2\dd t\le C .
\]

Differentiating \eqref{e:qdamp} twice in $y$ gives
\begin{equation}\label{e:bd3}
\begin{aligned}
 (Q_d)_{yyt}+\frac{P_b}{\nu}(Q_d)_{yy}
  ={}&\partial_y^2\Bigl[
      \frac R{V_b}(\vartheta_m)_y-\frac{P_b}{\nu}\upsilon_m
      -\frac{R\vartheta_m(V_c)_y}{V_bV_c}
      +(U_c)_t-\nu\partial_y\frac{(U_c)_y}{V_c}\Bigr]\\
     &-\frac{2(P_b)_y}{\nu}(Q_d)_y
      -\frac{(P_b)_{yy}}{\nu}Q_d .
\end{aligned}
\end{equation}
By the bounds already proved, \eqref{e:agm} and \eqref{eq:pr-gau}, the
right side of \eqref{e:bd3} has $L^2$ norm at most
$C(\Vert\partial_y^3\vartheta_m\Vert_{L^2}+\Vert\partial_yZ_m\Vert_{H^1}
+\Vert(Q_d)_y\Vert_{L^2}+\delta\tau^{-3/4})$. Multiplying \eqref{e:bd3}
by $(Q_d)_{yy}$ and using $P_b\ge\frac{R\underline k}{\overline k}$ and
$\Vert(Q_d)_{yy}(0)\Vert_{L^2}\le C$, we obtain
\[
 \sup_{0\le t\le T}\Vert(Q_d)_{yy}\Vert_{L^2}^2
 +\int_0^T\Vert(Q_d)_{yy}\Vert_{L^2}^2\dd t\le C .
\]
Differentiating
$\partial_y\log V_b=\frac{Q_d+\upsilon_m}{\nu}+\partial_y\log V_c$ twice
gives
\begin{equation}\label{e:V3}
 (V_b)_{yyy}
 =V_b\Bigl\{\frac{(Q_d)_{yy}+(\upsilon_m)_{yy}}\nu
                         +\partial_y^3\log V_c\Bigr\}
       +\frac{3(V_b)_y(V_b)_{yy}}{V_b}
       -\frac{2(V_b)_y^3}{V_b^2}.
\end{equation}
Subtracting the same identity for $V_c$, in which $Q_d+\upsilon_m$ is
replaced by $0$, and using
$|\partial_y^3\log V_c|+|(V_c)_y(V_c)_{yy}|+|(V_c)_y|^3\le
C\delta\tau^{-3/2}$ and $\sup_t\Vert(\varphi_m)_y\Vert_{L^\infty}\le C$,
we obtain
\[
 \Bigl|\partial_y^3\varphi_m
   -\frac{V_b}{\nu}\bigl((Q_d)_{yy}+(\upsilon_m)_{yy}\bigr)\Bigr|
 \le C\bigl(|(\varphi_m)_{yy}|+|(\varphi_m)_y|+\delta\tau^{-3/2}|\varphi_m|\bigr),
\]
and hence
\[
 \sup_{0\le t\le T}\Vert\partial_y^3\varphi_m\Vert_{L^2}^2
 +\int_0^T\Vert\partial_y^3\varphi_m\Vert_{L^2}^2\dd t\le C .
\]

For $0<t_0<T$ choose $r\in(\frac{t_0}2,t_0)$ with
\[
 \Vert\partial_y^3(\upsilon_m,\vartheta_m)(r)\Vert_{L^2}^2
 \le\frac2{t_0}\int_{t_0/2}^{t_0}
       \Vert\partial_y^3(\upsilon_m,\vartheta_m)\Vert_{L^2}^2\dd t
 \le\frac C{t_0}.
\]
Similarly to \eqref{e:ep2}, differentiating \eqref{e:bgpar} twice in
$y$, multiplying by $-\partial_y^4\upsilon_m$ and
$-\partial_y^4\vartheta_m$, and using the bounds just proved, together
with
$\Vert(U_b)_{yy}^2\Vert_{L^2}\le\Vert(U_b)_{yy}\Vert_{L^\infty}
\Vert(U_b)_{yy}\Vert_{L^2}$ and \eqref{e:agm} for the second derivative
of the heating term, gives
\[
 \frac{\dd}{\dd t}\bigl(\Vert\partial_y^3\upsilon_m\Vert_{L^2}^2
       +c_v\Vert\partial_y^3\vartheta_m\Vert_{L^2}^2\bigr)
 +C^{-1}\Vert\partial_y^4(\upsilon_m,\vartheta_m)\Vert_{L^2}^2
 \le C\bigl(\Vert\partial_yZ_m\Vert_{H^2}^2+\delta^2\tau^{-3/2}\bigr),
\]
whose right side has integral at most $C$ over $(0,T)$. Integration from
$r$ gives the last two terms of \eqref{e:bgch}, and the preceding
displays give the first two.

Let $\eta_n$ be a standard mollifier in $y$. Since
$\Vert\eta_n*B-B\Vert_{H^k}\le n^{-1}\Vert\partial_yB\Vert_{H^k}$, for
$n\ge n_0$ the data $\eta_n*(V_b,U_b,\Theta_b)(\cdot,0)$ differ from $B$
by functions in every $H^k$, satisfy \eqref{e:bgin} with $2C_M$, and
their $Z_m(\cdot,0)$ converge in $H^3\times H^2\times H^2$.
The same theorem of Huang and Wang~\cite{HuangWang2016}, applied as in
Lemma~\ref{l:bg-known}, gives solutions $W^n$ with these data; the
preceding positivity argument and \eqref{e:blow} give bounds on
$[0,T_e]$ independent of $n$. By Lemmas~\ref{l:a0},
\ref{l:au}, \eqref{e:bgch} and the continuation criterion of
Lemma~\ref{l:a0}, $W^n-B\in C([0,T_e];H^3)$. By Lemma~\ref{l:au}, with
$\int_0^{T_e}k\dd t$ bounded independently of $n$ by \eqref{e:blow} and
\eqref{e:agm}, $W^n\to(V_b,U_b,\Theta_b)$ in $C([0,T_e];L^2)$. Hence the
solution satisfies \eqref{e:bgch} with $T=T_e$, and the equations, the
energy identity, \eqref{e:bd3} and \eqref{e:V3} give \eqref{e:bgcl}.

The solution therefore satisfies the hypotheses of \eqref{e:bgch} on
every $[0,T]$ with $\underline k=c_b$, $\overline k=C_b$ and
$K_1=K(\infty)$, and \eqref{e:bgch}, whose constant does not depend on
$T$, gives the bounds of the lemma. For $t\ge t_0$,
$(\varphi_m)_t=(\upsilon_m)_y$ by \eqref{eq:bg-m} and \eqref{eq:pr-id},
and \eqref{e:bgpar}, \eqref{eq:pr-gau} and \eqref{e:put} give, for
$j=1,2$,
\[
 \Vert\partial_y^j(\upsilon_m,\vartheta_m)_t\Vert_{L^2}
 \le C_b\bigl(\Vert\partial_yZ_m\Vert_{H^2}
     +\Vert\partial_y^4(\upsilon_m,\vartheta_m)\Vert_{L^2}+\delta\tau^{-1}\bigr).
\]
With $\frac{\dd}{\dd t}\Vert\partial_y^3w\Vert_{L^2}^2
=-2\langle\partial_y^4w,\partial_y^2w_t\rangle$ for
$w=(\upsilon_m,\vartheta_m)$, the bounds of the lemma give
$\int_{t_0}^\infty(\Vert\partial_y^jZ_m\Vert_{L^2}^2
+|\frac{\dd}{\dd t}\Vert\partial_y^jZ_m\Vert_{L^2}^2|)\dd t<\infty$ for
$1\le j\le3$, which together with \eqref{e:blow} gives the limit of the
lemma.

By \eqref{eq:bg-m}, $X_t=\int_0^y(U_b)_y\dd z+U_b(0,t)=U_b$, and
$z=Y_0(x)$ gives $X(y,0)=\int_0^y\frac{\dd z}{\rho_0(X_0(z))}=X_0(y)$.
Since $X_y=V_b$ lies between positive constants, $X(\cdot,t)$ is a
global diffeomorphism of $\R$. The chain rule \eqref{e:chain} transforms
\eqref{eq:bg-m}, \eqref{eq:bg-mom} and \eqref{e:bgth} into
\eqref{eq:bg-eqs}, and \eqref{eq:bg-eul} has initial value $W_0$.
Since $X-\bar X=\int_0^y\varphi_m\dd z+\int_0^t\upsilon_m(0,r)\dd r$ and
$\partial_x\bar Y=\bar\rho$, $|\bar Y(X,t)-y|\le C_{b,T_e}(1+|y|^{1/2})$
on $[0,T_e]$; with
$(\rho,u,\theta)-\bar Q=(-\frac{\varphi_m}{V_bV_c},\upsilon_m,\vartheta_m)\circ Y
+(\frac1{V_c},U_c,\Theta_c)\circ Y-(\frac1{V_c},U_c,\Theta_c)\circ\bar Y$,
this bound, \eqref{e:bgcl}, \eqref{e:chain}, the mean value theorem and
\eqref{eq:pr-gau} give \eqref{e:brg}, with the time derivatives following
from \eqref{eq:bg-eqs}.
Uniqueness in the mass class follows from Lemma~\ref{l:au}. An Eulerian solution
$(\rho_c,u_c,\theta_c)$ of \eqref{eq:bg-eqs} with initial value $W_0$
solves \eqref{eq:bg-m}--\eqref{eq:bg-en} in the labels
$Y_c=\int_0^{x}\rho_c\dd z-\int_0^t(\rho_cu_c)(0,r)\dd r$, with
$Y_c(\cdot,0)=Y_0$, hence coincides with $(V_b,U_b,\Theta_b)$, and
$Y_c=Y$.

\end{proof}

\subsection{Auxiliary balance identities}\label{app:balances}

\begin{proof}[Proof of Lemma~\ref{l:osp}]
Mollification reduces the proof to smooth functions; the stated
integrability permits passage to the limit, including the endpoints.
Integrating
$\partial_tE_i+\partial_xF_i=S_i$ over $(x,\infty)$ and over
$(-\infty,x)$ gives \eqref{e:osp1}. Hence, by $\partial_xP_2=-E_2$ and
an integration by parts,
\[
 \frac{\dd I}{\dd t}
 =\int_\R\Bigl[(S_1-\partial_xF_1)P_2
   +E_1\Bigl(F_2+\int_x^\infty S_2\dd z\Bigr)\Bigr]\dd x
 =\int_\R\Bigl(E_1F_2-F_1E_2+S_1P_2+E_1\int_x^\infty S_2\dd z\Bigr)\dd x,
\]
and integration over $(t_1,t_2)$ gives \eqref{e:osp2}.
\end{proof}

\begin{proof}[Proof of Lemma~\ref{l:trc}]
Since $(E_h)_z=\frac{h}{\Lambda}$, $|h(z)|\le|z|$ and $0\le h'\le1$,
\[
\begin{gathered}
 |(E_h)_z|\le\frac{|h(z)|}{\Lambda_0},\qquad
 |(E_h)_{zz}|\le\frac{h'(z)}{\Lambda_0}+C|h(z)|,\\
 |(E_h)_{zw}|\le C|h(z)|,\qquad
 |(E_h)_w|+|(E_h)_{ww}|\le C\Bigl|\int_0^z|h(a)|\dd a\Bigr|,
\end{gathered}
\]
which gives $|S_h|\le C\mathcal M$. The function $E_h(u,w)$ belongs to
$C^1([s,T];L^1)$ with derivative $(E_h)_zu_t+(E_h)_ww_t$. Inserting the
two equations and using
\[
\begin{gathered}
 (E_h)_zF_y=\partial_y[(E_h)_zF]-[(E_h)_{zz}u_y+(E_h)_{zw}w_y]F,\\
 (E_h)_wf_y=\partial_y[(E_h)_wf]-[(E_h)_{zw}u_y+(E_h)_{ww}w_y]f,\\
 \alpha(E_h)_zu_y=\bigl(\alpha-\alpha_0\Lambda(u,w)\bigr)(E_h)_zu_y
   +\alpha_0\partial_y\Phi_h(u),
\end{gathered}
\]
we obtain the first identity in \eqref{e:trc1}; since $\Phi_h(u)$ and
$J_h$ lie in $W^{1,1}(\R)$ for each time, its integral over
$\R\times(t',t)$ gives the second.

Let $\mathcal M\in L^1(\R\times(s,\infty))$ and
$\Vert u(t)\Vert_{L^\infty}\to0$. For $h=g_\varepsilon$ and $t$ so large
that $\Vert u(t)\Vert_{L^\infty}\le\varepsilon$, $E_{g_\varepsilon}(u,w)$
equals $\int_0^u\frac{a\dd a}{\Lambda(a,w)}\ge\frac{u^2}{2\max_K\Lambda}$.
The second identity in \eqref{e:trc1} with $t'=s$ gives
\[
 \limsup_{t\to\infty}\frac{\Vert u(t)\Vert_{L^2}^2}{2\max_K\Lambda}
 \le\int_\R E_{g_\varepsilon}(u,w)(s)\dd y
   +\int_s^\infty\!\!\int_\R|S_{g_\varepsilon}|\dd y\dd r .
\]
As $\varepsilon\to0$, $|E_{g_\varepsilon}(u,w)|\le\frac{2\varepsilon|u|}{\Lambda_0}$
tends to zero with the majorant $\frac{u^2}{2\Lambda_0}$, and
$|S_{g_\varepsilon}|\le C\mathcal M$. Every term of $S_{g_\varepsilon}$
contains a factor $g_\varepsilon(u)$ or
$\int_0^u|g_\varepsilon(a)|\dd a$, bounded by $2\varepsilon$ and
$2\varepsilon|u|$, except
$-\frac{\chi_e(u/\varepsilon)}{\Lambda(u,w)}u_yF$; this term equals
$-\frac{u_yF}{\Lambda(u,w)}$ on $\{u=0\}$, where $u_y=0$ almost
everywhere, and tends to zero elsewhere. Dominated convergence shows
that the right side tends to zero.
\end{proof}

\subsection{A threshold controlled by the initial entropy}\label{app:bg-low}

We control the initial time interval by relative entropy, then apply
the theorem of Huang and Wang~\cite{HuangWang2016} at a fixed positive time.

\begin{proposition}[A threshold controlled by the initial entropy]
\label{l:bg-entropy}
Let
$(\varphi_0,\upsilon_0,\vartheta_0)
=(V_0-V_c(\cdot,0),U_0-U_c(\cdot,0),\Theta_0-\Theta_c(\cdot,0))$
belong to $H^1(\R)$, and put
\[
 E_0=\int_\R\left\{\frac{\upsilon_0^2}{2}
 +R\Theta_c(y,0)\Phi\left(\frac{V_0}{V_c(y,0)}\right)
 +c_v\Theta_c(y,0)\Phi\left(\frac{\Theta_0}{\Theta_c(y,0)}\right)\right\}\dd y,
 \qquad \Phi(z)=z-\log z-1.
\]
Suppose, for $A\ge1$, that
\[
 \Vert\varphi_0\Vert_{H^1}+E_0\le A,
 \qquad \inf_\R V_0,\ \inf_\R\Theta_0\ge A^{-1}.
\]
There is $\delta_{\rm ent}(A)>0$, depending only on $A$ and the fixed
data, and independent of $u_*$, such that $0<\delta\le\delta_{\rm ent}(A)$
gives a unique global solution with
$Z_m\in C([0,\infty);H^1)$ and with $V_b,\Theta_b$ bounded above
and below by positive constants. It satisfies \eqref{e:blow}.
The threshold requires no bound on
$\Vert\upsilon_0\Vert_{L^4}$, $\Vert\vartheta_0\Vert$ or
$\Vert\partial_y(\upsilon_0,\vartheta_0)\Vert$.
The resulting global $H^1$ bound may depend on the complete initial $H^1$ norm.
\end{proposition}
\begin{proof}
We first prove the estimates for smooth data and then pass to an
$H^1$ limit. Subtract $u_*$ and restrict initially
$\delta\le\min\{1,\theta_-/2\}$. On a local interval $[0,T]$, $T\le1$,
let $G(t)$ be the nonnegative entropy-error integral of
\cite[(3.6)]{PengShiWu2026}, with viscosity $\nu$ and profile residuals
$(U_c)_t-\nu\partial_y((U_c)_y/V_c)$ and $-\nu(U_c)_y^2/V_c$.
Bootstrap $G(t)\le2(1+E_0)$. The relative entropy identity bounds the
entropy and its dissipation by $C_A$. We spell out the dependence in the volume argument of
Appendix~\ref{app:bg-positive}. The bounds for $B_k$ use only the initial
volume box and the velocity $L^2$ bound. The entropy cell averages and
the identity for $F_k/Y_k$ then give $Y_k(t)\ge c_A$ for $t\le1$;
hence $V_b\ge c_A$. The temperature comparison gives
$\Theta_b\ge(A+C_At)^{-1}\ge c_A$. If
\[
 d(t)=\int\frac{|(\vartheta_m)_y|^2}{V_b\Theta_b^2}\dd y,
\]
the cell averages and the oscillation estimate for $\sqrt{\Theta_b}$
give
\[
 \Vert\Theta_b(t)\Vert_\infty
 \le C_A+C_A\Vert V_b(t)\Vert_\infty\bigl(d(t)+\delta^2\bigr).
\]
Here the reference derivative contributes $C_A\delta^2$ by the lower
bounds already obtained. The volume representation and Gronwall,
using $\int_0^T d(t)\dd t\le C_A$, give the upper volume bound.
Thus, before using any initial velocity $L^4$ or temperature $L^2$ bound,
we have
\[
 c_A\le V_b\le C_A,\qquad \Theta_b\ge c_A,
 \qquad \sup_{t\le T}\Vert\upsilon_m(t)\Vert^2
 +\int_0^T\!\int\left(
 \frac{|(\upsilon_m)_y|^2}{V_b\Theta_b}
 +\frac{|(\vartheta_m)_y|^2}{V_b\Theta_b^2}\right)\dd y\dd t\le C_A.
\]
Here and below $C_A$ depends only on the stated lower-order data;
all constants in the cutoff below are fixed from $A$ and the fixed data.
In particular, volume perturbations are uniformly bounded in $L^2$
and $L^\infty$. Choose a fixed number $M$ larger than twice every
temperature occurring in the reference profile. The entropy bound gives
\[
 b=U_b-u_*\in L^\infty(0,T;L^2),\qquad
 h=(\Theta_b-M)_+,\qquad
 \sup_{t\le T}\bigl(\Vert b(t)\Vert^2+\Vert h(t)\Vert_{L^1}\bigr)\le C_A.
\]
The equations for these actual fields are
\[
 b_t=\nu\partial_y\left(\frac{b_y}{V_b}\right)-(P_b)_y,
 \qquad
 c_v(\Theta_b)_t=\kappa\partial_y\left(\frac{(\Theta_b)_y}{V_b}\right)
       -P_bb_y+\frac{\nu b_y^2}{V_b}.
\]

We first show that the temperature and velocity acquire the missing
integrability at positive times. For a sufficiently large fixed $L$ put
\[
 F=\int\left(\frac{c_v}{2}h^2+b^2h+\frac L4b^4\right)\dd y,
 \qquad
 D=\int\bigl(h_y^2+h b_y^2+b^2b_y^2\bigr)\dd y.
\]
Testing the temperature equation by $h$ and the momentum equation by
$2bh$, and substituting the temperature equation in
$\int b^2h_t$, gives, with $\chi=\boldsymbol1_{\{\Theta_b>M\}}$,
\[
\begin{aligned}
 &\frac{\dd}{\dd t}\int\left(\frac{c_v}{2}h^2+b^2h\right)
 +\int\frac{\kappa h_y^2+\nu h b_y^2}{V_b}\\
 &\quad\le\int P_bb_yh+2\int P_bb h_y
 -2\left(\nu+\frac\kappa{c_v}\right)\int\frac{bb_yh_y}{V_b}
 -\frac1{c_v}\int\chi P_bb^2b_y
 +\frac\nu{c_v}\int\frac{\chi b^2b_y^2}{V_b}.
\end{aligned}
\]
The omitted term is nonpositive: regularizing $\chi$ produces
$-\frac\kappa{c_v}\int\chi'(\Theta_b)b^2(\Theta_b)_y^2/V_b$.
This also justifies the calculation at the truncation level.
The additional identity is
\[
 \frac14\frac{\dd}{\dd t}\int b^4
 +3\nu\int\frac{b^2b_y^2}{V_b}=3\int P_bb^2b_y.
\]
Since $P_b\le C_A(M+h)$, the pressure terms are bounded by small
multiples of the displayed dissipations and
\[
 C_A\int\bigl((M+h)^2h+(M+h)^2b^2\bigr)
 \le C_A(1+\Vert h\Vert_\infty^2).
\]
The mixed gradient term is bounded by a small multiple of
$\Vert h_y\Vert^2$ plus $C_A\int b^2b_y^2$.
Choosing $L$ first absorbs all such velocity terms. The one-dimensional
Nash inequalities, applied to $h$ and $b^2$, give
\[
 \Vert h\Vert_\infty^2\le C_A\Vert h_y\Vert^{4/3},\qquad
 F\le C_A\left(\Vert h_y\Vert^{2/3}
                  +\Vert(b^2)_y\Vert^{2/3}\right)\le C_AD^{1/3}.
\]
Consequently
\begin{equation}\label{e:bg-entropy-hot}
 F'+c_AD\le C_A,\qquad F'+c_AF^3\le C_A,
 \qquad F(t)\le C_A(1+t^{-1/2})\quad(0<t\le T).
\end{equation}
The scalar comparison in the last step is independent of the finite
initial value $F(0)$. Multiplication of the first inequality by
$t^\beta$, with the fixed choice $\beta=\frac34$, therefore yields
\begin{equation}\label{e:bg-entropy-weight}
 \int_0^T t^\beta D(t)\dd t\le C_{A,\beta}.
\end{equation}
Indeed, $t^\beta F(t)\to0$ at zero, and
$\int_0^1 t^{\beta-1}(1+t^{-1/2})\dd t<\infty$.

The density gradient is estimated without squaring its source.
On $\{\Theta_b\le M\}$ the entropy dissipation controls
$(\Theta_b)_y$ in $L^2$, up to the reference-profile derivative.
On its complement $(\Theta_b)_y=h_y$. Thus
\[
 \int_0^T t^\beta\Vert(\Theta_b)_y\Vert^2\dd t\le C_{A,\beta},
 \qquad \int_0^T\Vert(\Theta_b)_y\Vert\dd t\le C_{A,\beta}.
\]
Furthermore $\Vert\Theta_b\Vert_\infty\le M+C_AD^{1/3}$, and weighted
H\"older gives
\[
 \int_0^T\Vert\Theta_b\Vert_\infty\dd t
 \le C_A+C_A\left(\int_0^Tt^\beta D\dd t\right)^{1/3}
                  \left(\int_0^Tt^{-\beta/2}\dd t\right)^{2/3}
 \le C_{A,\beta}.
\]
For $Q=\nu(\log V_b)_y-b$, the exact equations give
\[
 Q_t+\frac{P_b}{\nu}Q
       =\frac R{V_b}(\Theta_b)_y-\frac{P_b}{\nu}b.
\]
Its damping coefficient is positive and $\Vert Q(0)\Vert\le C_A$.
The variation-of-constants formula and Minkowski's inequality imply
\[
 \sup_{t\le T}\Vert Q(t)\Vert
 \le C_A+C_A\int_0^T\left(
        \Vert(\Theta_b)_y\Vert+\Vert\Theta_b\Vert_\infty\Vert b\Vert\right)\dd t
 \le C_{A,\beta}.
\]
Hence $\sup_{t\le T}\Vert\varphi_m(t)\Vert_{H^1}\le C_A$.
This estimate uses $\beta<1$ and does not require an unweighted
integral of $\Vert(\Theta_b)_y\Vert^2$ independent of the initial data.

The entropy bootstrap closes directly. For $\Theta_b\ge c_A$ and
bounded positive $\Theta_c$, elementary one-variable inequalities give
\[
 \frac{\vartheta_m^2}{\Theta_b^2}
 +\frac{\vartheta_m^2}{\Theta_b}
 +\Phi\left(\frac{\Theta_c}{\Theta_b}\right)
 \le C_A\Phi\left(\frac{\Theta_b}{\Theta_c}\right).
\]
All the terms in $G$ are controlled by the entropy using these
inequalities, Cauchy's inequality and the volume box. For the only
term with an extra temperature factor, use
\[
 \int\Theta_b\varphi_m^2
 \le M\Vert\varphi_m\Vert^2+\Vert\varphi_m\Vert_\infty^2\Vert h\Vert_{L^1}
 \le C_A.
\]
The residual term with $\upsilon_m$ is an $L^2$ pairing; the factor
$|\vartheta_m|/\Theta_b$ in the other residual is uniformly bounded.
The profile derivatives are $O(\delta)$ on $[0,1]$, with bounded
Gaussian $L^1$ and $L^2$ norms. Thus $G(T)\le C_A\delta$.
Choosing $\delta\le\delta_a(A)$ strictly improves
$G(T)\le2(1+E_0)$.

We justify continuation and passage to $H^1$ data. Mollify the actual
initial fields, preserving their positive lower bounds, and write
$Z_m^n$ for the resulting smooth perturbations. They converge in
$H^1$, their complete initial $H^1$ norms are bounded by a finite
number $H$, and their low-order data satisfy the same assumptions
with $2A$ in place of $A$. The final threshold is chosen for this
enlarged bound. Constants $C_H$ may also depend on $A$ and the fixed
data. All bootstrap
and hot-energy bounds proved above are uniform in $n$.
The unweighted hot-energy inequality gives
$\sup F_n+\int_0^T D_n\le C_A(1+F_n(0))\le C_H$.
Suppress $n$. Testing the equation for $Q$ by $Q$ gives
\[
 \int_0^T\!\int\Theta_bQ^2\dd y\dd t
 \le C_A\left(\Vert Q(0)\Vert^2
  +\sup_{t\le T}\Vert Q(t)\Vert
      \int_0^T\Vert(\Theta_b)_y\Vert\dd t
  +\sup_{t\le T}\Vert b(t)\Vert^2
      \int_0^T\Vert\Theta_b\Vert_\infty\dd t\right)\le C_A.
\]
Hence $\int_0^T\!\int\Theta_b|\partial_y\varphi_m|^2\dd y\dd t\le C_A$.
Put
\[
 S_T=\sup_{\R\times[0,T]}\Theta_b,\qquad
 E_u(T)=\sup_{r\le T}\Vert\partial_y\upsilon_m(r)\Vert^2
            +\int_0^T\Vert\partial_y^2\upsilon_m\Vert^2\dd r,
\]
\[
 E_T(T)=\sup_{r\le T}\Vert\partial_y\vartheta_m(r)\Vert^2
            +\int_0^T\Vert\partial_y^2\vartheta_m\Vert^2\dd r.
\]
Testing the velocity and temperature perturbation equations by
$-\partial_y^2\upsilon_m$ and $-\partial_y^2\vartheta_m$ gives
\[
 E_u(T)\le C_H(1+S_T),\qquad
 E_T(T)\le C_H\bigl(1+S_T+E_u(T)^{3/2}\bigr).
\]
The pressure terms use the preceding weighted density-gradient bound
and the hot-energy bound; the coefficient terms use
$\sup\Vert\partial_y\varphi_m\Vert^2\le C_A$ and
$\Vert f\Vert_\infty^2\le2\Vert f\Vert\Vert\partial_yf\Vert$.
The viscous heating term satisfies
\[
 \int_0^T\Vert\partial_y\upsilon_m\Vert_{L^4}^4\dd t
 \le2\sup_{t\le T}\Vert\partial_y\upsilon_m(t)\Vert^2
   \left(\int_0^T\Vert\partial_y\upsilon_m\Vert^2\dd t\right)^{1/2}
   \left(\int_0^T\Vert\partial_y^2\upsilon_m\Vert^2\dd t\right)^{1/2}
 \le C_HE_u(T)^{3/2}.
\]
All profile terms have bounded integrals on $[0,1]$.
Since $S_T^2\le C_H(1+E_T(T)^{1/2})$, these estimates yield
$S_T^2\le C_H(1+S_T^{3/4})$ and hence, for $T\le1$,
\begin{equation}\label{e:bg-entropy-h1approx}
 \sup_{t\le T}\Vert Z_m^n(t)\Vert_{H^1}^2
 +\int_0^T\Vert\partial_y^2(U_b^n,\Theta_b^n)\Vert^2\dd t\le C_H,
 \qquad c_A\le V_b^n,\Theta_b^n\le C_H.
\end{equation}
No additional restriction on $\delta$ occurs here. For each fixed
$n$, the calculation proving \eqref{e:bgch} has a finite constant
depending on its smooth initial norm. In conjunction with
\eqref{e:bg-entropy-h1approx}, it prevents the breakdown in
Lemma~\ref{l:a0} before time $1$: if its maximal time $T_n<1$,
take $t_0=T_n/4$ in \eqref{e:bgch} to obtain a finite $H^3$ bound
up to $T_n$, a contradiction. This uses no uniform bound for
the approximating $H^3$ norms: the estimates used below are only
\eqref{e:bg-entropy-h1approx}.

Retain the positive diffusion terms in the difference estimate of
Lemma~\ref{l:au}. Its coefficient has an integral bounded by $C_H$,
by \eqref{e:bg-entropy-h1approx} and the one-dimensional Sobolev
inequality. The initial convergence consequently gives
\[
 Z_m^n\longrightarrow Z_m\quad\hbox{in }C([0,1];L^2),\qquad
 \partial_y(U_b^n,\Theta_b^n)\longrightarrow
       \partial_y(U_b,\Theta_b)\quad\hbox{in }L^2(0,1;L^2).
\]
Uniform $H^1$ bounds and
$\Vert f\Vert_\infty^2\le2\Vert f\Vert\Vert f_y\Vert$ also give
uniform convergence of the fields. Thus the reciprocal-volume
coefficients converge uniformly, all diffusion fluxes pass to the
limit, and
$((U_b^n)_y)^2/V_b^n\to(U_b)_y^2/V_b$
in $L^1(\R\times(0,1))$. This identifies the viscous heating term.

The equations give the strong initial $H^1$ trace as well. In fact,
\[
 \int_0^1\Vert(U_b^n)_y\Vert_{L^4}^4\dd t
 +\int_0^1\Vert(V_b^n)_y\,\partial_y(U_b^n,\Theta_b^n)\Vert^2\dd t
 \le C_H
\]
by \eqref{e:bg-entropy-h1approx} and the one-dimensional Sobolev
inequality. Hence $((U_b^n)_t,(\Theta_b^n)_t)$ are bounded in
$L^2(0,1;L^2)$; these bounds pass to the limit. Taking
$B=(V_c,U_c,\Theta_c)(\cdot,0)$, the interpolation embedding
$H^1(0,1;L^2)\cap L^2(0,1;H^2)\subset C([0,1];H^1)$ gives
$(U_b-B_u,\Theta_b-B_T)\in C([0,1];H^1)$.
Moreover $(V_b)_t=(U_b)_y\in L^2(0,1;H^1)$, so
$V_b-B_v\in C([0,1];H^1)$ directly. The $L^2$ convergence at zero
identifies these traces with the prescribed $H^1$ initial datum.
The limit is therefore in the uniqueness class of Lemma~\ref{l:au}.
All low-order bounds established above pass to it, so the good-time
restart of Huang and Wang~\cite{HuangWang2016} used below applies without any higher
initial regularity assumption.

On $[1/2,1]$, \eqref{e:bg-entropy-hot} and
\eqref{e:bg-entropy-weight} are unweighted bounds depending only on
$A$. Splitting into $\{\Theta_b\le2M\}$ and its complement gives
\[
 \int_{1/2}^1\left(\Vert b_y\Vert^2+\Vert(\Theta_b)_y\Vert^2\right)\dd t
 \le C_A.
\]
The entropy on $\{\Theta_b\le2M\}$ and the hot part of $F$ give
\[
 \Vert\vartheta_m(t)\Vert^2\le C_A(1+F(t)).
\]
Together with the uniform density $H^1$ bound and the bound for $F$,
this provides a time $s\in[1/2,1]$ at which the full perturbation
$H^1$ norm is at most $C_A$. The reference-profile difference between
times $s$ and zero is $O(\delta)$ in $H^1$.
Restarting \cite[Theorem~1.1 and Proposition~3.1]{HuangWang2016}
at $s$, with the finite-interval positivity verification in
Appendix~\ref{app:bg-positive}, therefore requires only
$\delta\le\delta_b(2A)$ after the enlargement used in the approximation.
The restarted and original reference
profiles have vanishing uniform difference and square-integrable
spatial-derivative differences. Taking
$\delta_{\rm ent}(A)=\min\{1,\theta_-/2,\delta_a(2A),\delta_b(2A)\}$
proves \eqref{e:blow} for the original profile.
\end{proof}

\section{Additional derivative estimates}\label{s:coef}

\paragraph{Use as a planar background.}
In the three-dimensional application, the one-dimensional notation
$x,\partial_x,D_t$ becomes $x_1,\partial_1,D_b$, and
$(\rho,u,\theta)$ becomes $(\rho_b,u_b,\theta_b)$; the velocity field is
$u_be_1$, where $e_1=(1,0,0)^\top$. All norms here remain one-dimensional.
The conserved masses, diffusion waves and characteristic budgets are
given in Lemmas~\ref{l:cmp}, \ref{l:pr}, \ref{l:asy} and~\ref{l:dep}.
The additional coefficient bounds are Lemma~\ref{l:mcoef} and
Corollary~\ref{l:coef}; the weighted pressure decomposition is
Lemma~\ref{l:prs}. Proposition~\ref{p:uniform} gives a common wave-strength
threshold depending only on the lower-order data and the fixed
parameters, including the prescribed Gaussian exponents. Finite
constants and entrance times may depend on the complete initial datum.

Throughout this appendix $(V_b,U_b,\Theta_b)$ is the solution of
Lemma~\ref{l:bg} with $0<\delta\le\delta_1(M)$, $m_\pm$ and $h_*$ are the
masses and the phase of \eqref{e:mass}, $t\ge S_0$ with $S_0$ from
Lemma~\ref{l:dep}, norms of functions of $x$ are over $\R_x$, and
$\delta_1(M)$ may be decreased and $S_0$ increased in each statement.

\begin{lemma}\label{l:mcoef}
Put $Z=(V_b,b_m,\Theta_b)$,
$(V_h,H_h,\Theta_h)=(V_c,H_c,\Theta_c)(y+h_*,t)$ and
\begin{equation}\label{e:cmp}
\begin{gathered}
 (v_a,b_a,E_a)=\widehat{\mathcal U}(y+h_*,t)
   +\sum_{\sigma=\pm1}\vartheta_\sigma r_\sigma,
 \qquad T_a=\frac1{c_v}\Bigl(E_a-\frac{b_a^2}2\Bigr),\qquad
 p_a=\frac{RT_a}{v_a},\\
 W=(W_v,W_b,W_T)=(V_b-v_a,\,b_m-b_a,\,\Theta_b-T_a).
\end{gathered}
\end{equation}
Then $Z$ and $(v_a,b_a,T_a)$ take values in
$\mathcal K=[\frac1{2\rho_-},\frac2{\rho_-}]\times[-1,1]
\times[\frac{\theta_-}2,2\theta_-]$. With
$\Lambda_a=\kappa\partial_y\frac{(T_a)_y}{v_a}+\frac{\nu(b_a)_y^2}{v_a}$,
the functions
\begin{equation}\label{e:cdft}
\begin{gathered}
 f_v=(v_a)_t-(b_a)_y,\qquad
 f_b=(b_a)_t+(p_a)_y-\nu\partial_y\frac{(b_a)_y}{v_a},\qquad
 f_T=c_v(T_a)_t+p_a(b_a)_y-\Lambda_a\\
 \text{satisfy}\quad\Vert\partial_y^j(f_v,f_b,f_T)\Vert_{L^2}
 \le C_b\tau^{-\frac54-\frac j2}\qquad(0\le j\le3).
\end{gathered}
\end{equation}
Moreover,
\begin{equation}\label{e:wval}
\begin{gathered}
 (W_v,W_b)=\bigl((R_U)_1,(R_U)_2\bigr),\qquad
 c_vW_T=(R_U)_3-\frac{b_m+b_a}2W_b,\\
 \sup_{t\ge S_0}\sqrt\tau\Vert W\Vert_{L^\infty}\le C_b,\qquad
 \lim_{t\to\infty}\bigl(\sqrt\tau\Vert W\Vert_{L^\infty}
       +\Vert(R_U)_y\Vert_{L^\infty}\bigr)=0,
\end{gathered}
\end{equation}
\begin{equation}\label{e:wtC}
 \lim_{S\to\infty}\Bigl\{\sup_{t\ge S}\tau\Vert W_y\Vert_{H^1}^2
 +\int_S^\infty\tau\Vert W_{yy}\Vert_{H^1}^2\dd t\Bigr\}=0,
\end{equation}
and the quantity in braces is finite for $S=S_0$. Put $P=P_b-p_a$,
$\Sigma=\log\frac{\Theta_b}{T_a}-\frac{\gamma-1}\gamma\log\frac{P_b}{p_a}$ and
$\Lambda_b=\kappa\partial_y\frac{(\Theta_b)_y}{V_b}+\frac{\nu(b_m)_y^2}{V_b}$.
Then
\begin{equation}\label{e:wdif}
\begin{gathered}
 P_t+\frac{\gamma P_b}{V_b}(W_b)_y=\Xi_P,\qquad (W_b)_t+P_y=\Xi_B,\qquad
 \Sigma_t=\Xi_\Sigma,\\
 \Xi_P=(\gamma-1)\Bigl(\frac{\Lambda_b}{V_b}-\frac{\Lambda_a}{v_a}\Bigr)
   -\gamma\Bigl(\frac{P_b}{V_b}-\frac{p_a}{v_a}\Bigr)(b_a)_y
   -\frac{(\gamma-1)f_T-p_af_v}{v_a},\\
 \Xi_B=\nu\partial_y\Bigl(\frac{(b_m)_y}{V_b}-\frac{(b_a)_y}{v_a}\Bigr)-f_b,
 \qquad
 \Xi_\Sigma=\frac1{c_p}\Bigl(\frac{\Lambda_b}{\Theta_b}-\frac{\Lambda_a}{T_a}\Bigr)
   -\frac{f_T}{c_pT_a}-\frac{(\gamma-1)f_v}{\gamma v_a},\\
 \Vert(\Xi_P,\Xi_B,\Xi_\Sigma)\Vert_{L^2}^2
 \le C_b\bigl(\Vert W_{yy}\Vert_{L^2}^2+\tau^{-1}\Vert W_y\Vert_{L^2}^2
   +\tau^{-5/2}\bigr),
\end{gathered}
\end{equation}
and
\begin{equation}\label{e:munw}
\begin{gathered}
 \lim_{t\to\infty}\Vert Z_y\Vert_{H^2}=0,\qquad
 \sup_{t\ge S_0}\sqrt\tau\bigl(\Vert b_m\Vert_{L^\infty}
       +\Vert Z_y\Vert_{L^\infty}\bigr)<\infty,\\
 \int_{S_0}^\infty\bigl(\Vert Z_{yy}\Vert_{H^1}^2
   +\Vert\partial_y^4(b_m,\Theta_b)\Vert_{L^2}^2+\Vert(b_m)_y\Vert_{L^2}^2
   +\Vert(P_b)_y\Vert_{H^2}^2+\Vert Z_y\Vert_{L^4}^4\bigr)\dd t<\infty.
\end{gathered}
\end{equation}
Finally, fix $0<c_0<\frac18\min\{c_-,c_+\}$, let $\sigma=\pm1$, and put
$\Omega_\sigma(t)=\{y:\sigma(y+h_*)\ge c_0\tau\}$,
$Z_{\sigma,*}=(v_\sigma,0,\theta_\sigma)$,
$Z_{\sigma,a}=Z_{\sigma,*}+\vartheta_\sigma(-1,s_\sigma,\frac{p_*}{c_v})
-(0,0,\frac{s_\sigma^2\vartheta_\sigma^2}{2c_v})$ and
\[
 \mathsf B=\begin{pmatrix}0&-1&0\\-\frac{P_b}{V_b}&0&\frac R{V_b}\\
   0&\frac{P_b}{c_v}&0\end{pmatrix},\quad
 \mathsf V=\operatorname{diag}\Bigl(0,\frac\nu{V_b},\frac\kappa{c_vV_b}\Bigr),\quad
 \mathsf Q=\Bigl(0,\ \nu(b_m)_y\partial_y\frac1{V_b},\
   \frac{\kappa(\Theta_b)_y}{c_v}\partial_y\frac1{V_b}
   +\frac{\nu(b_m)_y^2}{c_vV_b}\Bigr)^{\!\top}.
\]
For $\varsigma\in\R$ and smooth functions $a,f$ on $\mathcal K$, with $f_Z$
the row of partial derivatives of $f$, put
$L=[af_Z(\varsigma I-\mathsf B)](Z_{\sigma,*})$ and
\[
 F=L(Z-Z_{\sigma,*})+af_Z\mathsf VZ_y,\qquad
 H=\bigl(af_Z(\varsigma I-\mathsf B)-L\bigr)Z_y+af_Z\mathsf Q
   -(af_Z\mathsf V)_yZ_y,
\]
where $a,f_Z,\mathsf B,\mathsf V$ are evaluated at $Z$. Then, for $s\ge S_0$,
\begin{equation}\label{e:E10}
\begin{gathered}
 \sup_{t\ge s}\sqrt\tau\bigl(\Vert Z-Z_{\sigma,a}\Vert_{L^\infty(\Omega_\sigma(t))}
   +\Vert Z_y\Vert_{L^\infty(\Omega_\sigma(t))}\bigr)\le\omega_b(s),\qquad
 a(Z)(\partial_t+\varsigma\partial_y)f(Z)=\partial_yF+H,\\
 \sup_{t\ge S_0}\sqrt\tau\Vert F\Vert_{L^\infty(\Omega_\sigma(t))}\le C_b,\qquad
 \sup_{t\ge s}\tau\Vert H\Vert_{L^\infty(\Omega_\sigma(t))}\le\omega_b(s),
\end{gathered}
\end{equation}
with $C_b,\omega_b$ depending also on $a,f,\varsigma$; if $\varsigma=s_\sigma$,
then $\sup_{t\ge s}\sqrt\tau\Vert F\Vert_{L^\infty(\Omega_\sigma(t))}\le\omega_b(s)$.
\end{lemma}

\begin{proof}
By \eqref{eq:prof}, $V_c=\frac{R\Theta_c}{p_*}$ and
$|\Theta_c-\theta_-|\le\delta$, and Lemma~\ref{l:pr} gives
$|\partial_y^j\vartheta_\sigma|\le C_M\tau^{-\frac{j+1}2}
e^{-\frac{(y-s_\sigma t)^2}{8d_\sigma\tau}}$ for the masses of
\eqref{e:mass}. Since $R\Theta_h=p_*V_h$, $R=(\gamma-1)c_v$ and
$r_\sigma=(-1,s_\sigma,p_*)$,
\[
 p_a=p_*+\frac1{v_a}\Bigl(\gamma p_*\sum_\sigma\vartheta_\sigma
   -\frac{\gamma-1}2b_a^2\Bigr).
\]
With \eqref{eq:pr-gau} this gives, for $j\ge1$,
\begin{equation}\label{e:cmpb}
\begin{gathered}
 \Vert\partial_y^j(v_a,T_a)\Vert_{L^\infty}
 +\Vert\partial_y^{j-1}(b_a,p_a-p_*)\Vert_{L^\infty}\le C_M\tau^{-\frac j2},\\
 \Vert\partial_y^j(v_a,T_a)\Vert_{L^2}
 +\Vert\partial_y^{j-1}(b_a,p_a-p_*)\Vert_{L^2}\le C_M\tau^{\frac14-\frac j2},
\end{gathered}
\end{equation}
and, with $\lim_{t\to\infty}\Vert Z_m\Vert_{L^\infty}=0$ from
Lemma~\ref{l:bg}, the values in $\mathcal K$.

Put $\widehat{\mathcal U}_h=\widehat{\mathcal U}(y+h_*,t)$,
$\widehat{\mathcal U}_\sigma=(v_\sigma,0,c_v\theta_\sigma)$ and
$\mathcal U_a=(v_a,b_a,E_a)$, so that
$\mathcal U_a-\widehat{\mathcal U}_h=\sum_\sigma\vartheta_\sigma r_\sigma$.
For functions $v,b,E$ with $\theta=\frac1{c_v}(E-\frac{b^2}2)$ and
$p=\frac{R\theta}v$, put $\mathsf F(v,b,E)=(-b,p,pb)$ and let
$\mathsf M(v,b,E)$ be the matrix with
$\mathsf M(v,b,E)\partial_y(v,b,E)
=(0,\frac{\nu b_y}v,\frac{\kappa\theta_y+\nu bb_y}v)$. By \eqref{e:cdft},
$(f_v,f_b,f_T+b_af_b)=(\mathcal U_a)_t
+\partial_y\bigl(\mathsf F(\mathcal U_a)
-\mathsf M(\mathcal U_a)\partial_y\mathcal U_a\bigr)$. This,
\eqref{e:Uhat} at $y+h_*$, the equation for $\vartheta_\sigma$ in
Lemma~\ref{l:pr} and
\[
 D\mathsf F(\widehat{\mathcal U}_\sigma)
 =\begin{pmatrix}0&-1&0\\-\frac{p_*}{v_\sigma}&0&\frac{\gamma-1}{v_\sigma}\\
   0&p_*&0\end{pmatrix},\qquad
 D\mathsf F(\widehat{\mathcal U}_\sigma)r_\sigma
 =\Bigl(-s_\sigma,\frac{\gamma p_*}{v_\sigma},s_\sigma p_*\Bigr)=s_\sigma r_\sigma
\]
give
\[
\begin{gathered}
 (f_v,f_b,f_T+b_af_b)=\partial_y(0,\widehat R_2,\widehat R_3)(y+h_*,t)
 +\sum_\sigma r_\sigma\bigl(d_\sigma(\vartheta_\sigma)_{yy}
   -k_\sigma s_\sigma(\vartheta_\sigma^2)_y\bigr)+\partial_y\mathcal N,\\
 \begin{aligned}
 \mathcal N={}&\sum_\sigma\vartheta_\sigma\bigl(D\mathsf F(\widehat{\mathcal U}_h)
   -D\mathsf F(\widehat{\mathcal U}_\sigma)\bigr)r_\sigma
 +\mathsf F(\mathcal U_a)-\mathsf F(\widehat{\mathcal U}_h)
   -D\mathsf F(\widehat{\mathcal U}_h)(\mathcal U_a-\widehat{\mathcal U}_h)\\
 &-\mathsf M(\mathcal U_a)\partial_y(\mathcal U_a-\widehat{\mathcal U}_h)
 -\bigl(\mathsf M(\mathcal U_a)-\mathsf M(\widehat{\mathcal U}_h)\bigr)
   \partial_y\widehat{\mathcal U}_h .
 \end{aligned}
\end{gathered}
\]
In the first and the last term of $\mathcal N$, $\vartheta_\sigma$ multiplies
a bounded function that is at most $C\delta e^{-\frac{c_G(y+h_*)^2}\tau}$
where $\sigma(y+h_*)\ge0$, by \eqref{eq:ct-gau} and \eqref{eq:pr-gau}; since
$|s_\sigma t+h_*|\ge\frac{c_\sigma t}2$, their derivatives of order at most
four have $L^2$ norms at most $C_b\tau^{-3}$. The other two terms satisfy
$|\partial_y^j(\cdot)|\le C_M\tau^{-1-\frac j2}\sum_\sigma
e^{-\frac{(y-s_\sigma t)^2}{8d_\sigma\tau}}$, and, by \eqref{e:Uhat} and
\eqref{eq:pr-gau}, $|\partial_y^j\widehat R_2|+|\partial_y^j\widehat R_3|
\le C\delta\tau^{-1-\frac j2}e^{-\frac{c_Gy^2}\tau}$. With \eqref{e:cmpb} for
$b_af_b$, this proves the bound in \eqref{e:cdft}.

By \eqref{e:Um} and \eqref{e:cmp},
$W=Z_m-(v_a-V_c,\,b_a-H_c,\,T_a-\Theta_c)$. By \eqref{eq:pr-gau},
$\Vert(V_h-V_c,H_h-H_c,\Theta_h-\Theta_c)\Vert_{H^4}
\le|h_*|\Vert\partial_y(V_c,H_c,\Theta_c)\Vert_{H^4}\le C_b\tau^{-1/4}$
and $\Vert\partial_y(V_h-V_c,H_h-H_c,\Theta_h-\Theta_c)\Vert_{H^3}
\le|h_*|\Vert\partial_y^2(V_c,H_c,\Theta_c)\Vert_{H^3}\le C_b\tau^{-3/4}$.
With Lemmas~\ref{l:pr} and \ref{l:dep},
\begin{equation}\label{e:wbas}
 \lim_{t\to\infty}\Vert W\Vert_{H^3}=0,\qquad
 \int_{S_0}^\infty\Vert W_y\Vert_{H^2}^2\dd t<\infty .
\end{equation}
Since $Z=Z_m+(V_c,H_c,\Theta_c)$, Lemmas~\ref{l:bg} and \ref{l:dep},
\eqref{e:agm}, \eqref{e:bprod}, \eqref{eq:pr-gau} and
$P_b=\frac{R\Theta_b}{V_b}$ give the limit and the integral in
\eqref{e:munw}. By \eqref{e:CH} and \eqref{e:cmp},
$R_U=(V_b,b_m,E_m)-(v_a,b_a,E_a)$; with $c_v\Theta_b=E_m-\frac{b_m^2}2$ and
$c_vT_a=E_a-\frac{b_a^2}2$ this gives the identities in \eqref{e:wval}, and
Lemma~\ref{l:asy} gives the bounds for $W$. The limit of $(R_U)_y$ follows
from the identities and \eqref{e:wbas}.

By \eqref{e:bgth} and \eqref{e:pbt},
$(P_b)_t+\frac{\gamma P_b}{V_b}(b_m)_y=\frac{(\gamma-1)\Lambda_b}{V_b}$ and
$\partial_t(\log\Theta_b-\frac{\gamma-1}\gamma\log P_b)=\frac{\Lambda_b}{c_p\Theta_b}$. By
\eqref{e:cdft}, the same computation with $p_a=\frac{RT_a}{v_a}$ gives
\[
 (p_a)_t+\frac{\gamma p_a}{v_a}(b_a)_y
 =\frac{(\gamma-1)(\Lambda_a+f_T)-p_af_v}{v_a},\qquad
 \partial_t(\log T_a-\frac{\gamma-1}\gamma\log p_a)
 =\frac{\Lambda_a+f_T}{c_pT_a}+\frac{(\gamma-1)f_v}{\gamma v_a}.
\]
Subtracting these identities from the preceding ones, and the definition of
$f_b$ from \eqref{eq:bg-mom}, gives the equations in \eqref{e:wdif}. Put
$G=\frac{\gamma P_b}{V_b}$,
$F_h=\kappa\bigl(\frac{(\Theta_b)_y}{V_b}-\frac{(T_a)_y}{v_a}\bigr)
=\frac{\kappa(W_T)_y}{V_b}-\frac{\kappa W_v(T_a)_y}{V_bv_a}$ and
$J=\nu\bigl(\frac{(b_m)_y^2}{V_b}-\frac{(b_a)_y^2}{v_a}\bigr)$. Since
$\Lambda_b-\Lambda_a=\partial_yF_h+J$,
\begin{equation}\label{e:wsys}
\begin{gathered}
 \Xi_P=\frac{\gamma-1}{V_b}\partial_yF_h+\mathfrak R_P,\qquad
 \Xi_B=\nu\partial_y\frac{(W_b)_y}{V_b}+\mathfrak R_B,\qquad
 \Xi_\Sigma=\frac{\partial_yF_h}{c_p\Theta_b}+\mathfrak R_\Sigma,\\
 \mathfrak R_P=\frac{(\gamma-1)J}{V_b}-\frac{(\gamma-1)W_v\Lambda_a}{V_bv_a}
   -\Bigl(G-\frac{\gamma p_a}{v_a}\Bigr)(b_a)_y
   -\frac{(\gamma-1)f_T-p_af_v}{v_a},\\
 \mathfrak R_B=-\nu\partial_y\frac{W_v(b_a)_y}{V_bv_a}-f_b,\qquad
 \mathfrak R_\Sigma=\frac{J}{c_p\Theta_b}-\frac{W_T\Lambda_a}{c_p\Theta_bT_a}
   -\frac{f_T}{c_pT_a}-\frac{(\gamma-1)f_v}{\gamma v_a}.
\end{gathered}
\end{equation}

To prove \eqref{e:wtC}, set
\[
\begin{gathered}
 E_d=\frac12\sum_{j=1}^2\int\bigl(G^{j-1}|\partial_y^jP|^2
   +G^j|\partial_y^jW_b|^2+c_p\Theta_bG^j|\partial_y^j\Sigma|^2\bigr)\dd y,\qquad
 K_d=\int (W_b)_yP_{yy}\dd y,\\
 D_d=\Vert P_{yy}\Vert_{L^2}^2+\Vert((W_b)_{yy},(W_T)_{yy})\Vert_{H^1}^2,\qquad
 N_d=\tau^{-3/2}\Vert W\Vert_{L^\infty}^2+\Vert(f_v,f_b,f_T)\Vert_{H^2}^2 .
\end{gathered}
\]
In $\mathcal K$, $\Theta_b=T_ae^\Sigma(\frac{P_b}{p_a})^{\frac{\gamma-1}\gamma}$ and
$V_b=v_ae^\Sigma(\frac{P_b}{p_a})^{-\frac1\gamma}$. Differentiating these and the
definitions of $P,\Sigma$ twice, and using \eqref{e:cmpb}, \eqref{e:wval},
\eqref{e:wbas} and, at the end, a larger $S_0$, we obtain
\begin{equation}\label{e:wPS}
\begin{gathered}
 \frac1C\Vert W_y\Vert_{H^1}^2-C_b\tau^{-1/2}\Vert W\Vert_{L^\infty}^2
 \le E_d\le C\Vert W_y\Vert_{H^1}^2+C_b\tau^{-1/2}\Vert W\Vert_{L^\infty}^2,\\
 \Vert W_{yy}\Vert_{L^2}^2\le CD_d
   +C_b\bigl(\tau^{-1}+\Vert W_y\Vert_{L^\infty}^2\bigr)\Vert W_y\Vert_{H^1}^2
   +C_b\tau^{-3/2}\Vert W\Vert_{L^\infty}^2 .
\end{gathered}
\end{equation}
Apply $\partial_y^j$, $j=1,2$, to \eqref{e:wdif}, multiply by
$G^{j-1}\partial_y^jP$, $G^j\partial_y^jW_b$ and $c_p\Theta_bG^j\partial_y^j\Sigma$,
integrate and add. The terms $-\int(G^j)_y\partial_y^jP\partial_y^jW_b$ from the
hyperbolic part cancel $\int P_yG_y(W_b)_y$ at $j=1$ and $2\int GG_yP_{yy}(W_b)_{yy}$
at $j=2$, and $G^{j-1}\frac{\gamma-1}{V_b}=\frac{(\gamma-1)G^j}{\gamma P_b}$. With
$\mathfrak s_j=\frac{\gamma-1}{\gamma P_b}\partial_y^jP+\partial_y^j\Sigma
-\frac{\partial_y^jW_T}{\Theta_b}$ and
$\mathfrak h_j=\partial_y^jF_h-\frac{\kappa\partial_y^{j+1}W_T}{V_b}$, this
gives
\[
\begin{aligned}
 &E_d'+\sum_{j=1}^2\int G^j\Bigl(\frac{\nu|\partial_y^{j+1}W_b|^2}{V_b}
   +\frac{\kappa|\partial_y^{j+1}W_T|^2}{V_b\Theta_b}\Bigr)\dd y
 =-\int GG_{yy}P_{yy}(W_b)_y\dd y\\
 &\quad+\frac12\sum_{j=1}^2\int\bigl((G^{j-1})_t|\partial_y^jP|^2
   +(G^j)_t|\partial_y^jW_b|^2+(c_p\Theta_bG^j)_t|\partial_y^j\Sigma|^2\bigr)\dd y
 +\sum_{j=1}^2(I_{1,j}+I_{2,j}+I_{3,j}),
\end{aligned}
\]
where
\[
\begin{aligned}
 I_{1,j}={}&-\nu\int\Bigl(\frac{(G^j)_y}{V_b}\partial_y^jW_b\,\partial_y^{j+1}W_b
   +\partial_y(G^j\partial_y^jW_b)\Bigl[\partial_y^j,\frac1{V_b}\Bigr](W_b)_y\Bigr)\dd y,\\
 I_{2,j}={}&\int\Bigl(G^{j-1}\partial_y^jP\Bigl[\partial_y^j,\frac{\gamma-1}{V_b}\Bigr]
   +c_p\Theta_bG^j\partial_y^j\Sigma\Bigl[\partial_y^j,\frac1{c_p\Theta_b}\Bigr]\Bigr)
   \partial_yF_h\dd y\\
 &-\int\Bigl(\partial_y(G^j\mathfrak s_j)\,\partial_y^jF_h
   +\partial_y\Bigl(\frac{G^j\partial_y^jW_T}{\Theta_b}\Bigr)\mathfrak h_j
   +\frac\kappa{V_b}\Bigl(\frac{G^j}{\Theta_b}\Bigr)_y\partial_y^jW_T\,
     \partial_y^{j+1}W_T\Bigr)\dd y,\\
 I_{3,j}={}&\int\bigl(G^{j-1}\partial_y^jP\,\partial_y^j\mathfrak R_P
   +G^j\partial_y^jW_b\,\partial_y^j\mathfrak R_B
   +c_p\Theta_bG^j\partial_y^j\Sigma\,\partial_y^j\mathfrak R_\Sigma\bigr)\dd y .
\end{aligned}
\]
By the definitions of $\Sigma$ and $P$,
\[
\begin{gathered}
 \mathfrak s_2=-\frac{W_T(T_a)_{yy}}{\Theta_bT_a}+\frac{(\gamma-1)P(p_a)_{yy}}{\gamma P_bp_a}
   -\frac{(\Theta_b)_y^2}{\Theta_b^2}+\frac{(T_a)_y^2}{T_a^2}
   +\frac{\gamma-1}\gamma\Bigl(\frac{(P_b)_y^2}{P_b^2}-\frac{(p_a)_y^2}{p_a^2}\Bigr),\\
 \mathfrak h_j=\kappa\Bigl[\partial_y^j,\frac1{V_b}\Bigr](W_T)_y
   -\kappa\partial_y^j\frac{W_v(T_a)_y}{V_bv_a},
\end{gathered}
\]
so that $\partial_y(G^2\mathfrak s_2)$ contains no third derivative of $P$.
Since $G_{yy}-\frac{\gamma P_{yy}}{V_b}
=\gamma\bigl(\frac{(p_a)_{yy}}{V_b}+2(P_b)_y\partial_y\frac1{V_b}
+P_b\partial_y^2\frac1{V_b}\bigr)$,
\[
 \Bigl|\int GG_{yy}P_{yy}(W_b)_y\dd y\Bigr|
 \le C\Vert (W_b)_y\Vert_{L^\infty}\Vert P_{yy}\Vert_{L^2}
   \bigl(\Vert P_{yy}\Vert_{L^2}+\Vert W_{yy}\Vert_{L^2}+C_b\tau^{-3/4}
   +C_b\Vert W_y\Vert_{L^2}\bigr).
\]
By \eqref{e:bgth} and \eqref{e:pbt},
\[
 c_v(\Theta_b)_t=\Lambda_b-P_b(b_m)_y,\qquad
 G_t=-\frac{(\gamma+1)G}{V_b}(b_m)_y+\frac{\gamma(\gamma-1)}{V_b^2}\Lambda_b,
\]
with $(b_m)_y=(W_b)_y+(b_a)_y$ and
$\Lambda_b=\frac{\kappa(W_T)_{yy}}{V_b}+\frac{\kappa(T_a)_{yy}}{V_b}
-\frac{\kappa(V_b)_y(\Theta_b)_y}{V_b^2}+\frac{\nu(b_m)_y^2}{V_b}$. By
\eqref{e:cmpb}, each weight $\mathfrak m\in\{G^{j-1},G^j,c_p\Theta_bG^j\}$
satisfies $\mathfrak m_t=\mathfrak a(W_b)_y+\mathfrak b(W_T)_{yy}+\mathfrak r$
with $\mathfrak a,\mathfrak b$ smooth functions of $(V_b,\Theta_b)$ and
$|\mathfrak r|\le C_b\tau^{-1}+C|W_y|^2$. At $j=2$ the terms with
$\mathfrak a,\mathfrak b$ are at most
$C\Vert((W_b)_y,(W_T)_{yy})\Vert_{L^\infty}\Vert\partial_y^2(P,W_b,\Sigma)\Vert_{L^2}^2$;
at $j=1$, for $\mathfrak c\in\{\mathfrak a,\mathfrak b\}$ and $w\in\{P,W_b,\Sigma\}$,
\[
\begin{aligned}
 \int\mathfrak c(W_b)_yw_y^2\dd y
 &=-2\int\mathfrak cW_bw_yw_{yy}\dd y-\int\mathfrak c_yW_bw_y^2\dd y,\\
 \int\mathfrak c(W_T)_{yy}w_y^2\dd y
 &=-2\int\mathfrak c(W_T)_yw_yw_{yy}\dd y-\int\mathfrak c_y(W_T)_yw_y^2\dd y,
\end{aligned}
\]
and $\Vert W_b\Vert_{L^\infty}^2\le C_b\tau^{-1}$,
$\Vert(W_T)_y\Vert_{L^\infty}^2\le\Vert W_y\Vert_{H^1}^2$ by \eqref{e:wval} and
\eqref{e:agm}; $\Sigma_{yy}=\frac{(W_T)_{yy}}{\Theta_b}-\frac{(\gamma-1)P_{yy}}{\gamma P_b}
+\mathfrak s_2$. In $I_{3,1}$ the derivative on the defects is moved to the
other factor, and at $j=2$ the first term of $\mathfrak R_B$ is integrated by
parts once. Let $0<\varepsilon\le1$. By Cauchy's inequality, \eqref{e:agm},
\eqref{e:cmpb}, \eqref{e:wval}, \eqref{e:wPS} and
$\lim_{t\to\infty}\Vert W_y\Vert_{W^{1,\infty}}=0$ from \eqref{e:wbas}, after
increasing $S_0$,
\[
 E_d'+\frac1{C_1}\Vert((W_b)_{yy},(W_T)_{yy})\Vert_{H^1}^2
 \le\varepsilon\Vert P_{yy}\Vert_{L^2}^2
   +\frac{C_b}\varepsilon\Bigl(\bigl(\tau^{-1}+\Vert W_y\Vert_{H^2}^2\bigr)E_d
   +N_d\Bigr).
\]
By \eqref{e:wdif},
\[
 K_d'+\Vert P_{yy}\Vert_{L^2}^2
 =\int\bigl(G(W_b)_{yy}^2+G_y(W_b)_y(W_b)_{yy}+P_{yy}\partial_y\Xi_B
   -(W_b)_{yy}\partial_y\Xi_P\bigr)\dd y,
\]
and by \eqref{e:wsys} and the same bounds,
$K_d'+\frac12\Vert P_{yy}\Vert_{L^2}^2
\le C_2\Vert((W_b)_{yy},(W_T)_{yy})\Vert_{H^1}^2
+C_b((\tau^{-1}+\Vert W_y\Vert_{H^2}^2)E_d+N_d)$. Since $|K_d|\le CE_d$, we
choose $\eta>0$ with $2\eta C_1C_2\le1$ and
$\frac12E_d\le\mathcal E=E_d+\eta K_d\le2E_d$, and then
$\varepsilon=\frac\eta4$. Then, with $k=C_b\Vert W_y\Vert_{H^2}^2$,
\[
 \mathcal E'+\frac1CD_d\le\Bigl(\frac{C_b}\tau+k\Bigr)\mathcal E+C_bN_d,\qquad
 (\tau\mathcal E)'+\frac\tau CD_d\le C_b\mathcal E+k\tau\mathcal E+C_b\tau N_d .
\]
We use the following averaging estimate. Let $r_0\ge1$, let $E\ge0$ be
absolutely continuous on $[r_0,\infty)$, and let $k,g\ge0$ satisfy
$(\tau E)'\le k\tau E+g$. For $t\ge2r_0$ there is
$t_1\in[\frac t2,t]$ with $E(t_1)\le\frac2t\int_{t/2}^tE\dd r$; hence
$\tau_{t_1}E(t_1)\le4\int_{t/2}^tE\dd r$, and integration over $[t_1,t]$
gives
\begin{equation}\label{e:avst}
 \tau E(t)\le\exp\Bigl(\int_{r_0}^\infty k\dd r\Bigr)
 \int_{t/2}^\infty(4E+g)\dd r .
\end{equation}

By \eqref{e:wPS}, \eqref{e:wbas}, \eqref{e:wval} and \eqref{e:cdft},
$\int_{S_0}^\infty(k+\mathcal E+\tau N_d)\dd t<\infty$. Gronwall's inequality and
\eqref{e:avst}, with $g=C_b(\mathcal E+\tau N_d)$, give
\[
 \sup_{t\ge S_0}\tau\mathcal E+\int_{S_0}^\infty\tau D_d\dd t<\infty,
 \qquad \lim_{t\to\infty}\tau\mathcal E=0.
\]
Consequently,
\[
 \int_S^\infty\tau D_d\dd t\le C\tau_S\mathcal E(S)
 +C\int_S^\infty\bigl(C_b\mathcal E+k\tau\mathcal E+C_b\tau N_d\bigr)\dd t
 \longrightarrow0 .
\]
With \eqref{e:wPS}, \eqref{e:wval} and
\[
 \int_{S_0}^\infty\tau\Vert W_y\Vert_{L^\infty}^2\Vert W_y\Vert_{H^1}^2\dd t
 \le\sup_{t\ge S_0}\tau\Vert W_y\Vert_{H^1}^2
   \int_{S_0}^\infty\Vert W_y\Vert_{H^1}^2\dd t,
\]
this proves
\begin{equation}\label{e:wlow}
 \lim_{S\to\infty}\Bigl\{\sup_{t\ge S}\tau\Vert W_y\Vert_{H^1}^2
 +\int_S^\infty\tau\bigl(\Vert W_{yy}\Vert_{L^2}^2
   +\Vert((W_b)_{yyy},(W_T)_{yyy})\Vert_{L^2}^2\bigr)\dd t\Bigr\}=0,
\end{equation}
the quantity in braces being finite for $S=S_0$. Since $b_m=b_a+W_b$ and
$Z_y=W_y+\partial_y(v_a,b_a,T_a)$, \eqref{e:wlow}, \eqref{e:agm},
\eqref{e:wval} and \eqref{e:cmpb} give the second bound in \eqref{e:munw}.

Put $\Psi_v=\nu\partial_y\log\frac{V_b}{v_a}-W_b$. By \eqref{eq:bg-m},
$(v_a)_t=(b_a)_y+f_v$, the equation for $W_b$ in \eqref{e:wdif} and
$P_y=\frac{R(W_T)_y}{V_b}-\frac{RW_T(v_a)_y}{V_bv_a}-\frac{W_v(p_a)_y}{V_b}
-\frac{P_b}\nu(\Psi_v+W_b)$, which follows from $V_bP=RW_T-p_aW_v$,
\[
 (\Psi_v)_t+\frac{P_b}\nu\Psi_v
 =\frac R{V_b}(W_T)_y-\frac{P_b}\nu W_b-\frac{RW_T(v_a)_y}{V_bv_a}
   -\frac{W_v(p_a)_y}{V_b}+f_b-\nu\partial_y\frac{f_v}{v_a}.
\]
Applying $\partial_y^2$ gives
$(\Psi_v)_{yyt}+\frac{P_b}\nu(\Psi_v)_{yy}=N_v$ with
\[
 N_v=\partial_y^2\Bigl[\frac R{V_b}(W_T)_y-\frac{P_b}\nu W_b
   -\frac{RW_T(v_a)_y}{V_bv_a}-\frac{W_v(p_a)_y}{V_b}
   +f_b-\nu\partial_y\frac{f_v}{v_a}\Bigr]
   -\frac{2(P_b)_y}\nu(\Psi_v)_y-\frac{(P_b)_{yy}}\nu\Psi_v ,
\]
which contains no third derivative of $W_v$. Multiplying by
$(\Psi_v)_{yy}$ gives
$[\tau\Vert(\Psi_v)_{yy}\Vert_{L^2}^2]'
+\frac\tau C\Vert(\Psi_v)_{yy}\Vert_{L^2}^2\le C\tau\Vert N_v\Vert_{L^2}^2$.
By \eqref{e:wlow}, \eqref{e:wval}, \eqref{e:cmpb}, \eqref{e:munw} and
\eqref{e:cdft} with $j\le3$,
$\int_{S_0}^\infty\tau\Vert N_v\Vert_{L^2}^2\dd t<\infty$; for instance
$\tau\Vert(P_b)_y(\Psi_v)_y\Vert_{L^2}^2
\le\sup_{t\ge S_0}\tau\Vert(P_b)_y\Vert_{L^\infty}^2\Vert(\Psi_v)_y\Vert_{L^2}^2$
and $\tau\Vert(P_b)_{yy}\Psi_v\Vert_{L^2}^2
\le\sup_{t\ge S_0}\tau\Vert\Psi_v\Vert_{L^\infty}^2\Vert(P_b)_{yy}\Vert_{L^2}^2$.
Hence $\lim_{t\to\infty}\tau\Vert(\Psi_v)_{yy}\Vert_{L^2}^2=0$ and
$\int_S^\infty\tau\Vert(\Psi_v)_{yy}\Vert_{L^2}^2\dd t
\le C\tau_S\Vert(\Psi_v)_{yy}(S)\Vert_{L^2}^2
+C\int_S^\infty\tau\Vert N_v\Vert_{L^2}^2\dd t\to0$. Since
$\nu\partial_y^3\log\frac{V_b}{v_a}=(\Psi_v)_{yy}+(W_b)_{yy}$ and
\[
 \partial_y^3\log\frac{V_b}{v_a}=\frac{(W_v)_{yyy}}{V_b}
   -\frac{W_v(v_a)_{yyy}}{V_bv_a}
   -3\Bigl(\frac{(V_b)_y(V_b)_{yy}}{V_b^2}-\frac{(v_a)_y(v_a)_{yy}}{v_a^2}\Bigr)
   +2\Bigl(\frac{(V_b)_y^3}{V_b^3}-\frac{(v_a)_y^3}{v_a^3}\Bigr),
\]
\eqref{e:wlow}, \eqref{e:wval}, \eqref{e:cmpb} and \eqref{e:munw} give
$\int_S^\infty\tau\Vert(W_v)_{yyy}\Vert_{L^2}^2\dd t\to0$, which completes
\eqref{e:wtC}.

By \eqref{e:wsys}, \eqref{e:cmpb}, \eqref{e:wval} and \eqref{e:munw},
\[
 \Vert(\Xi_P,\Xi_B,\Xi_\Sigma)\Vert_{L^2}
 \le C\Vert W_{yy}\Vert_{L^2}+C_b\tau^{-1/2}\Vert W_y\Vert_{L^2}
   +C_b\tau^{-3/4}\Vert W\Vert_{L^\infty}+C\Vert(f_v,f_b,f_T)\Vert_{L^2},
\]
which with \eqref{e:wval} and \eqref{e:cdft} is the last line of
\eqref{e:wdif}.

For $t\ge S_0$ and $y\in\Omega_\sigma(t)$, $|y+h_*|\ge c_0\tau$ and
$|y-s_{-\sigma}t|\ge c_{-\sigma}t$. Since
$(v_a,b_a,T_a)-Z_{\sigma,a}=\bigl(V_h-v_\sigma-\vartheta_{-\sigma},\,
H_h+s_{-\sigma}\vartheta_{-\sigma},\,\Theta_h-\theta_\sigma
+\frac{p_*}{c_v}\vartheta_{-\sigma}
-\frac{(b_a-s_\sigma\vartheta_\sigma)(b_a+s_\sigma\vartheta_\sigma)}{2c_v}\bigr)$,
\eqref{eq:ct-gau}, \eqref{eq:pr-gau} and Lemma~\ref{l:pr} give, on
$\Omega_\sigma(t)$,
\begin{equation}\label{e:extb}
\begin{gathered}
 |\partial_y^j((v_a,b_a,T_a)-Z_{\sigma,a})|\le C_b\tau^{-3}\quad(0\le j\le3),\\
 |\partial_y^lZ_{\sigma,a}|\le C_M\tau^{-\frac{l+1}2}\quad(l=1,2),\qquad
 |Z_{\sigma,a}-Z_{\sigma,*}|\le C_M\tau^{-1/2}.
\end{gathered}
\end{equation}
With $Z-Z_{\sigma,a}=W+(v_a,b_a,T_a)-Z_{\sigma,a}$, \eqref{e:wval},
\eqref{e:wtC} and \eqref{e:agm}, this gives the first bound in
\eqref{e:E10}. By \eqref{eq:bg-m}, \eqref{eq:bg-mom} and \eqref{e:bgth},
$Z_t+\mathsf BZ_y=\mathsf VZ_{yy}+\mathsf Q$. Hence
$a(Z)(\partial_t+\varsigma\partial_y)f(Z)
=af_Z\bigl((\varsigma I-\mathsf B)Z_y+\mathsf VZ_{yy}+\mathsf Q\bigr)$, and
$LZ_y=\partial_y[L(Z-Z_{\sigma,*})]$ and
$af_Z\mathsf VZ_{yy}=\partial_y(af_Z\mathsf VZ_y)-(af_Z\mathsf V)_yZ_y$ give the
identity in \eqref{e:E10}. On $\Omega_\sigma(t)$,
$|af_Z(\varsigma I-\mathsf B)-L|\le C|Z-Z_{\sigma,*}|$,
$|\mathsf Q|+|(af_Z\mathsf V)_yZ_y|\le C|Z_y|^2$ and
$|Z-Z_{\sigma,*}|\le|Z-Z_{\sigma,a}|+C_M\tau^{-1/2}$, which give the bounds
for $F$ and $H$. If $\varsigma=s_\sigma$, then
$(s_\sigma I-\mathsf B(Z_{\sigma,*}))(-1,s_\sigma,\frac{p_*}{c_v})^\top=0$,
because $\frac{p_*}{v_\sigma}(1+\frac R{c_v})=\frac{\gamma p_*}{v_\sigma}
=s_\sigma^2$; so $|L(Z_{\sigma,a}-Z_{\sigma,*})|\le C_M\tau^{-1}$ and
$|F|\le C|Z-Z_{\sigma,a}|+C|Z_y|+C_M\tau^{-1}$ on $\Omega_\sigma(t)$.
\end{proof}

\begin{corollary}\label{l:coef}
Put
\begin{equation}\label{e:bc}
 h_b=\partial_xu,\qquad A_b=D_tu,\qquad B_b=D_t\theta,\qquad
 g_T=\frac{\partial_x\theta}{\theta}.
\end{equation}
Then
\begin{equation}\label{e:bal}
 D_t\rho=-\rho h_b,\qquad \rho A_b+\partial_xp=\nu\partial_xh_b,\qquad
 c_v\rho B_b=\kappa\partial_x^2\theta+\nu h_b^2-ph_b,
\end{equation}
and there is $C_0$, depending only on the fixed data, such that, for
$s\ge S_0$, with $V_h,\Theta_h,Z,\mathcal K$ as in Lemma~\ref{l:mcoef},
\begin{gather}
 \sup_{t\ge s}\Vert(\rho-\rho_-,u-u_*,\theta-\theta_-)\Vert_{L^\infty}
   \le C_0\delta+\omega_b(s),\label{e:E0}\\
 \sup_{t\ge s}\sqrt\tau\Vert u-u_*\Vert_{L^\infty}\le C_M+\omega_b(s),
   \label{e:cbg}
\end{gather}
\begin{equation}\label{e:E1}
\begin{gathered}
 \sup_{t\ge s}\tau\Vert\partial_x(\rho,\theta)\Vert_{L^\infty}^2
   \le C_0\delta^2+\omega_b(s),\qquad
 \partial_x(\rho,\theta)=\Bigl(-\frac{(V_h)_y}{V_b^3},
   \frac{(\Theta_h)_y}{V_b}\Bigr)\circ Y+k_1,\\
 \Bigl|\Bigl(-\frac{(V_h)_y}{V_b^3},\frac{(\Theta_h)_y}{V_b}\Bigr)\circ Y\Bigr|
   \le C_0\delta\tau^{-1/2},\qquad
 \sup_{t\ge s}\tau\Vert k_1\Vert_{L^\infty}^2
   +\int_s^\infty\bigl(\Vert k_1\Vert_{L^\infty}^2+\Vert k_1\Vert_{L^2}^2\bigr)\dd t
   \le\omega_b(s),
\end{gathered}
\end{equation}
\begin{gather}
 \sup_{t\ge s}\tau\Vert(h_b,\partial_xp)\Vert_{L^\infty}^2\le\omega_b(s),
   \label{e:E2}\\
 \int_s^\infty\Vert(h_b,\partial_xp,A_b,B_b,\partial_xh_b,\partial_x^2p,
   \partial_x^2\theta,\partial_xA_b,\partial_xB_b)\Vert_{L^\infty}^2\dd t
   \le\omega_b(s),\label{e:E3}
\end{gather}
\begin{equation}\label{e:curv}
\begin{gathered}
 \partial_x^2\rho=c_\rho+e_{\rho,2},\qquad
 \partial_x^2\theta=c_\theta+e_{\theta,2},\\
 c_\rho=\Bigl(-\frac{(V_h)_{yy}}{V_b^4}+\frac{3(V_h)_y^2}{V_b^5}\Bigr)\circ Y,
 \qquad
 c_\theta=\Bigl(\frac{(\Theta_h)_{yy}}{V_b^2}
   -\frac{(V_h)_y(\Theta_h)_y}{V_b^3}\Bigr)\circ Y,\\
 |c_\rho|+|c_\theta|\le\frac{C_0\delta}\tau,\qquad
 \int_s^\infty\bigl(\sqrt\tau\Vert(e_{\rho,2},e_{\theta,2})\Vert_{L^2}^2
   +\tau\Vert(e_{\rho,2},e_{\theta,2})\Vert_{L^\infty}^2\bigr)\dd t
   \le\omega_b(s),
\end{gathered}
\end{equation}
\begin{gather}
 \int_s^\infty\bigl(\Vert\partial_x^3\rho\Vert_{L^2}^2
   +\Vert(h_b,B_b)\Vert_{L^2}^2
   +\sqrt\tau\Vert(\partial_xh_b,\partial_x^2h_b,\partial_xA_b,\partial_xB_b)
     \Vert_{L^2}^2\bigr)\dd t\le\omega_b(s),\label{e:E5}\\
 \sup_{t\ge s}\bigl(\Vert\partial_x(\rho,u,\theta)\Vert_{H^2}
   +\Vert(A_b,B_b)\Vert_{H^1}\bigr)\le\omega_b(s),\label{e:E6}
\end{gather}
\begin{equation}\label{e:unwE}
\begin{gathered}
 \sup_{t\ge S_0}\sqrt\tau\bigl(\Vert u-u_*\Vert_{L^\infty}
   +\Vert\partial_x(\rho,u,\theta)\Vert_{L^\infty}\bigr)<\infty,\\
 \int_s^\infty\bigl(\Vert\partial_x^2(\rho,u,\theta)\Vert_{H^1}^2
   +\Vert\partial_x^4(u,\theta)\Vert_{L^2}^2+\Vert\partial_xu\Vert_{L^2}^2
   +\Vert\partial_xp\Vert_{H^1}^2\bigr)\dd t\le\omega_b(s),
\end{gathered}
\end{equation}
\begin{equation}\label{e:E8}
\begin{gathered}
 \sup_{t\ge S_0}\tau\Vert(u-u_*,\partial_x\rho,\partial_x\theta,
   \partial_xp,h_b)\Vert_{L^\infty}^2<\infty,\qquad
 \int_{S_0}^\infty\tau^a\Vert(A_b,B_b,h_b,\partial_xp)\Vert_{L^\infty}^2\dd t
   <\infty\quad(a<\tfrac23),\\
 \int_{S_0}^\infty\tau\Bigl(\Vert(\partial_xh_b,\partial_x^2p)\Vert_{L^\infty}^2
   +\Bigl\Vert\Bigl(\partial_x^2p,\partial_xA_b,\partial_xB_b,
     \partial_x\frac{\partial_x^2\theta}{\rho\theta}\Bigr)\Bigr\Vert_{L^2}^2
   \Bigr)\dd t<\infty,\\
 \int_{S_0}^\infty\tau^{1/4}\Bigl(\Bigl\Vert\Bigl(\partial_x^2\log\rho,
   \partial_x^2\log\theta,\frac{\partial_x^2\theta}{\rho\theta}\Bigr)
   \Bigr\Vert_{L^\infty}^2+\Vert\partial_x^3\log\rho\Vert_{L^2}^2\Bigr)\dd t
   <\infty,
\end{gathered}
\end{equation}
and, for smooth functions $f,g$ on $\mathcal K$, integers $n\ge0$,
$c_1>0$ and $\xi_\sigma=\frac{Y-s_\sigma t}{\sqrt\tau}$,
\begin{equation}\label{e:floc}
 \int_s^\infty\sqrt\tau\,\Bigl\Vert\bigl(g(Z)\circ Y\bigr)|\xi_\sigma|^n
   e^{-c_1\xi_\sigma^2}\bigl(|\partial_x(f(Z)\circ Y)|
   +|D_t(f(Z)\circ Y)|\bigr)\Bigr\Vert_{L^\infty}^2\dd t\le\omega_b(s),
\end{equation}
with $\omega_b$ depending also on $f,g,n,c_1$.
\end{corollary}

\begin{proof}
By \eqref{e:bc}, \eqref{eq:bg-eqs} is \eqref{e:bal}. For a function $F(y,t)$,
\eqref{e:chain} gives
\begin{equation}\label{e:cder}
\begin{gathered}
 [\partial_x(F\circ Y)]\circ X=\frac{F_y}{V_b},\qquad
 [D_t(F\circ Y)]\circ X=F_t,\\
 [\partial_x^2(F\circ Y)]\circ X
       =\frac{F_{yy}}{V_b^2}-\frac{(V_b)_yF_y}{V_b^3},\qquad
 \Vert F\circ Y\Vert_{L^2}^2=\int_\R V_bF^2\dd y.
\end{gathered}
\end{equation}
and, by induction, $[\partial_x^k(F\circ Y)]\circ X$ is a linear combination of
$V_b^{-k-l}\partial_y^iF\,\partial_y^{j_1}V_b\cdots\partial_y^{j_l}V_b$ with
$i,j_1,\dots,j_l\ge1$ and $i+j_1+\dots+j_l=k$. With
$(\rho,u,\theta)\circ X=(\frac1{V_b},u_*+b_m,\Theta_b)$, \eqref{eq:bg-mom} and
\eqref{e:bgth},
\[
 h_b\circ X=\frac{(b_m)_y}{V_b},\qquad
 A_b\circ X=\nu\partial_y\frac{(b_m)_y}{V_b}-(P_b)_y,\qquad
 c_vB_b\circ X=\Lambda_b-P_b(b_m)_y .
\]

Since $V_c=\frac{R\Theta_c}{p_*}$, $|\Theta_c-\theta_-|\le\delta$ and
$|H_c|\le C\delta$ by \eqref{eq:pr-gau}, Lemma~\ref{l:bg} gives \eqref{e:E0}.
By \eqref{e:CH}, $b_m=H_h+\sum_\sigma s_\sigma\vartheta_\sigma+(R_U)_2$, and
Lemma~\ref{l:pr}, \eqref{eq:pr-gau} and Lemma~\ref{l:asy} give \eqref{e:cbg}.

By \eqref{e:cder}, $k_1\circ X=\bigl(-\frac{(V_b-V_h)_y}{V_b^3},
\frac{(\Theta_b-\Theta_h)_y}{V_b}\bigr)$, where
$(V_b-V_h,\Theta_b-\Theta_h)=(W_v,W_T)
+\bigl(-\sum_\sigma\vartheta_\sigma,\frac{p_*}{c_v}\sum_\sigma\vartheta_\sigma
-\frac{b_a^2}{2c_v}\bigr)$. By \eqref{eq:pr-gau},
$|(V_h)_y|+|(\Theta_h)_y|\le C\delta\tau^{-1/2}$, and by \eqref{e:agm},
\eqref{e:cmpb} and Lemma~\ref{l:pr},
\[
 \sup_{t\ge s}\tau\Vert k_1\Vert_{L^\infty}^2
 +\int_s^\infty\bigl(\Vert k_1\Vert_{L^\infty}^2+\Vert k_1\Vert_{L^2}^2\bigr)\dd t
 \le C\sup_{t\ge s}\tau\Vert W_y\Vert_{H^1}^2
   +C\int_s^\infty\Vert W_y\Vert_{H^1}^2\dd t+C_M\tau_s^{-1/2},
\]
which tends to zero by \eqref{e:wtC} and \eqref{e:wbas}. This gives
\eqref{e:E1}. In the same way,
\[
\begin{gathered}
 e_{\rho,2}\circ X=-\frac{(V_b-V_h)_{yy}}{V_b^4}
   +\frac{3(V_b-V_h)_y(V_b+V_h)_y}{V_b^5},\\
 e_{\theta,2}\circ X=\frac{(\Theta_b-\Theta_h)_{yy}}{V_b^2}
   -\frac{(V_b-V_h)_y(\Theta_b)_y+(V_h)_y(\Theta_b-\Theta_h)_y}{V_b^3},
\end{gathered}
\]
$|c_\rho|+|c_\theta|\le\frac{C\delta}\tau$ by \eqref{eq:pr-gau}, and the
integral in \eqref{e:curv} follows from
$\tau\Vert(W_v,W_T)_{yy}\Vert_{L^\infty}^2\le\tau\Vert W_{yy}\Vert_{H^1}^2$,
\eqref{e:wtC}, \eqref{e:cmpb}, the bound for $k_1$ and
$\sup_{t\ge S_0}\tau\Vert\partial_y(V_b,V_h,\Theta_b)\Vert_{L^\infty}^2<\infty$
from \eqref{e:munw} and \eqref{eq:pr-gau}.

By \eqref{e:cder}, $\partial_xp\circ X=\frac{(P_b)_y}{V_b}$ with
$(P_b)_y=(p_a)_y+P_y$ and $V_bP=RW_T-p_aW_v$. Hence, by \eqref{e:cmpb},
\eqref{e:wval} and \eqref{e:munw},
$\tau\Vert(h_b,\partial_xp)\Vert_{L^\infty}^2
\le C_b\tau^{-1}+C\tau\Vert W_y\Vert_{L^\infty}^2$, and \eqref{e:agm} and
\eqref{e:wtC} give \eqref{e:E2}.

By \eqref{e:cder} and the formulas for $h_b,A_b,B_b$, each component of the
vector in \eqref{e:E3} and its derivative $\partial_x$, and
$\partial_x^k\rho$ ($k=2,3$), $\partial_x^ku$ ($1\le k\le4$),
$\partial_x^k\theta$ ($2\le k\le4$) and $\partial_x^kp$ ($k=1,2$), are,
composed with $X$, sums of monomials, each a smooth function of $Z$ times
derivatives of $Z$ and $P_b$, and every monomial contains one of
$\partial_y^i(V_b,P_b)$ ($2\le i\le3$), $\partial_y^i(b_m,\Theta_b)$
($2\le i\le4$), $(b_m)_y$, $(P_b)_y$, or two factors from $Z_y$. By
\eqref{e:munw}, whose limit bounds the remaining factors in $L^\infty$, the
squares of the $L^2$ norms of these monomials are integrable over
$(S_0,\infty)$. With $\Vert F\circ Y\Vert_{L^\infty}^2\le2\Vert F\Vert_{L^2}
\Vert F_y\Vert_{L^2}$ from \eqref{e:agm}, this proves \eqref{e:E3}, the
integral in \eqref{e:unwE} and the first two terms of \eqref{e:E5}. The first
limit in \eqref{e:munw} gives \eqref{e:E6}, and its second bound gives the
first line of \eqref{e:unwE}. In the last term of \eqref{e:E5}, replace
$(V_b,b_m,\Theta_b,P_b)$ by $(v_a,b_a,T_a,p_a)+(W_v,W_b,W_T,P)$: the terms
without $W$ have squared $L^2$ norms at most $C_b\tau^{-5/2}$ by
\eqref{e:cmpb}, and the others are bounded by \eqref{e:wtC}, \eqref{e:wval},
\eqref{e:munw} and \eqref{e:agm}.

For $w\in\{W_v,W_b,W_T\}$ and $\ell>0$, the mean value theorem and Cauchy's
inequality give $|w_y(y)|\le\frac{2\Vert w\Vert_{L^\infty}}\ell
+\ell^{1/2}\Vert w_{yy}\Vert_{L^2}$. Minimizing in $\ell$ and using
\eqref{e:wval} and H\"older's inequality,
\[
 \Vert W_y\Vert_{L^\infty}^2\le C\Vert W\Vert_{L^\infty}^{2/3}
   \Vert W_{yy}\Vert_{L^2}^{4/3},\qquad
 \int_{S_0}^\infty\tau^a\Vert W_y\Vert_{L^\infty}^2\dd t
 \le C_b\Bigl(\int_{S_0}^\infty\tau\Vert W_{yy}\Vert_{L^2}^2\dd t\Bigr)^{2/3}
   \Bigl(\int_{S_0}^\infty\tau^{3a-3}\dd t\Bigr)^{1/3},
\]
and the last integral is finite for $a<\frac23$. By the formulas for
$h_b,A_b,B_b,\partial_xp$, \eqref{e:cmpb}, \eqref{e:wval} and
\eqref{e:munw}, $|(A_b,B_b,h_b,\partial_xp)\circ X|
\le C(|W_{yy}|+|W_y|)+C_b\tau^{-1}$, and
$\tau^a\Vert W_{yy}\Vert_{L^\infty}^2\le\tau\Vert W_{yy}\Vert_{H^1}^2$ for
$a\le1$; this gives the second bound in the first line of \eqref{e:E8}, and
\eqref{e:unwE}, \eqref{e:E1} and \eqref{e:E2} give the first. The second and
third lines follow in the same way from the first line, \eqref{e:cder},
\eqref{e:agm}, \eqref{e:wtC}, \eqref{e:wval}, \eqref{e:cmpb}, \eqref{e:munw}
and \eqref{e:curv}; for instance $\tau\Vert\partial_x^2p\Vert_{L^\infty}^2
\le C\tau\Vert(P_b)_{yy}\Vert_{L^2}\Vert(P_b)_{yyy}\Vert_{L^2}
+C_b\Vert(P_b)_y\Vert_{L^\infty}^2$, where $(P_b)_{yyy}$ contains
$(W_v)_{yyy}$ and $(W_T)_{yyy}$ of \eqref{e:wtC}, and
$\tau^{1/4}\Vert\partial_x^2\log\theta\Vert_{L^\infty}^2
\le C\tau^{1/4}\bigl(\frac{\delta^2}{\tau^2}+\Vert e_{\theta,2}\Vert_{L^\infty}^2
+\Vert\partial_x\theta\Vert_{L^\infty}^4\bigr)$.

By \eqref{e:cder} and $Z_t=\mathsf VZ_{yy}-\mathsf BZ_y+\mathsf Q$,
$\bigl(|\partial_x(f(Z)\circ Y)|+|D_t(f(Z)\circ Y)|\bigr)\circ X
\le C(|Z_y|+|(b_m,\Theta_b)_{yy}|)$. If $t\ge S_0$ and
$y\notin\Omega_\sigma(t)$, then $|y-s_\sigma t|\ge\frac{c_\sigma t}2$, since
$c_0<\frac{c_\sigma}8$, so that
$(|\xi_\sigma|^ne^{-c_1\xi_\sigma^2})\circ X\le C_b\tau^{-3}$ there, while
$\sup_{t\ge S_0}\Vert Z_y\Vert_{W^{1,\infty}}<\infty$ by \eqref{e:munw}. On
$\Omega_\sigma(t)$, by \eqref{e:extb},
$|Z_y|+|(b_m,\Theta_b)_{yy}|\le|W_y|+|W_{yy}|+C_b\tau^{-1}$, and by
\eqref{e:agm}
\[
 \sqrt\tau\Vert(W_y,W_{yy})\Vert_{L^\infty}^2
 \le\Vert W_y\Vert_{L^2}^2+2\tau\Vert W_{yy}\Vert_{L^2}^2
   +\Vert W_{yyy}\Vert_{L^2}^2 .
\]
By \eqref{e:wbas} and \eqref{e:wtC}, this gives \eqref{e:floc}.
\end{proof}

\section{Scalar parabolic estimates}\label{app:sc}

We retain the kernel notation $I,D_0,H,K,h,X,G$ of
Section~\ref{s:kern}. After proving Lemma~\ref{l:kern}, we establish
the additional estimates used in Appendix~\ref{s:pres}.

\subsection{Proof of the scalar kernel bounds}

\begin{proof}[Proof of Lemma~\ref{l:kern}]
Let $M=\sup_{\R\times[s,t]}|b-\bar b|$ and $M_1=\sup_{\R\times[s,t]}|b_\xi|$.
Duhamel's formula for $\pa_t-D_0\pa_\xi^2+\bar b\pa_\xi$ gives
\begin{equation}\label{e:kvol}
 H(t,s;\xi,\eta)=G(h,X)-\int_s^t\!\!\int_\R
 G_x\bigl(t-q,\xi-z-\bar b(t-q)\bigr)(b-\bar b)(z,q)H(q,s;z,\eta)\dd z\dd q,
\end{equation}
where $G_x$ is the derivative of $G(h,x)$ in $x$. Let $H^{(0)}=G(h,X)$, and let
$H^{(n+1)}$ be the second term on the right of \eqref{e:kvol}, with its sign, with
$H$ replaced by $H^{(n)}$, so that $H=\sum_{n\ge0}H^{(n)}$. We use
\[
\begin{gathered}
 |G_x(h,x)|\le C_0h^{-\frac12}G(2h,x),\qquad G(h,x)\le\sqrt2\,G(2h,x),\qquad
 C_0=\Bigl(\frac2{eD_0}\Bigr)^{\frac12},\\
 \int_\R G(2h_1,x-z)G(2h_2,z)\dd z=G(2(h_1+h_2),x),\qquad
 \int_0^h(h-q)^{-\frac12}q^{\frac m2}\dd q
 =\frac{\Gamma(\frac12)\Gamma(\frac m2+1)}{\Gamma(\frac{m+3}2)}h^{\frac{m+1}2},
\end{gathered}
\]
the last with $m=n$ and $m=n-1$. By induction on $n$,
\[
 |H^{(n)}|\le\frac{\sqrt2\,(C_0\sqrt\pi M)^nh^{\frac n2}}{\Gamma(\frac n2+1)}G(2h,X),
 \qquad
 |\pa_\eta H^{(n)}|\le\frac{\sqrt\pi\,C_0(C_0\sqrt\pi M)^nh^{\frac{n-1}2}}
 {\Gamma(\frac{n+1}2)}G(2h,X).
\]
Differentiating \eqref{e:kvol} in $\xi$ and integrating by parts in $z$,
\[
 H_\xi=G_x(h,X)-\int_s^t\!\!\int_\R G_x\bigl(t-q,\xi-z-\bar b(t-q)\bigr)
 \bigl[(b-\bar b)H_z+b_zH\bigr](z,q;\eta)\dd z\dd q .
\]
The same Volterra iteration, with inhomogeneous term bounded by
$C(1+(q-s)M_1)(q-s)^{-\frac12}G(2(q-s),\cdot)$, gives
$\epsilon_1>0$, depending only
on $D_0$, such that $M\sqrt h\le\epsilon_1$ implies
\[
 |H-G(h,X)|\le CM\sqrt h\,G(2h,X),\qquad
 \sqrt h\,|H_\eta|\le CG(2h,X),\qquad
 \sqrt h\,|H_\xi|\le C(1+hM_1)G(2h,X),
\]
and, for $|X|\le2\sqrt h$,
$H\ge(4\pi D_0h)^{-\frac12}(e^{-\frac1{D_0}}-CM\sqrt h)
\ge c_2h^{-\frac12}$ with $c_2=\frac12(4\pi D_0)^{-\frac12}e^{-\frac1{D_0}}<1$.
Choose $\varpi\le1$ with $A_0\sqrt\varpi\le\epsilon_1$. If $h\le\varpi\tau_s$,
then \eqref{e:kH} gives $M\sqrt h\le A_0\sqrt{h/\tau_s}\le\epsilon_1$ and
$hM_1\le A_1$.

For the lower bound for all $X$, let $N=1+\lfloor\frac{X^2}h\rfloor$,
$t_i=s+\frac{ih}N$, $z_i=\eta+\frac iN(\xi-\eta)$ and
$B_i=[z_i-\frac12\sqrt{h/N},z_i+\frac12\sqrt{h/N}]$ for $0<i<N$, with
$B_0=\{\eta\}$, $B_N=\{\xi\}$. For $x\in B_i$, $x'\in B_{i+1}$,
$|x'-x-\bar b\frac hN|\le\frac{|X|}N+\sqrt{\frac hN}\le2\sqrt{\frac hN}$, so
$H(t_{i+1},t_i;x',x)\ge c_2\sqrt{N/h}$. Since $H\ge0$,
\[
 H(t,s;\xi,\eta)\ge\int_{B_1}\!\!\cdots\int_{B_{N-1}}\prod_{i=0}^{N-1}
 H(t_{i+1},t_i;\eta_{i+1},\eta_i)\dd\eta_1\cdots\dd\eta_{N-1}
 \ge c_2^N\sqrt{\frac Nh}\ge c_2h^{-\frac12}e^{-\log\frac1{c_2}\,\frac{X^2}h},
\]
with $\eta_0=\eta$, $\eta_N=\xi$. This proves the bounds in (a) for $H$, $H_\xi$ and
$H_\eta$. In $(s,\eta)$, $K$ solves $\pa_sK+D_0K_{\eta\eta}+\pa_\eta(bK)=0$ and
tends to $\delta(\eta-\xi)$ as $s\to t$. Hence $K(t,s;\xi,\eta)$ is the kernel
$H(h,0;\eta,\xi)$ for the drift $-b(\eta,t-q)$, $0\le q\le h$, with $-\bar b$ in
place of $\bar b$; this drift has the same $M$ and $M_1$, and
$\eta-\xi+\bar bh=-X$. Thus $K$, $K_\xi$, $K_\eta$ satisfy the bounds proved for
$H$, $H_\eta$, $H_\xi$, respectively. The mixed bounds follow from
\[
 H_{\xi\eta}(t,s;\xi,\eta)=\int_\R H_\xi(t,m;\xi,z)H_\eta(m,s;z,\eta)\dd z,
 \qquad m=\frac{s+t}2,
\]
and the same identity for $K$. The functions $K_\xi+K_\eta$ and $H_\xi+H_\eta$
vanish at $t=s$ and solve
$\pa_tw=D_0w_{\xi\xi}-bw_\xi-b_\xi K_\xi$ and
$\pa_tw=D_0w_{\xi\xi}-\pa_\xi(bw)-\pa_\xi(b_\xi H)$, respectively. By Duhamel's
formula,
\begin{equation}\label{e:kdef}
\begin{aligned}
 (K_\xi+K_\eta)(t,s;\xi,\eta)&=-\int_s^t\!\!\int_\R K(t,q;\xi,z)b_z(z,q)
 K_z(q,s;z,\eta)\dd z\dd q,\\
 (H_\xi+H_\eta)(t,s;\xi,\eta)&=\int_s^t\!\!\int_\R H_z(t,q;\xi,z)b_z(z,q)
 H(q,s;z,\eta)\dd z\dd q.
\end{aligned}
\end{equation}
The defect bounds follow from \eqref{e:kdef}, its $\xi$-derivative, the bounds
already proved, $\int_\R G(Ch_1,x-z)G(Ch_2,z)\dd z=G(C(h_1+h_2),x)$ and
$\int_s^t(t-q)^{-\frac12}(q-s)^{-\frac12}\dd q=\pi$. This proves (a). Only
$H_\xi$, $K_\eta$ and the mixed derivatives used $M_1$.

For (b), since $\pa_tX=-\bar b$ and $\pa_\xi X=1$, the equations of $H$ and $K$
give, for $\mu\in\R$,
\begin{equation}\label{e:kmom}
\begin{aligned}
 \frac{\dd}{\dd t}\int_\R e^{\mu X}H\dd\xi
 &=\int_\R\bigl(D_0\mu^2+\mu(b-\bar b)\bigr)e^{\mu X}H\dd\xi,\\
 \frac{\dd}{\dd t}\int_\R e^{\mu X}K\dd\xi
 &=\int_\R\bigl(D_0\mu^2+\mu(b-\bar b)+b_\xi\bigr)e^{\mu X}K\dd\xi .
\end{aligned}
\end{equation}
By \eqref{e:kH} and $\int_s^t\tau_q^{-\frac12}\dd q\le2\sqrt h$, the two integrals
are at most $e^{D_0\mu^2h+2|\mu|A_0\sqrt h}$ and
$e^{D_0\mu^2h+2|\mu|A_0\sqrt h}(\frac{\tau_t}{\tau_s})^{A_1}$; $\mu=0$ gives
$\int_\R K\dd\xi\le(\frac{\tau_t}{\tau_s})^{A_1}$. Taking $\mu=\pm\frac x{2D_0h}$
and using $\frac{A_0x}{D_0\sqrt h}\le\frac{x^2}{8D_0h}+\frac{2A_0^2}{D_0}$,
\[
 \int_{|X|\ge x}H\dd\xi\le2e^{\frac{2A_0^2}{D_0}-\frac{x^2}{8D_0h}},\qquad
 \int_{|X|\ge x}K\dd\xi\le2e^{\frac{2A_0^2}{D_0}-\frac{x^2}{8D_0h}}
 \Bigl(\frac{\tau_t}{\tau_s}\Bigr)^{A_1}\qquad(x\ge0).
\]
Let $\ell=\frac{\varpi h}2$. Since $\tau_{t-\ell}\ge\frac{\tau_t}2\ge\frac h2$,
(a) applies on $[t-\ell,t]$. Write
\[
 H(t,s;\xi,\eta)=\int_\R H(t,t-\ell;\xi,z)H(t-\ell,s;z,\eta)\dd z,\qquad
 Z=z-\eta-\bar b(h-\ell),
\]
so that $\xi-z-\bar b\ell=X-Z$. On $\{|Z|\ge\frac{|X|}2\}$ bound the first factor
by (a) and integrate the second by the tail bound; on
$\{|X-Z|\ge\frac{|X|}2\}$ bound the first factor by (a) and use
$\int_\R H(t-\ell,s;z,\eta)\dd z=1$. This gives
$H\le C\ell^{-\frac12}e^{-\frac{X^2}{Ch}}\le CG(Ch,X)$, and the same computation
with $K$ gives the bound for $K$. Then
$H_\xi(t,s;\xi,\eta)=\int_\R H_\xi(t,t-\ell;\xi,z)H(t-\ell,s;z,\eta)\dd z$, (a)
and $\int_\R G(Ch_1,x-z)G(Ch_2,z)\dd z=G(C(h_1+h_2),x)$ give
$|H_\xi|\le C\ell^{-\frac12}G(Ch,X)$, and likewise for $K_\xi$. Let
$\ell'=\min\{\frac h2,\varpi\tau_s\}$. By
\[
 H_\eta(t,s;\xi,\eta)=\int_\R H(t,s+\ell';\xi,z)H_\eta(s+\ell',s;z,\eta)\dd z,
\]
(a) on $[s,s+\ell']$, the bound for $H$ and the same convolution identity,
$|H_\eta|\le C\ell'^{-\frac12}G(Ch,X)$, and
$\ell'^{-\frac12}\le\sqrt2\,h^{-\frac12}+\varpi^{-\frac12}\tau_s^{-\frac12}$. The
same with $K$ gives the bound for $K_\eta$, and
$K_{\xi\eta}(t,s;\xi,\eta)=\int_\R K_\xi(t,m;\xi,z)K_\eta(m,s;z,\eta)\dd z$ gives
that for $K_{\xi\eta}$. The bounds for $H$ and $H_\eta$ used only the first
inequality in \eqref{e:kH}.

For (c), uniqueness in \eqref{e:Tdef} gives $u=T_{P,\Omega}f$ and
$\bar u=T_{P,\Omega}|f|$. Since $T_0\le\varpi\tau_r$, (a) applies throughout.
For $h\ge1$, (a),
$\sup_\xi\Vert\mathcal Sw(\xi,\cdot)\Vert\le\Vert\Psi\Vert_{L^1}\sup_\xi\Vert w(\xi,\cdot)\Vert$
and
$\Vert G(Ch,\xi-\cdot-\bar bh)\Vert_{L^2}\le Ch^{-\frac14}$ bound the three
integrands by $Ch^{-\frac34}$, $Ch^{-\frac54}$ and $Ch^{-1}$. For $h<1$,
$\mathcal S\pa_\xi K=\Psi'*K$ and $\mathcal S\pa_\xi\pa_\eta K=\Psi'*K_\eta$, so
(a) and $\int_\R K\dd\xi\le C$ bound them by $Ch^{-\frac14}$, $Ch^{-\frac34}$ and
$C$. Integration in $h$ gives the three integral bounds; the fourth bound is (a)
with $\int_\R G(Ch,X)\dd\eta=1$. Since $\frac{1-\hat\Psi(k)}k\in L^2(\R)$,
\[
 \Vert g-\mathcal Sg\Vert_{L^\infty}\le\frac1{2\pi}\int_\R
 \Bigl|\frac{1-\hat\Psi(k)}k\Bigr|\,|k\hat g(k)|\dd k
 \le\frac1{\sqrt{2\pi}}\Bigl\Vert\frac{1-\hat\Psi(k)}k\Bigr\Vert_{L^2}\Vert g_\xi\Vert_{L^2}.
\]

Let $C_1\ge1$ be a constant for which (a) holds with $C=C_1$,
and put $h_j=2^{-j}T_0$ and $I_j=(t_0-h_j,t_0-h_{j+1}]$ for $j\ge0$. For $s\in I_j$
the age $h=t_0-s$ lies in $[\frac{h_j}2,h_j)$, and
$C_1h\le C_1h_j+\frac h{C_1}\le(2C_1+2C_1^{-1})h$, so
\[
 \int_\R G(C_1h_j,\xi_0-\xi')G\bigl(C_1^{-1}h,\xi'-\eta-\bar bh\bigr)\dd\xi'
 =G\bigl(C_1h_j+C_1^{-1}h,\xi_0-\eta-\bar bh\bigr)
 \ge\frac{G(C_1h,\xi_0-\eta-\bar bh)}{\sqrt{2+2C_1^{-2}}}.
\]
By the upper bound for $K_\xi$ and the lower bound for $K$ in (a), and
$h^{-\frac12}\le\sqrt2\,h_j^{-\frac12}$, this gives
\begin{equation}\label{e:kcmp}
 |K_\xi(t_0,s;\xi_0,\eta)|
 \le Ch_j^{-\frac12}\int_\R G(C_1h_j,\xi_0-\xi')K(t_0,s;\xi',\eta)\dd\xi'
 \qquad(s\in I_j).
\end{equation}
Multiply \eqref{e:kcmp} by $|f(\eta,s)|$ and integrate over $\eta\in\R$ and
$s\in I_j$. Since $K\ge0$, $\bar u(\cdot,t_0)\le\varepsilon_0$ on
$(\xi_0-\rho,\xi_0+\rho)$, $\bar u\le F_0T_0$ by $\int_\R K\dd\eta=1$, and
$\int_{|x|\ge\rho}G(C_1h_j,x)\dd x\le\sqrt2\,e^{-\frac{\rho^2}{8D_0C_1T_0}}$, we
obtain, with $\varepsilon=\varepsilon_0+F_0T_0e^{-\frac{\rho^2}{8D_0C_1T_0}}$,
\[
 \Bigl|\int_{I_j}\!\int_\R K_\xi(t_0,s;\xi_0,\eta)f(\eta,s)\dd\eta\dd s\Bigr|
 \le Ch_j^{-\frac12}\int_\R G(C_1h_j,\xi_0-\xi')\bar u(\xi',t_0)\dd\xi'
 \le Ch_j^{-\frac12}\varepsilon .
\]
By (a), the same integral is also at most
$F_0\int_{I_j}\int_\R|K_\xi|\dd\eta\dd s\le CF_0\sqrt{h_j}$. Summing over $j$,
\[
 |u_\xi(\xi_0,t_0)|\le C\sum_{j\ge0}\min\Bigl\{F_0\sqrt{h_j},\frac{\varepsilon}{\sqrt{h_j}}\Bigr\}
 \le C\sum_{h_j\le\varepsilon/F_0}F_0\sqrt{h_j}
 +C\sum_{h_j>\varepsilon/F_0}\frac{\varepsilon}{\sqrt{h_j}}
 \le\frac{2C}{1-2^{-1/2}}\sqrt{F_0\varepsilon},
\]
which proves \eqref{e:kres}.

For the last assertion of (c), $w=\bar u-T_{P,J\times(r,t_0)}|f|$ solves $Pw=0$
in $J\times(r,t_0)$, $w(\cdot,r)=0$, and $w=\bar u$ on $\pa J\times(r,t_0)$. For
$\xi\in\pa J$, $r<s<t\le t_0$ and $|\eta-\xi_0|\le\rho$,
$|\xi-\eta-\bar b(t-s)|\ge\frac\rho2-|\bar b|T_0\ge\frac{3\rho}8$, so (a) gives
\[
 \bar u(\xi,t)\le F_0\int_r^t\!\!\int_{|\xi-\eta-\bar b(t-s)|\ge\frac{3\rho}8}
 C_1G\bigl(C_1(t-s),\xi-\eta-\bar b(t-s)\bigr)\dd\eta\dd s
 \le\sqrt2\,C_1F_0T_0e^{-\frac{9\rho^2}{512D_0C_1T_0}} .
\]
By the maximum principle $w$ is bounded by the right-hand side in
$J\times(r,t_0]$, which proves the assertion. If $b\equiv\bar b$, then $H$ and
$K$ are the Gaussians $G(h,X)$, so (a) holds for all $h$ and the proof of (c)
needs no restriction on $T_0$.
\end{proof}

\subsection{Additional kernel bounds}

\begin{corollary}\label{l:kernv}
\begin{enumerate}
\item[(i)] Let $K^L$ and $K^{\mathscr P}$ be the kernels $K$ for $I=[0,\infty)$,
$D_0=d(\rho_-^{-1})$ and the drifts $b$ of $L_\sigma$ and of $\mathscr P_\sigma$
in \eqref{e:nco}, and let $\bar b=s_\sigma(\rho_-^{-1})$. Then
$e^{-C_M}K^L\le K^{\mathscr P}\le e^{C_M}K^L$, and $K^{\mathscr P}$ satisfies the
pointwise bounds for $K$, $K_\xi$, $K_\eta$, $K_{\xi\eta}$ in
\cref{l:kern}\textup{(b)} with $C$ replaced by $C_M$ and $A_1$ by $C\delta$.
\item[(ii)] Let $v_0>0$ and $A_2\ge0$. Assume the first inequality of
\eqref{e:kH}, $|b_\xi|\le A_2$ on $\R\times I$, and $|b_\xi|\le A_1\tau^{-1}$ on
$\{(\xi,t)\in\R\times I:|\xi|\le v_0\tau\}$. Then the bounds for $H$, $H_\eta$ and
$H_\xi+H_\eta$ in \cref{l:kern}\textup{(a)}, \textup{(b)} hold, and
$|H_\xi|\le C(h^{-\frac12}+A_2^{\frac12})G(Ch,X)$, with $C$ depending only on
$D_0$ and an upper bound for $A_0$. There are $\varpi'\in(0,\varpi]$ and $C'$,
depending only on $D_0$, $v_0$, $|\bar b|$ and upper bounds for $A_0$ and $A_1$,
such that
\[
\begin{aligned}
 |H_\xi(t,s;\xi,\eta)|&\le C'h^{-\frac12}G(C'h,X)
 &&\bigl(|\xi|\le\tfrac{v_0}2\tau_t\bigr),\\
 |(H_\xi+H_\eta)(t,s;\xi,\eta)|&\le C'\sqrt h\,\bigl(A_1\tau_s^{-1}
 +A_2e^{-\frac{\tau_s}{C'}}\bigr)G(C'h,X)
 &&\bigl(h\le\varpi'\tau_s,\ |\eta|\le\tfrac{v_0}2\tau_s\bigr).
\end{aligned}
\]
\item[(iii)] Assume \eqref{e:kH}. There is $C_*\ge1$, depending only on $D_0$ and
upper bounds for $A_0$ and $A_1$, such that, if $0<w_2<w_1$ and $C_*w_1w_2\le w_1-w_2$,
then for $s<t$ in $I$, $S=\tau_s$ and
$v(\xi)=\int_\R H(t,s;\xi,\eta)g(\eta)\dd\eta$,
\[
\begin{aligned}
 \Bigl\Vert e^{\frac{w_2(\xi-\bar bh)^2}{S+h}}v\Bigr\Vert_{L^2}
 &\le C\Bigl\Vert e^{\frac{w_1\eta^2}S}g\Bigr\Vert_{L^2},\\
 \Bigl\Vert e^{\frac{w_2(\xi-\bar bh)^2}{S+h}}v_\xi\Bigr\Vert_{L^2}
 &\le C\Bigl(\Bigl\Vert e^{\frac{w_1\eta^2}S}g_\eta\Bigr\Vert_{L^2}
 +S^{-\frac12}\Bigl\Vert e^{\frac{w_1\eta^2}S}g\Bigr\Vert_{L^2}\Bigr).
\end{aligned}
\]
\item[(iv)] Let $s<T$, let $\mathcal D$ be a measurable $k\times k$ matrix function
on $\R\times(s,T)$ with $|\mathcal D|\le A$, and let
$Q\in C([s,T];L^2(\R;\R^k))\cap L^2(s,T;H^1(\R;\R^k))$ be a weak solution of
$Q_t=\pa_\xi(D_0Q_\xi-\mathcal DQ)$ with $Q(s)\in L^1$. Then
$Q\in C([s,T];L^1)$ and
$\Vert Q(t)\Vert_{L^1}\le2^{1+\frac{16A^2(t-s)}{\pi D_0}}\Vert Q(s)\Vert_{L^1}$.
\end{enumerate}
\end{corollary}

\begin{proof}
(i) By \cref{l:pr}, $L_\sigma(p_\sigma z)=p_\sigma(\mathscr P_\sigma z
-\alpha R^{\rm prof}_\sigma z)$, $e^{-C_M}\le p_\sigma\le e^{C_M}$ and
$\int_0^\infty\Vert R^{\rm prof}_\sigma\Vert_{L^\infty}\dd t\le C_M$. Let
$V(F_c(y,t),t)=\alpha R^{\rm prof}_\sigma(y,t)$, let $b$ be the drift of $L_\sigma$,
and let $\tilde K$ be the kernel of $\pa_t-D_0\pa_\xi^2+b\pa_\xi+V$. Thus
$z\mapsto p_\sigma z$ maps solutions of $\mathscr P_\sigma z=0$ to solutions of
$L_\sigma w+\alpha R^{\rm prof}_\sigma w=0$, and by \eqref{e:nco}
\[
 K^{\mathscr P}(t,s;F_c(y,t),F_c(z,s))=\frac{p_\sigma(z,s)}{p_\sigma(y,t)}\,
 \tilde K(t,s;F_c(y,t),F_c(z,s)).
\]
If $w\ge0$ solves $w_t-D_0w_{\xi\xi}+bw_\xi=0$, then
$(\pa_t-D_0\pa_\xi^2+b\pa_\xi+V)\bigl(e^{\pm\int_s^t\Vert V\Vert_{L^\infty}}w\bigr)
=(V\pm\Vert V\Vert_{L^\infty})e^{\pm\int_s^t\Vert V\Vert_{L^\infty}}w$, so the
maximum principle gives $e^{-C_M}K^L\le\tilde K\le e^{C_M}K^L$, and the comparison
of $K^{\mathscr P}$ with $K^L$ follows. By \eqref{e:kH} for $L_\sigma$ with
$A_0=A_1=C\delta$, \cref{l:kern}(b) bounds $K^{\mathscr P}$ by
$C_M(\frac{\tau_t}{\tau_s})^{C\delta}G(C_Mh,X)$. The derivative bounds follow from
this bound by the composition steps in the proof of \cref{l:kern}(b), with
\cref{l:kern}(a) for the drift of $\mathscr P_\sigma$, which satisfies
\eqref{e:kH} with $A_0=A_1=C_M$.

(ii) The bounds for $H$, $H_\eta$ and $H_\xi+H_\eta$ in \cref{l:kern}(a), (b)
used only the first inequality of \eqref{e:kH}. With $M,M_1,\epsilon_1$
as in the short-time proof of \cref{l:kern}(a),
$\ell\le\min\{\varpi\tau_{t-\ell},A_2^{-1}\}$ gives
$M\sqrt\ell\le\epsilon_1$ and $\ell M_1\le1$ on $[t-\ell,t]$. Hence $\sqrt\ell\,|H_\xi(t,t-\ell;\xi,z)|\le CG(C\ell,\xi-z-\bar b\ell)$.
The proof of \cref{l:kern}(b) with $\ell=\min\{\frac{\varpi h}2,A_2^{-1}\}$ and
$\ell^{-\frac12}\le(\frac{\varpi h}2)^{-\frac12}+A_2^{\frac12}$ gives the global
bound for $H_\xi$. For the central bound, let $\chi\in C_c^\infty(\R)$ with
$0\le\chi\le1$, $\chi=1$ on $[-\frac78,\frac78]$ and $\chi=0$ outside $[-1,1]$, and
\[
 b^p(\xi,t)=\bar b+\chi\Bigl(\frac\xi{v_0\tau}\Bigr)(b-\bar b),\qquad
 |b^p-\bar b|\le A_0\tau^{-\frac12},\qquad
 |b^p_\xi|\le\Bigl(A_1+\frac{\Vert\chi'\Vert_{L^\infty}A_0}{v_0}\Bigr)\tau^{-1}.
\]
Let $H^p$ be the kernel $H$ for $b^p$; \cref{l:kern} applies to it. Let
$\epsilon_2=\min\{\frac18,\frac{v_0}{8(1+|\bar b|)}\}$,
$\ell=\min\{\frac{\varpi h}2,\epsilon_2\tau_t\}$ and $t'=t-\ell$. Since
$v_t=D_0v_{\xi\xi}-\pa_\xi(b^pv)-\pa_\xi((b-b^p)v)$,
\[
 H_\xi(t,s;\xi,\eta)=\int_\R H^p_\xi(t,t';\xi,z)H(t',s;z,\eta)\dd z
 +\int_{t'}^t\!\!\int_\R H^p_{\xi z}(t,q;\xi,z)(b-b^p)(z,q)H(q,s;z,\eta)\dd z\dd q .
\]
The first term is at most $C\ell^{-\frac12}G(Ch,X)\le Ch^{-\frac12}G(Ch,X)$, by
\cref{l:kern}(a) for $H^p$, \cref{l:kern}(b) for $H$ and $h\le\tau_t$. In the second, $b-b^p=0$ for
$|z|\le\frac78v_0\tau_q$, and $\tau_q\ge\frac78\tau_t$; so for
$|\xi|\le\frac{v_0}2\tau_t$ the integrand vanishes unless
$|\xi-z-\bar b(t-q)|\ge\frac{17}{64}v_0\tau_t-|\bar b|\epsilon_2\tau_t
\ge\frac{v_0}8\tau_t$. There
$|H^p_{\xi z}(t,q;\xi,z)|\le C(t-q)^{-1}G(C(t-q),\cdot)$ and
$G(Cx,y)\le\sqrt2\,e^{-\frac{v_0^2\tau_t^2}{512D_0Cx}}G(2Cx,y)$ for
$|y|\ge\frac{v_0}8\tau_t$. With $|b-b^p|\le A_0\tau_q^{-\frac12}$,
$H(q,s)\le CG(C(q-s),\cdot)$ and
$\int_0^\ell x^{-1}e^{-\frac{v_0^2\tau_t^2}{512D_0Cx}}\dd x\le C$, the second term
is at most $C\tau_t^{-\frac12}G(Ch,X)\le Ch^{-\frac12}G(Ch,X)$. For the defect, let
$\varpi'=\min\{\varpi,\frac{v_0}{4(1+|\bar b|)}\}$ and split the integral in
\eqref{e:kdef} into $|z|\le v_0\tau_q$, where $|b_z|\le A_1\tau_s^{-1}$, and
$|z|>v_0\tau_q$, where $|b_z|\le A_2$ and
$|z-\eta-\bar b(q-s)|\ge\frac{v_0}2\tau_s-|\bar b|\varpi'\tau_s\ge\frac{v_0}4\tau_s$;
in the second part
$G(C(q-s),z-\eta-\bar b(q-s))\le\sqrt2\,e^{-\frac{v_0^2\tau_s}{128D_0C\varpi'}}
G(2C(q-s),z-\eta-\bar b(q-s))$. The bounds of \cref{l:kern}(a) for $H_z$ and $H$
then give the bound for $H_\xi+H_\eta$.

(iii) Let $C_1$ be a constant for which the bounds for $H$, $H_\xi$ and
$H_\xi+H_\eta$ in \cref{l:kern}(a), (b) hold with $C=C_1$, and
$C_*=\max\{1,8D_0C_1\}$. Since $\xi-\bar bh=X+\eta$,
\[
 \frac{w_2(\xi-\bar bh)^2}{S+h}
 \le\frac{w_1\eta^2+\frac{w_1w_2}{w_1-w_2}X^2}{S+h}
 \le\frac{w_1\eta^2}S+\frac{w_1w_2}{w_1-w_2}\,\frac{X^2}h
 \le\frac{w_1\eta^2}S+\frac{X^2}{8D_0C_1h}.
\]
Hence $e^{\frac{w_2(\xi-\bar bh)^2}{S+h}}H(t,s;\xi,\eta)e^{-\frac{w_1\eta^2}S}
\le\sqrt2\,C_1G(2C_1h,X)$, whose integrals in $\xi$ and in $\eta$ are at most
$\sqrt2\,C_1$; Schur's test gives the first bound. For $h\ge\varpi S$,
$|H_\xi|\le C(\varpi S)^{-\frac12}G(Ch,X)$ by \cref{l:kern}(b). For $h<\varpi S$,
$v_\xi=\int_\R Hg_\eta\dd\eta+\int_\R(H_\xi+H_\eta)g\dd\eta$ and
$|H_\xi+H_\eta|\le CA_1\sqrt{\varpi}\,S^{-\frac12}G(Ch,X)$ by \cref{l:kern}(a).
Schur's test, with the same splitting of the exponent, gives the second bound.

(iv) If $A=0$, the heat semigroup is contractive in $L^1$, so the assertion follows.
Let $A>0$, $s'\in[s,T)$ and $h_0=\frac{\pi D_0}{16A^2}$. Since
$\Vert G_x(h,\cdot)\Vert_{L^1}=(\pi D_0h)^{-\frac12}$, the map
\[
 W\mapsto G(t-s')*Q(s')-\int_{s'}^tG_x(t-q)*(\mathcal DW)(q)\dd q,
\]
where $G(h)*$ is convolution with $G(h,\cdot)$, is a contraction with constant
$2A\sqrt{h_0/(\pi D_0)}=\frac12$ on
$C([s',s'+h_0]\cap[s,T];L^1\cap L^2)$. Its fixed point $W$ satisfies
$\sup_t\Vert W(t)\Vert_{L^1}\le2\Vert Q(s')\Vert_{L^1}$. The weak solution $Q$
satisfies the same Duhamel formula, so $Q=W$ by the contraction in $C(L^2)$.
Applying this on $[s+jh_0,s+(j+1)h_0]$, $j\ge0$, proves (iv).
\end{proof}

\subsection{Uniqueness and reference-operator estimates}

Let $r<T<\infty$, let $d$ be measurable with $d_{\min}\le d\le d_{\max}$ for constants
$0<d_{\min}\le d_{\max}$, and let $B_0,B_1,B_2$ be bounded measurable matrix functions on
$\R\times(r,T)$. If $v\in C([r,T];L^2)$ with $v_y\in L^2((r,T)\times\R)$ is a weak solution of
\begin{equation}\label{e:uniq}
 v_t=\pa_y(dv_y-B_1v)+B_2v_y+B_0v,\qquad v(r)=0,
\end{equation}
then $v=0$. Indeed, the equation gives $v_t\in L^2(r,T;H^{-1})$, and testing by $v$ gives
\[
 \frac12\frac{\dd}{\dd t}\Vert v\Vert_{L^2}^2+\int_\R d|v_y|^2\dd y
 =\int_\R\bigl(v_y^\top B_1v+v^\top B_2v_y+v^\top B_0v\bigr)\dd y
 \le\frac{d_{\min}}2\Vert v_y\Vert_{L^2}^2+C\Vert v\Vert_{L^2}^2,
\]
with $C$ depending on $d_{\min}$ and the bounds of $B_0,B_1,B_2$; Gronwall's inequality gives
$v=0$.

In this appendix $\sigma=\pm1$, $d=d(V_c)$, $\varsigma=s_\sigma(V_c)$, and $L$ is the number of
\eqref{e:nfwdef}. For $P\in\{L_\sigma,\mathscr P_\sigma\}$, $s<t$ and $g\in L^1\cap L^2$,
$S_P(t,s)g$ is the value at time $t$ of the solution $u$ of $Pu=0$ on $\R\times(s,t)$ with
$u(s)=g$.

\begin{lemma}\label{l:c2}
For $h\in L^1(\R)\cap H^1(\R)$,
\begin{equation}\label{e:nash}
 \Vert h\Vert_{L^2}^6\le\Vert h\Vert_{L^1}^4\Vert h_y\Vert_{L^2}^2 .
\end{equation}
If $\tau_r\ge L^2$, then, for $t>s\ge r$ and $P\in\{L_\sigma,\mathscr P_\sigma\}$,
\begin{equation}\label{e:ck}
 \Vert S_P(t,s)\Vert_{L^1\to L^1}\le C_M,\qquad
 \Vert S_P(t,s)\Vert_{L^1\to L^\infty}\le C_M(t-s)^{-\frac12} .
\end{equation}
\end{lemma}
\begin{proof}
Since $h(y)^2=2\int_{-\infty}^yhh_y\dd z=-2\int_y^\infty hh_y\dd z$,
$\Vert h\Vert_{L^\infty}^2\le\Vert h\Vert_{L^2}\Vert h_y\Vert_{L^2}$; with
$\Vert h\Vert_{L^2}^2\le\Vert h\Vert_{L^\infty}\Vert h\Vert_{L^1}$ this gives \eqref{e:nash}.

Let $h$ solve $L_\sigma h=0$ on $\R\times(s,t)$ with $h(s)=g\in L^1\cap H^1$ and $g\ge0$; then
$h\ge0$. With $w_0,w_1$ of \eqref{e:nfwdef}, the two energy identities in the proof of
Lemma~\ref{l:nf}(v), with $F=G=0$, also hold for this nonzero initial datum. Together with
the test against $w_0$ and \eqref{e:nfwts}, they give
\[
\begin{gathered}
 \frac{\dd}{\dd t}\int_\R w_0h\dd y\le0,\qquad
 \frac12\frac{\dd}{\dd t}\int_\R w_0h^2\dd y+\int_\R w_0dh_y^2\dd y\le0,\\
 \frac{\dd}{\dd t}\int_\R w_1h_y^2\dd y\le0 .
\end{gathered}
\]
Constants $C$, $c_*$ below depend only
on the fixed data. By $C^{-1}\le w_i\le C$, the first inequality gives
$\Vert h(t')\Vert_{L^1}\le C\Vert g\Vert_{L^1}$ for $t'\ge s$, and the second, with
$E=\frac12\int w_0h^2$ and \eqref{e:nash}, gives $E'\le-c_*E^3\Vert g\Vert_{L^1}^{-4}$, hence
$E(t')\le C(t'-s)^{-\frac12}\Vert g\Vert_{L^1}^2$. Integrating the second inequality over
$(\frac{s+t}2,t)$ and using that $\int w_1h_y^2$ does not increase,
\[
 \frac{t-s}2\Vert h_y(t)\Vert_{L^2}^2\le C\int_{\frac{s+t}2}^t\Vert h_y\Vert_{L^2}^2\dd t'
 \le CE\Bigl(\frac{s+t}2\Bigr)\le C(t-s)^{-\frac12}\Vert g\Vert_{L^1}^2 ,
\]
and $\Vert h(t)\Vert_{L^\infty}^2\le\Vert h(t)\Vert_{L^2}\Vert h_y(t)\Vert_{L^2}\le C(t-s)^{-1}
\Vert g\Vert_{L^1}^2$. Applied to $g^+$ and $g^-$, and by density, this is \eqref{e:ck} for
$L_\sigma$, with $C$ in place of $C_M$.

If $\mathscr P_\sigma z=0$, then $\tilde h=p_\sigma z$ solves
$L_\sigma\tilde h+\alpha R^{\rm prof}_\sigma\tilde h=0$ by Lemma~\ref{l:pr}. With
$\epsilon(t)=\alpha\Vert R^{\rm prof}_\sigma(t)\Vert_{L^\infty}$, the function
$e^{\int_s^t\epsilon}S_{L_\sigma}(t,s)|\tilde h(s)|$ is a supersolution of this equation, so
comparison gives $|\tilde h(t)|\le e^{\int_s^t\epsilon}S_{L_\sigma}(t,s)|\tilde h(s)|$. By
Lemma~\ref{l:pr}, $e^{-C_M}\le p_\sigma\le e^{C_M}$ and
$\int_0^\infty\epsilon\dd t\le\frac{\alpha C\delta}{c_e}$; hence \eqref{e:ck} for
$\mathscr P_\sigma$ follows from that for $L_\sigma$.
\end{proof}

\begin{lemma}\label{l:c4}
Let $P\in\{L_\sigma,\mathscr P_\sigma\}$, let $b$ be its drift in \eqref{e:nco},
$\mathsf P=\pa_t-D_0\pa_\xi^2+b\pa_\xi$, and let $\varpi$ be the number of
Lemma~\ref{l:kern} for $b$.
\begin{enumerate}[label=(\roman*)]
\item Let $t_1<t_2$ with $t_2-t_1\le\varpi\tau_{t_1}$, and let
 $U\in C([t_1,t_2];H^1)\cap L^2(t_1,t_2;H^2)$ solve $\mathsf PU=F$ with $U$ and $F$ bounded.
 Then, for $t_1<t\le t_2$,
 \begin{equation}\label{e:cdh}
  \Vert U_\xi(t)\Vert_{L^\infty}\le C_M(t-t_1)^{-\frac12}\Vert U(t_1)\Vert_{L^\infty}
  +C_M\int_{t_1}^t(t-s)^{-\frac12}\Vert F(s)\Vert_{L^\infty}\dd s ,
 \end{equation}
 and the same bound holds for $\Vert u_y(t)\Vert_{L^\infty}$ if $u(y,t)=U(F_c(y,t),t)$ and
 $Pu=f$ with $f(y,t)=F(F_c(y,t),t)$.
\item Let $r\ge0$, $A,K\ge0$, and let $z\in C([r,T];H^3)\cap L^2(r,T;H^4)$ for every finite
 $T>r$ solve $\mathscr P_\sigma z=f$ on $\R\times[r,\infty)$ with
 $\sup_{t\ge r}\Vert z(t)\Vert_{L^\infty}\le A$ and
 $\Vert\pa_y^jf(t)\Vert_{L^\infty}\le K\tau^{-1-\frac j2}$ $(0\le j\le2,\ t\ge r)$. Then
 $\Vert\pa_y^jz(t)\Vert_{L^\infty}\le C_M(A+K)\tau^{-\frac j2}$ for $1\le j\le3$ and $t\ge2r+1$.
\end{enumerate}
\end{lemma}
\begin{proof}
(i) By \eqref{e:nco}, $b$ satisfies \eqref{e:kH} with $A_0=A_1=C_M$ on $[0,\infty)$, and
$\mathsf P$ has the kernel $K$ of Section~\ref{s:kern}. The function
\[
 \int_\R K(t,t_1;\cdot,\eta)U(\eta,t_1)\dd\eta
 +\int_{t_1}^t\!\!\int_\R K(t,s;\cdot,\eta)F(\eta,s)\dd\eta\dd s
\]
solves $\mathsf P\,\cdot=F$ in the class of $U$ and equals $U$ at $t_1$, so it is $U$ by
\eqref{e:uniq}. For $t_1\le s<t$, $t-s\le\varpi\tau_s$, and Lemma~\ref{l:kern}(a) gives
$|K_\xi(t,s;\xi,\eta)|\le C_M(t-s)^{-\frac12}G(C_M(t-s),X)$ with
$\int_\R G(C_M(t-s),X)\dd\eta=1$; this gives \eqref{e:cdh}. The second assertion follows from
$u_y=f_cU_\xi(F_c,t)$ and $f_c\le C$.

(ii) Let $t\ge2r+1$, $h=\frac{\varpi\tau_t}{1+\varpi}$ and $t_1=t-h$; then $t_1\ge r$,
$h=\varpi\tau_{t_1}$ and $\tau_{t_1}\le\tau_s\le\tau_t\le2\tau_{t_1}$ for $s\in[t_1,t]$. Put
$z(y,s)=U(F_c(y,s),s)$ and $f(y,s)=\hat f(F_c(y,s),s)$. By \eqref{e:nco}, $\mathsf PU=\hat f$,
$|\pa_\xi^j\hat f|\le CK\tau^{-1-\frac j2}$ $(j\le2)$, $|b_\xi|\le C_M\tau^{-1}$,
$|b_{\xi\xi}|\le C_M\tau^{-\frac32}$, and, $D_0$ being constant,
\[
 \mathsf PU_\xi=-b_\xi U_\xi+\hat f_\xi,\qquad
 \mathsf PU_{\xi\xi}=-2b_\xi U_{\xi\xi}-b_{\xi\xi}U_\xi+\hat f_{\xi\xi} .
\]
By \eqref{e:cdh} on $[t_1,s]$, $\Vert U_\xi(s)\Vert_{L^\infty}\le C_M(A+K)\tau_t^{-\frac12}$
for $s\in[t_1+\frac h2,t]$. Then the right side of the equation of $U_\xi$ is at most
$C_M(A+K)\tau_t^{-\frac32}$ there, and \eqref{e:cdh} on $[t_1+\frac h2,s]$ gives
$\Vert U_{\xi\xi}(s)\Vert_{L^\infty}\le C_M(A+K)\tau_t^{-1}$ for $s\in[t_1+\frac{3h}4,t]$;
in the same way, \eqref{e:cdh} for $U_{\xi\xi}$ on $[t_1+\frac{3h}4,t]$ gives
$\Vert U_{\xi\xi\xi}(t)\Vert_{L^\infty}\le C_M(A+K)\tau_t^{-\frac32}$. With
\[
 z_y=f_cU_\xi,\qquad z_{yy}=f_c^2U_{\xi\xi}+(f_c)_yU_\xi,\qquad
 z_{yyy}=f_c^3U_{\xi\xi\xi}+3f_c(f_c)_yU_{\xi\xi}+(f_c)_{yy}U_\xi
\]
and $|\pa_y^{j+1}f_c|\le C\delta\tau^{-\frac{j+1}2}$, this proves (ii).
\end{proof}

\section{Decomposition of the pressure}\label{s:pres}

Throughout this appendix $M\ge1$, \eqref{e:dat-aux} holds, $0<\delta\le\delta_1(M)$,
$(V_b,U_b,\Theta_b)$ is the solution of Lemma~\ref{l:bg}, $m_\pm$ and $h_*$ are
given by \eqref{e:mass}, $\mathcal B,\mathcal G,\mathcal F,\chi_m,\Psi_\sigma$ are
the functions of \eqref{e:pot}, and $y$ and $t$ are the mass coordinate and the
time of Section~\ref{s:mass}; $\tau=1+t$, $\tau_r=1+r$, the indices
$\sigma,\rho,j$ belong to $\{-1,+1\}$ and $o=-\sigma$. Vectors have the
components $(+,-)$ in this order, and the square of a vector or of a matrix is
the sum of the squares of its entries. Every restart time $r$ satisfies
$r\ge S_0$. The Gaussian exponents below and the gaps between them are fixed
before the datum; in each lemma $\delta_1(M)$, which then depends also on these
exponents, may be decreased, and $r$ may be increased.

\subsection{The pressure identity and the equation of the corrected potentials}

The first line of \eqref{e:potq} is the identity
\[
 P_b-p_*=\frac{\nu(b_m)_y}{V_b}-\mathcal B_t .
\]
For $a>0$, $0\le r\le T\le\infty$ and a function $f$ on $\R\times(r,T)$ put
\begin{equation}\label{e:nobs}
\begin{aligned}
 A_f^a(t)&=\sqrt\tau\int_\R e^{-\frac{a(y+h_*)^2}{2\tau}}|f(y,t)|^2\dd y,
 \qquad
 \mathscr O^a_{r,T}(f)=\int_r^TA_f^a(t)\dd t,\\
 \mathscr N^a_{r,T}(f)&=\int_r^T\Bigl[\tau^{\frac14}
   \bigl\Vert e^{-\frac{a(y+h_*)^2}{4\tau}}f_y\bigr\Vert_{L^2}
   +\tau^{-\frac14}\bigl\Vert e^{-\frac{a(y+h_*)^2}{4\tau}}f\bigr\Vert_{L^2}
   \Bigr]\dd t .
\end{aligned}
\end{equation}
Fix $a_s$ and $b_0$ with $0<a_s\le\frac{c_G}4$ and $0<b_0<\frac{a_s}{64}$, where
$c_G$ is the exponent of \eqref{eq:ct-gau} and \eqref{eq:pr-gau}; write
$A_f=A_f^{a_s}$, $\mathscr O_{r,T}=\mathscr O^{a_s}_{r,T}$,
$\mathscr N_{r,T}=\mathscr N^{a_s}_{r,T}$, and put
\begin{equation}\label{e:pwt}
\begin{gathered}
 z_j=\frac{y-s_jt}{\sqrt\tau},\qquad
 G^a=e^{-\frac{a(y+h_*)^2}{2\tau}},\qquad G=G^{a_s},\qquad
 G_0=e^{-\frac{b_0(y+h_*)^2}{\tau}},\qquad W_c=\tau^{-\frac12}G_0,\\
 \mathscr Q_r(v)=\int_r^\infty\!\!\int_\R W_c|v|^2\dd y\dd t,\qquad
 J_r(v)=\Vert v(r)\Vert_{L^2}^2+\mathscr Q_r(v)
   +\int_r^\infty\!\!\int_\R G_0|v_y|^2\dd y\dd t .
\end{gathered}
\end{equation}
Thus $G\le G_0^2\le G_0$ and $b_0<\frac{c_G}{256}$.

For a real measurable $2\times2$ matrix $\mathcal M$ define its residual
$\mathcal R$ by
\begin{equation}\label{e:gmat}
 \mathcal M_{\sigma\sigma}=-2k(V_c)a_\sigma+2k(V_c)a_o+\mathcal R_{\sigma\sigma},
 \qquad
 \mathcal M_{\sigma o}=2k(V_c)a_\sigma+2\ell_{\rm g}(V_c)a_o+\mathcal R_{\sigma o},
 \qquad \ell_{\rm g}(n)=\frac{3\gamma-5}{4n},
\end{equation}
with the waves of \eqref{e:CH}. With $\mathfrak b_\sigma$ of \eqref{e:pot} and
$L=\operatorname{diag}(L_+,L_-)$ put
\begin{equation}\label{e:gop}
\begin{gathered}
 S=\begin{pmatrix}\frac12&\mathfrak b_+(V_c)\\[2pt]
                  \frac12&\mathfrak b_-(V_c)\end{pmatrix},\qquad
 A_{\rm ph}=\begin{pmatrix}0&\frac1{V_c}\\[2pt]\gamma p_*&0\end{pmatrix},\qquad
 \Lambda=\begin{pmatrix}s_+(V_c)&0\\0&s_-(V_c)\end{pmatrix},\\
 d_0=\frac{(1,-1)}{2\mathfrak b_+(V_c)},\qquad
 (\mathcal C_{\mathcal M}\mathbf H)_\sigma
 =\Bigl[L_\sigma\mathfrak b_\sigma(V_c)
   -\sum_\rho\mathcal M_{\sigma\rho}\pa_y\mathfrak b_\rho(V_c)\Bigr]d_0\mathbf H
   -2d(V_c)\,\pa_y\mathfrak b_\sigma(V_c)\,\pa_y(d_0\mathbf H),\\
 \mathcal K_{\mathcal M}=L-\mathcal M\pa_y-\mathcal C_{\mathcal M},\qquad
 \mathcal K_{\mathcal M}(S\mathbf v)
 =S\bigl[\mathbf v_t-d(V_c)\mathbf v_{yy}
         +(A_{\rm ph}-S^{-1}\mathcal MS)\mathbf v_y\bigr]
\end{gathered}
\end{equation}
for every two-vector $\mathbf v$ with $\mathbf v_t$, $\mathbf v_{yy}$ locally square
integrable. Indeed, $\mathfrak b_\sigma(n)=\frac1{2s_\sigma(n)n}$ and
$s_\sigma(n)^2=\frac{\gamma p_*}n$ give
$s_\sigma(V_c)\mathfrak b_\sigma(V_c)=\frac1{2V_c}$ and
$\gamma p_*\mathfrak b_\sigma(V_c)=\frac{s_\sigma(V_c)}2$, that is,
$SA_{\rm ph}=\Lambda S$. Since $(S\mathbf v)_\sigma=\frac{v_1}2+\mathfrak b_\sigma(V_c)v_2$
and $L_\sigma(fg)=fL_\sigma g+gL_\sigma f-2d(V_c)f_yg_y$, with
$\mathfrak b=(\mathfrak b_+(V_c),\mathfrak b_-(V_c))^\top$ and
$L\mathfrak b=(L_+\mathfrak b_+(V_c),L_-\mathfrak b_-(V_c))^\top$,
\[
 L(S\mathbf v)-\mathcal M(S\mathbf v)_y
 =S\bigl[\mathbf v_t-d(V_c)\mathbf v_{yy}
       +(A_{\rm ph}-S^{-1}\mathcal MS)\mathbf v_y\bigr]
 +(L\mathfrak b-\mathcal M\mathfrak b_y)v_2-2d(V_c)\mathfrak b_y(v_2)_y ,
\]
and $\mathfrak b_-=-\mathfrak b_+$ gives $d_0S\mathbf v=v_2$; the last two terms
are $\mathcal C_{\mathcal M}(S\mathbf v)$.

Write $\mathfrak b_\sigma=\mathfrak b_\sigma(V_c)$ and $d=d(V_c)$. By
\eqref{eq:pr-gau} and \eqref{e:pt}, every derivative $\pa_y^k\pa_t^lV_c$ with
$k+l\ge1$ and $k+2l\le3$ is bounded by
$C\delta\tau^{-\frac{k+2l}2}e^{-\frac{c_Gy^2}{2\tau}}$, and
$y^2\ge\frac{(y+h_*)^2}2-h_*^2$. Expanding
$\pa_y(d_0\mathbf v)=d_0\mathbf v_y+(d_0)_y\mathbf v$ in \eqref{e:gop}, and using
$\mathfrak b_\sigma(n)=\frac\sigma{2\sqrt{\gamma p_*n}}$,
$s_\sigma(n)\mathfrak b_\sigma'(n)=-\frac1{4n^2}$ and
$(d_0)_\rho=\rho\sqrt{\gamma p_*V_c}$, we obtain, with $b_1=\frac{c_G}4$ and
for $\tau\ge h_*^2$,
\begin{equation}\label{e:pcnx}
\begin{gathered}
 \mathcal C_{\mathcal M}\mathbf v=\mathcal C^0\mathbf v+B^1\mathbf v_y+C^1\mathbf v,
 \qquad
 \mathcal C^0_{\sigma\rho}=s_\sigma(V_c)(\mathfrak b_\sigma)_y(d_0)_\rho
 =-\rho\frac{\sqrt{\gamma p_*}\,(V_c)_y}{4V_c^{3/2}},\qquad
 B^1_{\sigma\rho}=-2d(\mathfrak b_\sigma)_y(d_0)_\rho,\\
 C^1_{\sigma\rho}=\Bigl((\mathfrak b_\sigma)_t-d(\mathfrak b_\sigma)_{yy}
   -\sum_j\mathcal M_{\sigma j}(\mathfrak b_j)_y\Bigr)(d_0)_\rho
   -2d(\mathfrak b_\sigma)_y\pa_y(d_0)_\rho,\\
 |\pa_y^j\mathcal C^0|\le C\delta\tau^{-\frac{j+1}2}
   e^{-\frac{b_1(y+h_*)^2}\tau}\quad(0\le j\le2),\qquad
 |\mathcal C^0_t|\le C\delta\tau^{-\frac32}e^{-\frac{b_1(y+h_*)^2}\tau},\\
 |B^1|\le C\delta\tau^{-\frac12}e^{-\frac{b_1(y+h_*)^2}\tau},\qquad
 |C^1|\le C\delta\bigl(\tau^{-1}+\tau^{-\frac12}|\mathcal M|\bigr)
   e^{-\frac{b_1(y+h_*)^2}\tau}.
\end{gathered}
\end{equation}

With $V_h$ and $H_h$ of Lemma~\ref{l:mcoef} put $n_\chi=\chi_m-V_h$. With
$\mathfrak a_\sigma(n,U)$ of Lemma~\ref{l:pot} write
$\mathfrak a_\sigma=\mathfrak a_\sigma(V_c,\mathcal U)
=\frac{b_m}2+\mathfrak b_\sigma(V_c)\mathcal F$,
$\mathfrak a=(\mathfrak a_+,\mathfrak a_-)$, and
\[
 \mathfrak a^a_\sigma=\frac{H_h}2
   +\frac12\sum_j\Bigl(1+\frac{s_\sigma(V_c)}{s_j}\Bigr)a_j .
\]
For $n_0,n>0$ and $\zeta=(\zeta_+,\zeta_-)\in\R^2$ put
\begin{equation}\label{e:pj}
\begin{gathered}
 U^\sharp(n_0,n,\zeta)=\Bigl(v^\sharp,\ \zeta_++\zeta_-,\
   \frac{p_*v^\sharp+s_+(n_0)n_0(\zeta_+-\zeta_-)}{\gamma-1}\Bigr),
 \qquad v^\sharp=n-\frac{\zeta_+-\zeta_-}{s_+(n_0)},\\
 \mathfrak j^\sharp_\sigma(n_0,n,\zeta)
   =\mathfrak j_\sigma\bigl(n_0,U^\sharp(n_0,n,\zeta)\bigr),\qquad
 \mathcal M_{\sigma\rho}=\int_0^1\pa_{\zeta_\rho}\mathfrak j^\sharp_\sigma
   \bigl(V_c,\chi_m,\mathfrak a^a+\theta(\mathfrak a-\mathfrak a^a)\bigr)\dd\theta,
\end{gathered}
\end{equation}
with $\mathfrak j_\sigma(n,U)$ of Lemma~\ref{l:pot}. With
$\ell^0_\sigma,\Lambda_\sigma,U_*,R_\rho,R_0$ of Lemma~\ref{l:nf} put, for
$e\in\R^3$,
\begin{equation}\label{e:pcf}
\begin{gathered}
 K^0_\sigma(n)e=D_U\Lambda_\sigma(n,U_*(n))[e]R_0-s_\sigma(n)(\ell^0_\sigma)'(n)e,
 \qquad
 J^0_\rho(n)=\frac{(\ell^0_+(n)-\ell^0_-(n))R_\rho(n)}{2\mathfrak b_+(n)},\\
 (T_{\rm ref})_{\sigma\rho}=(V_c)_yK^0_\sigma(V_c)R_\rho(V_c)
   +s_\sigma(V_c)\pa_y[\mathfrak b_\sigma(V_c)]J^0_\rho(V_c),\qquad
 \widetilde{\mathcal M}=\mathcal M+T_{\rm ref},
\end{gathered}
\end{equation}
and
\begin{equation}\label{e:pB}
\begin{gathered}
 M_{\rm ref}=\bigl(\pa_{\zeta_\rho}\mathfrak j^\sharp_\sigma
   (V_c,V_h,\mathfrak a^a)\bigr)_{\sigma\rho},\qquad
 B_{\rm ref}=A_{\rm ph}-S^{-1}(M_{\rm ref}+T_{\rm ref})S,\\
 B_D=B_{\rm ref}-S^{-1}\operatorname{diag}(\mathcal M-M_{\rm ref})S,\qquad
 \mathcal K=\mathcal K_{\widetilde{\mathcal M}},\\
 B=B_D-S^{-1}\operatorname{offdiag}(\mathcal M-M_{\rm ref})S
   =A_{\rm ph}-S^{-1}\widetilde{\mathcal M}S,
\end{gathered}
\end{equation}
where $\operatorname{diag}$ and $\operatorname{offdiag}$ keep the diagonal and
the off-diagonal entries; $B_D$ keeps both off-diagonal entries of $M_{\rm ref}$
and all of $T_{\rm ref}$. By \eqref{e:gop},
$\mathcal K(S\mathbf v)=S[\mathbf v_t-d(V_c)\mathbf v_{yy}+B\mathbf v_y]$.

\begin{lemma}\label{l:am}
For every sufficiently late fixed $r$ the following holds on $\R\times[r,\infty)$.
The points $U^\sharp(V_c,\chi_m,\mathfrak a^a+\theta(\mathfrak a-\mathfrak a^a))$,
$0\le\theta\le1$, lie in a fixed neighbourhood of the curve $U_*$, and
$U^\sharp(V_c,\chi_m,\mathfrak a)=\mathcal U$. The matrices $\mathcal M$ and
$\widetilde{\mathcal M}$ have the form \eqref{e:gmat} with residuals
$\mathcal R$ and $\widetilde{\mathcal R}=\mathcal R+T_{\rm ref}$, and for
$\mathcal R'\in\{\mathcal R,\widetilde{\mathcal R}\}$
\begin{equation}\label{e:pmat}
 |\mathcal R'|\le C_M|R_U|+C_b\bigl(\tau^{-1}+\tau^{-\frac12}G_0\bigr),\qquad
 \sup_{t\ge r}\tau\Vert\widetilde{\mathcal M}(t)\Vert_{L^\infty}^2<\infty,\qquad
 \max_\sigma\Vert T_{L_\sigma,r}|\mathcal R'|^2\Vert_{L^\infty}\le\omega_b(r).
\end{equation}
Moreover
\begin{equation}\label{e:pT}
 |T_{\rm ref}|\le C\delta\tau^{-\frac12}e^{-\frac{c_G(y+h_*)^2}{2\tau}},
\end{equation}
and, with $z=\mathcal U-U_*(V_c)$ and the polynomials
$\mathcal P=(\mathcal P_+,\mathcal P_-)^\top$ of Lemma~\ref{l:nf},
\begin{equation}\label{e:pTo}
 \Bigl|(V_c)_yK^0_\sigma(V_c)z
   +s_\sigma(V_c)\pa_y[\mathfrak b_\sigma(V_c)]\,d_0\mathcal P(V_c,z)
   -\sum_\rho(T_{\rm ref})_{\sigma\rho}\mathfrak a_\rho\Bigr|
 \le C|(V_c)_y||z|^2 .
\end{equation}
\end{lemma}
\begin{proof}
By \eqref{e:pot}, $\chi_m-\frac{\mathcal F}{\gamma p_*}=V_b$,
$\mathfrak a_+-\mathfrak a_-=\frac{\mathcal F}{s_+(V_c)V_c}$ and
$\mathfrak a_++\mathfrak a_-=b_m$; hence
$U^\sharp(V_c,\chi_m,\mathfrak a)=(V_b,b_m,E_m)=\mathcal U$. In \eqref{e:pj},
$\frac{\zeta_+-\zeta_-}{s_+(n_0)}=\frac{\zeta_\sigma-\zeta_o}{s_\sigma(n_0)}$ and
$p(U^\sharp)-p_*=\bigl[s_\sigma(n_0)n_0(\zeta_\sigma-\zeta_o)
-\frac{\gamma-1}2(\zeta_\sigma+\zeta_o)^2\bigr]/v^\sharp$, so the three terms of
$\mathfrak j_\sigma$ have the quadratic parts
\[
 \frac{s_\sigma(n_0)}{2n_0}(n-n_0)(\zeta_\sigma-\zeta_o)
   -\frac{(\zeta_\sigma-\zeta_o)^2}{2n_0},\qquad
 \frac{(\gamma-1)(\zeta_\sigma+\zeta_o)^2}{4n_0},\qquad
 -\frac{(\gamma-1)(\zeta_\sigma^2-\zeta_o^2)}{2n_0},
\]
and, with $k=k(n_0)$ and $\ell_{\rm g}=\ell_{\rm g}(n_0)$,
\begin{equation}\label{e:jsq}
 \mathfrak j^\sharp_\sigma(n_0,n,\zeta)
 =-k\zeta_\sigma^2+2k\zeta_\sigma\zeta_o+\ell_{\rm g}\zeta_o^2
  +\frac{s_\sigma(n_0)}{2n_0}(n-n_0)(\zeta_\sigma-\zeta_o)
  +\mathfrak e_\sigma(n_0,n,\zeta),
\end{equation}
where $|\mathfrak e_\sigma|\le C(|n-n_0|+|\zeta|)^3$ and
$|D_{(n,\zeta)}\mathfrak e_\sigma|\le C(|n-n_0|+|\zeta|)^2$ for
$|n-n_0|+|\zeta|$ small, uniformly for $n_0\in[\frac{v_-}2,2v_-]$.

With $\ell_0$ and $\ell_\sigma(n)$ of Section~\ref{s:poteq},
$\chi_m=\ell_0\mathcal U$ and $\mathfrak a_\sigma=\ell_\sigma(V_c)(\mathcal U-U_*(V_c))$.
Put $\mathcal U_a=\widehat{\mathcal U}(y+h_*,t)+\sum_j\vartheta_jr_j$, so that
$\mathcal U-\mathcal U_a=R_U$. Since $R\Theta_c=p_*V_c$,
$\widehat{\mathcal U}(y+h_*,t)-U_*(V_c)=(V_h-V_c)R_0+(0,H_h,0)^\top$, and
$\ell_0R_0=1$, $\ell_0r_j=0$, $\ell_\sigma(n)R_0=0$,
$\ell_\sigma(n)r_j=\frac{s_j+s_\sigma(n)}2$; hence
$\ell_0\mathcal U_a=V_h$, $\ell_\sigma(V_c)(\mathcal U_a-U_*(V_c))=\mathfrak a^a_\sigma$,
$n_\chi=\ell_0R_U$ and $\mathfrak a_\sigma-\mathfrak a^a_\sigma=\ell_\sigma(V_c)R_U$.
By \eqref{eq:pr-gau}, Lemma~\ref{l:pr} and Lemma~\ref{l:asy}, for $\tau\ge h_*^2$,
\[
\begin{gathered}
 |n_\chi|+|\mathfrak a-\mathfrak a^a|\le C|R_U|,\qquad
 \sqrt\tau\Vert R_U(t)\Vert_{L^\infty}\le C_b,\qquad
 \lim_{t\to\infty}\sqrt\tau\Vert R_U(t)\Vert_{L^\infty}=0,\\
 |V_h-V_c|+|H_h|\le C_b\tau^{-\frac12}e^{-\frac{c_Gy^2}{2\tau}},\qquad
 \Bigl|\mathfrak a^a_\sigma-a_\sigma-\frac{H_h}2\Bigr|
 \le C_M\delta\tau^{-1}\sum_{j=\pm1}e^{-\frac{z_j^2}{16d_j}} .
\end{gathered}
\]
For the last bound,
$\mathfrak a^a_\sigma-a_\sigma-\frac{H_h}2=\frac{s_\sigma(V_c)-s_\sigma}{2s_\sigma}a_\sigma
+\frac{s_\sigma(V_c)+s_o}{2s_o}a_o$ with $s_\sigma=s_\sigma(v_\sigma)$ and
$s_o=-s_\sigma(v_o)$. For $j=\pm1$, \eqref{eq:ct-gau} gives
$|V_c-v_j|\le C\delta e^{-\frac{c_Gy^2}\tau}$ if $jy\ge0$, and
$\frac{c_Gy^2}\tau+\frac{(y-s_jt)^2}{16d_j\tau}\ge\frac{c_Gs_j^2t^2}{(1+16c_Gd_j)\tau}$,
while $|y-s_jt|\ge|s_j|t$ if $jy<0$; with the bound for $a_j$ in Lemma~\ref{l:pr}
this gives the claim. Thus the segment
of \eqref{e:pj} stays in a fixed neighbourhood of $(V_c,V_c,0)$, and
differentiating \eqref{e:jsq} in $\zeta$ along it gives, with $k,\ell_{\rm g}$
evaluated at $V_c$,
\[
\begin{gathered}
 \mathcal R_{\sigma\sigma}=\frac{s_\sigma(V_c)}{2V_c}(\chi_m-V_c)+\mathfrak e',\qquad
 \mathcal R_{\sigma o}=(k+\ell_{\rm g})H_h-\frac{s_\sigma(V_c)}{2V_c}(\chi_m-V_c)
   +\mathfrak e'',\\
 |\mathfrak e'|+|\mathfrak e''|\le C|R_U|+C\bigl(|V_h-V_c|+|a_+|+|a_-|+|H_h|\bigr)^2
   +C_M\tau^{-1};
\end{gathered}
\]
the terms with $H_h$ cancel in $\mathcal R_{\sigma\sigma}$. Since
$\chi_m-V_c=n_\chi+V_h-V_c$, $|a_\pm|^2\le C_M\tau^{-1}$ and
$e^{-\frac{c_Gy^2}{2\tau}}\le e^{\frac{c_G}2}G_0$ for $\tau\ge h_*^2$ (as
$b_0\le\frac{c_G}4$), this gives the first bound of \eqref{e:pmat} for
$\mathcal R$ and $\tau|\mathcal M|^2\le C_b$. By Lemma~\ref{l:sq},
Lemma~\ref{l:pos}(i) and the third bound of \eqref{e:gauss} for $\tau^{-1}G_0^2$
(with $v=0$, $h=h_*$),
\[
 \max_\sigma\Vert T_{L_\sigma,r}|\mathcal R|^2\Vert_{L^\infty}
 \le C_M\max_\sigma\Vert T_{L_\sigma,r}|R_U|^2\Vert_{L^\infty}
   +C_b\bigl(\tau_r^{-1}+\tau_r^{-\frac12}\bigr)\le\omega_b(r).
\]
The rows $K^0_\sigma$, $J^0_\rho$ and the vectors $R_\rho$ are smooth and
bounded on $[\frac{v_-}2,2v_-]$, so $|T_{\rm ref}|\le C|(V_c)_y|$; with
\eqref{eq:pr-gau} this gives \eqref{e:pT} for $\tau\ge h_*^2$, and
$\widetilde{\mathcal R}$ satisfies \eqref{e:pmat} as well.

For \eqref{e:pTo}, recall from Section~\ref{s:poteq} that
$\Lambda_\sigma(n,U)=O_\sigma(n,U)-\ell^0_\sigma(n)[s_\sigma(n)\mathbb I-\mathsf A(U)]$,
$O_\sigma(n,U)=\ell_\sigma(n)\mathsf B(U)-d(n)\ell_\sigma(n)$, and
$\ell^0_\sigma(n)R_0=f_0'(n)=\frac c{2n}$ by Lemma~\ref{l:nf}(i). Along
$U_*(n)+\epsilon R_0=U_*(n+\epsilon)$ the pressure equals $p_*$ and $b=0$, so
$\mathsf F$ is constant there and $D_U\mathsf A(U_*(n))[R_0]R_0=0$. Only the
entries $\frac\nu v$, $\frac bv(\nu-\frac\kappa{c_v})$, $\frac\kappa{c_vv}$ of
$\mathsf B$ depend on $U$, so
$D_U\mathsf B(U_*(n))[R_0]R_0=(0,0,-\frac{\kappa p_*}{c_v(\gamma-1)n^2})^\top$ and,
with $c_p=\gamma c_v$ and $s_\sigma(n)^2n=\gamma p_*$,
\[
 D_UO_\sigma(n,U_*(n))[R_0]R_0=-\frac{\kappa p_*}{2c_vs_\sigma(n)n^3}
 =-\frac{s_\sigma(n)f_0'(n)}n,\qquad
 (\ell^0_\sigma)'(n)R_0=f_0''(n)=-\frac{f_0'(n)}n .
\]
Since $D_U\Lambda_\sigma[e]=D_UO_\sigma[e]+\ell^0_\sigma D_U\mathsf A[e]$, this gives
\[
 K^0_\sigma(n)R_0=0,\qquad (\ell^0_+(n)-\ell^0_-(n))R_0=0 .
\]
The rows $\ell_+(n),\ell_-(n),\ell_0$ are dual to $R_+(n),R_-(n),R_0$, and
$\ell_0U_*(n)=n$; therefore
$z=\sum_\rho R_\rho(V_c)\mathfrak a_\rho+R_0(\chi_m-V_c)$. By the two
cancellations, $K^0_\sigma(V_c)z=\sum_\rho K^0_\sigma(V_c)R_\rho(V_c)\mathfrak a_\rho$
and, since $\mathcal P_\sigma(n,z)=\ell^0_\sigma(n)z+P_{2,\sigma}(n,z)$ with
$P_{2,\sigma}$ quadratic in $z$,
$d_0\mathcal P(V_c,z)=\sum_\rho J^0_\rho(V_c)\mathfrak a_\rho+O(|z|^2)$. With
$|\pa_y\mathfrak b_\sigma(V_c)|\le C|(V_c)_y|$ this is \eqref{e:pTo}.
\end{proof}

By \eqref{e:pmat}, for $\mathcal M=\widetilde{\mathcal M}$ the last bound of
\eqref{e:pcnx} reads $|C^1|\le C_b\delta\tau^{-1}e^{-\frac{b_1(y+h_*)^2}\tau}$.

For $r\ge S_0$ put $\mathcal E_{\sigma,r}=E^{\rm base}_\sigma-W_{\sigma,r}$, with
$E^{\rm base}_\sigma$ of Lemma~\ref{l:all} and $W_{\sigma,r}$ of \eqref{e:ubC}.
Since $\mathfrak b_-=-\mathfrak b_+$, $2f_0(V_c)=c\log\frac{V_c}{v_-}$,
$K_++K_-=B_L$ and $\bar A_++\bar A_-=A_++A_--B_L$, \eqref{e:pot} and
\eqref{e:plat} give
$\sum_\sigma(\mathcal E_{\sigma,r}+W_{\sigma,r})=\mathcal B-c\log\frac{V_c}{v_-}-A_+-A_-$;
with the identity $P_b-p_*=\frac{\nu(b_m)_y}{V_b}-\mathcal B_t$ this gives
\begin{equation}\label{e:bpot}
\begin{gathered}
 P_b-p_*=\frac{\nu(b_m)_y}{V_b}-\pa_t(\mathcal E_{+,r}+\mathcal E_{-,r})-j_r
   -c\,\pa_t\log V_c(y+h_*,t)-\pa_t(A_++A_-),\\
 j_r=\pa_t(W_{+,r}+W_{-,r})+c\,\pa_t\log\frac{V_c(y,t)}{V_c(y+h_*,t)} .
\end{gathered}
\end{equation}

Put $\mathcal E_r=(\mathcal E_{+,r},\mathcal E_{-,r})^\top$,
$u_r=(u_{+,r},u_{-,r})^\top$, $W_r=(W_{+,r},W_{-,r})^\top$,
$E^{\rm base}=(E^{\rm base}_+,E^{\rm base}_-)^\top$, $z=\mathcal U-U_*(V_c)$,
\begin{equation}\label{e:pn}
\begin{gathered}
 D^{\rm n}_{\sigma\rho}=\int_0^1\!\!\int_0^1\pa_n\pa_{\zeta_\rho}\mathfrak j^\sharp_\sigma
   (V_c,V_h+\theta n_\chi,\lambda\mathfrak a^a)\dd\lambda\dd\theta,\qquad
 C^{\rm n}_{\sigma j}=\frac12\sum_\rho\Bigl(1+\frac{s_\rho(V_c)}{s_j}\Bigr)
   D^{\rm n}_{\sigma\rho},\\
 f^{\rm n}_\sigma=\sum_ja_jC^{\rm n}_{\sigma j}n_\chi,\qquad
 f^{\rm o}_\sigma=\ell_{\rm g}(V_c)a_o^2,
\end{gathered}
\end{equation}
and, with $R_{\rm c}$ of \eqref{e:nfs} for the index $\sigma$ and
$\mathfrak b_\sigma=\mathfrak b_\sigma(V_c)$,
\begin{equation}\label{e:pfr}
\begin{gathered}
 Q_P=d_0(W_r-u_r),\qquad Q_A=d_0(\bar A_+,\bar A_-)^\top,\\
\begin{aligned}
 \mathfrak f_\sigma={}&\mathfrak j_\sigma(V_c,\mathcal U)+k(V_c)a_\sigma^2
   -R^{\rm prof}_\sigma+R_{\rm c}
   -\sum_\rho\widetilde{\mathcal M}_{\sigma\rho}\bigl(\mathfrak a_\rho-\pa_yf_0(V_c)
   -(W_{\rho,r})_y-a_\rho\bigr)\\
 &+\Bigl[L_\sigma\mathfrak b_\sigma-\sum_\rho\widetilde{\mathcal M}_{\sigma\rho}
   (\mathfrak b_\rho)_y\Bigr](Q_P+Q_A)
   -2d(V_c)(\mathfrak b_\sigma)_y\,\pa_y(Q_P+Q_A).
\end{aligned}
\end{gathered}
\end{equation}

\begin{lemma}\label{l:peq}
For every sufficiently late fixed $r$,
\begin{equation}\label{e:pceq}
 \mathcal K\mathcal E_r=f^{\rm n}+f^{\rm o}+f^{\rm s}+g
   +\mathcal C_{\widetilde{\mathcal M}}u_r,\qquad
 \mathcal E_r(r)=E^{\rm base}(r)-\mathcal P(V_c,z)(r),
\end{equation}
where $f^{\rm n}+f^{\rm o}+f^{\rm s}+g=\mathfrak f$,
$f^{\rm s}_\sigma(y,t)=\tau^{-1}F_\sigma(\frac y{\sqrt\tau})$ with a smooth
function $F_\sigma$ determined by the profile, $|F_\sigma^{(j)}(\xi)|\le
C\delta e^{-c_G\xi^2}$ for $0\le j\le3$, and
\[
 g=\widetilde{\mathcal M}(W_r)_y+K_cz_y+g_c,\qquad
 |K_c|\le C_b\tau^{-\frac12}e^{-\frac{c_G(y+h_*)^2}{2\tau}},\qquad
 |g_c|\le C_b\tau^{-\frac32}(G^*_0+G^*_++G^*_-),
\]
with $G^*_0=e^{-\frac{b_*(y+h_*)^2}\tau}$, $G^*_j=e^{-b_*z_j^2}$ and
$b_*=\min\{\frac{c_G}8,\frac1{32d_-},\frac1{32d_+}\}$.
Moreover $|C^{\rm n}|\le C$, $\Vert g_c(t)\Vert_{L^2}\le C_b\tau^{-\frac54}$ and
$\int_r^\infty(\Vert(W_r)_y\Vert_{L^2}^2+\Vert z_y\Vert_{L^2}^2)\dd t<\infty$.
\end{lemma}
\begin{proof}
Write $d=d(V_c)$ and $\Delta\mathcal G=\mathcal G-G_c$. By Lemma~\ref{l:Deq} with
$\zeta\equiv1$ and $W=W_{\sigma,r}$, where $D_\sigma=\mathcal E_{\sigma,r}+\bar A_\sigma$,
by $L_\sigma\bar A_\sigma=L_\sigma A_\sigma=R^{\rm prof}_\sigma-k(V_c)a_\sigma^2$
(Lemma~\ref{l:pr}), $(\Delta\mathcal G)_y=\mathcal F$, $(\bar A_\sigma)_y=a_\sigma$ and
$(\mathfrak b_\sigma)_y=\mathfrak b_\sigma'(V_c)(V_c)_y$,
\[
\begin{gathered}
 L_\sigma\mathcal E_{\sigma,r}=\mathfrak j_\sigma(V_c,\mathcal U)+k(V_c)a_\sigma^2
   -R^{\rm prof}_\sigma+R_{\rm c}+(L_\sigma\mathfrak b_\sigma)\Delta\mathcal G
   -2d(\mathfrak b_\sigma)_y(\Delta\mathcal G)_y,\\
 \mathfrak a_\sigma=(\mathcal E_{\sigma,r})_y+a_\sigma+(W_{\sigma,r})_y
   +\pa_yf_0(V_c)-(\mathfrak b_\sigma)_y\Delta\mathcal G .
\end{gathered}
\]
Since $\mathfrak b_-=-\mathfrak b_+$, $\Psi_+-\Psi_-=2\mathfrak b_+\mathcal G$ and
$K_+-K_-=2\mathfrak b_+G_c$, so that $\Delta\mathcal G=d_0\mathcal E_r+d_0u_r+Q_P+Q_A$.
Add and subtract $\sum_\rho\widetilde{\mathcal M}_{\sigma\rho}\mathfrak a_\rho$ in the
equation, insert the formula for $\mathfrak a_\rho$ into the added sum and that for
$\Delta\mathcal G$ into the last two terms. By \eqref{e:gop}, the terms containing
$\mathcal E_r$ are
$(\widetilde{\mathcal M}(\mathcal E_r)_y+\mathcal C_{\widetilde{\mathcal M}}\mathcal E_r)_\sigma$,
those containing $u_r$ are $(\mathcal C_{\widetilde{\mathcal M}}u_r)_\sigma$, and the
others are $\mathfrak f_\sigma$. This is the equation in \eqref{e:pceq}; the trace
follows from the definition of $\mathcal E_{\sigma,r}$ and $u_{\sigma,r}(r)=0$.

It remains to split $\mathfrak f$. Let $\tau_r\ge h_*^2$. By
Lemma~\ref{l:asy}, \eqref{eq:pr-gau} and Lemma~\ref{l:pr},
$|z|\le C_b\tau^{-\frac12}$. Let $g_0$ be $e^{-\frac{c_Gy^2}{2\tau}}$ or
$e^{-\frac{c_G(y+h_*)^2}{2\tau}}$ and $g_j=e^{-\frac{z_j^2}{16d_j}}$. By the argument
for the last bound in the proof of Lemma~\ref{l:am}, $g_0g_j+g_+g_-$ is at most
$C_be^{-c_bt}(G^*_0+G^*_++G^*_-)$ with $c_b>0$, and $g_0\le C_bG^*_0$; below,
$R^{\rm g}$ denotes any function bounded by
$C_b\tau^{-\frac32}(G^*_0+G^*_++G^*_-)$, and $\mathfrak f\equiv\mathfrak f'$ means
$\mathfrak f-\mathfrak f'=R^{\rm g}$.

At $\zeta=0$, $U^\sharp(n_0,n,0)=U_*(n)$, so that $p(U^\sharp)=p_*$, the second
component of $U^\sharp$ vanishes and $v^\sharp=n$; since
$\pa_{\zeta_\rho}p(U^\sharp)=\frac{\rho s_+(n_0)n_0}n$ at $\zeta=0$, the formula for
$\mathfrak j_\sigma$ in Section~\ref{s:poteq} gives
\[
 \mathfrak j^\sharp_\sigma(n_0,n,0)=0,\qquad
 \pa_{\zeta_\rho}\mathfrak j^\sharp_\sigma(n_0,n,0)=\frac{\rho s_+(n_0)(n-n_0)}{2n},
 \qquad
 \pa_n(\pa_{\zeta_+}+\pa_{\zeta_-})\mathfrak j^\sharp_\sigma(n_0,n,0)=0 .
\]
The first identity and \eqref{e:pj} give the two mean value formulas
\[
 \mathfrak j_\sigma(V_c,\mathcal U)=\sum_\rho\mathcal M_{\sigma\rho}
   (\mathfrak a_\rho-\mathfrak a^a_\rho)+\mathfrak j^\sharp_\sigma(V_c,\chi_m,\mathfrak a^a),
 \qquad
 \mathfrak j^\sharp_\sigma(V_c,\chi_m,\mathfrak a^a)-\mathfrak j^\sharp_\sigma(V_c,V_h,\mathfrak a^a)
 =n_\chi\sum_\rho D^{\rm n}_{\sigma\rho}\mathfrak a^a_\rho,
\]
and $n_\chi\sum_\rho D^{\rm n}_{\sigma\rho}\mathfrak a^a_\rho
=f^{\rm n}_\sigma+\frac12H_hn_\chi\sum_\rho D^{\rm n}_{\sigma\rho}$. The segments in
\eqref{e:pn} lie in the neighbourhood of Lemma~\ref{l:am}, so $|D^{\rm n}|+|C^{\rm n}|\le C$,
and the third identity gives $|\sum_\rho D^{\rm n}_{\sigma\rho}|\le C|\mathfrak a^a|$.
With $|n_\chi|\le C|R_U|\le C_b\tau^{-\frac12}$ and
$|\mathfrak a^a|\le C_M\tau^{-\frac12}$ this gives $H_hn_\chi\sum_\rho D^{\rm n}_{\sigma\rho}=R^{\rm g}$.

By \eqref{e:jsq} with $n_0=V_c$, $n=V_h$, $\zeta=\mathfrak a^a$, the bounds in the
proof of Lemma~\ref{l:am}, $k(n)+\ell_{\rm g}(n)=\frac{\gamma-1}n$ and
$|H_h^2-H_c^2|\le C_b\tau^{-\frac32}e^{-\frac{c_Gy^2}{2\tau}}$,
\[
 \mathfrak j^\sharp_\sigma(V_c,V_h,\mathfrak a^a)+k(V_c)a_\sigma^2
 \equiv\ell_{\rm g}(V_c)a_o^2+\frac{(\gamma-1)H_c^2}{4V_c};
\]
the products $a_\sigma a_o$, $H_ha_o$ and $(V_h-V_c)(a_\sigma-a_o)$ have different
centres, and the cubic remainder is $R^{\rm g}$.

For $R_{\rm c}$, write $Q=Q_\sigma$, $F_2=F_{2,\sigma}$, $\mathsf B=\mathsf B(\mathcal U)$,
evaluate $Q$, $F_2$, $\mathcal R_\sigma$ and their derivatives at $(V_c,\mathcal U)$
with the $n$-derivatives at fixed $U$, and use $\mathcal U_y=z_y+R_0(V_c)_y$. Then
\eqref{e:nfs}, $(F_{\rm nf})_n=(F_2)_n-Q_{Un}[(\mathsf B-d\mathbb I)\mathcal U_y]
+d'(V_c)Q_U\mathcal U_y$ and
$L_\sigma V_c=(V_c)_t-d(V_c)_{yy}+s_\sigma(V_c)(V_c)_y$ give
\[
\begin{gathered}
 R_{\rm c}=A^{\rm c}z_y+dQ_n(V_c)_{yy}-Q_n(V_c)_t
   +(V_c)_y\bigl[\mathcal R_\sigma(V_c,\mathcal U)R_0-(F_2)_n-s_\sigma(V_c)Q_n\bigr]
   +(V_c)_y^2D^{\rm c}+\pa_yF_{\rm c},\\
 A^{\rm c}e=(V_c)_y\bigl\{Q_{UU}[e,\mathsf BR_0]+Q_{UU}[R_0,\mathsf Be]
   +Q_{Un}[(\mathsf B+d\mathbb I)e]-d'(V_c)Q_Ue\bigr\},\\
 D^{\rm c}=Q_{UU}[R_0,\mathsf BR_0]+Q_{Un}[(\mathsf B+d\mathbb I)R_0]
   -d'(V_c)Q_UR_0+dQ_{nn}.
\end{gathered}
\]
Let $H_P(n)$, $H_F(n)$ be the Hessians of $P_{2,\sigma}(n,\cdot)$, $F_{2,\sigma}(n,\cdot)$.
Since $\ell^0_\sigma R_0=f_0'$ (Lemma~\ref{l:nf}(i)),
$Q_n=(\ell^0_\sigma)'z-H_P[R_0,z]+O(|z|^2)$ and $-(F_2)_n=H_F[R_0,z]+O(|z|^2)$;
applied to $R_0$, with $\mathsf A(U_*)R_0=0$, Lemma~\ref{l:nf}(i) gives
$D_U\Lambda_\sigma(n,U_*)[z]R_0=s_\sigma(n)H_P[z,R_0]+H_F[z,R_0]$. Thus the part of
$(V_c)_y[-(F_2)_n-s_\sigma(V_c)Q_n]$ linear in $z$ is $(V_c)_yK^0_\sigma(V_c)z$, and
$|\mathcal R_\sigma(V_c,\mathcal U)|\le C|z|^2$. Put
\[
 D_{00,\sigma}(n)=D^{\rm c}\big|_{(n,U_*(n))},\qquad
 J_{0,\sigma}(n)=\ell^0_\sigma(n)\bigl(\mathsf B(U_*(n))-d(n)\mathbb I\bigr)R_0 ;
\]
since $Q_U(n,U_*(n))=\ell^0_\sigma(n)$, $-J_{0,\sigma}(V_c)(V_c)_y$ is the value of
$F_{\rm c}$ at $z=0$, $\widetilde z=-(0,H_c,\frac{H_c^2}2)$, apart from $F_2$. The
remaining terms are products of $z,z_y,\widetilde z_y$ with derivatives of $V_c$ and
$H_c$; by \eqref{eq:pr-gau}, \eqref{e:pt} and $|z|\le C_b\tau^{-\frac12}$,
\[
 R_{\rm c}=(V_c)_yK^0_\sigma(V_c)z+D_{00,\sigma}(V_c)(V_c)_y^2
   -\pa_y\bigl(J_{0,\sigma}(V_c)(V_c)_y\bigr)+K^{\rm R}z_y+R^{\rm g},
 \qquad |K^{\rm R}|\le C\bigl(|(V_c)_y|+|H_c|\bigr).
\]

Since $\pa_yf_0(V_c)=\frac{H_c}2$ and $\widetilde{\mathcal M}=\mathcal M+T_{\rm ref}$,
\[
\begin{aligned}
 &-\sum_\rho\widetilde{\mathcal M}_{\sigma\rho}\bigl(\mathfrak a_\rho-\pa_yf_0(V_c)
   -(W_{\rho,r})_y-a_\rho\bigr)\\
 &\quad=-\sum_\rho\mathcal M_{\sigma\rho}\Bigl(\mathfrak a_\rho-a_\rho-\frac{H_c}2\Bigr)
   -\sum_\rho(T_{\rm ref})_{\sigma\rho}\Bigl(\mathfrak a_\rho-a_\rho-\frac{H_c}2\Bigr)
   +(\widetilde{\mathcal M}(W_r)_y)_\sigma .
\end{aligned}
\]
In the first sum write $\mathfrak a_\rho=(\mathfrak a_\rho-\mathfrak a^a_\rho)+\mathfrak a^a_\rho$;
the part with $\mathfrak a_\rho-\mathfrak a^a_\rho$ cancels the first mean value
formula, and $|\mathcal M||\mathfrak a^a_\rho-a_\rho-\frac{H_c}2|=R^{\rm g}$ by
$\tau|\mathcal M|^2\le C_b$, the bound for $\mathfrak a^a-a-\frac{H_h}2$ in the proof of
Lemma~\ref{l:am} and $|H_h-H_c|\le C_b\tau^{-1}e^{-\frac{c_Gy^2}{2\tau}}$. In the second sum,
$T_{\rm ref}a=R^{\rm g}$ by \eqref{e:pT}. In the bracket term of \eqref{e:pfr},
$L_\sigma\mathfrak b_\sigma=s_\sigma(V_c)(\mathfrak b_\sigma)_y+(\mathfrak b_\sigma)_t
-d(\mathfrak b_\sigma)_{yy}$ and $Q_P=d_0\mathcal P(V_c,z)$; with $|Q_P|\le C|z|$,
\[
 \Bigl[L_\sigma\mathfrak b_\sigma-\sum_\rho\widetilde{\mathcal M}_{\sigma\rho}
   (\mathfrak b_\rho)_y\Bigr]Q_P-2d(\mathfrak b_\sigma)_y\pa_yQ_P
 \equiv s_\sigma(V_c)(\mathfrak b_\sigma)_yd_0\mathcal P(V_c,z)
   -2d(\mathfrak b_\sigma)_yd_0D_z\mathcal P(V_c,z)z_y .
\]
By \eqref{e:pTo},
$(V_c)_yK^0_\sigma(V_c)z+s_\sigma(V_c)(\mathfrak b_\sigma)_yd_0\mathcal P(V_c,z)
-\sum_\rho(T_{\rm ref})_{\sigma\rho}\mathfrak a_\rho=R^{\rm g}$. By
\eqref{e:CH} and $F_d(\xi)\le e^{-\frac{\xi^2}{4d}}$ for $\xi\le0$,
$|\bar A_+|\le C_Me^{-\frac{z_+^2}{4d_+}}$ for $y\le s_+t$ and
$|\bar A_-|\le C_Me^{-\frac{z_-^2}{4d_-}}$ for $y\ge s_-t$, while the coefficients of
$Q_A$ and $(Q_A)_y$ in $\mathfrak f$ carry a factor $(V_c)_y$; so the terms with $Q_A$
are $R^{\rm g}$. So are $R^{\rm prof}_\sigma$, by its formula in Lemma~\ref{l:pr} and
the argument for $\mathfrak a^a$, and $H_hn_\chi\sum_\rho D^{\rm n}_{\sigma\rho}$. Collecting,
$\mathfrak f=f^{\rm n}+f^{\rm o}+f^{\rm s}+g$ with
\[
\begin{gathered}
 f^{\rm s}_\sigma=\frac{(\gamma-1)H_c^2}{4V_c}+D_{00,\sigma}(V_c)(V_c)_y^2
   -\pa_y\bigl(J_{0,\sigma}(V_c)(V_c)_y\bigr)
   +\frac{H_c}2\sum_\rho(T_{\rm ref})_{\sigma\rho},\\
 (K_c)_\sigma=K^{\rm R}-2d(\mathfrak b_\sigma)_yd_0D_z\mathcal P(V_c,z),
\end{gathered}
\]
and $g_c$ the sum of the terms denoted $R^{\rm g}$. Since $V_c$, $H_c$ and
$T_{\rm ref}$ are functions of $\frac y{\sqrt\tau}$ times $1$, $\tau^{-\frac12}$
and $\tau^{-\frac12}$, $f^{\rm s}_\sigma=\tau^{-1}F_\sigma(\frac y{\sqrt\tau})$, and
\eqref{eq:pr-gau}, applied to at most five derivatives of the profile, gives the
bound for $F_\sigma^{(j)}$. The bound for $K_c$ follows from \eqref{eq:pr-gau}
and $\tau\ge h_*^2$, and $\Vert R^{\rm g}\Vert_{L^2}\le C_b\tau^{-\frac54}$.
Finally, $z=(V_b-V_c,b_m,c_v(\Theta_b-\Theta_c)+\frac{b_m^2}2)$, so
Lemma~\ref{l:dep} and \eqref{eq:pr-gau} give $\int_r^\infty\Vert z_y\Vert_{L^2}^2\dd t<\infty$;
with $|\pa_y\mathcal P(V_c,z)|\le C(|z_y|+|z||(V_c)_y|)$ and \eqref{e:nfe} this gives
the integral of $\Vert(W_r)_y\Vert_{L^2}^2$.
\end{proof}

\subsection{Weighted scalar estimates}

Let $P=\pa_t-d\pa_y^2+B\pa_y$ with $d,B$ as in \eqref{e:Tdef}, let
$P^*w=-w_t-(dw)_{yy}-(Bw)_y$, and let $w>0$ and $\Upsilon\ge0$ be bounded and
smooth enough for the integrations by parts below. With
$\omega=w(1+\Upsilon)^{-\frac12}$ and $z_\Upsilon=\frac{\Upsilon_y}{1+\Upsilon}$,
\begin{equation}\label{e:rcw}
 \frac12\frac{\dd}{\dd t}\int_\R\omega u^2\dd y
 +\int_\R d\omega\Bigl(u_y^2-z_\Upsilon uu_y+\frac38z_\Upsilon^2u^2\Bigr)\dd y
 +\int_\R\Bigl[\frac{w\,P\Upsilon}{4(1+\Upsilon)^{\frac32}}
   +\frac{P^*w}{2(1+\Upsilon)^{\frac12}}\Bigr]u^2\dd y
 =\int_\R\omega u\,Pu\dd y ,
\end{equation}
and $x^2-xz+\frac38z^2\ge\frac{x^2}4+\frac{z^2}{24}$. Indeed, two integrations by
parts give $\frac12\frac{\dd}{\dd t}\int\omega u^2+\int d\omega u_y^2
+\frac12\int(P^*\omega)u^2=\int\omega uPu$, and, with
$\phi(\Upsilon)=(1+\Upsilon)^{-\frac12}$,
\[
 P^*(w\phi)=\phi P^*w-\phi'wP\Upsilon-2\phi'(dw\Upsilon_y)_y-\phi''dw\Upsilon_y^2 ;
\]
one more integration by parts in the term with $(dw\Upsilon_y)_y$ gives
\eqref{e:rcw}. The inequality is
$(\frac{\sqrt3}2x-\frac z{\sqrt3})^2\ge0$.

Put $b_{\rm ac}=\frac1{16d(\rho_-^{-1})}$. Then $b_{\rm ac}\le\frac1{8d_\pm}$, and
by Lemma~\ref{l:pr}
\[
 |a_j|^2+\tau^{\frac12}|(a_j)_y|\le C_M\tau^{-1}e^{-2b_{\rm ac}z_j^2}\qquad(j=\pm1).
\]
For $\sigma=\pm1$ and $b>0$ put, with $\lambda(u)=u+4d_\sigma b(1-u)$,
\begin{equation}\label{e:qresp}
\begin{gathered}
 \Phi_b(\xi)=\int_0^1\frac{e^{-\frac{b\xi^2}{\lambda(u)}}}{\sqrt{u\lambda(u)}}\dd u,
 \qquad
 P_\sigma\bigl[\Phi_b(z_\sigma)\bigr]=\tau^{-1}e^{-bz_\sigma^2},\\
 \bigl\Vert T_{\mathscr L,r}\bigl(\tau^{-1}e^{-bz_\sigma^2}\bigr)\bigr\Vert_{L^\infty}\le C
 \quad(\mathscr L\in\{L_\sigma,\mathscr P_\sigma\}),
\end{gathered}
\end{equation}
with $C$ depending only on $b$, $M$ and the fixed data. Indeed, in
$x=y-s_\sigma t$, $P_\sigma=\pa_t-d_\sigma\pa_x^2$ and
$P_\sigma[\Phi(\frac x{\sqrt\tau})]=\tau^{-1}(\mathfrak L\Phi)(z_\sigma)$ with
$\mathfrak L\Phi=-d_\sigma\Phi''-\frac\xi2\Phi'$. For
$\phi_\lambda(\xi)=\lambda^{-\frac12}e^{-\frac{b\xi^2}\lambda}$,
$\mathfrak L\phi_\lambda=\phi_\lambda[\frac{2bd_\sigma}\lambda
+\frac{b\xi^2}{\lambda^2}(\lambda-4bd_\sigma)]$, and $\lambda-4bd_\sigma=u\lambda'$;
hence $\pa_u[u^{\frac12}\phi_{\lambda(u)}]=u^{-\frac12}\mathfrak L\phi_{\lambda(u)}$ and
$\mathfrak L\Phi_b=\phi_1=e^{-b\xi^2}$. Moreover
$0\le\Phi_b\le2\max\{1,(4d_\sigma b)^{-\frac12}\}$ and
$|\Phi_b'|+|\Phi_b''|\le C(1+\xi^2)e^{-b\xi^2/\max\{1,4d_\sigma b\}}$. By
\eqref{eq:ct-gau}, $|V_c-v_\sigma|\le C\delta e^{-\frac{c_Gy^2}\tau}$ for
$\sigma y\ge0$, while $|y-s_\sigma t|\ge|s_\sigma|t$ for $\sigma y<0$; therefore
\[
 |(L_\sigma-P_\sigma)\Phi_b(z_\sigma)|
 =\Bigl|-(d(V_c)-d_\sigma)\frac{\Phi_b''(z_\sigma)}\tau
   +(s_\sigma(V_c)-s_\sigma)\frac{\Phi_b'(z_\sigma)}{\sqrt\tau}\Bigr|
 \le C\delta\tau^{-\frac12}e^{-c_bt}
\]
with $c_b>0$ depending only on $b$ and the fixed data. The function
$\Phi_b(z_\sigma)+C\delta\int_r^t(1+t')^{-\frac12}e^{-c_bt'}\dd t'$ is a bounded
supersolution, and comparison gives the bound for $L_\sigma$;
Lemma~\ref{l:pos}(v) gives the one for $\mathscr P_\sigma$.

For $a>0$ write $W^a=\tau^{-\frac12}G^a$, so that $W_c=W^{2b_0}$, and let
$\mathcal M^{\rm d}$ denote the diagonal part of a matrix $\mathcal M$.

\begin{lemma}\label{l:geo}
Let $0<a\le\frac{c_G}2$ and let $\mathcal M$ be a real measurable matrix on
$\R\times[S_0,\infty)$ of the form \eqref{e:gmat} with
$\sup\tau|\mathcal R|^2<\infty$ and
$\max_\sigma\Vert T_{L_\sigma,r}|\mathcal R_{\sigma\cdot}|^2\Vert_{L^\infty}\to0$ as
$r\to\infty$, where $|\mathcal R_{\sigma\cdot}|^2=\sum_\rho|\mathcal R_{\sigma\rho}|^2$.
For every sufficiently late fixed $r$ the following holds, with $C$ depending only
on $M$, $a$ and the fixed data. Let $\mathfrak v_\sigma\ge0$ be measurable with
$\Vert T_{L_\sigma,r}\mathfrak v_\sigma\Vert_{L^\infty}\le1$ $(\sigma=\pm1)$, let
$g\in H^1$, let $r<T\le\infty$, and let $F=F_1+f$ with $F\in L^2((r,T')\times\R)$
for $T'<T$ finite and
\[
 \mathcal B_F=\Bigl(\int_r^T\Vert F_1\Vert_{L^2}\dd t\Bigr)^2
   +\sum_{\sigma=\pm1}\int_r^T\!\!\int_\R\frac{|f_\sigma|^2}{W^a+\mathfrak v_\sigma}\dd y\dd t
 <\infty .
\]
Then $\mathcal K_{\mathcal M}\mathbf H=F$, $\mathbf H(r)=g$, has a unique solution
on $(r,T)$ with $\mathbf H\in C([r,T'];H^1)$ and
$\mathbf H_t,\mathbf H_{yy}\in L^2((r,T')\times\R)$ for $T'<T$ finite, and
\[
 \sup_{r\le t<T}\Vert\mathbf H(t)\Vert_{L^2}^2
 +\int_r^T\Bigl(\Vert\mathbf H_y\Vert_{L^2}^2
   +\sum_{\sigma=\pm1}\int_\R(W^a+\mathfrak v_\sigma)|\mathbf H_\sigma|^2\dd y\Bigr)\dd t
 \le C\bigl(\Vert g\Vert_{L^2}^2+\mathcal B_F\bigr).
\]
The same holds with $\mathcal K_{\mathcal M}$ replaced by $L-\mathcal M\pa_y$; by
$L-\mathcal M^{\rm d}\pa_y$, where only the diagonal entries of $\mathcal M$ need the
form \eqref{e:gmat} and the bounds for $\mathcal R_{\sigma\sigma}$; and by $L$, for every
$r$ for which \eqref{e:nfwdef} applies and with $C$ depending only on $a$ and the fixed
data.
\end{lemma}
\begin{proof}
Write $d=d(V_c)$, $k=k(V_c)$, $\ell_{\rm g}=\ell_{\rm g}(V_c)$, and let
$\tau_r\ge h_*^2$ with \eqref{e:nfwdef} applicable. Put
$Q_\sigma=\tau^{-1}e^{-b_{\rm ac}z_\sigma^2}$,
$W^{\rm o}_\sigma=\tau^{-\frac12}e^{-b_{\rm ac}z_o^2}$,
$\eta_r=\max_\sigma\Vert T_{L_\sigma,r}|\mathcal R_{\sigma\cdot}|^2\Vert_{L^\infty}+\tau_r^{-\frac12}$
and
\[
 \Upsilon_\sigma=T_{L_\sigma,r}\Bigl(\frac{|\mathcal R_{\sigma\cdot}|^2}{\eta_r}+Q_\sigma
   +W^{\rm o}_\sigma+W^a+\mathfrak v_\sigma\Bigr).
\]
By \eqref{e:qresp}, and by \eqref{e:gauss} with $(v,h)=(s_o,0)$ and $(0,h_*)$, for
which $\sigma(s_\sigma(V_c)-v)\ge\min\{c_-,c_+\}$,
$0\le\Upsilon_\sigma\le C_\Upsilon$ with $C_\Upsilon$ depending only on $M$, $a$ and
the fixed data. Let $p_\sigma=e^{-\alpha A_\sigma}$ and $h_\sigma=p_\sigma\mathbf H_\sigma$;
by Lemma~\ref{l:pr}, $e^{-C_M}\le p_\sigma\le e^{C_M}$, $(p_\sigma)_y=-\alpha a_\sigma p_\sigma$
and $L_\sigma(p_\sigma z)=p_\sigma(\mathscr P_\sigma z-\alpha R^{\rm prof}_\sigma z)$. With
$B_{\sigma\rho}=\mathcal M_{\sigma\rho}+2ka_\sigma\mathbf 1_{\rho=\sigma}$, the equation
$\mathcal K_{\mathcal M}\mathbf H=F$ becomes
\[
 L_\sigma h_\sigma=\sum_\rho\frac{p_\sigma}{p_\rho}B_{\sigma\rho}
   \bigl((h_\rho)_y+\alpha a_\rho h_\rho\bigr)
   +p_\sigma(\mathcal C_{\mathcal M}\mathbf H+F)_\sigma-\alpha R^{\rm prof}_\sigma h_\sigma,
\]
and by \eqref{e:gmat}, $B_{\sigma\sigma}=2ka_o+\mathcal R_{\sigma\sigma}$,
$B_{\sigma o}=2ka_\sigma+2\ell_{\rm g}a_o+\mathcal R_{\sigma o}$.

Apply \eqref{e:rcw} with $P=L_\sigma$, $w=w_0$ of \eqref{e:nfwdef} for the index
$\sigma$, $\Upsilon=\Upsilon_\sigma$ and $u=h_\sigma$. By \eqref{e:nfwts},
$P^*w_0\ge c_1g_r$, and $L_\sigma\Upsilon_\sigma$ is the source of $\Upsilon_\sigma$.
Summing over $\sigma$ gives, for $E=\frac12\sum_\sigma\int\omega_\sigma h_\sigma^2$
and with $c_3>0$ depending only on $C_\Upsilon$ and the fixed data,
\[
\begin{gathered}
 E'+c_3\mathcal D\le\sum_\sigma\int\omega_\sigma h_\sigma L_\sigma h_\sigma,\\
 \mathcal D=\sum_\sigma\int\Bigl[(h_\sigma)_y^2+\Bigl(z_{\Upsilon_\sigma}^2+g_r
   +\frac{|\mathcal R_{\sigma\cdot}|^2}{\eta_r}+Q_\sigma+W^{\rm o}_\sigma+W^a+\mathfrak v_\sigma\Bigr)
   h_\sigma^2\Bigr].
\end{gathered}
\]
In the right side, $|(\omega_\sigma)_y|\le C(|(w_0)_y|+|z_{\Upsilon_\sigma}|)$ with
$|(w_0)_y|^2\le Cg_r$, and $|a_\rho|^2\le C_MQ_\rho$. The terms are bounded as
follows, with $\varepsilon>0$ fixed small and $k_r\ge0$, $\int_r^\infty k_r\le C_M$.

The residual terms of row $\sigma$ are at most
$\varepsilon\mathcal D+C_\varepsilon C_M\eta_r\int\frac{|\mathcal R_{\sigma\cdot}|^2}{\eta_r}h_\sigma^2$.
The terms with $a_o$ in row $\sigma$ are at most
$\varepsilon\mathcal D+C_\varepsilon\int(a_o^2+|a_oa_\sigma|)h_\sigma^2
+C\int a_o^2|h_\sigma h_o|$; here $a_o^2\le C_M\tau^{-\frac12}W^{\rm o}_\sigma$,
$\Vert a_oa_\sigma(t)\Vert_{L^\infty}\le C_Me^{-c_et}$ and
$\tau^{-1}e^{-2b_{\rm ac}z_o^2}|h_\sigma h_o|\le\tau^{-\frac14}(W^{\rm o}_\sigma h_\sigma^2+Q_oh_o^2)$.
The entry $2ka_\sigma$ of $B_{\sigma o}$ is integrated by parts together with its
conjugation factor; since $(\frac{p_\sigma}{p_o})_y=\alpha(a_o-a_\sigma)\frac{p_\sigma}{p_o}$,
\[
\begin{aligned}
 \int\omega_\sigma h_\sigma\frac{p_\sigma}{p_o}2ka_\sigma\bigl((h_o)_y+\alpha a_oh_o\bigr)
 ={}&-\int\frac{p_\sigma}{p_o}2k\omega_\sigma a_\sigma(h_\sigma)_yh_o\\
 &-\int\frac{p_\sigma}{p_o}\bigl[2k(\omega_\sigma)_ya_\sigma+2k_y\omega_\sigma a_\sigma
   +2k\omega_\sigma(a_\sigma)_y-2\alpha k\omega_\sigma a_\sigma^2\bigr]h_\sigma h_o ,
\end{aligned}
\]
which is at most $\varepsilon\mathcal D+C_\varepsilon\int a_\sigma^2h_o^2
+C_M\int\tau^{-1}e^{-2b_{\rm ac}z_\sigma^2}|h_\sigma h_o|+k_rE$, because
$\Vert(V_c)_ya_\sigma(t)\Vert_{L^\infty}\le C\delta e^{-c_et}$; moreover
$a_\sigma^2\le C_M\tau^{-\frac12}W^{\rm o}_o$ and
$\tau^{-1}e^{-2b_{\rm ac}z_\sigma^2}|h_\sigma h_o|\le\tau^{-\frac14}(Q_\sigma h_\sigma^2+W^{\rm o}_oh_o^2)$.
The term with $R^{\rm prof}_\sigma$ is at most $k_rE$ by Lemma~\ref{l:pr}. By the
computation for \eqref{e:pcnx}, $|\mathcal C^0|+|B^1|\le C\delta W^a$ and
$|C^1|\le C\delta(\tau^{-\frac12}+|\mathcal M|)W^a$ for $a\le\frac{c_G}2$, and
$|\mathcal M|\le C_M\tau^{-\frac12}+|\mathcal R|$; with
$|\mathbf H_y|\le C_M(|h_y|+|a||h|)$ the connection term is at most
$C_M\delta\mathcal D+C_M\bigl(1+(\sup\tau|\mathcal R|^2)^{\frac12}\bigr)\tau_r^{-\frac12}\mathcal D$.
Finally
$\sum_\sigma\int\omega_\sigma h_\sigma p_\sigma F_\sigma
\le C\sqrt E\Vert F_1\Vert_{L^2}+\frac{c_3}4\mathcal D
+C\sum_\sigma\int\frac{|f_\sigma|^2}{W^a+\mathfrak v_\sigma}$.

Fix $\varepsilon$, then $\delta_1(M)$ with $C_M\delta_1(M)\le\frac{c_3}8$, then $r$ with
$C_\varepsilon C_M\eta_r+C_M\tau_r^{-\frac14}
+C_M(1+(\sup\tau|\mathcal R|^2)^{\frac12})\tau_r^{-\frac12}\le\frac{c_3}8$. Then
\[
 E'+\frac{c_3}2\mathcal D\le k_rE+C\sqrt E\Vert F_1\Vert_{L^2}
   +C\sum_\sigma\int\frac{|f_\sigma|^2}{W^a+\mathfrak v_\sigma},
\]
and the integrating factor $e^{-\int_r^tk_r}$ and the supremum of $\sqrt E$ give
$\sup E+\int_r^T\mathcal D\le C(\Vert g\Vert_{L^2}^2+\mathcal B_F)$. Since
$(\mathbf H_\sigma)_y=p_\sigma^{-1}((h_\sigma)_y+\alpha a_\sigma h_\sigma)$ and
$a_\sigma^2\le C_MQ_\sigma$, this is the estimate. For $L-\mathcal M\pa_y$ the
connection term is absent, for $L-\mathcal M^{\rm d}\pa_y$ also $B_{\sigma o}=0$ and
only $\mathcal R_{\sigma\sigma}$ enters, and for $L$ the test is done with
$\mathbf H_\sigma$ itself and
$\Upsilon_\sigma=T_{L_\sigma,r}(W^a+\mathfrak v_\sigma)$, without conjugation and
residual.

On $(r,T')$, $F\in L^2$, and Galerkin approximation with the tests $\mathbf v$ and
$-\mathbf v_{yy}$, where $\mathbf v=S^{-1}\mathbf H$ solves
$\mathbf v_t-d\mathbf v_{yy}+(A_{\rm ph}-S^{-1}\mathcal MS)\mathbf v_y=S^{-1}F$ by
\eqref{e:gop} (for the three other operators $\mathbf H$ is tested directly), gives
$\frac{\dd}{\dd t}\Vert\mathbf v\Vert_{H^1}^2+c_*\Vert\mathbf v_y\Vert_{H^1}^2
\le C(\Vert\mathbf v\Vert_{H^1}^2+\Vert F\Vert_{L^2}^2)$ with $c_*>0$ and $C$ depending on
$\sup|\mathcal M|$; this gives existence, uniqueness and the stated regularity, and
the weighted identities above hold by approximation of $\Upsilon_\sigma$ and spatial
cutoffs.
\end{proof}

\begin{lemma}\label{l:flx}
Let $a>0$, $\sigma=\pm1$, and let $v\in C([r,T];L^2)\cap L^2(r,T;H^1)$ for every
finite $T>r$ be a weak solution of
\[
 v_t=\pa_y\bigl(d(V_c)v_y-s_\sigma(V_c)v+f_{\rm flx}\bigr),\qquad
 \int_r^\infty\bigl(\Vert v\Vert_{L^2}^2+\tau\Vert f_{\rm flx}\Vert_{L^2}^2\bigr)\dd t<\infty .
\]
For every sufficiently late fixed $r$ and every $s\ge r$, with $C$ depending only on
$a$ and the fixed data,
\[
\begin{gathered}
 \sup_{t\ge s}\tau\Vert v(t)\Vert_{L^2}^2+\int_s^\infty\tau\Vert v_y\Vert_{L^2}^2\dd t
   +\mathscr O^a_{s,\infty}(v)
 \le C\Bigl[\tau_s\Vert v(s)\Vert_{L^2}^2
   +\int_s^\infty\bigl(\Vert v\Vert_{L^2}^2+\tau\Vert f_{\rm flx}\Vert_{L^2}^2\bigr)\dd t\Bigr],\\
 \tau_t\Vert v(t)\Vert_{L^2}^2\le C\int_{t/2}^\infty
   \bigl(\Vert v\Vert_{L^2}^2+\tau\Vert f_{\rm flx}\Vert_{L^2}^2\bigr)\dd t'\qquad(t\ge2r).
\end{gathered}
\]
If $v=u_y$ and $L_\sigma u=f_{\rm flx}$, then
$\mathscr O^a_{s,\infty}(u_t)\le C\int_s^\infty\tau(\Vert v_y\Vert_{L^2}^2
+\Vert f_{\rm flx}\Vert_{L^2}^2)\dd t+C\mathscr O^a_{s,\infty}(v)$.
\end{lemma}
\begin{proof}
Write $d=d(V_c)$. The equation reads
$Pv+\pa_y[s_\sigma(V_c)]v=\pa_yf_{\rm flx}$ with
$P=\pa_t-d\pa_y^2+(s_\sigma(V_c)-d_y)\pa_y$. Let
$\Upsilon=T_{L_\sigma,r}W^a$, which is bounded by \eqref{e:gauss}, and apply
\eqref{e:rcw} with $w=w_1$ of \eqref{e:nfwdef}. Here
$P\Upsilon=L_\sigma\Upsilon-d_y\Upsilon_y$,
$P^*w_1+2\pa_y[s_\sigma(V_c)]w_1
=-[(w_1)_t+(w_1s_\sigma(V_c))_y-2w_1\pa_y[s_\sigma(V_c)]+((w_1)_yd)_y]\ge c_1g_r$ by
\eqref{e:nfwts}, and
\[
 -\int(\omega v)_yf_{\rm flx}\dd y
 =-\int\Bigl(\omega v_y+\frac{(w_1)_y}{\sqrt{1+\Upsilon}}v
   -\frac12\omega z_\Upsilon v\Bigr)f_{\rm flx}\dd y .
\]
With $|(w_1)_y|^2\le Cg_r$ and
$\frac14|d_y\omega z_\Upsilon|v^2\le\varepsilon d\omega z_\Upsilon^2v^2+C\delta^2\tau^{-1}v^2$,
this gives, for $E_v=\frac12\int\omega v^2\asymp\Vert v\Vert_{L^2}^2$,
\begin{equation}\label{e:flxe}
 E_v'+c_*\Bigl(\Vert v_y\Vert_{L^2}^2+\int_\R W^av^2\dd y\Bigr)
 \le C\delta^2\tau^{-1}\Vert v\Vert_{L^2}^2+C\Vert f_{\rm flx}\Vert_{L^2}^2,
\end{equation}
with $c_*>0$ depending only on $a$ and the fixed data. Multiplying by $\tau$,
\[
 (\tau E_v)'+c_*\tau\Bigl(\Vert v_y\Vert_{L^2}^2+\int_\R W^av^2\dd y\Bigr)
 \le C\bigl(\Vert v\Vert_{L^2}^2+\tau\Vert f_{\rm flx}\Vert_{L^2}^2\bigr);
\]
integration from $s$ gives the first line, since $\tau W^a=\sqrt\tau G^a$, and
\eqref{e:avst} with $k=0$ gives the second. If $v=u_y$, then
$u_t=dv_y-s_\sigma(V_c)v+f_{\rm flx}$, which gives the last bound. The tests are
justified by time averaging and spatial cutoffs.
\end{proof}

\begin{lemma}\label{l:bulk}
Let $0<a\le c_G$. For every sufficiently late fixed $r$,
\[
 \mathscr O^a_{r,\infty}\bigl(\pa_y(P_b-p_a)\bigr)
 +\mathscr O^a_{r,\infty}\bigl(\pa_y(b_m-b_a)\bigr)<\infty .
\]
\end{lemma}
\begin{proof}
We estimate the two acoustic derivatives together. Use $P=P_b-p_a$ and
$\Xi_P,\Xi_B$ of \eqref{e:wdif}, and put $g_P=\frac{\gamma P_b}{V_b}$.
Since $P$ is a smooth function of the exact and approximate states, vanishing
when $W=0$, \eqref{e:wbas}, \eqref{e:wval}, \eqref{e:cmpb}, \eqref{e:wtC} and
\eqref{e:wdif} give
\[
 \sup_{t\ge r}\tau\Vert(P_y,(W_b)_y)\Vert_{L^2}^2
 +\int_r^\infty\bigl(\Vert(P_y,(W_b)_y)\Vert_{L^2}^2+\tau\Vert(P_{yy},(W_b)_{yy})\Vert_{L^2}^2
   +\tau\Vert(\Xi_P,\Xi_B)\Vert_{L^2}^2\bigr)\dd t<\infty ;
\]
With $\xi=\frac{y+h_*}{\sqrt\tau}$, set
\[
 I=\tau\int_0^{\xi}e^{-\frac{az^2}2}\dd z,\qquad \mathcal J(t)=\int_\R IP_y(W_b)_y\dd y .
\]
Then $I_y=\sqrt\tau G^a$, $|I|\le C\tau$ and $|I_t|\le C$. Differentiating
$P_t+g_P(W_b)_y=\Xi_P$ and $(W_b)_t+P_y=\Xi_B$ in $y$ and integrating by parts,
\[
\begin{aligned}
 \mathcal J'={}&\frac12\int I_y(P_y^2+g_P(W_b)_y^2)\dd y-\frac12\int I(g_P)_y(W_b)_y^2\dd y
   +\int I_tP_y(W_b)_y\dd y\\
 &-\int(I_y(W_b)_y+I(W_b)_{yy})\Xi_P\dd y-\int(I_yP_y+IP_{yy})\Xi_B\dd y .
\end{aligned}
\]
The first term is at least $c_*(A^a_{P_y}+A^a_{(W_b)_y})$, where $c_*>0$ depends only on $a$ and the fixed data. Write
$g_P=\frac{\gamma p_*}{V_h}+R_G$; by Lemma~\ref{l:rate} and \eqref{eq:pr-gau},
$\Vert R_G\Vert_{L^\infty}\le C_b\tau^{-\frac12}$, and
$|I\pa_y\frac{\gamma p_*}{V_h}|\le C\delta\sqrt\tau G^{2c_G}\le C\delta\sqrt\tau G^a$.
Moreover
\[
 -\frac12\int I(R_G)_y(W_b)_y^2\dd y=\frac12\int I_yR_G(W_b)_y^2\dd y+\int IR_G(W_b)_y(W_b)_{yy}\dd y,
\]
where the first term is at most $C_b\tau^{-\frac12}A^a_{(W_b)_y}$ and the second at most
$C_b\Vert (W_b)_y\Vert_{L^2}\sqrt\tau\Vert (W_b)_{yy}\Vert_{L^2}$. The $\Xi$-terms are at most
$\frac{c_*}4(A^a_{P_y}+A^a_{(W_b)_y})+C\tau(\Vert(P_{yy},(W_b)_{yy})\Vert_{L^2}^2
+\Vert(\Xi_P,\Xi_B)\Vert_{L^2}^2)$, and $|\int I_tP_y(W_b)_y|\le C\Vert(P_y,(W_b)_y)\Vert_{L^2}^2$.
For $C\delta\le\frac{c_*}4$ and $C_b\tau_r^{-\frac12}\le\frac{c_*}4$, integration over $(r,T)$
and $|\mathcal J|\le C\tau\Vert(P_y,(W_b)_y)\Vert_{L^2}^2$ give
$\mathscr O^a_{r,T}(P_y)+\mathscr O^a_{r,T}((W_b)_y)\le C_b$ uniformly in $T$.
\end{proof}

\begin{lemma}\label{l:bcor}
For every sufficiently late fixed $r$ and $\sigma=\pm1$,
\begin{equation}\label{e:bfsq}
 \int_r^\infty\tau\bigl(\Vert\pa_yF_{\rm bulk}\Vert_{L^2}^2
   +\Vert G_{\rm bulk}\Vert_{L^2}^2\bigr)\dd t<\infty ,
\end{equation}
$J_r(u_{\sigma,r})<\infty$, and the quantity
\[
 \sup_{t\ge S}\tau\Vert(u_{\sigma,r})_y\Vert_{L^2}^2
 +\int_S^\infty\tau\Vert(u_{\sigma,r})_{yy}\Vert_{L^2}^2\dd t
 +\mathscr O_{S,\infty}((u_{\sigma,r})_y)+\mathscr O_{S,\infty}((u_{\sigma,r})_t)
 +\mathscr O_{S,\infty}((W_{\sigma,r})_t)+\mathscr O_{S,\infty}(j_r)
\]
is finite for $S=r$ and tends to zero as $S\to\infty$; here $F_{\rm bulk},G_{\rm bulk}$
are those of \eqref{e:nfs} for the index $\sigma$.
\end{lemma}
\begin{proof}
With the identity for $z$ in the proof of Lemma~\ref{l:am}, \eqref{e:wtC},
\eqref{e:munw}, Lemmas~\ref{l:dep}, \ref{l:pr}, \ref{l:asy} and \eqref{eq:pr-gau} give
\[
 \Vert(z,\widetilde z)\Vert_{L^\infty}\le C_b\tau^{-\frac12},\qquad
 \sup_{t\ge r}\tau\Vert(z_y,\widetilde z_y)\Vert_{H^1}^2
 +\int_r^\infty\bigl(\Vert(z_y,\widetilde z_y)\Vert_{L^2}^2
   +\tau\Vert(z_{yy},\widetilde z_{yy})\Vert_{L^2}^2\bigr)\dd t<\infty ;
\]
the profile terms of $z$ have $\Vert\pa_y^2(\cdot)\Vert_{L^2}^2\le C_b\tau^{-\frac52}$,
since $V_h-V_c$ is a difference of translates. By \eqref{e:nfs} and Lemma~\ref{l:nf}(i),
\[
 |G_{\rm bulk}|\le C(|z|^2|z_y|+|z_y|^2),\qquad
 |\pa_yF_{\rm bulk}|\le C\bigl(|\widetilde z||\widetilde z_y|+|(V_c)_y||\widetilde z|^2
   +|z_{yy}|+(|z_y|+|(V_c)_y|)|z_y|\bigr).
\]
Since $\sqrt\tau\Vert(z,\widetilde z,z_y,(V_c)_y)\Vert_{L^\infty}\le C_b$,
these products satisfy the single bound
\[
 \tau\bigl(\Vert\pa_yF_{\rm bulk}\Vert_{L^2}^2+\Vert G_{\rm bulk}\Vert_{L^2}^2\bigr)
 \le C_b\Bigl(\tau\Vert z_{yy}\Vert_{L^2}^2
   +\Vert(z_y,\widetilde z_y)\Vert_{L^2}^2+\tau^{-3/2}\Bigr).
\]
Its right side is integrable, proving \eqref{e:bfsq}.

Write $u=u_{\sigma,r}$ and $F=F_{\rm bulk}$, $G=G_{\rm bulk}$. Apply \eqref{e:rcw} with
$P=L_\sigma$, $w=w_0$ and $\Upsilon=T_{L_\sigma,r}W_c$. Since
$-\int(\omega u)_yF=-\int(\omega u_y+\frac{(w_0)_y}{\sqrt{1+\Upsilon}}u
-\frac12\omega z_\Upsilon u)F$, the work of $\pa_yF+G$ is at most
$\frac{c_*}2\int(u_y^2+(g_r+z_\Upsilon^2)u^2)+C\Vert F\Vert_{L^2}^2+C\Vert u\Vert_{L^2}\Vert G\Vert_{L^2}$,
and Lemma~\ref{l:nf}(iii), (v) give $J_r(u)<\infty$. The function $v=u_y$ solves
$v_t=\pa_y(d(V_c)v_y-s_\sigma(V_c)v+\pa_yF+G)$, $v(r)=0$; Lemma~\ref{l:flx} with
$a=a_s$, \eqref{e:nfe} and \eqref{e:bfsq} give the assertions for $u$.

Lemma~\ref{l:bulk} and \eqref{e:wdif} give
$\mathscr O_{r,\infty}(P_t)+\mathscr O_{r,\infty}((W_b)_t)+\mathscr O_{r,\infty}(\Sigma_t)<\infty$,
since $\int_r^\infty\sqrt\tau\Vert\Xi\Vert_{L^2}^2\dd t<\infty$.
The smooth inverse map $(P,W_b,\Sigma)\mapsto W$ and \eqref{e:wval}, written as
$(R_U)_3=c_vW_T+b_aW_b+\frac12W_b^2$, give
\[
\begin{aligned}
 &A_{W_t}+A_{z_t}+A_{\pa_t[\mathcal P_\sigma(V_c,z)]}\\
 &\qquad\le C_b\Bigl(A_{P_t}+A_{(W_b)_t}+A_{\Sigma_t}+\tau^{-2}+\tau e^{-c_bt}\Bigr).
\end{aligned}
\]
Indeed, the inverse vanishes at $(P,W_b,\Sigma)=0$, so its time dependence through
$(v_a,b_a,T_a)$ contributes at most $C_b\tau^{-1}|W|$. Moreover,
$z=R_U+(v_a,b_a,E_a)-U_*(V_c)$ has only translated profile terms near the contact
and acoustic Gaussians away from it. The chain rule for $\mathcal P_\sigma$ adds at most
$C|z||(V_c)_t|$. Thus the preceding observations and
\eqref{e:ubC} give $\mathscr O_{r,\infty}((W_{\sigma,r})_t)<\infty$. Finally
\[
 c\,\pa_t\log\frac{V_c(y,t)}{V_c(y+h_*,t)}
 =-c\int_0^{h_*}\pa_y\pa_t\log V_c(y+z',t)\dd z',\qquad
 \Bigl\Vert\pa_t\log\frac{V_c}{V_c(\cdot+h_*,\cdot)}\Bigr\Vert_{L^2}^2\le C_b\tau^{-\frac52},
\]
so \eqref{e:bpot} gives $\mathscr O_{r,\infty}(j_r)<\infty$. The tails tend to zero as
tails of convergent integrals, and the supremum by the second bound of
Lemma~\ref{l:flx}.
\end{proof}

\subsection{Source estimates and solvability}

\begin{lemma}\label{l:emi}
For every sufficiently late fixed $r$,
\[
 \int_r^\infty\tau^{-\frac12}\sum_{j=\pm1}\int_\R e^{-b_{\rm ac}z_j^2}n_\chi^2\dd y\dd t<\infty,
\]
and every bounded measurable family $C_{\sigma j}$ satisfies
\[
 \sum_{\sigma,j=\pm1}\int_r^\infty\sqrt\tau\int_\R e^{b_{\rm ac}z_j^2}
   |a_jC_{\sigma j}n_\chi|^2\dd y\dd t<\infty .
\]
\end{lemma}
\begin{proof}
By \eqref{e:pot}, $\chi_m=V_b+\frac{\mathcal F}{\gamma p_*}$, and $R=(\gamma-1)c_v$ gives
$R\Theta_b=p_*\chi_m+\frac{\gamma-1}\gamma\mathcal F-\frac{\gamma-1}2b_m^2$. Hence \eqref{eq:bg-m}, the
second line of \eqref{e:potq}, $\frac{(\gamma-1)\kappa}{\gamma R}=c$ and $(V_c)_t=\pa_y\frac{c(V_c)_y}{V_c}$
(by \eqref{eq:prof} and \eqref{eq:pr-id}) give, with $D=\frac c{V_b}$,
\[
\begin{gathered}
 (n_\chi)_t-\pa_y\bigl(D(n_\chi)_y\bigr)=\pa_yJ,\\
 J=c\Bigl(\frac1{V_b}-\frac1{V_h}\Bigr)(V_h)_y+\frac{(\gamma-1)c\mathcal F_y}{\gamma p_*V_b}
   +\frac{(\gamma-1)(\nu-\gamma c)b_m(b_m)_y}{\gamma p_*V_b}-\frac{\gamma-1}{\gamma p_*}(P_b-p_*)b_m .
\end{gathered}
\]
For $f\in H^1$ with $Q_y=f$ and $Q$ bounded,
$\int f^4=-3\int Qf^2f_y\le3\Vert Q\Vert_{L^\infty}\Vert f\Vert_{L^4}^2\Vert f_y\Vert_{L^2}$,
hence $\int f^4\le9\Vert Q\Vert_{L^\infty}^2\Vert f_y\Vert_{L^2}^2$. Applied to
$(f,Q)=(b_m,\mathcal B)$ and $(\mathcal F,\mathcal G)$, with \eqref{e:bdpot},
\eqref{e:munw} and $\mathcal F=V_b(P_b-p_*)+\frac{\gamma-1}2b_m^2$, this gives
\[
 \int_r^\infty\!\!\int_\R\bigl(b_m^4+\mathcal F^4+|(P_b-p_*)b_m|^2+|P_b-p_*|^4\bigr)
 \dd y\dd t<\infty .
\]
For $\tau\ge h_*^2$, the first term of $J$ is bounded by
$C\delta\tau^{-\frac12}e^{-\frac{c_Gy^2}{2\tau}}(|V_b-V_c|+|V_c-V_h|)$; the Gaussian
integral of Lemma~\ref{l:dep} and \eqref{eq:pr-gau} make it square integrable, and
\eqref{e:munw}, Lemma~\ref{l:dep} and Lemma~\ref{l:rate} do so for the other terms. Thus
$J\in L^2((r,\infty)\times\R)$, and testing by $n_\chi$ gives
$\sup_{t\ge r}\Vert n_\chi\Vert_{L^2}^2+\int_r^\infty\Vert(n_\chi)_y\Vert_{L^2}^2\le C_b$,
since $n_\chi=\ell_0R_U$ and $R_U(r)\in L^1\cap L^\infty$.

Let $r$ be so large that $\Vert D_y(t)\Vert_{L^\infty}\le\frac14\min\{c_-,c_+\}$ for
$t\ge r$ (by \eqref{e:munw}). For $P=\pa_t-\pa_y(D\pa_y)=\pa_t-D\pa_y^2-D_y\pa_y$ and
$j=\pm1$, $-j(-D_y-s_j)\ge\frac34c_j$, and \eqref{e:gauss} with $(\sigma,v,h)=(-j,s_j,0)$
bounds
\[
 \Upsilon=\sum_{j=\pm1}T_{P,r}\bigl(\tau^{-\frac12}e^{-b_{\rm ac}z_j^2}\bigr)
\]
by a constant depending only on the fixed data. Apply \eqref{e:rcw} with $w=1$, for
which $P^*w=0$, and $u=n_\chi$: since $(\omega n_\chi)_y=\omega((n_\chi)_y-\frac12z_\Upsilon n_\chi)$,
\[
 \frac12\frac{\dd}{\dd t}\int\omega n_\chi^2
 +\int D\omega\Bigl(\frac{(n_\chi)_y^2}8+\frac{z_\Upsilon^2n_\chi^2}{48}\Bigr)
 +\int\frac{P\Upsilon}{4(1+\Upsilon)^{\frac32}}n_\chi^2\le C\Vert J\Vert_{L^2}^2,
\]
and $P\Upsilon=\tau^{-\frac12}\sum_je^{-b_{\rm ac}z_j^2}$ gives the first assertion. The
second follows from $\sqrt\tau e^{b_{\rm ac}z_j^2}|a_j|^2\le C_M\tau^{-\frac12}e^{-b_{\rm ac}z_j^2}$.
\end{proof}

\begin{lemma}\label{l:solv}
Let $1<p_w<\frac32$, $C_A\ge1$ and $b_A\ge\frac{a_s}2$. For every sufficiently late
fixed $r$ the following holds. Let $g\in H^1$ and $f=f^{\rm a}+f^{\rm b}+A\varphi$, where:
\begin{enumerate}[label=(\alph*)]
\item $f^{\rm a}_\sigma=\sum_{j=\pm1}f_{\sigma j}$ with
\begin{equation}\label{e:psrc}
 \sum_{\sigma,j=\pm1}\int_r^\infty\sqrt\tau\int_\R(1+z_j^2)^{p_w}|f_{\sigma j}|^2\dd y\dd t<\infty ;
\end{equation}
\item $f^{\rm b}=f^{\rm g}+f^{\rm w}+\sum_iK_i\pa_yw_i$ (a finite sum), with
\[
 \int_r^\infty\bigl(\Vert f^{\rm g}\Vert_{L^2}+\tau\Vert f^{\rm g}\Vert_{L^2}^2\bigr)\dd t
 +\int_r^\infty\!\!\int_\R\Bigl(\frac1{W_c}+\tau\Bigr)|f^{\rm w}|^2\dd y\dd t
 +\sum_i\int_r^\infty\Vert\pa_yw_i\Vert_{L^2}^2\dd t<\infty,
\]
and matrices $K_i$ with
$|K_i|\le C_K\bigl(|a_+|+|a_-|+|R_U|+\tau^{-1}+\tau^{-\frac12}G_0^{\frac12}\bigr)$;
\item $|\pa_y^jA|\le C_A\delta\tau^{-\frac{j+1}2}e^{-\frac{b_A(y+h_*)^2}\tau}$ $(0\le j\le2)$,
$|A_t|\le C_A\delta\tau^{-\frac32}e^{-\frac{b_A(y+h_*)^2}\tau}$, and
$\varphi\in C([r,T];L^2)$ with $J_r(\varphi)<\infty$ and $\mathscr O_{r,T}(\varphi_t)<\infty$
for every finite $T$.
\end{enumerate}
Then $\mathcal KU=f$, $U(r)=g$, with $\mathcal K$ of \eqref{e:pB}, has a unique solution
with $U\in C([r,T];H^1)$ and $U_t,U_{yy}\in L^2((r,T)\times\R)$ for finite $T$, and
\begin{equation}\label{e:solv}
 J_r(U)<\infty,\qquad
 \mathscr O_{r,T}(U_t)\le C_1+C_0\delta^2\mathscr O_{r,T}(\varphi_t)\qquad(r<T<\infty),
\end{equation}
where $C_0$ depends only on $C_A$, $b_A$, the exponents of \eqref{e:pwt} and the fixed
data, and $C_1$ does not depend on $T$. If $g=0$ and $f=A\varphi$, then $J_r(U)\le C_{b,r}\delta^2J_r(\varphi)$ and
$\mathscr O_{r,T}(U_t)\le C_0\delta^2\mathscr O_{r,T}(\varphi_t)+C_{b,r}\delta^2J_r(\varphi)$.
\end{lemma}
\begin{proof}
Write $d=d(V_c)$ and $k=k(V_c)$; let
$\tau_r\ge h_*^2$ with \eqref{e:nfwdef} applicable.

\emph{Step 0.} Put $Q_\sigma=\tau^{-1}e^{-b_{\rm ac}z_\sigma^2}$ and
$N_\sigma=\tau^{-\frac12}(1+z_o^2)^{-p_w}$,
\[
 \eta_r=\max_\sigma\bigl\Vert T_{L_\sigma,r}\bigl(|\widetilde{\mathcal R}|^2+|R_U|^2
   +\tau^{-2}\bigr)\bigr\Vert_{L^\infty}+\tau_r^{-\frac12},\qquad
 \mathfrak v_\sigma=c_{\mathfrak v}\Bigl[Q_\sigma+N_\sigma
   +\frac{|\widetilde{\mathcal R}|^2+|R_U|^2+\tau^{-2}}{\eta_r}\Bigr].
\]
By \eqref{e:pmat}, Lemma~\ref{l:sq} and Lemma~\ref{l:pos}(i), $\eta_r\to0$; by
\eqref{e:qresp} and the first bound of \eqref{e:gauss} with $(v,h)=(s_o,0)$ there is
$c_{\mathfrak v}>0$, depending only on $M$ and the fixed data, with
$\Vert T_{L_\sigma,r}\mathfrak v_\sigma\Vert_{L^\infty}\le1$. Since $|a_\sigma|^2\le C_MQ_\sigma$,
$|a_o|^2\le C_MN_\sigma$ and $\tau^{-1}G_0=\tau^{-\frac12}W_c$, every $K$ as in (b) has
$|K|^2\le C_{K,b}(W_c+\mathfrak v_\sigma)$ in row $\sigma$ and $\tau|K|^2\le C_{K,b}$; by
\eqref{e:pmat}, $K=\widetilde{\mathcal M}$ is such a matrix. Lemma~\ref{l:geo} applies
to $\widetilde{\mathcal M}$ by Lemma~\ref{l:am}, with $a=2b_0$, $W^{2b_0}=W_c$.

\emph{Step 1: $f=(f_{\sigma\sigma})_\sigma$.} Let $(L_\sigma-\widetilde{\mathcal M}_{\sigma\sigma}\pa_y)Z_\sigma
=f_{\sigma\sigma}$, $Z_\sigma(r)=0$. Put $w_\sigma=-2ka_o-\widetilde{\mathcal R}_{\sigma\sigma}$,
so that $\widetilde{\mathcal M}_{\sigma\sigma}=-2ka_\sigma-w_\sigma$, and $v=p_\sigma Z_\sigma$
with $p_\sigma=e^{-\alpha A_\sigma}$. By Lemma~\ref{l:pr},
\[
 L_\sigma v+w_\sigma v_y=p_\sigma f_{\sigma\sigma}+c_\sigma v,\qquad
 c_\sigma=-\alpha w_\sigma a_\sigma-\alpha R^{\rm prof}_\sigma,
\]
and by Lemma~\ref{l:asy}, \eqref{e:pmat} and Lemma~\ref{l:pr},
$\sup_y\tau|c_\sigma|\to0$ as $t\to\infty$; by Lemma~\ref{l:opp} and \eqref{e:pmat},
$\eta'_r=\Vert T_{L_\sigma,r}w_\sigma^2\Vert_{L^\infty}+\tau_r^{-\frac12}\to0$. In
$\xi=F_c(y,t)$ of \eqref{e:nco}, with $(d,B)=(d(V_c),s_\sigma(V_c))$, $L_\sigma$ becomes
$\pa_t-D_0\pa_\xi^2+b\pa_\xi$ with $|b-\bar b|\le C\delta\tau^{-\frac12}$,
$|b_\xi|\le C\delta\tau^{-1}$ and $\bar b=s_\sigma(\rho_-^{-1})$. Put $\eta=\xi-\bar bt$ and
$B_c=b-\bar b$. Since $f_cs_\sigma(V_c)=\bar b$ and, by \eqref{eq:ct-gau},
$|F_c(y,t)-f_c(v_\sigma)y|\le C\delta(\sqrt\tau+|y-s_\sigma t|)$,
\[
 C^{-1}(1+z_\sigma^2)\le1+\frac{\eta^2}\tau\le C(1+z_\sigma^2).
\]
For $\mathcal W=\tau^{-\frac12}(1+\frac{\varepsilon_w\eta^2}\tau)^{p_w}$ and
$q=\frac\eta{\sqrt\tau}$,
\[
 \frac{\mathcal W_t+D_0\mathcal W_{\eta\eta}}{\mathcal W}
 =\frac1\tau\Bigl[-\frac12-\frac{p_w\varepsilon_wq^2}{1+\varepsilon_wq^2}
   +\frac{2D_0p_w\varepsilon_w}{1+\varepsilon_wq^2}
   +\frac{4D_0p_w(p_w-1)\varepsilon_w^2q^2}{(1+\varepsilon_wq^2)^2}\Bigr]
 \le\frac1\tau\Bigl[-\frac12+\varepsilon_wD_0p_w(p_w+1)\Bigr],
\]
and $|\pa_\eta\log\mathcal W|\le p_w\sqrt{\varepsilon_w}\tau^{-\frac12}$ gives
$|(B_c\mathcal W)_\eta|\le C\delta\tau^{-1}\mathcal W$. Fix $\varepsilon_w$ with
$\varepsilon_wD_0p_w(p_w+1)\le\frac18$; for $C\delta\le\frac18$,
\[
 \mathcal W_t+D_0\mathcal W_{\eta\eta}+(B_c\mathcal W)_\eta\le-\frac1{4\tau}\mathcal W .
\]
Apply \eqref{e:rcw} in $(\eta,t)$ with $P=\pa_t-D_0\pa_\eta^2+B_c\pa_\eta$,
$w=\mathcal W$, $\Upsilon=\frac1{\eta'_r}T_{L_\sigma,r}w_\sigma^2$ and $u=v$; the
equation reads $Pv+f_cw_\sigma v_\eta=p_\sigma f_{\sigma\sigma}+c_\sigma v$. The drift
work is at most $\frac{D_0}8\int\omega v_\eta^2+C\eta'_r\int\omega\frac{w_\sigma^2}{\eta'_r}v^2$,
the term with $c_\sigma$ at most $o(1)\int\tau^{-1}\mathcal Wv^2$, and the source work at most
$\frac1{32}\int\tau^{-1}\mathcal Wv^2+C\int\tau\mathcal W|f_{\sigma\sigma}|^2$. For late $r$
this gives, after return to $y$ and $Z_\sigma=p_\sigma^{-1}v$,
\begin{equation}\label{e:pown}
\begin{aligned}
 &\sup_{t\ge r}\tau^{-\frac12}\int_\R(1+z_\sigma^2)^{p_w}|Z_\sigma|^2\dd y
 +\int_r^\infty\!\!\int_\R(1+z_\sigma^2)^{p_w}
   \bigl(\tau^{-\frac12}|(Z_\sigma)_y|^2+\tau^{-\frac32}|Z_\sigma|^2\bigr)\dd y\dd t\\
 &\quad\le C_b\int_r^\infty\!\!\int_\R\sqrt\tau(1+z_\sigma^2)^{p_w}|f_{\sigma\sigma}|^2\dd y\dd t .
\end{aligned}
\end{equation}
Since $\tau G_0\le C_b(1+z_\sigma^2)^{p_w}$, this bounds $J_r(Z)$ and
$\mathscr O^{4b_0}_{r,\infty}(Z_y)$.

\emph{Step 2: $f=(f_{\sigma o})_\sigma$.} Let $\mathscr P_\sigma V_\sigma=h_\sigma$,
$V_\sigma(r)=0$, with $h_\sigma=f_{\sigma o}+\widetilde{\mathcal M}_{\sigma o}(Z_o)_y$. By
\eqref{e:pmat} and \eqref{e:pown} for the index $o$,
$\int_r^\infty\!\int\sqrt\tau(1+z_o^2)^{p_w}|h_\sigma|^2<\infty$. Since
$\mathscr P_\sigma=L_\sigma-\mathcal M_{\sigma\sigma}\pa_y$ with $\mathcal M_{\sigma\sigma}=-2ka_\sigma$,
whose residual $-2ka_o$ has $\Vert T_{L_\sigma,r}(a_o^2)\Vert_{L^\infty}\le C_M\tau_r^{-\frac12}$
by Lemma~\ref{l:opp}, Lemma~\ref{l:geo} for $L-\mathcal M^{\rm d}\pa_y$ with
$\mathfrak v_\sigma=c_{\mathfrak v}N_\sigma$ gives $J_r(V)+\int_r^\infty\Vert V_y\Vert_{L^2}^2<\infty$.

For the contact derivative split $h_\sigma=h^{\rm c}_\sigma+h^{\rm e}_\sigma$,
$h^{\rm c}_\sigma=h_\sigma\mathbf 1_{\{|y-s_ot|\le\varepsilon_ct\}}$, with $\varepsilon_c>0$
small, and $V=V^{\rm c}+V^{\rm e}$ accordingly. On the support of $h^{\rm e}_\sigma$,
$1+z_o^2\ge c_*\tau$, so
$\int_r^\infty\tau^{p_w+\frac12}\Vert h^{\rm e}_\sigma\Vert_{L^2}^2<\infty$; hence
$\int_r^\infty(\Vert h^{\rm e}_\sigma\Vert_{L^2}+\tau\Vert h^{\rm e}_\sigma\Vert_{L^2}^2)<\infty$.
Lemma~\ref{l:geo} with $F_1=h^{\rm e}$ bounds $\int_r^\infty\Vert V^{\rm e}_y\Vert_{L^2}^2$, and
Lemma~\ref{l:flx} for $v=(V^{\rm e}_\sigma)_y$, with
$f_{\rm flx}=h^{\rm e}_\sigma-2ka_\sigma v$, gives $\mathscr O^a_{r,\infty}(V^{\rm e}_y)<\infty$
for every $a>0$.

By \eqref{e:nco} and \eqref{e:xder},
\[
\begin{gathered}
 (V^{\rm c}_\sigma)_y(y,t)=\int_r^t\int_\R\mathsf k(t,y;s,z)h^{\rm c}_\sigma(z,s)\dd z\dd s,\\
 \mathsf k=f_c(y,t)f_c(z,s)K^{\mathscr P}_\xi\bigl(t,s;F_c(y,t),F_c(z,s)\bigr),
\end{gathered}
\]
and by Lemma~\ref{l:kernv}(i) and Lemma~\ref{l:kern}(b),
\[
 |\mathsf k|\le C_M\Bigl(\frac{\tau_t}{\tau_s}\Bigr)^{C\delta}(t-s)^{-\frac12}G\bigl(C_M(t-s),X\bigr),\qquad
 X=F_c(y,t)-F_c(z,s)-\bar b(t-s) .
\]
For $|z-s_os|\le\varepsilon_cs$, $f_c(v_o)s_o=-\bar b$ and \eqref{eq:ct-gau} give
$X=F_c(y,t)-\bar b(t-2s)-f_c(v_o)(z-s_os)-E_1$ with $|E_1|\le C\delta\sqrt{\tau_s}$. On
$\{|y+h_*|\ge\varepsilon_ct\}$, $\sqrt\tau G_0^2\le e^{-b_0\varepsilon_c^2t/2}$, so this part of
$\mathscr O^{4b_0}_{r,\infty}(V^{\rm c}_y)$ is bounded by
$\int_r^\infty\Vert V^{\rm c}_y\Vert_{L^2}^2$. On $\{|y+h_*|\le\varepsilon_ct\}$ put
$P_o(z,s)=(1+z_o(z,s)^2)^{p_w}$ and
\[
 \mathcal S(t,y;s,z)=\mathbf 1_{\{|z-s_os|\le\varepsilon_cs\}}\tau_t^{\frac14}G_0(y,t)|\mathsf k(t,y;s,z)|
   \tau_s^{-\frac14}P_o(z,s)^{-\frac12};
\]
then $\tau_t^{\frac14}G_0|(V^{\rm c}_\sigma)_y|\le\int\!\!\int\mathcal S\,
\tau_s^{\frac14}P_o^{\frac12}|h^{\rm c}_\sigma|\dd z\dd s$, and by the Schur test it suffices
to bound $\sup_{t,y}\iint\mathcal S\dd z\dd s$ and $\sup_{s,z}\iint\mathcal S\dd y\dd t$.

If $\frac t8\le s\le\frac{7t}8$, then $t-s\asymp\tau_s\asymp\tau_t$. With $u=z_o(z,s)$,
$\dd z=\sqrt{\tau_s}\dd u$, the Gaussian $G(C_M(t-s),X)$ has width of order one in $u$ and
$|E_1|\le C\delta\sqrt{\tau_s}$, so
\[
 \int_\R\mathcal S\dd z\le C_M\tau_t^{-\frac12}G_0(y,t)
   \Bigl(1+\frac{(F_c(y,t)-\bar b(t-2s))^2}{\tau_t}\Bigr)^{-\frac{p_w}2},
\]
and the substitution $u=\frac{F_c(y,t)-\bar bt+2\bar bs}{\sqrt{\tau_t}}$ bounds the
integral over $s$ by $C_M\int_\R(1+u^2)^{-\frac{p_w}2}\dd u$, finite as $p_w>1$. For
fixed $(s,z)$, $\pa_tX=(F_c)_t-\bar b$ with $|(F_c)_t|\le C\delta\tau^{-\frac12}$, so
$\int G(C_M(t-s),X)\dd t\le C_M$, and $G_0(y,t)\le e^{-\frac{b_0(y+h_*)^2}{8\tau_s}}$ for $t\le8s$;
hence $\iint\mathcal S\dd y\dd t\le C_M\int\tau_s^{-\frac12}e^{-\frac{b_0(y+h_*)^2}{8\tau_s}}\dd y\le C_M$.

If $s<\frac t8$ or $s>\frac{7t}8$, then $|\bar b(t-2s)|\ge\frac34|\bar b|t$, while
$|F_c(y,t)|\le C(\varepsilon_ct+|h_*|)$ and $|f_c(v_o)(z-s_os)|\le C\varepsilon_ct$; for
$\varepsilon_c$ small and $t\ge C_b$, $|X|\ge\frac{|\bar b|}2t$. With $t-s\le t$ this gives
$|\mathsf k|\le C_Me^{-c_1t}G(2C_M(t-s),X)$, $c_1>0$, and since $|\pa_zX|$ and $|\pa_yX|$ are
bounded above and below, both Schur integrals over this region are at most $C_M$; the
times $t\le C_b$ contribute at most $C_b$. Hence
$\mathscr O^{4b_0}_{r,\infty}(V^{\rm c}_y)\le C_M\iint\sqrt\tau(1+z_o^2)^{p_w}|h^{\rm c}_\sigma|^2
+C\int_r^\infty\Vert V^{\rm c}_y\Vert_{L^2}^2<\infty$.

Now let $u\in\{Z_\sigma,V_\sigma\}$, so that $u_t=du_{yy}-\beta_uu_y+f_u$ with
$\beta_u=s_\sigma(V_c)-\widetilde{\mathcal M}_{\sigma\sigma}$, $f_u=f_{\sigma\sigma}$, or
$\beta_u=s_\sigma(V_c)+2ka_\sigma$, $f_u=h_\sigma$. Then $v=u_y$ solves
$v_t=\pa_y(dv_y-\beta_uv+f_u)$, $v(r)=0$. Testing by $\sqrt\tau Gv$, with
$|\pa_t(\sqrt\tau G)|+|\pa_y(\sqrt\tau G)|\le C(\sqrt\tau G\cdot G^{4b_0})^{\frac12}$ and
$\sqrt\tau G\le\sqrt\tau G^{4b_0}$, gives
\[
 \int_r^\infty\!\!\int_\R\sqrt\tau G|u_{yy}|^2\le C\Bigl(\mathscr O^{4b_0}_{r,\infty}(u_y)
   +\int_r^\infty\sqrt\tau\Vert f_u\Vert_{L^2}^2\dd t\Bigr)<\infty,
\]
and the equation gives $\mathscr O_{r,\infty}(u_t)<\infty$.

\emph{Step 3: (a) with trace.} With
$\mathcal A=\widetilde{\mathcal M}+2k\operatorname{diag}(a_+,a_-)$, let
$(L-\widetilde{\mathcal M}\pa_y)Y=\mathcal AV_y$, $Y(r)=g$. Then $U^{\rm a}=Z+V+Y$ solves
$(L-\widetilde{\mathcal M}\pa_y)U^{\rm a}=f^{\rm a}$, $U^{\rm a}(r)=g$. In row $\sigma$,
$\sum_\rho|\mathcal A_{\sigma\rho}|^2\le C(|a_+|^2+|a_-|^2+|\widetilde{\mathcal R}|^2)
\le C_Mc_{\mathfrak v}^{-1}\mathfrak v_\sigma$, so Lemma~\ref{l:geo} for $L-\widetilde{\mathcal M}\pa_y$
gives $\sup\Vert Y\Vert_{L^2}^2+\int_r^\infty\Vert Y_y\Vert_{L^2}^2+J_r(Y)
\le C(\Vert g\Vert_{L^2}^2+C_M\int_r^\infty\Vert V_y\Vert_{L^2}^2)$. The components of
$q=Y_y$ solve $q_t=\pa_y(dq_y-\Lambda q+\widetilde{\mathcal M}q+\mathcal AV_y)$ with
$\int_r^\infty\tau\Vert\widetilde{\mathcal M}q+\mathcal AV_y\Vert_{L^2}^2\le C_b\int_r^\infty
(\Vert q\Vert_{L^2}^2+\Vert V_y\Vert_{L^2}^2)$ and $q(r)=g_y$; Lemma~\ref{l:flx} and
$Y_t=dq_y-\Lambda q+\widetilde{\mathcal M}q+\mathcal AV_y$ give
$\mathscr O_{r,\infty}(Y_t)<\infty$. Thus $J_r(U^{\rm a})+\mathscr O_{r,\infty}(U^{\rm a}_t)<\infty$.

\emph{Step 4: (b).} Let $(L-\widetilde{\mathcal M}\pa_y)U^{\rm b}=f^{\rm b}$, $U^{\rm b}(r)=0$.
By Step 0, Lemma~\ref{l:geo} with $F_1=f^{\rm g}$ and $f=f^{\rm w}+\sum_iK_i\pa_yw_i$ gives
$\sup\Vert U^{\rm b}\Vert_{L^2}^2+\int_r^\infty\Vert U^{\rm b}_y\Vert_{L^2}^2+J_r(U^{\rm b})\le C_b\mathfrak N$,
where
\[
 \mathfrak N=\Bigl(\int_r^\infty\Vert f^{\rm g}\Vert_{L^2}\Bigr)^2
   +\int_r^\infty\tau\Vert f^{\rm g}\Vert_{L^2}^2
   +\int_r^\infty\!\!\int_\R\Bigl(\frac1{W_c}+\tau\Bigr)|f^{\rm w}|^2
   +\sum_i\int_r^\infty\Vert\pa_yw_i\Vert_{L^2}^2 ;
\]
and Lemma~\ref{l:flx} for $U^{\rm b}_y$, with
$f_{\rm flx}=(\widetilde{\mathcal M}U^{\rm b}_y+f^{\rm b})_\sigma$ and
$\int_r^\infty\tau\Vert f^{\rm b}\Vert_{L^2}^2\le C_b\mathfrak N$, gives
$\mathscr O_{r,\infty}(U^{\rm b}_t)\le C_{b,r}\mathfrak N$.

\emph{Step 5: (c).} Let $(L-\widetilde{\mathcal M}\pa_y)U^{\rm c}=A\varphi$, $U^{\rm c}(r)=0$.
Since $b_0\le b_A$, $\frac{|A\varphi|^2}{W_c}\le C_A^2\delta^2W_c|\varphi|^2$, and
Lemma~\ref{l:geo} gives $J_r(U^{\rm c})\le C_b\delta^2\mathscr Q_r(\varphi)$. Let
$L_\sigma Z^0_\sigma=(A\varphi)_\sigma$, $Z^0(r)=0$; Lemma~\ref{l:geo} for $L$ gives
$\sup\Vert Z^0\Vert_{L^2}^2+\int_r^\infty\Vert Z^0_y\Vert_{L^2}^2+\mathscr Q_r(Z^0)
\le C\delta^2\mathscr Q_r(\varphi)$ with $C$ depending only on $C_A$, $b_A$, the fixed Gaussian exponents and the fixed data. If $D$, $B$, $e$ are matrix functions and $z$, $q$, $f_2$ vector functions of
matching dimensions with $z_t-Dz_{yy}+Bz_y=(B-e)q+f_2$, then differentiating in $y$ and
writing $z_y=(z_y-q)+q$ gives
\begin{equation}\label{e:cnr}
 (\pa_t-\pa_yD\pa_y+\pa_yB)(z_y-q)=\pa_y(Dq_y-eq+f_2)-q_t .
\end{equation}
Apply \eqref{e:cnr} with $z=Z^0_\sigma$, $D=d$, $B=s_\sigma(V_c)$, $e=f_2=0$ and
$q=q_\sigma=\frac{(A\varphi)_\sigma}{s_\sigma(V_c)}$, and put $h=(Z^0_\sigma)_y-q_\sigma$,
$x=\sigma(y+h_*)$, $\mathsf s=\sigma s_\sigma(V_c)\ge\underline s=\min\{c_-,c_+\}$. Then
$h_t-\pa_x(dh_x)+\pa_x(\mathsf sh)=\pa_x(dq_x)-q_t$ and $h(r)=-q(r)$. Fix
$0<c_0<\frac{\underline s}4$, $0<\kappa_0<\min\{\frac{a_sc_0}8,\frac{\underline s}{16\sup d}\}$
and a smooth convex $\varrho$ with $\varrho=0$ on $(-\infty,0]$ and $\varrho'=1$ on
$[1,\infty)$, and put
\[
 \phi=\frac x{\sqrt\tau}+(\kappa_0-\tau^{-\frac12})\varrho(x-c_0\tau),\qquad w=\tau e^{-\phi}.
\]
For $\tau$ large, $\tau^{-\frac12}\le\phi_x\le\kappa_0$, $\phi_{xx}\ge0$, and, by
$\varrho(x-c_0\tau)\ge x-c_0\tau-1$ for $x\ge c_0\tau$, $\phi_t+\mathsf s\phi_x\ge(\mathsf s-c_0)\phi_x$.
Since $|d_x|+|\mathsf s_x|\le C\delta\tau^{-\frac12}\le C\delta\phi_x$ and
$d\phi_x\le\frac{\underline s}{16}$,
\[
 \frac{w_t+(dw_x)_x+\mathsf sw_x-\mathsf s_xw}w
 =\frac1\tau-(\phi_t+\mathsf s\phi_x)+d\phi_x^2-d\phi_{xx}-d_x\phi_x-\mathsf s_x
 \le-\frac{\underline s}2\tau^{-\frac12}
\]
for $C\delta\le\frac{\underline s}8$. Moreover $w\ge e^{-\frac1{a_s}}\tau G$ (for
$x\le c_0\tau$ by $\xi\le\frac{a_s\xi^2}2+\frac1{2a_s}$, for $x>c_0\tau$ by
$\kappa_0x\le\frac{a_s\xi^2}8$), and $(1+\sqrt\tau\phi_x^2)w\le C\tau e^{\frac{b_0x^2}\tau}$.
Testing by $wh$ gives
\[
 \frac12\frac{\dd}{\dd t}\int wh^2+\int dwh_x^2+\frac{\underline s}4\int\tau^{-\frac12}wh^2
 \le-\int(wh)_xdq_x-\int whq_t,
\]
and with $|w_x|=\phi_xw$,
\[
 \mathscr O_{r,T}(h)+\int_r^T\!\!\int wh_x^2
 \le C\int w(r)q(r)^2+C\int_r^T\!\!\int\Bigl[\tau e^{\frac{b_0x^2}\tau}|q_x|^2
   +\sqrt\tau w|q_t|^2\Bigr].
\]
By (c) and $2b_A-b_0\ge\max\{b_0,\frac{a_s}2\}$,
$\tau e^{\frac{b_0x^2}\tau}|q_x|^2\le C\delta^2(G_0|\varphi_y|^2+W_c|\varphi|^2)$,
$\sqrt\tau w|q_t|^2\le C\delta^2(\sqrt\tau G|\varphi_t|^2+W_c|\varphi|^2)$ and
$w(r)q(r)^2\le C\delta^2|\varphi(r)|^2$, so
$\mathscr O_{r,T}(h)+\iint wh_x^2\le C\delta^2(\mathscr O_{r,T}(\varphi_t)+J_r(\varphi))$ with $C$
depending only on $C_A$, $b_A$, the fixed Gaussian exponents and the fixed data; this holds first for smooth $\varphi$
with compact support and then by completion. Since
$(Z^0_\sigma)_t=d(q_y+h_y)-s_\sigma(V_c)h$ and $\sqrt\tau G\le Cw$,
\[
 \mathscr O_{r,T}(Z^0_t)\le C_0\delta^2\bigl(\mathscr O_{r,T}(\varphi_t)+J_r(\varphi)\bigr).
\]
The difference $U^{\rm c}-Z^0$ solves
$(L-\widetilde{\mathcal M}\pa_y)(U^{\rm c}-Z^0)=\widetilde{\mathcal M}Z^0_y$ with zero trace, so
Step 4 with $K=\widetilde{\mathcal M}$, $w=Z^0$ gives
$\mathscr O_{r,\infty}((U^{\rm c}-Z^0)_t)\le C_{b,r}\delta^2\mathscr Q_r(\varphi)$. Hence
$\mathscr O_{r,T}(U^{\rm c}_t)\le2C_0\delta^2\mathscr O_{r,T}(\varphi_t)+C_{b,r}\delta^2J_r(\varphi)$.

\emph{Step 6: connection.} Let $Q=U^{\rm a}+U^{\rm b}+U^{\rm c}$, so that
$(L-\widetilde{\mathcal M}\pa_y)Q=f$, $Q(r)=g$, and let $\mathcal KX=\mathcal C_{\widetilde{\mathcal M}}Q$,
$X(r)=0$; then $U=Q+X$ solves $\mathcal KU=f$, $U(r)=g$. By \eqref{e:pcnx},
$\frac{|\mathcal C_{\widetilde{\mathcal M}}Q|^2}{W_c}\le C_b\delta^2(G_0|Q_y|^2+W_c|Q|^2)$, and
Lemma~\ref{l:geo} gives $J_r(X)\le C_b\delta^2J_r(Q)$. Moreover
$(L-\widetilde{\mathcal M}\pa_y)X=\mathcal C_{\widetilde{\mathcal M}}U$; write $X=X^0+X^1$ with zero
traces and sources $\mathcal C^0U$ and $B^1U_y+C^1U$. By \eqref{e:pcnx}, $A=\mathcal C^0$
satisfies (c) with $b_A=\frac{c_G}4\ge\frac{a_s}2$ and a fixed $C_A$, and $B^1U_y+C^1U$
belongs to (b) with $f^{\rm g}=0$ and
$\mathfrak N\le C_b\delta^2J_r(U)$. Step 5 and Step 4 give, on every finite strip, where
$\mathscr O_{r,T}(X_t)<\infty$,
\[
 \mathscr O_{r,T}(X_t)\le C_0'\delta^2\bigl(\mathscr O_{r,T}(Q_t)+\mathscr O_{r,T}(X_t)\bigr)
   +C_{b,r}\delta^2J_r(U),
\]
with $C_0'$ depending only on the exponents of \eqref{e:pwt} and the fixed data. For
$C_0'\delta_1(M)^2\le\frac12$ this gives
$\mathscr O_{r,T}(X_t)\le2C_0'\delta^2\mathscr O_{r,T}(Q_t)+C_{b,r}\delta^2J_r(U)$, and
Steps 3--5 give \eqref{e:solv}; for $g=0$, $f=A\varphi$ only Steps 5 and 6 enter, which
gives the last assertion. Uniqueness and the regularity on finite strips follow from
Lemma~\ref{l:geo}, since each source lies in $L^2((r,T)\times\R)$.
\end{proof}

\subsection{The reference field}\label{s:ref}

Write $\tilde s_\sigma=s_\sigma(V_c)+2k(V_c)a_\sigma$, so that
$\mathscr P_\sigma=\pa_t-d(V_c)\pa_y^2+\tilde s_\sigma\pa_y$, and $\underline s=\min\{c_-,c_+\}$;
by Lemma~\ref{l:cen}, $a_0<\frac{\underline s}{32}$.

\begin{lemma}\label{l:c5}
Let $\sigma=\pm1$, $e\in\{0,s_o\}$, $b_f>0$, $K\ge0$, and let $f$ be a function on
$\R\times[0,\infty)$ with $f,\pa_yf,\pa_y^2f,\pa_tf$ continuous and
\begin{equation}\label{e:csrc}
 |\pa_y^jf|\le K\tau^{-1-\frac j2}e^{-\frac{b_f(y-et)^2}\tau}\quad(0\le j\le2),\qquad
 |(\pa_t+e\pa_y)f|\le K\tau^{-2}e^{-\frac{b_f(y-et)^2}\tau}.
\end{equation}
The problem $\mathscr P_\sigma w=f$, $w(0)=0$, has a unique solution with $w\in C([0,T];H^1)$
and $w_t,w_{yy}\in L^2((0,T)\times\R)$ for finite $T$. There are $t_0$ and $C$, depending only on
$M$, $b_f$ and the fixed data, such that $\sigma(\tilde s_\sigma-e)\ge\frac{\underline s}2$ for
$t\ge t_0$,
\begin{equation}\label{e:cval}
\begin{gathered}
 \sup_{t\ge0}\bigl(\Vert w\Vert_{L^\infty}+\Vert w_y\Vert_{L^\infty}+\Vert w_{yy}\Vert_{L^\infty}
   +\Vert w_t\Vert_{L^\infty}\bigr)\le CK,\qquad
 \Vert w(t)\Vert_{L^\infty}\le CK\tau^{-\frac14}\quad(t\ge t_0),\\
 \tau^{\frac12}|w|+\tau|w_y|+\tau^{\frac32}\bigl(|w_{yy}|+|w_t|\bigr)\le CK
 \qquad(|y|\le a_0t,\ t\ge t_0),
\end{gathered}
\end{equation}
and, for $1<p_w<\frac32$ and $h=w_y-\frac f{\tilde s_\sigma-e}$, with $C_{p_w}$ depending also
on $p_w$,
\begin{equation}\label{e:cout}
 \int_{t_0}^\infty\!\!\int_\R\tau^{-\frac12}(1+z_\sigma^2)^{p_w}|h|^2\dd y\dd t\le C_{p_w}K^2 .
\end{equation}
The source $f^{\rm o}_\sigma=\ell_{\rm g}(V_c)a_o^2$ of \eqref{e:pn} satisfies \eqref{e:csrc}
with $e=s_o$, $K=C_M$ and $b_f=\frac1{8d_o}$.
\end{lemma}
\begin{proof}
Write $D=d(V_c)$; constants $C$, $c_*$ depend only on $M$, $b_f$ and the fixed data, and
$C_T$ also on $T$.

Galerkin approximation with the tests $(-1)^k\pa_y^{2k}w$, $0\le k\le3$, gives a solution
with $\sup_{t\le T}\Vert w\Vert_{H^3}^2+\int_0^T\Vert w\Vert_{H^4}^2\dd t\le
C_T\int_0^T\Vert f\Vert_{H^2}^2\dd t\le C_TK^2$; uniqueness is \eqref{e:uniq}. The derivatives
$v_1=w_y$ and $v_2=w_{yy}$ vanish at $t=0$ and solve
\[
\begin{aligned}
 (v_1)_t-D(v_1)_{yy}+(\tilde s_\sigma-D_y)(v_1)_y+(\tilde s_\sigma)_yv_1&=f_y,\\
 (v_2)_t-D(v_2)_{yy}+(\tilde s_\sigma-2D_y)(v_2)_y+\bigl(2(\tilde s_\sigma)_y-D_{yy}\bigr)v_2
 &=f_{yy}-(\tilde s_\sigma)_{yy}v_1 .
\end{aligned}
\]
Fix $T$ and $\Lambda\ge\max\{8\sup D,\frac2{b_f}\}$; by $(y-et)^2\ge\frac{y^2}2-e^2t^2$,
$|\pa_y^jf|\le C_TKe^{-\frac{y^2}{\Lambda\tau}}$ on $[0,T]$. For
$H=e^{A_Tt-\frac{y^2}{\Lambda\tau}}$ and bounded $\widetilde B,\widetilde q$,
\[
 \frac{(\pa_t-D\pa_y^2+\widetilde B\pa_y+\widetilde q)H}H
 =A_T+\frac{2D}{\Lambda\tau}+\Bigl(\frac1\Lambda-\frac{4D}{\Lambda^2}\Bigr)\frac{y^2}{\tau^2}
   -\frac{2\widetilde By}{\Lambda\tau}+\widetilde q
 \ge A_T-\frac{4\sup\widetilde B^2}\Lambda-\sup|\widetilde q|,
\]
since $\frac1\Lambda-\frac{4D}{\Lambda^2}\ge\frac1{2\Lambda}$ and
$\frac{2|\widetilde By|}{\Lambda\tau}\le\frac{y^2}{4\Lambda\tau^2}+\frac{4\widetilde B^2}\Lambda$.
Take $A_T$ so that the right side is at least $1$ for the three operators above. Comparison
of $\pm w$, then of $\pm v_1$, then of $\pm v_2$ with $C_TKH$ (the source of $v_2$ contains
$v_1$) and the equation for $w_t$ give, with $b_T=\frac1{\Lambda\tau_T}$,
\begin{equation}\label{e:cfin}
 |w|+|w_y|+|w_{yy}|+|w_t|\le C_TKe^{-b_Ty^2}\qquad(0\le t\le T).
\end{equation}

Let $t_1$ be so large that $\tau_{t_1}\ge L^2$, with $L$ of \eqref{e:nfwdef}, and
$2k(V_c)|a_\sigma|\le\frac{\underline s}2$ for $t\ge t_1$; it will be enlarged finitely many
times, depending only on $M$, $b_f$ and the fixed data. Since $\sigma s_\sigma(V_c)\ge\underline s$
and $-\sigma e\ge0$, $\sigma(\tilde s_\sigma-e)\ge\frac{\underline s}2$ for $t\ge t_1$. By the
maximum principle on $\R\times(t_1,T)$ and the third bound of \eqref{e:gauss} with $(v,h)=(e,0)$
and $\frac{\underline s}2$ in place of $b_0$,
\[
 \Vert w(t)\Vert_{L^\infty}\le\Vert w(t_1)\Vert_{L^\infty}
 +\Bigl\Vert T_{\mathscr P_\sigma,t_1}\bigl(K\tau^{-1}e^{-\frac{b_f(y-et)^2}\tau}\bigr)\Bigr\Vert_{L^\infty}
 \le CK\qquad(t\ge t_1),
\]
where \eqref{e:cfin} with $T=t_1$ was used. Lemma~\ref{l:c4}(ii) with $r=t_1$ and $A=CK$ gives
\begin{equation}\label{e:cdec}
 \Vert\pa_y^jw(t)\Vert_{L^\infty}\le CK\tau^{-\frac j2}\qquad(1\le j\le3,\ t\ge2t_1+1),
\end{equation}
and with \eqref{e:cfin} for $T=2t_1+1$ and the equation this is the first bound of
\eqref{e:cval}. By \eqref{e:ck}, $\Vert f(t)\Vert_{L^1}\le CK\tau^{-\frac12}$ and \eqref{e:cfin}
for $T=t_1$, for $t\ge t_1$,
\[
 \Vert w(t)\Vert_{L^1}\le\Vert S_{\mathscr P_\sigma}(t,t_1)w(t_1)\Vert_{L^1}
 +\int_{t_1}^t\Vert S_{\mathscr P_\sigma}(t,s)f(s)\Vert_{L^1}\dd s\le CK\sqrt\tau .
\]
For $\tau\ge\max\{4,\tau_{t_1}^2\}$ let $\ell=\sqrt\tau-1\ge t_1$, so that $t-\ell\ge\frac\tau2$
and $\tau_\ell=\sqrt\tau$. The response to $f$ on $(\ell,t)$ is bounded by
$CK\tau_\ell^{-\frac12}$ by the third bound of \eqref{e:gauss} on $\R\times(\ell,\infty)$, and
\eqref{e:ck} gives
\[
 \Vert w(t)\Vert_{L^\infty}\le\Vert S_{\mathscr P_\sigma}(t,\ell)w(\ell)\Vert_{L^\infty}
 +CK\tau_\ell^{-\frac12}
 \le CK\bigl[(t-\ell)^{-\frac12}\tau_\ell^{\frac12}+\tau_\ell^{-\frac12}\bigr]\le CK\tau^{-\frac14};
\]
for $t\ge t_1$ with $\tau<\max\{4,\tau_{t_1}^2\}$ this follows from $\Vert w(t)\Vert_{L^\infty}\le CK$.

For the central bounds put $x=\sigma y$, $e'=\sigma e\le0$, $\zeta=\frac{x-e't}{\sqrt\tau}$, and
write, as functions of $(x,t)$, $D=d(V_c)$, $B=\sigma\tilde s_\sigma$, $\mathsf s=\sigma s_\sigma(V_c)$
and $q=\frac f{B-e'}$; then $w_x-q=\sigma h$. Since $w_t-Dw_{xx}+Bw_x=(B-e')q$,
\eqref{e:cnr} with $(z,D,B,e,f_2)=(w,D,B,e',0)$ gives
\[
 (\sigma h)_t-\bigl(D(\sigma h)_x\bigr)_x+(B\sigma h)_x=-R,\qquad R=(q_t+e'q_x)-(Dq_x)_x .
\]
On $\{x\le2a_0t\}$ the function $B-\mathsf s=2\sigma k(V_c)a_\sigma$ and its derivatives of order
at most two are bounded by $Ce^{-c_*t}$, since $a_\sigma$ is centred at $x=c_\sigma t$; moreover
$(\pa_t+e'\pa_x)f=(\pa_t+e\pa_y)f$, $|\mathsf s_x|+\tau^{\frac12}|\mathsf s_{xx}|\le C\delta\tau^{-\frac12}$,
$|\mathsf s_t|\le C\delta\tau^{-1}$, and, if $e'<0$,
$|\mathsf s_x||f|\le CK\tau^{-\frac32}e^{-c_*t}e^{-\frac{b_f}2\zeta^2}$, the centres $0$ and $e't$
being separated. With \eqref{e:csrc} this gives, for $b=\frac{b_f}2$ and $t\ge t_1$,
\[
 |q|+\tau^{\frac12}|q_x|+\tau|q_{xx}|\le CK\tau^{-1}e^{-b\zeta^2},\qquad
 |R|\le CK\tau^{-2}e^{-b\zeta^2}\qquad(x\le2a_0t).
\]

Let $g=g_{t_1}$ of Lemma~\ref{l:nf}(iii), as a function of $x$ (it is even), $K_0\ge1$, and
\[
 m=\mathsf s^{-1}\exp\Bigl(K_0\int_{-\infty}^xg(x')\dd x'\Bigr).
\]
By Lemma~\ref{l:nf}(iii), $C^{-1}\le m\le Ce^{\pi C_1K_0\delta}$ and $(\mathsf sm)_x=K_0g\,\mathsf sm$;
with $|\mathsf s_t|+|\mathsf s_{xx}|+\mathsf s_x^2+|D_x\mathsf s_x|\le Cg$,
$|D_x|+|\mathsf s_x|\le C\delta\tau^{-\frac12}$, $|g'|\le\tau_{t_1}^{-\frac12}g$ and
$g^2\le C_1\delta\tau_{t_1}^{-1}g$,
\[
 \frac{m_t-(Dm_x)_x+(\mathsf sm)_x}m=\mathsf sK_0g+R_m,\qquad
 |R_m|\le C_2\bigl(1+K_0\tau_{t_1}^{-\frac12}+K_0^2\delta\tau_{t_1}^{-1}\bigr)g,
\]
with $C_2$ depending only on the fixed data. Fix $K_0=\frac{3C_2}{\underline s}$ and enlarge $t_1$ so
that $K_0\tau_{t_1}^{-\frac12}+K_0^2\delta\tau_{t_1}^{-1}\le2$; then
$\frac{m_t-(Dm_x)_x+(\mathsf sm)_x}m\ge0$. Put
$\widetilde B=B-D_x-2D\frac{m_x}m$, $\epsilon(t)=\sup_{x\le2a_0t}|\frac{((B-\mathsf s)m)_x}m|\le Ce^{-c_*t}$
and $E(t)=\exp\int_{t_1}^t\epsilon$. Then $\tilde u=\frac{\sigma h}{mE}$ solves
\[
 \tilde u_t-D\tilde u_{xx}+\widetilde B\tilde u_x+\widetilde C\tilde u=-\frac R{mE},\qquad
 \widetilde C=\frac{m_t-(Dm_x)_x+(Bm)_x}m+\epsilon\ge0,
\]
and, since $|\frac{m_x}m|\le K_0C_1\delta+C\delta$, $\widetilde B\ge\frac34\underline s$ on
$\{x\le2a_0t,\ t\ge t_1\}$ for $\delta$ small and $t_1$ large.

Let $\lambda\ge\sup D$ with $\frac1{4\lambda}<\frac b2$, $b'=\frac{\underline s}8$,
$\mathcal D=\tau-\frac x{b'}$, $\Phi=F_\lambda$ of \eqref{e:CH}, $\mathcal W=\mathcal D^{-\frac32}\Phi(\zeta)$
and $P_m=\pa_t-D\pa_x^2+\widetilde B\pa_x$. From $\Phi''=-\frac\zeta{2\lambda}\Phi'$,
$\mathcal D_t=1$, $\mathcal D_x=-\frac1{b'}$, $\zeta_t=-\frac{e'}{\sqrt\tau}-\frac\zeta{2\tau}$ and
$\zeta_x=\tau^{-\frac12}$,
\[
\begin{aligned}
 P_m\mathcal W={}&\mathcal D^{-\frac32}\Phi'\Bigl[\frac{\widetilde B-e'}{\sqrt\tau}
   +\Bigl(\frac D{2\lambda}-\frac12\Bigr)\frac\zeta\tau\Bigr]
 +\frac32\mathcal D^{-\frac52}\Bigl(\frac{\widetilde B}{b'}-1\Bigr)\Phi\\
 &-\frac{15D}{4b'^2}\mathcal D^{-\frac72}\Phi-\frac{3D}{b'\sqrt\tau}\mathcal D^{-\frac52}\Phi' .
\end{aligned}
\]
On $\{x\le2a_0t\}$, $\mathcal D\ge\frac\tau2$ (as $\frac{2a_0}{b'}<\frac12$),
$\frac{\widetilde B}{b'}\ge6$, $(\frac D{2\lambda}-\frac12)\zeta\ge0$ for $\zeta\le0$, and
$\frac\zeta{2\tau}\le\frac{a_0+|e'|/2}{\sqrt\tau}$ for $\zeta>0$. Enlarge $t_1$ so that
$\tau_{t_1}\ge\frac{24\sup D}{b'\underline s}+\frac{\sup D}{b'^2}$; then, for $t\ge t_1$,
\[
 P_m\mathcal W\ge\frac{\widetilde B-a_0}{2\sqrt\tau}\mathcal D^{-\frac32}\Phi'(\zeta)
 \ge\frac{\underline s}{4\sqrt\tau}\mathcal D^{-\frac32}\Phi'(\zeta)\ge c_3\tau^{-2}e^{-b\zeta^2},
\]
where the last step uses $\frac{\mathcal D}\tau\le1+\frac{|e'|+|\zeta|}{b'}$ (from
$|x|\le(|e'|+|\zeta|)\tau$) and $\frac1{4\lambda}<\frac b2$, and $c_3>0$ depends also on
$\lambda$. Moreover, for $\varepsilon=\frac{\underline s}{2\sup D}$,
$P_me^{\varepsilon(x-2a_0t)}=\varepsilon(\widetilde B-2a_0-\varepsilon D)e^{\varepsilon(x-2a_0t)}\ge0$.
Enlarge $\lambda$ so that $\frac1{4\lambda\tau_{t_1}}<b_{t_1}$, with $b_{t_1}$ of \eqref{e:cfin}
for $T=t_1$; since $F_\lambda(\zeta)\ge\frac{c_*}{1+|\zeta|}e^{-\frac{\zeta^2}{4\lambda}}$ for $\zeta\le0$,
\eqref{e:cfin} gives $|h(\cdot,t_1)|\le CK\mathcal W(\cdot,t_1)$ on $x\le2a_0t_1$, and on
$x=2a_0t$, $|h|\le CK$ by the first bound of \eqref{e:cval}. Comparison for $\tilde u$ on
$\{-n<x<2a_0t,\ t_1<t<T\}$, with $\sup_{t_1\le t\le T}|h(-n,t)|\to0$ as $n\to\infty$, gives
\begin{equation}\label{e:cinc}
 |h|\le CK\Bigl[\mathcal D^{-\frac32}F_\lambda(\zeta)+e^{\varepsilon(x-2a_0t)}\Bigr]
 \qquad(x\le2a_0t,\ t\ge t_1).
\end{equation}

On $x\le\frac32a_0t$, $\mathcal D\ge\frac\tau2$ and $e^{\varepsilon(x-2a_0t)}\le e^{-\frac12\varepsilon a_0t}$,
so $|h|\le CK\tau^{-\frac32}$ there. Since $w(\cdot,t)$ vanishes at $-\infty$,
$\int_{-\infty}^{a_0t}\mathcal D^{-\frac32}\dd x\le2b'(\frac\tau2)^{-\frac12}$,
$\int_{-\infty}^{a_0t}e^{\varepsilon(x-2a_0t)}\dd x=\frac1\varepsilon e^{-\varepsilon a_0t}$ and
$\int_\R|q|\dd x\le CK\tau^{-\frac12}$ give $|w|\le CK\tau^{-\frac12}$ on $x\le a_0t$, and
$|w_x|\le|q|+|h|\le CK\tau^{-1}$ on $x\le\frac32a_0t$. By Taylor's formula,
$|h_x(x)|\le|h(x)|+|h(x+1)|+\frac12\sup_{[x,x+1]}|h_{xx}|$, where
$|h_{xx}|\le|w_{xxx}|+|q_{xx}|\le CK\tau^{-\frac32}$ by \eqref{e:cdec}; hence
$|h_x|\le CK\tau^{-\frac32}$ on $x\le a_0t$ for $a_0t\ge2$. Finally $w_{xx}=q_x+(\sigma h)_x$ and
\[
 w_t=Dw_{xx}-Bw_x+(B-e')q=Dq_x+D(\sigma h)_x-B\sigma h-e'q,
\]
where $e'q=0$ for $e'=0$ and $|q|\le CKe^{-c_*t}$ on $|x|\le a_0t$ for $e'<0$. This proves the
last line of \eqref{e:cval} for $t\ge t_0$, with $t_0$ the largest of the thresholds used.

For $f^{\rm o}_\sigma$: by Lemma~\ref{l:pr},
$(\pa_t+s_o\pa_y)a_o=d_o(a_o)_{yy}-2k_oa_o(a_o)_y$ (the equation of $\vartheta_o$ times $s_o$) and
$|\pa_y^ja_o|\le C_M\tau^{-\frac{j+1}2}e^{-\frac{z_o^2}{8d_o}}$, so $a_o^2$ satisfies
\eqref{e:csrc} with $e=s_o$, $K=C_M$, $b_f=\frac1{4d_o}$. The derivatives of $\ell_{\rm g}(V_c)$
of orders $1,2$ and $(\pa_t+s_o\pa_y)\ell_{\rm g}(V_c)$ are bounded by
$C\delta\tau^{-\frac12}e^{-\frac{c_Gy^2}{2\tau}}$, and
$e^{-\frac{c_Gy^2}{2\tau}-\frac{z_o^2}{8d_o}}\le e^{-c_*t}$ by $y^2+(y-s_ot)^2\ge\frac{s_o^2t^2}2$;
this gives the assertion for $f^{\rm o}_\sigma$.

It remains to prove \eqref{e:cout}; since $s_\sigma=\sigma c_\sigma$, in $x$,
$1+z_\sigma^2=1+\frac{\eta^2}\tau$ with $\eta=x-c_\sigma t$. Let $t\ge t_0$. On $-(|e'|+1)\tau\le x\le\frac{a_0}2t$,
$|h|^2\le CK^2\tau^{-3}+CK^2e^{-c_*t}$ and $1+\frac{\eta^2}\tau\le C\tau$; on $x<-(|e'|+1)\tau$,
$\zeta\le-\sqrt\tau$, $1+\frac{\eta^2}\tau\le C(1+\zeta^2)$ and $F_\lambda(\zeta)\le e^{-\frac{\zeta^2}{4\lambda}}$.
Thus \eqref{e:cinc} gives
\[
 \int_{-\infty}^{\frac{a_0}2t}\Bigl(1+\frac{\eta^2}\tau\Bigr)^{p_w}|h|^2\dd x
 \le CK^2\bigl(\tau^{p_w-2}+e^{-c_*t}\bigr).
\]
Let $v=p_\sigma w$, so that $L_\sigma v=p_\sigma f-\alpha R^{\rm prof}_\sigma v$ by Lemma~\ref{l:pr},
and let $\chi\in C^\infty(\R)$ with $\chi=0$ on $(-\infty,\frac{a_0}4]$ and $\chi=1$ on
$[\frac{a_0}2,\infty)$. With $\chi=\chi(\frac xt)$ and $L_\sigma=\pa_t-D\pa_x^2+\mathsf s\pa_x$
in $x$,
\[
 L_\sigma(\chi v)=\chi\bigl(p_\sigma f-\alpha R^{\rm prof}_\sigma v\bigr)
 +(\chi_t-D\chi_{xx}+\mathsf s\chi_x)v-2D\chi_xv_x .
\]
On $x\ge\frac{a_0}4t$, $|f|\le CKe^{-c_*t}$, and $|R^{\rm prof}_\sigma|\le C\delta e^{-c_et}$; the
derivatives of $\chi$ are supported in $\frac{a_0}4t\le x\le\frac{a_0}2t$ and bounded by $Ct^{-1}$,
and there $|v|+\tau^{\frac12}|v_x|\le CK\tau^{-\frac12}$ by \eqref{e:cval} and
$|a_\sigma|\le Ce^{-c_*t}$. Hence $|L_\sigma(\chi v)|\le CK\tau^{-\frac32}$ for $t\ge t_0$, and
\eqref{e:cdh} for $L_\sigma$ on $[t-\frac\varpi2\tau_t,t]$, with
$\Vert\chi v\Vert_{L^\infty}\le CK\tau^{-\frac14}$, gives, for $\varphi=v_x$ and $t\ge t_3=2t_0+1$,
\[
 |\varphi|\le CK\bigl[\tau^{-\frac12}\tau^{-\frac14}+\tau^{\frac12}\tau^{-\frac32}\bigr]\le CK\tau^{-\frac34}
 \qquad\bigl(x\ge\tfrac{a_0}2t\bigr).
\]
On $x>\frac{a_0}2t$, differentiation of the equation of $v$ gives
\[
\begin{gathered}
 \varphi_t-d_\sigma\varphi_{xx}+c_\sigma\varphi_x=E,\\
 E=(D-d_\sigma)\varphi_{xx}+(D_x-\mathsf s+c_\sigma)\varphi_x-\mathsf s_x\varphi
 +\pa_x\bigl(p_\sigma f-\alpha R^{\rm prof}_\sigma v\bigr),
\end{gathered}
\]
and, by \eqref{eq:ct-gau}, \eqref{e:cval}, \eqref{e:cdec}, \eqref{e:csrc} and Lemma~\ref{l:pr},
$|E|\le CKe^{-c_*(t+\frac{x^2}\tau)}$ there for $t\ge t_3$. Moreover
$\varphi=p_\sigma(\sigma h+q-\sigma\alpha a_\sigma w)$ gives $|\varphi(\frac{a_0}2t,t)|\le CK\tau^{-\frac32}$,
and \eqref{e:cfin} gives $|\varphi(x,t_3)|\le CKe^{-c_*x^2}$.

Put $\bar c=c_\sigma-\frac{a_0}2$, $\xi=x-\frac{a_0}2t$,
$\phi_u(X)=(4\pi d_\sigma u)^{-\frac12}e^{-\frac{X^2}{4d_\sigma u}}$, and, for $\xi>0$, $u>0$,
\[
 \Pi(\xi,u)=F_{d_\sigma}\Bigl(\frac{\bar cu-\xi}{\sqrt u}\Bigr)
 +e^{\frac{\bar c\xi}{d_\sigma}}F_{d_\sigma}\Bigl(-\frac{\xi+\bar cu}{\sqrt u}\Bigr),\qquad
 P(\xi,u)=\pa_u\Pi(\xi,u)=\frac\xi u\,\phi_u(\xi-\bar cu),
\]
where $P$ is computed with $e^{\frac{\bar c\xi}{d_\sigma}}\phi_u(\xi+\bar cu)=\phi_u(\xi-\bar cu)$.
The first term of $\Pi$ solves $\Pi_u-d_\sigma\Pi_{\xi\xi}+\bar c\Pi_\xi=0$, and the second is
$e^{\frac{\bar c\xi}{d_\sigma}}$ times a solution of the same equation with $-\bar c$ in place of
$\bar c$, hence also a solution; $\Pi(0,u)=1$, $\Pi(\xi,0^+)=0$ and $\Pi(\xi,\infty)=1$. Thus
$\int_0^\infty P(\xi,u)\dd u=1$, $P\ge0$, and, for $\beta\in C^1([t_3,\infty))$,
\[
 \int_{t_3}^tP(\xi,t-s)\beta(s)\dd s=\beta(t_3)\Pi(\xi,t-t_3)+\int_{t_3}^t\Pi(\xi,t-s)\beta'(s)\dd s
\]
solves $\pa_t-d_\sigma\pa_\xi^2+\bar c\pa_\xi=0$, which is $\pa_t-d_\sigma\pa_x^2+c_\sigma\pa_x$ in
$(x,t)$, with the value $\beta(t)$ at $\xi=0$ and $0$ at $t=t_3$.
Let $\Theta_2$ be this function for $\beta(s)=CK\tau_s^{-\frac32}$, and
\[
 \Theta_1=\int_\R\phi_{t-t_3}\bigl(x-x'-c_\sigma(t-t_3)\bigr)\bigl|\varphi(x',t_3)\bigr|\dd x',\qquad
 \Theta_3=\int_{t_3}^t\!\!\int_\R\phi_{t-s}\bigl(x-x'-c_\sigma(t-s)\bigr)\bigl|E(x',s)\bigr|\dd x'\dd s,
\]
with $\varphi(\cdot,t_3)$ and $E(\cdot,s)$ extended by zero to $x'<\frac{a_0}2t_3$ and
$x'<\frac{a_0}2s$. Then $\Theta=\Theta_1+\Theta_2+\Theta_3\ge0$,
$(\pa_t-d_\sigma\pa_x^2+c_\sigma\pa_x)\Theta=|E|$ on $\{x>\frac{a_0}2t\}$, and $\Theta\ge|\varphi|$
at $t=t_3$ and at $x=\frac{a_0}2t$; comparison on $\{\frac{a_0}2t<x<n\}$, with
$\sup_{t_3\le t\le T}|\varphi(n,t)|\to0$, gives $|\varphi|\le\Theta$.

Since $\xi-\bar c(t-s)=\eta+\bar cs$ and $t-s\le\tau$,
$P(\xi,t-s)\le\sqrt2P_2(\xi,t-s)e^{-\frac{(\eta+\bar cs)^2}{8d_\sigma\tau}}$, where $P_2$ is $P$
with $2d_\sigma$ in place of $d_\sigma$, and $\int_0^\infty P_2\dd u=1$. If $\eta\ge0$, or
$s\le\sqrt\tau$, or $\eta<0$ and $s\le\frac{|\eta|}{2\bar c}$, then
$(\eta+\bar cs)^2\ge\frac{\eta^2}4-\bar c^2\tau$; otherwise
$\tau_s^{-\frac32}\le C(\sqrt\tau+|\eta|)^{-\frac32}$. Hence
$\Theta_2\le CK[e^{-\frac{\eta^2}{C\tau}}+(\sqrt\tau+|\eta|)^{-\frac32}\mathbf 1_{\{\eta<0\}}]$.
Completing squares with $\eta=(x-x'-c_\sigma(t-s))+(x'-c_\sigma s)$ and
$\frac{(x'-c_\sigma s)^2}\tau\le C(\frac{x'^2}{\tau_s}+s)$ gives $\Theta_1+\Theta_3\le CKe^{-\frac{\eta^2}{C\tau}}$.
With $|\varphi|\le CK\tau^{-\frac34}$, on $x\ge\frac{a_0}2t$, $t\ge t_3$,
\[
 |\varphi|\le CK\Bigl[\min\bigl\{\tau^{-\frac34},e^{-\frac{\eta^2}{C\tau}}\bigr\}
 +(\sqrt\tau+|\eta|)^{-\frac32}\mathbf 1_{\{\eta<0\}}\Bigr].
\]
Splitting at $|\eta|=\sqrt{2C\tau\log(2+\tau)}$ and using $p_w>1$ and $\dd x=\dd\eta$,
\[
\begin{aligned}
 \int_{\frac{a_0}2t}^\infty\Bigl(1+\frac{\eta^2}\tau\Bigr)^{p_w}|\varphi|^2\dd x
 &\le CK^2\Bigl[\frac{\log(2+\tau)^{p_w+\frac12}}\tau
 +\int_{-c_\sigma t}^0\frac{(\tau+\eta^2)^{p_w}}{\tau^{p_w}(\sqrt\tau+|\eta|)^3}\dd\eta\Bigr]\\
 &\le CK^2\Bigl[\frac{\log(2+\tau)^{p_w+\frac12}}\tau+\tau^{p_w-2}\Bigr].
\end{aligned}
\]
On $x\ge\frac{a_0}2t$, $\sigma h=p_\sigma^{-1}\varphi+\sigma\alpha a_\sigma w-q$, with
$|a_\sigma w|\le CK\tau^{-\frac34}e^{-\frac{\eta^2}{8d_\sigma\tau}}$ and
$|q|\le CK\tau^{-1}e^{-c_*(t+\frac{x^2}\tau)}$; so $h$ satisfies the same bound. Together with
the incoming bound, \eqref{e:cfin} on $[t_0,t_3]$, and $p_w<\frac32$, multiplication by
$\tau^{-\frac12}$ and integration in $t$ give \eqref{e:cout}.
\end{proof}

\begin{lemma}\label{l:ref}
Let $f^{\rm o}$, $f^{\rm s}$ be those of Lemma~\ref{l:peq} and $\mathcal C^0$ that of
\eqref{e:pcnx}. The system
\[
 \mathscr P_\sigma Z^{\rm ref}_\sigma=f^{\rm o}_\sigma+f^{\rm s}_\sigma
 +\sum_{\rho=\pm1}\mathcal C^0_{\sigma\rho}Z^{\rm ref}_\rho\quad(\sigma=\pm1),\qquad Z^{\rm ref}(0)=0,
\]
has a unique solution with $Z^{\rm ref}\in C([0,T];H^1)$ and
$Z^{\rm ref}_t,Z^{\rm ref}_{yy}\in L^2((0,T)\times\R)$ for finite $T$. There are $t_0$ and $C$,
depending only on $M$ and the fixed data, such that, with $Z=Z^{\rm ref}$,
\[
\begin{gathered}
 \sup_{t\ge0}\bigl(\Vert Z\Vert_{L^\infty}+\Vert Z_y\Vert_{L^\infty}+\Vert Z_{yy}\Vert_{L^\infty}
   +\Vert Z_t\Vert_{L^\infty}\bigr)\le C,\\
 \tau^{\frac12}|Z|+\tau|Z_y|+\tau^{\frac32}\bigl(|Z_{yy}|+|Z_t|\bigr)\le C
 \qquad(|y|\le a_0t,\ t\ge t_0),
\end{gathered}
\]
$\sigma\tilde s_\sigma\ge\frac{\underline s}2$ and $\sigma(\tilde s_\sigma-s_o)\ge\frac{\underline s}2$
for $t\ge t_0$, $f^{\rm c}=f^{\rm s}+\mathcal C^0Z$ satisfies \eqref{e:csrc} with $e=0$, $K=C$ and
$b_f=\frac{c_G}4$, $\mathscr O_{t_0,\infty}(Z_t)<\infty$, and, for $1<p_w<\frac32$,
\[
 h_\sigma=Z_{\sigma,y}-\frac{f^{\rm o}_\sigma}{\tilde s_\sigma-s_o}-\frac{f^{\rm c}_\sigma}{\tilde s_\sigma}
 \quad\text{satisfies}\quad
 \sum_{\sigma=\pm1}\int_{t_0}^\infty\!\!\int_\R\tau^{-\frac12}(1+z_\sigma^2)^{p_w}|h_\sigma|^2\dd y\dd t
 \le C_{p_w}.
\]
\end{lemma}
\begin{proof}
Constants $C$ depend only on $M$ and the fixed data. By \eqref{e:pcnx}, $\mathcal C^0$ is a smooth
function of $V_c$ times $(V_c)_y$, so the bounds for the derivatives of $V_c$ stated before
\eqref{e:pcnx} give $|\pa_y^j\mathcal C^0|\le C\delta\tau^{-\frac{j+1}2}e^{-\frac{c_Gy^2}{2\tau}}$
$(0\le j\le2)$ and $|\mathcal C^0_t|\le C\delta\tau^{-\frac32}e^{-\frac{c_Gy^2}{2\tau}}$. Since
$\pa_y^jf^{\rm s}_\sigma=\tau^{-1-\frac j2}F^{(j)}_\sigma(\xi)$ and
$\pa_tf^{\rm s}_\sigma=-\tau^{-2}(F_\sigma+\frac\xi2F_\sigma')(\xi)$ with $\xi=\frac y{\sqrt\tau}$,
$f^{\rm s}$ satisfies \eqref{e:csrc} with $e=0$, $K=C\delta$, $b_f=\frac{c_G}4$, and
$f^{\rm o}$ does with $e=s_o$, $K=C_M$, $b_f=\frac1{8d_o}$ by Lemma~\ref{l:c5}. Let $t_0$ be the
largest of the numbers $t_0$ of Lemma~\ref{l:c5} for $\sigma=\pm1$ and these exponents, and put
\[
 \Vert Z\Vert_X=\sup_{t\ge0}\bigl(\Vert Z\Vert_{L^\infty}+\Vert Z_y\Vert_{L^\infty}
   +\Vert Z_{yy}\Vert_{L^\infty}+\Vert Z_t\Vert_{L^\infty}\bigr)
 +\sup_{t\ge t_0,\,|y|\le a_0t}\bigl[\tau^{\frac12}|Z|+\tau|Z_y|+\tau^{\frac32}(|Z_{yy}|+|Z_t|)\bigr].
\]
If $\Vert Z\Vert_X<\infty$, then $\mathcal C^0Z$ satisfies \eqref{e:csrc} with $e=0$,
$b_f=\frac{c_G}4$ and $K=C\delta\Vert Z\Vert_X$: on $|y|\le a_0t$, $t\ge t_0$, every term of
$\pa_y^j(\mathcal C^0Z)$ $(j\le2)$ and of $\pa_t(\mathcal C^0Z)$ is a derivative of $\mathcal C^0$
times a derivative of $Z$ with the decay recorded in $\Vert Z\Vert_X$, and elsewhere
$e^{-\frac{c_Gy^2}{4\tau}}\le C\tau^{-2}$.

By linearity and \eqref{e:cval}, solving the scalar equations with sources
$f^{\rm o}_\sigma+f^{\rm s}_\sigma+(\mathcal C^0Z)_\sigma$ and zero initial values
defines a map on $X$, with value at zero bounded by $C$ and Lipschitz constant
$C_5\delta$. For $C_5\delta\le\frac12$, the contraction mapping theorem gives a fixed
point with $\Vert Z\Vert_X\le2C$; Lemma~\ref{l:c5}, applied to this fixed source, gives
its stated regularity on finite time intervals. Uniqueness is \eqref{e:uniq} with
$B_1=0$, $B_2=-\operatorname{diag}(\tilde s_+,\tilde s_-)-\pa_y[d(V_c)]I$ and $B_0=\mathcal C^0$. Since $f^{\rm c}$ is now a fixed source, uniqueness in Lemma~\ref{l:c5}
shows that $Z_\sigma$ is the sum of the solutions of Lemma~\ref{l:c5} with the sources
$f^{\rm o}_\sigma$ and $f^{\rm c}_\sigma$, and \eqref{e:cout} for both gives the bound for
$h_\sigma$. To estimate $\mathscr O_{t_0,\infty}(Z_t)$, let $t\ge t_0$ with
$a_0t\ge2|h_*|$. The central bound on $|y|\le a_0t$, the
boundedness of $Z_t$ and $G\le e^{-\frac{a_sy^2}{8\tau}}$ on $|y|>a_0t$ give
$A_{Z_t}(t)\le C\tau^{-2}+C\tau e^{-c_*t}$, with $c_*>0$ depending also on $a_s$; hence
$\mathscr O_{t_0,\infty}(Z_t)<\infty$.
\end{proof}

\subsection{The homogeneous field and the pressure}\label{s:hom}

Fix $1<p_w<\frac32$ and $r_0\ge t_0$, with $t_0$ of Lemma~\ref{l:ref}, so late that
$\tau_{r_0}\ge h_*^2$ and Lemmas~\ref{l:am}--\ref{l:solv} hold with $r=r_0$ (Lemma~\ref{l:bulk}
with $a=a_s$); write $u=u_{r_0}$, $W=W_{r_0}$ and $\mathcal E=\mathcal E_{r_0}$.

\begin{lemma}\label{l:for}
The problem
\begin{equation}\label{e:frc}
 \mathcal K\mathcal V=\mathfrak f+\mathcal C_{\widetilde{\mathcal M}}u,\qquad
 \mathcal V(r_0)=-\mathcal P(V_c,z)(r_0),
\end{equation}
has a unique solution with $\mathcal V\in C([r_0,T];H^1)$ and
$\mathcal V_t,\mathcal V_{yy}\in L^2((r_0,T)\times\R)$ for finite $T$, and
$\mathscr O_{r_0,\infty}(\mathcal V_t)<\infty$. Put
\begin{equation}\label{e:hqe}
\begin{gathered}
 \mathcal H=\mathcal E-\mathcal V,\qquad \psi_h=S^{-1}\mathcal H,\qquad \mathfrak q_h=\pa_y\psi_h;\\
 \text{then}\quad \mathcal K\mathcal H=0,\qquad \mathcal H_++\mathcal H_-=\psi_{h,1},\qquad
 (\psi_h)_t-d(V_c)(\psi_h)_{yy}+B(\psi_h)_y=0,\\
 (\mathfrak q_h)_t=\pa_y\bigl[d(V_c)(\mathfrak q_h)_y-B\mathfrak q_h\bigr],
\end{gathered}
\end{equation}
with $B$ of \eqref{e:pB}. Moreover
$\mathfrak q_h(r_0)=\pa_y[S^{-1}E^{\rm base}(r_0)]\in L^1\cap L^2$; for finite $T$,
$\mathfrak q_h\in C([r_0,T];L^1\cap L^2)$ and $(\mathfrak q_h)_y\in L^2((r_0,T)\times\R)$, and
$\mathfrak q_h$ is the only solution of the last equation in \eqref{e:hqe} in this class with its
value at $r_0$. In particular, for $r\ge r_0$, $\mathfrak q_h$ is on $[r,\infty)$ the solution in
this class with the value $\mathfrak q_h(r)\in L^1\cap L^2$ at $r$.
\end{lemma}
\begin{proof}
Write $Z=Z^{\rm ref}$, $d=d(V_c)$, $\mathcal A=\widetilde{\mathcal M}+2k(V_c)\operatorname{diag}(a_+,a_-)$,
$q^{\rm o}_\sigma=\frac{f^{\rm o}_\sigma}{\tilde s_\sigma-s_o}$ and
$q^{\rm c}_\sigma=\frac{f^{\rm c}_\sigma}{\tilde s_\sigma}$. Since
$\mathcal K=L-\widetilde{\mathcal M}\pa_y-\mathcal C_{\widetilde{\mathcal M}}$ by \eqref{e:gop} and
\eqref{e:pB}, $\mathscr P_\sigma=L_\sigma+2k(V_c)a_\sigma\pa_y$ and
$\mathcal C_{\widetilde{\mathcal M}}=\mathcal C^0+B^1\pa_y+C^1$ by \eqref{e:pcnx},
Lemma~\ref{l:ref} gives
\[
 \mathcal KZ=f^{\rm o}+f^{\rm s}-g_Z,\qquad g_Z=\mathcal AZ_y+B^1Z_y+C^1Z,\qquad
 Z_y=q^{\rm o}+q^{\rm c}+h .
\]
Let $U^{\rm f}$ be the solution of Lemma~\ref{l:solv} with $r=r_0$ of
\begin{equation}\label{e:frU}
 \mathcal KU^{\rm f}=f^\star,\qquad f^\star=f^{\rm n}+g+g_Z+\mathcal C_{\widetilde{\mathcal M}}u,\qquad
 U^{\rm f}(r_0)=-\mathcal P(V_c,z)(r_0)-Z(r_0).
\end{equation}
The terms of $f^\star$ lie in the classes of Lemma~\ref{l:solv}. Class (a):
$f^{\rm n}_\sigma=\sum_ja_jC^{\rm n}_{\sigma j}n_\chi$, by Lemma~\ref{l:emi} and
$(1+z_j^2)^{p_w}\le Ce^{b_{\rm ac}z_j^2}$, and $(\mathcal Ah)_\sigma=\sum_j\mathcal A_{\sigma j}h_j$,
by Lemma~\ref{l:ref} and $\sup\tau|\mathcal A|^2<\infty$ (\eqref{e:pmat} and Lemma~\ref{l:pr}).
Class (b): $\widetilde{\mathcal M}W_y$, $K_cz_y$ and $B^1u_y$ are terms $K_i\pa_yw_i$, by
Lemma~\ref{l:peq}, \eqref{e:nfe}, \eqref{e:pmat} and
$|K_c|+|B^1|\le C_b\tau^{-\frac12}G_0^{\frac12}$. The terms $g_c$, $\mathcal Aq^{\rm o}$,
$\mathcal Aq^{\rm c}$ and $B^1Z_y+C^1Z$ belong to $f^{\rm g}$: $\Vert g_c\Vert_{L^2}\le C_b\tau^{-\frac54}$
by Lemma~\ref{l:peq}, and each of the other three is bounded by
$C_b\tau^{-\frac32}e^{-\frac{(y-y_0)^2}{C\tau}}$ with $y_0=s_ot$, $0$, $-h_*$ respectively (for
$B^1Z_y+C^1Z$ by the central bounds of Lemma~\ref{l:ref} and \eqref{e:pcnx}), so that its
$L^2$ norm is at most $C_b\tau^{-\frac54}$. The term $C^1u$ belongs to $f^{\rm w}$, since
$|C^1|\le C_b\delta\tau^{-1}e^{-\frac{b_1(y+h_*)^2}\tau}$ and $b_1\ge b_0$ give
$(\frac1{W_c}+\tau)|C^1|^2\le C_bW_c$, and $J_{r_0}(u)<\infty$. Class (c):
$\mathcal C^0u$, with $A=\mathcal C^0$, $b_A=b_1=\frac{c_G}4\ge\frac{a_s}2$ and $\varphi=u$, by
Lemma~\ref{l:bcor}. The trace lies in $H^1$ by Lemma~\ref{l:ref}, $|\mathcal P(V_c,z)|\le C|z|$ and
$|\pa_y\mathcal P(V_c,z)|\le C(|z_y|+|(V_c)_y||z|)$. Hence \eqref{e:solv}, with
$\mathscr O_{r_0,\infty}(u_t)<\infty$ from Lemma~\ref{l:bcor}, gives
$J_{r_0}(U^{\rm f})+\mathscr O_{r_0,\infty}(U^{\rm f}_t)<\infty$.

By Lemma~\ref{l:peq}, $f^{\rm n}+f^{\rm o}+f^{\rm s}+g=\mathfrak f$, so $\mathcal V=Z+U^{\rm f}$
solves \eqref{e:frc}, and
$\mathscr O_{r_0,\infty}(\mathcal V_t)\le2\mathscr O_{r_0,\infty}(Z_t)+2\mathscr O_{r_0,\infty}(U^{\rm f}_t)<\infty$.
Uniqueness follows from Lemma~\ref{l:solv}.

By \eqref{e:pceq} with $r=r_0$, $\mathcal K\mathcal H=0$ and
$\mathcal H(r_0)=E^{\rm base}(r_0)$. By \eqref{e:pB} with $\mathbf v=\psi_h$ this is the equation
for $\psi_h$ in \eqref{e:hqe}, and its $y$-derivative is the equation for $\mathfrak q_h$. Since
$\mathfrak b_-=-\mathfrak b_+$,
$S^{-1}=\begin{pmatrix}1&1\\ \frac1{2\mathfrak b_+}&-\frac1{2\mathfrak b_+}\end{pmatrix}$, so that
$\mathcal H_++\mathcal H_-=\psi_{h,1}$; with $\Xi=S^{-1}E^{\rm base}$, Lemma~\ref{l:all},
\eqref{e:pot}, $2f_0(V_c)=c\log\frac{V_c}{v_-}$, $K_++K_-=B_L$, $K_+-K_-=2\mathfrak b_+G_c$,
$\bar A_++\bar A_-=A_++A_--B_L$, $(\mathcal G-G_c)_y=\mathcal F$ and $Q_A$ of \eqref{e:pfr} give
\begin{equation}\label{e:hpri}
\begin{aligned}
 \Xi_1&=\mathcal B-c\log\frac{V_c}{v_-}-A_+-A_-,&\qquad (\Xi_1)_y&=b_m-H_c-a_+-a_-,\\
 \Xi_2&=\mathcal G-G_c-Q_A,&\qquad (\Xi_2)_y&=\mathcal F-(Q_A)_y .
\end{aligned}
\end{equation}
At $t=r_0$, $b_m,\mathcal F\in L^1\cap L^2$ by Lemma~\ref{l:cmp} and the boundedness of the
solution, $H_c$, $a_\pm$ and $(\mathfrak b_+)_y$ belong to $L^1\cap L^2$, and $\bar A_\pm$ are
bounded; with
$(Q_A)_y=\frac{a_+-a_-}{2\mathfrak b_+}-\frac{(\mathfrak b_+)_y}{2\mathfrak b_+^2}(\bar A_+-\bar A_-)$
this gives $\mathfrak q_h(r_0)=\pa_y\Xi(r_0)\in L^1\cap L^2$. On $[r_0,T]$,
$\mathfrak q_h=\pa_y\Xi-\pa_y(S^{-1}W)-\pa_y(S^{-1}\mathcal V)$; \eqref{e:hpri}, \eqref{e:munw},
$W=\mathcal P(V_c,z)+u$, Lemma~\ref{l:nf}(v), the regularity of $\mathcal V$ and the bounds for
$z$ in the proof of Lemma~\ref{l:bcor} give $\mathfrak q_h\in C([r_0,T];L^2)$ and
$(\mathfrak q_h)_y\in L^2((r_0,T)\times\R)$.

By \eqref{e:nco} with the matrix $B$, $Q(F_c(y,t),t)=\frac{\mathfrak q_h(y,t)}{f_c(y,t)}$ solves
$Q_t=\pa_\xi(D_0Q_\xi-bQ)$ with $b=f_cB+((F_c)_t-d(f_c)_y)I$ bounded on $\R\times[r_0,T]$, and
$\int_\R|Q(\xi,t)|\dd\xi=\Vert\mathfrak q_h(t)\Vert_{L^1}$; Lemma~\ref{l:kernv}(iv) gives
$\mathfrak q_h\in C([r_0,T];L^1)$. Uniqueness in this class is \eqref{e:uniq} with $B_1=B$ and
$B_0=B_2=0$, on $[r_0,T]$ and on $[r,T]$.
\end{proof}

\begin{lemma}\label{l:hq}
Let $a_s$, $b_0$ be the exponents of \eqref{e:pwt}, $b_a$ the function of
\eqref{e:cmp}, $R_M=S^{-1}\widetilde{\mathcal M}S$,
and put
\begin{equation}\label{e:pdecq}
\begin{aligned}
 g_{\rm ref}&=\frac{\nu(b_a)_y}{V_b}-c\,\pa_t\log\frac{V_c(y+h_*,t)}{v_-}-\pa_t(A_-+A_+),\\
 R_p&=\frac{\nu(b_m-b_a)_y}{V_b}-j_{r_0}-\pa_t(\mathcal V_++\mathcal V_-)
   -d(V_c)(\mathfrak q_{h,1})_y-(R_M\mathfrak q_h)_1 .
\end{aligned}
\end{equation}
Then, on $\R\times[r_0,\infty)$,
\begin{equation}\label{e:pq}
 P_b-p_*=g_{\rm ref}+R_p+\frac{\mathfrak q_{h,2}}{V_c},\qquad
 \limsup_{t\to\infty}\tau A_{g_{\rm ref}}(t)\le C_P\delta^2,\qquad
 \mathscr O_{r_0,\infty}(R_p)<\infty,
\end{equation}
where $C_P$ depends only on the fixed data, and not on $M$, $a_s$, $b_0$, $r_0$ or the initial datum.
\end{lemma}
\begin{proof}
Write $d=d(V_c)$ and $\Phi=S^{-1}U^{\rm f}$. Since $B=A_{\rm ph}-R_M$ and
$(A_{\rm ph}\mathfrak q_h)_1=\frac{\mathfrak q_{h,2}}{V_c}$, the first row of
\eqref{e:hqe} gives
\[
 -\pa_t\psi_{h,1}=\frac{\mathfrak q_{h,2}}{V_c}-d(\mathfrak q_{h,1})_y-(R_M\mathfrak q_h)_1 .
\]
Substituting this, $\mathcal E=\mathcal H+\mathcal V$ and
$\mathcal H_++\mathcal H_-=\psi_{h,1}$ into \eqref{e:bpot} with $r=r_0$ gives the
first identity in \eqref{e:pq}.

With $V_h$, $H_h$ of Lemma~\ref{l:mcoef} and $b_a=H_h+a_-+a_+$ (by \eqref{e:cmp},
$\widehat{\mathcal U}=(V_c,H_c,c_v\Theta_c)$, $r_\sigma=(-1,s_\sigma,p_*)$ and
$a_\sigma=s_\sigma\vartheta_\sigma$), $g_{\rm ref}=g^{(1)}+g^{(2)}$,
\[
 g^{(1)}=\frac{\nu(H_h)_y}{V_b}-c\,\pa_t\log\frac{V_h}{v_-},\qquad
 g^{(2)}=\frac{\nu(a_-+a_+)_y}{V_b}-\pa_t(A_-+A_+).
\]
By the bounds for the derivatives of $V_c$ stated before \eqref{e:pcnx},
$H_c=\frac{c(V_c)_y}{V_c}$ and $V_b\ge\frac1{2\rho_-}$ (Lemma~\ref{l:mcoef}),
$|g^{(1)}|\le C\delta\tau^{-1}e^{-\frac{c_G(y+h_*)^2}{2\tau}}$ with $C$ depending only on the fixed
data, so $G\le1$ gives
\[
 \tau A_{g^{(1)}}(t)\le C^2\delta^2\tau^{-\frac12}\int_\R e^{-\frac{c_G(y+h_*)^2}\tau}\dd y
 =C_P\delta^2,\qquad C_P=C^2\sqrt{\frac\pi{c_G}} .
\]
By \eqref{e:CH}, $A_j$ is a function of $z_j$, so $\pa_tA_j=-(s_j+\frac{y-s_jt}{2\tau})a_j$;
with Lemma~\ref{l:pr} and $(y+h_*)^2+(y-s_jt)^2\ge\frac{(s_jt+h_*)^2}2$ this gives
$A_{g^{(2)}}(t)\le C_be^{-c_bt}$, and
$\tau A_{g_{\rm ref}}\le\tau A_{g^{(1)}}+2(\tau A_{g^{(1)}}\,\tau A_{g^{(2)}})^{\frac12}+\tau A_{g^{(2)}}$
gives the bound for $\limsup\tau A_{g_{\rm ref}}$.

By Lemma~\ref{l:bulk}, Lemma~\ref{l:bcor} and Lemma~\ref{l:for}, the first three terms of $R_p$
have finite $\mathscr O_{r_0,\infty}$. By \eqref{e:pmat}, $|R_M|\le C_b\tau^{-\frac12}$, and
$G\le G_0$ gives $A_{R_M\mathfrak q_h}\le C_b\int_\R W_c|\mathfrak q_h|^2\dd y$; it remains to
prove
\begin{equation}\label{e:hqfin}
 \mathscr Q_{r_0}(\mathfrak q_h)+\mathscr O_{r_0,\infty}\bigl((\mathfrak q_h)_y\bigr)<\infty .
\end{equation}
By $\mathcal V=Z^{\rm ref}+U^{\rm f}$ and $W=\mathcal P(V_c,z)+u$,
\[
 \psi_h=\Xi-S^{-1}W-S^{-1}Z^{\rm ref}-\Phi,
\]
and \eqref{e:hqfin} follows from the bounds below for
$\mathscr Q_{r_0}(\pa_y\Upsilon)+\mathscr O_{r_0,\infty}(\pa_y^2\Upsilon)$,
$\Upsilon\in\{\Xi,S^{-1}W,S^{-1}Z^{\rm ref},\Phi\}$. Here
$|\pa_y^jS^{-1}|\le C\delta\tau^{-\frac j2}e^{-\frac{c_Gy^2}{2\tau}}$ $(j=1,2)$ and
$|\pa_tS^{-1}|\le C\delta\tau^{-1}$; a function $f$ bounded by $C_b\tau^{-\frac32}$ has
$A_f\le C_b\tau^{-2}$, and one bounded by $C_be^{-c_bt}$ or by
$C_M\tau^{-\frac12}e^{-\frac{(y-s_jt)^2}{C\tau}}$ has $A_f+\int_\R W_c|f|^2\dd y\le C_be^{-c_bt}$.

For $\Upsilon=\Xi$: by \eqref{e:hpri}, $b_a=H_h+a_++a_-$, $\mathcal F=V_b(P_b-p_*)+\frac{\gamma-1}2b_m^2$
and $(Q_A)_y$ above,
\[
\begin{gathered}
 (\Xi_1)_y=(b_m-H_c)-a_+-a_-,\qquad (\Xi_1)_{yy}=(b_m-b_a)_y+(H_h-H_c)_y,\\
 (\Xi_2)_y=V_b(P_b-p_*)+\frac{\gamma-1}2b_m^2-(Q_A)_y,\\
 (\Xi_2)_{yy}=V_b(P_b-p_a)_y+V_b(p_a)_y+(V_b)_y(P_b-p_*)+(\gamma-1)b_m(b_m)_y-(Q_A)_{yy}.
\end{gathered}
\]
For $\tau\ge1$, $(x+a)^2\ge\frac{x^2}2-a^2$, $\frac{b_0}2x^2-2|x|\ge-\frac2{b_0}$ and
$\operatorname{sech}^2x\ge e^{-2|x|}$ give $G_0\le e^{b_0h_*^2+\frac2{b_0}}\operatorname{sech}^2\frac y{\sqrt\tau}$;
with $\upsilon_m=b_m-H_c$ of \eqref{e:Um}, \eqref{e:hcur} and the last quantity in
Lemma~\ref{l:dep},
\[
 \mathscr Q_{r_0}(P_b-p_*)+\mathscr Q_{r_0}(b_m-H_c)
 \le C_*e^{b_0h_*^2+\frac2{b_0}}\int_{r_0}^\infty\!\!\int_\R(q_+^2+q_-^2)\tau^{-\frac12}
 \operatorname{sech}^2\frac y{\sqrt\tau}\dd y\dd t<\infty .
\]
By \eqref{e:munw}, $b_m^2+(V_b)_y^2\le C_b\tau^{-1}$ and $\int_{r_0}^\infty\Vert(b_m)_y\Vert_{L^2}^2<\infty$;
so $\int W_cb_m^4\le C_b\tau^{-1}\int W_c(|b_m-H_c|^2+H_c^2)$,
$A_{(V_b)_y(P_b-p_*)}\le C_b\int W_c|P_b-p_*|^2$ and
$A_{b_m(b_m)_y}\le C_b\tau^{-\frac12}\Vert(b_m)_y\Vert_{L^2}^2$, with
$\int W_cH_c^2\le C\delta^2\tau^{-1}$. By Lemma~\ref{l:bulk},
$\mathscr O((b_m-b_a)_y)+\mathscr O((P_b-p_a)_y)<\infty$. By
$p_a-p_*=\frac{\gamma p_*\sum_j\vartheta_j-\frac{\gamma-1}2b_a^2}{v_a}$ (proof of
Lemma~\ref{l:mcoef}) and \eqref{e:cmpb}, $(p_a)_y$ is the sum of a function bounded by
$C_M\tau^{-\frac32}$ and of functions bounded by $C_M\tau^{-\frac12}e^{-\frac{(y-s_jt)^2}{C\tau}}$;
$|(H_h-H_c)_y|\le|h_*|\sup|(H_c)_{yy}|\le C_b\tau^{-\frac32}$. The functions $a_\pm$, $(a_\pm)_y$ are
bounded by $C_M\tau^{-\frac12}e^{-\frac{(y-s_\pm t)^2}{8d_\pm\tau}}$, and so are the terms of
$(Q_A)_y$, $(Q_A)_{yy}$ that contain $a_\pm$ or $(a_\pm)_y$; each other term is a derivative of
$V_c$ times $\bar A_\pm$, hence bounded by $C_be^{-c_bt}$, since
$|\bar A_+|\le C_Me^{-\frac{z_+^2}{4d_+}}$ for $y\le s_+t$ and
$|\bar A_-|\le C_Me^{-\frac{z_-^2}{4d_-}}$ for $y\ge s_-t$. This gives the bound for $\Xi$.

For $\Upsilon=S^{-1}W$: \eqref{e:ubC}, \eqref{e:nfe}, Lemma~\ref{l:dep} and
\eqref{eq:pr-gau} give $\sup_{t\ge r_0}\Vert W(t)\Vert_{L^2}<\infty$, and
Lemma~\ref{l:peq} gives $\int_{r_0}^\infty\Vert W_y\Vert_{L^2}^2\dd t<\infty$.
The bounds for $z$ at the start of the proof of Lemma~\ref{l:bcor} and its
weighted $u_{yy}$ estimate give, by the product rule,
\[
 \int_{r_0}^\infty\tau\Vert W_{yy}\Vert_{L^2}^2\dd t
 \le C_b\int_{r_0}^\infty\Bigl[\tau\Vert(z_{yy},u_{yy})\Vert_{L^2}^2
   +\Vert z_y\Vert_{L^2}^2+\tau^{-3/2}\Bigr]\dd t<\infty .
\]
The bounds for the derivatives of $S^{-1}$ recorded above now give
\[
 \int_\R W_c|\pa_y(S^{-1}W)|^2\dd y+A_{\pa_y^2(S^{-1}W)}
 \le C_b\Bigl[\sqrt\tau\Vert W_{yy}\Vert_{L^2}^2
   +\tau^{-\frac12}\Vert W_y\Vert_{L^2}^2+\tau^{-\frac32}\Vert W\Vert_{L^2}^2\Bigr],
\]
whose right side is integrable on $(r_0,\infty)$.

For $\Upsilon=S^{-1}Z^{\rm ref}$: by Lemma~\ref{l:ref}, $|\pa_y(S^{-1}Z^{\rm ref})|\le C\tau^{-1}$
and $|\pa_y^2(S^{-1}Z^{\rm ref})|\le C\tau^{-\frac32}$ on $|y|\le a_0t$ for $t\ge t_0$, and both are
bounded; with $W_c+\sqrt\tau G\le C\sqrt\tau e^{-\frac{b_0y^2}{4\tau}}$ on $|y|\ge2|h_*|$, the two
quantities are at most $C\tau^{-2}+C_b\tau e^{-c_bt}$.

For $\Upsilon=\Phi$: by \eqref{e:gop} and $B=A_{\rm ph}-R_M$,
$\Phi_t-d\Phi_{yy}+(A_{\rm ph}-R_M)\Phi_y=S^{-1}f^\star$. Let
$H_0=\operatorname{diag}(V_c,\frac1{\gamma p_*})$, so that
$H_0A_{\rm ph}=\begin{pmatrix}0&1\\1&0\end{pmatrix}$, and $w=\sqrt\tau G$. Multiplication of the
equation by $w\Phi_{yy}^\top H_0$ and
$\int w\Phi_{yy}^\top H_0A_{\rm ph}\Phi_y=-\frac12\int w_y\Phi_y^\top H_0A_{\rm ph}\Phi_y$ give
\[
 \int_\R dw\,\Phi_{yy}^\top H_0\Phi_{yy}\dd y
 =\int_\R w\Phi_{yy}^\top H_0\bigl(\Phi_t-R_M\Phi_y-S^{-1}f^\star\bigr)\dd y
 -\frac12\int_\R w_y\Phi_y^\top H_0A_{\rm ph}\Phi_y\dd y .
\]
Since $dH_0\ge c_*$, $|w_y|\le CG_0$ (as $b_0<\frac{a_s}4$) and $w|R_M|^2\le C_bG_0$, Cauchy's
inequality gives
\[
 A_{\Phi_{yy}}\le C_b\bigl(A_{\Phi_t}+A_{f^\star}\bigr)+C_b\int_\R G_0|\Phi_y|^2\dd y .
\]
Here $\mathscr O(\Phi_t)\le C\mathscr O(U^{\rm f}_t)+C\mathscr Q(U^{\rm f})$,
$\int\!\!\int G_0|\Phi_y|^2+\mathscr Q(\Phi_y)\le CJ_{r_0}(U^{\rm f})$, and
$\mathscr O_{r_0,\infty}(f^\star)<\infty$, because every source of Lemma~\ref{l:solv} has finite
$\mathscr O_{r_0,\infty}$: in (a) by $G\le1$, in (b) by $\sqrt\tau G\le\tau$ and
$\tau|K_i|^2\le C_{K,b}$ (Step 0 of the proof of Lemma~\ref{l:solv}), and in (c) by
$\sqrt\tau G|A|^2\le C\delta^2W_c$. The test is first made
on finite strips, where $\Phi\in L^2(r_0,T;H^2)$; hence
$\mathscr Q_{r_0}(\Phi_y)+\mathscr O_{r_0,\infty}(\Phi_{yy})<\infty$, which completes the proof of
\eqref{e:hqfin}.
\end{proof}

\subsection{Estimates for propagated data}

Put $a_\xi=\frac{a_s}{16}$ and, for $s<T\le\infty$ and a function $v$ on
$\R\times(s,T)$,
\begin{equation}\label{e:xnorm}
 \mathscr N^\xi_{s,T}(v)=\int_s^T\Bigl[\tau^{\frac14}
   \bigl\Vert e^{-\frac{a_\xi\xi^2}{4\tau}}v_\xi\bigr\Vert_{L^2}
   +\tau^{-\frac14}\bigl\Vert e^{-\frac{a_\xi\xi^2}{4\tau}}v\bigr\Vert_{L^2}
   \Bigr]\dd t .
\end{equation}
Since $|V_c-v_-|\le C\delta$, $\frac12\le\rho_-V_c\le2$; hence
$2^{-\frac12}\le f_c\le2^{\frac12}$ and $|F_c(y,t)|\le\sqrt2|y|$. If
$\tau\ge h_*^2$, then $a_\xi F_c^2\le\frac{a_s}8y^2\le a_s(y+h_*)^2+a_sh_*^2$
and
$e^{-\frac{a_s(y+h_*)^2}{4\tau}}\le e^{\frac{a_s}4}e^{-\frac{a_\xi F_c(y,t)^2}{4\tau}}$.
With \eqref{e:xder}, $|(f_c)_y|\le C\delta\tau^{-\frac12}$ and
$\dd y=\frac{\dd\xi}{f_c}$, this gives, for $w(y,t)=f_c(y,t)v(F_c(y,t),t)$,
\begin{equation}\label{e:xmass}
 \mathscr N_{s,T}(w)\le C\mathscr N^\xi_{s,T}(v)\qquad(\tau_s\ge h_*^2).
\end{equation}

\begin{lemma}\label{l:dr}
Let $\bar b\ne0$, $0<v_0\le\frac{|\bar b|}4$, $A_0,A_1,A_2\ge0$ and
$I=[r,\infty)$, let $b$ satisfy the hypotheses of \cref{l:kernv}\textup{(ii)},
and let $H$ be the kernel of $v_t=D_0v_{\xi\xi}-\pa_\xi(bv)$. For $s\ge r$ put
$(H(s)z)(\xi,t)=\int_\R H(t,s;\xi,\eta)z(\eta)\dd\eta$. Then, with $\mathscr N^\xi$ as in
\eqref{e:xnorm},
\begin{equation}\label{e:birth}
\begin{aligned}
 \mathscr N^\xi_{s,\infty}(H(s)z)&\le C\tau_s^{\frac38}\Vert z\Vert_{L^1},\\
 \mathscr N^\xi_{s,\infty}(H(s)z)&\le C\Vert z\Vert_{L^1}
   &&(z=0\text{ on }|\eta|<v_0\tau_s),\\
 \mathscr N^\xi_{s,\infty}(H(s)z)&\le C\tau_s^{\frac58}\Vert z\Vert_{L^2}
   &&(z=0\text{ on }|\eta|>2v_0\tau_s),
\end{aligned}
\end{equation}
\begin{equation}\label{e:bcol}
 \mathscr N^\xi_{s,\infty}\bigl(H(\cdot,s;\cdot,\eta)\bigr)
 \le C\tau_s^{\frac38}\Bigl(1+\frac{|\eta|}{\sqrt{\tau_s}}\Bigr)^{-\frac34}
 \qquad(|\eta|\le2v_0\tau_s),
\end{equation}
and, for $w>0$,
\begin{equation}\label{e:dsmb}
 \mathscr N^\xi_{s,\infty}(H(s)z)\le C_0\Bigl[\tau_s^{\frac14}
   \bigl\Vert e^{\frac{w\eta^2}{\tau_s}}z\bigr\Vert_{L^2}
   +\tau_s^{\frac34}\bigl\Vert e^{\frac{w\eta^2}{\tau_s}}z_\eta\bigr\Vert_{L^2}\Bigr]
   +C\bigl\Vert e^{\frac{w\eta^2}{\tau_s}}z\bigr\Vert_{L^2},
\end{equation}
where $C_0$ depends only on $D_0$, $\bar b$, $a_s$, $w$ and an upper bound for
$A_0$, and $C$ depends in addition on $v_0$ and upper bounds for $A_1$ and
$A_2$.
\end{lemma}

\begin{proof}
Write $S=\tau_s$, $h=t-s$ and $X=\xi-\eta-\bar bh$. Let $C_1\ge1$, depending
only on $D_0$ and an upper bound for $A_0$, be such that
$H\le C_1G(C_1h,X)$, $|H_\eta|\le C_1(h^{-\frac12}+S^{-\frac12})G(C_1h,X)$ and
$|H_\xi|\le C_1(h^{-\frac12}+A_2^{\frac12})G(C_1h,X)$ (\cref{l:kern}(b),
\cref{l:kernv}(ii)), and let $\varpi'$ and $C'\ge C_1$ be the constants of
\cref{l:kernv}(ii); $G(C_1h,X)\le(\frac{C'}{C_1})^{\frac12}G(C'h,X)$. For
$\lambda\ge1$ put $c_1=\bar b^2(4a_\xi^{-1}+8D_0\lambda+2w^{-1})^{-1}$ and
$c_2=(4a_\xi^{-1}+8D_0\lambda)^{-1}$. Since $\xi-X-\eta=\bar bh$, Cauchy's
inequality gives
\[
 \frac{a_\xi\xi^2}{4\tau_t}+\frac{X^2}{8D_0\lambda h}+\frac{w\eta^2}{2S}
 \ge\frac{\bar b^2h^2}{4a_\xi^{-1}\tau_t+8D_0\lambda h+2w^{-1}S}\ge\frac{c_1h^2}{\tau_t},
 \qquad
 \frac{a_\xi\xi^2}{2\tau_t}+\frac{X^2}{4D_0\lambda h}\ge\frac{2c_2(\eta+\bar bh)^2}{\tau_t}.
\]
With $G(\lambda h,X)=\sqrt2e^{-\frac{X^2}{8D_0\lambda h}}G(2\lambda h,X)$,
$G(\lambda h,X)^2=(4\pi D_0\lambda h)^{-\frac12}e^{-\frac{X^2}{4D_0\lambda h}}G(\lambda h,X)$
and Schur's test this gives
\begin{equation}\label{e:gsch}
 \Bigl\Vert e^{-\frac{a_\xi\xi^2}{4\tau_t}}\int_\R G(\lambda h,X)|g(\eta)|\dd\eta
 \Bigr\Vert_{L^2_\xi}
 \le\sqrt2e^{-\frac{c_1h^2}{\tau_t}}\bigl\Vert e^{\frac{w\eta^2}{2S}}g\bigr\Vert_{L^2},
 \qquad
 \bigl\Vert e^{-\frac{a_\xi\xi^2}{4\tau_t}}G(\lambda h,X)\bigr\Vert_{L^2_\xi}
 \le\frac{e^{-\frac{c_2(\eta+\bar bh)^2}{\tau_t}}}{(4\pi D_0h)^{\frac14}} .
\end{equation}

On $|\xi|\le\frac{v_0}2\tau_t$, $|H_\xi|\le C'h^{-\frac12}G(C'h,X)$; on
$|\xi|>\frac{v_0}2\tau_t$, $e^{-\frac{a_\xi\xi^2}{4\tau_t}}\le
e^{-\frac{a_\xi v_0^2\tau_t}{16}}$ and
$\Vert G(C_1h,X)\Vert_{L^2_\xi}\le(4\pi D_0h)^{-\frac14}$. By \eqref{e:gsch} with
$\lambda=C'$, the integrand of $\mathscr N^\xi_{s,\infty}(H(\cdot,s;\cdot,\eta))$
is at most
\[
 C\bigl(\tau_t^{\frac14}h^{-\frac34}+\tau_t^{-\frac14}h^{-\frac14}\bigr)
   e^{-\frac{c_2(\eta+\bar bh)^2}{\tau_t}}
 +C\tau_t^{\frac14}\bigl(h^{-\frac34}+h^{-\frac14}\bigr)
   e^{-\frac{a_\xi v_0^2\tau_t}{16}} .
\]
The time integral of the second term is at most
$Ce^{-\frac{a_\xi v_0^2S}{32}}$. In the first,
$\tau_t^{-\frac14}h^{-\frac14}\le\tau_t^{\frac14}h^{-\frac34}$, and for $h\ge S$,
$\tau_t\le2h$, $\tau_t^{\frac14}h^{-\frac34}\le2h^{-\frac12}$ and
\[
 \int_0^\infty h^{-\frac12}e^{-\frac{c_*(\eta+\bar bh)^2}{h}}\dd h
 =\frac1{|\bar b|}\sqrt{\frac\pi{c_*}}\,e^{-2c_*(\eta\bar b+|\eta\bar b|)}
 \le\frac1{|\bar b|}\sqrt{\frac\pi{c_*}}\qquad(c_*>0).
\]
For $h\le S$, $\tau_t\le2S$; $\int_0^{\sqrt S}h^{-\frac34}\dd h=4S^{\frac18}$,
and for $\sqrt S<h\le S$, $h^{-\frac34}\le S^{-\frac38}$ and
$\int_\R e^{-\frac{c_2(\eta+\bar bh)^2}{2S}}\dd h=\frac{\sqrt{2\pi S}}{|\bar b|\sqrt{c_2}}$.
This proves the first line of \eqref{e:birth} for a unit mass at $\eta$, and
hence, by Minkowski's inequality, for $z\in L^1$. If $|\eta|\ge v_0S$, then
$|\eta+\bar bh|\ge\frac{v_0S}2$ for $h\le\frac{|\eta|}{2|\bar b|}$, and
$h^{-\frac34}\le(\frac{2|\bar b|}{v_0S})^{\frac34}$ for
$h>\frac{|\eta|}{2|\bar b|}$; the range $h\le S$ then contributes at most $C$,
and Minkowski's inequality gives the second line. If
$\sqrt S\le|\eta|\le2v_0S$, then $\frac{|\eta|}{2|\bar b|}\le\frac S4$, and the
same splitting at $h=\frac{|\eta|}{2|\bar b|}$ bounds the range $h\le S$ by
\[
 CS^{\frac14}|\eta|^{\frac14}e^{-\frac{c_2\eta^2}{8S}}
   +CS^{\frac14}|\eta|^{-\frac34}\sqrt S
 \le CS^{\frac38}\Bigl(1+\frac{|\eta|}{\sqrt S}\Bigr)^{-\frac34};
\]
for $h\ge S$, $|\eta+\bar bh|\ge|\bar b|h-2v_0S\ge\frac{|\bar b|h}2$ gives
$Ce^{-\frac{c_2\bar b^2S}{16}}$, and $e^{-c_*S}\le
C(1+\frac{|\eta|}{\sqrt S})^{-\frac34}$ for $|\eta|\le2v_0S$. With the first
line for $|\eta|<\sqrt S$, this is \eqref{e:bcol}. By Minkowski's inequality,
\eqref{e:bcol}, the Cauchy--Schwarz inequality and
$\int_\R(1+\frac{|\eta|}{\sqrt S})^{-\frac32}\dd\eta=4\sqrt S$, the third line
of \eqref{e:birth} follows.

For \eqref{e:dsmb} let $v=H(s)z$. Integration by parts in $\eta$ gives
\[
 v_\xi=\int_\R Hz_\eta\dd\eta+\int_{|\eta|\le\frac{v_0}2S}(H_\xi+H_\eta)z\dd\eta
   +\int_{|\eta|>\frac{v_0}2S}(H_\xi+H_\eta)z\dd\eta
 \qquad(h\le\varpi'S),
\]
and $v_\xi=\int_\R H_\xi z\dd\eta$ for $h>\varpi'S$. By \eqref{e:gsch} with
$\lambda=C_1$ and
\[
 \int_0^\infty\tau_t^{\frac14}e^{-\frac{c_1h^2}{\tau_t}}\dd h\le CS^{\frac34},
 \qquad
 \int_0^\infty\tau_t^{-\frac14}e^{-\frac{c_1h^2}{\tau_t}}\dd h\le CS^{\frac14},
\]
the first integral and $v$ contribute the first term of \eqref{e:dsmb}, with
$C_0$ depending only on $C_1$ and $c_1$. By \cref{l:kernv}(ii),
$|H_\xi+H_\eta|\le C'\sqrt h(A_1S^{-1}+A_2e^{-\frac S{C'}})G(C'h,X)$ in the
second integral, and
$|H_\xi|+|H_\eta|\le C(h^{-\frac12}+A_2^{\frac12}+1)G(C'h,X)$ in the third and
for $h>\varpi'S$. Apply \eqref{e:gsch} with $\lambda=C'$, using
$\int_0^S\sqrt h\,e^{-\frac{c_1h^2}{2S}}\dd h\le CS^{\frac34}$,
$\Vert e^{\frac{w\eta^2}{2S}}g\Vert_{L^2}\le e^{-\frac{wv_0^2S}8}
\Vert e^{\frac{w\eta^2}S}g\Vert_{L^2}$ for $g=0$ on $|\eta|\le\frac{v_0}2S$,
and $\frac{h^2}{\tau_t}\ge\frac{\varpi'h}{1+\varpi'}$ for $h\ge\varpi'S$. These
contributions are at most
\[
 C\bigl(A_1+A_2Se^{-\frac S{C'}}+S^{\frac54}e^{-\frac{wv_0^2S}8}+1\bigr)
 \bigl\Vert e^{\frac{w\eta^2}S}z\bigr\Vert_{L^2},
\]
which proves \eqref{e:dsmb}.
\end{proof}

\subsection{Central derivatives}

\begin{lemma}\label{l:ga}
There is $C_b$ such that $\tau|(R_U)_y(y,t)|\le C_b$ for $t\ge S_0$ and
$|y|\le\frac{c_-}{1024}t$.
\end{lemma}

\begin{proof}
By Lemma~\ref{l:bg}, \eqref{e:cmp} and \eqref{e:wval}, $(R_U)_y$ is continuous
and bounded on $\R\times[0,T_e]$ for every $T_e<\infty$, and
$\Vert(R_U)_y(t)\Vert_{L^\infty}\to0$. Suppose that the assertion fails. Then
there are $(z_n^0,t_n^0)$ with $t_n^0\to\infty$,
$|z_n^0|\le\frac{c_-}{1024}t_n^0$ and $t_n^0|(R_U)_y(z_n^0,t_n^0)|\to\infty$.
Lemma~\ref{l:dbl} with $f=(R_U)_y$, $p=1$, $c_1=\bar c$ and
$c_2=\frac{c_-}{1024}$ gives $(z_n,t_n)$ with
\begin{equation}\label{e:gbox}
 \beta_n=|(R_U)_y(z_n,t_n)|\to0,\qquad t_n\beta_n\to\infty,\qquad
 \sup_{S_n}|(R_U)_y|\le2\beta_n,\qquad |y|\le\frac{c_-}{512}t\ \text{ on }S_n,
\end{equation}
where $S_n,J_n,I_n$ and $r_n=t_n-H_n$ are given by \eqref{e:win} with
$H_n=\beta_n^{-1}$. In this proof $n$ is large and $o(\cdot)$ refers to
$n\to\infty$. On $S_n$, $\tau^{-1}\le\frac2{t_n}=o(\beta_n)$, hence
$4H_n\le\varpi t_n$. Let $Q,K,h,A,F$ be as in \eqref{e:bmat}, put
$z=\mathcal U-U_*(V_c)$, and recall
$(\varphi_m)_y=V_bQ+\frac{H_c}c\varphi_m$,
$(\vartheta_m)_y=\Theta_bK+\frac{H_c}c\vartheta_m$, $(\upsilon_m)_y=V_bh$.
By \eqref{e:cmp}, Lemma~\ref{l:pr}, \eqref{eq:pr-gau} and
$|h_*|\le C_M\delta^{-1}$,
$\Vert\partial_y^i[(v_a,b_a,T_a)-(V_c,H_c,\Theta_c)]\Vert_{L^2}\le
C_b\tau^{-\frac i2-\frac14}$ and
$|\partial_y^i[(v_a,b_a,T_a)-(V_c,H_c,\Theta_c)]|\le C_b\tau^{-\frac{i+1}2}$
for $1\le i\le3$. With $W=Z_m-[(v_a,b_a,T_a)-(V_c,H_c,\Theta_c)]$,
\eqref{e:wval}, Lemma~\ref{l:asy}, \eqref{e:ZR}, \eqref{e:wtC} and
Lemma~\ref{l:dep}, this gives
\begin{equation}\label{e:gseed}
\begin{gathered}
 |R_U|+|Z_m|+|z|+|b_m|\le C_b\tau^{-\frac12},\qquad
 |\partial_yZ_m|+|z_y|+|(Q,K,h)|\le C\beta_n\qquad\text{on }S_n,\\
 \sup_{t\ge t_n-2H_n}\tau\bigl(\Vert\partial_yZ_m\Vert_{L^2}^2
   +\Vert(Q,K,h)\Vert_{L^2}^2\bigr)\to0,\qquad
 \int_{t_n-2H_n}^{t_n}\Vert(q,k,\ell)\Vert_{L^2(J_n)}^2\dd t=o(t_n^{-1}),
\end{gathered}
\end{equation}
where $(q,k,\ell)=\partial_y(Q,K,h)$; for the last bound,
$\int_{t_n-2H_n}^{t_n}\Vert\partial_y^2W\Vert_{L^2}^2\dd t\le\frac2{t_n}
\int_{t_n-2H_n}^\infty\tau\Vert W_{yy}\Vert_{L^2}^2\dd t$ and
\begin{equation}\label{e:gconv}
\begin{gathered}
 \partial_y^2\varphi_m=V_b\Bigl(q+Q^2+\frac{2H_cQ}c\Bigr)
   +\Bigl(\partial_y\frac{H_c}c+\frac{H_c^2}{c^2}\Bigr)\varphi_m,\qquad
 \partial_y^2\upsilon_m=V_b\ell+V_b\Bigl(Q+\frac{H_c}c\Bigr)h,\\
 \partial_y^2\vartheta_m=\Theta_b\Bigl(k+K^2+\frac{2H_cK}c\Bigr)
   +\Bigl(\partial_y\frac{H_c}c+\frac{H_c^2}{c^2}\Bigr)\vartheta_m .
\end{gathered}
\end{equation}

By \eqref{e:bmat} and $A+\frac{(H_c)_y}{V_c}=\frac{(H_c)_y}{V_b}$,
$h_t=\frac1{V_b}\partial_y(\upsilon_m)_t-h(h+\frac{(H_c)_y}{V_b})$
and $\partial_y[P_b(K-Q)]=P_b(k-q)+P_b(K-Q)^2$; hence
\begin{equation}\label{e:ghsys}
\begin{gathered}
 q_t=\ell_y+A_{yy},\qquad k_t=-(\gamma-1)\ell_y+F_{yy},\\
 \ell_t=\partial_y\Bigl\{\frac1{V_b}\bigl[\nu\ell_y-P_b(k-q)-P_b(K-Q)^2+S_y\bigr]\Bigr\}
   -\Bigl(2h+\frac{(H_c)_y}{V_b}\Bigr)\ell-\Bigl(\frac{(H_c)_y}{V_b}\Bigr)_yh,
\end{gathered}
\end{equation}
with $S=\nu A_y+(\nu-c)\partial_y\frac{(H_c)_y}{V_c}$. Let $\psi_n$ and
$\lambda_n=\frac{(\psi_n)_y}{\psi_n}$ be as in the proof of Lemma~\ref{l:win},
and put $w=\frac{\psi_nV_c}{V_b^2}$, $m=\frac{V_b^2}{\Theta_b}$,
\[
 E_2=\frac12\int_\R w\bigl(Rq^2+c_vk^2+m\ell^2\bigr)\dd y,\qquad
 D_2=\int_\R w\Bigl(\frac{\nu V_b}{\Theta_b}\ell_y^2+\frac\kappa{V_b}k_y^2\Bigr)\dd y .
\]
By \eqref{e:bmat}, $\frac{w_y}w=\lambda_n-2Q-\frac{H_c}c$,
$\frac{(wm)_y}{wm}=\lambda_n-K$,
$\frac{w_t}w=\frac{(H_c)_y}{V_c}-2h-\frac{2(H_c)_y}{V_b}$ and
$\frac{(wm)_t}{wm}=(\gamma-1)h-F$; moreover $\frac{mP_b}{V_b}=R$,
$|w\lambda_n|\le CH_n^{-1}$ and $w\lambda_n^2\le CH_n^{-2}$. Multiply the
equations \eqref{e:ghsys} by $Rwq$, $c_vwk$ and $wm\ell$, integrate over $\R$
and integrate by parts. The terms with $q\ell_y$, $k\ell_y$ and
$(k-q)\ell_y$ cancel, and, with
$\mathcal T=F_y-\frac{\kappa k_y}{c_vV_b}$,
\begin{equation}\label{e:ghe}
\begin{aligned}
 E_2'+D_2={}&R\int_\R w(\lambda_n-K)(k-q)\ell\dd y
   +R\int_\R w\bigl[\ell_y+(\lambda_n-K)\ell\bigr](K-Q)^2\dd y\\
 &-\nu\int_\R\frac{wV_b}{\Theta_b}(\lambda_n-K)\ell\ell_y\dd y
   -\int_\R\frac{wV_b}{\Theta_b}\bigl[\ell_y+(\lambda_n-K)\ell\bigr]S_y\dd y\\
 &-\kappa\int_\R\frac w{V_b}\Bigl(\lambda_n-2Q-\frac{H_c}c\Bigr)kk_y\dd y
   -c_v\int_\R w\Bigl[k_y+\Bigl(\lambda_n-2Q-\frac{H_c}c\Bigr)k\Bigr]\mathcal T\dd y\\
 &+R\int_\R wqA_{yy}\dd y-\int_\R wm\Bigl(\frac{(H_c)_y}{V_b}\Bigr)_yh\ell\dd y\\
 &+\frac12\int_\R w\Bigl(\frac{(H_c)_y}{V_c}-2h-\frac{2(H_c)_y}{V_b}\Bigr)
   (Rq^2+c_vk^2)\dd y\\
 &+\int_\R wm\Bigl(\frac{\gamma-5}2h-\frac{(H_c)_y}{V_b}-\frac F2\Bigr)\ell^2\dd y .
\end{aligned}
\end{equation}
In the last integral, $wm\frac{\kappa k}{c_vV_b}=\frac\kappa{c_v}\frac{wV_b}{\Theta_b}k$
and $\int_\R\frac{wV_b}{\Theta_b}k\ell^2\dd y=-\int_\R\frac{wV_b}{\Theta_b}
K\bigl[2\ell\ell_y+(\lambda_n-Q-K-\frac{H_c}c)\ell^2\bigr]\dd y$. By the formula for
$A_y$ in the proof of Lemma~\ref{l:win},
\[
\begin{gathered}
 A_{yy}=-\frac{(H_c)_y}{V_b}q+\mathcal A,\qquad
 S_y=-\frac{\nu(H_c)_y}{V_b}q+\nu\mathcal A+(\nu-c)\partial_y^2\frac{(H_c)_y}{V_c},\\
 \mathcal A=-\Bigl(\frac{(H_c)_y}{V_b}\Bigr)_yQ
   -\Bigl((H_c)_{yy}-\frac{H_c(H_c)_y}c\Bigr)\frac Q{V_b}\\
   +\Bigl[(H_c)_{yyy}-\frac{2H_c(H_c)_{yy}}c
     +\Bigl(\frac{H_c^2}{c^2}-\partial_y\frac{H_c}c\Bigr)(H_c)_y\Bigr]
     \Bigl(\frac1{V_b}-\frac1{V_c}\Bigr),
\end{gathered}
\]
and differentiation of $F$ in \eqref{e:bmat} gives $\mathcal T=\mathcal T_1+\mathcal T_0$,
\[
 \mathcal T_1=\frac\kappa{c_vV_b}\Bigl[(k-q)K+(K-Q)k-\Bigl(Q+\frac{H_c}c\Bigr)k\Bigr]
   +\frac{2\nu V_bh\ell}{c_v\Theta_b}+\frac{\gamma H_c(k-q)}{V_b}
   +\frac{2\nu(H_c)_y\ell}{c_v\Theta_b},
\]
where $\mathcal T_0$ collects the other terms of $F_y$. By \eqref{e:gseed},
\eqref{eq:pr-gau} and $|\frac1{V_b}-\frac1{V_c}|\le C|Z_m|$, on $S_n$,
$|\mathcal T_1|\le C(\beta_n+\tau^{-\frac12})|(q,k,\ell)|$,
$|F-\frac{\kappa k}{c_vV_b}|\le C(\beta_n^2+\tau^{-1})$, and, with norms over
$J_n$ except for $Z_m$,
\[
\begin{gathered}
 \Vert\mathcal A\Vert_{L^2}^2+\Bigl\Vert\Bigl(\frac{(H_c)_y}{V_b}\Bigr)_yh\Bigr\Vert_{L^2}^2
 \le C(\tau^{-3}+\beta_n^2\tau^{-2})\Vert(Q,K,h)\Vert_{L^2}^2
   +C\tau^{-4}\Vert Z_m\Vert_{L^2(\R)}^2,\\
 \Bigl\Vert\partial_y^2\frac{(H_c)_y}{V_c}\Bigr\Vert_{L^2}^2\le C\tau^{-\frac72},\qquad
 \Vert\mathcal T_0\Vert_{L^2}^2\le C(\beta_n^4+\tau^{-2})\Vert(Q,K,h)\Vert_{L^2}^2
   +C\tau^{-3}\Vert Z_m\Vert_{L^2(\R)}^2+C\tau^{-\frac92}.
\end{gathered}
\]
Young's inequality in each integral of \eqref{e:ghe}, with the integration by
parts above in the last one and with the weight $\beta_n$ in $R\int_\R wq\mathcal A\dd y$
and in the $h\ell$ term, the bounds for $w\lambda_n$ and these estimates give,
for $t_n-2H_n\le t\le t_n$, since $\beta_nt_n\to\infty$ and $H_n=\beta_n^{-1}$,
\begin{equation}\label{e:ghineq}
 E_2'+\frac12D_2\le C(\beta_n+\tau^{-1})E_2
   +C\beta_n\Vert(q,k,\ell)\Vert_{L^2(J_n)}^2+C\beta_n^2\Vert(Q,K,h)\Vert_{L^2}^2
   +Ct_n^{-3}\Vert Z_m\Vert_{L^2}^2+Ct_n^{-\frac72}.
\end{equation}
By \eqref{e:gseed} there is $s_n\in[t_n-2H_n,r_n]$ with
$E_2(s_n)\le\frac C{H_n}\int_{t_n-2H_n}^{t_n}\Vert(q,k,\ell)\Vert_{L^2(J_n)}^2\dd t
=o(\beta_nt_n^{-1})$. On $[s_n,t_n]$, $\int C(\beta_n+\tau^{-1})\dd t\le C$, and the
time integral of the last four terms of \eqref{e:ghineq} is
$o(\beta_nt_n^{-1})+C\beta_n^2H_no(t_n^{-1})+CH_nt_n^{-3}=o(\beta_nt_n^{-1})$, by
\eqref{e:gseed} and Lemma~\ref{l:dep}. Gronwall's inequality, $w\ge C^{-1}$ on
$I_n$, \eqref{e:gconv}, its $y$-derivative and \eqref{e:gseed} give
\begin{equation}\label{e:ghess}
 \sup_{r_n\le t\le t_n}\bigl(\Vert\partial_y^2Z_m\Vert_{L^2(I_n)}^2
   +\Vert(R_U)_{yy}\Vert_{L^2(I_n)}^2+\Vert z_{yy}\Vert_{L^2(I_n)}^2\bigr)
 +\int_{r_n}^{t_n}\Vert\partial_y^3(\upsilon_m,\vartheta_m)\Vert_{L^2(I_n)}^2\dd t
 =o(\beta_nt_n^{-1}),
\end{equation}
where $R_U$ is converted by \eqref{e:wval},
$(R_U)_3=c_vW_T+b_aW_b+\frac12W_b^2$, and $z$ by
$z=(\varphi_m,b_m,c_v\vartheta_m+\frac12b_m^2)$.

For $\mathscr L\in\{L_\sigma,P_V\}$ let $T=T_{\mathscr L,\R\times(r_n,t_n)}$, and let
$\mathcal S_n$, $\xi_n$, $\tilde u$, $\tilde f$ be as in the proof of
Lemma~\ref{l:win}; Lemma~\ref{l:kern}(a) applies on $[r_n,t_n]$ with
$A_0=A_1=C\delta$. Let $u=Tf$ and $|f|\le F_0$, $|f_y|\le F_1$ where stated.
Then, with $C_j$ depending only on $j$ and the fixed data,
\begin{equation}\label{e:gker}
\begin{gathered}
 |\mathcal S_n[u_{yy}(\cdot,t_n)]|\le C\sup_{r_n<s<t_n}\bigl(\Vert f_1(s)\Vert_{L^2}
   +\Vert f_2(s)\Vert_{L^1}\bigr)\quad\text{if }f=f_1+f_2,\\
 |u_{yy}(y,t_n)|\le C\sqrt{H_n}F_1+C\delta\bigl(H_n\tau_{r_n}^{-1}
   +\sqrt{H_n}\tau_{r_n}^{-\frac12}\bigr)F_0,\\
 |\mathcal S_n[u_{yy}(\cdot,t_n)]|\le C_jH_n^{-j}\sup_{r_n<s<t_n}\bigl(
   \Vert f^0(s)\Vert_{L^1}+\Vert f^1(s)\Vert_{L^\infty}\bigr)
   \quad\text{if }f=f^0+\partial_yf^1,\\
 |v_{yy}(y,t_n)|\le CH_n^{-1}\Vert v(r_n)\Vert_{L^\infty}\quad\text{if }\mathscr Lv=0
   \text{ on }\R\times(r_n,t_n),\qquad
 |\mathcal S_n(g_{yy})|\le C_jH_n^{-j}\Vert g\Vert_{L^\infty},
\end{gathered}
\end{equation}
where in the third line $f^0,f^1$ vanish unless $2\bar cH_n\le|y-z_n|\le4\bar cH_n$,
and in the last $g=0$ on $(z_n-2\bar cH_n,z_n+2\bar cH_n)$. Indeed, by
\cref{l:kern}(a), $K_{\xi\xi}=-K_{\xi\eta}+(K_\xi+K_\eta)_\xi$ and
$\Vert b_\xi\Vert_{L^\infty}\le C\delta\tau_{r_n}^{-1}\le Ch^{-1}$,
\[
 |K_\xi|\le\frac{CG(Ch,X)}{\sqrt h},\quad |K_{\xi\xi}|\le\frac{CG(Ch,X)}h,\quad
 |K_{\xi\xi}+K_{\xi\eta}|\le\frac{C\delta G(Ch,X)}{\tau_{r_n}}
 \quad(h=t_n-s\le H_n).
\]
By \eqref{e:xder}, $u_{yy}=f_c^2\tilde u_{\xi\xi}+(f_c)_y\tilde u_\xi$, and the
definition of $\mathcal S_n$ gives
$\mathcal S_n[u_{yy}]=\mathcal S[f_c\tilde u_{\xi\xi}+\frac{(f_c)_y}{f_c}\tilde u_\xi](\xi_n)$,
with $f_c,(f_c)_y$ at $t_n$ as functions of $\xi$. Integrating by parts in the
difference $\mathcal S[f_c\tilde u_{\xi\xi}](\xi_n)-f_c(z_n,t_n)\mathcal S\tilde u_{\xi\xi}(\xi_n)$
and using $|\partial_\xi f_c|\le C\delta\tau^{-\frac12}$ and
$\int_\R(|\Psi|+|x\Psi'(x)|)^2\dd x<\infty$,
\[
 \bigl|\mathcal S_n[u_{yy}(\cdot,t_n)]-f_c(z_n,t_n)\mathcal S\tilde u_{\xi\xi}(\xi_n,t_n)\bigr|
 \le C\delta t_n^{-\frac12}\Vert u_y(t_n)\Vert_{L^2},
\]
and
$\mathcal S\tilde u_{\xi\xi}(\xi_n,t_n)=\int_{r_n}^{t_n}\!\int_\R\bigl(\int_\R
\Psi''(\xi_n-\xi)K(t_n,s;\xi,\eta)\dd\xi\bigr)\tilde f(\eta,s)\dd\eta\dd s$. The inner
integral has $L^\infty_\eta$ and $L^2_\eta$ norms at most $C$ and $Ch^{-\frac14}$
for $h\le1$, and, by the bound for $K_{\xi\xi}$, at most $Ch^{-\frac32}$ and
$Ch^{-\frac54}$ for $h>1$; with
$\Vert u_y(t_n)\Vert_{L^2}\le C\sup_s(H_n^{\frac12}\Vert f_1(s)\Vert_{L^2}
+H_n^{\frac14}\Vert f_2(s)\Vert_{L^1})$ and $H_n\le t_n$ this is the first line.
The second follows from
$\tilde u_{\xi\xi}=\int\!\!\int K_\xi\tilde f_\eta\dd\eta\dd s
+\int\!\!\int(K_{\xi\xi}+K_{\xi\eta})\tilde f\dd\eta\dd s$,
$|\tilde f_\eta|\le CF_1$ and $|(f_c)_y\tilde u_\xi|\le C\delta\tau^{-\frac12}
\sqrt{H_n}F_0$, the fourth from the bounds for $K_\xi$ and $K_{\xi\xi}$. The
third and the last follow as the third and the last lines of \eqref{e:Sn},
after two integrations by parts that put the derivatives on the kernel of
$\mathcal S$ and on $f_c$.

Fix $\sigma$. Let $\zeta_n,\eta_n,u_n,W_n,E_n$, $\mathfrak a$ and $Z$ be as in the
proof of Lemma~\ref{l:win}(b), now with the windows above, so that
\eqref{e:Erow} holds on $\Omega_n$, and $\sup_{t\ge r_n}\Vert E_n(t)\Vert_{L^\infty}\to0$
by Lemmas~\ref{l:all} and \ref{l:nf}(iii) and \eqref{e:nfe}. By \eqref{e:nfs},
\eqref{e:ghess} and \eqref{e:gseed},
\[
 \sup_{r_n<t<t_n}\bigl(\Vert\partial_y(\zeta_nF_{\rm bulk})\Vert_{L^2}
   +\Vert\zeta_nG_{\rm bulk}\Vert_{L^1}
   +\Vert\partial_y^2[\mathcal P_\sigma(V_c,z)]\Vert_{L^2(I_n)}\bigr)=o(\beta_n),
\]
so the first line of \eqref{e:gker} gives $|\mathcal S_n[\partial_y^2u_n(\cdot,t_n)]|=o(\beta_n)$.
Write $\eta_nE_n$ on $\R\times[r_n,t_n]$ as in the proof of
Lemma~\ref{l:win}(b). The initial value contributes $CH_n^{-1}\Vert E_n(r_n)\Vert_{L^\infty}
=o(\beta_n)$ to $\mathcal S_n[\partial_y^2(\eta_nE_n)(\cdot,t_n)]$, and the cutoff terms
$C_jH_n^{-j}$, by \eqref{e:gker}. By $-\mathfrak j_\sigma=\mathcal H_\sigma(V_c,\mathfrak a,Z)$
and \eqref{e:hex}, $|\mathfrak j_\sigma|\le C(|\mathfrak a|^2+|Z|^2)$ and
$|\partial_y\mathfrak j_\sigma|\le C(|\mathfrak a|+|Z|)(|\mathfrak a_y|+|Z_y|)
+C(|\mathfrak a|^2+|Z|^2)|(V_c)_y|$; by \eqref{e:rows}, \eqref{e:gseed} and
Lemma~\ref{l:pr}, and since $q_{-\sigma}$ and $e_m$ are smooth functions of
$(V_c,H_c,Z_m)$ vanishing at $Z_m=0$ (proof of Lemma~\ref{l:win}(a)),
$|\mathfrak a|+|Z|\le C_b\tau^{-\frac12}$ and $|\mathfrak a_y|+|Z_y|\le C\beta_n$ on $S_n$. Substituting
$\mathcal U_y=z_y+U_*'(V_c)(V_c)_y$ in \eqref{e:nfs} and expanding as in the proof of
\eqref{e:nfr} gives $R_{\rm c}=C_{\rm c}z_y+f_{\rm c}$, with
$|C_{\rm c}|\le C(|(V_c)_y|+|H_c|)$ and $f_{\rm c}$ a sum of products of smooth
functions of $(V_c,z,H_c)$ with $(V_c)_y$, $(V_c)_t$, $(V_c)_{yy}$, $(V_c)_y^2$ or
$(H_c)_y$, where the coefficient of $(V_c)_y$ vanishes at $(z,H_c)=0$. Hence the three terms
$\eta_n(\mathfrak j_\sigma+k(V_c)a_\sigma^2)$, $2\eta_nd(V_c)\mathfrak b_\sigma'(V_c)(V_c)_y\mathcal F$
and $\eta_nf_{\rm c}$ satisfy, on $\R\times(r_n,t_n)$,
\[
 |f|\le C_b\tau^{-1},\qquad |f_y|\le C_b\bigl(\beta_n\tau^{-\frac12}+\tau^{-\frac32}
   +H_n^{-1}\tau^{-1}\bigr),
\]
and the second line of \eqref{e:gker} bounds their contribution by
$C_b(\beta_n^{\frac12}t_n^{-\frac12}+\beta_n^{-\frac12}t_n^{-\frac32}
+\beta_n^{-1}t_n^{-2})=o(\beta_n)$. Since
$\Vert\eta_nR^{\rm prof}_\sigma\Vert_{L^2}\le C\delta e^{-c_et}\sqrt{H_n}$ and
$\Vert\eta_nC_{\rm c}z_y\Vert_{L^2}\le C\delta\tau^{-\frac12}\Vert z_y\Vert_{L^2}=o(t_n^{-1})$,
the first line of \eqref{e:gker} shows that $\eta_nR^{\rm prof}_\sigma$ and
$\eta_nC_{\rm c}z_y$ contribute $o(\beta_n)$.

Let $V_n=T(\eta_n\mathcal A_\sigma\Delta\mathcal G)$ with $\mathscr L=L_\sigma$. In the
coordinate $\xi$, $w_1=\partial_\xi\tilde V_n$ and $w_2=\partial_\xi w_1$ vanish at $r_n$
and satisfy
\[
 (\partial_t-D_0\partial_\xi^2+b\partial_\xi+b_\xi)w_1=\tilde f_\xi,\qquad
 (\partial_t-D_0\partial_\xi^2+b\partial_\xi+2b_\xi)w_2=\tilde f_{\xi\xi}-b_{\xi\xi}w_1,
\]
with $f=\eta_n\mathcal A_\sigma\Delta\mathcal G$. Since $|b_\xi|\le C\delta\tau^{-1}$,
$|b_{\xi\xi}|\le C\delta\tau^{-\frac32}$, $\tau_{t_n}\le2\tau_{r_n}$,
$f_{yy}=f_c^2\tilde f_{\xi\xi}+(f_c)_y\tilde f_\xi$ and $\partial_y^2V_n=f_c^2w_2+(f_c)_yw_1$,
comparison gives
$|\partial_y^2V_n(t_n)|\le8T|f_{yy}|(t_n)+C\delta t_n^{-\frac12}\sup_{r_n<s<t_n}T|f_y|(s)$.
Here
\[
 f_{yy}=\eta_n\bigl[(\partial_y^2\mathcal A_\sigma)\Delta\mathcal G
   +2(\partial_y\mathcal A_\sigma)\mathcal F+\mathcal A_\sigma\mathcal F_y\bigr]
   +2(\eta_n)_y\bigl[(\partial_y\mathcal A_\sigma)\Delta\mathcal G+\mathcal A_\sigma\mathcal F\bigr]
   +(\eta_n)_{yy}\mathcal A_\sigma\Delta\mathcal G,
\]
and, by \eqref{e:FRU} and \eqref{e:gbox},
$|\mathcal F_y|\le(\gamma-1+p_*)\cdot2\beta_n+\gamma p_*\sum_\rho|(\vartheta_\rho)_y|$ on
$\operatorname{supp}\eta_n$, where $|(\vartheta_\rho)_y|\le C_Me^{-c_Mt}$ by
Lemma~\ref{l:pr}, since $|y-s_\rho t|\ge\frac{c_\rho}2t$ there. By \eqref{e:Aresp},
Lemma~\ref{l:bd} ($|\Delta\mathcal G|\le C_b$), $|\mathcal F|\le C_b\tau^{-\frac12}$ and
the bound for $T|f_y|$ in the proof of Lemma~\ref{l:win}(b),
\[
 |\partial_y^2V_n(t_n)|\le16C_0\delta(\gamma-1+p_*)\beta_n+o(\beta_n),
\]
with $C_0$ the constant of \eqref{e:Aresp} for $i=0$, which depends only on the
fixed data. Altogether, with $C_0'$ depending only on the fixed data,
\[
 \bigl|\mathcal S_n[\partial_y^2(\eta_nE_n)(\cdot,t_n)]\bigr|\le C_0'\delta\beta_n+o(\beta_n).
\]
Differentiating the identity for $\ell_\sigma(V_c)R_U$ in the proof of
Lemma~\ref{l:win}(b) gives
\[
\begin{aligned}
 \ell_\sigma(V_c)(R_U)_y={}&\partial_y^2E_n+\partial_y^2W_n-\partial_y[\ell_\sigma(V_c)]R_U
   -\frac12\bigl((H_c)_y(y+h_*,t)-(H_c)_y\bigr)\\
 &-\partial_y\Bigl[\mathfrak b_\sigma'(V_c)(V_c)_y\Delta\mathcal G
   +\sum_{\rho=\pm1}a_\rho\ell_\sigma(V_c)\bigl(R_\rho(v_\rho)-R_\rho(V_c)\bigr)\Bigr].
\end{aligned}
\]
By \eqref{eq:pr-gau}, Lemma~\ref{l:pr}, \eqref{e:bdpot} and $(\Delta\mathcal G)_y=\mathcal F$,
the last three terms and $\partial_y[\ell_\sigma(V_c)]R_U$ are at most
$C_b\tau^{-1}+C_Me^{-c_Mt}$ on $\operatorname{supp}\eta_n$. Let
$g=\eta_n\ell_\sigma(V_c)(R_U)_y$ at $t=t_n$. By
\[
 \eta\,\partial_y^2X=\partial_y^2(\eta X)-2\partial_y(\eta_yX)+\eta_{yy}X,\qquad
 \partial_y^2X=\eta\,\partial_y^2X+\partial_y^2((1-\eta)X)+2\partial_y(\eta_yX)-\eta_{yy}X,
\]
with $(\eta,X)=(\eta_n,E_n)$, $(\eta_n,W_n)$ and $(\zeta_n,\mathcal P_\sigma(V_c,z))$, the
bounds for $\partial_y^2(\eta_nE_n)$, $\partial_y^2u_n$ and
$\partial_y^2[\mathcal P_\sigma(V_c,z)]$ above, the last lines of \eqref{e:gker} and
\eqref{e:Sn}, and $|\mathcal S_ng|\le C\min\{\Vert g\Vert_{L^\infty},\Vert g\Vert_{L^2}\}$,
$|\mathcal S_ng|\le C_0'\delta\beta_n+o(\beta_n)$. Since
$\Vert g_y\Vert_{L^2}+t_n^{-\frac12}\Vert g\Vert_{L^2}=o(\beta_n)$ by \eqref{e:ghess} and
\eqref{e:gseed}, the fifth line of \eqref{e:Sn} gives
\[
 |\ell_\sigma(V_c)(R_U)_y(z_n,t_n)|\le\sqrt2C_0'\delta\beta_n+o(\beta_n).
\]

For the third row let $u_N,C_N,D_N,e_N,\Delta_d,\Delta_b$ be as in the proof of
Lemma~\ref{l:win}(b), so that \eqref{e:Nrow} holds on $\Omega_n$. By \eqref{e:npr} and
\eqref{e:gseed},
$\int_{r_n}^{t_n}\Vert\zeta_nF_N\Vert_{L^2}^2\dd t=o(H_nt_n^{-1})$ and
$\int_{r_n}^{t_n}\Vert(\zeta_nF_N)_y\Vert_{L^2}^2\dd t=o(t_n^{-1})$, and by \eqref{e:ghess}
(the highest term being a bounded multiple of $\partial_y^3\upsilon_m$)
$\int_{r_n}^{t_n}\Vert(\zeta_nF_N)_{yy}\Vert_{L^2}^2\dd t=o(\beta_nt_n^{-1})$. The two
lines of \eqref{e:pchie} give, as in the proof of Lemma~\ref{l:win}(b),
\[
 \sup_{r_n\le t\le t_n}\Vert u_N\Vert_{L^2}^2+\int_{r_n}^{t_n}\Vert(u_N)_y\Vert_{L^2}^2\dd t
 =o(H_nt_n^{-1}),\qquad
 \sup_{r_n\le t\le t_n}\Vert(u_N)_y\Vert_{L^2}^2+\int_{r_n}^{t_n}\Vert(u_N)_{yy}\Vert_{L^2}^2\dd t
 =o(t_n^{-1}),
\]
hence $\sup_{r_n\le t\le t_n}\Vert D_N(t)\Vert_{L^\infty}=o(1)$ as there. Since
$P_\chi u_N=\partial_y(\zeta_nF_N)$,
\[
 \frac12\frac{\dd}{\dd t}\Vert(u_N)_{yy}\Vert_{L^2}^2+\int_\R d_0(u_N)_{yyy}^2\dd y
 =\int_\R\bigl[-(d_0)_y(u_N)_{yy}+(b_0)_y(u_N)_y+b_0(u_N)_{yy}
   -(\zeta_nF_N)_{yy}\bigr](u_N)_{yyy}\dd y .
\]
With $(d_0)_y=-d_0(g_\chi+\frac{H_c}c)$, $\int_{r_n}^{t_n}\Vert g_\chi\Vert_{L^\infty}^2\dd t\le C_b$,
$|b_0|\le C\delta\tau^{-\frac12}$ and $|(b_0)_y|\le C\delta\tau^{-1}+C\delta\tau^{-\frac12}|g_\chi|$,
Young's and Gronwall's inequalities give
$\sup_{r_n\le t\le t_n}\Vert(u_N)_{yy}\Vert_{L^2}^2=o(\beta_nt_n^{-1})$. Hence, by
\eqref{e:ghess} and $|\partial_y^2(H_c\mathcal B)|\le|(H_c)_{yy}\mathcal B|
+C_b\delta\tau^{-\frac32}+C\delta\tau^{-\frac12}|(b_m)_y|$,
$\sup_{r_n\le t\le t_n}\Vert(C_N)_{yy}\Vert_{L^2(I_n)}=o(\beta_n)$, and, by \eqref{e:wtC},
$\Vert(D_N)_{yy}\Vert_{L^2}\le\Vert(e_N)_y\Vert_{L^2}+\Vert(C_N)_{yy}\Vert_{L^2}=o(t_n^{-\frac12})$.
In \eqref{e:npr} substitute $(V_b)_y=(V_c)_y+(\varphi_m)_y$, $(b_m)_y=(H_c)_y+(\upsilon_m)_y$,
$(\chi_m)_y=(V_c)_y+(\chi_m-V_c)_y$ and
$g_\chi=\frac{(\chi_m-V_c)_y}{\chi_m}-\frac{(\chi_m-V_c)(V_c)_y}{V_c\chi_m}$, and write
$R_N=R_N^0+R_N^1+R_N^2$, where $R_N^1$ collects the terms with one factor
$(\varphi_m)_y$, $(\upsilon_m)_y$ or $(\chi_m-V_c)_y$ and
$R_N^2=-\frac{(\gamma-1)cd_0(\chi_m-V_c)_y(\upsilon_m)_y}{\gamma p_*\chi_m}$. On $S_n$, by
\eqref{e:gseed} and \eqref{eq:pr-gau}, $|R_N^0|\le C_b\tau^{-1}$ and
$|\partial_yR_N^0|\le C_b(\beta_n\tau^{-\frac12}+\tau^{-\frac32})$, the coefficients in
$R_N^1$ are at most $C_b\tau^{-\frac12}$, so that $\Vert\eta_nR_N^1\Vert_{L^2}=o(t_n^{-1})$,
and $\Vert R_N^2\Vert_{L^1}\le C\Vert\partial_y(\mathcal U-\widehat{\mathcal U})\Vert_{L^2}^2=o(t_n^{-1})$;
the term $(1-\frac{V_c}{\chi_m})h_*(V_c)_t$ is at most $C_b\tau^{-\frac32}$ with
derivative at most $C_b(\beta_n\tau^{-1}+\tau^{-2})$; and
$\Vert\Delta_d(D_N)_{yy}+\Delta_b(D_N)_y\Vert_{L^2}\le C_b\tau^{-\frac12}
\Vert(D_N)_{yy}\Vert_{L^2}+C_b\tau^{-1}\Vert(D_N)_y\Vert_{L^2}=o(t_n^{-1})$, since
$\partial_y(\Delta_d(D_N)_y)+(\Delta_b-(\Delta_d)_y)(D_N)_y=\Delta_d(D_N)_{yy}+\Delta_b(D_N)_y$.
Writing $\eta_nD_N$ as $\eta_nE_n$, now with $\mathscr L=P_V$, \eqref{e:gker} gives
$|\mathcal S_n[\partial_y^2(\eta_nD_N)(\cdot,t_n)]|=o(\beta_n)$. By the identity for
$\ell_0R_U$ in the proof of Lemma~\ref{l:win}(b),
\[
 \ell_0(R_U)_y=(D_N)_{yy}-(C_N)_{yy}-\bigl((V_c)_y(y+h_*,t)-(V_c)_y-h_*(V_c)_{yy}\bigr),
\]
where the last term is at most $\frac{h_*^2}2\Vert\partial_y^3V_c\Vert_{L^\infty}\le C_b\tau^{-\frac32}$,
and \eqref{e:Sn} applied to $g=\eta_n\ell_0(R_U)_y$ at $t=t_n$ gives
$|\ell_0(R_U)_y(z_n,t_n)|=o(\beta_n)$.

By \eqref{e:rinv} with $n=V_c(z_n,t_n)$, the bounds for $\ell_\pm(V_c)(R_U)_y$ and
$\ell_0(R_U)_y$ at $(z_n,t_n)$ give $\beta_n\le C_1'\delta\beta_n+o(\beta_n)$, with
$C_1'$ depending only on the fixed data. Decreasing $\delta_1(M)$ so that
$2C_1'\delta_1(M)\le1$, this is a contradiction.
\end{proof}

\subsection{The diagonal system}

For late $r$ put, on $\R\times[r,\infty)$, with $B_{\rm ref}=(B_{{\rm ref},ik})$
and $B_D$ of \eqref{e:pB},
\begin{equation}\label{e:frm}
\begin{gathered}
 \lambda_\pm=\frac{\operatorname{tr}B_{\rm ref}\pm\sqrt{(B_{{\rm ref},11}
   -B_{{\rm ref},22})^2+4B_{{\rm ref},12}B_{{\rm ref},21}}}2,\qquad
 \mathcal R_j=\begin{pmatrix}1\\[2pt]
   \frac{\lambda_j-B_{{\rm ref},11}}{B_{{\rm ref},12}}\end{pmatrix},\\
 \mathcal R=(\mathcal R_+,\mathcal R_-),\qquad \mathcal P=S\mathcal R,\qquad
 \Lambda_{\rm ref}=\operatorname{diag}(\lambda_+,\lambda_-),
\end{gathered}
\end{equation}
\begin{equation}\label{e:drf}
 \Delta=\operatorname{diag}(\mathcal M-M_{\rm ref}),\qquad
 \mathcal Z=\mathcal P^{-1}\Delta\mathcal P,\qquad
 b^D_j=\lambda_j-\mathcal Z_{jj}.
\end{equation}
By \eqref{e:gop} and \eqref{e:pB}, $SB_{\rm ref}S^{-1}=\Lambda-M_{\rm ref}-T_{\rm ref}$.
By \eqref{e:jsq}, $\pa_\zeta\mathfrak j^\sharp_\sigma(n_0,n_0,0)=0$, so $M_{\rm ref}$ is a
smooth function of $(V_c,V_h-V_c,\mathfrak a^a)$ vanishing at $(V_c,0,0)$; with
$V_h=V_c(y+h_*,t)$, $\delta|h_*|\le C_M$, $\tau_r\ge h_*^2$, \eqref{eq:pr-gau},
Lemma~\ref{l:pr} and \eqref{e:pT}, $M_{\rm ref}+T_{\rm ref}$ and its derivatives
$\pa_y^i\pa_t^l$, $i+2l\le2$, are bounded by $C_M\tau^{-\frac{1+i+2l}2}$ times Gaussians
of width $\sqrt\tau$ centred at $0$, $-h_*$ and $s_\pm t$. Hence $B_{{\rm ref},12}\ge\frac1{2V_c}$,
$\lambda_\pm$ are real and $|\lambda_j-s_j(V_c)|\le C_M\tau^{-\frac12}$, $B_{\rm ref}\mathcal R
=\mathcal R\Lambda_{\rm ref}$, and $\mathcal P$ is a smooth function of
$(V_c,M_{\rm ref}+T_{\rm ref})$ with $\mathcal P=I$ at $M_{\rm ref}+T_{\rm ref}=0$, since the
columns $(1,\pm\frac1{2\mathfrak b_+(V_c)})^\top$ of $S^{-1}$ are the eigenvectors of $A_{\rm ph}$
with first component $1$; so $|\mathcal P-I|\le C_M\tau^{-\frac12}$. By \eqref{e:pB},
$B_D\mathcal R=\mathcal R(\Lambda_{\rm ref}-\mathcal Z)$. By \eqref{e:pj}, $\mathcal M-M_{\rm ref}$
is a smooth function of $(V_c,V_h,\mathfrak a^a,\chi_m,\mathfrak a)$ vanishing at
$(\chi_m,\mathfrak a)=(V_h,\mathfrak a^a)$; with $\chi_m-V_h=\ell_0R_U$,
$\mathfrak a_\sigma-\mathfrak a^a_\sigma=\ell_\sigma(V_c)R_U$ (proof of Lemma~\ref{l:am}),
Lemmas~\ref{l:asy}, \ref{l:ga} and \eqref{e:wval},
\begin{equation}\label{e:dDel}
 |\Delta|\le C|R_U|\le C_b\tau^{-\frac12},\qquad
 |\Delta_y|\le C|(R_U)_y|+C_b\tau^{-1}\le C_b,\qquad
 |\Delta_y|\le C_b\tau^{-1}\quad\Bigl(|y|\le\frac{c_-}{1024}t\Bigr).
\end{equation}

Let $Q$ be a weak solution of $Q_t=\pa_y(dQ_y-B_DQ)$, $d=d(V_c)$, with
$Q\in C([r,T];L^2)$, $Q_y\in L^2((r,T)\times\R)$, and put $u=\mathcal R^{-1}Q$,
$K^{\rm fr}_y=\mathcal R^{-1}\mathcal R_y$, $K^{\rm fr}_t=\mathcal R^{-1}\mathcal R_t$. Since
$(\mathcal Ru)_y=\mathcal R(u_y+K^{\rm fr}_yu)$ and $B_D\mathcal R=\mathcal R(\Lambda_{\rm ref}
-\mathcal Z)$, the equation is
$\mathcal Ru_t+\mathcal R_tu=\pa_y\{\mathcal R[du_y+dK^{\rm fr}_yu-(\Lambda_{\rm ref}-\mathcal Z)u]\}$;
multiplying by $\mathcal R^{-1}$ and using
$dK^{\rm fr}_yu_y=\pa_y(dK^{\rm fr}_yu)-\pa_y(dK^{\rm fr}_y)u$,
\begin{equation}\label{e:dsc}
 (u_j)_t=\pa_y\bigl\{d(u_j)_y-b^D_ju_j\bigr\}+g_j,
\end{equation}
\begin{equation}\label{e:dsrc}
\begin{gathered}
 g=\pa_y(\mathsf Eu)+\mathsf Fu,\qquad
 \mathsf E=2dK^{\rm fr}_y+\operatorname{offdiag}\mathcal Z,\\
 \mathsf F=d(K^{\rm fr}_y)^2-\pa_y(dK^{\rm fr}_y)-K^{\rm fr}_y(\Lambda_{\rm ref}-\mathcal Z)
   -K^{\rm fr}_t .
\end{gathered}
\end{equation}
Let $\beta_j$ be the drift $b$ of \eqref{e:nco} with $B=b^D_j$, and
$z_j(F_c(y,t),t)=\frac{u_j(y,t)}{f_c(y,t)}$. By \eqref{e:nco},
\[
 (z_j)_t=\pa_\xi\bigl(D_0(z_j)_\xi-\beta_jz_j\bigr)+\pa_\xi(\tilde{\mathsf E}z)_j
   +(\tilde{\mathsf F}z)_j,\qquad
 \tilde{\mathsf E}(F_c,t)=f_c\mathsf E,\quad \tilde{\mathsf F}(F_c,t)=\mathsf F,
\]
and $\Vert z(t)\Vert_{L^1_\xi}=\Vert u(t)\Vert_{L^1}$. Put $v_0=2^{-14}c_-$ and, for a
fixed $b_w>0$ depending only on the fixed data,
\[
 g_0=e^{-\frac{a_s\xi^2}{2\tau}},\qquad g_\pm=e^{-\frac{b_w(\xi\mp c_-t)^2}\tau},
 \qquad c_-=s_+(\rho_-^{-1}),
\]
so that $f_cs_j(V_c)=jc_-$. Write $\tilde{\mathsf E}=\mathsf E^c+\mathsf E^r+J$ and
$\tilde{\mathsf F}=\mathsf F^c+\mathsf F^r$, where $J(F_c,t)=f_c\operatorname{offdiag}\mathcal Z$
and $\mathsf E^c(F_c,t)$, $\mathsf F^c(F_c,t)$ are $2f_cdK^{\rm fr}_y$ and
$-K^{\rm fr}_y\Lambda_{\rm ref}-K^{\rm fr}_t$ computed with $a_+=a_-=0$ in $M_{\rm ref}$.
For late $r$ and $j,k=\pm1$,
\begin{equation}\label{e:dfco}
\begin{gathered}
 |\beta_j-jc_-|\le C_M\tau^{-\frac12},\qquad |\pa_\xi\beta_j|\le C_b,\qquad
 |\pa_\xi\beta_j|+\tau^{\frac12}|\pa_\xi J|\le C_b\tau^{-1}\quad(|\xi|\le v_0\tau),\\
 |\mathsf E^c|+|\mathsf F^c|+\tau^{\frac12}\bigl(|\pa_\xi\mathsf E^c|+|\pa_\xi\mathsf F^c|\bigr)
   \le C_M\delta\tau^{-\frac12}g_0,\qquad
 |\mathsf F^r_{ik}|\le C_b\bigl[\tau^{-1}(g_0+g_{-k})+\tau^{-\frac32}\bigr],\\
 |\mathsf E^r|+|J|+\tau^{-\frac14}\bigl(\Vert\mathsf E^r\Vert_{L^2}+\Vert J\Vert_{L^2}\bigr)
   +\tau^{\frac14}\Vert\pa_\xi\mathsf E^r\Vert_{L^2}\le C_b\tau^{-1},\qquad
 \Vert\pa_\xi J\Vert_{L^2}\le C_b\tau^{-1},
\end{gathered}
\end{equation}
where $\mathsf F^r_{ik}$ is the entry acting on $z_k$. Indeed,
$\beta_j-jc_-=f_c(\lambda_j-s_j(V_c)-\mathcal Z_{jj})+(F_c)_t-d(f_c)_y$; the bounds after
\eqref{e:nco}, $|\mathcal Z|\le C|\Delta|$, $|\pa_y\mathcal Z|\le C(|\Delta_y|+|\mathcal P_y||\Delta|)$
and \eqref{e:dDel} give the first line, since $|\xi|\le v_0\tau$ implies
$|y|\le\sqrt2|\xi|\le\frac{c_-}{1024}t$. In the coefficients $\mathsf E^c,\mathsf F^c$ only $V_c$,
$V_h-V_c$, $H_h$ and $T_{\rm ref}$ occur; $\pa_y$ of $V_c$ costs $C\delta\tau^{-\frac12}$, of the
other three $C_M\tau^{-1}\le C_M\delta\tau^{-\frac12}$ for $\tau\ge\delta^{-2}$, and the Gaussians
are bounded by $Cg_0$ because $\xi^2\le2y^2$, $\tau\ge h_*^2$ and $\frac{c_G}8\ge\frac{a_s}2$.
The other parts of $K^{\rm fr}_y$, $K^{\rm fr}_t$ contain a derivative of $a_\pm$ or a
product of two coefficients of size $C_M\tau^{-\frac12}$. A Gaussian of width $\sqrt\tau$
has $L^2$ norm of order $\tau^{1/4}$, and each spatial derivative contributes
$\tau^{-1/2}$; this gives the bounds for $\mathsf E^r$. In column $k$ of
$K^{\rm fr}_y\Lambda_{\rm ref}+K^{\rm fr}_t$, the derivative of $\mathcal R_k$ in its argument
$a_k$ multiplies $(\pa_t+\lambda_k\pa_y)a_k=(\lambda_k-s_k)(a_k)_y+d_k(a_k)_{yy}-2k_ka_k(a_k)_y$
(Lemma~\ref{l:pr}), which is at most $C_M\tau^{-\frac32}$ because $\lambda_k-s_k$ is at most
$C_M\tau^{-\frac12}$ near $y=s_kt$ and the profile is exponentially small elsewhere; the
argument $a_{-k}$ gives $(\pa_t+\lambda_k\pa_y)a_{-k}$, bounded by $C_M\tau^{-1}g_{-k}$. This is
the third bound. Every entry of $\operatorname{offdiag}\mathcal Z$ contains $\mathcal P-I$, so
$|J|\le C_M\tau^{-\frac12}|\Delta|$ with Gaussian envelopes, and $|\pa_\xi J|\le
C_M(\tau^{-1}|\Delta|+\tau^{-\frac12}|\Delta_y|)$; with \eqref{e:dDel} and
$\Vert(R_U)_y\Vert_{L^2}\le C_b\tau^{-\frac12}$ (\eqref{e:wval}, \eqref{e:wtC}) this gives the
global bounds for $J$. On $|\xi|\le v_0\tau$, the sharper bound
$|\Delta_y|\le C_b\tau^{-1}$ in \eqref{e:dDel} gives $|J_\xi|\le C_b\tau^{-3/2}$.
Finally, $\beta_j$ is a smooth function of $\mathcal U$ and of the reference
fields, and by Lemma~\ref{l:bg} $\mathcal U$ and $\mathcal U_y$ are bounded and uniformly
H\"older continuous on $\R\times J'$ for compact $J'\subset[r,\infty)$ (the H\"older bound in
$t$ follows from $\Vert f\Vert_{L^\infty}^2\le\Vert f\Vert_{L^2}\Vert f_y\Vert_{L^2}$ with
$f=\mathcal U_y(t)-\mathcal U_y(t')$ and $\int\Vert\pa_y\pa_t\mathcal U\Vert_{L^2}^2<\infty$ on
$J'$). Thus $\beta_j$ satisfies the hypotheses of \cref{l:kernv}(ii) on $I=[r,\infty)$ with
$\bar b=jc_-$, $A_0=C_M$, $A_1=A_2=C_b$ and $v_0$, and $v_0\le\frac{c_-}4$. Let $H_j$ be its
kernel $H$.

\begin{lemma}\label{l:dp}
For every sufficiently late fixed $r$ and every $p\in L^1(\R;\R^2)\cap L^2(\R;\R^2)$, the
problem
\[
 \pa_tQ_D=\pa_y\bigl\{d(V_c)(Q_D)_y-B_DQ_D\bigr\},\qquad Q_D(r)=p,
\]
has a unique solution with $Q_D\in C([r,T];L^1\cap L^2)$ and
$(Q_D)_y\in L^2((r,T)\times\R)$ for every finite $T$, and, with
$L_p=\Vert p\Vert_{L^1}+\Vert p\Vert_{L^2}$,
\begin{equation}\label{e:dpb}
\begin{gathered}
 \sup_{t\ge r}\Vert Q_D(t)\Vert_{L^1}\le C_{b,r}\Vert p\Vert_{L^1},\qquad
 \Vert Q_D(t)\Vert_{L^2}\le C_{b,r}L_p\tau^{-\frac14},\\
 \int_r^\infty\tau^{\frac38}\Vert(Q_D)_y\Vert_{L^2}^2\dd t+\mathscr Q_r(Q_D)
 \le C_{b,r}L_p^2,
\end{gathered}
\end{equation}
\begin{equation}\label{e:dpn}
 \mathscr N_{r,\infty}\Bigl(\frac{(Q_D)_2}{V_c}\Bigr)\le C_{b,r}L_p .
\end{equation}
\end{lemma}

\begin{proof}
Since $|B_D|\le C$, the energy identity
$\frac12\frac{\dd}{\dd t}\Vert Q\Vert_{L^2}^2+\int_\R d|Q_y|^2\dd y=\int_\R Q_y^\top B_DQ\dd y$
and Galerkin approximation give a solution with $Q_D\in C([r,T];L^2)$,
$(Q_D)_y\in L^2((r,T)\times\R)$; it is unique by \eqref{e:uniq}. In $\xi=F_c(y,t)$,
$\frac{Q_D}{f_c}$ solves the equation of \cref{l:kernv}(iv) with the bounded matrix
$f_cB_D+((F_c)_t-d(f_c)_y)I$, so $Q_D\in C([r,T];L^1)$. Let $u$, $g$, $z$ be as above with
$Q=Q_D$.

For $r\le T<\infty$ put
\[
 \mathsf X_T=\sup_{r\le t\le T}\Vert z(t)\Vert_{L^1},\qquad
 \mathsf O_T=\sum_j\int_r^T\!\!\int_\R\tau^{-\frac12}(g_0+g_{-j})|z_j|\dd\xi\dd t,
\]
which are finite. Write $h=t-s$ and $X=\xi-\eta-jc_-h$. By \cref{l:kernv}(ii) and
\cref{l:kern}(b), $H_j\le CG(Ch,X)$ and $|\pa_\eta H_j|\le C(h^{-\frac12}+\tau_s^{-\frac12})G(Ch,X)$
with $C=C_M$; since $\int_\R H_j\dd\xi=1$, $\int_\R|\pa_\eta H_j|\dd\xi\le C_M h^{-\frac12}$
for $h\le\tau_s$ and $\le C_M\tau_s^{-\frac12}$ for $h\ge\tau_s$. Gaussian convolution gives
$\int_\R\tau_t^{-\frac12}(g_0+g_{-j})(\xi,t)G(Ch,X)\dd\xi\le C\tau_t^{-\frac12}
(e^{-\frac{c_*(\eta+jc_-h)^2}{\tau_t}}+e^{-\frac{c_*(\eta+jc_-h+jc_-t)^2}{\tau_t}})$, and both
centres move with speed at least $c_-$ in $h$; hence
\begin{equation}\label{e:docc}
 \sup_{\eta}\int_s^\infty\!\!\int_\R\tau^{-\frac12}(g_0+g_{-j})H_j\dd\xi\dd t\le C_M,\qquad
 \sup_{\eta}\int_s^\infty\!\!\int_\R\tau^{-\frac12}(g_0+g_{-j})|\pa_\eta H_j|\dd\xi\dd t
 \le C_M\tau_s^{-\frac14} .
\end{equation}
For the second bound, the ages $h\le\sqrt{\tau_s}$ contribute
$C_M\tau_s^{-\frac12}\int_0^{\sqrt{\tau_s}}h^{-\frac12}\dd h$, the ages
$\sqrt{\tau_s}<h\le\tau_s$ contribute $\tau_s^{-\frac14}$ times the first bound, and for
$h\ge\tau_s$ one uses $\pa_\eta H_j(t,s)=\int_\R H_j(t,m;\cdot,z)\pa_\eta H_j(m,s;z,\cdot)\dd z$,
$m=s+\frac{\tau_s}2$, $\int_\R|\pa_\eta H_j(m,s;z,\eta)|\dd z\le C_M\tau_s^{-\frac12}$ and the
first bound. By Duhamel's formula for the equation of $z_j$ (which holds by \eqref{e:uniq}),
\[
 z_j(t)=\int_\R H_j(t,r;\cdot,\eta)z_j(\eta,r)\dd\eta
 +\int_r^t\!\!\int_\R\bigl[H_j(\tilde{\mathsf F}z)_j-\pa_\eta H_j(\tilde{\mathsf E}z)_j\bigr]
 (\eta,s)\dd\eta\dd s .
\]
By \eqref{e:dfco},
\[
\begin{gathered}
 \int_r^T\Vert\mathsf F^cz\Vert_{L^1}\dd t\le C_M\delta\mathsf O_T,\qquad
 \int_r^T\Vert\mathsf F^rz\Vert_{L^1}\dd t\le C_b\tau_r^{-\frac12}(\mathsf O_T+\mathsf X_T),\\
 \Vert\mathsf E^cz(s)\Vert_{L^1}\le C_M\delta\tau_s^{-\frac12}
   \min\Bigl\{\Vert z(s)\Vert_{L^1},\int_\R g_0|z(s)|\dd\xi\Bigr\},\qquad
 \Vert(\mathsf E^r+J)z(s)\Vert_{L^1}\le C_b\tau_s^{-1}\Vert z(s)\Vert_{L^1}.
\end{gathered}
\]
With $\int_{t-\tau_s\le s\le t}h^{-\frac12}\tau_s^{-\frac12}\dd s\le C$ and
$\int_r^t(h^{-\frac12}\mathbf1_{h\le\tau_s}+\tau_s^{-\frac12})\tau_s^{-1}\dd s\le C\tau_r^{-\frac12}$,
the $L^1$ norm of the Duhamel formula and the occupation bounds \eqref{e:docc} give
\[
 \mathsf X_T+\mathsf O_T\le C_M\Vert z(r)\Vert_{L^1}
 +\bigl(C_M\delta+C_b\tau_r^{-\frac14}\bigr)(\mathsf X_T+\mathsf O_T).
\]
Decreasing $\delta_1(M)$ and increasing $r$ so that the bracket is at most $\frac12$, we get
$\mathsf X_\infty+\mathsf O_\infty\le C_M\Vert p\Vert_{L^1}$, and the first bound in \eqref{e:dpb}.

We use the following identities. If $\mathcal K_{\mathcal M}\mathbf H=F$ with $\mathbf H$ as in Lemma~\ref{l:geo}, then
$\mathbf v=S^{-1}\mathbf H$, $\mathbf q=\mathbf v_y$ and $\mathbf z=S\mathbf q$ satisfy,
with $J_S=S_yS^{-1}$ and $T_S=S_tS^{-1}$,
\begin{equation}\label{e:exder}
\begin{aligned}
 \mathbf q_t&=\pa_y\bigl[d(V_c)\mathbf q_y-A_{\rm ph}\mathbf q
   +S^{-1}\mathcal MS\mathbf q+S^{-1}F\bigr],\\
 \mathbf z_t-d\mathbf z_{yy}+(\Lambda-d_yI)\mathbf z_y+\Lambda_y\mathbf z
 &=\pa_y\bigl[(\mathcal M-dJ_S)\mathbf z+F\bigr]-dJ_S\mathbf z_y\\
 &\quad+(T_S+J_S\Lambda-J_S\mathcal M+dJ_S^2)\mathbf z-J_SF,
\end{aligned}
\end{equation}
where $d=d(V_c)$. The first line is the $y$-derivative of the equation for
$\mathbf v$; for the second, $\mathbf q=S^{-1}\mathbf z$,
$\mathbf q_t=S^{-1}(\mathbf z_t-T_S\mathbf z)$,
$\mathbf q_y=S^{-1}(\mathbf z_y-J_S\mathbf z)$ and $A_{\rm ph}S^{-1}=S^{-1}\Lambda$
turn the first line into $\mathbf z_t-T_S\mathbf z=Y_y-J_SY$ with
$Y=d(\mathbf z_y-J_S\mathbf z)-\Lambda\mathbf z+\mathcal M\mathbf z+F$.

Put $x=y+h_*$, $\mathbf z=SQ_D$ and $\mathcal M_D=M_{\rm ref}+T_{\rm ref}+\Delta$,
with $J_S,T_S$ as in \eqref{e:exder}. By \eqref{e:pB}, $B_D=A_{\rm ph}-S^{-1}\mathcal M_DS$.
The computation after \eqref{e:exder}, with $F=0$ and $\mathcal M_D$ in place of $\mathcal M$, gives
\[
 \mathbf z_t=\pa_y\bigl[d\mathbf z_y-\operatorname{diag}(\tilde s_j)\mathbf z\bigr]+f,\qquad
 f=\pa_y\bigl[(\mathcal N-dJ_S)\mathbf z\bigr]-dJ_S\mathbf z_y
   +(T_S+J_S\Lambda-J_S\mathcal M_D+dJ_S^2)\mathbf z,
\]
with $\mathcal N=\mathcal M_D+2k(V_c)\operatorname{diag}(a_+,a_-)$. By \eqref{e:pmat} and
$|\mathcal M-M_{\rm ref}|\le C|R_U|$, $\mathcal M_D=\widetilde{\mathcal M}
-\operatorname{offdiag}(\mathcal M-M_{\rm ref})$ has the form \eqref{e:gmat} with a residual
$\mathfrak r$, $|\mathfrak r|\le C_M|R_U|+C_b(\tau^{-1}+\tau^{-\frac12}G_0)$, so that
$\mathcal N_{jj}=2k(V_c)a_{-j}+\mathfrak r_{jj}$ and
$\mathcal N_{j,-j}=2k(V_c)a_j+2\ell_{\rm g}(V_c)a_{-j}+\mathfrak r_{j,-j}$. Since $S$ depends
only on $V_c$, $|J_S|+\tau^{\frac12}|T_S|\le C\delta W_c$ for $\tau\ge h_*^2$.

Let $h_j=1-2k_j\tau(a_j)_y$. By \eqref{e:CH}, $h_j=1+2d_j\tau\pa_y^2\log p_j$ is a function
of $z_j$, and $p_j(y,t)=\int_\R\Pi(y-s_jt-x)G_j(x)\dd x$, with $\Pi=1$ on $(-\infty,0)$,
$\Pi=e^{-\alpha s_jm_j}$ on $(0,\infty)$ and $G_j$ the Gaussian density of variance
$2d_j\tau$. Two differentiations give $h_j=\frac{\operatorname{Var}\mu}{2d_j\tau}$ for the
probability measure $\mu\propto\Pi(y-s_jt-x)G_j(x)\dd x$, whose density with respect to $G_j$
lies between $e^{-\alpha|s_jm_j|}$ and $e^{\alpha|s_jm_j|}$; hence $C_M^{-1}\le h_j\le C_M$. By
$P_jA_j=-k_ja_j^2$ (Lemma~\ref{l:pr}),
$(a_j)_t=\pa_y[d_j(a_j)_y-s_ja_j-k_ja_j^2]$, and a direct computation gives
\[
\begin{gathered}
 (h_j)_t=\pa_y\bigl[d_j(h_j)_y-(s_j+2k_ja_j)h_j\bigr],\\
 \mathsf d_j:=(h_j)_t-\pa_y\bigl[d(h_j)_y-\tilde s_jh_j\bigr]
 =\pa_y\bigl[(d_j-d)(h_j)_y+(\tilde s_j-s_j-2k_ja_j)h_j\bigr].
\end{gathered}
\]
Fix $\theta_c=\frac{c_-}{64}$. On $|x|\ge\frac{\theta_c}2\tau$ either $V_c-v_j$ and its
derivatives, or $h_j-1$ and its derivatives, are at most $C_be^{-c_bt}$; so
$|\mathsf d_j|\le C_be^{-c_bt}$ there. With $b_{\rm ac}$ of Appendix~\ref{s:pres}, put
$Q_j=\tau^{-1}e^{-b_{\rm ac}z_j^2}$ and $W_j=\tau^{-\frac12}e^{-b_{\rm ac}z_j^2}$; then
$|(h_j)_y|^2\le C_MQ_j$, $|a_j|^2+|(a_j)_y|\le C_MQ_j$ and
$|d_y|^2\le C\delta^2\tau^{-\frac12}W_c$.

For $0\le\theta\le1$ let $P_{j,\theta}=\pa_t-d\pa_y^2
+[\tilde s_j-\theta(d_y+2d\frac{(h_j)_y}{h_j})]\pa_y$, so that $P_{j,0}=\mathscr P_j$. Then
\begin{equation}\label{e:dres}
 \bigl\Vert T_{P_{j,\theta},r}(Q_j+W_{-j}+W_c)\bigr\Vert_{L^\infty}\le C_M .
\end{equation}
For $W_{-j}$ and $W_c$ this is \eqref{e:gauss} with $\sigma=j$ and $v=s_{-j}$, respectively
$v=0$, $h=h_*$, since the drift of $P_{j,\theta}$ is $s_j(V_c)+O(C_M\tau^{-\frac12})$. For
$Q_j$, let $P^e_{j,\theta}$ be $P_{j,\theta}$ with $d$, $s_j(V_c)$, $k(V_c)$ replaced by $d_j$,
$s_j$, $k_j$ and $d_y$ by $0$. Since $2k_ja_j=-2d_j\frac{(p_j)_y}{p_j}$ and $p_j$, $h_j$ are
functions of $z_j$, $P^e_{j,\theta}[\Phi(z_j)]=\tau^{-1}[-d_j\Phi''
-(\frac z2+2d_j\frac{p_j'}{p_j}+2\theta d_j\frac{h_j'}{h_j})\Phi'](z_j)$, which is $Q_j$ for
\[
 \Phi'(z)=-\frac{e^{-\frac{z^2}{4d_j}}}{d_jp_j^2h_j^{2\theta}}\int_0^z
   e^{\frac{u^2}{4d_j}}p_j^2h_j^{2\theta}e^{-b_{\rm ac}u^2}\dd u .
\]
Since
$b_{\rm ac}\le\frac1{8d_j}$, $|\Phi'|+|\Phi''|\le C_M(1+|z|)e^{-b_{\rm ac}z^2}$; the
coefficient errors are at most $C_M\delta\tau^{-\frac12}e^{-c_bt}$ times these, as in the proof
of \eqref{e:qresp}, and $\Phi(z_j)+C_M+C_M\delta\int_r^t\tau^{-\frac12}e^{-c_bt'}\dd t'$ is a
bounded nonnegative supersolution. Next, if $P'=P+\mathsf w\pa_y$, $d\ge d_{\min}$, $f\ge0$,
$v=T_{P',r}f$ and $\Vert T_{P,r}\mathsf w^2\Vert_{L^\infty}\le\frac{d_{\min}}{16}$, then
$\Vert v\Vert_{L^\infty}\le2\Vert T_{P,r}f\Vert_{L^\infty}$: on finite domains,
$Pv=f-\mathsf wv_y$, $P(v^2)=2vPv-2dv_y^2$ and the Cauchy--Schwarz inequality for the
positive operator $T_{P,r}$ give, with $M_v=\Vert v\Vert_{L^\infty}$,
$J_v^2=\Vert T_{P,r}v_y^2\Vert_{L^\infty}$ and $F=\Vert T_{P,r}f\Vert_{L^\infty}$,
\[
 M_v\le F+\tfrac14\sqrt{d_{\min}}\,J_v,\qquad
 d_{\min}J_v^2\le M_v\bigl(F+\tfrac14\sqrt{d_{\min}}\,J_v\bigr),
\]
so $J_v\le\frac{4F}{3\sqrt{d_{\min}}}$ and $M_v\le2F$. Apply this $N$ times with
$P=P_{j,\frac iN}$ and $\mathsf w=-\frac1N(d_y+2d\frac{(h_j)_y}{h_j})$; by \eqref{e:dres},
$\Vert T_{P,r}\mathsf w^2\Vert_{L^\infty}\le\frac{C_M}{N^2}\le\frac{d_{\min}}{16}$ for
$N=N(M)$. By Lemma~\ref{l:sq}, \eqref{e:pmat} and \cref{l:kernv}(i),
$\Vert T_{\mathscr P_j,r}|\mathfrak r|^2\Vert_{L^\infty}\le C_M\omega_b(r)$; hence
\[
 \eta_r:=\max_j\bigl\Vert T_{P_{j,1},r}|\mathfrak r|^2\bigr\Vert_{L^\infty}+\tau_r^{-1}
 \le2^NC_M\omega_b(r)+\tau_r^{-1}\to0 .
\]
Put $\Upsilon_j=T_{P_{j,1},r}(\eta_r^{-1}|\mathfrak r|^2+Q_j+W_{-j})$, so that
$0\le\Upsilon_j\le\Upsilon_*=1+2C_M$.

Let $w$, $h>0$ satisfy $w_t=\pa_y(dw_y-bw)+f_w$ and $h_t=\pa_y(dh_y-bh)+\mathsf d$, let
$\mathsf z=\frac wh$, $P=\pa_t-d\pa_y^2+(b-d_y-2d\frac{h_y}h)\pa_y$, $P\Upsilon=V$ with
$\Upsilon\ge0$, $\phi=(1+\Upsilon)^{-\frac\chi2}$, $z_\Upsilon=\frac{\Upsilon_y}{1+\Upsilon}$, and
let $a>0$ be a weight and $0\le\chi\le1$ a constant. Then $P\mathsf z=\frac{f_w-\mathsf d\mathsf z}h$,
$P^*(ah)=-h[a_t+(da_y)_y+ba_y]-a\mathsf d$, and the computation of \eqref{e:rcw} with the
exponent $\frac\chi2$ in place of $\frac12$ gives
\begin{equation}\label{e:dden}
\begin{aligned}
 &\frac12\frac{\dd}{\dd t}\int_\R ah\phi\mathsf z^2\dd y
 +\int_\R adh\phi\Bigl[\mathsf z_y^2-\chi z_\Upsilon\mathsf z\mathsf z_y
   +\frac{\chi(\chi+2)}8z_\Upsilon^2\mathsf z^2\Bigr]\dd y
 +\frac\chi4\int_\R\frac{ah\phi V}{1+\Upsilon}\mathsf z^2\dd y\\
 &\quad=\frac12\int_\R h\phi\bigl[a_t+da_{yy}+(b+d_y)a_y\bigr]\mathsf z^2\dd y
 -\frac12\int_\R a\phi\,\mathsf d\,\mathsf z^2\dd y+\int_\R a\phi\mathsf zf_w\dd y ,
\end{aligned}
\end{equation}
and the bracket is at least $\frac14\mathsf z_y^2+\frac\chi{24}z_\Upsilon^2\mathsf z^2$. If
$\chi$ is a function with $0\le\chi\le1$, the right side acquires terms bounded by
$C(\Upsilon_*)(|\chi_t|+|\chi_y|+|\chi_{yy}|+\frac{\chi_y^2}\chi)\mathsf z^2
+C(\Upsilon_*)|\chi_y||\mathsf z||\mathsf z_y|+C(\Upsilon_*)|\chi_y||z_\Upsilon|\mathsf z^2$.

Let $\zeta,\chi$ be smooth functions of $\frac x\tau$ with $0\le\zeta,\chi\le1$, $\zeta=1$
for $|x|\le\frac{\theta_c}2\tau$, $\zeta=0$ for $|x|\ge\theta_c\tau$, $\chi=0$ for
$|x|\le\theta_c\tau$, $\chi=1$ for $|x|\ge2\theta_c\tau$, and
$\frac{\zeta_y^2}{\zeta(1-\zeta)}+\frac{\chi_y^2}\chi\le C\tau^{-2}$. Let $\psi_c$ be smooth,
bounded and nondecreasing, $\psi_c(0)=0$, $\psi_c'\ge1$ on $[-3\theta_c,3\theta_c]$ and
$\psi_c'=0$ outside $[-4\theta_c,4\theta_c]$, and for $\mu_c\ge1$ and $\epsilon_c>0$ put
\[
 w^c_j=e^{-\mu_c\psi_c(\frac{jx}\tau)},\qquad
 h_c=-\epsilon_c\int_0^{\frac x{\sqrt\tau}}e^{-b_0z^2}\dd z,\qquad
 e_j=2V_c\bigl(1+js_+(V_c)h_c\bigr),
\]
$F_j=(1+\Upsilon_j)^{-\frac\chi2}$, $q_j=\frac{\mathbf z_j}{h_j}$ and
\[
 E=\frac12\sum_j\int_\R e_jw^c_j\bigl[\zeta\mathbf z_j^2+(1-\zeta)h_jF_jq_j^2\bigr]\dd y .
\]
Choose $\epsilon_c$ so small that $|s_+(V_c)h_c|\le\frac12$; then
$C_M^{-1}\Vert Q_D\Vert_{L^2}^2\le E\le C_M\Vert Q_D\Vert_{L^2}^2$.
Apply \eqref{e:dden} to $\mathbf z_j$, with $(h,\chi,a)=(1,0,\zeta e_jw^c_j)$, hence
$\mathsf d=(\tilde s_j)_y$, and with $(h,\Upsilon,a)=(h_j,\Upsilon_j,(1-\zeta)e_jw^c_j)$,
$P=P_{j,1}$, $b=\tilde s_j$, $f_w=f_j$. The following bounds hold for late $r$.
\begin{enumerate}
\item[(a)] Since $\psi_c'$ vanishes for $|x|\ge4\theta_c\tau$ and $j\tilde s_j\ge\frac{c_-}2$,
$(\pa_t+\tilde s_j\pa_y)\log w^c_j=-\frac{\mu_c}\tau\psi_c'(\frac{jx}\tau)(j\tilde s_j-\frac{jx}\tau)
\le-\frac{\mu_cc_-}{4\tau}\psi_c'$. By $(h_c)_y=-\epsilon_cW_c$ and
$(h_c)_t=\frac{\epsilon_cx}{2\tau}W_c$, $(\pa_t+\tilde s_j\pa_y)\log e_j\le
-\frac{\epsilon_c}2s_+(V_c)^2W_c+C\delta W_c$ on $|x|\le\theta_c\tau$, and
$|(\tilde s_j)_y|\le C\delta W_c$ there. The remaining terms of
$\frac12[a_t+da_{yy}+(\tilde s_j+d_y)a_y]$, including those with $\zeta$ and $\chi$, are at
most $C_{\mu_c}[\tau^{-1}(1+\frac{|x|}{\sqrt\tau})^2e^{-\frac{b_0x^2}\tau}+\tau^{-2}+e^{-c_bt}]$ and
$\frac{C(\Upsilon_*)}\tau$ on $\theta_c\tau\le|x|\le2\theta_c\tau$; the $\zeta$-terms of the two
identities cancel up to $C_Me^{-c_bt}\tau^{-1}$, because $F_j=1$ and $|h_j-1|\le C_Me^{-c_bt}$
where $\zeta_y\ne0$ or $\zeta_t\ne0$. Splitting at $|x|=\tau^{\frac12+\frac1{24}}$,
$\tau^{-1}(1+\frac{|x|}{\sqrt\tau})^2e^{-\frac{b_0x^2}\tau}\le\varepsilon W_c+k_\varepsilon(t)$
with $\int_r^\infty k_\varepsilon<\infty$ for every $\varepsilon>0$ and late $r$.
\item[(b)] The potential term is at least
$c_M\chi(\eta_r^{-1}|\mathfrak r|^2+Q_j+W_{-j})q_j^2$, and $|\mathsf d_j|\le C_be^{-c_bt}$ on
$\operatorname{supp}(1-\zeta)$.
\item[(c)] On $\operatorname{supp}\zeta$, $|a_\pm|\le C_Me^{-c_bt}$; there the work
$\int\zeta e_jw^c_j\mathbf z_jf_j$ is at most
$\varepsilon\int\zeta|\pa_y\mathbf z|^2+C\delta\int W_c|\mathbf z|^2$ plus the integral of
$C_\varepsilon\zeta[\omega_b(r)\tau^{-1}+C_b\tau^{-\frac12}W_c+C_{\mu_c}\tau^{-\frac32}]|\mathbf z|^2$,
with $\omega_b(r)=\sup_{t\ge r}\tau\Vert R_U\Vert_{L^\infty}^2$ (Lemma~\ref{l:asy}) and
$\tau^{-1}G_0^2\le\tau^{-\frac12}W_c$. On $\operatorname{supp}(1-\zeta)$ the coefficients
$J_S$, $T_S$ are at most $C\delta e^{-c_bt}$, and in the work
$\int a^E_jF_jq_j\pa_y(\mathcal N\mathbf z)_j=-\int(a^E_jF_jq_j)_y(\mathcal N\mathbf z)_j$,
$a^E_j=(1-\zeta)e_jw^c_j$, $\mathbf z_\rho=h_\rho q_\rho$, the parts with $(a^E_j)_y$ and
$\chi_y$ cost $C_{\mu_c}\tau^{-\frac32}|q|^2$. The
residual entries cost $\varepsilon$ times the gradient and $\chi z_{\Upsilon_j}^2q_j^2$ terms plus
$C_\varepsilon|\mathfrak r|^2|q|^2$, paid by $\eta_r$ times the potential where $\chi=1$ and by
$\frac{\mu_c}\tau$ where $\chi<1$; the opposite profile $2ka_{-j}$ and the receiving profile
$2ka_j$ cost $CQ_{-j}q_j^2$ and $CQ_jq_{-j}^2$, paid by the potentials of the other family since
$Q_i=\tau^{-\frac12}W_i$. In the emitting entry,
\[
 -\int_\R(a^E_jF_jq_j)_y\,2\ell_{\rm g}a_{-j}h_{-j}q_{-j}\dd y
 =\int_\R a^E_jF_jq_j\bigl[2\ell_{\rm g}a_{-j}h_{-j}(q_{-j})_y
   +(2\ell_{\rm g}a_{-j}h_{-j})_yq_{-j}\bigr]\dd y,
\]
the first term costs $\varepsilon\Vert(1-\zeta)^{\frac12}(q_{-j})_y\Vert^2+C_\varepsilon\int
(1-\zeta)Q_{-j}q_j^2$ and the second $\varepsilon\int(1-\zeta)Q_{-j}q_{-j}^2+C_\varepsilon\int
(1-\zeta)Q_{-j}q_j^2$, since $|(\ell_{\rm g}a_{-j}h_{-j})_y|\le C_MQ_{-j}$.
\end{enumerate}
Choose $\epsilon_c$, then $\delta$ with $C\delta\le c_*\epsilon_c$, then $\mu_c$ large
depending on $\Upsilon_*$, then $\varepsilon$, then $r$. Moreover
$(Q_D)_y=S^{-1}(\mathbf z_y-J_S\mathbf z)$, $\mathbf z_j=h_jq_j$ and $|(h_j)_y|^2\le C_MQ_j$,
where $Q_j\le C_Me^{-c_bt}$ on $\{\chi<1\}$. Summing the identities therefore gives
\begin{equation}\label{e:dener}
 E'+c_b\Bigl(\Vert(Q_D)_y\Vert_{L^2}^2+\int_\R W_c|Q_D|^2\dd y\Bigr)\le k_b(t)E,
 \qquad k_b\ge0,\quad \int_r^\infty k_b\dd t<\infty .
\end{equation}

By \eqref{e:nash} for each component,
$\Vert Q_D\Vert_{L^2}^6\le4\Vert Q_D\Vert_{L^1}^4\Vert(Q_D)_y\Vert_{L^2}^2$,
so $\tilde E=e^{-\int_r^tk_b}E$ satisfies $\tilde E'\le-c_b'\Vert p\Vert_{L^1}^{-4}\tilde E^3$
and $\tilde E\le\min\{\tilde E(r),C_{b,r}\Vert p\Vert_{L^1}^2(t-r)^{-\frac12}\}$; this gives the
second bound in \eqref{e:dpb}. Multiplying \eqref{e:dener} by $\tau^{\frac38}$,
\[
 \int_r^\infty\tau^{\frac38}\Vert(Q_D)_y\Vert_{L^2}^2\dd t
 \le C_{b,r}L_p^2+C_b\int_r^\infty\bigl(\tau^{-\frac58}+\tau^{\frac38}k_b\bigr)E\dd t
 \le C_{b,r}L_p^2,
\]
and the integral of \eqref{e:dener} bounds $\mathscr Q_r(Q_D)$; this proves \eqref{e:dpb}.

For \eqref{e:dpn} let $T<\infty$; $\mathscr N^\xi_{r,T}(z)$ is finite. Write the flux in the
Duhamel formula for $z_j$ as the source $\tilde{\mathsf E}_\xi z+\tilde{\mathsf E}z_\xi$, and
use Minkowski's inequality over the birth time $s$ with \eqref{e:birth} and \eqref{e:dsmb}
($w=\frac{a_s}4$, so that $e^{\frac{w\xi^2}\tau}g_0\le e^{-\frac{a_\xi\xi^2}{4\tau}}$). By
\eqref{e:dfco}: the trace costs $C\tau_r^{\frac38}\Vert z(r)\Vert_{L^1}$; $\mathsf F^cz$
costs $(C_0C_M\delta+C_b\delta\tau_r^{-\frac14})\mathscr N^\xi_{r,T}(z)$ by \eqref{e:dsmb};
$\mathsf E^cz_\xi$ and $\mathsf E^c_\xi z$ cost $C_b\delta\tau_r^{-\frac18}\mathscr N^\xi_{r,T}(z)$
by the first line of \eqref{e:birth}, since
$\Vert g_0e^{\frac{a_\xi\xi^2}{4\tau}}\Vert_{L^2}\le C\tau^{\frac14}$; in $\mathsf F^rz$, the
part with $g_0$ costs $C_b\tau_r^{-\frac18}\mathsf O_\infty$, the part with $g_{-k}$ on
$|\eta|\ge v_0\tau_s$ costs $C_b\tau_r^{-\frac12}\mathsf O_\infty$ by the second line of
\eqref{e:birth}, and the rest $C_b\tau_r^{-\frac18}\mathsf X_\infty$, since $g_{-k}\le Ce^{-c_*t}$
on $|\eta|\le v_0\tau$; $(\mathsf E^r+J)z_\xi$ costs
$C_b\int_r^\infty\tau^{-\frac38}\Vert z_\xi\Vert_{L^2}\le C_{b,r}L_p$ by \eqref{e:dpb};
$\mathsf E^r_\xi z$ costs $C_{b,r}L_p\int_r^\infty\tau^{-\frac98}$; $J_\xi z$ costs
$C_b\mathsf X_\infty\int_r^\infty\tau^{-\frac98}$ on $|\eta|\le v_0\tau$ and
$C_{b,r}L_p\int_r^\infty\tau^{-\frac54}$ elsewhere. Hence
\[
 \mathscr N^\xi_{r,T}(z)\le C_{b,r}L_p
 +\bigl(C_0C_M\delta+C_b\tau_r^{-\frac18}\bigr)\mathscr N^\xi_{r,T}(z),
\]
where $C_0$ depends only on $M$, $a_s$ and the fixed data. Decrease $\delta_1(M)$ so that
$C_0C_M\delta\le\frac14$ and increase $r$; then $\mathscr N^\xi_{r,\infty}(z)\le C_{b,r}L_p$.
Since $\frac{(Q_D)_2}{V_c}=\sum_k\frac{\mathcal R_{2k}}{V_c}f_cz_k(F_c,t)$ and
$|\frac{\mathcal R_{2k}}{V_c}|+\tau^{\frac12}|\pa_y\frac{\mathcal R_{2k}}{V_c}|\le C$,
\eqref{e:xmass} gives \eqref{e:dpn}.
\end{proof}

\begin{lemma}\label{l:dk}
Let $r$, $p$ and $Q_D$ be as in Lemma~\ref{l:dp}, $u=\mathcal R^{-1}Q_D$, $g$ as in
\eqref{e:dsrc}, and let $H^D_j(t,s)$ be the evolution of $w_t=\pa_y\{d(V_c)w_y-b^D_jw\}$ from $s$
to $t$. Then the histories
\begin{equation}\label{e:dhist}
\begin{gathered}
 N_j(t)=H^D_j(t,r)|u_j(r)|+\int_r^tH^D_j(t,s)|g_j(s)|\dd s\ge|u_j(t)|,\qquad
 \sum_jm^N_j\le C_{b,r}L_p,\\
 \Vert N_j(t)\Vert_{L^1}=\Vert u_j(r)\Vert_{L^1}+\int_r^t\Vert g_j(s)\Vert_{L^1}\dd s
 \le m^N_j=\Vert u_j(r)\Vert_{L^1}+\int_r^\infty\Vert g_j(s)\Vert_{L^1}\dd s,
\end{gathered}
\end{equation}
and
\begin{equation}\label{e:dgap}
 j\bigl[b^D_j-s_{-j}(V_c)\bigr]\ge c_-,\qquad
 \int_r^\infty\tau^{\frac38}\Vert u_y\Vert_{L^2}^2\dd t\le C_{b,r}L_p^2 .
\end{equation}
\end{lemma}

\begin{proof}
By \eqref{e:xder}, $H^D_j(t,s)$ has the kernel $f_c(y,t)H_j(t,s;F_c(y,t),F_c(z,s))\ge0$ and
preserves the integral. Since $|\tilde{\mathsf E}|\le C$, $\int_r^\infty\Vert g\Vert_{L^1}$ is
at most $\int_r^\infty(\Vert\tilde{\mathsf F}z\Vert_{L^1}+\Vert\tilde{\mathsf E}_\xi z\Vert_{L^1}
+\Vert\tilde{\mathsf E}z_\xi\Vert_{L^1})\dd t$; by \eqref{e:dfco} as in the proof of
\eqref{e:dpn}, without the time factors in \eqref{e:birth}, this is at most
$C_{b,r}(\mathsf X_\infty+\mathsf O_\infty+\mathscr N^\xi_{r,\infty}(z)+L_p)\le C_{b,r}L_p$; for
instance $\int_r^\infty\Vert\mathsf E^cz_\xi\Vert_{L^1}\le C_M\delta\tau_r^{-\frac12}
\mathscr N^\xi_{r,\infty}(z)$. Duhamel's formula for \eqref{e:dsc} (which holds by \eqref{e:uniq})
and positivity give $|u_j|\le N_j$, and conservation of the integral gives the identity for
$\Vert N_j(t)\Vert_{L^1}$. By \eqref{e:frm}, $j[b^D_j-s_{-j}(V_c)]\ge2|s_j(V_c)|-C_M\tau^{-\frac12}
\ge c_-$ for late $r$. Finally $u_y=\mathcal R^{-1}(Q_D)_y-K^{\rm fr}_yu$ with
$|K^{\rm fr}_y|\le C_M\tau^{-\frac12}$, and \eqref{e:dpb} gives the last bound.
\end{proof}

\subsection{Exterior sources}

\begin{lemma}\label{l:ex}
Let $c_{\rm gap},\theta>0$. For every sufficiently late fixed $r$ the following holds, with
$C_b$ depending also on $c_{\rm gap}$ and $\theta$. For $j=\pm1$ let $b_j$ be a scalar function
on $\R\times[r,\infty)$ such that the drift of \eqref{e:nco} with $B=b_j$ satisfies the
hypotheses stated before \eqref{e:xder}, with $j[b_j-s_{-j}(V_c)]\ge c_{\rm gap}$, let
$H^b_j(t,s)$ be the evolution of $w_t=\pa_y\{d(V_c)w_y-b_jw\}$ from $s$ to $t$, and let
$0\le N_j(r)\in L^1(\R)$, $0\le n_j\in L^1((r,\infty)\times\R)$,
\begin{equation}\label{e:exhist}
 N_j(t)=H^b_j(t,r)N_j(r)+\int_r^tH^b_j(t,s)n_j(s)\dd s,\qquad
 m^N_j=\Vert N_j(r)\Vert_{L^1}+\int_r^\infty\Vert n_j(s)\Vert_{L^1}\dd s ,
\end{equation}
let $F^o,F^g\in L^2((r,T)\times\R;\R^2)$ for every finite $T$ vanish for $|y|<\theta\tau$,
with $F^o_{-j}=0$ if $m^N_j=0$, and put
\begin{equation}\label{e:exbud}
 B_F(T)=\sum_{j:\,m^N_j>0}m^N_j\int_r^T\!\!\int_\R\frac{|F^o_{-j}|^2}{N_j}\dd y\dd t
 +\Bigl(\int_r^T\Vert F^g\Vert_{L^2}\dd t\Bigr)^2 ,
\end{equation}
where the quotient is $0$ where $F^o_{-j}$ and $N_j$ vanish and $\infty$ where only $N_j$
vanishes. Then $\mathcal K\mathbf H=F^o+F^g$, $\mathbf H(r)=0$, has a unique solution on
$(r,\infty)$ with $\mathbf H\in C([r,T];H^1)$ and $\mathbf H_t,\mathbf H_{yy}\in
L^2((r,T)\times\R)$ for finite $T$, and
\[
 \mathscr O_{r,T}\bigl(\pa_y(S^{-1}\mathbf H)\bigr)\le C_bB_F(T)\qquad(r<T<\infty).
\]
\end{lemma}

\begin{proof}
By Lemma~\ref{l:am}, \eqref{e:pmat} and Lemma~\ref{l:asy}, $\widetilde{\mathcal M}$ satisfies
the hypotheses of Lemma~\ref{l:geo}, and $\sup_{t\ge r}\tau\Vert\widetilde{\mathcal M}
\Vert_{L^\infty}^2<\infty$. For finite $T$, Lemma~\ref{l:geo} with $F_1=F^o+F^g$ and $f=0$
gives the solution on $(r,T)$; by uniqueness these solutions agree on common intervals.

Let $m^N_j>0$. If $N_j(r),n_j\in C_c^\infty$, then $N_j$ is a classical solution of
$(N_j)_t=\pa_y\{d(V_c)(N_j)_y-b_jN_j\}+n_j$, and $N_j$, $(N_j)_y$ decay as
$|y|\to\infty$. Integrating this equation over $(y,\infty)$ for $j=+1$ and over
$(-\infty,y)$ for $j=-1$ shows that
\[
 \mathfrak m_+(y,t)=\int_y^\infty N_+(z,t)\dd z,\qquad
 \mathfrak m_-(y,t)=\int_{-\infty}^yN_-(z,t)\dd z,
\]
for which $(\mathfrak m_j)_y=-jN_j$, satisfy
\[
 L_{-j}\mathfrak m_j=j\bigl[b_j-s_{-j}(V_c)\bigr]N_j+\int_{I_j(y)}n_j(z,t)\dd z
 \ge c_{\rm gap}N_j,\qquad I_+(y)=(y,\infty),\quad I_-(y)=(-\infty,y).
\]
Since $0\le\mathfrak m_j\le m^N_j$, comparison on finite domains gives
\begin{equation}\label{e:exob}
 0\le T_{L_{-j},r}\frac{N_j}{m^N_j}\le\frac1{c_{\rm gap}} .
\end{equation}
For general data let $N^{(k)}_j$ be given by \eqref{e:exhist} with $C_c^\infty$ data
converging in $L^1$. Then $N^{(k)}_j\to N_j$ in $C([r,T];L^1)$, on every finite domain
$\Omega$ the functions $\min\{N^{(k)}_j,n\}$ converge to $\min\{N_j,n\}$ in $L^2(\Omega)$,
and \eqref{e:exob} passes to the limit by \eqref{e:Tdef}.

Put $\mathfrak v_{-j}=c_{\rm gap}\frac{N_j}{m^N_j}$ if $m^N_j>0$ and $\mathfrak v_{-j}=0$
otherwise. By \eqref{e:exob}, $\Vert T_{L_{-j},r}\mathfrak v_{-j}\Vert_{L^\infty}\le1$, and,
for $m^N_j>0$,
\[
 \int_r^T\!\!\int_\R\frac{|F^o_{-j}|^2}{W^{a_s}+\mathfrak v_{-j}}\dd y\dd t
 \le\frac{m^N_j}{c_{\rm gap}}\int_r^T\!\!\int_\R\frac{|F^o_{-j}|^2}{N_j}\dd y\dd t .
\]
Lemma~\ref{l:geo} with $a=a_s$, $g=0$, $F_1=F^g$ and $f=F^o$ therefore gives
$\int_r^T(\Vert\mathbf H_y\Vert_{L^2}^2+\int_\R W^{a_s}|\mathbf H|^2\dd y)\dd t\le C_bB_F(T)$.
By the bounds for the derivatives of $V_c$ stated before \eqref{e:pcnx},
$y^2\ge\frac{(y+h_*)^2}2-h_*^2$, $\tau\ge h_*^2$ and $a_s\le\frac{c_G}4$, the matrices
$J_S=S_yS^{-1}$ and $T_S=S_tS^{-1}$ of \eqref{e:exder} satisfy
\begin{equation}\label{e:exS}
 |S_y|+|(S^{-1})_y|+|J_S|+|\pa_yd(V_c)|+|\Lambda_y|+\tau^{\frac12}|T_S|
 \le C\delta W^{a_s} .
\end{equation}
Hence $\mathbf q=\pa_y(S^{-1}\mathbf H)=S^{-1}\mathbf H_y+(S^{-1})_y\mathbf H$ and
$\mathbf z=S\mathbf q$ satisfy
\begin{equation}\label{e:expr}
 \int_r^T\Vert\mathbf z\Vert_{L^2}^2\dd t\le C\int_r^T\Vert\mathbf q\Vert_{L^2}^2\dd t
 \le C_bB_F(T).
\end{equation}

Let $\zeta_0\in C_c^\infty((-\frac\theta2,\frac\theta2);[0,1])$ equal $1$ on
$[-\frac\theta4,\frac\theta4]$, and put $\zeta=\zeta_0(\frac y\tau)$, $\mathbf u=\zeta\mathbf z$,
$d=d(V_c)$ and $\mathcal T_\sigma w=w_t-\pa_y(dw_y-s_\sigma(V_c)w)$; row $\sigma$ of the
left side of the second line of \eqref{e:exder} is $\mathcal T_\sigma z_\sigma$. Since
$\zeta F=\zeta_yF=0$ for $F=F^o+F^g$, multiplication of that line by $\zeta$ gives
\[
\begin{gathered}
 \mathcal T_\sigma u_\sigma=\pa_y\bigl[(\widetilde{\mathcal M}-dJ_S)\mathbf u\bigr]_\sigma
 +\bigl[-dJ_S\mathbf u_y+(T_S+J_S\Lambda-J_S\widetilde{\mathcal M}+dJ_S^2)\mathbf u
 +\mathbf e\bigr]_\sigma,\\
 \mathbf e=\bigl(\zeta_t-d\zeta_{yy}+(\Lambda-d_yI)\zeta_y\bigr)\mathbf z
 -2d\zeta_y\mathbf z_y-\zeta_y(\widetilde{\mathcal M}-2dJ_S)\mathbf z .
\end{gathered}
\]
For each $\sigma$ apply the computation for \eqref{e:flxe}, with $a=a_s$ and
$f_{\rm flx}=[(\widetilde{\mathcal M}-dJ_S)\mathbf u]_\sigma$, to $u_\sigma$, add the work of the
bracket, and sum over $\sigma$; let $E$ be the sum of the energies $E_v$ of that computation.
With $|\zeta_t|+|\zeta_y|\le C_\theta\tau^{-1}$, $|\zeta_{yy}|\le C_\theta\tau^{-2}$,
$\zeta\mathbf z_y=\mathbf u_y-\zeta_y\mathbf z$, \eqref{e:exS}, $\tau|\widetilde{\mathcal M}|^2\le
C_b$ and $D_{\mathbf u}=\Vert\mathbf u_y\Vert_{L^2}^2+\int_\R W^{a_s}|\mathbf u|^2\dd y$, this
gives, for $\varepsilon>0$,
\[
 E'+c_*D_{\mathbf u}\le\bigl(C\delta+C_b\tau^{-\frac12}+\varepsilon\bigr)D_{\mathbf u}
 +C_{b,\theta,\varepsilon}\tau^{-1}\Vert\mathbf z\Vert_{L^2}^2 .
\]
Decrease $\delta_1(M)$, fix $\varepsilon$ and increase $r$. Multiplication by $\tau$,
$E\le C\Vert\mathbf z\Vert_{L^2}^2$, $E(r)=0$, $\tau W^{a_s}=\sqrt\tau G$ and \eqref{e:expr}
give $\int_r^T\sqrt\tau\int_\R\zeta^2G|\mathbf z|^2\dd y\dd t\le C_bB_F(T)$. Let $r$ be so late
that $|h_*|\le\frac{\theta\tau_r}8$. Where $\zeta<1$, $|y+h_*|\ge\frac{\theta\tau}8$, so
$\sqrt\tau G\le\sqrt\tau e^{-\frac{a_s\theta^2\tau}{128}}\le C_\theta$; with
$|\mathbf q|\le C|\mathbf z|$ and \eqref{e:expr} this proves the assertion. The energy
identities hold by the regularity of $\mathbf H$, as in Lemma~\ref{l:flx}.
\end{proof}

\subsection{The homogeneous field}

\begin{lemma}\label{l:wmx}
For every sufficiently late fixed $r$ and every $p\in L^1(\R;\R^2)\cap L^2(\R;\R^2)$ the
problem
\begin{equation}\label{e:wh}
 \mathfrak q_t=\pa_y\bigl\{d(V_c)\mathfrak q_y-B\mathfrak q\bigr\},\qquad \mathfrak q(r)=p,
\end{equation}
with $B$ of \eqref{e:pB}, has a unique solution with $\mathfrak q\in C([r,T];L^2)$ and
$\mathfrak q_y\in L^2((r,T)\times\R)$ for finite $T$, and there are functions $d_1,d_2$ on
$\R\times(r,\infty)$, independent of $T$, with
\begin{equation}\label{e:wmx}
 \frac{\mathfrak q_2}{V_c}=d_1+d_2,\qquad
 \mathscr N_{r,\infty}(d_1)\le C_{b,r}L_p,\qquad
 \mathscr O_{r,T}(d_2)\le C_{b,r}L_p^2\Bigl[1+\int_r^T\frac{\varepsilon_b(t)^2}\tau\dd t\Bigr]
 \quad(T>r),
\end{equation}
where $L_p=\Vert p\Vert_{L^1}+\Vert p\Vert_{L^2}$,
$\varepsilon_b(t)=\sqrt\tau\Vert R_U(t)\Vert_{L^\infty}$, and $C_{b,r}$ does not depend on $p$
and $T$.
\end{lemma}

\begin{proof}
Write $d=d(V_c)$, $\mathscr N=\mathscr N_{r,\infty}$, $\mathscr O=\mathscr O_{r,\infty}$,
$\mathcal T_\beta w=w_t-\pa_y(dw_y-\beta w)$ for a scalar or diagonal $\beta$, and
\[
 X(H)=\sup_{t\ge r}\Vert H(t)\Vert_{L^2}^2+\int_r^\infty\Bigl(\Vert H_y\Vert_{L^2}^2
 +\int_\R W_c|H|^2\dd y\Bigr)\dd t .
\]
Let $Q_D$, $u=\mathcal R^{-1}Q_D$, $g$, $H^D_j$, $N_j$ and $m^N_j$ be as in Lemmas~\ref{l:dp}
and \ref{l:dk}, with $Q_D(r)=p$. By \eqref{e:dpb}, \eqref{e:dhist}, \eqref{e:dgap} and
$|\mathcal R^{-1}|\le C$,
\[
\begin{gathered}
 \Vert u(t)\Vert_{L^2}\le C_{b,r}L_p\tau^{-\frac14},\qquad
 \int_r^\infty\Bigl(\int_\R W_c|u|^2\dd y+\tau^{\frac38}\Vert u_y\Vert_{L^2}^2\Bigr)\dd t
 \le C_{b,r}L_p^2,\\
 |u_j|\le N_j,\qquad
 \int_r^\infty\Vert g_j\Vert_{L^1}\dd t\le m^N_j\le C_{b,r}L_p .
\end{gathered}
\]

\emph{Step 1.} Put $O=\operatorname{offdiag}(\mathcal M-M_{\rm ref})$ with $\mathcal M$ of
\eqref{e:pj}, $C_{\sigma o}=O_{\sigma o}$, and $\chi_c=\chi_0(\frac{F_c}{v_0\tau})$ with $v_0$ of
\eqref{e:dfco} and $\chi_0\in C_c^\infty((-2,2);[0,1])$ equal to $1$ on $[-1,1]$. As for
\eqref{e:dDel}, $|O|\le C|R_U|\le C\varepsilon_b\tau^{-\frac12}$. By Lemma~\ref{l:geo}, the
problems
\[
 \mathcal K\mathcal W_c=\chi_cOSQ_D,\qquad \mathcal K\mathcal W_e=(1-\chi_c)OSQ_D,\qquad
 \mathcal W_c(r)=\mathcal W_e(r)=0
\]
have unique solutions in the class of Lemma~\ref{l:ex}. By \eqref{e:pB}, $B_D-B=S^{-1}OS$ and
$\mathcal K(S\mathbf v)=S[\mathbf v_t-d\mathbf v_{yy}+B\mathbf v_y]$; for
$\mathbf v=S^{-1}(\mathcal W_c+\mathcal W_e)$ this gives $\mathbf v_t-d\mathbf v_{yy}+B\mathbf v_y
=(B_D-B)Q_D$, and differentiation gives $(\mathbf v_y)_t=\pa_y\{d(\mathbf v_y)_y-B\mathbf v_y
+(B_D-B)Q_D\}$. Adding the equation of $Q_D$ shows that
\begin{equation}\label{e:wident}
 \mathfrak q=Q_D+\pa_y(S^{-1}\mathcal W_c)+\pa_y(S^{-1}\mathcal W_e)
\end{equation}
solves \eqref{e:wh} in the stated class; by \eqref{e:uniq} it is the only solution.

\emph{Step 2.} Since $SQ_D=\mathcal Pu$, $(1-\chi_c)OSQ_D=F^o+F^g$ with
$F^o_\sigma=(1-\chi_c)C_{\sigma o}\mathcal P_{oo}u_o$ and
$F^g_\sigma=(1-\chi_c)C_{\sigma o}\mathcal P_{o\sigma}u_\sigma$. By $|u_j|\le N_j$ and
$|\mathcal P-I|\le C_M\tau^{-\frac12}$,
\[
 \frac{|F^o_{-j}|^2}{N_j}\le C\frac{\varepsilon_b^2}\tau N_j,\qquad
 \Vert F^g\Vert_{L^2}\le C_{b,r}L_p\varepsilon_b\tau^{-\frac54},
\]
and $F^o_{-j}=0$ where $N_j=0$; with $\Vert N_j(t)\Vert_{L^1}\le m^N_j$, the budget
\eqref{e:exbud} of these sources satisfies
\begin{equation}\label{e:webud}
 B_F(T)\le C_{b,r}L_p^2\Bigl[1+\int_r^T\frac{\varepsilon_b^2}\tau\dd t\Bigr].
\end{equation}
Since $|F_c|\le\sqrt2|y|$, $F^o$ and $F^g$ vanish for $|y|<\frac{v_0}{\sqrt2}\tau$. The
histories \eqref{e:dhist} are those of \eqref{e:exhist} with $b_j=b^D_j$,
$N_j(r)=|u_j(r)|$ and $n_j=|g_j|$; the drift $\beta_j$ of $b^D_j$ satisfies the hypotheses
stated before \eqref{e:xder}, and $c_{\rm gap}=c_-$ by \eqref{e:dgap}. Lemma~\ref{l:ex} with
$\theta=\frac{v_0}{\sqrt2}$ and \eqref{e:webud} give
\[
 \mathscr O_{r,T}\bigl(\pa_y(S^{-1}\mathcal W_e)\bigr)
 \le C_{b,r}L_p^2\Bigl[1+\int_r^T\frac{\varepsilon_b^2}\tau\dd t\Bigr]\qquad(T>r).
\]

\emph{Step 3.} Put $\mathcal M^Z=\operatorname{diag}(s_\sigma(V_c)-b^D_\sigma)$,
$\mathcal L=L-\mathcal M^Z\pa_y$ and $\mathcal K_0=L-\widetilde{\mathcal M}\pa_y$, and write
$k=k(V_c)$. By \eqref{e:frm}, $\lambda_j$ is an eigenvalue of $\Lambda-E$ with
$E=M_{\rm ref}+T_{\rm ref}$, so $|\lambda_j-s_j(V_c)+E_{jj}|\le C|E|^2\le C_M\tau^{-1}$. By
\eqref{e:jsq} at $(V_c,V_h,\mathfrak a^a)$,
$E_{jj}=-2k\mathfrak a^a_j+2k\mathfrak a^a_{-j}+\frac{s_j(V_c)}{2V_c}(V_h-V_c)
+(T_{\rm ref})_{jj}+O(C_M\tau^{-1})$, and
$\mathfrak a^a_\sigma-a_\sigma=\frac{H_h}2+\frac12(\frac{s_\sigma(V_c)}{s_\sigma}-1)a_\sigma
+\frac12(1+\frac{s_\sigma(V_c)}{s_o})a_o$, where the brackets vanish at $V_c=v_\sigma$,
respectively $V_c=v_o$; since $|V_c-v_\sigma|\le C\delta e^{-\frac{c_Gy^2}\tau}$ for
$\sigma y\ge0$ and $|y-s_\sigma t|\ge|s_\sigma|t$ for $\sigma y<0$, the last two terms are at
most $C_Me^{-c_bt}$. With $|H_h|+|V_h-V_c|+|T_{\rm ref}|\le C_M\tau^{-\frac12}G_0$
(\eqref{eq:pr-gau}, \eqref{e:pT}, $\delta|h_*|\le C_M$, $\tau\ge h_*^2$) and
$|\mathcal Z|\le C|R_U|$, \eqref{e:drf} gives
\begin{equation}\label{e:wbD}
 b^D_\sigma=s_\sigma(V_c)+2ka_\sigma-2ka_o+r_\sigma,\qquad
 |r_\sigma|\le C|R_U|+C_M\bigl(\tau^{-1}+\tau^{-\frac12}G_0\bigr).
\end{equation}
Thus $\mathcal M^Z$ has the diagonal form \eqref{e:gmat} with residual $-r_\sigma$, and by
Lemmas~\ref{l:asy}, \ref{l:sq} and \eqref{e:gauss}, $\sup\tau|r_\sigma|^2<\infty$ and
$\Vert T_{L_\sigma,r}|r_\sigma|^2\Vert_{L^\infty}\to0$. Hence, for
$\mathscr K\in\{\mathcal K,\mathcal K_0,\mathcal L,L\}$, potentials $\mathfrak v_\sigma\ge0$ with
$\Vert T_{L_\sigma,r}\mathfrak v_\sigma\Vert_{L^\infty}\le1$ and $\mathscr KH=F_1+f$, $H(r)=0$,
Lemma~\ref{l:geo} with $a=2b_0$ gives
\begin{equation}\label{e:wE1}
 X(H)\le C_b\Bigl[\Bigl(\int_r^\infty\Vert F_1\Vert_{L^2}\dd t\Bigr)^2
 +\sum_\sigma\int_r^\infty\!\!\int_\R\frac{|f_\sigma|^2}{W_c+\mathfrak v_\sigma}\dd y\dd t\Bigr].
\end{equation}
For $\mathscr K\in\{\mathcal K_0,\mathcal L,L\}$, row $\sigma$ of $H_y$ solves
$v_t=\pa_y(dv_y-s_\sigma(V_c)v+f_{\rm flx})$ with $f_{\rm flx}=(\mathcal M^\star H_y+F_1+f)_\sigma$,
$\mathcal M^\star=\widetilde{\mathcal M},\mathcal M^Z,0$, and $\tau|\mathcal M^\star|^2\le C_b$;
Lemma~\ref{l:flx} with $s=r$ gives
\begin{equation}\label{e:wE2}
 \mathscr O(H_y)\le C_b\Bigl[X(H)+\int_r^\infty\tau\Vert F_1+f\Vert_{L^2}^2\dd t\Bigr].
\end{equation}
If $f=Kw_y$ with $|K|\le C_b(|a_+|+|a_-|+|R_U|+\tau^{-1}+\tau^{-\frac12}G_0^{\frac12})$ and
$\mathfrak v_\sigma$ is the potential of Step~0 of the proof of Lemma~\ref{l:solv}, then
$|K_{\sigma\cdot}|^2\le C_b(W_c+\mathfrak v_\sigma)$ and $\tau|K|^2\le C_b$ (loc.\ cit.), so
\eqref{e:wE1} and \eqref{e:wE2} give $X(H)+\mathscr O(H_y)\le C_b\int_r^\infty\Vert w_y
\Vert_{L^2}^2\dd t$. Finally, if $w_t=\pa_y(dw_y-\beta w+f)$, $w(r)=0$, with
$\beta=\operatorname{diag}(b^D_\sigma)$ or $\beta=\Lambda$, then $w=H_y$ for $\mathcal LH=f$,
respectively $LH=f$, $H(r)=0$, by \eqref{e:uniq}, and
\begin{equation}\label{e:wE3}
 \mathscr O(w)\le C_b\Bigl[\Bigl(\int_r^\infty\Vert f\Vert_{L^2}\dd t\Bigr)^2
 +\int_r^\infty\tau\Vert f\Vert_{L^2}^2\dd t\Bigr].
\end{equation}

\emph{Step 4.} Let $H^0_j(t,s)$ be the evolution of $w_t=\pa_y\{dw_y-s_j(V_c)w\}$. By
\eqref{e:nco} its drift in $\xi=F_c$ is $jc_-+(F_c)_t-d(f_c)_y$, which satisfies \eqref{e:kH}
and the hypotheses of \cref{l:kernv}\textup{(ii)} with $\bar b=jc_-$, $A_0=A_1=A_2=C\delta$
and $v_0$. For a function $f$ on $\R$ at time $s$ define $f^\xi$ by
$f^\xi(F_c(y,s))=\frac{f(y)}{f_c(y,s)}$; for a matrix function $A$ and a two-vector $v$ put
$(T_Av)_j(t)=\int_r^tH^0_j(t,s)(Av)_j(s)\dd s$, and let
$G_*=e^{-\frac{c_G(y+h_*)^2}{4\tau}}$. Fix $w=\min\{C_*^{-1},\frac{3c_G}{64}\}$ with $C_*$ of
\cref{l:kernv}\textup{(iii)} for these drifts. Then, for $r<T\le\infty$,
\begin{equation}\label{e:wbody}
 \mathscr O_{r,T}\Bigl(\int_r^tH^0_j(t,s)f(s)\dd s\Bigr)
 \le C\int_r^T\tau_s^{\frac32}\bigl\Vert e^{\frac{w\eta^2}{\tau_s}}f^\xi(s)\bigr\Vert_{L^2}^2\dd s,
\end{equation}
and, if $|A|+\tau^{\frac12}|A_y|\le K\delta\tau^{-\frac12}G_*$,
\begin{equation}\label{e:wcont}
 \mathscr N(T_Av)\le C_3\delta\,\mathscr N(v),\qquad
 \mathscr O_{r,T}(T_Av)\le C_3^2\delta^2\mathscr O_{r,T}(v),
\end{equation}
with $C$ depending only on the fixed data and the fixed Gaussian exponents,
and $C_3$ also on $K$. Indeed, let
$s<t$, $h=t-s$ and $v=H^0_j(t,s)f$. \Cref{l:kernv}\textup{(iii)} with
$(w_1,w_2)=(w,\frac w2)$, for which $C_*w_1w_2\le w_1-w_2$, gives
$\Vert e^{\frac{w(\xi-jc_-h)^2}{2\tau_t}}v^\xi\Vert_{L^2}\le C\Vert e^{\frac{w\eta^2}{\tau_s}}
f^\xi\Vert_{L^2}$. As for \eqref{e:xmass}, $G\le Ce^{-\frac{a_\xi F_c^2}{2\tau}}$, so
$A_v(t)\le C\sqrt{\tau_t}\Vert e^{-\frac{a_\xi\xi^2}{4\tau_t}}v^\xi\Vert_{L^2}^2$, and
$\frac{a_\xi\xi^2}4+\frac w2(\xi-jc_-h)^2\ge c_wh^2$ with $c_w=\frac{a_\xi wc_-^2}{2a_\xi+4w}$.
Hence $\sqrt{A_v(t)}\le Ck(t,s)\tau_s^{\frac34}\Vert e^{\frac{w\eta^2}{\tau_s}}f^\xi
\Vert_{L^2}$ with $k(t,s)=\tau_t^{\frac14}\tau_s^{-\frac34}e^{-\frac{c_wh^2}{\tau_t}}$, and
$\sup_t\int_r^tk\dd s+\sup_s\int_s^\infty k\dd t\le C$ (split at $h=\tau_s$); Minkowski's
inequality and Schur's test give \eqref{e:wbody}. For $f=(Av)_j$,
$\pa_\eta f^\xi=\frac1{f_c}\pa_y\frac{Av}{f_c}$, $|(f_c)_y|\le C\delta\tau^{-\frac12}$, and
$\frac{wF_c^2}\tau-\frac{c_G(y+h_*)^2}{4\tau}\le4w-\frac{a_s(y+h_*)^2}{4\tau}$ since
$F_c^2\le4(y+h_*)^2+4h_*^2$ and $4w\le\frac{c_G-a_s}4$; hence
\[
\begin{gathered}
 \tau_s^{\frac14}\bigl\Vert e^{\frac{w\eta^2}{\tau_s}}f^\xi\bigr\Vert_{L^2}
 +\tau_s^{\frac34}\bigl\Vert e^{\frac{w\eta^2}{\tau_s}}\pa_\eta f^\xi\bigr\Vert_{L^2}
 \le CK\delta\bigl[\tau_s^{-\frac14}\Vert G^{\frac12}v\Vert_{L^2}
 +\tau_s^{\frac14}\Vert G^{\frac12}v_y\Vert_{L^2}\bigr],\\
 \tau_s^{\frac32}\bigl\Vert e^{\frac{w\eta^2}{\tau_s}}f^\xi\bigr\Vert_{L^2}^2
 \le CK^2\delta^2A_v(s).
\end{gathered}
\]
\eqref{e:dsmb} for $H^0_j$, Minkowski's inequality in $s$ and \eqref{e:xmass} give the first
bound in \eqref{e:wcont}, and \eqref{e:wbody} gives the second.

\emph{Step 5.} Since $SQ_D=\mathcal Pu$, the rows of the central source are
\begin{equation}\label{e:wcrow}
\begin{gathered}
 (\chi_cOSQ_D)_\sigma=\bigl(s_\sigma(V_c)-s_o(V_c)\bigr)\Gamma_\sigma u_o+F_{\times,\sigma},
 \qquad F_{\times,\sigma}=\chi_cC_{\sigma o}\mathcal P_{o\sigma}u_\sigma,\\
 \Gamma_\sigma=\mathfrak p_\sigma\chi_cC_{\sigma o},\qquad
 \mathfrak p_\sigma=\frac{\mathcal P_{oo}}{s_\sigma(V_c)-s_o(V_c)} .
\end{gathered}
\end{equation}
They vanish outside $\{|F_c|\le2v_0\tau\}$, which lies in $\{|y|\le\frac{c_-}{1024}t\}$ since
$|y|\le\sqrt2|F_c|$. By $|\mathcal P-I|\le C_M\tau^{-\frac12}$, $|F_\times|\le C_b\tau^{-1}|u|$,
so $\int_r^\infty\Vert F_\times\Vert_{L^2}\dd t\le C_{b,r}L_p$ and
$\int_r^\infty\tau\Vert F_\times\Vert_{L^2}^2\dd t\le C_{b,r}L_p^2$.

By \eqref{e:nco}, $\pa_tN^\xi_j=\pa_\xi(D_0\pa_\xi N^\xi_j-\beta_jN^\xi_j)+|g_j|^\xi$, and
$\int_\R N^\xi_j\dd\xi=\Vert N_j\Vert_{L^1}\le m^N_j$. With $x=j\xi$, $\epsilon=\tau^{-\frac12}$,
$f_\epsilon(z)=(\epsilon+\sqrt{\epsilon^2+z^2})^{-\frac34}$ and
$\chi_1\in C_c^\infty((-3v_0,3v_0);[0,1])$ equal to $1$ on $[-2v_0,2v_0]$, put
\[
 \Phi(x,t)=\int_{-\infty}^{x/\tau}f_\epsilon\chi_1\dd z,\qquad
 W_{\rm occ}=\frac1\tau f_\epsilon\Bigl(\frac x\tau\Bigr)\chi_1\Bigl(\frac x\tau\Bigr).
\]
Then $0\le\Phi\le C$, $\Phi_x=W_{\rm occ}$, $\Phi_t\ge-\frac x\tau W_{\rm occ}$ (since
$f_\epsilon$ increases in $t$), and $|\Phi_{xx}|\le C\tau^{-\frac12}W_{\rm occ}+C\tau^{-2}$
(since $|f_\epsilon'|\le\frac3{4\epsilon}f_\epsilon$). As $j\beta_j\ge\frac{c_-}2$ by
\eqref{e:dfco}, $\frac x\tau\le3v_0\le\frac{c_-}8$ on $\operatorname{supp}W_{\rm occ}$ and the
birth term is nonnegative,
\[
 \frac{\dd}{\dd t}\int_\R\Phi N^\xi_j\dd\xi
 \ge\int_\R\bigl(\Phi_t+j\beta_j\Phi_x-D_0|\Phi_{xx}|\bigr)N^\xi_j\dd\xi
 \ge\frac{c_-}4\int_\R W_{\rm occ}N^\xi_j\dd\xi-C\tau^{-2}m^N_j
\]
for late $r$. Integration in $t$ and $W_{\rm occ}\ge\tau^{-\frac58}(2+\frac{|\xi|}{\sqrt\tau})^{-\frac34}$
for $|\xi|\le2v_0\tau$ give
\begin{equation}\label{e:wocc}
 \int_r^\infty\!\!\int_{|\xi|\le2v_0\tau}\tau^{-\frac58}\Bigl(1+\frac{|\xi|}{\sqrt\tau}
 \Bigr)^{-\frac34}N^\xi_j\dd\xi\dd t\le Cm^N_j,\qquad
 \int_r^\infty\!\!\int_{|F_c|\le2v_0\tau}\tau^{-1}N_j\dd y\dd t\le Cm^N_j .
\end{equation}

Let $\mathcal U_a=\widehat{\mathcal U}(y+h_*,t)+\sum_j\vartheta_jr_j$, so that
$\mathcal U-\mathcal U_a=R_U$, and $R_a=(\mathcal U_a)_t+\mathsf A(\mathcal U_a)(\mathcal U_a)_y
-\pa_y(\mathsf B(\mathcal U_a)(\mathcal U_a)_y)$ with $\mathsf A,\mathsf B$ of \eqref{e:sfun}.
By \eqref{e:cmp} and \eqref{e:cdft}, $R_a=(f_v,f_b,f_T+b_af_b)$, so
$\Vert R_a\Vert_{L^\infty}^2\le\Vert R_a\Vert_{L^2}\Vert(R_a)_y\Vert_{L^2}\le C_b\tau^{-3}$. By
\eqref{e:nfflx},
\[
\begin{gathered}
 (R_U)_t=\pa_yJ_U+H_U,\qquad
 J_U=\mathsf B(\mathcal U)(R_U)_y+[\mathsf B(\mathcal U)-\mathsf B(\mathcal U_a)](\mathcal U_a)_y,\\
 H_U=-\mathsf A(\mathcal U)(R_U)_y-[\mathsf A(\mathcal U)-\mathsf A(\mathcal U_a)](\mathcal U_a)_y
 -R_a,
\end{gathered}
\]
and on $|y|\le\frac{c_-}{1024}t$ Lemma~\ref{l:ga}, $|R_U|\le C_b\tau^{-\frac12}$ and
$|\pa_y\mathcal U_a|+\tau^{\frac12}|\pa_t\mathcal U_a|\le C_M\tau^{-\frac12}$ give
$|J_U|+|H_U|\le C_b\tau^{-1}$. By \eqref{e:pj}, \eqref{e:pB} and the proof of
Lemma~\ref{l:am}, $C_{\sigma o}=\mathscr C(V_c,\mathcal U_a,R_U)$ with $\mathscr C$ smooth and
$\mathscr C(n,U,0)=0$; hence, with $D_o=\pa_t+s_o(V_c)\pa_y$,
\[
 D_oC_{\sigma o}=\pa_yJ_C+H_C,\quad J_C=\mathscr C_RJ_U,\quad
 H_C=\mathscr C_RH_U-(\mathscr C_R)_yJ_U+\mathscr C_n(V_c)_t+\mathscr C_U(\mathcal U_a)_t
 +s_o(V_c)(C_{\sigma o})_y,
\]
and $|J_C|+|H_C|+|(C_{\sigma o})_y|\le C_b\tau^{-1}$ on $|y|\le\frac{c_-}{1024}t$, since
$|\mathscr C_n|+|\mathscr C_U|\le C|R_U|$. With $|(\chi_c)_y|+|D_o\chi_c|\le C\tau^{-1}$ and
$|(\mathfrak p_\sigma)_y|+|D_o\mathfrak p_\sigma|\le C_M\tau^{-\frac12}$ this gives
\begin{equation}\label{e:wcurr}
\begin{gathered}
 D_o\Gamma_\sigma=\pa_yJ_\Gamma+H_\Gamma,\qquad J_\Gamma=\mathfrak p_\sigma\chi_cJ_C,\\
 H_\Gamma=\mathfrak p_\sigma\bigl[\chi_cH_C-(\chi_c)_yJ_C+(D_o\chi_c)C_{\sigma o}\bigr]
 +(D_o\mathfrak p_\sigma)\chi_cC_{\sigma o}-(\mathfrak p_\sigma)_y\chi_cJ_C,\\
 |\Gamma_\sigma|\le C_b\tau^{-\frac12},\qquad
 |(\Gamma_\sigma)_y|+|J_\Gamma|+|H_\Gamma|\le C_b\tau^{-1}.
\end{gathered}
\end{equation}

For $h_\sigma=\Gamma_\sigma u_o$ and $m_o=s_o(V_c)-b^D_o$, the identity
$\mathcal T_s(\Gamma u)=\Gamma\mathcal T_su+uD_o\Gamma-\pa_y(d\Gamma_yu)-d\Gamma_yu_y$ for
$s=s_o(V_c)$, $\mathcal T_su_o=\pa_y(m_ou_o)+g_o$ (by \eqref{e:dsc}) and \eqref{e:wcurr} give
\begin{equation}\label{e:wnf}
\begin{gathered}
 \mathcal T_{s_o(V_c)}h_\sigma=\pa_yF_\Gamma+B_y+B_v+\Gamma_\sigma g_o,\qquad
 F_\Gamma=\bigl(J_\Gamma-d(\Gamma_\sigma)_y+m_o\Gamma_\sigma\bigr)u_o,\\
 B_y=-\bigl(J_\Gamma+d(\Gamma_\sigma)_y\bigr)(u_o)_y,\qquad
 B_v=\bigl(H_\Gamma-m_o(\Gamma_\sigma)_y\bigr)u_o .
\end{gathered}
\end{equation}
Let $\mathcal LZ^0=((s_\sigma(V_c)-s_o(V_c))h_\sigma)_\sigma$ and $\mathcal LZ^\times=F_\times$,
with zero traces, and $Z=Z^0+Z^\times$; by \eqref{e:wcrow}, $\mathcal LZ=\chi_cOSQ_D$. By
\eqref{e:cnr} with $z=Z^0_\sigma$, $D=d$, $B=b^D_\sigma$, $e=s_o(V_c)$, $q=h_\sigma$ and
$f_2=-(b^D_\sigma-s_\sigma(V_c))h_\sigma$, by \eqref{e:wnf}, and by Duhamel's formula
(\eqref{e:uniq}),
\begin{equation}\label{e:wscal}
\begin{aligned}
 (Z^0_\sigma)_y={}&h_\sigma-H^D_\sigma(t,r)h_\sigma(r)
 -\int_r^tH^D_\sigma(t,s)\pa_y\bigl[F_\Gamma+(b^D_\sigma-s_\sigma(V_c))h_\sigma\bigr](s)\dd s\\
 &-\int_r^tH^D_\sigma(t,s)\bigl(B_y+B_v+\Gamma_\sigma g_o\bigr)(s)\dd s .
\end{aligned}
\end{equation}
By \eqref{e:wcurr} and $|m_o|+|b^D_\sigma-s_\sigma(V_c)|\le C_M\tau^{-\frac12}$, the fields
$h_\sigma$, $F_\Gamma+(b^D_\sigma-s_\sigma(V_c))h_\sigma$, $B_v$ and $B_y$ are at most
$C_b\tau^{-\frac12}|u|$, $C_b\tau^{-1}|u|$, $C_b\tau^{-1}N_o$ and $C_b\tau^{-1}|u_y|$, and vanish
outside $\{|F_c|\le2v_0\tau\}$. Hence $\mathscr O(h_\sigma)\le C_b\iint W_c|u|^2$; the flux term
has $\mathscr O\le C_{b,r}L_p^2$ by \eqref{e:wE3}; by the third line of \eqref{e:birth} and by
\eqref{e:xmass}, $\mathscr N(H^D_\sigma(\cdot,r)h_\sigma(r))\le C_b\tau_r^{\frac58}
\Vert h_\sigma(r)\Vert_{L^2}\le C_{b,r}L_p$, and the $B_y$ term has
$\mathscr N\le C_b\int_r^\infty\tau^{\frac58}\Vert B_y\Vert_{L^2}\le C_b\int_r^\infty
\tau^{-\frac38}\Vert u_y\Vert_{L^2}\le C_{b,r}L_p$; by \eqref{e:bcol} and \eqref{e:wocc} the
$B_v$ term has $\mathscr N\le C_bm^N_o$; by the first line of \eqref{e:birth} the last term
has $\mathscr N\le C_b\int_r^\infty\tau^{\frac38}\Vert\Gamma_\sigma g_o\Vert_{L^1}\le C_bm^N_o$;
and \eqref{e:wE1}, \eqref{e:wE2} give $X(Z^\times)+\mathscr O(Z^\times_y)\le C_{b,r}L_p^2$. Thus
$Z_y=z_1+z_2$, where $z_1$ is the sum of the terms of \eqref{e:wscal} with $h_\sigma(r)$, $B_y+B_v$ and
$\Gamma_\sigma g_o$, and
\[
 \mathscr N(z_1)+\sqrt{\mathscr O(z_2)}\le C_{b,r}L_p,\qquad
 \mathscr O_{r,T}(z_1)<\infty\quad(T<\infty),
\]
the last since $h_\sigma(r)$ and the sources of $z_1$ are square integrable on
$(r,T)\times\R$.

\emph{Step 6.} Let $\mathfrak v_\sigma=c_-\frac{N_o}{m^N_o}$ if $m^N_o>0$ and $\mathfrak v_\sigma=0$
otherwise. By Step~2, \eqref{e:exob} gives $\Vert T_{L_\sigma,r}\mathfrak v_\sigma\Vert_{L^\infty}
\le1$, and $|u_o|\le N_o$ and \eqref{e:wocc} give
\[
 \int_r^\infty\!\!\int_\R\frac{|(s_\sigma(V_c)-s_o(V_c))h_\sigma|^2}{\mathfrak v_\sigma}
 \dd y\dd t\le C_bm^N_o\int_r^\infty\!\!\int_{|F_c|\le2v_0\tau}\tau^{-1}N_o\dd y\dd t
 \le C_b(m^N_o)^2 .
\]
\eqref{e:wE1} with $F_1=F_\times$ and $f=((s_\sigma(V_c)-s_o(V_c))h_\sigma)_\sigma$, for
$\mathscr K=\mathcal K$ and $\mathscr K=\mathcal L$, gives
\begin{equation}\label{e:wceng}
 X(\mathcal W_c)+X(Z)\le C_{b,r}L_p^2 .
\end{equation}
Let $LZ_0=\mathcal C_{\widetilde{\mathcal M}}\mathcal W_c$ and
$\mathcal K_0Z_1=(\widetilde{\mathcal M}-\mathcal M^Z)Z_y+\widetilde{\mathcal M}(Z_0)_y$, with zero
traces. By \eqref{e:pcnx},
$|\mathcal C^0|+|B^1|+\tau^{\frac12}|C^1|\le C_b\tau^{-\frac12}G_*$, and
$\frac{\tau^{-1}G_*^2}{W_c}\le W_c\le1$; so \eqref{e:wE1} gives $X(Z_0)\le C_bX(\mathcal W_c)$.
By \eqref{e:gmat}, \eqref{e:pmat} and \eqref{e:wbD}, both matrices in the equation of $Z_1$
have the bound of $K$ in Step~3, so $X(Z_1)+\mathscr O((Z_1)_y)\le C_b\int_r^\infty
(\Vert Z_y\Vert_{L^2}^2+\Vert(Z_0)_y\Vert_{L^2}^2)\dd t\le C_{b,r}L_p^2$. Since
$\mathcal K_0Z=\mathcal LZ-(\widetilde{\mathcal M}-\mathcal M^Z)Z_y$,
$\mathcal K_0Z_0=\mathcal C_{\widetilde{\mathcal M}}\mathcal W_c-\widetilde{\mathcal M}(Z_0)_y$ and
$\mathcal K_0\mathcal W_c=\chi_cOSQ_D+\mathcal C_{\widetilde{\mathcal M}}\mathcal W_c$, uniqueness
for $\mathcal K_0$ in Lemma~\ref{l:geo} gives $\mathcal W_c=Z+Z_0+Z_1$.

Split $Z_0=V_0+Z_g$ by $LV_0=\mathcal C^0\mathcal W_c$ and $LZ_g=B^1(\mathcal W_c)_y+C^1\mathcal W_c$, with
zero traces. The last source $f_g$ has
$\iint\frac{|f_g|^2}{W_c}+\int_r^\infty\tau\Vert f_g\Vert_{L^2}^2\le C_bX(\mathcal W_c)$, so
$\mathscr O((Z_g)_y)\le C_{b,r}L_p^2$ by \eqref{e:wE1}, \eqref{e:wE2}. Put
$\mathsf a=\Lambda^{-1}\mathcal C^0$ and $\varphi=(V_0)_y-\mathsf a\mathcal W_c$. By
\eqref{e:cnr} with $z=V_0$, $D=d$, $B=\Lambda$, $e=0$, $q=\mathsf a\mathcal W_c$ and $f_2=0$, by
$(\mathcal W_c)_t=d(\mathcal W_c)_{yy}-\Lambda (\mathcal W_c)_y+\widetilde{\mathcal M}(\mathcal W_c)_y+\mathcal C_{\widetilde{\mathcal M}}
\mathcal W_c+\chi_cOSQ_D$ (\eqref{e:gop}, \eqref{e:pB}), and by
$\pa_y(d(\mathsf a\mathcal W_c)_y)-\mathsf ad(\mathcal W_c)_{yy}=\pa_y(d\mathsf a_y\mathcal W_c)
+(d\mathsf a)_y(\mathcal W_c)_y$,
\begin{equation}\label{e:wcaus}
\begin{aligned}
 \mathcal T_\Lambda\varphi={}&\mathsf a\Lambda (\mathcal W_c)_y+\bigl[(d\mathsf a)_y
 -\mathsf a\widetilde{\mathcal M}\bigr](\mathcal W_c)_y+\pa_y(d\mathsf a_y\mathcal W_c)\\
 &-\mathsf a\,\mathcal C_{\widetilde{\mathcal M}}\mathcal W_c-\mathsf a_t\mathcal W_c
 -\mathsf a\chi_cOSQ_D,\qquad \varphi(r)=0 .
\end{aligned}
\end{equation}
By \eqref{e:pcnx}, $|\mathsf a|+\tau^{\frac12}|\mathsf a_y|+\tau|\mathsf a_t|\le C\delta\tau^{-\frac12}G_*$. The nondivergence terms other than $\mathsf a\Lambda (\mathcal W_c)_y$ are
therefore at most $C_b\tau^{-1}G_*(|(\mathcal W_c)_y|+|\mathcal W_c|+|u|)$, and by \eqref{e:wbody}, with
$e^{\frac{wF_c^2}\tau}G_*\le e^{4w}G^{\frac12}$ (Step~4), the $\mathscr O$ of their response is
at most $C_b[\int_r^\infty\Vert (\mathcal W_c)_y\Vert_{L^2}^2\dd t+\iint W_c(|\mathcal W_c|^2+|u|^2)]\le
C_{b,r}L_p^2$. For the flux $d\mathsf a_y\mathcal W_c$,
$\int_r^\infty\Vert d\mathsf a_y\mathcal W_c\Vert_{L^2}\dd t\le
C\delta(\int_r^\infty\tau^{-\frac32}\dd t)^{\frac12}(\iint W_c|\mathcal W_c|^2)^{\frac12}$ and
$\int_r^\infty\tau\Vert d\mathsf a_y\mathcal W_c\Vert_{L^2}^2\dd t\le C\delta^2\iint W_c|\mathcal W_c|^2$, so
\eqref{e:wE3} bounds its response; moreover $\mathscr O(\mathsf a\mathcal W_c)\le C\delta^2
X(\mathcal W_c)$. With $(\mathcal W_c)_y=Z_y+(V_0)_y+(Z_g)_y+(Z_1)_y$ and Duhamel's formula for
\eqref{e:wcaus},
\begin{equation}\label{e:wser}
 (\mathcal W_c)_y=Z_y+T_{\mathsf a\Lambda}(\mathcal W_c)_y+q_g,\qquad \mathscr O(q_g)\le C_{b,r}L_p^2 .
\end{equation}
The coefficient $\mathsf a\Lambda=\Lambda^{-1}\mathcal C^0\Lambda$ satisfies the hypothesis of
\eqref{e:wcont} with $K$ depending only on the fixed data; decrease $\delta_1(M)$ so that
$C_3\delta\le\frac12$. Then $Q_1=\sum_{n\ge0}T_{\mathsf a\Lambda}^nz_1$ and
$Q_2=\sum_{n\ge0}T_{\mathsf a\Lambda}^n(z_2+q_g)$ converge in $\mathscr N$ and in
$\sqrt{\mathscr O}$, respectively, with $\mathscr N(Q_1)+\sqrt{\mathscr O(Q_2)}\le C_{b,r}L_p$,
and on every finite strip both converge in $\sqrt{\mathscr O_{r,T}}$, since
$\mathscr O_{r,T}(z_1)<\infty$. Hence $Q_1+Q_2=Z_y+q_g+T_{\mathsf a\Lambda}(Q_1+Q_2)$, and by \eqref{e:wser}
$(\mathcal W_c)_y-Q_1-Q_2=T_{\mathsf a\Lambda}((\mathcal W_c)_y-Q_1-Q_2)$, where $\mathscr O_{r,T}((\mathcal W_c)_y-Q_1-Q_2)<\infty$;
by \eqref{e:wcont}, $\mathscr O_{r,T}((\mathcal W_c)_y-Q_1-Q_2)\le\frac14\mathscr O_{r,T}((\mathcal W_c)_y-Q_1-Q_2)$, so
$(\mathcal W_c)_y=Q_1+Q_2$.

\emph{Step 7.} By \eqref{e:wident} and $\pa_y(S^{-1}\mathcal W_c)=S^{-1}(\mathcal W_c)_y+(S^{-1})_y
\mathcal W_c$, $\frac{\mathfrak q_2}{V_c}=d_1+d_2$ with
\[
 d_1=\frac{(Q_D)_2+(S^{-1}Q_1)_2}{V_c},\qquad
 d_2=\frac1{V_c}\bigl[S^{-1}Q_2+(S^{-1})_y\mathcal W_c+\pa_y(S^{-1}\mathcal W_e)\bigr]_2 .
\]
Multiplication by a function $m$ with $|m|+\tau^{\frac12}|m_y|\le C$ changes $\mathscr N$ by at
most a constant factor, and \eqref{e:dpn} gives $\mathscr N(d_1)\le C_{b,r}L_p$. By
\eqref{e:exS}, $\mathscr O((S^{-1})_y\mathcal W_c)\le C\delta^2\iint W_c|\mathcal W_c|^2$; with
\eqref{e:wceng} and Step~2 this gives the bound for $\mathscr O_{r,T}(d_2)$. The fields
$d_1,d_2$ do not depend on $T$.
\end{proof}

\subsection{The decomposition of the pressure}

\begin{lemma}\label{l:prs}
Let $0<\delta\le\delta_1(M)$, let $a_s$ and $b_0$ be the exponents of \eqref{e:pwt}, and let $C_P$ be
the constant of Lemma~\ref{l:hq}, which depends only on the fixed data. There is $r\ge r_0$,
fixed after the initial datum, such that
on $\R\times[r,\infty)$
\begin{equation}\label{e:pspl}
 P_b-p_*=g_1+g_2,\qquad g_1=d_1,\qquad g_2=g_{\rm ref}+d_2+R_p,
\end{equation}
with $g_{\rm ref}$, $R_p$ of \eqref{e:pdecq} and $d_1,d_2$ of Lemma~\ref{l:wmx} for
$p=\mathfrak q_h(r)$; moreover $\mathscr N_{r,\infty}(g_1)<\infty$, and for every $\xi>0$ there
is $C_{b,r,\xi}$ with
\begin{equation}\label{e:prlog}
 \mathscr O_{r,T}(g_2)\le C_{b,r,\xi}+(4C_P\delta^2+\xi)\log\frac{2+T}{2+r}
 \qquad(T\ge r).
\end{equation}
\end{lemma}

\begin{proof}
Choose $r\ge r_0$ so that Lemma~\ref{l:wmx} applies. By Lemma~\ref{l:for},
$p=\mathfrak q_h(r)\in L^1\cap L^2$, and uniqueness in \eqref{e:wh} identifies its solution
with $\mathfrak q_h$. Thus \eqref{e:pq} and \eqref{e:wmx} give \eqref{e:pspl} and
$\mathscr N_{r,\infty}(g_1)<\infty$.

Fix $\xi>0$. By \eqref{e:pq} and Lemma~\ref{l:asy}, choose $S\ge r$ such that
\[
 \tau A_{g_{\rm ref}}\le C_P\delta^2+\frac\xi4,\qquad
 4C_{b,r}L_p^2\varepsilon_b^2\le\frac\xi2\qquad(t\ge S).
\]
where $C_{b,r}$ is the constant in \eqref{e:wmx}.
The inequality $|x+y+z|^2\le2|x|^2+4|y|^2+4|z|^2$, \eqref{e:wmx} and
$\mathscr O_{r_0,\infty}(R_p)<\infty$ yield, for $T\ge S$,
\[
 \mathscr O_{r,T}(g_2)\le C_{b,r,\xi}+(2C_P\delta^2+\xi)\log\frac{1+T}{1+S}.
\]
Here the constant absorbs the integrals on $[r,S]$, which are finite by
\eqref{e:pdecq} and \eqref{e:wmx}. Since
$\log\frac{1+T}{1+S}\le\log2+\log\frac{2+T}{2+r}$, this proves \eqref{e:prlog};
enlarging the same constant covers $r\le T<S$.
\end{proof}

\subsection{Uniformity of the wave-strength threshold}

\begin{proposition}[Uniform wave-strength threshold]\label{p:uniform}
Under the assumptions of Theorem~\ref{thm:planar}, the estimates in
Lemma~\ref{l:mcoef}, Corollary~\ref{l:coef} and Lemma~\ref{l:prs}
hold after decreasing $\delta_0(A)$ using only $A$, the fixed data
and the fixed Gaussian exponents.
The persistent contact and absorption coefficients have the same
dependence; $C_P$ in \eqref{e:prlog} depends only on the
fixed data. Higher derivative bounds, finite source budgets and
entrance times may depend on the complete datum. The Gaussian exponents
$a_s,b_0$ in \eqref{e:pwt} are fixed before the datum, and $C_P$ is
independent of them.
For the main theorem they can be chosen using only the physical data.
\end{proposition}
\begin{proof}
Here $C_A$ may depend on $A$, the fixed data and the fixed Gaussian
exponents; $C_b$ and the entrance times may depend on the complete
datum. Lemma~\ref{l:cmp} and \eqref{e:dat} give
$|m_-|+|m_+|+\delta|h_*|\le C_A$.
In Lemma~\ref{l:mcoef}, the reference bounds
\eqref{e:cmpb} use only this mass bound. Its differentiated
energy estimates use finite $C_b$ bounds and vanishing tails in
\eqref{e:wbas} and \eqref{e:wtC}; they impose no further high-data
wave restriction. Corollary~\ref{l:coef} consequently retains the
fixed-data contact terms in \eqref{e:E0}, \eqref{e:E1} and
\eqref{e:curv}, the bound $C_A+\omega_b(s)$ in \eqref{e:cbg},
and the remaining vanishing tails. Increase the entrance time so
that the solution lies in the fixed positive box and the bounds
in \eqref{e:E6} are at most one.

The pressure estimates have the same separation. In
Lemma~\ref{l:geo}, the conjugations depend only on the acoustic masses;
the persistent connection absorption is $C_A\delta$, and its residual
responses tend to zero as the entrance time increases. The bounds
for $p_j^{\pm1}$, $h_j^{\pm1}$, the response in \eqref{e:dres} and
the number of drift-comparison steps are also controlled by $C_A$.
The reference-field contraction in Lemma~\ref{l:ref} has the same
dependence. In \eqref{e:dfco}, choose $r$ after fixing the datum so that
\[
 \tau_r\ge\max\{h_*^2,\delta^{-2}\},\qquad
 \sup_{t\ge r}\sqrt\tau\Vert R_U(t)\Vert_{L^\infty}\le1.
\]
Then $|\beta_j-jc_-|\le C_A\tau^{-1/2}$. The derivative bounds for
this drift may still be $C_b$: the bounds for $H,H_\eta$ in
Lemma~\ref{l:kernv}(ii) and the leading constant in \eqref{e:dsmb}
are independent of those derivative bounds. The two absorptions in
Lemma~\ref{l:dp} therefore have coefficients
$C_A\delta+C_b\tau_r^{-1/4}$ and
$C_A\delta+C_b\tau_r^{-1/8}$.
The other persistent pressure restrictions in Lemmas~\ref{l:solv},
\ref{l:ga} and~\ref{l:wmx} depend only on the fixed Gaussian exponents
and the fixed data. All finite
source budgets may depend on the complete initial datum. Finally,
Lemma~\ref{l:hq} obtains $C_P$ from the contact Gaussian alone;
the translation disappears in its spatial integral and the acoustic
terms vanish in its limiting estimate. Thus \eqref{e:prlog} has the
same fixed-data $C_P$, uniformly for
$\theta_+/\theta_-\in[1/2,3/2]$.
\end{proof}

\section*{Acknowledgments}
The authors acknowledge the use of generative AI tools in developing and
checking mathematical arguments and in drafting and revising the manuscript.

\begingroup
\small
\urlstyle{same}

\endgroup
\end{document}